\documentclass[12pt]{amsbook}
\usepackage{amssymb}
\usepackage[all]{xy}

\usepackage[colorlinks=true,linkcolor=magenta,citecolor=magenta]{hyperref}
\makeindex

\newtheorem{theorem}{Theorem}[chapter]
\newtheorem{addendum}[theorem]{Addendum}
\newtheorem{answer}[theorem]{Answer}
\newtheorem{cat}[theorem]{Cat}
\newtheorem{claim}[theorem]{Claim}
\newtheorem{comment}[theorem]{Comment}
\newtheorem{conclusion}[theorem]{Conclusion}

\newtheorem{convention}[theorem]{Convention}
\newtheorem{definition}[theorem]{Definition}

\newtheorem{exercise}[theorem]{Exercise}
\newtheorem{fact}[theorem]{Fact}
\newtheorem{method}[theorem]{Method}
\newtheorem{observation}[theorem]{Observation}

\newtheorem{principle}[theorem]{Principle}

\newtheorem{proposition}[theorem]{Proposition}
\newtheorem{question}[theorem]{Question}
\newtheorem{questions}[theorem]{Questions}
\newtheorem{remark}[theorem]{Remark}

\newtheorem{strategy}[theorem]{Strategy}

\begin{document}

\title{Invitation to finite groups}

\author{Teo Banica}
\address{Department of Mathematics, University of Cergy-Pontoise, F-95000 Cergy-Pontoise, France. {\tt teo.banica@gmail.com}}


\begin{abstract}
This is an introduction to the finite groups, with focus on the groups of permutations and reflections, and more generally, on the finite groups of unitary matrices. We first discuss the basics of group theory, featuring the cyclic, dihedral and symmetric groups, and the structure result for finite abelian groups. Then we study the complex reflection groups, with general theory and examples, classification results, and with a look into braid groups too.  We then go into the study of representation theory, and of more advanced aspects, such as Tannakian duality, Brauer theorems and Clebsch-Gordan rules. Finally, we discuss, using representation theory methods, a number of advanced analytic aspects, for the most in relation with questions coming from probability.
\end{abstract}

\maketitle

\chapter*{Preface}

There is a lot of symmetry in the real world, surrounding us. Minerals, plants, animals, we all have interesting symmetry features, witnessing for some built-in symmetry, in the various laws of mathematics, physics, chemistry and biology, having produced us.

\bigskip

Mathematically speaking, understanding this symmetry is a key problem. Have a look for instance at the snowflake pictured below, isn't this beautiful, by all possible beauty standards in this world, and wouldn't you like to know more about its symmetry:
$$
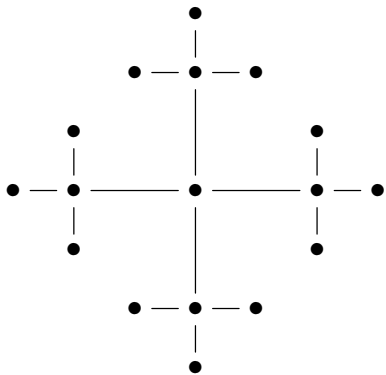
\xymatrix@R=10pt@C=10pt{
&&&\bullet\ar@{-}[d]\\
&&\bullet\ar@{-}[r]&\bullet\ar@{-}[dd]\ar@{-}[r]&\bullet\\
&\bullet\ar@{-}[d]&&&&\bullet\ar@{-}[d]\\
\bullet\ar@{-}[r]&\bullet\ar@{-}[rr]&&\bullet\ar@{-}[rr]\ar@{-}[dd]&&\bullet\ar@{-}[r]&\bullet\\
&\bullet\ar@{-}[u]&&&&\bullet\ar@{-}[u]\\
&&\bullet\ar@{-}[r]&\bullet\ar@{-}[d]\ar@{-}[r]&\bullet\\
&&&\bullet
}$$

Actually, understanding why snowflakes are made like this is a quite difficult question, requiring you to know well all basic mathematics, all basic physics, including quantum mechanics, and then a bit of quantum chemistry, and some advanced thermodynamics too. So, perhaps not the easiest example to start with. Maybe, for later.

\bigskip

More modestly, what we can do, as mathematicians, is to have at least a good understanding of abstract, mathematical symmetry. And here, things are quite straightforward. Symmetries are encoded by mathematical objects called groups, and the simplest such groups are those which are finite. So, as a reasonable objective, we can try to understand the finite groups. And for the laws of nature, and snowflakes, these can come later.

\bigskip

This book is an introduction to the finite groups, with focus on the groups of permutations and reflections, and more generally, on the finite groups of unitary matrices. The text is organized quite symmetrically, in 4 parts, each having 4 chapters, each having 4 sections, plus an informal exercise section at the end, as follows:

\bigskip

Part I discusses the basics of group theory, all standard techniques that you need to know, notably in relation with morphisms and basic operations, and with our main examples here, source of inspiration for everything, being the cyclic groups $\mathbb Z_N$, the dihedral groups $D_N$, and the symmetric groups $S_N$. We also discuss, at the end, more advanced topics, including group actions, Sylow theorems, the structure theorem for finite abelian groups, and with a look into Pontrjagin duality and Fourier transforms too.

\bigskip

Part II goes into more advanced theory and examples, with the aim of discussing the main connections between groups and geometry, and mathematics in general. We first discuss, in some detail, the complex reflection groups $\mathbb Z_s\wr S_N$ and their subgroups. Then we go on a topological excursion, featuring the braid groups $B_k$, and related knot invariants, followed by a geometric excursion in small dimensions, featuring the Platonic solids and their symmetries. Finally, we have a look at the Galois groups.

\bigskip

Part III goes into representation theory, with the fundamental results of Peter-Weyl explained, and then with a discussion of more advanced aspects, including tensor categories, Tannaka duality, Brauer theorems and Clebsch-Gordan rules. We then go on an even more advanced discussion, centered on the notion of easiness, and with this often spilling into the general compact group case. Finally, we discuss quantum extensions of this, notably with the basic theory of the quantum permutation groups $S_N^+$ explained.

\bigskip

Part IV discusses, using representation theory techniques, a number of more advanced analytic aspects, for the most in relation with questions coming from probability. We first discuss the computation of asymptotic character laws for the main groups, featuring the normal, Poisson and Bessel laws. Then we get into some related technical questions, in relation with liberation and free probability, and with Gram determinants as well. We end with an introduction to the subfactor and planar algebra theory of Jones.

\bigskip

In the hope that you will find this book useful. Most of the material is of course very old and standard, but the choice of topics to be discussed, and the order of presentation, was of course a tricky business, and by the way thanks here to my cats, for some help with this. Also, and in the hope that you will appreciate this too, I tried to present everything with simple notations and proofs, even if some notations are not always standard.

\bigskip

\

{\em Cergy, November 2025}

\smallskip

{\em Teo Banica}

\baselineskip=15.95pt
\tableofcontents
\baselineskip=14pt

\part{Finite groups}

\ \vskip50mm

\begin{center}
{\em And I miss you

Like the deserts miss the rain

And I miss you

Like the deserts miss the rain}
\end{center}

\chapter{Group theory}

\section*{1a. Group theory}

The symmetries of nature are encoded by groups, and with groups being something very simple, namely some sets, with a composition operation, which must satisfy what we can expect from a multiplication. The precise definition of groups is as follows:

\index{multiplication}
\index{group}
\index{abelian group}
\index{associativity}
\index{finite group}
\index{group law}

\begin{definition}
A group is a set $G$ endowed with a multiplication operation 
$$(g,h)\to gh$$
which must satisfy the following conditions:
\begin{enumerate}
\item Associativity: we have $(gh)k=g(hk)$, for any $g,h,k\in G$.

\item Unit: there is an element $1\in G$ such that $g1=1g=g$, for any $g\in G$.

\item Inverses: for any $g\in G$ there is $g^{-1}\in G$ such that $gg^{-1}=g^{-1}g=1$.
\end{enumerate}
\end{definition}

As a first observation, this is an abstract definition, making use of some abstract notations, that we can fine-tune and adapt to our situation, if needed. For instance, in order to avoid certain confusions, we can use a dot for the multiplication operation:
$$(g,h)\to g\cdot h$$

This being said, we can use in fact any symbol for the multiplication operation, as for instance $(g,h)\to g\;\!\heartsuit\;\! h$. The same goes for the unit, which can be denoted as $\bullet\in G$ if we want to, and for the inverses as well, which can be denoted $\!\not\! g$, if that is convenient. More on such alternative conventions in a moment, once we will have some examples.

\bigskip

There are many examples of groups, with typically all the basic systems of numbers that we know being groups. Here are some standard illustrations, for this fact:

\index{groups of numbers}

\begin{theorem}
We have the following groups, and non-groups:
\begin{enumerate}
\item $(\mathbb Z,+)$ is a group.

\item $(\mathbb Q,+)$, $(\mathbb R,+)$, $(\mathbb C,+)$ are groups as well.

\item $(\mathbb N,+)$ is not a group.

\item $(\mathbb Q^*,\cdot\,)$ is a group.

\item $(\mathbb R^*,\cdot\,)$, $(\mathbb C^*,\cdot\,)$ are groups as well.

\item $(\mathbb N^*,\cdot\,)$, $(\mathbb Z^*,\cdot\,)$ are not groups.
\end{enumerate}
\end{theorem}

\begin{proof}
All this is clear from the definition of the groups, as follows:

\medskip

(1) The group axioms are indeed satisfied for $(\mathbb Z,+)$, with the group operation being the sum $(g,h)\to g+h$, the unit element being 0, and the inverse map being $g\to-g$. Indeed, the axioms correspond to the following formulae, which are all trivial:
$$(g+h)+k=g+(h+k)$$
$$g+0=0+g=g$$
$$g+(-g)=(-g)+g=0$$

Observe now that $(\mathbb Z,+)$ has the following special property:
$$g+h=h+g$$

However, we have not included this property in our axioms from Definition 1.1, which with the notations there would read $gh=hg$, because we would like our notion of group to cover as well examples which do not satisfy this, $gh\neq hg$. More on these later.

\medskip

(2) Once again, the group axioms are satisfied for $(\mathbb Q,+),(\mathbb R,+),(\mathbb C,+)$, for the same reasons as for $(\mathbb Z,+)$, and with the remark that for $\mathbb Q$ we are using here the fact that the sum of any two rational numbers is a rational number, coming from:
$$\frac{a}{b}+\frac{c}{d}=\frac{ad+bc}{bd}$$

Observe that these latter groups all have the property $g+h=h+g$.

\medskip

(3) In $(\mathbb N,+)$ the first two group axioms are satisfied, for the same reasons as for $(\mathbb Z,+)$. However, we do not have inverses, so we do not have a group:
$$-1\notin\mathbb N$$

As a comment here, such beasts, which satisfy the first two group axioms, but not the third one, are called semigroups. Thus, $(\mathbb N,+)$ is a semigroup, which is not a group.

\medskip

(4) The group axioms are indeed satisfied for $(\mathbb Q^*,\cdot\,)$, with the product $gh$ being the usual product, 1 being the usual 1, and $g^{-1}$ being the usual $g^{-1}$. Observe that we must remove indeed the element $0\in\mathbb Q$, because in a group, any element must be invertible:
$$0\notin\mathbb Q^*$$

Finally, observe that $(\mathbb Q^*,\cdot\,)$ has the commutativity property $gh=hg$.

\medskip

(5) Once again, the group axioms are satisfied for $(\mathbb R^*,\cdot\,),(\mathbb C^*,\cdot\,)$, for the same reasons as for $(\mathbb Q^*,\cdot\,)$, and with the remark that for $\mathbb C^*$ we are using here the fact that the nonzero complex numbers can be inverted, with this coming from the following formula:
$$\frac{1}{a+ib}=\frac{a-ib}{a^2+b^2}$$

Observe that these latter groups both have the property $gh=hg$.

\medskip

(6) In what regards $(\mathbb N^*,\cdot\,),(\mathbb Z^*,\cdot\,)$, here the first two group axioms are satisfied, but not the third one, for instance due to the fact that the element 2 has no inverse:
$$\frac{1}{2}\notin\mathbb Z^*$$

Thus, both $(\mathbb N^*,\cdot\,),(\mathbb Z^*,\cdot\,)$ are semigroups, which are not groups.
\end{proof}

What is next? More examples, of course. In order to construct further examples of groups, from those that we already have, we can use the following simple fact:

\begin{proposition}
Given a group $(G,\cdot\,)$, a subset $H\subset G$ is a subgroup when:
\begin{enumerate}
\item $H$ is stable under multiplication, $g,h\in H\implies gh\in H$.

\item $H$ contains the unit element, $1\in H$.

\item $H$ is stable under taking inverses, $g\in H\implies g^{-1}\in H$.
\end{enumerate}
\end{proposition}

\begin{proof}
This is indeed something self-explanatory, the idea being that, once the above conditions are satisfied, $(H,\cdot\,)$ is a group, with the group axioms being all trivial.
\end{proof}

So, let us see how this technology applies to the various groups from Theorem 1.2. We first have there the group $(\mathbb Z,+)$, and it is easy to construct subgroups, such as:
$$(2\mathbb Z,+)\quad,\quad (3\mathbb Z,+)\quad,\quad (4\mathbb Z,+)\quad,\quad \ldots$$

However, this is something which is not very satisfying, because save for a rescaling by $2,3,4,\ldots\,$, we seem to always obtain the same group, namely $(\mathbb Z,+)$ itself. In order to further comment on this phenomenon, let us introduce the following notion:

\begin{definition}
Two groups are called isomorphic, and we write $G\simeq H$, or even $G=H$, when we have a bijection as follows, called group isomorphism,
$$f:G\simeq H$$
which transforms the group operations of $G$ into those of $H$.
\end{definition}

In other words, a group isomorphism must satisfy the following conditions:
$$f(gh)=f(g)f(h)$$
$$f(1)=1$$
$$f(g^{-1})=f(g)^{-1}$$

Now with this notion in hand, let us go back to the phenomenon discovered above, regarding the subgroups of $(\mathbb Z,+)$. We have the following result, regarding these:

\begin{theorem}
The subgroups of $(\mathbb Z,+)$ are as follows:
\begin{enumerate}
\item These are exactly the groups $(n\mathbb Z,+)$, with $n\in\mathbb N$.

\item And, save for $\{0\}$, any such group is isomorphic to $(\mathbb Z,+)$ itself.
\end{enumerate}
\end{theorem}

\begin{proof}
This follows indeed from some basic arithmetic, the idea being as follows:

\medskip

(1) To start with, given an integer $n\in\mathbb N$, the subset $n\mathbb Z\subset\mathbb Z$ is indeed a subgroup, in the sense of Proposition 1.3, as shown by the following formulae:
$$na+nb=n(a+b)$$
$$0=n0$$
$$-na=n(-a)$$

Conversely now, consider a subgroup $G\subset\mathbb Z$. When $G$ has only one element, this element must be the unit $0$, and our group follows to be covered by our construction:
$$G=\{0\}=0\mathbb Z$$

In the case $|G|>1$ now, consider the smallest $n>0$ belonging to $G$. By using the addition and inverse operations, we obtain successively:
$$n\in G\implies\pm n\in G\implies\pm 2n\in G\implies\pm3n\in G\implies\ldots$$

Since we must have $0\in G$ too, we conclude that we have an inclusion as follows:
$$n\mathbb Z\subset G$$

Our claim is that this inclusion must be an equality. Indeed, assume that this does not hold, and pick $m\in G-n\mathbb Z$, that we can assume by inverting to be positive, $m>0$. Now let us divide $m$ by $n$. This produces a formula as follows, with $0<r<n$:
$$m=cn+r$$

But from $m\in G$ and $n\mathbb Z\subset G$ we obtain $r\in G$, and together with $0<r<n$ this contradicts our choice above of $n$, as being the smallest positive integer belonging to $G$. Thus, contradiction, and we conclude that we have $G=n\mathbb Z$, as desired.

\medskip

(2) This is clear, because given $n\in\mathbb N$, we can consider the following map:
$$f:\mathbb Z\to n\mathbb Z\quad,\quad a\to na$$

This map is then bijective, assuming $n>0$ as in the statement, and the three group isomorphism conditions given after Definition 1.4 are then satisfied, as follows:
$$f(a+b)=na+nb=f(a)+f(b)$$
$$f(1_{\mathbb Z})=n=1_{n\mathbb Z}$$
$$f(-a)=-na=-f(a)$$

Thus, our map $f$ above is a group isomorphism $\mathbb Z\simeq n\mathbb Z$, as desired.
\end{proof}

Very nice all this, we have now a Theorem that we can be proud of, but in what regards the question that we started with, namely finding interesting new groups, coming as subgroups of the groups in Theorem 1.2, what we did so far is rather a failure.

\bigskip

Nevermind. Getting to item (2) in Theorem 1.2, we have there the groups $(\mathbb Q,+)$, $(\mathbb R,+)$, $(\mathbb C,+)$, and these groups do have interesting subgroups, as shown by:

\begin{proposition}
The following is a subgroup of $(\mathbb R,+)$,
$$G=\left\{a+b\sqrt{2}\,\Big|a,b\in\mathbb Z\right\}$$
which is not isomorphic to $(\mathbb Z,+)$, $(\mathbb Q,+)$, $(\mathbb R,+)$ or $(\mathbb C,+)$.
\end{proposition}

\begin{proof}
Here the fact that we have indeed a subgroup is clear from definitions, and in what regards the fact that this group is indeed something new, this is quite routine, that we will leave here as an exercise, the idea with this being as follows:

\medskip

(1) With $\sqrt{2}$ being irrational, it is clear that $G$ looks a bit like $\mathbb Q$, that is, dense and messy inside $\mathbb R$, which makes any isomorphism with $\mathbb Z$ hopeless. In practice, proving that we have $G\not\simeq\mathbb Z$ is not difficult, by contradiction, and we will leave this as an exercise.

\medskip

(2) Getting now to $\mathbb Q$, this is a more serious contender for an isomorphism, in view of how $G$ looks like. But here, again, proving that we have $G\not\simeq\mathbb Q$ is not very difficult, by contradiction, and we will leave this again as an instructive exercise.

\medskip

(3) In what regards now $\mathbb R$, this is excluded right away, by the fact that $G$ is countable, while $\mathbb R$ is uncountable. And the same goes for $\mathbb C$, which is uncountable too.
\end{proof}

Quite nice all this, and it looks that with the subgroups $G\subset\mathbb Q$, we have found a goldmine for groups, with other available such beasts including:
$$\left\{a\sqrt{2}+b\sqrt{3}\,\Big|a,b\in\mathbb Z\right\}\quad,\quad \left\{a\sqrt{2}+b\sqrt{3}+c\sqrt{7}\,\Big|a,b,c\in\mathbb Z\right\}\quad,\quad \ldots$$

However, this is not exactly true, as we will soon discover. Indeed, getting back now to the group theory generalities, we have the following result, which can help us:

\begin{theorem}
Given two groups $G,H$, their product $G\times H$, with operations
$$(g,h)(g',h')=(gg',hh')$$
$$1=(1,1)$$
$$(g,h)^{-1}=(g^{-1},h^{-1})$$
is a group too.
\end{theorem}

\begin{proof}
This is indeed something self-explanatory, with the verification of the group axioms being something trivial, that we will leave as an exercise. In fact, we have called this statement Theorem instead of Proposition only due to its theoretical importance.
\end{proof}

Now with the above product notion in hand, we can go back to the mysterious group constructed in Proposition 1.6, and reveal its quite trivial nature, as follows:

\begin{proposition}
We have an isomorphism of groups, as follows,
$$\left\{a+b\sqrt{2}\,\Big|a,b\in\mathbb Z\right\}\simeq\mathbb Z^2$$
given by $a+b\sqrt{2}\to(a,b)$.
\end{proposition}

\begin{proof}
It is better here to construct instead the isomorphism in the other sense. So, consider the following map, going from right to left:
$$(a,b)\to a+b\sqrt{2}$$

This map is then sujective, and injective too, because $\sqrt{2}$ is irrational. Thus we have a bijection, and this bijection is the group isomorphism we are looking for.
\end{proof}

Damn, I would say, in relation with our initial objective, namely finding interesting subgroups of the groups in Theorem 1.2. And even worse is true, because with some help with advanced arithmetic, we can trivialize other such beasts too, as for instance:
$$\left\{a+b\sqrt{2}+c\sqrt{3}+d\pi+e\sqrt[7]{5}\,\Big|a,b,c,d,e\in\mathbb Z\right\}\simeq\mathbb Z^5$$

Nevermind. Getting back now to Theorem 1.2, we still have there available the groups $(\mathbb Q^*,\cdot\,)$, $(\mathbb R^*,\cdot\,)$, $(\mathbb C^*,\cdot\,)$, whose subgroups we can investigate. And leaving here aside $(\mathbb Q^*,\cdot\,)$ and $(\mathbb R^*,\cdot\,)$, which are smaller, and focusing instead of $(\mathbb C^*,\cdot\,)$, we are led to:

\begin{theorem}
The unit circle in the complex plane
$$\mathbb T=\left\{z\in\mathbb C\Big||z|=1\right\}$$
is a subgroup of $(\mathbb C^*,\cdot\,)$, not isomorphic as continuous group to our previous groups.
\end{theorem}

\begin{proof}
Many things can be said here, the idea being as follows:

\medskip

(1) To start with, the fact that $\mathbb T$ is indeed a multiplicative group is clear, coming from the formula $|zz'|=|z|\cdot|z'|$, and with the inverses being given by $z^{-1}=\bar{z}$.

\medskip

(2) Next, you might wonder why we chose that $\mathbb T$ symbol for our group. In answer, if we look at the various powers of $\mathbb T$, these are all tori, as follows:
$$\mathbb T^2\subset\mathbb C^2\simeq\mathbb R^4$$
$$\mathbb T^3\subset\mathbb C^3\simeq\mathbb R^6$$
$$\ldots$$

Thus, $\mathbb T$ itself is best thought of as being a torus, with $\mathbb T$ standing for this.

\medskip

(3) Getting now to the non-isomorphism with our previous groups, as a first observation, we have $(-1)^2=1$ inside $\mathbb T$, which in additive notation would translate into $2(-1)=0$, and this rules out the isomorphism with $(\mathbb Z,+)$, $(\mathbb Q,+)$, $(\mathbb R,+)$, $(\mathbb C,+)$.

\medskip

(4) Our group $\mathbb T$ being uncountable, it is not isomorphic to $(\mathbb Q^*,\cdot\,)$ either.

\medskip

(5) With $w=e^{2\pi i/3}$ we have $w^3=1$, which rules out the isomorphism with $(\mathbb R^*,\cdot\,)$.

\medskip

(6) Regarding now the non-isomorphism $\mathbb C^*\not\simeq\mathbb T$, this is something more tricky, because there is no obvious algebraic reason for this. As a solution, we certainly have no isomorphism of continuous groups $\mathbb C^*\simeq\mathbb T$, because $\mathbb T$ is compact, while $\mathbb C^*$ is not.

\medskip

(7) As a related comment here, observe that we have $\mathbb C^*\simeq\mathbb R_+^*\times\mathbb T$, with this coming from $z=re^{it}$, and with this leading us into the study of the subgroup $\mathbb R_+^*\subset(\mathbb R^*,\cdot\,)$, which looks new too. But this latter group is in fact not new, because we have a group isomorphism $(\mathbb R,+)\simeq(\mathbb R_+^*,\cdot\,)$ given by the usual exponential function, $t\to e^t$.

\medskip

(8) Finally, talking products, $\mathbb T$ cannot be isomorphic either, as a continuous group, to any product of the groups constructed before, due to its compacity, as in (6).
\end{proof}

Along the same lines,  subgroups of $\mathbb C^*$, we have as well the following result:

\begin{theorem}
The set of $N$-th roots of unity in the complex plane,
$$\mathbb Z_N=\left\{1,w,\ldots,w^{N-1}\Big| w=e^{2\pi i/N}\right\}$$
is a subgroup of $\mathbb T$, called cyclic group. This group is isomorphic via $w^a\to a$ to
$$\mathbb Z_N=\{0,1,\ldots,N-1\}$$
the group of remainders modulo $N$, with the usual addition of these remainders.
\end{theorem}

\begin{proof}
There are several things going on here, the idea being that we have two groups, which are isomorphic, and in view of this isomorphism, we can use the same notation for both, namely $\mathbb Z_N$. In practice now, the proof goes as follows:

\medskip

(1) In what regards the first $\mathbb Z_N$, this is indeed a subgroup of $\mathbb T$, with its multiplication, unit and inverse operations being as follows, all exponents being modulo $N$:
$$w^aw^b=w^{a+b}$$
$$1=w^0$$
$$(w^a)^{-1}=w^{-a}$$

(2) In what regards the second $\mathbb Z_N$, this is also a group, with its multiplication, unit and inverse operations being as follows, with all numbers being taken modulo $N$:
$$(a,b)\to a+b$$
$${\rm unit}=0$$
$$a\to-a$$

(3) Finally, the two groups that we have are indeed isomorphic, with the isomorphism in one sense being given by $a\to w^a$, and in the other sense, by $w^a\to a$, as stated.
\end{proof}

To summarize now, by the end of the day, we have as main examples of groups those in Theorem 1.2, Theorem 1.9 and Theorem 1.10, and their products.

\bigskip

However, the story is not over here, because in what regards the products, interesting things can happen. As an example, here are an isomorphism and a non-isomorphism, which are quite easy to establish, and that we will leave as an instructive exercise:
$$\mathbb Z_2\times\mathbb Z_3\simeq\mathbb Z_6\quad,\quad\mathbb Z_2\times\mathbb Z_2\not\simeq\mathbb Z_4$$

We will be back to such things later in this book, first in chapter 2 when discussing more in detail the cyclic groups, and then with a ultimate discussion on the subject in chapter 4, where we will investigate, in general, the groups satisfying $gh=hg$.

\section*{1b. Matrix groups} 

Getting now to more complicated examples of groups, we have many of them coming from space geometry and linear algebra, as shown by the following result:

\index{groups of matrices}
\index{special linear group}
\index{general linear group}
\index{orthogonal group}
\index{unitary group}
\index{special orthogonal group}
\index{special unitary group}

\begin{theorem}
We have the following groups:
\begin{enumerate}
\item $(\mathbb R^N,+)$ and $(\mathbb C^N,+)$.

\item $(M_N(\mathbb R),+)$ and $(M_N(\mathbb C),+)$.

\item $(GL_N(\mathbb R),\cdot\,)$ and $(GL_N(\mathbb C),\cdot\,)$, the invertible matrices.

\item $(SL_N(\mathbb R),\cdot\,)$ and $(SL_N(\mathbb C),\cdot\,)$, with S standing for ``special'', meaning $\det=1$.

\item $(O_N,\cdot\,)$ and $(U_N,\cdot\,)$, the orthogonal and unitary matrices.

\item $(SO_N,\cdot\,)$ and $(SU_N,\cdot\,)$, with S standing as above for $\det=1$.
\end{enumerate}
\end{theorem}

\begin{proof}
All this is clear from definitions, and from various standard results from linear algebra, that you certainly know well, the idea with all this being as follows:

\medskip

(1) The axioms are indeed clearly satisfied for $\mathbb R^N,\mathbb C^N$, with the sum being the usual sum of vectors, $-v$ being the usual $-v$, and the null vector $0$ being the unit. 

\medskip

(2) Once again, the axioms are clearly satisfied for $M_N(\mathbb R),M_N(\mathbb C)$, with the sum being the usual sum of matrices, $-A$ being the usual $-A$, and the null matrix $0$ being the unit. Observe that what we have here is in fact a particular case of (1), because any $N\times N$ matrix can be regarded as a $N^2\times1$ vector, and so at the group level we have:
$$(M_N(\mathbb R),+)\simeq(\mathbb R^{N^2},+)$$ 
$$(M_N(\mathbb C),+)\simeq(\mathbb C^{N^2},+)$$

(3) Regarding now $GL_N(\mathbb R),GL_N(\mathbb C)$, these are groups because the product of invertible matrices is invertible, according to the following formula:
$$(AB)^{-1}=B^{-1}A^{-1}$$

Observe that at $N=1$ we obtain the groups $(\mathbb R^*,\cdot),(\mathbb C^*,\cdot)$.

\medskip

(4) The sets $SL_N(\mathbb R),SL_N(\mathbb C)$ formed by the real and complex matrices of determinant 1 are subgroups of the groups in (3), because of the following formula, which shows that the matrices satisfying $\det A=1$ are stable under multiplication:
$$\det(AB)=\det(A)\det(B)$$

(5) Regarding now $O_N,U_N$, here the group property is clear too from definitions, and is best seen by using the associated linear maps, because the composition of two isometries is an isometry. Equivalently, assuming $U^*=U^{-1}$ and $V^*=V^{-1}$, we have:
\begin{eqnarray*}
(UV)^*
&=&V^*U^*\\
&=&V^{-1}U^{-1}\\
&=&(UV)^{-1}
\end{eqnarray*}

Observe that at $N=1$ we obtain the groups $\mathbb Z_2$, $\mathbb T$.

\medskip

(6) The sets of matrices $SO_N,SU_N$ in the statement are obtained by intersecting the groups in (4) and (5), and so they are groups indeed:
$$SO_N=O_N\cap SL_N(\mathbb R)\quad,\quad 
SU_N=U_N\cap SL_N(\mathbb C)$$

Thus, all the sets in the statement are indeed groups, as claimed.
\end{proof}

Getting back now to the general case, that of the arbitrary groups, can we have some theory going on, in the context of Definition 1.1? And with this question in mind, a look at the various examples of groups that we have, from Theorems 1.2, 1.9, 1.10 and 1.11, suggests having as starting point for our considerations the following fact:

\begin{fact}
The groups can be divided into classes, as follows:
\begin{enumerate}
\item Groups satisfying $gh=hg$, and groups not satisfying this.

\item Groups which are finite, and groups which are infinite.
\end{enumerate}
\end{fact}

In practice now, let us first have a look at the condition $gh=hg$, called commutativity of the multiplication, which seems to best distinguish between the simple and the complicated groups. In relation with it, we have the following convention:

\begin{convention}
When the group multiplication is commutative we call $G$ abelian, en hommage to Abel, and we usually denote its multiplication additively,
$$(g+h)\to g+h$$
with the group axioms taking the following form:
\begin{enumerate}
\item Associativity: we have $(g+h)+k=g+(h+k)$, for any $g,h,k\in G$.

\item Unit: there is an element $0\in G$ such that $g+0=0+g=g$, for any $g\in G$.

\item Inverses: for any $g\in G$ there is $-g\in G$ such that $g+(-g)=(-g)+g=0$.
\end{enumerate}
\end{convention}

However, this will be not a general rule, and rather the converse is true, in the sense that if a group is denoted as above, this means that the group must be abelian. 

\bigskip

At the level of the basic examples of groups that we have, from Theorem 1.2 and Theorem 1.11, somehow confirming our $\cdot$ vs $+$ rule, the situation is as follows:

\begin{theorem}
The following happen:
\begin{enumerate}
\item $(\mathbb Z,+)$, $(\mathbb Q,+)$, $(\mathbb R,+)$, $(\mathbb C,+)$ are abelian.

\item $(\mathbb Q^*,\cdot)$, $(\mathbb R^*,\cdot)$, $(\mathbb C^*,\cdot)$ are abelian too.

\item $(\mathbb R^N,+)$, $(\mathbb C^N,+)$ and $(M_N(\mathbb R),+)$, $(M_N(\mathbb C),+)$ are abelian too.

\item $(GL_N(\mathbb R),\cdot)$, $(GL_N(\mathbb C),\cdot)$, $(SL_N(\mathbb R),\cdot)$, $(SL_N(\mathbb C),\cdot)$ with $N\geq2$ are not abelian.

\item $(SO_2,\cdot)$ is abelian, while $(SO_N,\cdot)$ with $N\geq3$ is not abelian.

\item $(O_N,\cdot)$, $(U_N,\cdot)$, $(SU_N,\cdot)$ with $N\geq2$ are not abelian. 
\end{enumerate}
\end{theorem}

\begin{proof}
There are 20 groups in this statement, and as a matter of cutting from complexity, and doing all this quickly, we can proceed as follows:

\medskip

(1) We have inclusions of sets as follows, which are additive group inclusions:
$$\mathbb Z\subset\mathbb Q\subset\mathbb R\subset\mathbb C$$

Now since the group on the right is abelian, due to $x+y=y+x$ for any $x,y\in\mathbb C$, all its subgroups are abelian too. Thus, all the above groups are indeed abelian.

\medskip

(2) This is similar, using the following inclusions, and $xy=yx$, for any $x,y\in\mathbb C$:
$$\mathbb Q^*\subset\mathbb R^*\subset\mathbb C^*$$

(3) In what regards the first two groups, the proof is similar, by using the following inclusion, and the fact that we have $x+y=y+x$ for any two vectors $x,y\in\mathbb C^N$:
$$\mathbb R^N\subset\mathbb C^N$$

As for the last two groups, abelianity follows from the following isomorphisms: 
$$M_N(\mathbb R)\simeq\mathbb R^{N^2}\quad,\quad 
M_N(\mathbb C)\simeq\mathbb C^{N^2}$$

(4) We have inclusions as follows, between the groups in the statement:
$$\xymatrix@R=55pt@C=40pt{
GL_N(\mathbb R)\ar[r]&GL_N(\mathbb C)\\
SL_N(\mathbb R)\ar[r]\ar[u]&SL_N(\mathbb C)\ar[u]}$$

As a first observation, at $N=1$ all these groups are abelian, by using the same argument as before, along with the following isomorphism, and (2):
$$GL_1(\mathbb C)=\mathbb C^*$$

In order to prove now the non-abelianity at $N\geq2$, again by using the same argument, it is enough to prove that the smallest group, namely $SL_N(\mathbb R)$, is not abelian. But here, observe that we have inclusions as follows, obtained by adding a 1 entry, diagonally:
$$SL_2(\mathbb R)\subset SL_3(\mathbb R)\subset SL_4(\mathbb R)\subset\ldots$$ 

Thus, it is enough to prove that the group $SL_2(\mathbb R)$ is not abelian, and this can be done as follows, by using a basic Jordan block, and its transpose:
$$\begin{pmatrix}1&1\\0&1\end{pmatrix}
\begin{pmatrix}1&0\\1&1\end{pmatrix}
=\begin{pmatrix}2&1\\1&1\end{pmatrix}
\neq
\begin{pmatrix}1&1\\1&2\end{pmatrix}
=\begin{pmatrix}1&0\\1&1\end{pmatrix}
\begin{pmatrix}1&1\\0&1\end{pmatrix}$$

(5,6) The temptation is high here to use the same method as in (4), based on the following diagram of inclusions, between our various groups of isometries:
$$\xymatrix@R=60pt@C=55pt{
O_N\ar[r]&U_N\\
SO_N\ar[r]\ar[u]&SU_N\ar[u]}$$

At $N=1$ all these groups are abelian, by using the same argument as before, along with the following isomorphism, with $\mathbb T\subset\mathbb C^*$ being the unit circle in the plane:
$$U_1=\mathbb T$$

At $N\geq 2$ now, you would say that it is enough to prove that $SO_2$ is not abelian, which would basically terminate the proof, for all groups in question. However, this is wrong, because $SO_2$ is the group of rotations in the plane, which do commute with each other. Thus, we have some work to do here, and this can be done as follows:

\medskip

(7) To start with, as just mentioned, $SO_2$ is the group of rotations in the plane, which do commute with each other. In case you are not very familiar with this, let us pick $U=\binom{a\ b}{c\ d}$ having determinant one, $ad-bc=1$. The condition $U^t=U^{-1}$ reads:
$$\begin{pmatrix}a&c\\b&d\end{pmatrix}
=\begin{pmatrix}d&-b\\-c&a\end{pmatrix}$$

Thus we must have $U=\binom{a\ -c}{c\ \ a}$. Moreover, the determinant one condition becomes $a^2+c^2=1$, so $a=\cos t$, $c=\sin t$, and we conclude that $U$ is the rotation of angle $t$:
$$U=\begin{pmatrix}\cos t&-\sin t\\ \sin t&\cos t\end{pmatrix}$$

Summarizing, claim proved, and with these rotations commuting, $SO_2$ is abelian.

\medskip

(8) Next, let us study $O_2$. Although we will not formally need this now, it is worth working out the structure of $O_2$ first. Since $V^t=V^{-1}$ implies $\det V=\pm1$, we have two cases, with the case $\det V=1$ being studied above, and with the case $\det V=-1$ leading, via a similar computation, to the symmetry with respect to $Ox$ rotated by $t/2$:
$$V=\begin{pmatrix}\cos t&\sin t\\ \sin t&-\cos t\end{pmatrix}$$

Getting now to what we wanted to do here, namely non-abelianity of $O_2$, this is best seen by using a basic symmetry, and a basic rotation, as follows:
$$\begin{pmatrix}0&1\\1&0\end{pmatrix}
\begin{pmatrix}0&1\\-1&0\end{pmatrix}
=\begin{pmatrix}-1&0\\0&1\end{pmatrix}
\neq
\begin{pmatrix}1&0\\0&-1\end{pmatrix}
=\begin{pmatrix}0&1\\-1&0\end{pmatrix}
\begin{pmatrix}0&1\\1&0\end{pmatrix}$$

Now since we have inclusions $O_2\subset O_3\subset\ldots\,$, all these groups are non-abelian.

\medskip

(9) Let us study now $SU_2$. For this purpose, let us pick $U=\binom{a\ b}{c\ d}$ having determinant one, $ad-bc=1$. The complex isometry condition $U^*=U^{-1}$ reads:
$$\begin{pmatrix}\bar{a}&\bar{c}\\\bar{b}&\bar{d}\end{pmatrix}
=\begin{pmatrix}d&-b\\-c&a\end{pmatrix}$$

We conclude, a bit as in (7), that the solutions are as follows, with $|a|^2+|b|^2=1$:
$$U=\begin{pmatrix}a& b\\-\bar{b}&\bar{a}\end{pmatrix}$$

What is next? A key observation of Pauli, namely that $SU_2$ is spanned by:
$$c_1=\begin{pmatrix}1&0\\ 0&1\end{pmatrix}\quad,\quad
c_2=\begin{pmatrix}i&0\\ 0&-i\end{pmatrix}\quad,\quad
c_3=\begin{pmatrix}0&1\\ -1&0\end{pmatrix}\quad,\quad 
c_4=\begin{pmatrix}0&i\\ i&0\end{pmatrix}$$

And the point is that all this physics trickery helps, because we get, right away:
$$\begin{pmatrix}i&0\\0&-i\end{pmatrix}
\begin{pmatrix}0&1\\-1&0\end{pmatrix}
=\begin{pmatrix}0&i\\i&0\end{pmatrix}
\neq
\begin{pmatrix}0&-i\\-i&0\end{pmatrix}
=\begin{pmatrix}0&1\\-1&0\end{pmatrix}
\begin{pmatrix}i&0\\0&-i\end{pmatrix}$$

Now since we have inclusions $SU_2\subset SU_3\subset\ldots\,$, all the groups are non-abelian. As for the groups $U_2\subset U_3\subset\ldots\,$, which all contain $SU_2$, these are non-abelian either.

\medskip

(10) Still with me I hope, after all these computations, and coming as good news, we only have the groups $SO_3\subset SO_4\subset\ldots\,$ left. And here, if we manage to prove that $SO_3$ is not abelian, done and we can go get a beer. So, let us do this quickly, as follows:
\begin{eqnarray*}
\begin{pmatrix}0&1&0\\1&0&0\\0&0&-1\end{pmatrix}
\begin{pmatrix}0&1&0\\-1&0&0\\0&0&1\end{pmatrix}
&=&\begin{pmatrix}-1&0&0\\0&1&0\\0&0&-1\end{pmatrix}\\
&\neq&\begin{pmatrix}1&0&0\\0&-1&0\\0&0&-1\end{pmatrix}\\
&=&\begin{pmatrix}0&1&0\\-1&0&0\\0&0&1\end{pmatrix}
\begin{pmatrix}0&1&0\\1&0&0\\0&0&-1\end{pmatrix}
\end{eqnarray*}

Of course, you might wonder where the above matrices come from. In answer, I used those in (8) via a trick, namely the fact that we have an embedding $O_2\subset SO_3$.
\end{proof}

Thinking a bit at this, with a bit of all-around algebraic know-how, that you will certainly master too, as you read this book, Theorem 1.14 was in fact something quite easy. In fact, thinking some more, here is a new version of Theorem 1.14, and of pretty much everything that we learned so far in this book, coming with a compact proof:

\begin{theorem}
We have groups of invertible matrices as follows,
$$\xymatrix@R=22pt@C=5pt{
&GL_N(\mathbb R)\ar[rr]&&GL_N(\mathbb C)\\
O_N\ar[rr]\ar[ur]&&U_N\ar[ur]\\
&SL_N(\mathbb R)\ar[rr]\ar[uu]&&SL_N(\mathbb C)\ar[uu]\\
SO_N\ar[uu]\ar[ur]\ar[rr]&&SU_N\ar[uu]\ar[ur]
}$$
the abelian ones being those at $N=1$, plus $SO_2$. 
\end{theorem}

\begin{proof}
This is something that we already know, from the above, but which is worth a new discussion, at a more advanced level, as follows:

\medskip

(1) To start with, we have indeed groups as in the statement, with this coming from basic results from linear algebra, as explained in the proof of Theorem 1.11.

\medskip

(2) At $N=1$, our diagram is as follows, consisting of abelian subgroups of $\mathbb C^*$:
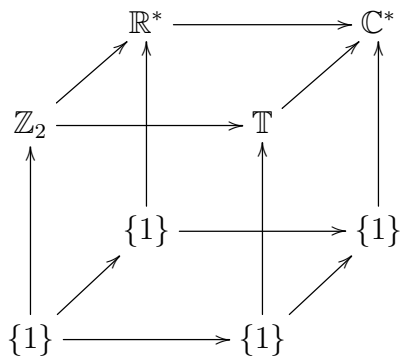
$$\xymatrix@R=22pt@C=19pt{
&\mathbb R^*\ar[rr]&&\mathbb C^*\\
\mathbb Z_2\ar[rr]\ar[ur]&&\mathbb T\ar[ur]\\
&\{1\}\ar[rr]\ar[uu]&&\{1\}\ar[uu]\\
\{1\}\ar[uu]\ar[ur]\ar[rr]&&\{1\}\ar[uu]\ar[ur]
}$$

(3) Still talking $N=1$, observe that we have in the above diagram all the multiplicative number groups that we know, except for $\mathbb Q^*$. It is possible to further improve this, by talking about $GL_N(F)$ and $SL_N(F)$, over an arbitrary field $F$, but in what regards the corresponding orthogonal and unitary groups, things here are more complicated. So, we will not do so, and declare ourselves satisfied with the examples that we have.

\medskip

(4) At $N=2$ now, the diagram in the statement is as follows:
$$
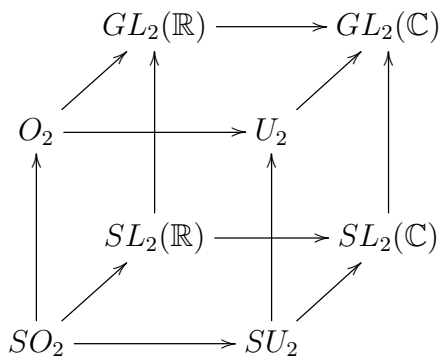
\xymatrix@R=22pt@C=7pt{
&GL_2(\mathbb R)\ar[rr]&&GL_2(\mathbb C)\\
O_2\ar[rr]\ar[ur]&&U_2\ar[ur]\\
&SL_2(\mathbb R)\ar[rr]\ar[uu]&&SL_2(\mathbb C)\ar[uu]\\
SO_2\ar[uu]\ar[ur]\ar[rr]&&SU_2\ar[uu]\ar[ur]
}$$

And here, assuming a bit of linear algebra know-how, $SO_2$ is abelian, as being the group of rotations in the plane, then $SL_2(\mathbb R),O_2,SU_2$ are not abelian, due respectively to the  presence of Jordan blocks, symmetries and rotations, and Pauli matrices, and finally the remaining groups, which are bigger, follow to be not abelian either.

\medskip

(5) Finally, at $N\geq3$ all the groups under consideration contain $SO_3$, which in turn contains $O_2$, which is not abelian, so these groups are not abelian.
\end{proof}

Which sounds quite nice and clean, good and valuable 1 page of mathematics that we learned, and we will stop with our simplifications here.

\section*{1c. General results}

What we did in the above with matrices, which correspond to the linear transformations of $\mathbb R^N$ or $\mathbb C^N$, can be partly extended to more general transformations, not necessarily linear. In order to discuss this, let us start with the following basic fact:

\begin{proposition}
Given any invertible transformation $T$ of the plane, or of the space, or of anything else, the set of transformations
$$G=\left\{T^k\Big|k\in\mathbb Z\right\}$$
is a group, which is abelian, with the usual composition of transformations.
\end{proposition}

\begin{proof}
This is indeed clear from the definition of the groups, as follows:

\medskip

(1) We have $T^kT^l=T^{k+l}$, which implies associativity, $(T^kT^l)T^m=T^k(T^lT^m)$.

\medskip

(2) Our group has indeed a unit element, which is the transformation $T^0=id$.

\medskip

(3) We have inverses as well, given by $(T^k)^{-1}=T^{-k}$, for any $k\in\mathbb Z$.

\medskip

Finally, the fact that our group is indeed abelian is clear too, from $T^kT^l=T^{k+l}$.
\end{proof}

We can further complicate the above construction by adding to the picture a second transformation $S$. We are led in this way to the following statement:

\begin{proposition}
Given two invertible transformations $S,T$ the set
$$G=\left\{\ldots S^{k_i}T^{l_i}S^{k_{i+1}}T^{l_{i+1}}\ldots\Big|k_i,l_i\in\mathbb Z\right\}$$
is a group, with this being the smallest group containing $S,T$. Moreover:
\begin{enumerate}
\item When $S=T$, we recover the group $G=\{T^k|k\in\mathbb Z\}$ constructed before.

\item In fact, when $S=T^l$, we recover the previous group $G=\{T^k|k\in\mathbb Z\}$.

\item When $ST=TS$ this group is abelian, and given by $G=\{S^kT^l|k,l\in\mathbb Z\}$.
\end{enumerate}
\end{proposition}

\begin{proof}
This is indeed something self-explanatory and elementary, with the group property coming as in the proof of Proposition 1.16, and with (1,2,3) being clear.
\end{proof}

More generally now, we have the following result:

\begin{theorem}
We can talk, in the obvious way, by considering all possible compositions, about the group generated by $n$ invertible transformations $T_1,\ldots,T_n$,
$$G=\Big<T_1,\ldots,T_n\Big>$$
and this group is abelian at $n=1$. However, in the general case, unless we have 
$$T_iT_j=T_iT_j$$
for any $i,j$, this group will not be abelian.
\end{theorem}

\begin{proof}
Again, this is something self-explanatory and elementary, with the group property coming as in the proof of Proposition 1.16, and with the rest being clear.
\end{proof}

With this discussed, and thinking some more, we do not actually need transformations $T_1,\ldots,T_n$ in order to construct a group, because we can use abstract variables as well. Many things can be said here, notably with the following result:

\begin{theorem}
We can talk about groups defined with generators and relations,
$$G=\left<g_1,\ldots,g_n\Big|\mathcal R\right>$$
and in particular about free groups, with generators having no relations between them:
$$F_n=\left<g_1,\ldots,g_n\Big|\emptyset\right>$$
Moreover, we have models for the free groups, coming from topology.
\end{theorem}

\begin{proof}
Again, this is something elementary, save for the last part, about the models, which requires some topology know-how. Observe also that we have $F_1=\mathbb Z$. More generally, the free abelian groups, which are by definition the quotients of the free groups $F_n$ by the commutation relations $gh=hg$ between generators, are the groups $\mathbb Z^n$.
\end{proof}

Getting now to the modeling part, let us start with the following notion:

\begin{definition}
A topological space $X$ is called path connected when any two points $x,y\in X$ can be connected by a path. That is, given any two points $x,y\in X$, we can find a continuous function $f:[0,1]\to X$ such that $f(0)=x$ and $f(1)=y$.
\end{definition}

The problem is now, given a connected space $X$, how to count its ``holes''. And this is quite subtle problem, because as examples of such spaces we have:

\medskip

(1) The sphere, the donut, the double-holed donut, the triple-holed donut, and so on. These spaces are quite simple, and intuition suggests to declare that the number of holes of the $N$-holed donut is, and you guessed it right, $N$.

\medskip

(2) However, we have as well as example the empty sphere, I mean just the crust of the sphere, and while this obviously falls into the class of ``one-holed spaces'', this is not the same thing as a donut, its hole being of different nature.

\medskip

(3) As another example, consider again the sphere, but this time with two tunnels drilled into it, in the shape of a cross. Whether that missing cross should account for 1 hole, or for 2 holes, or for something in between, I will leave it up to you. 

\medskip

Summarizing, things are quite tricky, suggesting that the ``number of holes'' of a topological space $X$ is not an actual number, but rather something more complicated. 

\bigskip

Now with this in mind, let us formulate the following definition:

\index{homotopy group}
\index{hole}
\index{loop}
\index{topological space}
\index{null loop}

\begin{definition}
The homotopy group $\pi_1(X)$ of a connected space $X$ is the group of loops based at a given point $*\in X$, with the following conventions,
\begin{enumerate}
\item Two such loops are identified when one can pass continuously from one loop to the other, via a family of loops indexed by $t\in[0,1]$,

\item The composition of two such loops is the obvious one, namely is the loop obtained by following the first loop, then the second loop, 

\item The unit loop is the null loop at $*$, which stays there, and the inverse of a given loop is the loop itself, followed backwards,
\end{enumerate}
with the remark that the group $\pi_1(X)$ defined in this way does not depend on the choice of the given point $*\in X$, where the loops are based.
\end{definition}

This definition is obviously something non-trivial, based on some preliminary thinking on the subject, the technical details being as follows:

\medskip

-- The fact that the set $\pi_1(X)$ defined as above is indeed a group is obvious, with all the group axioms being clear from definitions.

\medskip

-- Obvious as well is the fact that, since $X$ is assumed to be connected, this group does not depend on the choice of the given point $*\in X$, where the loops are based.

\medskip

As basic examples now, for spaces having ``no holes'', such as $\mathbb R$ itself, or $\mathbb R^N$, and so on, we have $\pi_1=\{1\}$. In fact, having no holes can only mean, by definition, $\pi_1=\{1\}$:

\begin{definition}
A space is called simply connected when:
$$\pi_1=\{1\}$$
That is, any loop inside our space must be contractible.
\end{definition}

As further illustrations now for Definition 1.21, here are a few basic computations of homotopy groups, and with (4) proving the last assertion in Theorem 1.19:

\index{free group}

\begin{theorem}
We have the following computations of homotopy groups:
\begin{enumerate}
\item For the circle, we have $\pi_1=\mathbb Z$.

\item For the torus, we have $\pi_1=\mathbb Z\times\mathbb Z$.

\item For the disk minus $2$ points, we have $\pi_1=F_2$.

\item In fact, for the disk minus $N$ points, we have $\pi_1=F_N$.
\end{enumerate}
\end{theorem}

\begin{proof}
These results are all standard, as follows:

\medskip

(1) The first assertion is clear, because a loop on the circle must wind $n\in\mathbb Z$ times around the center, and this parameter $n\in\mathbb Z$ uniquely determines the loop, up to the identification in Definition 1.21. Thus, the homotopy group of the circle is the group of such parameters $n\in\mathbb Z$, which is of course the group $\mathbb Z$ itself.

\medskip

(2) In what regards now the second assertion, the torus being a product of two circles, we are led to the conclusion that its homotopy group must be some kind of product of $\mathbb Z$ with itself. But pictures show that the two standard generators of $\mathbb Z$, and so the two copies of $\mathbb Z$ themselves, commute, $gh=hg$, so we obtain the product of $\mathbb Z$ with itself, subject to commutation, which is the usual product $\mathbb Z\times\mathbb Z$:
$$\left<g,h\Big|gh=hg\right>=\mathbb Z\times\mathbb Z$$

It is actually instructive here to work out explicitly the proof of the commutation relation. We can use the usual drawing convention for the torus, namely:
$$
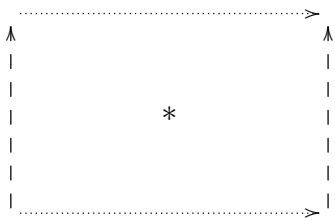
\xymatrix@R=28pt@C=50pt{
\ar@.[rr]&&\\
&\ast\\
\ar@{-->}[uu]\ar@.[rr]&&\ar@{-->}[uu]}$$

The standard generators $g,h$ of the homotopy group are then as follows:
$$
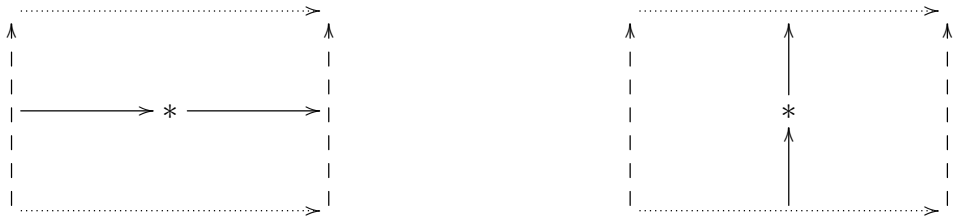
\xymatrix@R=28pt@C=50pt{
\ar@.[rr]&&&&\ar@.[rr]&&\\
\ar[r]&\ast\ar[r]&&&&\ast\ar[u]&\\
\ar@{-->}[uu]\ar@.[rr]&&\ar@{-->}[uu]&&\ar@{-->}[uu]\ar@.[rr]&\ar[u]&\ar@{-->}[uu]}$$

Regarding now the two compositions $gh,hg$, these are as follows:
$$
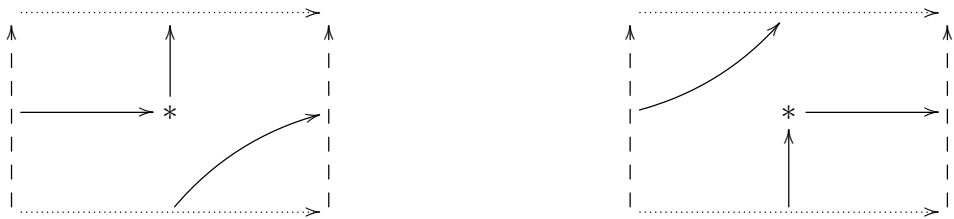
\xymatrix@R=28pt@C=50pt{
\ar@.[rr]&&&&\ar@.[rr]&&\\
\ar[r]&\ast\ar[u]&&&\ar@/_/[ur]&\ast\ar[r]&\\
\ar@{-->}[uu]\ar@.[rr]&\ar@/^/[ur]&\ar@{-->}[uu]&&\ar@{-->}[uu]\ar@.[rr]&\ar[u]&\ar@{-->}[uu]}$$

But these two pictures coincide, up to homotopy, with the following picture:
$$
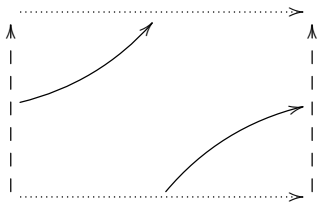
\xymatrix@R=28pt@C=50pt{
\ar@.[rr]&&\\
\ar@/_/[ur]&&\\
\ar@{-->}[uu]\ar@.[rr]&\ar@/^/[ur]&\ar@{-->}[uu]}$$

Thus we have indeed $gh=hg$, as desired, which gives the formula in (2).

\medskip

(3) Regarding now the disk minus $2$ points, the result here is quite clear, because the homotopy group is generated by the 2 loops around the 2 missing points, and these 2 loops are obviously free, algebrically speaking. Thus, we obtain a free product of the group $\mathbb Z$ with itself, which is by definition the free group on 2 generators $F_2$.

\medskip

(4) This is again clear, because the homotopy group is generated here by the $N$ loops around the $N$ missing points, which are free, algebrically speaking. Thus, we obtain a $N$-fold free product of $\mathbb Z$ with itself, which is the free group on $N$ generators $F_N$.
\end{proof}

Many other things can be said about the homotopy groups, and related algebraic invariants. We will be back to topology later in this book, on several occasions.

\section*{1d. Further examples}

Back now to the matrix groups, at a more specialized level, we have the groups $B_N,C_N$, consisting of the orthogonal and unitary bistochastic matrices. Let us start with:

\index{bistochastic matrix}
\index{row-stochastic}
\index{column-stochastic}

\begin{definition}
A square matrix $U\in M_N(\mathbb C)$ is called bistochastic when
$$\begin{matrix}
U_{11}&\ldots&U_{1N}&\to&1\\
\vdots&&\vdots\\
U_{N1}&\ldots&U_{NN}&\to&1\\
\downarrow&&\downarrow\\
1&&1
\end{matrix}$$
that is, when all its rows and columns sum up to $1$.
\end{definition}

In what follows we will be interested in the orthogonal and unitary bistochastic matrices, which are quite interesting objects. These matrices form groups, as follows:

\index{bistochastic group}
\index{discrete Fourier transform}

\begin{theorem}
The real and complex bistochastic groups, which are the sets
$$B_N\subset O_N\quad,\quad 
C_N\subset U_N$$
consisting of matrices which are bistochastic, are isomorphic to $O_{N-1}$, $U_{N-1}$.
\end{theorem}

\begin{proof}
There are several things going on here, the idea being as follows:

\medskip

(1) To start with, the bistochasticity condition in Definition 1.24 can be reformulated in linear algebra terms as follows, with $\xi\in\mathbb C^N$ being the all-one vector:
$$U\xi=U^t\xi=\xi$$

But this makes it clear that if $U,V$ are bistochastic, so is their product $UV$, and we conclude that both the above sets $B_N\subset O_N$ and $C_N\subset U_N$ are indeed groups.

\medskip

(2) Next, let us pick a matrix $F\in U_N$ satisfying the following condition, where $e_0,\ldots,e_{N-1}$ is the standard basis of $\mathbb C^N$, and where $\xi$ is the all-one vector:
$$Fe_0=\frac{1}{\sqrt{N}}\xi$$ 

Observe that such matrices exist indeed, one standard example being the Fourier matrix $F=\frac{1}{\sqrt{N}}(w^{ij})_{ij}$, with $w=e^{2\pi i/N}$. We have then the following computation:
\begin{eqnarray*}
U\xi=\xi
&\iff&UFe_0=Fe_0\\
&\iff&F^*UFe_0=e_0\\
&\iff&F^*UF=diag(1,V)
\end{eqnarray*}

Thus we have isomorphisms as in the statement, given by $V_{ij}\to(F^*UF)_{ij}$.
\end{proof}

As yet another basic example of a continuous group, we have the symplectic group, $Sp_N\subset U_N$ with $N\in2\mathbb N$, which appears as a modification of $O_N$. Let us begin with:

\index{super-space}
\index{super-identity}

\begin{definition}
The ``super-space'' $\bar{\mathbb C}^N$ is the usual space $\mathbb C^N$, with its standard basis $\{e_1,\ldots,e_N\}$, with a chosen sign $\varepsilon=\pm 1$, and a chosen involution on the indices:
$$i\to\bar{i}$$
The ``super-identity'' matrix is $J_{ij}=\delta_{i\bar{j}}$ for $i\leq j$ and $J_{ij}=\varepsilon\delta_{i\bar{j}}$ for $i\geq j$.
\end{definition}

Up to a permutation of the indices, we have a decomposition $N=2p+q$, such that the involution is, in standard permutation notation:
$$(12)\ldots (2p-1,2p)(2p+1)\ldots (q)$$

Thus, up to a base change, the super-identity is as follows, where $N=2p+q$ and $\varepsilon=\pm 1$, with the $1_q$ block at right disappearing if $\varepsilon=-1$:
$$J=\begin{pmatrix}
0&1\ \ \ \\
\varepsilon 1&0_{(0)}\\
&&\ddots\\
&&&0&1\ \ \ \\
&&&\varepsilon 1&0_{(p)}\\
&&&&&1_{(1)}\\
&&&&&&\ddots\\
&&&&&&&1_{(q)}
\end{pmatrix}$$

In the case $\varepsilon=1$, the super-identity is the following matrix:
$$J_+(p,q)=\begin{pmatrix}
0&1\ \ \ \\
1&0_{(1)}\\
&&\ddots\\
&&&0&1\ \ \ \\
&&&1&0_{(p)}\\
&&&&&1_{(1)}\\
&&&&&&\ddots\\
&&&&&&&1_{(q)}
\end{pmatrix}$$

In the case $\varepsilon=-1$ now, the diagonal terms vanish, and the super-identity is:
$$J_-(p,0)=\begin{pmatrix}
0&1\ \ \ \\
-1&0_{(1)}\\
&&\ddots\\
&&&0&1\ \ \ \\
&&&-1&0_{(p)}
\end{pmatrix}$$

With the above notions in hand, we have the following result:

\index{super-orthogonal group}
\index{symplectic group}

\begin{theorem}
The super-orthogonal group, which is by definition
$$\bar{O}_N=\left\{U\in U_N\Big|U=J\bar{U}J^{-1}\right\}$$
with $J$ being the super-identity matrix, is as follows:
\begin{enumerate}
\item At $\varepsilon=1$ we have $\bar{O}_N=O_N$.

\item At $\varepsilon=-1$ we have $\bar{O}_N=Sp_N$.
\end{enumerate}
\end{theorem}

\begin{proof}
These is something quite tricky, the idea being as follows:

\medskip

(1) At $\varepsilon=1$, consider the root of unity $w=e^{\pi i/4}$, and let us set:
$$K=\frac{1}{\sqrt{2}}\begin{pmatrix}w&w^7\\ w^3&w^5\end{pmatrix}$$

This matrix $K$ is then unitary, and we have the following formula:
$$K\begin{pmatrix}0&1\\1&0\end{pmatrix}K^t=1$$

Thus the following matrix is unitary as well, and satisfies $CJC^t=1$:
$$C=\begin{pmatrix}K^{(1)}\\&\ddots\\&&K^{(p)}\\&&&1_q\end{pmatrix}$$

Now in terms of $V=CUC^*$, the relations $U=J\bar{U}J^{-1}=$ unitary simply read:
$$V=\bar{V}={\rm unitary}$$

We conclude that we have an isomorphism $\bar{O}_N=O_N$ as in the statement. 

\medskip

(2) At $\varepsilon=-1$, this depends a bit on what you call symplectic group $Sp_N$, and for our purposes here, we will take the above formula $Sp_N=\bar{O}_N$ as a definition for it.
\end{proof}

We can say more about the symplectic group $Sp_N$, as follows:

\begin{theorem}
The symplectic group $Sp_N\subset U_N$, which is by definition
$$Sp_N=\left\{U\in U_N\Big|U=J\bar{U}J^{-1}\right\}$$
with $J$ being as above, consists of the $SU_2$ patterned matrices,
$$U=\begin{pmatrix}
a&b&\ldots\\
-\bar{b}&\bar{a}\\
\vdots&&\ddots
\end{pmatrix}$$
which are unitary, $U\in U_N$. In particular, we have $Sp_2=SU_2$.
\end{theorem}

\begin{proof}
At $N=2$, to start with, given a matrix $U=\binom{a\ b}{c\ d}$, the condition $U=J\bar{U}J^{-1}$ reformulates as follows, which gives $d=\bar{a}$ and $c=-\bar{b}$, as desired:
$$\begin{pmatrix}a&b\\c&d\end{pmatrix}
\begin{pmatrix}0&1\\-1&0\end{pmatrix}
=\begin{pmatrix}0&1\\-1&0\end{pmatrix}
\begin{pmatrix}\bar{a}&\bar{b}\\\bar{c}&\bar{d}\end{pmatrix}
\iff
\begin{pmatrix}-b&a\\-d&c\end{pmatrix}
=\begin{pmatrix}\bar{c}&\bar{d}\\-\bar{a}&-\bar{b}\end{pmatrix}$$

In the general case, $N\in2\mathbb N$, the proof is similar, with the condition $U=J\bar{U}J^{-1}$ corresponding precisely to the fact that $U$ must be $SU_2$-patterned, as stated.
\end{proof}

We will be back later to the symplectic groups, towards the end of the present book, with more results about them. In the meantime, have a look at the mechanics book of Arnold \cite{ar1}, which explains what the symplectic groups and geometry are good for.

\section*{1e. Exercises}

We had many examples of groups in this chapter, sometimes barely investigated. Of course no worries, we will be back to this, but here are some exercises about this:

\begin{exercise}
Explore the isomorphisms between products of cylic groups.
\end{exercise}

\begin{exercise}
Review some more, if needed, the basic theory of $SO_2,O_2,SU_2$.
\end{exercise}

\begin{exercise}
Experiment as well with $SO_3$, how can you parametrize it?
\end{exercise}

\begin{exercise}
Further explore the groups constructed with generators and relations.
\end{exercise}

\begin{exercise}
Find a family of explicit transformations $T_1,\ldots,T_n$ which are free.
\end{exercise}

\begin{exercise}
Learn more about homotopy groups, and topology in general.
\end{exercise}

\begin{exercise}
Learn more about the bistochastic matrices and groups.
\end{exercise}

\begin{exercise}
Learn also about advanced mechanics, and symplectic groups.
\end{exercise}

As bonus exercise, learn if needed some more linear algebra. This is very related to group theory, and the more linear algebra you know, the better that will be.

\chapter{Finite groups}

\section*{2a. Cyclic groups}

We have seen so far that we have many examples of groups, related to pretty much everything mathematics. In order to have now some theory going, we obviously have to impose some conditions on the groups that we consider. We will impose, for the purposes of this chapter, and in fact for most of the present book, the finiteness condition:

\begin{definition}
A group $(G,\cdot\,)$ is called finite when it is finite as a set:
$$|G|<\infty$$
In this case its cardinality $|G|\in\mathbb N$ is called its order.
\end{definition}

In order to study the finite groups, let us begin with a summary of the examples that we know, from chapter 1. And here, browsing through the material from chapter 1, we can see that the vast majority of the groups considered there were infinite.

\bigskip

Of course these infinite groups $G$ potentially have many interesting finite subgroups $H\subset G$, and studying these will be a main theme in this book. However, when it comes to explicit examples that we already have, these basically reduce to the cyclic groups:

\index{cyclic group}
\index{roots of unity}

\begin{definition}
The cyclic group $\mathbb Z_N$ is defined as follows:
\begin{enumerate}
\item As the additive group of remainders modulo $N$.

\item As the multiplicative group of the $N$-th roots of unity.
\end{enumerate}
This group is finite, having $N$ elements, and is abelian.
\end{definition}

Observe that (1,2) are indeed equivalent, because if we set $w=e^{2\pi i/N}$, then any remainder modulo $N$ defines a $N$-th root of unity, according to the following formula:
$$k\to w^k$$

We obtain in this way all the $N$-roots of unity, and so our correspondence is bijective. Moreover, our correspondence transforms the sum of remainders modulo $N$ into the multiplication of the $N$-th roots of unity, due to the following formula:
$$w^kw^l=w^{k+l}$$

Thus, the groups defined in (1,2) above are indeed isomorphic, via $k\to w^k$, and we agree to denote by $\mathbb Z_N$ the corresponding group, and call it cyclic group.

\bigskip

We will use in what follows both the above conventions for $\mathbb Z_N$, which are both useful, depending on the context. As a general comment about this, let us record:

\begin{comment}
Both the above conventions for $\mathbb Z_N$ are useful. The additive one
$$\mathbb Z_N=\{0,1,2,\ldots,N-1\}$$
is good for doing quick algebra, while the multiplicative one, with $\mathbb Z_N$ being
$$\xymatrix@R=25pt@C=4pt{
&w^2\ar@{.}@/^/[rr]\ar@{.}@/_/[dl]&&w\ar@{.}@/^/[dr]\\
w^3&&\ \ \ast&&1\\
&\ar@{.}@/_/[rr]\ar@{.}@/^/[ul]&&w^{N-1}\ar@{.}@/_/[ur]}$$
with $w=e^{2\pi i/N}$, is obviously ``cyclic'', and brings geometric understanding.
\end{comment}

Observe now that the cyclic groups $\mathbb Z_N$ are indeed abelian, as said above. We can construct further abelian groups by taking products of such cyclic groups, as follows:

\begin{theorem}
The following groups are all finite, and abelian,
$$G=\mathbb Z_{N_1}\times\ldots\times\mathbb Z_{N_k}$$
for any choice of the numbers $N_1,\ldots,N_k\in\mathbb N$.
\end{theorem}

\begin{proof}
This is something trivial, coming from the fact that a product of abelian groups must be abelian too. We have called this Theorem instead of Proposition only because its theoretical importance. Indeed, we will see later, in chapter 4 below, that any finite abelian group must appear as above, as a product of cyclic groups.
\end{proof}

As just mentioned, we will talk more about this later, in chapter 4. In the meantime, however, let us examine a bit the groups appearing in Theorem 2.4. In the simplest case, that of a product of two cyclic groups, we have the following useful result:

\begin{theorem}
Given two integers satisfying $(M,N)=1$, we have:
$$\mathbb Z_M\times\mathbb Z_N\simeq\mathbb Z_{MN}$$
In the case $(M,N)>1$ this fails, and the group on the left is not cyclic. 
\end{theorem}

\begin{proof}
This follows from some basic arithmetic, the idea being as follows:

\medskip

(1) In order to establish the isomorphism in the statement, consider the following map, which is obviously well-defined, for any two positive integers $M,N\in\mathbb N$:
$$f:\mathbb Z_M\times\mathbb Z_N\to\mathbb Z_{MN}\quad,\quad f(a,b)=Na+Mb$$

To be more precise, the fact that $f$ is well-defined comes from the following facts:
$$a=a'(M)\implies Na=Na'(MN)$$
$$b=b'(N)\implies Mb=Mb'(MN)$$

In order to prove now that $f$ is bijective, when $(M,N)=1$, we can invoke a standard result from basic arithmetic, coming as a consequence of the division algorithm for the integers, stating that for $(M,N)=1$, we can always find $p,q\in\mathbb Z$ such that:
$$Mp+Nq=1$$ 

Indeed, this shows that we have $f(q,p)=1$, and by further multiplying everything by a given $c\in\mathbb Z_{MN}$, taken arbitrary, we can have $c\in Im(f)$, as follows:
\begin{eqnarray*}
f(cq,cp)
&=&Ncq+Mcp\\
&=&(Nq+Mp)c\\
&=&c
\end{eqnarray*}

Thus $f$ is surjective, and since its domain and range have the same cardinality, it is bijective. Regarding now the group isomorphism conditions, we first have:
\begin{eqnarray*}
f(a+a',b+b')
&=&N(a+a')+M(b+b')\\
&=&(Na+Mb)+(Na'+Mb')\\
&=&f(a,b)+f(a',b')
\end{eqnarray*}

The condition regarding the units is also satisfied, coming from:
$$f(0,0)=N0+M0=0$$

Finally, the condition regarding the inverses holds too, as shown by:
$$f(-a,-b)=-Na-Mb=-f(a,b)$$

Thus our map $f:\mathbb Z_M\times\mathbb Z_N\to\mathbb Z_{MN}$ is indeed an isomorphism, as desired.
 
\medskip

(2) With this discussed, let us prove now the second assertion, stating that when assuming $(M,N)>1$, we have a non-isomorphism, as follows: 
$$\mathbb Z_M\times\mathbb Z_N\not\simeq\mathbb Z_{MN}$$

But, how to prove such things? Well, we must trick. Assume by contradiction that we have isomorphism, and with the isomorphism, from right to left, being as follows:
$$f:\mathbb Z_{MN}\to \mathbb Z_M\times\mathbb Z_N\quad,\quad f(1)=(a,b)$$

Now if we set $d=[M,N]<MN$, then we have, with $d$ terms in all sums:
\begin{eqnarray*}
f(d)
&=&f(1+\ldots+1)\\
&=&f(1)+\ldots+f(1)\\
&=&(a,b)+\ldots+(a,b)\\
&=&(da,db)\\
&=&(0,0)
\end{eqnarray*}

But this is a contradiction, because the elements $d\neq0$ and $0$ have the same image under $f$, and so $f$ cannot be injective. Thus, non-isomorphism proved, and the very last assertion follows too, because $\mathbb Z_M\times\mathbb Z_N$ having $MN$ elements, and being not isomorphic to the only cyclic group having $MN$ elements, namely $\mathbb Z_{MN}$, it cannot be cyclic.
\end{proof}

Getting back now to general theory, the above proof was quite instructive, containing a number of new ideas, that we would like to formalize next. To start with, inspired by the various maps $f$ used above, which are not necessarily bijective, let us formulate the following definition, fine-tuning our notion of group isomorphism, from chapter 1:

\begin{definition}
A group morphism is a map $f:G\to K$ satisfying:
$$f(gh)=f(g)f(h)$$
$$f(1)=1$$
$$f(g^{-1})=f(g)^{-1}$$
When $f$ is bijective, we say that it is a group isomorphism.
\end{definition}

Many things can be said about group morphisms, and there are countless interesting examples too. Staying theoretical, for the moment, a first result, we have:

\begin{theorem}
Given a group morphism $f:G\to K$, the following happen:
\begin{enumerate}
\item Its kernel $\ker f=\{g\in G|f(g)=1\}$ is a subgroup of $G$.

\item $f$ is injective precisely when $\ker f=\{1\}$.

\item Its image $Im(f)=\{f(g)|g\in G\}$ is a subgroup of $K$.

\item $f$ is surjective precisely when $Im(f)=K$.
\end{enumerate}
\end{theorem}

\begin{proof}
These assertions are all elementary, the idea being as follows:

\medskip

(1) The fact that $\ker f$ is stable under multiplication comes from:
\begin{eqnarray*}
g,h\in\ker f
&\implies&f(g)=1,f(h)=1\\
&\implies&f(gh)=1\\
&\implies&gh\in\ker f
\end{eqnarray*}

Next, the fact that $\ker f$ contains the unit is trivial, coming from:
$$f(1)=1\implies 1\in\ker f$$

Finally, the fact that $\ker f$ is stable under taking inverses comes from:
\begin{eqnarray*}
g\in\ker f
&\implies&f(g)=1\\
&\implies&f(g^{-1})=1\\
&\implies&g^{-1}\in\ker f
\end{eqnarray*}

We conclude that $\ker f\subset G$ is indeed a subgroup, as stated.

\medskip

(2) This is something more subtle, which can be proved as follows:
\begin{eqnarray*}
f\ {\rm injective}
&\iff&[g\neq h\implies f(g)\neq f(h)]\\
&\iff&[gh^{-1}\neq 1\implies f(g)f(h)^{-1}\neq1]\\
&\iff&[gh^{-1}\neq 1\implies f(gh^{-1})\neq1]\\
&\iff&[k\neq1\implies f(k)\neq 1]\\
&\iff&[k\neq 1\implies k\notin\ker f]\\
&\iff&\ker f=\{1\} 
\end{eqnarray*}

(3) The fact that $Im(f)$ is stable under multiplication is clear, coming from:
\begin{eqnarray*}
g,h\in Im(f)
&\implies&g=f(g'),h=f(h')\\
&\implies&gh=f(g'h')\\
&\implies&gh\in Im(f)
\end{eqnarray*}

Next, the fact that $Im(f)$ contains the unit is trivial, coming from:
$$f(1)=1\implies 1\in Im(f)$$

Finally, the fact that $Im(f)$ is stable under taking inverses comes from:
\begin{eqnarray*}
g\in Im(f)
&\implies&g=f(h)\\
&\implies&g^{-1}=f(h^{-1})\\
&\implies&g^{-1}\in Im(f)
\end{eqnarray*}

We conclude that $Im(f)\subset K$ is indeed a subgroup, as stated.

\medskip

(4) Finally, the last assertion is simply the definition of surjectivity.
\end{proof}

The above result is something very useful, and as a first illustration for it, we can reformulate the essentials of Theorem 2.5 and its proof in the following way:

\begin{theorem}
We have a group morphism as follows,
$$f:\mathbb Z_M\times\mathbb Z_N\to\mathbb Z_{MN}\quad,\quad f(a,b)=Na+Mb$$
which is an isomorphism precisely when $(M,N)=1$.
\end{theorem}

\begin{proof}
We already know this from Theorem 2.5 and its proof, but let us see how a brand new proof goes, using Theorem 2.7. And here, things are quick, as follows:

\medskip

(1) To start with, it is clear that $f$ is a group morphism.

\medskip

(2) Since $|\mathbb Z_M\times\mathbb Z_N|=|\mathbb Z_{MN}|$, our map $f$ is bijective when it is injective.

\medskip

(3) When $(M,N)=1$ our map is injective, and so bijective, as shown by:
\begin{eqnarray*}
(a,b)\in\ker f
&\iff&Na+Mb=0(MN)\\
&\iff&a=0,b=0
\end{eqnarray*}

(4) When $(M,N)=d>1$ our map is not injective, and so is not bijective either, as shown by the following computation, with $M=dm$, $N=dn$:
$$f(m,n)=Nm+Mn=0+0=0$$

Thus, we are led to the conclusions in the statement.
\end{proof}

We can further improve the above result, in the following way:

\begin{theorem}
We have a group morphism as follows,
$$f:\mathbb Z_M\times\mathbb Z_N\to\mathbb Z_{MN}\quad,\quad f(a,b)=Na+Mb$$
whose kernel and image are both cyclic groups, as follows,
$$\ker f\simeq\mathbb Z_{(M,N)}\quad,\quad 
Im(f)=\mathbb Z_{[M,N]}$$
and which is an isomorphism precisely when $(M,N)=1$.
\end{theorem}

\begin{proof}
This is again something straightforward, as follows:

\medskip

(1) As before, it is clear that we have a group morphism $f$ as above. 

\medskip

(2) Regarding now the computation of the kernel of $f$, consider the greatest common divisor $d=(M,N)$, and let us write $M=dm$, $N=dn$. We have then:
\begin{eqnarray*}
\ker f
&=&\left\{(a,b)\Big|Na+Mb=0(MN)\right\}\\
&=&\left\{(a,b)\Big|dna+dmb=0(d^2mn)\right\}\\
&=&\left\{(a,b)\Big|na+mb=0(dmn)\right\}
\end{eqnarray*}

Now since $m,n$ are prime to each other, we have the following implications:
$$na+mb=0(dmn)\implies na+mb=0(m)\implies m|a$$
$$na+mb=0(dmn)\implies na+mb=0(n)\implies n|b$$

Thus we must have $a=mp$ and $b=nq$, for some numbers $p,q$, which gives:
\begin{eqnarray*}
\ker f
&=&\left\{(mp,nq)\Big|mnp+mnq=0(dmn)\right\}\\
&=&\left\{(mp,nq)\Big|p+q=0(d)\right\}
\end{eqnarray*}

(3) Now let us examine what we found. The last equation, $p+q=0(d)$, uniquely determines $p$ in terms of $q$, and vice versa, with the solution to it being:
$$p=-q\in\mathbb Z_d$$

Thus, we have an identification $\ker f\simeq\mathbb Z_d$, given from right to left by:
$$p\to(mp,-np)$$

Moreover, since this latter identification is clearly a group morphism, we conclude that we have a group isomorphism as in the statement, namely:
$$\ker f\simeq\mathbb Z_d=\mathbb Z_{(M,N)}$$

(4) Regarding now the computation of the image, with $d=(M,N)$ and $M=dm$, $N=dn$ as before, observe that the formula of our morphism $f$ is as follows:
\begin{eqnarray*}
f(a,b)
&=&Na+Mb\\
&=&dna+dmb\\
&=&d(na+mb)
\end{eqnarray*}

We conclude that we have the following formula, with the embedding in the middle being the standard one, obtained by multiplying the group elements by $d$:
$$f(a,b)\in\mathbb Z_{dmn}\subset\mathbb Z_{d^2mn}=\mathbb Z_{MN}$$ 

Moreover, due to $(m,n)=1$, the above quantities $na+mb$ in the formula of $f(a,b)$ can take any values, and we conclude that we have, as claimed:
$$Im(f)=\mathbb Z_{dmn}=\mathbb Z_{[M,N]}$$

(5) Finally, the last assertion is clear from what we found.
\end{proof}

The above result was something quite advanced, but is not the end of the story with this, as we will soon discover. Also, more generally, we have the question of investigating what happens for the arbitrary products $\mathbb Z_{N_1}\times\ldots\times\mathbb Z_{N_k}$. We will be back to this later, and in the meantime, we will leave some exploration here as an exercise.

\section*{2b. Kernels, images}

Very nice all the above, good things that we are learning, and with everything ultimately coming from Definition 2.6 and Theorem 2.7, via some simple arithmetic. Now in view of this, let us get back to generalities, and formulate the following question:

\begin{question}
What can we say about an arbitrary group morphism
$$f:G\to K$$
generalizing $f$ injective when $\ker f=\{1\}$, and $f$ surjective when $Im(f)=K$?
\end{question}

In answer to this, some experiments first. We already have some substantial data coming from our various computations with cyclic groups, and their products, with the summary of what we know here, along with a little more, being as follows:

\begin{proposition}
We have group morphisms as follows, all surjective:
\begin{enumerate}
\item $g:\mathbb Z\to\mathbb Z_N$ given by $a\to a$, with $\ker g=N\mathbb Z$.

\item $h:\mathbb Z_{MN}\to\mathbb Z_N$ given by $a\to Ma$, with $\ker h=\mathbb Z_M$.

\item $k:\mathbb Z_M\times\mathbb Z_N\to\mathbb Z_{[M,N]}$ given by $(a,b)\to Na+Mb$, with $\ker k=\mathbb Z_{(M,N)}$.
\end{enumerate}
\end{proposition}

\begin{proof}
This is something quite self-explanatory, the idea being as follows:

\medskip

(1) We can certainly talk about the group morphism $g:\mathbb Z\to\mathbb Z_N$ given by $a\to a$, which is obviously surjective, and whose kernel is obviously $\ker g=N\mathbb Z$.

\medskip

(2) Similar situation here, we can certainly talk about the morphism $h:\mathbb Z_{MN}\to\mathbb Z_N$ given by $a\to Ma$, which is obviously surjective, and whose kernel is $\ker h=\mathbb Z_M$.

\medskip

(3) This is the only non-trivial statement, coming from Theorem 2.9, and with the convention that for the image, we use the embedding $\mathbb Z_{[M,N]}\subset\mathbb Z_{MN}$ discussed there.
\end{proof}

Before going further, with some conjectures inspired by Proposition 2.11, let us make a bit more effort, and enlarge as much as possible our list of examples. And here, browsing through the material from chapter 1, we are led to the following statement, involving various basic groups, finite or infinite, and covering Proposition 2.11:

\begin{theorem}
We have group morphisms as follows, all surjective:
\begin{enumerate}
\item $f:\mathbb R\to\mathbb T$ given by $t\to e^{2\pi it}$, with $\ker f=\mathbb Z$.

\item $g:\mathbb Z\to\mathbb Z_N$ given by $k\to e^{2\pi ik/N}$, with $\ker g=N\mathbb Z$.

\item $h:\mathbb Z_{MN}\to\mathbb Z_N$ given by $w\to w^M$, with $\ker h=\mathbb Z_M$.

\item $k:\mathbb Z_M\times\mathbb Z_N\to\mathbb Z_{[M,N]}$ given by $(v,w)\to v^Nw^M$, with $\ker k=\mathbb Z_{(M,N)}$.

\item $p:H\times K\to H$ given by $(h,k)\to h$, with $\ker p=K$.

\item $q:H\times K\to K$ given by $(h,k)\to k$, with $\ker q=H$.
\end{enumerate}
\end{theorem}

\begin{proof}
This is something quite self-explanatory, and with all the cyclic groups in question being written here multiplicatively, the details being as follows:

\medskip

(1) The fact that we have indeed a group morphism is clear, coming from:
\begin{eqnarray*}
f(s+t)
&=&e^{2\pi i(s+t)}\\
&=&e^{2\pi s}e^{2\pi it}\\
&=&f(s)f(t)
\end{eqnarray*}

As for the formula of the kernel, this is something clear too, coming from:
$$e^{2\pi it}=1\iff t\in\mathbb Z$$

(2) This is again clear, with $g:\mathbb Z\to\mathbb Z_N$ appearing as the restriction of $f:\mathbb R\to\mathbb T$, save for a rescaling by $1/N$ on the real line, and with $\ker g=N\mathbb Z$ coming from: 
$$e^{2\pi ik/N}=1\iff k\in N\mathbb Z$$

Alternatively, in additive notation for $\mathbb Z_N$, the morphism $g$ appears as follows, as in Proposition 2.11, and with this making again clear the formula $\ker g=N\mathbb Z$:
$$g:\mathbb Z\to\mathbb Z_N\quad,\quad k\to k$$

(3) This is again something that we know from Proposition 2.11, in additive notation, where our morphism takes the following form, making clear the formula $\ker h=\mathbb Z_M$:
$$h:\mathbb Z_{MN}\to\mathbb Z_N\quad,\quad k\to Mk$$

However, since things in multiplicative notation are nice too, and more geometric, making the connection with the fundamental wrapping operation in (1,2), let us see how the multiplicative proof goes too. Consider the following map from $\mathbb Z_{MN}$ to itself:
$$h:\mathbb Z_{MN}\to \mathbb Z_{MN}\quad,\quad w\to w^M$$

This map is then certainly well-defined, and is a group morphism, as shown by:
$$h(vw)=(vw)^N=v^Nw^N=h(v)h(w)$$

Now observe that the image of $h$ is the group of the $N$-th roots of unity. Thus, by restricting the image, we obtain the morphism in the statement, namely:
$$h:\mathbb Z_{MN}\to\mathbb Z_N\quad,\quad w\to w^M$$

As for the computation of the kernel, this is geometrically clear too, as follows:
$$\ker h
=\left\{w\in\mathbb Z_{MN}\Big|w^M=1\right\}
=\mathbb Z_M$$

(4) Again this is something that we know from Proposition 2.11, or rather from Theorem 2.9, with everything coming, in additive notation, from the following morphism:
$$k:\mathbb Z_M\times\mathbb Z_N\to\mathbb Z_{MN}\quad,\quad (a,b)\to Na+Mb$$

Now observe that in multiplicative notation, this morphism is as follows:
$$k:\mathbb Z_M\times\mathbb Z_N\to\mathbb Z_{MN}\quad,\quad (v,w)\to v^Nw^M$$

As for the computation of the kernel and image, that we already know in additive notation, from Theorem 2.9, we will leave the details here as an instructive exercise.

\medskip

(5) This is something clear, because given a group product $H\times K$, the projection map $p:H\times K\to H$, given by $(h,k)\to h$, is a group morphism, with $\ker p=K$. 

\medskip

(6) This is clear too, because the other projection map, namely $q:H\times K\to K$ given by $(h,k)\to k$, is a group morphism too, with $\ker q=H$.
\end{proof}

The problem is now, what to conjecture, based on the above? We first have:

\begin{observation}
In the context of a surjective group morphism
$$f:G\to K$$
the groups $G,K$ and $H=\ker f$ seem to be related, but not via the formula
$$G=H\times K$$
as basic examples and intuition might suggest, at a first glance.
\end{observation}

To be more precise here, in the context of Theorem 2.12, and going backwards, the formula $G=H\times K$ certainly holds in the context of (5,6), which is nice. However, when skipping (4), which looks more complicated, and getting to (3), our previous study of products of cyclic groups, quickly reported in (4), shows that $\mathbb Z_{MN}\simeq\mathbb Z_M\times\mathbb Z_N$ fails, in general. Also, in the context of (2), we certainly do not have $\mathbb Z\simeq\mathbb Z\times\mathbb Z_N$. And things become even more clear in the context of (1), where we certainly have:
$$\mathbb R\not\simeq\mathbb Z\times\mathbb T$$

In answer now to Question 2.10, and coming with some inspiration from the above examples and counterexamples, we have the following key result:

\begin{theorem}
Given an abelian group $G$ and a subgroup $H\subset G$, we can talk about the quotient group $G/H=\left\{g+H|g\in G\right\}$, which is abelian too, with operation:
$$(g+H)+(k+H)=(g+k)+H$$
With this notion in hand, given a morphism of abelian groups $f:G\to K$ we have:
$$G/\ker f=Im(f)$$
As illustrations for this latter formula, we have the following isomorphisms:
$$\mathbb R/\mathbb Z=\mathbb T\quad,\quad 
\mathbb Z/N\mathbb Z=\mathbb Z_N\quad,\quad 
\mathbb Z_{MN}/\mathbb Z_M=\mathbb Z_N$$
Moreover, with some skill, all this can be extended to the non-abelian case too.
\end{theorem}

\begin{proof}
Many things going on here, the idea being as follows:

\medskip

(1) To start with, all this is a bit abstract, but if you have a simpler solution to Question 2.10, I would be interested. Getting to the proof now, we can certainly talk about $G/H$ as above. Regarding now the group operation, we must first prove that this is well-defined. So, assume that we have group elements $g,k,s,t$ such that:
$$g+H=s+H\quad,\quad k+H=t+H$$

But this means $g=s+e$ and $k=t+f$, for certain $e,f\in H$, and with this, we have the following computation, which shows that our operation is indeed well-defined:
\begin{eqnarray*}
(g+k)+H
&=&(s+e+t+f)+H\\
&=&(s+t)+(e+f+H)\\
&=&(s+t)+H
\end{eqnarray*}

Regarding now the proof of associativity, this goes as follows:
\begin{eqnarray*}
\left[(g+H)+(h+H)\right]+(k+H)
&=&(g+h+H)+(k+H)\\
&=&g+h+k+H\\
&=&(g+H)+(h+k+H)\\
&=&(g+H)+\left[(h+H)+(k+H)\right]
\end{eqnarray*}

Regarding now the unit and inverses, our claim is that these are as follows:
$$0=H\quad,\quad -(g+H)=(-g)+H$$

Indeed, in what regards the unit, this comes from the following computation:
$$(g+H)+H=H+(g+H)=g+H$$

As for the formula for inverses, the verification here goes as follows:
$$(g+H)+(-g+H)=(-g+H)+(g+H)=H$$

Thus, we have our abstract quotient group $G/H$, no question about this.

\medskip

(2) Getting now to the second assertion, by restricting the image we can assume that we are dealing with a surjective morphism $f:G\to K$, and we want to prove that:
$$G/\ker f=K$$

But this is something straightforward, with the isomorphism being as follows:
$$g+\ker f\to f(g)$$

To be more precise, the various verifications to be performed go as follows:

\medskip

-- This map is well-defined. Indeed, this comes from the following computation:
\begin{eqnarray*}
g+\ker f=h+\ker f
&\implies&g-h\in\ker f\\
&\implies&f(g-h)=0\\
&\implies&f(g)=f(h)
\end{eqnarray*}

\medskip

-- This map, say $\varphi$, is a group morphism. Indeed, this can be checked as follows:
\begin{eqnarray*}
\varphi((g+\ker f)+(h+\ker f))
&=&\varphi(g+h+\ker f)\\
&=&f(g+h)\\
&=&f(g)+f(h)\\
&=&\varphi(g+\ker f)+\varphi(h+\ker f)
\end{eqnarray*}

-- This map is injective. Inded, its kernel can be computed as follows:
\begin{eqnarray*}
\ker(\varphi)
&=&\left\{g+\ker f\Big|f(g)=0\right\}\\
&=&\left\{g+\ker f\Big|g\in\ker f\right\}\\
&=&\ker f+\ker f\\
&=&\ker f\\
&=&0
\end{eqnarray*}

-- This map is surjective. Indeed, this is something which is obvious. Thus, done with all the verifications, and the proof of our second assertion is now complete.

\medskip

(3) Next, the various illustrations in the statement come from the various surjective morphisms from the beginning of Theorem 2.12, which give, as stated:
$$\mathbb R/\mathbb Z=\mathbb T\quad,\quad 
\mathbb Z/N\mathbb Z=\mathbb Z_N\quad,\quad 
\mathbb Z_{MN}/\mathbb Z_M=\mathbb Z_N$$

Be said in passing, as a consequence of this, if you ever come upon $\mathbb Z/N\mathbb Z$, in certain algebra classes or books, now you know, that colleagues of mine are talking about $\mathbb Z_N$.

\medskip

(4) Now that we are done with this, let us comment as well about the remaining morphisms from Theorem 2.12. The last two morphisms there give the following isomorphisms, which are of course something true, but not very interesting:
$$(H\times K)/K=H\quad,\quad (H\times K)/H=K$$ 

Thus, we are left with the morphism in (4) there, which was as follows:
$$f:\mathbb Z_M\times\mathbb Z_N\to\mathbb Z_{[M,N]}\quad,\quad\ker f=\mathbb Z_{(M,N)}$$

By applying now what we found above, we are led to the following conclusion:
$$(\mathbb Z_M\times\mathbb Z_N)/\mathbb Z_{(M,N)}=\mathbb Z_{[M,N]}$$

(5) Which looks very nice, but as a matter of better sleeping at night, let us do some doublechecks, and make sure that this is true indeed. As a first observation, the cardinalities match, due to the following well-known formula from arithmetic:
$$MN=(M,N)[M,N]$$ 

Next, as basic particular cases for our result we have $N=1$, where our formula is $\mathbb Z_M/\mathbb Z_1=\mathbb Z_M$, which is true, and then, more generally, the case $(M,N)=1$, where our formula is $\mathbb Z_M\times\mathbb Z_N=\mathbb Z_{MN}$, again true, as we know well since Theorem 2.5. 

\medskip

(6) Next, yet another particular case is $M=N$, where our formula reads:
$$(\mathbb Z_N\times\mathbb Z_N)/\mathbb Z_N=\mathbb Z_N$$

But this is something which is clear too, coming from the morphism $f:\mathbb Z_N\times\mathbb Z_N\to\mathbb Z_N$ given by $f(a,b)=a+b$, which is clearly surjective, and has $\mathbb Z_N$ as kernel.

\medskip

(7) Summarizing, our formula holds indeed, for trivial reasons, in all the cases examined above. As for a non-trivial application, avoiding all the above, and some other coincidences that might happen, at the level of cardinalities, this can only be:
$$(\mathbb Z_4\times\mathbb Z_6)/\mathbb Z_2=\mathbb Z_{12}$$

And we will leave some meditation about this, as an interesting exercise.

\medskip

(8) Finally, in what regards the last assertion, in the non-abelian case we can certainly set $G/H=\left\{gH|g\in G\right\}$, but we will run into troubles if we attempt to define:
$$(gH)(kH)=(gk)H$$

Indeed, asssume that we have group elements $g,k,s,t$ such that:
$$gH=sH\quad,\quad kH=tH$$

But this means $g=se$ and $k=tf$, for certain $e,f\in H$, which leads to the following wrong computation, when attempting to prove that our operation is well-defined:
\begin{eqnarray*}
(gk)H
&=&(setf)H\\
&=&(stef)H\quad\leftarrow {\rm wrong}\\
&=&(st)(efH)\\
&=&(st)H
\end{eqnarray*}

Which does not look easy to understand or fix, and we will leave the discussion here, which is something more technical, for a bit later, at the end of this chapter.
\end{proof}

As a conclusion now, to what we learned so far in this chapter, the cyclic groups $\mathbb Z_N$ are the core objects of finite group theory, and the source of inspiration for anything basic that can be said, in general. However, and here comes the point, the cyclic groups and their quotients remain something abelian, $gh=hg$, and so quite particular, and getting beyond abelianity, with a number of new examples of groups, will be our next task. 

\section*{2c. Dihedral groups}

Passed the cyclic group $\mathbb Z_N$, as a second basic example of a finite group, we have the symmetric group $S_N$. This is something very familiar, appearing as follows:

\index{symmetric group}
\index{permutation group}
\index{signature}
\index{permutation}

\begin{theorem}
The permutations of $\{1,\ldots,N\}$, which are the bijections
$$\sigma:\{1,\ldots,N\}\to\{1,\ldots,N\}$$
form a group $S_N$, called symmetric group, having $N!$ elements.
\end{theorem}

\begin{proof}
It is indeed clear from definitions that $S_N$ is a group. Regarding now the element count, in order to construct a permutation $\sigma\in S_N$, we have:

\medskip

-- $N$ choices for the value of $\sigma(N)$.

\medskip

-- $(N-1)$ choices for the value of $\sigma(N-1)$.

\medskip

-- $(N-2)$ choices for the value of $\sigma(N-2)$.

\smallskip

$\vdots$

\smallskip

-- and so on, up to 1 choice for the value of $\sigma(1)$.

\medskip

Thus, we have indeed $N!$ choices for a permutation $\sigma\in S_N$, as claimed.
\end{proof}

In practice, many interesting things can be said about the symmetric groups $S_N$, ranging from very basic to fairly advanced. Here is a bit more about them:

\begin{theorem}
The symmetric groups $S_N$ have the following properties:
\begin{enumerate}
\item $S_1=\{1\}$.

\item $S_2=\mathbb Z_2$.

\item $S_3$ has $6$ elements, has $\mathbb Z_2$ and $\mathbb Z_3$ as subgroups, and is not abelian.

\item $S_3$ is best viewed as being the symmetry group of an equilateral triangle.

\item $S_N$ has $N!$ elements, and is not abelian at $N\geq3$.
\end{enumerate}
\end{theorem}

\begin{proof}
Many things going on here, mostly about $S_3$, the idea being as follows:

\medskip

(1) This is clear.

\medskip

(2) This is clear again, because the group $S_2$ consists of the identity $id$, and of the permutation $\tau=(1\leftrightarrow 2)$, which satisfies $\tau^2=id$, so that we have:
$$S_2=\left\{id,\tau\Big|\tau^2=id\right\}$$

But this is very similar to the structure of $\mathbb Z_2$, which written multiplicatively, is:
$$\mathbb Z_2=\left\{1,-1\Big|(-1)^2=1\right\}$$

Thus we have indeed an isomorphism $S_2\simeq\mathbb Z_2$, constructed as follows:
$$id\to 1\quad,\quad \tau\to -1$$

(3) Regarding now $S_3$, we know from Theorem 2.15 that this has $3!=6$ elements, that we can represent if we want using our various conventions above. And with this:

\medskip

-- First, we can see that we have $\mathbb Z_2\subset S_3$, and in fact that we have 3 natural inclusions of type $\mathbb Z_2\subset S_3$, coming from the 3 permutations of type $i\leftrightarrow j$.

\medskip

-- Then, we have an inclusion $\mathbb Z_3\subset S_3$, and in fact that we have 2 natural inclusions of type $\mathbb Z_3\subset S_3$, coming from the 2 permutations of type $i\to j\to k\to i$.

\medskip

-- Finally, in what regards the non-abelianity, this comes from the fact that the permutation $1\leftrightarrow 2$ does not commute with the permutation $2\leftrightarrow 3$, as shown by:
$$(1\leftrightarrow 2)(2\leftrightarrow 3)\quad:\quad1\to 2\to 3$$
$$(2\leftrightarrow 3)(1\leftrightarrow 2)\quad:\quad1\to 1\to 2$$

(4) This is obviously something quite informal, but which is good to know, when dealing with $S_3$. To be more precise, the 6 symmetries of $\triangle$ are as follows:

\medskip

-- First we have the rotations of angles $0^\circ$, $120^\circ$, $240^\circ$.

\medskip

-- And then we have the symmetries with respect to the 3 altitudes.

\medskip

So, this is the picture of $S_3$ which is best used, in practice, and as a first observation, I would say that all that we did in (3), namely the fact that $S_3$ contains three copies of $\mathbb Z_2$ and two copies of $\mathbb Z_3$, and is not abelian, is plainly clear with this picture. And here, exercise of course for you to review the proof of (3), with this idea in mind. 

\medskip

(5) In what regards now the groups $S_N$ at higher $N$, as a first observation, we have some obvious embeddings as follows, which are group morphisms:
$$S_1\subset S_2\subset S_3\subset S_4\subset\ldots$$

Now $S_3$ being not abelian, so must be the bigger groups $S_N$ with $N\geq 3$.
\end{proof}

The symmetric groups $S_N$ are key objects of group theory. We will be back to them on many occasions, and notably in chapter 3 below, with a systematic study of them. 

\bigskip

Passed $\mathbb Z_N$ and $S_N$, as a third interesting example of a finite group, which is something more advanced, we have the dihedral group $D_N$, which appears as follows:

\index{dihedral group}
\index{regular polygon}
\index{symmetry group}

\begin{definition}
The dihedral group $D_N$ is the symmetry group of 
$$
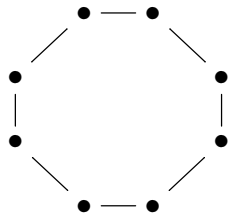
\xymatrix@R=12pt@C=13pt{
&\bullet\ar@{-}[r]\ar@{-}[dl]&\bullet\ar@{-}[dr]\\
\bullet\ar@{-}[d]&&&\bullet\ar@{-}[d]\\
\bullet\ar@{-}[dr]&&&\bullet\ar@{-}[dl]\\
&\bullet\ar@{-}[r]&\bullet}$$
that is, of the regular polygon having $N$ vertices.
\end{definition}

In order to understand how this works, here are the basic examples of regular $N$-gons, at small values of the parameter $N\in\mathbb N$, along with their symmetry groups:

\bigskip

\underline{$N=2$}. Here the $N$-gon is just a segment, and its symmetries are obviously the identity $id$, plus the symmetry $\tau$ with respect to the middle of the segment:
$$
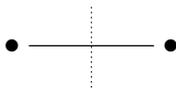
\xymatrix@R=10pt@C=20pt{
&\ar@{.}[dd]\\
\bullet\ar@{-}[rr]&&\bullet\\
&}$$

Thus we have $D_2=\{id,\tau\}$, which in group theory terms means $D_2=\mathbb Z_2$.

\bigskip

\underline{$N=3$}. Here the $N$-gon is an equilateral triangle, and we have 6 symmetries, the rotations of angles $0^\circ$, $120^\circ$, $240^\circ$, and the symmetries with respect to the altitudes: 
$$
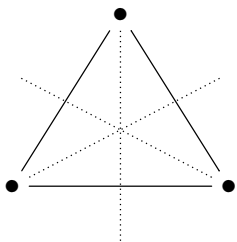
\xymatrix@R=13pt@C=28pt{
&\bullet\ar@{-}[dddr]\ar@{-}[dddl]\ar@{.}[dddd]\\
\ar@{.}[ddrr]&&\ar@{.}[ddll]\\
\\
\bullet\ar@{-}[rr]&&\bullet\\
&
}$$

Alternatively, we can say that the symmetries are all the $3!=6$ possible permutations of the vertices, and so that in group theory terms, we have $D_3=S_3$.

\bigskip

\underline{$N=4$}. Here the $N$-gon is a square, and as symmetries we have 4 rotations, of angles $0^\circ,90^\circ,180^\circ,270^\circ$, as well as 4 symmetries, with respect to the 4 symmetry axes, which are the 2 diagonals, and the 2 segments joining the midpoints of opposite sides:
$$
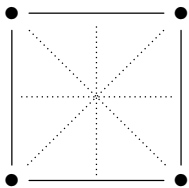
\xymatrix@R=22pt@C=22pt{
\bullet\ar@{-}[dd]\ar@{.}[ddrr]\ar@{-}[rr]&\ar@{.}[dd]&\bullet\ar@{-}[dd]\ar@{.}[ddll]\\
\ar@{.}[rr]&&\\
\bullet\ar@{-}[rr]&&\bullet
}$$

Thus, we obtain as symmetry group some sort of product between $\mathbb Z_4$ and $\mathbb Z_2$. Observe however that this product is not the usual one, our group being not abelian.

\bigskip

\underline{$N=5$}. Here the $N$-gon is a regular pentagon, and as symmetries we have 5 rotations, of angles $0^\circ,72^\circ,144^\circ,216^\circ,288^\circ$, as well as 5 symmetries, with respect to the 5 symmetry axes, which join the vertices to the midpoints of the opposite sides:
$$
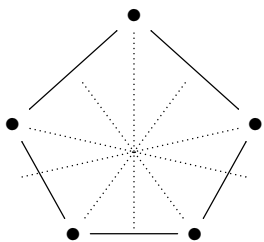
\xymatrix@R=11pt@C=10pt{
&&\bullet\ar@{-}[ddrr]\ar@{-}[ddll]\ar@{.}[dddd]\\
&&&&\\
\bullet\ar@{-}[ddr]\ar@{.}[drrrr]&&&&\bullet\ar@{-}[ddl]\ar@{.}[dllll]\\
&&&&\\
&\bullet\ar@{-}[rr]\ar@{.}[uuurr]&&\bullet\ar@{.}[uuull]&&
}$$

\underline{$N=6$}. Here the $N$-gon is a regular hexagon, and we have 6 rotations, of angles $0^\circ,60^\circ,120^\circ,180^\circ,240^\circ,300^\circ$, and 6 symmetries, with respect to the 6 symmetry axes, which are the 3 diagonals, and the 3 segments joining the midpoints of opposite sides:
$$
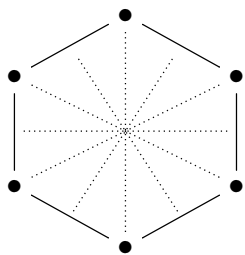
\xymatrix@R=2pt@C=11pt{
&&\bullet\ar@{-}[ddrr]\ar@{-}[ddll]\ar@{.}[dddddddd]\\
&\ar@{.}[ddddddrr]&&\ar@{.}[ddddddll]\\
\bullet\ar@{-}[dddd]\ar@{.}[ddddrrrr]&&&&\bullet\ar@{-}[dddd]\ar@{.}[ddddllll]\\
&&&&\\
\ar@{.}[rrrr]&&&&\\
&&&&\\
\bullet\ar@{-}[ddrr]&&&&\bullet\ar@{-}[ddll]\\
&&&&\\
&&\bullet
}$$

We can see from the above that the various dihedral groups $D_N$ have many common features, and that there are some differences as well, basically coming from the parity of the parameter $N\in\mathbb N$, which changes the structure of the symmetries, inside $D_N$. 

\bigskip

Our purpose in what follows will be that of better understanding the dihedral groups $D_N$, and if possible, with our study being uniform with respect to the parameter $N\in\mathbb N$. Let us start with the following basic fact, inspired by the above observations:

\index{rotations}
\index{symmetries}
\index{dihedral group}

\begin{proposition}
The dihedral group $D_N$ has $2N$ elements, as follows:
\begin{enumerate}
\item We have $N$ rotations $R_1,\ldots,R_N$, with $R_k$ being the rotation of angle $2k\pi/N$. When labeling the vertices of the $N$-gon $1,\ldots,N$, we have $R_k:i\to k+i$.

\item We have $N$ symmetries $S_1,\ldots,S_N$, with $S_k$ being the symmetry with respect to the $Ox$ axis rotated by $k\pi/N$. The symmetry formula is $S_k:i\to k-i$.
\end{enumerate}
\end{proposition}

\begin{proof}
This is clear, indeed. To be more precise, $D_N$ consists of:

\medskip

(1) The $N$ rotations, those of angles $2k\pi/N$ with $k=1,\ldots,N$. But these are exactly the rotations $R_1,\ldots,R_N$ from the statement.

\medskip

(2) The $N$ symmetries with respect to the $N$ possible symmetry axes. But these are exactly the symmetries $S_1,\ldots,S_N$ from the statement.
\end{proof}

With the above description of $D_N$ in hand, we can forget if we want about geometry and the regular $N$-gon, and talk about $D_N$ abstractly, as follows:

\index{multiplication table}

\begin{theorem}
The dihedral group $D_N$ is the group having $2N$ elements, $R_1,\ldots,R_N$ and $S_1,\ldots,S_N$, called rotations and symmetries, which multiply as follows,
$$R_kR_l=R_{k+l}\quad,\quad 
R_kS_l=S_{k+l}$$
$$S_kR_l=S_{k-l}\quad,\quad 
S_kS_l=R_{k-l}$$
with all the indices being taken modulo $N$.
\end{theorem}

\begin{proof}
With notations from Proposition 2.18, the various compositions between rotations and symmetries can be computed as follows:
$$R_kR_l\ :\ i\to l+i\to k+l+i$$
$$R_kS_l\ :\ i\to l-i\to k+l-i$$
$$S_kR_l\ :\ i\to l+i\to k-l-i$$
$$S_kS_l\ :\ i\to l-i\to k-l+i$$

But these are exactly the formulae for $R_{k+l},S_{k+l},S_{k-l},R_{k-l}$, as stated. Now since a group is uniquely determined by its multiplication rules, this gives the result.
\end{proof}

The above result is something very nice, which looks quite final, and we can even write a nice multiplication table for the dihedral group $D_N$, based on it, as follows: 
$$\xymatrix@R=0pt@C=15pt{
&R_1&R_2&\ldots&R_N&S_1&S_2&\ldots&S_N\\
R_1&R_2&R_3&\ldots&R_1&S_2&S_3&\ldots&S_1\\
\vdots\\
R_N&R_1&R_2&\ldots&R_N&S_1&S_2&\ldots&S_N\\
S_1&S_N&S_{N-1}&\ldots&S_1&R_N&R_{N-1}&\ldots&R_1\\
\vdots\\
S_N&S_{N-1}&S_{N-2}&\ldots&S_N&R_{N-1}&R_{N-2}&\ldots&R_N}$$

However, this is not the end of the story which the dihedral group $D_N$. Observe indeed that this group has the same cardinality as the following group:
$$E_N=\mathbb Z_N\times\mathbb Z_2$$

We obviously don't have $D_N\simeq E_N$, because $D_N$ is not abelian, while $E_N$ is. So, our next goal will be that of proving that $D_N$ appears by ``twisting'' $E_N$. Let us start with:

\begin{proposition}
The group $E_N=\mathbb Z_N\times\mathbb Z_2$ is the group having $2N$ elements, $r_1,\ldots,r_N$ and $s_1,\ldots,s_N$, which multiply according to the following rules,
$$r_kr_l=r_{k+l}\quad,\quad 
r_ks_l=s_{k+l}$$
$$s_kr_l=s_{k+l}\quad,\quad 
s_ks_l=r_{k+l}$$
with all the indices being taken modulo $N$.
\end{proposition}

\begin{proof}
With the notation $\mathbb Z_2=\{1,\tau\}$, the elements of the product group $E_N=\mathbb Z_N\times\mathbb Z_2$ can be labeled $r_1,\ldots,r_N$ and $s_1,\ldots,s_N$, as follows:
$$r_k=(k,1)\quad,\quad
s_k=(k,\tau)$$

These elements multiply then according to the formulae in the statement. Now since a group is uniquely determined by its multiplication rules, this gives the result.
\end{proof}

Let us compare now Theorem 2.19 and Proposition 2.20. In order to formally obtain $D_N$ from $E_N$, we must twist some of the multiplication rules of $E_N$, namely:
$$s_kr_l=s_{k+l}\to s_{k-l}$$
$$s_ks_l=r_{k+l}\to r_{k-l}$$

Informally, this amounts in following the rule ``$\tau$ switches the sign of what comes afterwards", and we are led in this way to the following definition:

\index{crossed product}

\begin{definition}
Given groups $H,K$, with an action $K\curvearrowright H$, the crossed product
$$G=H\rtimes K$$
is the set $H\times K$, with multiplication $(g,s)(h,t)=(gh^s,st)$.
\end{definition}

It is routine to check that $G$ is indeed a group. Observe that when the action is trivial, $h^s=h$ for any $h\in H$ and $s\in K$, we obtain the usual product $H\times K$. 

\bigskip

Now with this technology in hand, by getting back to the dihedral group $D_N$, we can improve Theorem 2.13, into a final result on the subject, as follows:

\index{dihedral group}
\index{crossed product decomposition}

\begin{theorem}
We have a crossed product decomposition as follows,
$$D_N=\mathbb Z_N\rtimes\mathbb Z_2$$
with $\mathbb Z_2=\{1,\tau\}$ acting on $\mathbb Z_N$ via switching signs, $k^\tau=-k$.
\end{theorem}

\begin{proof}
We have an action $\mathbb Z_2\curvearrowright\mathbb Z_N$ given by the formula in the statement, namely $k^\tau=-k$, so we can consider the corresponding crossed product group:
$$L_N=\mathbb Z_N\rtimes\mathbb Z_2$$

In order to understand the structure of $L_N$, we follow Proposition 2.20. The elements of $L_N$ can indeed be labeled $\rho_1,\ldots,\rho_N$ and $\sigma_1,\ldots,\sigma_N$, as follows:
$$\rho_k=(k,1)\quad,\quad 
\sigma_k=(k,\tau)$$

Now when computing the products of such elements, we basically obtain the formulae in Proposition 2.20, perturbed as in Definition 2.21. To be more precise, we have:
$$\rho_k\rho_l=\rho_{k+l}\quad,\quad 
\rho_k\sigma_l=\sigma_{k+l}$$
$$\sigma_k\rho_l=\sigma_{k+l}\quad,\quad 
\sigma_k\sigma_l=\rho_{k+l}$$

But these are exactly the multiplication formulae for $D_N$, from Theorem 2.19. Thus, we have an isomorphism $D_N\simeq L_N$ given by $R_k\to\rho_k$ and $S_k\to\sigma_k$, as desired.
\end{proof}

\section*{2d. Normality, order}

Now that we are a bit familiar with the non-abelian groups, let us go back to a question that we have open, namely the extension of Theorem 2.14. Let us start with:

\index{subgroup}
\index{cosets}
\index{left cosets}
\index{right cosets}
\index{order of group}
\index{Lagrange theorem}

\begin{theorem}[Lagrange]
Given a group $G$ and a subgroup $H\subset G$, the sets
$$G/H=\{gH\Big|g\in G\}\quad,\quad H\backslash G=\{Hg\Big|g\in G\}$$
consist of partitions of $G$ into subsets of size $|H|$, and we have the following formula:
$$|G/H|=|H\backslash G|=\frac{|G|}{|H|}$$
In particular, the order of the subgroup divides the order of the group, $|H|\ \big|\ |G|$.
\end{theorem}

\begin{proof}
The partition claim for the set $G/H$ constructed in the statement can be deduced as follows, and the proof for $H\backslash G$ is similar:
$$gH\cap kH\neq\emptyset\iff g^{-1}k\in H\iff gH=kH$$

But with this in hand, the cardinality formulae are all clear. 
\end{proof}

As a continuation of the above, which is something fundamental, we have:

\index{normal subgroup}

\begin{theorem}
Given a subgroup $H\subset G$ which is normal, in the sense that 
$$gH=Hg\quad,\quad\forall g\in G$$
the space $G/H=H\backslash G$ is a group, with multiplication $(gH)(kH)=gkH$.
\end{theorem}

\begin{proof}
We must first prove that the group operation in the statement is indeed well-defined. So, assume that we have group elements $g,k,s,t$ such that:
$$gH=sH\quad,\quad kH=tH$$

We have then the following computation, as desired, using the normality condition:
\begin{eqnarray*}
gkH
&=&gtH\\
&=&gHt\\
&=&sHt\\
&=&stH
\end{eqnarray*}

Regarding now the proof of associativity of our operation, this goes as follows:
\begin{eqnarray*}
\left[(gH)(hH)\right](kH)
&=&(ghH)(kH)\\
&=&ghkH\\
&=&(gH)(hkH)\\
&=&(gH)\left[(hH)(kH)\right]
\end{eqnarray*}

Regarding now the unit and inverses, our claim is that these are as follows:
$$1=H\quad,\quad (g+H)^{-1}=g^{-1}H$$

Indeed, in what regards the unit, this comes from the following computation, which uses the normality condition $gH=Hg$ in order to compute $H(gH)$:
$$(gH)H=H(gH)=gH$$

As for the formula for inverses, the verification here goes as follows, once again by using the normality condition, which produces $gg^{-1}=g^{-1}g=1$ quantities:
$$(gH)(g^{-1}H)=(g^{-1}H)(gH)=H$$

Thus, we have our abstract quotient group $G/H$, as desired.
\end{proof}

The above notion of normality is something quite tricky, for understanding, and for memorizing too, and in practice, the best way to remember it is as follows:

\begin{method}
The normality condition $gH=Hg$ is best remembered as
$$h\in H\implies ghg^{-1}\in H$$
with the keyword for this being ``invariance under conjugation''.
\end{method}

To be more precise here, you only have to remember that the normality condition on $H\subset G$ has something to do with conjugation. Indeed, with this in mind, the relevant condition can only be that $H$ itself must be invariant under conjugation by elements $g\in G$, which corresponds to the formula above, which in turn reads $gHg^{-1}=H$, and so normality, $gH=Hg$. Or at least, that's the method that I use, and it always works.

\bigskip

Next, we have the following result, generalizing the second part of Theorem 2.14:

\index{isomorphism theorem}

\begin{theorem}
Given an arbitrary morphism of groups,
$$f:G\to K$$
its kernel $\ker f\subset G$ is normal, and we have the following isomorphism,
$$G/\ker f=Im(f)$$
extending what we already knew, in the abelian case.
\end{theorem}

\begin{proof}
To start with, the kernel of $f$ is indeed normal, as shown by:
\begin{eqnarray*}
h\in\ker f
&\implies&f(h)=1\\
&\implies&f(ghg^{-1})=1\\
&\implies&ghg^{-1}\in \ker f
\end{eqnarray*}

As for the rest of the proof, this is straightforward, exactly as in the abelian case, and we will leave the verification that it is indeed so as an easy exercise.
\end{proof}

With this discussed, time now for some examples. In relation with the dihedral groups, which are basically the only non-abelian groups that we know well, we have:

\begin{proposition}
The following happen, in relation with $D_N=\mathbb Z_N\rtimes\mathbb Z_2$:
\begin{enumerate}
\item $\mathbb Z_N\subset D_N$ is normal, and $D_N/\mathbb Z_N=\mathbb Z_2$.

\item The various copies $\mathbb Z_2\subset D_N$ are in general not normal,

\item Save for $\mathbb Z_2=\{1,R_{N/2}\}$ when $N$ is even, for which $D_N/\mathbb Z_2=D_{N/2}$.
\end{enumerate}
\end{proposition}

\begin{proof}
Let us recall that $D_N$ is the group having $2N$ elements, $R_1,\ldots,R_N$ and $S_1,\ldots,S_N$, which multiply as follows, with all the indices being modulo $N$:
$$R_kR_l=R_{k+l}\quad,\quad 
R_kS_l=S_{k+l}$$
$$S_kR_l=S_{k-l}\quad,\quad 
S_kS_l=R_{k-l}$$

(1) The subgroup $\mathbb Z_N=\{R_1,\ldots,R_N\}$ is indeed normal, because the conjugates of rotations by either symmetries or rotations are again rotations, as shown by:
$$S_kR_lS_k^{-1}=S_{k-l}S_k=R_{-l}$$
$$R_kR_lR_k^{-1}=R_{k+l}R_{-l}=R_l$$

As for the quotient $D_N/\mathbb Z_N$, this having 2 elements, it must be $\mathbb Z_2$.

\medskip

(2) We have $N$ copies $\mathbb Z_2\subset D_N$, namely $\mathbb Z_2=\{1,S_k\}$, with $S_k$ being the $N$ symmetries, and none of these is normal, at $N>2$, due to the following conjugation formulae:
$$S_lS_kS_l^{-1}=R_{l-k}S_l=S_{2l-k}$$
$$R_lS_kR_l^{-1}=S_{l+k}R_{-l}=S_{2l+k}$$

(3) Thus, we are left with examining $\mathbb Z_2=\{1,R_{N/2}\}$ when $N$ is even, and with this covering the case $N=2$ too. And the subgroup $\mathbb Z_2=\{1,R_{N/2}\}$ is indeed normal, since $R_{N/2}$ commutes with everything, and the quotient is easily seen to be $D_N/\mathbb Z_2=D_{N/2}$.
\end{proof}

The above result is quite interesting, and we have the following extension of it:

\begin{proposition}
The following happen, regarding an arbitrary group $G$:
\begin{enumerate}
\item Any index $2$ subgroup $H\subset G$ is normal, and $G/H=\mathbb Z_2$.

\item A subgroup $\mathbb Z_2\subset G$ is normal when it is central, $\mathbb Z_2\subset Z(G)$.
\end{enumerate}
\end{proposition}

\begin{proof}
This is something a bit more abstract, the idea being as follows:

\medskip

(1) Assuming that $H\subset G$ has index $[G:H]=2$, meaning $|G|=2|H|$, let us try to prove $gH=Hg$. But for $g\in H$ this is clear, and for $g\notin H$ this is clear too, due to:
$$gH=Hg=G-H$$

As for the quotient $G/H$, this having 2 elements, it must be $\mathbb Z_2$.

\medskip

(2) With $\mathbb Z_2=\{1,h\}$, the normality condition on $\mathbb Z_2\subset G$ reads $ghg^{-1}=h$, so $gh=hg$, for any $g\in G$. Thus, $h$ must belong to the center $Z(G)=\{h\in G|gh=hg\}$.
\end{proof}

Moving on, as another continuation of Theorem 2.23, which is something fundamental too, we can talk about the order of group elements, inside any finite group, as follows:

\index{order of element}

\begin{theorem}
Given a finite group $G$, any $g\in G$ generates a cyclic subgroup
$$<g>=\{1,g,g^2,\ldots,g^{k-1}\}$$
with $k=ord(g)$ being the smallest number $k\in\mathbb N$ satisfying $g^k=1$. Also, we have
$$ord(g)\ \big|\ |G|$$
that is, the order of any group element divides the order of the group.
\end{theorem}

\begin{proof}
In order to prove the first assertion, let $g\in G$, and consider the semigroup $<g>\subset G$ formed by the sequence of powers of $g$:
$$<g>=\{1,g,g^2,g^3,\ldots\}\subset G$$

Since $G$ was assumed to be finite, the sequence of powers must cycle, $g^n=g^m$ for some $n<m$. But this shows that $g^k=1$, with $k=m-n$, which gives:
$$<g>=\{1,g,g^2,\ldots,g^{k-1}\}$$

Moreover, we can choose the number $k\in\mathbb N$ to be minimal with this property, and with this choice, we have a set without repetitions. Thus $<g>\subset G$ is indeed a group, and more specifically a cyclic group, whose order is as follows:
$$|<g>|=k=ord(g)$$

Thus, we proved the first assertion, and with this in hand, the second assertion, namely $ord(g)\,|\,|G|$, follows from Theorem 2.23, applied to the subgroup $<g>\subset G$.
\end{proof}

There are many concrete illustrations for the above result, and we will be back to this in the next chapter, when talking about the order of permutations $\sigma\in S_N$.

\bigskip

In the meantime, as a key application, we have the following result:

\index{Fermat theorem}

\begin{theorem}
We have the following congruence, for any prime $p$,
$$a^p=a(p)$$
called Fermat's little theorem.
\end{theorem}

\begin{proof}
There is a simple and a complicated proof for this, as follows:

\medskip

(1) The complicated proof, which is however quite straightforward, using nothing or almost, is by recurrence on $a\in\mathbb N$, using the following computation:
\begin{eqnarray*}
(a+1)^p
&=&\sum_{k=0}^p\binom{p}{k}a^k\\
&=&a^p+1(p)\\
&=&a+1(p)
\end{eqnarray*}

Here we have used the fact that all non-trivial binomial coefficients $\binom{p}{k}$ are multiples of $p$, as shown by a close inspection of these binomial coeffients, given by:
$$\binom{p}{k}=\frac{p(p-1)\ldots(p-k+1)}{k!}$$

Thus, we have the result for any $a\in\mathbb N$, and with the case $p=2$ being trivial, we can assume $p\geq3$, and here by using $a\to-a$ we get it for any $a\in\mathbb Z$, as desired.

\medskip

(2) As for the simple proof, we can argue here that $(\mathbb Z_p^*\,,\cdot\,)$ being obviously a group, having $p-1$ elements, its elements must satisfy $a^{p-1}=1$. Just like that.
\end{proof}

\section*{2e. Exercises}

Welcome to finite groups, and as exercises about them, we have:

\begin{exercise}
When is a product of cyclic groups $\mathbb Z_{N_1}\times\ldots\times\mathbb Z_{N_k}$ cyclic?
\end{exercise}

\begin{exercise}
Experiment some more with $(\mathbb Z_M\times\mathbb Z_N)/\mathbb Z_{(M,N)}=\mathbb Z_{[M,N]}$.
\end{exercise}

\begin{exercise}
Find other interesting examples of $G/H\simeq K$ isomorphisms.
\end{exercise}

\begin{exercise}
Learn more about permutations, including various notations for them.
\end{exercise}

\begin{exercise}
Write down the multiplication table for the groups $S_3$ and $S_4$.
\end{exercise}

\begin{exercise}
Find other interesting examples of $G=H\rtimes K$ decompositions.
\end{exercise}

\begin{exercise}
How many symmetries does the unit cube in $\mathbb R^3$ have?
\end{exercise}

\begin{exercise}
Check, with full details, $G/H\simeq K$, in the non-abelian case. 
\end{exercise}

As bonus exercise, compute orders of group elements $g\in G$, as many as you can.

\chapter{Permutations}

\section*{3a. Cayley theorem}

Let us go back now to the symmetric groups, which are fundamental objects in group theory, as we will soon discover. These groups are constructed as follows:

\index{symmetric group}
\index{permutation group}
\index{signature}

\begin{definition}
The symmetric group $S_N$ is the group of the permutations
$$\sigma:\{1,\ldots,N\}\to\{1,\ldots,N\}$$
of $\{1,\ldots,N\}$. This group has $N!$ elements, and is not abelian at $N\geq3$.
\end{definition}

To be more precise, these are things that we know since chapter 2, with the fact that $S_N$ is indeed a group being clear, and with the last assertion being something clear too. In addition to this, we also know that at $N\leq3$ we have the following formulae:
$$S_1=\{1\}\quad,\quad S_2=\mathbb Z_2\quad,\quad S_3=D_3$$

In practice, many interesting things can be said about $S_N$. To start with, there are many possible notations for the permutations, with the basic ones being as follows:

\bigskip

(1) The simplest way to encode a permutation, which is something very intuitive, consists in writing the numbers $1,\ldots,N$, and below them, their permuted versions:
$$\sigma=\begin{pmatrix}
1&2&3&4&5\\
2&1&4&5&3
\end{pmatrix}$$

(2) As a version of this notation, since the upper numbers are uniquely determined, we can erase them, if we want, and we get to a simpler notation, as follows:
$$\sigma=\begin{pmatrix}
2&1&4&5&3
\end{pmatrix}$$

(3) Another method, which takes a lot of writing, but can be useful for instance when composing, comes by denoting the permutations as diagrams, as follows:
$$\xymatrix@R=4mm@C=4mm{
&1\ar[ddr]&2\ar[ddl]&3\ar[ddrr]&4\ar[ddl]&5\ar[ddl]\\
\sigma=\\
&2&1&4&5&3}$$

(4) Next, the point is that, the upper and lower numbers being uniquely determined, we can erase them, if we want, and we get to a simpler notation, as follows:
$$\xymatrix@R=5mm@C=5mm{
&\ar[ddr]&\ar[ddl]&\ar[ddrr]&\ar[ddl]&\ar[ddl]\\
\sigma=\\
&&&&&}$$

(5) Moreover, we can erase if we want the tips of the arrows too, provided that we always remember that our convention is that things act from top to bottom:
$$\xymatrix@R=5mm@C=5mm{
&\ar@{-}[ddr]&\ar@{-}[ddl]&\ar@{-}[ddrr]&\ar@{-}[ddl]&\ar@{-}[ddl]\\
\sigma=\\
&&&&&}$$

(6) And we will end this discussion with this latter notation, which in practice is the best one, allowing among others quickly composing the permutations. Just remember, as said above, that in standard mathematics, things act from top to bottom. As opposed sometimes to physics, where traditional Feynman diagrams act from left to right.

\bigskip

At the level of the general theory now, by using the symmetric groups, we have the following fundamental result regarding the finite groups, due to Cayley:

\index{Cayley embedding}
\index{permutation group}

\begin{theorem}
Given a finite group $G$, we have an embedding as follows,
$$G\subset S_N\quad,\quad g\to(h\to gh)$$
with $N=|G|$. Thus, any finite group is a permutation group. 
\end{theorem}

\begin{proof}
Given a group element $g\in G$, we can associate to it the following map:
$$\sigma_g:G\to G\quad,\quad 
h\to gh$$

Since $gh=gh'$ implies $h=h'$, this map is bijective, and so is a permutation of $G$, viewed as a set. Thus, with $N=|G|$, we can view this map as a usual permutation, $\sigma_G\in S_N$. Summarizing, we have constructed so far a map as follows:
$$G\to S_N\quad,\quad 
g\to\sigma_g$$

Our first claim is that this is a group morphism. Indeed, this follows from:
$$\sigma_g\sigma_h(k)
=\sigma_g(hk)
=ghk
=\sigma_{gh}(k)$$

It remains to prove that this group morphism is injective. But this follows from:
\begin{eqnarray*}
g\neq h
&\implies&\sigma_g(1)\neq\sigma_h(1)\\
&\implies&\sigma_g\neq\sigma_h
\end{eqnarray*}

Thus, we are led to the conclusion in the statement.
\end{proof}

Observe that in the above statement the embedding $G\subset S_N$ that we constructed depends on a particular writing $G=\{g_1,\ldots,g_N\}$, which is needed in order to identify the permutations of $G$ with the elements of the symmetric group $S_N$. This is not very good, in practice, and as an illustration, for the basic examples of groups that we know, the Cayley theorem provides us with embeddings as follows:
$$\mathbb Z_N\subset S_N\quad,\quad 
D_N\subset S_{2N}\quad,\quad 
S_N\subset S_{N!}$$

And here the first embedding is the good one, the second one is not the best possible one, but can be useful, and the third embedding is useless. Thus, as a conclusion, the Cayley theorem remains something quite theoretical. We will be back to this later, in Part III, with a systematic study of the ``representation'' problem for finite groups.

\bigskip

Getting back now to our main series of finite groups, $\mathbb Z_N\subset D_N\subset S_N$, these are of course permutation groups, according to the above. However, and perhaps even more interestingly, these are as well subgroups of the orthogonal group $O_N$:
$$\mathbb Z_N\subset D_N\subset S_N\subset O_N$$

In order to explain this, we first have the following key result:

\index{permutation matrix}
\index{binary matrix}

\begin{theorem}
We have a group embedding as follows, obtained by regarding $S_N$ as the permutation group of the $N$ coordinate axes of $\mathbb R^N$,
$$S_N\subset O_N$$
which makes $\sigma\in S_N$ correspond to the matrix having $1$ on row $\sigma(j)$ and column $j$, for any $j$, and having $0$ entries elsewhere.
\end{theorem}

\begin{proof}
This is something quite fundamental, the idea being as follows:

\medskip

(1) To start with, we can certainly regard $S_N$ as being the permutation group of the $N$ coordinate axes of $\mathbb R^N$. Now since these permutations of the $N$ coordinate axes of $\mathbb R^N$ are isometries, this provides us with a group embedding $S_N\subset O_N$, as stated.

\medskip

(2) Regarding now the formula of this embedding, we have by definition:
$$\sigma(e_j)=e_{\sigma(j)}$$

Thus, the permutation matrix corresponding to $\sigma$ is given by:
$$\sigma_{ij}=
\begin{cases}
1&{\rm if}\ \sigma(j)=i\\
0&{\rm otherwise}
\end{cases}$$

We are theferore led to the conclusion in the statement.
\end{proof}

In practice now, in order to get familiar with the permutation matrices, nothing beats working out some examples. For $S_2=\{1,\tau\}$ the situation is very simple, as follows:
$$1=\begin{pmatrix}1&0\\0&1\end{pmatrix}\qquad,\qquad
\tau=\begin{pmatrix}0&1\\1&0\end{pmatrix}$$

Regarding now $S_3$, with the $(ijk)$ notation for its elements, we have:
$$(123)=\begin{pmatrix}1&0&0\\0&1&0\\0&0&1\end{pmatrix}\quad,\quad
(213)=\begin{pmatrix}0&1&0\\1&0&0\\0&0&1\end{pmatrix}\quad,\quad
(132)=\begin{pmatrix}1&0&0\\0&0&1\\0&1&0\end{pmatrix}$$
$$(321)=\begin{pmatrix}0&0&1\\0&1&0\\1&0&0\end{pmatrix}\quad,\quad
(231)=\begin{pmatrix}0&0&1\\1&0&0\\0&1&0\end{pmatrix}\quad,\quad
(312)=\begin{pmatrix}0&1&0\\0&0&1\\1&0&0\end{pmatrix}$$

And so on, you get the idea with this. The point now is that we can combine the above result with the Cayley theorem, and we obtain in this way:

\index{finite group}
\index{Cayley embedding}
\index{permutation group}

\begin{theorem}
Given a finite group $G$, we have an embedding as follows,
$$G\subset O_N\quad,\quad g\to(e_h\to e_{gh})$$
with $N=|G|$. Thus, any finite group is an orthogonal matrix group.
\end{theorem}

\begin{proof}
The Cayley theorem gives an embedding as follows:
$$G\subset S_N\quad,\quad g\to(h\to gh)$$

On the other hand, Theorem 3.3 provides us with an embedding as follows:
$$S_N\subset O_N\quad,\quad 
\sigma\to(e_i\to e_{\sigma(i)})$$

Thus, we are led to the conclusion in the statement.
\end{proof}

The same remarks as for the Cayley theorem apply. First, the embedding $G\subset O_N$ that we constructed depends on a particular writing $G=\{g_1,\ldots,g_N\}$. And also, for the basic examples of groups that we know, the embeddings that we obtain are as follows:
$$\mathbb Z_N\subset O_N\quad,\quad 
D_N\subset O_{2N}\quad,\quad 
S_N\subset O_{N!}$$

And here, as before with Cayley, the first embedding is the good one, the second one is not the best possible one, but can be useful, and the third embedding is useless. As already mentioned before, we intend to come back to such questions later, in Part III, with a systematic study of the ``representation problem'' for finite groups.

\bigskip

In any case, all this is not very good, and in order to advance, it is better to forget about the Cayley theorem, and build on Theorem 3.3 instead. So, let us start with the following result, standing as a useful complement to Theorem 3.3:

\begin{theorem}
The embedding $S_N\subset O_N$ is given by the formula
$$\sigma_{ij}=\delta_{i\sigma(j)}$$
with a Kronecker symbol on the right, and has the following properties:
\begin{enumerate}
\item Products of permutation matrices correspond to compositions of permutations.

\item Transposing permutation matrices corresponds to inverting permutations.

\item The permutation matrices are diagonalizable, over $\mathbb C$.

\item Their eigenvalues are roots of unity of order $N$.

\item The trace $Tr(\sigma)$ equals the number of fixed points of $\sigma$.

\item As for the determinant $\det(\sigma)$, this must be $\pm1$.
\end{enumerate}
\end{theorem}

\begin{proof}
We know that $S_N\subset O_N$ is given by $\sigma(e_j)=e_{\sigma(j)}$, and this gives right away the formula in the statement,  $\sigma_{ij}=\delta_{i\sigma(j)}$. As for the other assertions:

\medskip

(1) This comes from the fact that $S_N\subset O_N$ is a group morphism.

\medskip

(2) This comes too from our group morphism $S_N\subset O_N$, via the well-known inversion formula $U^t=U^{-1}$ for the orthogonal matrices.

\medskip

(3) This is something heavier, coming from the spectral theorem in linear algebra, which tells us that any $U\in O_N$ is diagonalizable over $\mathbb C$. In case you wonder if $\mathbb C$ is really needed, here is a quick computation, for the permutation $\sigma=(231)\in S_3$:
$$\det(x-\sigma)=\begin{vmatrix}x&0&-1\\-1&x&0\\0&-1&x\end{vmatrix}=x^3-1$$

(4) This comes again from linear algebra, and more specifically, from:
$$\sigma v=\lambda v\implies \sigma^Nv=\lambda^Nv\implies v=\lambda^Nv\implies\lambda^N=1$$

There are in fact many other things that can be said here, about both eigenvalues and eigenvectors, and with the spectral theorem evoked above not really needed, because the permutation matrices diagonalize quite explicitly. We will be back to this later.

\medskip

(5) This is something elementary, coming from the following computation:
$$Tr(\sigma)=\sum_i\sigma_{ii}=\sum_i\delta_{i\sigma(i)}=\#\left\{i\Big|\sigma(i)=i\right\}$$

(6) This is clear indeed from the definition of $\det$, and more on this later.
\end{proof}

Getting back now to the various groups that we know, we have here:

\begin{theorem}
We have the following finite groups of matrices:
\begin{enumerate}
\item $\mathbb Z_N\subset O_N$, the cyclic permutation matrices.

\item $D_N\subset O_N$, the dihedral permutation matrices.

\item $S_N\subset O_N$, the permutation matrices.
\end{enumerate}
\end{theorem}

\begin{proof}
This is something self-explanatory, the idea being that Theorem 3.3 provides us with embeddings as follows, given by the permutation matrices:
$$\mathbb Z_N\subset D_N\subset S_N\subset O_N$$

In practice now, the groups in the statement appear as follows:

\medskip

(1) The cyclic permutation matrices are by definition the matrices as follows, with 0 entries elsewhere, and form a group, which is isomorphic to the cyclic group $\mathbb Z_N$:
$$U=\begin{pmatrix}
&&&1&&&\\
&&&&1&&\\
&&&&&\ddots&\\
&&&&&&1\\
1&&&&&\\
&\ddots&&&&&\\
&&1&&&
\end{pmatrix}$$

(2) The dihedral matrices are the above cyclic permutation matrices, plus some suitable symmetry permutation matrices, and form a group which is isomorphic to $D_N$.

\medskip

(3) The permutation matrices, which by Theorem 3.3 form a group which is isomorphic to $S_N$, are the $0-1$ matrices having exactly one 1 on each row and column.
\end{proof}

All the above is quite interesting, and having a look now at what we basically know, namely Theorem 3.5 and Theorem 3.6, the following question appears:

\begin{question}
Is there a way of merging our main results, say by computing the laws of the spectral invariants from Theorem 3.5, over the groups in Theorem 3.6?
\end{question}

And good question this is. We will be back to it, at the end of this chapter.

\section*{3b. Cycles, signature}

We would like to discuss now some useful decomposition results, for the permutations. We can indeed talk about the cycle decomposition of permutations, as follows:

\begin{theorem}
Each permutation $\sigma\in S_N$ can be uniquely decomposed as 
$$\sigma=\sigma_1\ldots\sigma_k$$
product of disjoint cycles, and once this decomposition known, we have the formula
$$ord(\sigma)=\Big[|\sigma_1|,\ldots,|\sigma_k|\Big]$$
with $[.]$ standing for the least common multiple, and $|.|$ for the cycle length.  
\end{theorem}

\begin{proof}
The first assertion is something self-explanatory, but let us discuss this in detail, since we will systematically need all this, in what follows:

\medskip

(1) Regarding cycles, these are by definition the permutations acting as follows, on certain distinct elements $i_1,\ldots,i_k\in\{1,\ldots N\}$, and acting elsewhere by the identity:
$$i_1\to i_2\to i_3\to\ldots\to i_s\to i_1$$

The number $s$ is called length of the cycle $\sigma$, and we use the following notation:
$$|\sigma|=s$$

We will often assume $s\geq2$, as being part of the definition of the cycles. For the purposes of the present theorem, we can indeed make this assumption for what follows, by excluding the identity permutation from our cycle decomposition study.

\medskip

(2) The simplest examples of cycles are those having length 2. Such cycles are called transpositions, and their action is by exchanging two elements, as follows:
$$
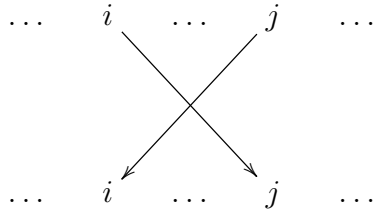
\xymatrix@R=50pt@C=15pt{
\ldots&i\ar[drr]&\ldots&j\ar[dll]&\ldots\\
\ldots&i&\ldots&j&\ldots
}$$

(3) Next, we have the cycles of lenght 3, which can appear in two configurations. Indeed, given $i<j<k$, a first possibility for our 3-cycle on $\{i,j,k\}$ is as follows:  
$$
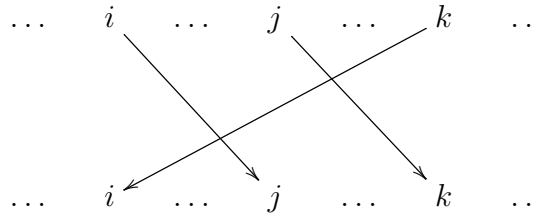
\xymatrix@R=50pt@C=15pt{
\ldots&i\ar[drr]&\ldots&j\ar[drr]&\ldots&k\ar[dllll]&\ldots\\
\ldots&i&\ldots&j&\ldots&k&\ldots
}$$

But then, we have as well a second possibility, as follows:
$$
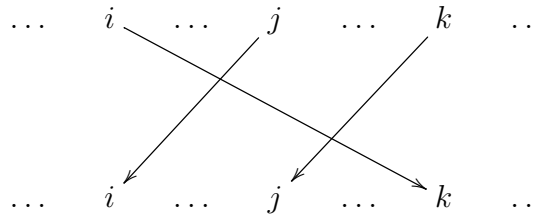
\xymatrix@R=50pt@C=15pt{
\ldots&i\ar[drrrr]&\ldots&j\ar[dll]&\ldots&k\ar[dll]&\ldots\\
\ldots&i&\ldots&j&\ldots&k&\ldots
}$$

(4) And so on, you get the point, with the number of possibilities for a $s$-cycle on a certain set of $s$ distinct elements $i_1<\ldots<i_s$ being $(s-1)!$, for obvious reasons.

\medskip

(5) Getting now to what the theorem says, the decomposition result there is something obvious, simply by picking $i\in\{1,\ldots,N\}$, and then iterating our permutation $\sigma$ on it. Indeed, we will obtain in this way certain elements, as follows:
$$i,\sigma(i),\sigma^2(i),\sigma^3(i),\ldots$$

But these will be the elements of the cycle passing through $i$, so done with that cycle. And then, once this cycle passing through $i$ determined, we can simply erase it, and repeat the procedure, by recurrence, until we completely decompose $\sigma$.

\medskip

(6) Finally, in what regards the last assertion, assume indeed that we have a permutation $\sigma\in S_N$ decomposed as a product of disjoint cycles, as follows:
$$\sigma=\sigma_1\ldots\sigma_k$$

Since the cycles $\sigma_1,\ldots,\sigma_k$ commute, the powers of $\sigma$ are given by:
$$\sigma^n=\sigma_1^n\ldots\sigma_k^n$$

Now when looking at the order, we are led to the following formula:
$$ord(\sigma)=\Big[|\sigma_1|,\ldots,|\sigma_k|\Big]$$

Thus, we are led to the conclusions in the statement.
\end{proof}

At a more advanced level now, we have the following result, that you surely know from linear algebra, and more specifically, from the theory of the determinant:

\index{signature}
\index{number of inversions}
\index{transpositions}
\index{odd cycles}
\index{crossings}

\begin{theorem}
The permutations have a signature function
$$\varepsilon:S_N\to\{\pm1\}$$
which can be defined in the following equivalent ways:
\begin{enumerate}
\item As $(-1)^c$, where $c$ is the number of inversions.

\item As $(-1)^t$, where $t$ is the number of transpositions.

\item As $(-1)^o$, where $o$ is the number of odd cycles.

\item As $(-1)^x$, where $x$ is the number of crossings.

\item As the sign of the corresponding permuted basis of $\mathbb R^N$.
\end{enumerate}
\end{theorem}

\begin{proof}
We have explain what the numbers $c,t,o,x$ appearing in (1-4) exactly are, then why they are well-defined modulo 2, then why they are equal to each other, and finally why the constructions (1-4) yield the same sign as (5). Let us begin with the first two steps, namely precise definition of the numbers $c,t,o,x$, modulo 2:

\medskip

(1) The idea here is that given any two numbers $i<j$ among $1,\ldots,N$, the permutation  can either keep them in the same order, $\sigma(i)<\sigma(j)$, or invert them:
$$\sigma(j)>\sigma(i)$$

Now by making $i<j$ vary over all pairs of numbers in $1,\ldots,N$, we can count the number of inversions, and call it $c$. This is an integer, $c\in\mathbb N$, which is well-defined.

\medskip

(2) Here the idea, which is something quite intuitive, is that any permutation appears as a product of switches, also called transpositions: 
$$i\leftrightarrow j$$

The decomposition as a product of transpositions is not unique, but the number $t$ of the needed transpositions is unique, when considered modulo 2. This follows for instance from the equivalence of (2) with (1,3,4,5), explained below.

\medskip

(3) Here the point is that any permutation decomposes, in a unique way, as a product of cycles, which are by definition permutations of the following type:
$$i_1\to i_2\to i_3\to\ldots\ldots\to i_k\to i_1$$

Some of these cycles have even length, and some others have odd length. By counting those having odd length, we obtain a well-defined number $o\in\mathbb N$.

\medskip

(4) Here the method is that of drawing the permutation, as we usually do, and by avoiding triple crossings, and then counting the number of crossings. This number $x$ depends on the way we draw the permutations, but modulo 2, we always get the same number. Indeed, this follows from the fact that we can continuously pass from a drawing to each other, and that when doing so, the number of crossings can only jump by $\pm2$.

\medskip

Summarizing, we have 4 different definitions for the signature of the permutations, which all make sense, constructed according to (1-4) above. Regarding now the fact that we always obtain the same number, this can be established as follows:

\medskip

(1)=(2) This is clear, because any transposition inverts once, modulo 2.

\medskip

(1)=(3) This is clear as well, because the odd cycles invert once, modulo 2.

\medskip

(1)=(4) This comes from the fact that the crossings correspond to inversions.

\medskip

(2)=(3) This follows by decomposing the cycles into transpositions.

\medskip

(2)=(4) This comes from the fact that the crossings correspond to transpositions.

\medskip

(3)=(4) This follows by drawing a product of cycles, and counting the crossings.

\medskip

Finally, in what regards the equivalence of all these constructions with (5), here simplest is to use (2). Indeed, we already know that the sign of a system of vectors switches when interchanging two vectors, and so the equivalence between (2,5) is clear. 
\end{proof}

As already mentioned, and as you surely know too, the permutations and their signature are key ingredients in linear algebra, in the theory of the determinant. It is tempting to take a break at this point from group theory, and talk a bit about this.

\bigskip

To start with, here is the definition of the determinant, I mean the correct one:

\begin{definition}
The determinant of a square real matrix is the signed volume of the parallelepiped formed by its column vectors. That is, we have
$${\rm det}(v_1\ldots v_N)=\pm vol<v_1,\ldots,v_N>$$
with the sign being $+$ when $\{v_1,\ldots,v_N\}\subset\mathbb R^N$ is positively oriented, in the sense that one can continuously pass from it to the standard basis of $\mathbb R^N$, and is $-$ otherwise.
\end{definition}

In case this sounds new to you, no fear, and as an illustration for this, here is the computation of the determinant for the $2\times2$ matrices, which is quite easy:

\begin{theorem}
In $2$ dimensions we have the following formula,
$$\det\begin{pmatrix}a&b\\ c&d\end{pmatrix}=ad-bc$$
with $\det:M_2(\mathbb R)\to\mathbb R$ being the function constructed above.
\end{theorem}

\begin{proof}
We must show that the area of the parallelogram formed by $\binom{a}{c},\binom{b}{d}$ equals $|ad-bc|$. We can assume $a,b,c,d>0$ for simplifying, the proof in general being similar. Moreover, by switching if needed the vectors $\binom{a}{c},\binom{b}{d}$, we can assume that we have:
$$\frac{a}{c}>\frac{b}{d}$$

Now let us slide the upper side of the parallelogram downwards left, as follows:
$$
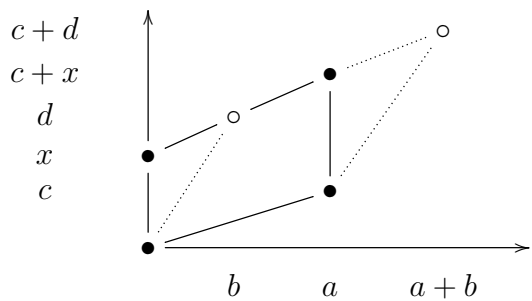
\xymatrix@R=1pt@C=16pt{
&&&\\
c+d&&&&\circ\\
c+x&&&\bullet\ar@{.}[ur]&\\
d&&\circ\ar@{-}[ur]&&&&\\
x&\bullet\ar@{-}[ur]\ar[uuuu]&&&\\
c&&&\bullet\ar@{-}[uuu]\ar@{.}[uuuur]&\\
&&&&\\
&\bullet\ar@{-}[uurr]\ar@{-}[uuu]\ar[rrrr]\ar@{.}[uuuur]&&&&\\
&&\ b\ &\ a\ &a+b}$$

Next, let us slide the new right side of the parallelogram downwards, as follows:
$$
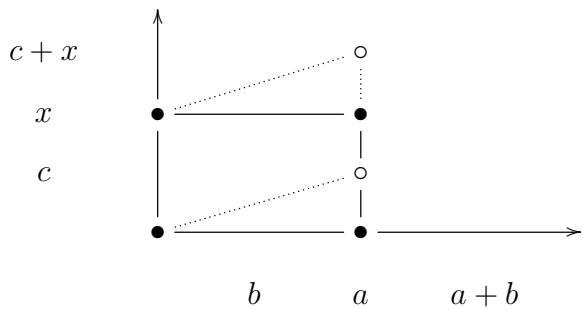
\xymatrix@R=10pt@C=20pt{
&&&\\
c+x&&&\circ&\\
x&\bullet\ar@{.}[urr]\ar[uu]\ar@{-}[rr]&&\bullet\ar@{.}[u]&\\
c&&&\circ\ar@{-}[u]&\\
&\bullet\ar@{.}[urr]\ar@{-}[uu]\ar@{-}[rr]&&\bullet\ar@{-}[u]\ar[rr]&&\\
&&\ b\ &\ a\ &a+b}$$

Let us compute now the area. Since our two sliding operations have not changed the area of the original parallelogram, this area is given by:
$$A=ax$$

In order to compute now the quantity $x$, observe that in the context of the first move we have two similar triangles, which lead by Thales to the following equation:
$$\frac{d-x}{b}=\frac{c}{a}$$

Thus we have $x=d-bc/a$, so the area is $A=ad-bc$, as desired.
\end{proof}

In the general case now, by similarly playing with the Thales theorem and other elementary geometry tools, we are led to some rules for computing the determinants, that you surely know well. And the point is that these rules eventually lead to:

\index{determinant formula}
\index{Sarrus formula}
\index{volume of parallelepiped}
\index{determinant}
\index{linear algebra}

\begin{theorem}
We have the following formula for the determinant,
$$\det A=\sum_{\sigma\in S_N}\varepsilon(\sigma)A_{1\sigma(1)}\ldots A_{N\sigma(N)}$$
with the signature function being the one introduced above.
\end{theorem}

\begin{proof}
This follows by recurrence over $N\in\mathbb N$, as follows:

\medskip

(1) When developing the determinant over the first column, we obtain a signed sum of $N$ determinants of size $(N-1)\times(N-1)$. But each of these determinants can be computed by developing over the first column too, and so on, and we are led to the conclusion that we have a formula as in the statement, with $\varepsilon(\sigma)\in\{-1,1\}$ being certain coefficients.

\medskip

(2) But these latter coefficients $\varepsilon(\sigma)\in\{-1,1\}$ can only be the signatures of the corresponding permutations $\sigma\in S_N$, with this being something that can be viewed again by recurrence, with either of the definitions (1-5) in Theorem 3.9 for the signature.
\end{proof}

As an illustration, in 2 dimensions we recover the usual formula, as follows:
$$\begin{vmatrix}a&b\\ c&d\end{vmatrix}
=\varepsilon(|\,|)\cdot ad+\varepsilon(\slash\hskip-2mm\backslash)\cdot cb
=ad-bc$$

In 3 dimensions now, we recover the well-known Sarrus formula, as follows:
$$\begin{vmatrix}a&b&c\\ d&e&f\\ g&h&i\end{vmatrix}=aei+bfg+cdh-ceg-bdi-afh$$

Observe that the triangles in the Sarrus formula correspond to the permutations of $\{1,2,3\}$, and their signs correspond to the signatures of these permutations:
\begin{eqnarray*}
\det
&=&\begin{pmatrix}*&&\\ &*&\\ &&*\end{pmatrix}
+\begin{pmatrix}&*&\\ &&*\\ *&&\end{pmatrix}
+\begin{pmatrix}&&*\\ *&&\\ &*&\end{pmatrix}\\
&-&\begin{pmatrix}&&*\\ &*&\\ *&&\end{pmatrix}
+\begin{pmatrix}&*&\\ *&&\\ &&*\end{pmatrix}
+\begin{pmatrix}*&&\\ &&*\\ &*&\end{pmatrix}
\end{eqnarray*}

Getting back now to groups, by using the above technology, we can formulate:

\begin{theorem}
The following set is a group, called alternating group,
$$A_N=\left\{\sigma\in S_N\Big|\varepsilon(\sigma)=1\right\}$$
and we have $A_N=S_N\cap SO_N$, as subgroups of $O_N$.
\end{theorem}

\begin{proof}
We have several assertions here, the idea being as follows:

\medskip

(1) In what regards the first assertion, we can say that the alternating group $A_N$ appears as the kernel of the signature map $\varepsilon:S_N\to\mathbb Z_2$. Moreover, we have:
$$S_N/A_N=\mathbb Z_2$$

(2) Consider now the standard embedding $S_N\subset O_N$, obtained by permuting the coordinate axes of $\mathbb R^N$, which is given by the permutation matrices. The determinant of a permutation $\sigma\in S_N$ is then its signature, and this gives $A_N=S_N\cap SO_N$, as stated.
\end{proof}

As yet another application of our technology in this section, in relation with some questions that we had open, in the context of Theorem 3.5, let us record as well:

\begin{theorem}
Given a permutation $\sigma\in S_N$, with cycle decomposition
$$\sigma=\sigma_1\ldots\sigma_k$$
the corresponding permutation matrix $U\in O_N$ can be diagonalized as follows:
\begin{enumerate}
\item $U$ is block-diagonal, with respect to the cycle decomposition of $\sigma$.

\item The eigenvalues corresponding to $\sigma_s$ are the $s$-th roots of unity.

\item The corresponding eigenvectors are formed by $s$-th roots of unity too.
\end{enumerate}
\end{theorem}

\begin{proof}
This is indeed something self-explanatory, obtained by merging what we have in Theorems 3.5 and 3.8, and we will leave the details here, as an instructive exercise.
\end{proof}

\section*{3c. Graph symmetries}

Getting now to the various subgroups $G\subset S_N$, many interesting things can be said about the symmetry groups of finite graphs. Let us start our study here with:

\begin{theorem}
Given a finite graph $X$, with vertices denoted $1,\ldots,N$, the symmetries of $X$, which are the permutations $\sigma\in S_N$ leaving invariant the edges, 
$$i-j\implies\sigma(i)-\sigma(j)$$
form a subgroup of the symmetric group, as follows, called symmetry group of $X$:
$$G(X)\subset S_N$$
As basic examples, for the empty graph, or for the simplex, we have $G(X)=S_N$.
\end{theorem}

\begin{proof}
Here the first assertion, regarding the group property of $G(X)$, is clear from definitions, because the symmetries of $X$ are stable under composition. The second assertion, regarding the empty graph and the simplex, is clear as well.
\end{proof}

Let us work out now some further examples. As a first result here, we have:

\begin{theorem}
The symmetry group of the $N$-gon graph
$$
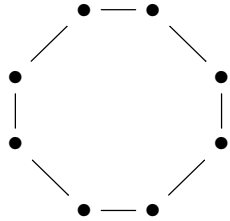
\xymatrix@R=13pt@C=13pt{
&\bullet\ar@{-}[r]\ar@{-}[dl]&\bullet\ar@{-}[dr]\\
\bullet\ar@{-}[d]&&&\bullet\ar@{-}[d]\\
\bullet\ar@{-}[dr]&&&\bullet\ar@{-}[dl]\\
&\bullet\ar@{-}[r]&\bullet}$$
is the dihedral group $D_N$.
\end{theorem}

\begin{proof}
This is something that we know well from chapter 2, the point being that the $N$-gon graph has the same symmetries as the regular $N$-gon, metric space.
\end{proof}

In order to reach now to some general theory, let us formulate:

\begin{definition}
A graph is called circulant if, when drawn with its $N$ vertices on a circle, it is invariant under rotations. Equivalently, we must have $\mathbb Z_N\subset G(X)$.
\end{definition}

Generally speaking, these graphs are the most ``beautiful'' ones, no matter your aesthetics. As a basic example here, we have the wheel graph, which is as follows:
$$
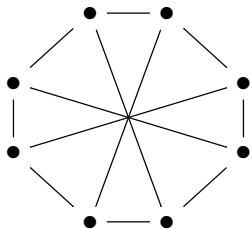
\xymatrix@R=14pt@C=16pt{
&\bullet\ar@{-}[r]\ar@{-}[dddr]\ar@{-}[dl]&\bullet\ar@{-}[dr]\ar@{-}[dddl]\\
\bullet\ar@{-}[d]\ar@{-}[drrr]&&&\bullet\ar@{-}[d]\\
\bullet\ar@{-}[dr]\ar@{-}[urrr]&&&\bullet\ar@{-}[dl]\\
&\bullet\ar@{-}[r]&\bullet}$$

Here is another circulant graph, again with 8 vertices, again with the picture suggesting the name ``circulant'', and of course, again being a beautiful graph as well:
$$
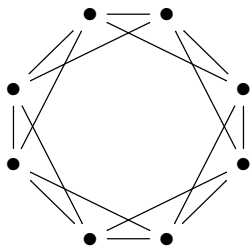
\xymatrix@R=16pt@C=16pt{
&\bullet\ar@{-}[r]\ar@{-}[drr]\ar@{-}[dl]&\bullet\ar@{-}[dr]\ar@{-}[dll]\\
\bullet\ar@{-}[d]\ar@{-}[ddr]&&&\bullet\ar@{-}[d]\\
\bullet\ar@{-}[dr]\ar@{-}[uur]\ar@{-}[drr]&&&\bullet\ar@{-}[dl]\ar@{-}[uul]\\
&\bullet\ar@{-}[r]\ar@{-}[urr]&\bullet\ar@{-}[uur]}$$

Getting back now to Definition 3.17 as stated, it is understood there that the vertices there are labeled in advance $1,\ldots,N$. In general, when the vertices are not labeled in advance, the convention is that the graph is called circulant when it is possible to label the vertices  $1,\ldots,N$, as for the graph to become circulant in the above sense.

\bigskip

In order to understand this, let us pick a graph which is obviously circulant, such as the wheel graph above, and mess up the labeling of the vertices, see what we get. For this purpose, let us put some random labels $1,2,\dots,8$ on our wheel graph, say as follows:
$$
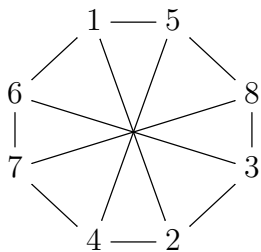
\xymatrix@R=13pt@C=17pt{
&1\ar@{-}[r]\ar@{-}[dddr]\ar@{-}[dl]&5\ar@{-}[dr]\ar@{-}[dddl]\\
6\ar@{-}[d]\ar@{-}[drrr]&&&8\ar@{-}[d]\\
7\ar@{-}[dr]\ar@{-}[urrr]&&&3\ar@{-}[dl]\\
&4\ar@{-}[r]&2}$$

Now let us redraw this graph, with the vertices $1,2,\ldots,8$ ordered on a circle, equally spaced, as Definition 3.17 requires. We get something not very beautiful, as follows:
$$
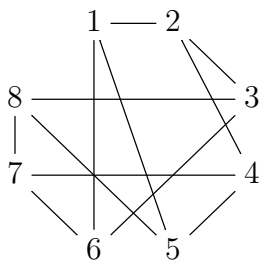
\xymatrix@R=14pt@C=17pt{
&1\ar@{-}[r]\ar@{-}[dddr]\ar@{-}[ddd]&2\ar@{-}[dr]\ar@{-}[ddr]\\
8\ar@{-}[d]&&&3\ar@{-}[lll]\ar@{-}[ddll]\\
7\ar@{-}[dr]&&&4\ar@{-}[dl]\ar@{-}[lll]\\
&6&5\ar@{-}[uull]}$$

So, here is the point. This graph, regarded as a graph with vertices labeled $1,2,\ldots,8$ is obviously not circulant, in the sense of Definition 3.17. However, when removing the labels, this graph does become circulant, as per our conventions above. 

\bigskip

All this might seem a bit confusing, when first seen, and you are probably in this situation, so here is a precise statement in this sense, coming with a full proof:

\begin{proposition}
The following graph is circulant, despite its bad look,
$$
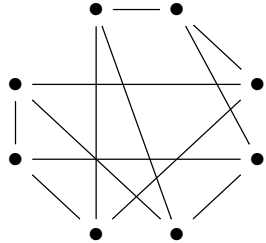
\xymatrix@R=16pt@C=17.5pt{
&\bullet\ar@{-}[r]\ar@{-}[dddr]\ar@{-}[ddd]&\bullet\ar@{-}[dr]\ar@{-}[ddr]\\
\bullet\ar@{-}[d]&&&\bullet\ar@{-}[lll]\ar@{-}[ddll]\\
\bullet\ar@{-}[dr]&&&\bullet\ar@{-}[dl]\ar@{-}[lll]\\
&\bullet&\bullet\ar@{-}[uull]}$$
in the sense that it can be put in circulant form, with a suitable labeling of the vertices.
\end{proposition}

\begin{proof}
As already mentioned, this normally follows from the above discussion, but let us prove this as well directly. The idea is as follows:

\medskip

(1) Let us label the vertices of our graph as follows, and I will explain in moment where this tricky labeling choice comes from:
$$
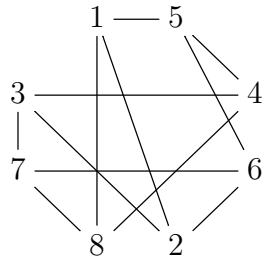
\xymatrix@R=14pt@C=17pt{
&1\ar@{-}[r]\ar@{-}[dddr]\ar@{-}[ddd]&5\ar@{-}[dr]\ar@{-}[ddr]\\
3\ar@{-}[d]&&&4\ar@{-}[lll]\ar@{-}[ddll]\\
7\ar@{-}[dr]&&&6\ar@{-}[dl]\ar@{-}[lll]\\
&8&2\ar@{-}[uull]}$$

Now let us redraw this graph, with the vertices $1,2,\ldots,8$ ordered on a circle, equally spaced, as Definition 3.17 requires. We get something very nice, as follows:
$$
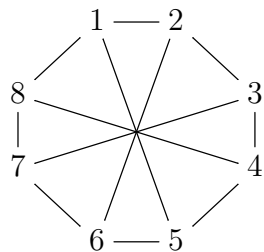
\xymatrix@R=13pt@C=17pt{
&1\ar@{-}[r]\ar@{-}[dddr]\ar@{-}[dl]&2\ar@{-}[dr]\ar@{-}[dddl]\\
8\ar@{-}[d]\ar@{-}[drrr]&&&3\ar@{-}[d]\\
7\ar@{-}[dr]\ar@{-}[urrr]&&&4\ar@{-}[dl]\\
&6\ar@{-}[r]&5}$$

Thus, our original graph was indeed circulant, as stated.

\medskip

(2) In order for everything to be fully clarified, we still must explain where the tricky labeling choice in (1) comes from. For this purpose, let us recall where the graph in the statement came from. Well, this graph was obtained by messing up the labeling of the vertices of the wheel graph, by using the following permutation: 
$$\sigma=\begin{pmatrix}
1&2&3&4&5&6&7&8\\
1&5&8&3&2&4&7&6
\end{pmatrix}$$

The point now is that, if we want to unmess our graph, we must use the inverse of the above permutation, obtained by reading things upside-down, which is given by:
$$\sigma^{-1}=\begin{pmatrix}
1&2&3&4&5&6&7&8\\
1&5&4&6&2&8&7&3
\end{pmatrix}$$

Thus, we must label our vertices $1,5,4,6,2,8,7,3$, precisely as done in (1).
\end{proof}

Moving on, with some general theory, we first have a basic observation, as follows:

\begin{proposition}
A graph $X$, with vertices labeled $1,2,\ldots,N$, is circulant precisely when for any two vertices $i,j\in X$ there is a permutation $\sigma\in S_N$ such that:
\begin{enumerate}
\item $\sigma$ maps one vertex to another, $\sigma(i)=j$.

\item $\sigma$ is cyclic, $\sigma(k)=k+s$ modulo $N$, for some $s$.

\item $\sigma$ leaves invariant the edges, $k-l\iff\sigma(k)-\sigma(l)$.
\end{enumerate}
\end{proposition}

\begin{proof}
This is obvious from definitions, and with the remark that the number $s$ appearing in (2) is uniquely determined by (1), as being $s=j-i$, modulo $N$.
\end{proof}

The point now is that, with this picture of the circulant graphs in mind, it is quite clear that if we remove the assumption (2), that our permutation is cyclic, we will reach to a quite interesting class of graphs, generalizing them. So, let us formulate:

\begin{definition}
A graph $X$, with vertices labeled $1,2,\ldots,N$, is called transitive when for any two vertices $i,j\in X$ there is a permutation $\sigma\in S_N$ such that:
\begin{enumerate}
\item $\sigma$ maps one vertex to another, $\sigma(i)=j$.

\item $\sigma$ leaves invariant the edges, $k-l\iff\sigma(k)-\sigma(l)$.
\end{enumerate}
\end{definition}

In short, what we did here is to copy the statement of Proposition 3.19, with the assumption (2) there removed, and call this a Definition. In view of this, obviously, any circulant graph is transitive. But, do we have other interesting examples?

\bigskip

As a first piece of answer to this question, which is very encouraging, we have:

\begin{theorem}
The cube graph, namely
$$
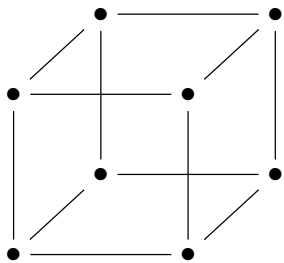
\xymatrix@R=18pt@C=20pt{
&\bullet\ar@{-}[rr]&&\bullet\\
\bullet\ar@{-}[rr]\ar@{-}[ur]&&\bullet\ar@{-}[ur]\\
&\bullet\ar@{-}[rr]\ar@{-}[uu]&&\bullet\ar@{-}[uu]\\
\bullet\ar@{-}[uu]\ar@{-}[ur]\ar@{-}[rr]&&\bullet\ar@{-}[uu]\ar@{-}[ur]
}$$
is transitive, but not circulant.
\end{theorem}

\begin{proof}
The fact that the cube is transitive is clear, because given any two vertices $i,j\in X$, we can certainly rotate the cube in 3D, as to have $i\to j$. As for the fact that the cube is not circulant, this is something more tricky, as follows:

\medskip

(1) As a first observation, when trying to draw the cube on a circle, in a somewhat nice and intuitive way, as to have it circulant, we reach to the following picture:
$$
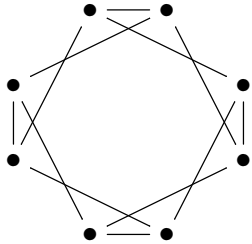
\xymatrix@R=16pt@C=16pt{
&\bullet\ar@{-}[r]\ar@{-}[drr]&\bullet\ar@{-}[dll]\\
\bullet\ar@{-}[d]\ar@{-}[ddr]&&&\bullet\ar@{-}[d]\\
\bullet\ar@{-}[uur]\ar@{-}[drr]&&&\bullet\ar@{-}[uul]\\
&\bullet\ar@{-}[r]\ar@{-}[urr]&\bullet\ar@{-}[uur]}$$

Thus, our cube is indeed not circulant, or at least not in an obvious way. 

\medskip

(2) However, this does not stand for a proof, and the problem of abstractly proving that the cube is not circulant remains. Normally this can be done by attempting to label the vertices in a circulant way. Indeed, up to some discussion here, that we will leave as an instructive exercise, we can always assume that $1,2$ are connected by an edge:
$$1-2$$

(3) But with this in hand, we can now start labeling the vertices of the cube, in a circulant way. Since $1-2$ implies via our circulant graph assumption $2-3$, $3-4$, and so on, in order to start our labeling, we must pick one vertex, and then follow a path on the cube, emanating from there. But, by some obvious symmetry reasons, this means that we can always assume that our first three vertices $1,2,3$ are as follows:
$$
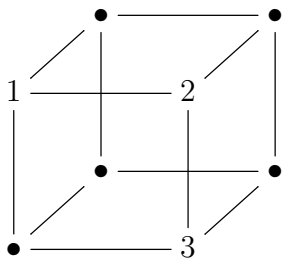
\xymatrix@R=16pt@C=20pt{
&\bullet\ar@{-}[rr]&&\bullet\\
1\ar@{-}[rr]\ar@{-}[ur]&&2\ar@{-}[ur]\\
&\bullet\ar@{-}[rr]\ar@{-}[uu]&&\bullet\ar@{-}[uu]\\
\bullet\ar@{-}[uu]\ar@{-}[ur]\ar@{-}[rr]&&3\ar@{-}[uu]\ar@{-}[ur]
}$$

(4) So, the question comes now, where the vertex 4 can be, as for all this to lead, in the end, to a circulant graph. And the point is that, among the two possible choices for the vertex 4, as new neighbors of 3, none works. Thus, our cube is indeed not circulant, and we will leave the remaining details here as an instructive exercise.
\end{proof}

Getting back to the basics, was our very first result, Theorem 3.16, good or bad news? And here, I don't know about you, but personally I feel quite frustrated by the fact that the computation there leads to $D_N$, instead to $\mathbb Z_N$ itself. So, let us fix this. It is quite obvious that the construction in Theorem 3.15 will work perfectly well for the oriented graphs, or for the colored graphs, and in view of this, let us formulate:

\begin{definition}
Given a colored oriented graph $X$, with vertices denoted $1,\ldots,N$, the symmetries of $X$, which are the permutations $\sigma\in S_N$ leaving invariant the edges, 
$$i-j\implies\sigma(i)-\sigma(j)$$
with their orientations and colors, form a subgroup of the symmetric group
$$G(X)\subset S_N$$
called symmetry group of $X$.
\end{definition}

Here, as before with the construction in Theorem 3.15, the fact that we obtain indeed a group is clear from definitions. Now with this convention in hand, we have:

\begin{theorem}
The symmetry group of the oriented $N$-gon
$$
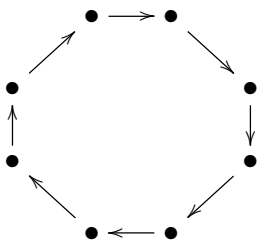
\xymatrix@R=15pt@C=17pt{
&\bullet\ar[r]&\bullet\ar[dr]\\
\bullet\ar[ur]&&&\bullet\ar[d]\\
\bullet\ar[u]&&&\bullet\ar[dl]\\
&\bullet\ar[ul]&\bullet\ar[l]}$$
is the cyclic group $\mathbb Z_N$.
\end{theorem}

\begin{proof}
This is clear from definitions, because once we choose a vertex $i$ and denote its image by $\sigma(i)=i+k$, the permutation $\sigma\in S_N$ leaving invariant the edges, with their orientation, must map $\sigma(i+1)=i+k+1$, $\sigma(i+2)=i+k+2$ and so on, and so must be an element of the cyclic group, in remainder modulo $N$ notation $\sigma=k\in\mathbb Z_N$.
\end{proof}

With this done, and the authority of $\mathbb Z_N$ restored, let us work out some general properties of the construction $X\to G(X)$. In order to do this, we can use:

\begin{theorem}
Having a group action on a colored graph $G\curvearrowright X$ is the same as saying that the action of $G$ leaves invariant the adjacency matrix $d$:
$$d_{ij}=d_{g(i)g(j)}\quad,\quad\forall g\in G$$
Equivalently, the action must preserve the spectral projections of $d$, as follows:
$$d=\sum_\lambda\lambda P_\lambda\implies (P_\lambda)_{ij}=(P_\lambda)_{g(i)g(j)}$$
Moreover, when allowing self-edges, each $P_\lambda$ will correspond to a colored graph $X_\lambda$.
\end{theorem}

\begin{proof}
This is indeed something self-explanatory, with the first assertion coming from definitions, and with the second one being clear too, based on the well-known fact that each real symmetric matrix, and in particular $d$, must be diagonalizable.
\end{proof}

Next, we have the following result, which is again something standard:

\begin{theorem}
Having an action on a colored graph $G\curvearrowright X$ is the same as saying that the action leaves invariant the color components of $X$. Equivalently, with
$$d=\sum_{c\in C}cd_c$$
being the color decomposition of the adjacency matrix, with color components
$$(d_c)_{ij}=\begin{cases}
1&{\rm if}\ d_{ij}=c\\
0&{\rm otherwise}
\end{cases}$$
the action must leave invariant all these color components $d_c$. Thus, the symmetry group $G(X)\subset S_N$ is the subgroup which preserves all these matrices $d_c$.
\end{theorem}

\begin{proof}
This is again something self-explanatory, and elementary.
\end{proof}

The point now is that, by combining Theorems 3.24 and 3.25, both trivialities, we are led to the following enigmatic statement, which all of the sudden wakes us up:

\begin{theorem}
Given an adjacency matrix of a graph $X$, which can be taken in a colored graph sense, $d\in M_N(\mathbb C)$, or even binary as usual, 
$$d\in M_N(0,1)$$
a group action $G\curvearrowright X$ must preserve all ``spectral-color'' components of this matrix, obtained by successively applying the spectral decomposition, and color decomposition.
\end{theorem}

\begin{proof}
This is clear indeed by combining Theorem 3.24 and Theorem 3.25, and with the remark that, indeed, even for a usual binary matrix $d\in M_N(0,1)$ this leads to something non-trivial, because the spectral components of this matrix are no longer binary, and so all of the sudden, we are into colors and everything.
\end{proof}

With the above result in hand, which is something quite unexpected, we are led into a quite interesting linear algebra question, which is surely new for you, namely:

\begin{question}
What are the spectral-color components of a matrix $d\in M_N(\mathbb C)$, or even of a usual binary matrix $d\in M_N(0,1)$?
\end{question}

This question is something non-trivial, and we will be back to it on several occasions, and notably at the end of this book, when talking planar algebras in the sense of Jones \cite{jo4}, which provide the good framework for the study of such questions.

\section*{3d. Derangements}

We still have work to do with permutations, in relation with Question 3.7. That does not look trivial, so let us focus on the simplest possible situation, as follows:

\begin{question}
What is the law of the random variable
$$\chi:S_N\to\mathbb N$$
counting the fixed points, or equivalently, the trace, via $S_N\subset O_N$?
\end{question}

In answer to this, let us first study the permutations $\sigma\in S_N$ having no fixed points at all, called derangements. Regarding them, we have the following remarkable result:

\index{random permutation}
\index{derangement}
\index{fixed points}

\begin{theorem}
The probability for a random $\sigma\in S_N$ to be a derangement is:
$$P=1-\frac{1}{1!}+\frac{1}{2!}-\ldots+(-1)^N\frac{1}{N!}$$
Thus we have the following asymptotic formula, in the $N\to\infty$ limit,
$$P\simeq\frac{1}{e}$$
with $e=2.7182\ldots$ being the usual constant from analysis.
\end{theorem}

\begin{proof}
This is something very classical, which is best viewed by using the inclusion-exclusion principle. Consider indeed the following subsets of $S_N$:
$$S_N^i=\left\{\sigma\in S_N\Big|\sigma(i)=i\right\}$$

By inclusion-exclusion, the probability that we are interested in is given by:
\begin{eqnarray*}
P
&=&\frac{1}{N!}\left(|S_N|-\sum_i|S_N^i|+\sum_{i<j}|S_N^i\cap S_N^j|-\ldots+(-1)^N\sum_{i_1<\ldots<i_N}|S_N^{i_1}\cap\ldots\cap S_N^{i_N}|\right)\\
&=&\frac{1}{N!}\sum_{k=0}^N(-1)^k\sum_{i_1<\ldots<i_k}|S_N^{i_1}\cap\ldots\cap S_N^{i_k}|\\
&=&\frac{1}{N!}\sum_{k=0}^N(-1)^k\sum_{i_1<\ldots<i_k}(N-k)!\\
&=&\frac{1}{N!}\sum_{k=0}^N(-1)^k\binom{N}{k}(N-k)!\\
&=&\sum_{k=0}^N\frac{(-1)^k}{k!}
\end{eqnarray*}

Since at the end we have the standard expansion of $1/e$, we obtain the result.
\end{proof}

More generally now, we have the following result, improving the above:

\begin{theorem}
The probability for a random permutation $\sigma\in S_N$ to have exactly $k$ fixed points, with $k\in\mathbb N$, is given by the following formula:
$$P=\frac{1}{k!}\left(1-\frac{1}{1!}+\frac{1}{2!}-\ldots+(-1)^{N-1}\frac{1}{(N-1)!}+(-1)^N\frac{1}{N!}\right)$$
Thus we have $P\simeq1/(ek!)$, in the $N\to\infty$ limit.
\end{theorem}

\begin{proof}
We already know, from Theorem 3.29, that this formula holds at $k=0$. In the general case now, we have to count the permutations $\sigma\in S_N$ having exactly $k$ points. Since having such a permutation amounts in choosing $k$ points among $1,\ldots,N$, and then permuting the $N-k$ points left, without fixed points allowed, we have:
\begin{eqnarray*}
\#\left\{\sigma\in S_N\Big|\chi(\sigma)=k\right\}
&=&\binom{N}{k}\#\left\{\sigma\in S_{N-k}\Big|\chi(\sigma)=0\right\}\\
&=&\frac{N!}{k!(N-k)!}\#\left\{\sigma\in S_{N-k}\Big|\chi(\sigma)=0\right\}\\
&=&N!\times\frac{1}{k!}\times\frac{\#\left\{\sigma\in S_{N-k}\Big|\chi(\sigma)=0\right\}}{(N-k)!}
\end{eqnarray*}

Now by dividing everything by $N!$, we obtain from this the following formula:
$$\frac{\#\left\{\sigma\in S_N\Big|\chi(\sigma)=k\right\}}{N!}=\frac{1}{k!}\times\frac{\#\left\{\sigma\in S_{N-k}\Big|\chi(\sigma)=0\right\}}{(N-k)!}$$

By using now the computation at $k=0$, that we already have, from Theorem 3.29, it follows that with $N\to\infty$ we have the following estimate:
$$P(\chi=k)
\simeq\frac{1}{k!}\cdot P(\chi=0)
\simeq\frac{1}{k!}\cdot\frac{1}{e}$$

Thus, we are led to the conclusion in the statement.
\end{proof}

In order to interpret now what we found, let us recall the following key definition:

\index{Poisson law}

\begin{definition}
The Poisson law of parameter $1$ is the following measure,
$$p_1=\frac{1}{e}\sum_k\frac{\delta_k}{k!}$$
and the Poisson law of parameter $t>0$ is the following measure,
$$p_t=e^{-t}\sum_k\frac{t^k}{k!}\,\delta_k$$
with the letter ``p'' standing for Poisson.
\end{definition}

We are using here some simplified notations for these laws. Observe that our laws have indeed mass 1, as they should, due to the following key formula:
$$e^t=\sum_k\frac{t^k}{k!}$$

We will see in the moment why these measures appear a bit everywhere, in various discrete probability contexts, the reasons for this coming from the Poisson Limit Theorem (PLT). In order to discuss this, we will need the following key result:

\index{Fourier transform}

\begin{theorem}
The Fourier transform of $p_t$ is given by
$$F_{p_t}(y)=\exp\left((e^{iy}-1)t\right)$$
for any $t>0$. In particular we have $p_s*p_t=p_{s+t}$, for any $s,t>0$.
\end{theorem}

\begin{proof}
We have indeed the following computation, for the Fourier transform:
\begin{eqnarray*}
F_{p_t}(y)
&=&e^{-t}\sum_k\frac{t^k}{k!}F_{\delta_k}(y)\\
&=&e^{-t}\sum_k\frac{t^k}{k!}\,e^{iky}\\
&=&e^{-t}\sum_k\frac{(e^{iy}t)^k}{k!}\\
&=&\exp(-t)\exp(e^{iy}t)\\
&=&\exp\left((e^{iy}-1)t\right)
\end{eqnarray*}

As for the second assertion, this comes from this formula, via the well-known fact that the logarithm of the Fourier transform linearizes the convolution. 
\end{proof}

Good news, we can establish now the Poisson Limit Theorem, as follows:

\index{PLT}
\index{Poisson Limit Theorem}
\index{Bernoulli laws}

\begin{theorem}[PLT]
We have the following convergence, in moments,
$$\left(\left(1-\frac{t}{n}\right)\delta_0+\frac{t}{n}\,\delta_1\right)^{*n}\to p_t$$
for any $t>0$.
\end{theorem}

\begin{proof}
Let us denote by $\nu_n$ the measure under the convolution sign, namely:
$$\nu_n=\left(1-\frac{t}{n}\right)\delta_0+\frac{t}{n}\,\delta_1$$

We have the following computation, for the Fourier transform of the limit: 
\begin{eqnarray*}
F_{\delta_r}(y)=e^{iry}
&\implies&F_{\nu_n}(y)=\left(1-\frac{t}{n}\right)+\frac{t}{n}\,e^{iy}\\
&\implies&F_{\nu_n^{*n}}(y)=\left(\left(1-\frac{t}{n}\right)+\frac{t}{n}\,e^{iy}\right)^n\\
&\implies&F_{\nu_n^{*n}}(y)=\left(1+\frac{(e^{iy}-1)t}{n}\right)^n\\
&\implies&F(y)=\exp\left((e^{iy}-1)t\right)
\end{eqnarray*}

Thus, we obtain indeed the Fourier transform of $p_t$, as desired.
\end{proof}

Now back to the permutation groups, we first have the following result:

\begin{theorem}
The number of fixed points, viewed as random variable,
$$\chi:S_N\to\mathbb N$$
follows the Poisson law $p_1$, in the $N\to\infty$ limit.
\end{theorem}

\begin{proof}
This is indeed a fancy reformulation of what we found before, in Theorem 3.30, by using the probabilistic notions from Definition 3.31.
\end{proof}

As a natural question now, that you might have, can we recover as well the parametric Poisson laws $p_t$ via permutations? In answer, yes, the result being as follows:

\begin{theorem}
Given a number $t\in(0,1]$, the number of fixed points of permutations $\sigma\in S_N$ among $\{1,\ldots,[tN]\}$, viewed as random variable
$$\chi_t:S_N\to\mathbb N$$
follows the Poisson law $p_t$, in the $N\to\infty$ limit.
\end{theorem}

\begin{proof}
As before in the proof of Theorem 3.29, we get by inclusion-exclusion:
\begin{eqnarray*}
P(\chi_t=0)
&=&\frac{1}{N!}\sum_{r=0}^{[tN]}(-1)^r\sum_{k_1<\ldots<k_r<[tN]}|S_N^{k_1}\cap\ldots\cap S_N^{k_r}|\\
&=&\frac{1}{N!}\sum_{r=0}^{[tN]}(-1)^r\binom{[tN]}{r}(N-r)!\\
&=&\sum_{r=0}^{[tN]}\frac{(-1)^r}{r!}\cdot\frac{[tN]!(N-r)!}{N!([tN]-r)!}
\end{eqnarray*}

Now with $N\to\infty$, we obtain from this the following estimate:
$$P(\chi_t=0)
\simeq\sum_{r=0}^{[tN]}\frac{(-1)^r}{r!}\cdot t^r
\simeq e^{-t}$$

More generally, by counting the permutations $\sigma\in S_N$ having exactly $r$ fixed points among $1,\ldots,[tN]$, as in the proof of Theorem 3.30, we obtain:
$$P(\chi_t=r)\simeq\frac{t^r}{r!e^t}$$

Thus, we obtain in the limit a Poisson law of parameter $t$, as stated.
\end{proof}

Many other things can be said, as a continuation of this. We will be back to this, on several occasions, in what follows, and notably in Part IV of the present book.

\section*{3e. Exercises}

This was a quite exciting chapter, we are now already quite deep into finite group theory, and its applications. As exercises on all this, we have:

\begin{exercise}
Try to decompose $S_4$, as a product of smaller groups.
\end{exercise}

\begin{exercise}
Clarify what we said, regarding the dihedral permutation matrices.
\end{exercise}

\begin{exercise}
Learn the theory of the determinant, constructed as signed volume.
\end{exercise}

\begin{exercise}
Find the structure of the alternating group $A_N$, at small values of $N$.
\end{exercise}

\begin{exercise}
Learn more, from probabilists, about the Poisson limits.
\end{exercise}

\begin{exercise}
Study the moments of the Poisson law $p_1$. Can you compute these?
\end{exercise}

\begin{exercise}
More generally, study the moments of $p_t$. Can you compute these?
\end{exercise}

\begin{exercise}
What can you say about derangements belonging to $A_N$?
\end{exercise}

As bonus exercise, experiment some more with the various possible notations for the permutations, and find the one which suits you best. Having this is quite useful.

\chapter{Abelian groups}

\section*{4a. Group actions}

Time to end this introductory Part I, and we still have many things to be discussed, including among others the structure theorem for the finite abelian groups, namely:
$$G=\mathbb Z_{N_1}\times\ldots\times\mathbb Z_{N_k}$$

We will do this slowly, by discussing all sorts of interesting things, and with this abelian group theorem, and its applications, coming towards the end. Getting started now, the key to everything more advanced, in group theory, is the following notion:

\begin{definition}
We say that a group acts on a set, $G\curvearrowright X$, when we have a map
$$G\times X\to X\quad,\quad (g,x)\to g(x)$$
satisfying the conditions $g(h(x))=(gh)(x)$, and $1(x)=x$. 
\end{definition}

Observe that nothing in particular is said about $g^{-1}(x)$, but of course with this point $y=g^{-1}(x)$ being subject to the condition $g(y)=x$. In relation with this, we have:

\begin{proposition}
Given an action $G\curvearrowright X$, the maps
$$x\to g(x)$$
must be bijective, for any $g\in G$.
\end{proposition}

\begin{proof}
The injectivity of the map in the statement comes from:
$$g(x)=g(y)\implies g^{-1}(g(x))=g^{-1}(g(y))\implies x=y$$

As for the surjectivity, this comes from the above observation, namely:
$$x\in X\implies g(g^{-1}(x))=x$$

Thus, we are led to the conclusion in the statement.
\end{proof}

Regarding now the basic examples of actions, all the finite groups that we know well, namely $\mathbb Z_N,D_N,S_N$, come by definition with actions on $\{1,\ldots,N\}$:
$$\mathbb Z_N,\ D_N,\ S_N\curvearrowright\{1,\ldots,N\}$$

In fact, given a subgroup $G\subset S_N$, or just a group morphism $G\to S_N$, we have:
$$G\curvearrowright\{1,\ldots,N\}$$

Now by thinking a bit, it is pretty much clear that this latter example is the generic one. In fact, this is true for infinite sets $X$ too, the result being as follows:

\begin{theorem}
Having an action $G\curvearrowright X$ as axiomatized before, namely
$$G\times X\to X\quad,\quad (g,x)\to g(x)$$
with $g(h(x))=(gh)(x)$, and $1(x)=x$, is the same as having a group morphism
$$G\to S_X\quad,\quad g\to(x\to g(x))$$
with $S_X$ being the group of bijections $\varphi:X\to X$.
\end{theorem}

\begin{proof}
This is indeed something self-explanatory, based on Proposition 4.2:

\medskip

(1) Given an action $G\curvearrowright X$, we know from Proposition 4.2 that the maps $x\to g(x)$ are bijections, so our action is encoded by the following map:
$$G\to S_X\quad,\quad g\to(x\to g(x))$$

But the action axioms, namely $g(h(x))=(gh)(x)$ and $1(x)=x$, show that this latter map $G\to S_X$ must be a group morphism, as desired.

\medskip

(2) Conversely, given a group morphism $G\to S_X$, we can denote it if we want as $g\to(x\to g(x))$, and then define a map in the following way:
$$G\times X\to X\quad,\quad (g,x)\to g(x)$$

Now since $G\to S_X$ was a group morphism, this map is an action, as desired.
\end{proof}

In order to further comment on this, making the link with the material from chapter 3, let us formulate the following definition, in the context of Definition 4.1:

\begin{definition}
Given an action $G\curvearrowright X$, the following conditions are equivalent,
\begin{enumerate}
\item $g\neq h\implies\exists x, g(x)\neq h(x)$.

\item $g\neq 1\implies\exists x, g(x)\neq x$.
\end{enumerate}
and if these conditions are satisfied, we say that $G\curvearrowright X$ is faithful.
\end{definition}

To be more precise, the equivalence between (1) and (2) is indeed clear, coming via $(g^{-1}h)(x)\neq x$, and in what regards the terminology, the ``faithfulness'' in question is with respect to $G$ itself, whose elements $g$ must be reminded by the maps $x\to g(x)$. Alternatively, the corresponding group morphism $G\to S_X$ must be injective.

\bigskip

We can now reformulate the essentials of Theorem 4.3, as follows:

\begin{theorem}
The following happen:
\begin{enumerate}
\item The faithful actions $G\curvearrowright X$ come from the subgroups $G\subset S_X$.

\item The faithful actions $G\curvearrowright\{1,\ldots,N\}$ come from the subgroups $G\subset S_N$.

\item The faithful actions on finite sets come from permutation groups.
\end{enumerate}
\end{theorem}

\begin{proof}
This is basically a reformulation of Theorem 4.3, as follows:

\medskip

(1) This is indeed a reformulation of Theorem 4.3, in the faithful case, and with the remark that we are not losing much when doing so, because an arbitrary action $G\curvearrowright X$, or morphism $G\to S_X$, must appear as follows, with $G\to H$ being a quotient map:
$$G\to H\subset S_X$$

(2) This follows from (1), with $X=\{1,\ldots,N\}$.

\medskip

(3) Assume indeed that we have an action $G\curvearrowright X$, with $|X|=N<\infty$. We can then label $X=\{1,\ldots,N\}$, and we get in this way in the context of (2), as desired.
\end{proof}

With this discussed, you might wonder at this point what is the deal with our notion of action $G\curvearrowright X$, because when restricting the attention to the most interesting case, faithful actions on finite sets, we are back to the permutation groups from chapter 3.

\bigskip

Good point, and in answer, there are many interesting group actions, quite often a bit abstract, which are better understood in terms of the formalism in Definition 4.1, than in terms of associated permutation groups. Here are a few illustrations, for this:

\begin{theorem}
Given a group $G$, we have actions as follows,
\begin{enumerate}
\item An action of $G$ on itself, $G\curvearrowright G$, by translation, $g(h)=gh$. 

\item Another action of $G$ on itself, $G\curvearrowright G$, by conjugation, $g(h)=ghg^{-1}$. 

\item An action $G\curvearrowright G/H$, given by $g(kH)=gkH$, for any subgroup $H\subset G$. 
\end{enumerate}
with $(1)$ producing the Cayley theorem, and with $(2,3)$ remaining to be explored.
\end{theorem}

\begin{proof}
This is indeed something intuitive and elementary, as follows:

\medskip

(1) With $g(h)=gh$ the first action axiom holds indeed, as shown by:
$$g(h(k))=g(hk)=ghk=gh(k)$$

We also have $1(k)=k$, so we have an action. Now with this action being clearly faithful, we obtain from it the Cayley theorem, $G\subset S_G$. Which is very nice.

\medskip

(2) With $g(h)=ghg^{-1}$ the first action axiom holds indeed, as shown by:
$$g(h(k))=g(hkh^{-1})=ghkh^{-1}g^{-1}=gh(k)$$

We also have $1(k)=k$, so we have an action. Finally, with our action being in general not faithful, for instance when $G$ is abelian, we will leave its study for later.

\medskip

(3) With $g(kH)=gkH$ the first action axiom holds indeed, as shown by:
$$g(h(kH))=g(hkH)=ghkH=gh(kH)$$

We also have $1(H)=H$, so we have an action. But with this action being in general not faithful, for instance when $H=G$, we will leave its study for later.
\end{proof}

Observe now that (3) above generalizes (1), which can be obtained with $H=\{1\}$. Which is quite exciting, because based on this, we can declare that the future of the Cayley theorem can only lie in studying group actions, such as those in (2,3).

\bigskip

Summarizing, we have now some good motivations for the study of the group actions $G\curvearrowright X$, as axiomatized in Definition 4.1, and with Theorems 4.3 and 4.5 and permutation groups always in mind. Which is good to know, and doing so will be our plan.

\bigskip

With this in mind, as our next observation regarding the group actions, we have:

\begin{proposition}
Given an action $G\curvearrowright X$, we have an equivalence relation
$$x\sim y\iff \exists g,\, g(x)=y$$
whose orbits are the orbits of the action, namely the sets
$$G(x)=\left\{g(x)\Big|g\in G\right\}$$
which orbits therefore follow to either coincide, or to be disjoint.
\end{proposition}

\begin{proof}
We have to check 3 axioms. The reflexivity axiom $x\sim x$ comes from:
$$1(x)=x$$

The symmetry axiom $x\sim y\implies y\sim x$ is something clear too, coming from:
$$g(x)=y\implies g^{-1}(y)=x$$

As for the transitivity axiom $x\sim y$, $y\sim z$ $\implies$ $x\sim z$, this comes from:
$$g(x)=y,\, h(y)=z\implies (hg)(x)=z$$

Finally, the last assertion, regarding the orbits, is clear from this.
\end{proof}

In order to further study the group actions, we can use the following key notion:

\index{stabilizer}

\begin{proposition}
Given an action $G\curvearrowright X$, and a point $x\in X$, the set
$$G_x=\left\{g\in G\Big|g(x)=x\right\}$$
is a subgroup of $G$, called stabilizer of $x$.
\end{proposition}

\begin{proof}
Again, we have to check 3 axioms. First, we have $1\in G_x$, due to:
$$1(x)=x$$

Next, the fact that we have $g\in G_x\implies g^{-1}\in G_x$ is clear too, coming from:
$$g(x)=x\implies g^{-1}(x)=x$$

Finally, the condition $g,h\in G_x\implies gh\in G_x$ is clear too, coming from:
$$g(x)=x,h(x)=x\implies gh(x)=g(x)=x$$

Thus, we are led to the conclusion in the statement.
\end{proof}

Observe the similarity between Propositions 4.7 and 4.8, and their proofs. In fact, we have the following key result, making the link between orbits and stabilizers:

\index{orbit}

\begin{theorem}
Given an action $G\curvearrowright X$ and a point $x\in X$, we have
$$|G(x)|=[G:G_x]$$
That is, the cardinality of the orbit is the index of the stabilizer.
\end{theorem}

\begin{proof}
In order to prove this, we must construct a bijection, as follows:
$$\varphi:G/G_x\to G(x)$$

But the formula of $\varphi$ can only be something straightforward, as follows:
$$\varphi(gG_x)=g(x)$$

So, let us see if this works. To start with, $\varphi$ is well-defined and injective, due to:
\begin{eqnarray*}
gG_x=hG_x
&\iff&g^{-1}h\in G_x\\
&\iff&g^{-1}h(x)=x\\
&\iff&g(x)=h(x)
\end{eqnarray*}

But $\varphi$ is clearly surjective too, and we therefore obtain the result.
\end{proof}

As a first illustration, for the action $G\curvearrowright G$ by translation, $g(h)=gh$, from Theorem 4.6 (1), the orbits are $G(h)=G$, the stabilizers are $G_h=\{1\}$, and we obtain:
$$|G|=[G:\{1\}]$$

Nevermind. As a second illustration for this, for the action $G\curvearrowright G$ by conjugation, $g(h)=ghg^{-1}$, from Theorem 4.6 (2), the orbits are the conjugacy classes:
$$G(h)=\left\{ ghg^{-1}\Big|g\in G\right\}$$

As for the stabilizers, these are the following subgroups, called centralizers:
$$G_h=\left\{g\in G\Big|ghg^{-1}=h\right\}$$

Which sounds quite interesting, the above orbits and stabilizers being something non-trivial, so let us record the application of Theorem 4.9 in this case as a theorem:

\index{conjugacy class}
\index{centralizer}

\begin{theorem}
We have the following formula, for any $h\in G$,
$$|G(h)|=[G:G_h]$$
telling us that the cardinality of the conjugacy class is the index of the centralizer.
\end{theorem}

\begin{proof}
This follows indeed from Theorem 4.9, and from the above discussion.
\end{proof}

But you might wonder, in relation with this, what is conjugacy good for? In answer, conjugacy can help in connection with many advanced questions, as shown by:

\index{central function}
\index{group algebra}
\index{character}

\begin{theorem}
Given a finite group $G$, the following happen:
\begin{enumerate}
\item The central functions, $f(gh)=f(hg)$, are those constant on conjugacy classes.

\item Any group morphism, or character, $\chi:G\to\mathbb T$ is such a central function.

\item In fact, the character $\chi_u(g)=Tr(u(g))$ of any morphism $u:G\to U_N$ is central.

\item Also, the center of the group algebra $\mathbb C[G]$ consists of the central functions.
\end{enumerate}
\end{theorem}

\begin{proof}
Many things going on here, with all this being quite advanced, and with the idea being as follows, coming with the promise to come back to all this, later:

\medskip

(1) We have indeed the following equivalence, coming by setting $k=hg$:
$$f(gh)=f(hg)\iff f(gkg^{-1})=f(k)$$

(2) This is clear indeed from the fact that the target group $\mathbb T$ is abelian:
$$\chi(gh)=\chi(g)\chi(h)=\chi(h)\chi(g)=\chi(hg)$$

(3) This assertion, which generalizes (1), can be proved as follows:
$$\chi_u(gh)=Tr(u(g)u(h))=Tr(u(h)u(g))=\chi_u(hg)$$

\medskip

(4) Let us define indeed the group algebra of our group $G$ as being the following complex vector space, with multiplication given by $g\cdot h=gh$ and linearity:
$$\mathbb C[G]=\left\{\sum_g\lambda_gg\Big|\lambda_g\in\mathbb C\right\}$$

Getting now to the computation of the center, this goes as follows:
\begin{eqnarray*}
\lambda\in Z(\mathbb C[G])
&\iff&\lambda h=h\lambda,\, \forall h\\
&\iff&\sum_g\lambda_ggh=\sum_g\lambda_ghg\\
&\iff&\sum_k\lambda_{kh^{-1}}k=\sum_k\lambda_{h^{-1}k}k\\
&\iff&\lambda_{kh^{-1}}=\lambda_{h^{-1}k}
\end{eqnarray*}

Thus, our algebra element $\lambda$ must be a central function on $G$, as stated.
\end{proof}

Moving on, as our last trick, we still have left applying Theorem 4.9 to the situation in Theorem 4.6 (3). So, given a subgroup $H\subset G$, consider the following action:
$$G\curvearrowright G/H\quad,\quad g(kH)=gkH$$

But in this case, a bit as before in the context of Theorem 4.6 (1), which corresponds to the case $H=\{1\}$, the orbits are $G(kH)=G/H$, independently of $k\in G$, the stabilizers are $G_{kH}=H$, also independently of $k\in G$, and the formula in Theorem 4.9 reads:
$$|G/H|=[G:H]$$

Nevermind. In order to recycle however this study, or at least part of it, and get something out of this, let us formulate, in the general context of Definition 4.1:

\begin{definition}
A group action $G\curvearrowright X$ is called transitive when:
$$\forall x,y\in X,\ \exists g\in G,\ g(x)=y$$
Equivalently, our action must have only $1$ orbit, namely the whole $X$.
\end{definition}

As a first observation, we already met such things in chapter 3, when discussing the transitive graphs, which require the symmetry group to be transitive, in the above sense. More on this later, and in the meantime, let us formulate the following result, making the link with Theorem 4.9, and which is something very useful, in practice:

\index{transitive action}

\begin{theorem}
For a transitive group action $G\curvearrowright X$ we have
$$X\simeq G/H$$
with $H=G_x$ being the stabilizer of any point $x\in X$.
\end{theorem}

\begin{proof}
This follows indeed from Theorem 4.9, or rather from its proof, but in view of future applications, let us have this worked out independently, with full details. Our claim is that the identification in the statement can be constructed as follows:
$$\varphi:G/H\to X\quad,\quad \varphi(gH)=g(x)$$

So, let us see if this works. To start with, $\varphi$ is well-defined and injective, due to:
\begin{eqnarray*}
gH=kH
&\iff&g^{-1}k\in H\\
&\iff&g^{-1}k(x)=x\\
&\iff&g(x)=k(x)
\end{eqnarray*}

But $\varphi$ is surjective too, by transitivity, and this gives the result.
\end{proof}

Finally, still in relation with transitivity, at a more advanced level, we have:

\index{orbital}
\index{higher orbital}

\begin{theorem}
Given a subgroup $G\subset S_N$ we can talk about
\begin{enumerate}
\item Orbits, coming via $x\sim y$ when $\exists g$, $g(x)=y$.

\item Orbitals, coming via $(x,z)\sim(y,t)$ when $\exists g$, $g(x)=y$, $g(z)=t$.

\item Higher orbitals, coming via $(x_i)\sim(y_i)$ when $\exists g$, $g(x_i)=y_i$, $\forall i$.
\end{enumerate}
and the corresponding notions of transitivity refine each other.
\end{theorem}

\begin{proof}
This is something quite self-explanatory, with the fact that $(x_i)\sim(y_i)$ as in the statement is indeed an equivalence relation coming from definitions. Alternatively, we can say that these equivalence relations come from the following group actions:
$$G\curvearrowright\{1,\ldots,N\}^k\quad,\quad g(x_1,\ldots, x_k)=(g(x_1),\ldots,g(x_k))$$

As for the last assertion, observe for instance that $\mathbb Z_N$ is transitive, but not doubly transitive. Many other things can be said here, ranging from elementary to fairly advanced, and we will leave some exploration here as an instructive exercise.
\end{proof}

\section*{4b. Sylow theorems}

Getting now towards the finite abelian groups, remember that we want to prove that any such group decomposes as a product of cyclic groups, as follows:
$$G=\mathbb Z_{N_1}\times\ldots\times\mathbb Z_{N_k}$$

There are many ways of proving this theorem, with some of them being rather quick and abstract, while some other being slower, but more useful in practice. We will follow here a slower, but quite efficient way. Our guiding principle will be as follows:

\index{local-global principle}

\begin{principle}
The finite abelian groups $G$ are best understood by splitting them as direct products, as follows, with $p$ varying over the prime numbers,
$$G=\prod_pG_p$$
and then working out the cyclic group decomposition of each of the components,
$$G_p=\mathbb Z_{p^{s_1}}\times\ldots\times\mathbb Z_{p^{s_r}}$$ 
as to reach in the end the conclusion that $G$ is a direct product of cyclic groups, and with each such cyclic component having as order a prime power, $q=p^s$.
\end{principle}

In short, we would like to use an approach inspired by arithmetic, and more specifically, by the usual factorization of the integers. In case you are a bit familiar with more advanced arithmetic, you will recognize here the mighty ``local-global'' principle from there.

\bigskip

But probably enough talking, let us get to work. Inspired by the above, let us formulate the following definition, which is something of independent interest, going beyond the abelian case, and that will keep us busy for the remainder of this section:

\index{Sylow subgroup}
\index{p-group}

\begin{definition}
Let $p\geq2$ be a prime number.
\begin{enumerate}
\item A $p$-group is a group $H$ satisfying $|H|=p^s$, with $s\geq1$.

\item A $p$-subgroup is a subgroup $H\subset G$ which is a $p$-group, $|H|=p^s$.

\item When $|G|=p^sn$ with $(n,p)=1$, we say that $H\subset G$ is a $p$-Sylow subgroup.
\end{enumerate}
\end{definition}

Our goal in this section will be that of proving the Sylow theorems, stating that any finite group $G$, with $p|\,|G|$, has $p$-Sylow subgroups, which moreover are conjugate, and with their number being $=1(p)$. Which is something very useful in group theory, especially when it comes to study an abstract group $G$, based on the knowledge of $N=|G|$.

\bigskip

As an illustration, assuming that $G$ is abelian, with $p|\,|G|$, it follows either from definitions, or from the above-mentioned Sylow theorems, that $G$ has exactly one $p$-Sylow subgroup, namely the group $G_p\subset G$ from Principle 4.15. More on this later.

\bigskip

In practice now, the Sylow theorems are something quite tricky, whose proof takes some time. Let us start our study with the following key result, due to Cauchy:

\index{Cauchy theorem}

\begin{theorem}[Cauchy]
Given a finite group $G$, and a prime number satisfying
$$p\ |\ |G|$$
$G$ has an element of order $p$. Equivalently, $G$ has a subgroup of order $p$.
\end{theorem}

\begin{proof}
This is something quite tricky, the idea being as follows:

\medskip

(1) We must find $g\neq1$ with $g^p=1$. In order to do so, consider the following set:
$$X=\left\{(g_1,\ldots,g_p)\in G^p\Big|g_1\ldots g_p=1\right\}$$

We have then an obvious action $\mathbb Z_p\curvearrowright X$, by rotation, as follows:
$$k(g_1,\ldots,g_p)=(g_{k+1},\ldots,g_{k+p})$$

Now let us decompose $X$ into orbits. This gives the following formula, with $F\subset X$ being the fixed points, and with the sum being over the non-trivial orbits $O$:
$$|X|=|F|+\sum_{|O|\geq2}|O|$$

(2) Next, let us look at this equality modulo $p$. To start with, we have:
$$|X|=|G|^{p-1}=0(p)$$

Also, in what regards the fixed points, we can say here that we have:
$$(1,\ldots,1)\in F\implies |F|\geq1$$

Finally, by Theorem 4.9 the size of any orbit must divide $|\mathbb Z_p|=p$, and so:
$$|O|\geq2\implies |O|=p$$

Now by putting everything together, modulo our $p\geq2$, we conclude that:
$$|F|\geq2$$

(3) But this is exactly what we need, because the fixed points are precisely the elements $(g,\ldots,g)\in G^p$ with $g^p=1$. Thus, we have found $g\neq1$ with $g^p=1$, as desired.

\medskip

(4) Finally, the last assertion is clear, because an element $g\in G$ of order $p$ must generate a subgroup $<g>=\mathbb Z_p$, and also, conversely, a subgroup $H\subset G$ of order $p$ must be of the form $H=\mathbb Z_p$, with its standard generator having order $p$. 
\end{proof}

In relation now with the various notions from Definition 4.18, we have:

\begin{remark}
The Cauchy theorem can be further interpreted as saying that, under the assumption $p|\,|G|$ there, the group $G$ must have a $p$-subgroup.
\end{remark}

To be more precise, in one sense this is clear, because the subgroup $\mathbb Z_p\subset G$ coming from Cauchy is a $p$-subgroup. As for the converse, assuming that $H\subset G$ is a $p$-subgroup, if we pick any $g\in H$, $g\neq1$ then we have the following implication, proving Cauchy:
$$ord(g)=p^r\implies ord(g^{p^{r-1}})=p$$

Getting now to the Sylow theorems, these are a continuation of the Cauchy theorem. We first have the following result, which improves the Cauchy theorem:

\index{Sylow theorem}

\begin{theorem}[Sylow 1]
Given a finite group $G$, and a prime number satisfying
$$p\ |\ |G|$$
$G$ has a $p$-Sylow subgroup.
\end{theorem}

\begin{proof}
We use the same idea as for Cauchy, namely a suitable group action:
\medskip

(1) Let us write $|G|=p^sn$ with $(n,p)=1$, and consider the following set:
$$X=\left\{S\subset G\Big||S|=p^s\right\}$$

We have then an action $G\curvearrowright X$ by left multiplication, as follows:
$$g(S)=gS$$

Now let us decompose $X$ into orbits. In view of Theorem 4.9, we can write the resulting counting formula as follows, with the sum being over points in different orbits:
$$|X|=\sum_S[G:G_S]$$

(2) Next, we would like to look at this equality, modulo $p$. For this purpose, observe first that the cardinality of the set $X$ is given by the following formula:
$$|X|=\binom{p^sn}{p^s}=\frac{p^sn(p^sn-1)\ldots(p^sn-p^s+1)}{p^s\ldots1}$$

Now let us compute the exponent of $p$ in the numerator. This exponent appears from the $p^s$ multiples of $p$ which are present, followed by the $p^{s-1}$ multiples of $p^2$ which are present, which must be counted again, followed by the $p^{s-2}$ multiples of $p^3$ which are present, which must be counted again too, and so on, up to the $1$ multiple of $p^s$ which is present. In other words, this exponent is given by the following formula:
$$e=p^s+p^{s-1}+\ldots+p+1$$

Now since this exponent does not depend on $n$, when looking at the denonimator, which corresponds to $n=1$, we will obtain exactly the same exponent. Thus, the various $p$ factors in our fraction will cancel each other, and we conclude that we have:
$$|X|\neq0(p)$$

(3) But with this, we can go back to the counting formula in (1), with the conclusion that we have at least one point $S\in X$ such that the following happens:
$$[G:G_S]\neq0(p)$$ 

Our claim now is that $G_S\subset G$ is the desired Sylow subgroup. Indeed:

\medskip

-- From $|G|=p^sn$ and $[G:G_S]\neq0(p)$ we deduce that we have $p^s|\,|G_s|$, and so:
$$p^s\leq|G_S|$$

\medskip

-- On the other hand, $G_S=\{g\in G|gS=S\}$, so if we pick $h\in S$, we obtain:
$$|G_S|=|G_Sh|\leq|S|=p^s$$

We conclude that $|G_S|=p^s$, so we have our Sylow subgroup, as desired.
\end{proof}

The second Sylow theorem, which actually improves the first one, is as follows:

\index{Sylow theorem}

\begin{theorem}[Sylow 2]
Given a finite group $G$, and a prime number satisfying $p|\,|G|$, the various $p$-Sylow subgroups $K\subset G$ are as follows:
\begin{enumerate}
\item Any $p$-subgroup $H\subset G$ is contained in a $p$-Sylow subgroup.

\item These $p$-Sylow subgroups $K\subset G$ are pairwise conjugate.
\end{enumerate}
\end{theorem}

\begin{proof}
This is something more technical, the idea being as follows:

\medskip

(1) To start with, the first assertion generalizes Sylow 1, and does not seem to be related much to the second assertion. However, this is wrong, because when it comes to proofs, this first assertion appears as a ``bonus'' to the proof of the second assertion.

\medskip

(2) In practice now, we know from Sylow 1 that we have a Sylow subgroup $K\subset G$. Our claim, which will prove both the assertions in the statement, is that if $H\subset G$ is a $p$-subgroup, then we can find an element $g\in G$ such that the following happens:
$$H\subset gKg^{-1}$$

Observe that this would prove indeed the first assertion, with $gKg^{-1}\subset G$ being obviously a Sylow subgroup, due to $|gKg^{-1}|=|K|$, and would prove the second assertion as well, simply by taking the $p$-subgroup $H\subset G$ to be an arbitrary Sylow subgroup.

\medskip

(3) Getting to work now, in order to prove our claim, on the menu, group actions as usual. So, given a Sylow subgroup $K\subset G$ as above, consider the following set:
$$G/K=\left\{gK\Big|g\in G\right\}$$

Next, given $H\subset G$ as above, we have an action $H\curvearrowright G/K$ by multiplication:
$$h(gK)=hgK$$

(4) Now let us decompose $G/K$ into orbits, and write the counting formula as follows:
$$|G/K|=\sum_O|O|$$

Since $K\subset G$ was assumed Sylow, $[G:K]\neq0(p)$, so we can find $O$ such that:
$$|O|\neq0(p)$$

(5) So, let us try to understand what $|O|\neq0(p)$ means. Since $O$ was an orbit coming from an action of a $p$-group, namely $H\subset G$, by Theorem 4.9 we conclude that $|O|$ must be a power of $p$. Thus, the condition $|O|\neq0(p)$ found above simply reads:
$$|O|=1$$

(6) Next, let us see what $|O|=1$ tells us. We must have, for a certain $g\in G$:
$$O=gK$$

And the claim is that this $g\in G$ is exactly the one we are looking for, realising:
$$H\subset gKg^{-1}$$

Indeed, our orbit $O=gK$ being a fixed point under the action of $H$, we must have $HO=O$, which amounts in saying that the following happens:
$$h\in H\implies hg\in gK$$

But the condition on the right reads $h\in gKg^{-1}$, so $H\subset gKg^{-1}$, as desired.
\end{proof}

As a next question now, which is important in relation with concrete applications, as we will soon discover, we would like to count the $p$-Sylow subgroups $K\subset G$. The answer to this question comes from the third Sylow theorem, which is as follows:

\index{Sylow theorem}

\begin{theorem}[Sylow 3]
Given a finite group $G$, and a prime number satisfying $p|\,|G|$, the number $m_p$ of $p$-Sylow subgroups $K\subset G$ satisfies:
\begin{enumerate}
\item $m_p|\,|G|$.

\item $m_p=1(p)$.

\item $m_p=[G:N_K]$, with $N_K=\{g\in G|gK=Kg\}$.
\end{enumerate}
\end{theorem}

\begin{proof}
This is something a bit more technical, the idea being as follows:

\medskip

(1) To start with, the last assertion makes reference to the notion of normalizer of a subgroup $K\subset G$, given by the formula above. Observe that $N_K$ is indeed a group, and in addition $K\subset N_K$ is normal, and $N_K\subset G$ is maximal with this property. Many other things can be said about normalizers, but for our purposes here, this will do.

\medskip

(2) Getting now to the proof, this will come, as usual, via the study of some suitable group actions. Consider the set formed by all $p$-Sylow subgroups $K_i\subset G$:
$$X=\{K_1,\ldots,K_m\}$$

We have then a group action $G\curvearrowright X$ by conjugation, as follows:
$$g(K_i)=gK_ig^{-1}$$

By Sylow 2 this action is transitive, so if we pick any $K\in X$, we have:
$$m=[G:G_K]$$

But the stabilizer being $G_K=\{g\in G|gKg^{-1}=K\}$, this is the same as the normalizer $N_K=\{g\in G|gK=Kg\}$ in the statement, and so we have proved the third assertion:
$$m=[G:N_K]$$

Observe that the first assertion, namely $m|\, |G|$, follows as well, from this.

\medskip

(3) Summarizing, it remains to prove $m=1(p)$. For this purpose, we will use the same action as above, on the set $X=\{K_i\}$, but restricted to a certain $K\in X$. That is, we pick a Sylow subgroup $K\in X$, and consider the action $K\curvearrowright X$ by conjugation:
$$g(K_i)=gK_ig^{-1}$$

Now let us decompose $X$ into orbits, and write the counting formula as follows, with $F\subset X$ being the fixed points, and with the sum being over the non-trivial orbits $O$:
$$|X|=|F|+\sum_{|O|\geq2}|O|$$

Our claim is that the formula $m=1(p)$ will come from this, by analysis modulo $p$.

\medskip

(4) In order to prove this, let us look at the fixed points $F\subset X$. The point $K\in X$ is certainly fixed, $K\in F$, and our claim is that this is the only fixed point:
$$|F|=1$$

Indeed, assume that we have a second fixed point, $L\in X$. The fixed point condition means $gLg^{-1}=L$ for any $g\in K$, which is the same as saying that we have:
$$K\subset N_L$$

In order to prove $K=L$, pick $g\in K$. We have then $g\in N_L$, so we can consider the element $gL\in N_L/L$. Now observe that the order of this latter element satisfies:
$$ord(gL)|ord(g)$$

On the other hand, we have the following formulae, both coming from definitions:
$$ord(gL)|\ |N_L/L|\neq0(p)$$
$$ord(g)|\ |K|=p^s$$

We conclude that $ord(gL)=1$, so $gL=L$, which gives $g\in L$, as desired.

\medskip

(5) Good news, we can now finish the proof. Indeed, the formula in (3) becomes:
$$|X|=1+\sum_{|O|\geq2}|O|$$

Now since by Theorem 4.9 the size of each orbit $O$ must appear as the index of certain stabilizer $G_L\subset K$, we conclude that we have the following implication:
$$|O|\geq 2\implies |O|=0(p)$$

Thus, our counting formula above reads $|X|=1(p)$, as desired.
\end{proof}

Still with me I hope, after all this algebra. Getting now to the applications of this, there are many of them, with the general principle being as follows:

\begin{principle}
When faced with an abstract group $G$, having order
$$|G|=p_1^{a_1}\ldots p_k^{a_k}$$
with $r=\sum_ia_i$ small, the Sylow theorems are what you need.
\end{principle}

So, let us explore now a bit this general principle, which is something fundamental, in advanced group theory. The situation here, quickly examined, is as follows:

\bigskip

-- At $r=1$, to start with, we know that each group having prime order, $|G|=p$, must be a cyclic group, $G=\mathbb Z_p$, and with this being something elementary.

\medskip

-- At $r=2$ now, we first have the situation $|G|=p^2$, with $p$ prime. But here, it is not difficult to show that $G$ must be abelian, and we will leave this as an exercise.

\medskip

-- Still at $r=2$, next we have the situation $|G|=2p$, with $p>2$ prime. But here, using Cauchy, the solutions are $G=\mathbb Z_{2p}$ and $G=D_p$, again exercise for you.

\medskip

-- And then, still at $r=2$, we have the remaining situation, $|G|=pq$ with $p>q>2$ primes. And here, many things can happen, and you need the Sylow theorems.

\bigskip

So, these are the first few observations that can be made, and as a target for a first application of the Sylow theorems, we have the case $|G|=15$. And here, we have:

\begin{theorem}
Any group of order $15$ must be cyclic, $G\simeq\mathbb Z_{15}$.
\end{theorem}

\begin{proof}
Assume that $G$ has order $|G|=15$. Then, the following happen:

\medskip

(1) According to the Sylow 3 theorem, the number $m_3$ of 3-Sylow subgroups, and the number $m_5$ of 5-Sylow subgroups, can be computed as follows:
$$m_3|5,\ m_3=1(3)\implies m_3=1$$
$$m_5|3,\ m_5=1(5)\implies m_5=1$$

(2) Thus, we have a unique 3-Sylow subgroup $H\subset G$, which is the same as saying that we have a unique copy $\mathbb Z_3\subset G$, which in addition, by Sylow 2, must be normal. Similarly, we have a unique 5-Sylow subgroup $K\subset G$, which is the same as saying that we have a unique copy $\mathbb Z_5\subset G$, which in addition, by Sylow 2, must be normal.

\medskip

(3) Summarizing, we have proved that $G$ has normal subgroups $\mathbb Z_3\subset G$ and $\mathbb Z_5\subset G$. The point now is that, by normality, we get a group morphism as follows:
$$\mathbb Z_5\to Aut(\mathbb Z_3)\quad,\quad g\to(h\to ghg^{-1})$$

But $Aut(\mathbb Z_3)=\mathbb Z_2$, so the above group morphism must be trivial. In other words, our subgroups $\mathbb Z_3\subset G$ and $\mathbb Z_5\subset G$ commute, and it follows that these subgroups generate inside $G$ a copy of $\mathbb Z_3\times\mathbb Z_5=\mathbb Z_{15}$. Thus, we obtain $G=\mathbb Z_{15}$, as desired.
\end{proof}

The above result looks quite encouraging, suggesting that the study of the general case $|G|=pq$ with $p>q>2$ primes might be similar. However, this is not exactly the case, because a close inspection of the above proof shows that this relies on some special arithmetic properties of the numbers $5>3$ used there as input, exploited in both steps (1,3) there. We will leave some further exploration here as an interesting exercise.

\section*{4c. Abelian groups} 

With the above ingredients in hand, we can go back to the abelian groups. We recall that we want to prove that the finite abelian groups appear as products of cyclic groups, as follows, and in addition, with the orders $N_i$ being prime powers:
$$G=\mathbb Z_{N_1}\times\ldots\times\mathbb Z_{N_k}$$

Let us start with the following result, making the link with Sylow subgroups:

\begin{theorem}
Given a finite abelian group $G$, and $p|\,|G|$, the set
$$G_p=\left\{g\in G\Big|\exists k\in\mathbb N,g^{p^k}=1\right\}$$
is a subgroup. In fact, this is the unique $p$-Sylow subgroup.
\end{theorem}

\begin{proof}
This is something quite elementary, not necessarily requiring our knowledge of Sylow theory, with the pedestrian proof for this being as follows:

\medskip

(1) To start with, the fact that the set in the statement $G_p\subset G$ is a subgroup is clear, coming from the following computation, valid inside any abelian group:
$$g^a=1,h^b=1\implies (gh)^{ab}=g^ah^b=1$$

Indeed, given two elements $g,h\in G$, having as orders powers of $p$, this computation shows that  $gh\in G$ has as order a certain power of $p$ too, as desired.

\medskip

(2) Next, observe that, by the Cauchy theorem, $G_p$ must be a $p$-group. With the remark that we are only using here Cauchy for abelian groups, whose proof is something quite elementary, simpler than that of Cauchy in general, say exercise for you.

\medskip

(3) Our next claim is that $G_p\subset G$ is the maximal $p$-subgroup. Indeed, given a $p$-subgroup $H\subset G$, and $h\in H$, from $ord(h)|\,|H|$ we get $h\in G_p$, as desired.

\medskip

(4) Finally, assuming $|G|=p^kn$ with $(n,p)=1$, we must prove the Sylow condition $|G_p|=p^k$. But this is best seen by contradiction. Indeed, assuming $p|\,|G/G_p|$, by Cauchy we would have a certain non-trivial element $hG_p\in G/G_p$ of order $p$. But this means $h\notin G_p$, $h^p\in G_p$, which in turn reads $h\notin G_p$, $h\in G_p$, which is contradictory.

\medskip

(5) Summarizing, we managed to prove our theorem, based only on the Cauchy theorem for abelian groups, which is something elementary. Alternatively, this follows from the Sylow theorems too. Indeed, in the abelian case, these theorems show that we have a unique $p$-Sylow subgroup $H\subset G$, which by (3) follows to be $H=G_p$, as desired.
\end{proof}

As a continuation of this, we have the following key result:

\begin{theorem}
Given a finite abelian group $G$, we have:
$$G=\prod_pG_p$$
That is, $G$ appears as the product of its Sylow subgroups.
\end{theorem}

\begin{proof}
By using the fact that our group $G$ is abelian, we have a group morphism as follows, with the order of the factors when computing $\prod_pg_p$ being irrelevant:
$$\prod_pG_p\to G\quad,\quad (g_p)\to\prod_pg_p$$

(1) Our first claim is that this morphism is injective. Indeed, let us consider an element in its kernel, which amounts in having an equation of the following type:
$$g_1\ldots g_k=1$$

Now since the elements $g_1$ and $g_2\ldots g_k$, which are inverse to each other, must have the same order, and the order of $g_1$ is a certain prime power, and that of $g_2\ldots g_k$ is not divisible by that prime, we conclude that the kernel is trivial, as claimed. 

\medskip

(2) It remains to prove that our morphism is surjective. Normally this can be done in the pedestrian way, by picking $g\in G$, writing its order as $ord(g)=p_1^{a_1}\ldots p_k^{a_k}$, and doing some arithmetic in order to reach to a writing of type $g=g_1\ldots g_k$, with $g_i\in G_{p_i}$.

\medskip

(3) However, there is better than this, using our previous knowledge of Sylow subgroups. Indeed, as mentioned in Theorem 4.24, the components $G_p\subset G$ are the various Sylow subgroups of $G$. Now by using this, we see right away that we have:
$$|G|=\prod_p|G_p|$$

Thus, by getting back now to our group morphism above, this is injective, and the cardinalities of its domain and range match, so this is an isomorphism, as desired.
\end{proof}

Getting now to what we wanted to do, in this section, namely structure theorem for the abelian groups, Theorem 4.25 does half of the job. For the other half, we must decompose the various components $G_p$, and this can be done as follows:

\index{abelian p-group}

\begin{theorem}
The abelian $p$-groups decompose as follows:
$$G=\mathbb Z_{p^{r_1}}\times\ldots\times\mathbb Z_{p^{r_s}}$$
That is, the abelian $p$-groups are the products of cyclic $p$-groups.
\end{theorem}

\begin{proof}
We can do this by recurrence on $|G|$, as follows:

\medskip

(1) Let us pick $g\in G$ of maximal order, say $ord(g)=p^k$, and consider the subgroup $H=<g>$ that it generates, inside $G$. By recurrence, the quotient group $G/H$ must decompose as follows, with the components $C_i$ being cyclic groups:
$$G/H=C_1\times\ldots\times C_n$$

Our goal will be that of producing, out of this, an isomorphism as follows:
$$G=H\times C_1\times\ldots\times C_n$$

(2) Let us start by fixing some notation. The subgroups $C_i\subset G/H$ appearing above being cyclic, we can denote them as $C_i=\{z_i^aH\}$, with $z_iH\in C_i$ being some chosen generators for them. And with this, the isomorphism that we have is:
$$\varphi:C_1\times\ldots\times C_n\to G/H\quad,\quad 
(z_1^{a_1}H,\ldots,z_n^{a_n}H)\to z_1^{a_1}\ldots z_n^{a_n}H$$

Our more precise claim now, which will prove the result, is that, with a suitable choice of the generators $z_iH\in C_i$, we can lift this into an isomorphism as follows:
$$\psi:H\times C_1\times\ldots\times C_n\to G\quad,\quad 
(g^a,z_1^{a_1}H,\ldots,z_n^{a_n}H)\to g^az_1^{a_1}\ldots z_n^{a_n}$$

(3) In order to do this, let us look at one of the components, $C=C_i$. If we pick an arbitrary generator $zH\in C$, with $z\in G$, the following happens, trivially:
$$ord(zH)|ord(z)$$

And our claim now, which will provide us with what is needed in (2), is that we can always arrange for our generator $zH\in C$, with $z\in G$, as to have equality:
$$ord(zH)=ord(z)$$

(4) Summarizing, we have eventually found something concrete to prove, in relation with what we want to do, so let us prove this. Let us start with an arbitrary generator $xH\in C$, with $x\in G$. Consider the two orders mentioned in (3), namely:
$$p^r=ord(xH)\quad,\quad p^s=ord(x)\quad,\quad r\leq s$$

Our goal will be that of suitably modifying our generator $xH$, as to have $r=s$.

\medskip

(5) In order to do so, let us look at the following group element $y\in G$:
$$y=x^{p^r}\quad,\quad ord(y)=p^{s-r}$$

Since $ord(xH)=p^r$ we have $ord(yH)=1$, which means $y\in H$. Now since $H=<g>$ was the group generated by $g$, we can write $y$ as follows, with $(n,p)=1$:
$$y=g^{np^t}$$

Now recall that $g\in G$ was chosen of maximal order $p^k$. Thus, we have:
$$ord(y)=p^{k-t}$$

We conclude that we have $s-r=k-t$. Now consider the following element:
$$z=xg^{-np^{t-r}}$$

Our claim is that this is the element $z\in G$ that we were looking for, in (3).

\medskip

(6) Indeed, we first have the following computation, which gives $ord(z)\leq p^r$:
$$z^{p^r}=x^{p^r}g^{-np^t}=y\cdot y^{-1}=1$$

Also, $zH=xH=C$, and so $ord(zH)=|C|=p^r$. Thus we have, as desired:
$$ord(zH)=ord(z)=p^r$$

(7) Time for the endgame. Let us go back to the isomorphism in (2), which was as follows, and with the generators $z_iH\in C_i$ with $z_i\in G$ being chosen as above:
$$\varphi:C_1\times\ldots\times C_n\to G/H\quad,\quad 
(z_1^{a_1}H,\ldots,z_n^{a_n}H)\to z_1^{a_1}\ldots z_n^{a_n}H$$

Our claim is that this lifts into an isomorphism as follows:
$$\psi:H\times C_1\times\ldots\times C_n\to G\quad,\quad 
(g^a,z_1^{a_1}H,\ldots,z_n^{a_n}H)\to g^az_1^{a_1}\ldots z_n^{a_n}$$

(8) Indeed, this latter map is well-defined, due to $ord(z_iH)=ord(z_i)$. It is also clear that $\psi$ is a group morphism. Also, since $\varphi$ is surjective, so must be $\psi$. Finally, since the cardinalities of the domain and range match, $\psi$ must be an isomorphism, as desired.
\end{proof}

Time now to put everything together. We obtain the following remarkable result:

\index{finite abelian group}

\begin{theorem}
The finite abelian groups are the products of cyclic groups:
$$G=\mathbb Z_{N_1}\times\ldots\times\mathbb Z_{N_k}$$
Moreover, we can choose the numbers $N_i$ to be prime powers.
\end{theorem}

\begin{proof}
This follows indeed by putting together all the above, and more specifically, by combining Theorem 4.25 and Theorem 4.26. As further remarks on this:

\medskip

(1) Observe that this generalizes both Theorem 4.25 and Theorem 4.26. Thus, in practice, this theorem is all you need to know, regarding the finite abelian groups. 

\medskip

(2) In fact, what you need to know is just the first part of the present theorem, because the second part is easy to recover, thanks to the following elementary isomorphisms:
$$\mathbb Z_{p_1^{a_1}\ldots p_k^{a_k}}=\mathbb Z_{p_1^{a_1}}\times\ldots\times\mathbb Z_{p_k^{a_k}}$$

(3) Finally, there is a uniqueness assertion too, stating that with $G$ fully split, with $N_i$ prime powers, the cyclic components will be unique, up to permutation. But this is something which follows from what we have, coming as a continuation of our study of cyclic groups from chapter 2, and we will leave the details here as an exercise.
\end{proof}

Very nice all this, and you might say, good algebra that we learned here, and time to stop our class, and go get a beer. Well, in answer, while Theorem 4.27 certainly closes the discussion, in what regards the finite abelian groups, when it comes to the arbitrary abelian groups, not necessarily finite, several interesting new phenomena appear.

\bigskip

This being said, with this book being officially dedicated to the finite groups, we will not talk much about this, unless really needed in our study of finite groups. So, in waiting for such a situation, let us just make a brief survey, of what is known about the abelian groups, in general. We first have the following key result, extending Theorem 4.27:

\index{finitely generated group}

\begin{theorem}
The finitely generated abelian groups are the products of type
$$G=\mathbb Z_{N_1}\times\ldots\times\mathbb Z_{N_k}$$
with $N_i\in\mathbb N\cup\{\infty\}$, with the convention $\mathbb Z_\infty=\mathbb Z$.
\end{theorem}

\begin{proof}
As mentioned, we won't talk much about this. The idea is as follows:

\medskip

(1) Given a finite abelian group $G$, we can consider its torsion subgroup $T\subset G$, consisting of elements of finite order, and with $F=G/T$, we then have:
$$G=T\times F$$

(2) Regarding now the $T$ component, this is a finite abelian group, and Theorem 4.27 applies to it. We obtain a decomposition as follows, with $N_i<\infty$:
$$T=\mathbb Z_{N_1}\times\ldots\times\mathbb Z_{N_r}$$

(3) As for the $F$ component, this is free in some technical sense, a bit like vector spaces are, and by using standard vector space techniques we reach to the conclusion that $F$ appears as follows, a bit like the vector spaces must be powers of the ground field:
$$F=\mathbb Z^s$$

(4) Now by putting everything together, we are led to the conclusion in the statement. So, this was for the story, and for more, have a look at Lang \cite{la1}.
\end{proof}

Which is again very nice, and the temptation is high to stop here, and go get that above-mentioned beer. However, not so quick, because we still have some interesting abelian groups, such as our beloved unit circle $\mathbb T$, not covered by the above.

\bigskip

So, what to do with these? In answer, we have the following result, also extending Theorem 4.27, which looks a bit like Theorem 4.28, but is something more advanced:

\index{compact abelian group}

\begin{theorem}
The compact smooth abelian groups are the products of type
$$G=\mathbb Z_{N_1}\times\ldots\times\mathbb Z_{N_k}$$
with $N_i\in\mathbb N\cup\{\infty\}$, but this time with the convention $\mathbb Z_\infty=\mathbb T$.
\end{theorem}

\begin{proof}
To start with, the statement makes sense indeed, with ``smooth'' meaning here having no singularities, that is, covering things like points too. As for the proof, there is certainly some magic trickery relating this to Theorem 4.28. More later.
\end{proof}

And with this, time to end this section. We have learned many interesting things, with those regarding the finite groups coming with full proofs, and with those regarding more general groups being still in need of some explanations. Among others, we would like to understand the mysterious duality $\mathbb Z\leftrightarrow\mathbb T$, relating Theorems 4.28 and 4.29.

\section*{4d. Duality, Fourier}

Welcome back. We would like to discuss now a certain duality for the finite abelian groups, which is in fact a self-duality, $G\leftrightarrow G$, and with this being something non-trivial, with for instance $\mathbb Z_N\leftrightarrow\mathbb Z_N$ producing, with $N\to\infty$, that mysterious duality $\mathbb Z\leftrightarrow\mathbb T$.

\bigskip

Let us start with the following general fact, regarding the arbitrary finite groups:

\index{representation}
\index{character}
\index{dual group}
\index{group of characters}

\begin{theorem}
The characters of a finite group $G$, which are the morphisms
$$\chi:G\to\mathbb T$$
form a finite abelian group $\widehat{G}$, having the following properties:
\begin{enumerate}
\item $\widehat{G}=\widehat{G}_{ab}$, where $G_{ab}=G/[G,G]$ is the abelianization of $G$.

\item We have a morphism $G\to\widehat{\widehat{G}}$, producing a morphism $G_{ab}\to\widehat{\widehat{G}}$.

\item The dual of a product is the product of duals, $\widehat{G\times H}=\widehat{G}\times\widehat{H}$.
\end{enumerate}
\end{theorem}

\begin{proof}
There are many things to be proved here, the idea being as follows:

\medskip

(1) As a first observation, making the link with some previous considerations, from Theorem 4.11, we can talk about the unitary representations $u:G\to U_N$ of any group $G$, and with the character of such a representation being the following function: 
$$\chi_u:G\to\mathbb C\quad,\quad\chi_u(g)=Tr(u(g))$$

Now observe that at $N=1$ our representation is a group morphism $u:G\to\mathbb T$, whose character equals itself, $\chi_u=u$. Thus, the present theorem is about doing representation theory in $N=1$ dimensions, which is something interesting, and good to know.

\medskip

(2) Getting now to the proof, our first claim is that $\widehat{G}$ is indeed a group, with the group operation being the pointwise multiplication of characters, as follows:
$$(\chi\rho)(g)=\chi(g)\rho(g)$$

Indeed, if $\chi,\rho$ are characters, so is $\chi\rho$, and so the multiplication is well-defined on $\widehat{G}$. Regarding the unit, this is the trivial character, constructed as follows:
$$1:G\to\mathbb T\quad,\quad 
g\to1$$ 

Finally, we have inverses, with the inverse of $\chi:G\to\mathbb T$ being its conjugate:
$$\bar{\chi}:G\to\mathbb T\quad,\quad 
g\to\overline{\chi(g)}$$

(3) Next, the fact that $\widehat{G}$ is finite comes from the following implication, which shows that for any $g\in G$, there are finitely many choices for the value $\chi(g)\in\mathbb T$:
$$g^k=1\implies\chi(g)^k=1$$

Finally, the fact that $\widehat{G}$ is indeed abelian is clear from definitions.

\medskip

(4) Let us prove now the equality $\widehat{G}=\widehat{G}_{ab}$. We recall that given a group $G$, its commutator subgroup $[G,G]\subset G$ is constructed as follows:
$$[G,G]=\left\{ghg^{-1}h^{-1}\Big|g,h\in G\right\}$$

This subgroup is then normal, and we can define the abelianization of $G$ as being $G_{ab}=G/[G,G]$. Now observe that for any character $\chi:G\to\mathbb T$, we have:
$$\chi(ghg^{-1}h^{-1})=1$$

Thus, our character factorizes as follows, into a character of the group $G_{ab}$:
$$\chi:G\to G_{ab}\to\mathbb T$$

Summarizing, we have constructed an identification $\widehat{G}=\widehat{G}_{ab}$, as claimed.

\medskip

(5) Next, we have a group morphism as follows, obtained by evaluation:
$$I:G\to\widehat{\widehat{G}}\quad,\quad g\to[\chi\to\chi(g)]$$

Since $\chi(ghg^{-1}h^{-1})=1$ we have $[G,G]\subset\ker I$, so $I$ factorizes as follows, as claimed:
$$I:G\to G_{ab}\to \widehat{\widehat{G}}$$

(6) Finally, in what regards the last assertion, namely $\widehat{G\times H}=\widehat{G}\times\widehat{H}$, observe that a character of a product of groups $\chi:G\times H\to\mathbb T$ must satisfy:
$$\chi(g,h)
=\chi(g,1)\chi(1,h)
=\chi_{|G}(g)\chi_{|H}(h)$$

Thus $\chi$ must appear as the product of its restrictions $\chi_{|G},\chi_{|H}$, which must be both characters, and this gives the identification in the statement.
\end{proof}

When the starting finite group $G$ is abelian itself, which is more or less the case of interest, due to the formula $\widehat{G}=\widehat{G}_{ab}$ established above, we have:

\index{cyclic group}
\index{product of cyclic groups}
\index{self-dual group}

\begin{theorem}
The character group operation $G\to\widehat{G}$ for the finite abelian groups, called Pontrjagin duality, has the following properties:
\begin{enumerate}
\item The dual of a cyclic group is the group itself, $\widehat{\mathbb Z}_N=\mathbb Z_N$.

\item The dual of a product is the product of duals, $\widehat{G\times H}=\widehat{G}\times\widehat{H}$.

\item Any product of cyclic groups $G=\mathbb Z_{N_1}\times\ldots\times\mathbb Z_{N_k}$ is self-dual, $G=\widehat{G}$.
\end{enumerate}
\end{theorem}

\begin{proof}
We have several things to be proved, the idea being as follows:

\medskip

(1) By multiplicativity, a character $\chi:\mathbb Z_N\to\mathbb T$ is uniquely determined by its value $z=\chi(g)$ on the standard generator $g\in\mathbb Z_N$. But this value must satisfy:
$$z^N=1$$

Conversely, any $z\in\mathbb Z_N$ defines a character $\chi:\mathbb Z_N\to\mathbb T$, by setting, for any $r\in\mathbb N$:
$$\chi(g^r)=z^r$$

Thus we have an identification of abelian groups $\widehat{\mathbb Z}_N=\mathbb Z_N$, as claimed.

\medskip

(2) This is something clear, that we already know, from Theorem 4.30.

\medskip

(3) This follows from (1) and (2). Alternatively, any character $\chi:G\to\mathbb T$ is uniquely determined by its values $\chi(g_1),\ldots,\chi(g_k)$ on the standard generators of $\mathbb Z_{N_1},\ldots,\mathbb Z_{N_k}$, which must belong to $\mathbb Z_{N_1},\ldots,\mathbb Z_{N_k}\subset\mathbb T$, and this gives $\widehat{G}=G$, as claimed.
\end{proof}

Observe that Theorem 4.31 and its proof use nothing, or almost. At a more advanced level now, making the link with Theorem 4.27, and with analysis too, we have:

\index{finite abelian group}
\index{Fourier transform}
\index{dual group}

\begin{theorem}
Given a finite abelian group $G$, the following happen:
\begin{enumerate}

\item The Pontrjagin duality for it is a self-duality, $G\simeq \widehat{G}$.

\item In particular, the Pontrjagin duality is indeed a duality, $G=\widehat{\widehat{G}}$.

\item We have an isomorphism of commutative algebras $\mathbb C[G]\simeq C(\widehat{G})$.
\end{enumerate}
\end{theorem}

\begin{proof}
These assertions are more advanced, the idea being as follows:

\medskip

(1) This comes from Theorem 4.27, and Theorem 4.31 (3).

\medskip

(2) This assertion, further building on Theorem 4.30 (2), is clear from (1).

\medskip

(3) Consider indeed the group algebra $\mathbb C[G]$, that we already met in Theorem 4.11. We have then a morphism of commutative algebras, as follows:
$$\varphi:\mathbb C[G]\to C(\widehat{G})\quad,\quad \sum_h\lambda_gg\to\left[\chi\to\sum_g\lambda_g\chi(g)\right]$$

But for $G=\mathbb Z_N$ this morphism is easily seen to be an isomorphism, and then by taking products and using Theorem 4.27, we obtain the result in general. We will actually come back to this in a moment, with some further details, and explicit formulae as well.
\end{proof}

We can feel that all this is related to Fourier analysis, and we have:

\begin{fact}
The following happen, regarding the locally compact abelian groups:
\begin{enumerate}
\item What we did in the finite case, namely group characters, and construction and basic properties of the dual, can be extended to them.

\item As basic examples of this, besides what we have in the finite case, and notably $\widehat{\mathbb Z}_N=\mathbb Z_N$, we have $\widehat{\mathbb Z}=\mathbb T$, $\widehat{\mathbb T}=\mathbb Z$, and also $\widehat{\mathbb R}=\mathbb R$.

\item With some care for analytic aspects, $\mathbb C[G]\simeq C(\widehat{G})$ remains true in this setting, and in the case $G=\mathbb R$, this isomorphism is the Fourier transform.
\end{enumerate}
\end{fact}

Obviously, all this is a bit heavy, and rather going into analysis, which is not our main topic for this book, but you get the point, we have 3 types of Fourier analysis in life:

\medskip

--  We first have the ``standard'' Fourier analysis, the one for real functions, that you surely know from advanced calculus, corresponding to the group $G=\mathbb R$.

\medskip

-- Then we have another type of Fourier analysis, that you surely met in advanced calculus too, called the ``Fourier series'' one, corresponding to $G=\mathbb Z,\mathbb T$.

\medskip

-- And finally we have the ``discrete'' version of Fourier analysis that we started to learn here, over the cyclic group $G=\mathbb Z_N$, and other finite abelian groups.

\medskip

In practice, all this is a bit complicated, and back now to the finite abelian groups, let us work out a softer version of all the above, which is what is really needed, in practice, when doing discrete Fourier analysis. The result here is as follows:

\index{Fourier matrix}
\index{Gauss sum}
\index{all-one matrix}

\begin{theorem}
Given a finite abelian group $G$, with dual group $\widehat{G}=\{\chi:G\to\mathbb T\}$, consider the corresponding Fourier coupling, namely:
$$\mathcal F_G:G\times\widehat{G}\to\mathbb T\quad,\quad 
(i,\chi)\to\chi(i)$$
\begin{enumerate}
\item Via the standard isomorphism $G\simeq\widehat{G}$, this Fourier coupling can be regarded as a square matrix, $F_G\in M_G(\mathbb T)$, which is a complex Hadamard matrix.

\item For $G=\mathbb Z_N$ we obtain in this way, via the identification $\mathbb Z_N=\{1,\ldots,N\}$, the standard Fourier matrix, $F_N=(w^{ij})_{ij}$ with $w=e^{2\pi i/N}$.

\item In general, when using a decomposition $G=\mathbb Z_{N_1}\times\ldots\times\mathbb Z_{N_k}$, the corresponding Fourier matrix is given by $F_G=F_{N_1}\otimes\ldots\otimes F_{N_k}$.
\end{enumerate}
Moreover, the matrix of the isomorphism $\mathbb C[G]\to C(\widehat{G})$ is this Fourier matrix $F_G$.
\end{theorem}

\begin{proof}
This is something very standard, the idea being as follows:

\medskip

(1) With the identification $G\simeq\widehat{G}$ made, our matrix is given by $(F_G)_{i\chi}=\chi(i)$, and the Hadamard condition, stating that the rows are pairwise orthogonal, comes from:
$$<R_i,R_j>
=\sum_\chi\chi(i)\overline{\chi(j)}
=\sum_\chi\chi(i-j)
=|G|\cdot\delta_{ij}$$

(2) This follows indeed from the elementary fact that, via the standard identifications $\mathbb Z_N=\widehat{\mathbb Z}_N=\{1,\ldots,N\}$, the Fourier coupling here is as follows, with $w=e^{2\pi i/N}$:
$$(i,j)\to w^{ij}$$

Thus, we have $F_N=(w^{ij})_{ij}$ with $w=e^{2\pi i/N}$, as stated. Let us record as well the matrix form of $F_N$. With indices $i,j\in\{0,1,\ldots,N-1\}$, this is as follows:
$$F_N=
\begin{pmatrix}
1&1&1&\ldots&1\\
1&w&w^2&\ldots&w^{N-1}\\
1&w^2&w^4&\ldots&w^{2(N-1)}\\
\vdots&\vdots&\vdots&&\vdots\\
1&w^{N-1}&w^{2(N-1)}&\ldots&w^{(N-1)^2}
\end{pmatrix}$$

And interesting matrix this is. For instance any circulant matrix is diagonalized by it, with a basic example here being the all-one matrix $\mathbb I_N$, which diagonalizes as:
$$\mathbb I_N=\frac{1}{N}\,F_N\,diag(N,0,\ldots,0)\,F_N^*$$

As another interesting feature of $F_N$, making the link with arithmetic, and many things to be learned here too, its trace is the standard quadratic Gauss sum:
$$Tr(F_N)=\sum_kw^{k^2}$$

(3) This comes indeed from the following formula, which is clear from definitions:
$$F_{H\times K}=F_H\otimes F_K$$

(4) As for the last assertion, making the link with Theorem 4.32 and Fact 4.33, and with Fourier analysis in general, this comes again from definitions.
\end{proof}

\section*{4e. Exercises}

This was a quite subtle algebra and analysis chapter, and as exercises, we have:

\begin{exercise}
Learn more about conjugation, centralizers, and central functions.
\end{exercise}

\begin{exercise}
Learn more about orbits, orbitals, and higher orbitals.
\end{exercise}

\begin{exercise}
Prove the Cauchy theorem for the abelian groups.
\end{exercise}

\begin{exercise}
Prove that a group with $|G|=p^2$ must be abelian.
\end{exercise}

\begin{exercise}
Prove that $|G|=2p$ implies $G=\mathbb Z_{2p}$ or $G=D_{2p}$.
\end{exercise}

\begin{exercise}
Further explore the groups with $|G|=pq$, with $p>q>2$.
\end{exercise}

\begin{exercise}
Learn, as much as you can, about abelian groups and duality.
\end{exercise}

\begin{exercise}
Learn also more about the complex Hadamard matrices.
\end{exercise}

As bonus exercise, read some functional analysis, which is related to all this.

\part{Rotation groups}

\ \vskip50mm

\begin{center}
{\em Red, red wine

Goes to my head

Makes me forget that I

Still need her so}
\end{center}

\chapter{Reflection groups}

\section*{5a. Cube symmetries}

We have seen a lot of interesting theory in Part I, but at the level of concrete examples of finite groups, I mean those that we can really understand, geometrically, we are still left with three of them, namely the cyclic group $\mathbb Z_N$, the dihedral group $D_N$, and the symmetric group $S_N$. These groups are all understood in terms of their action on $\{1,\ldots,N\}$, which becomes increasingly finer, according to the following embeddings:
$$\mathbb Z_N\subset D_N\subset S_N$$

Equivalently, and getting a bit more geometric, we can equally view these groups as permuting the standard basis $\{e_1,\ldots,e_N\}$ of $\mathbb R^N$, and with this being something useful, both in order to understand the geometry of $\mathbb R^N$, and of the groups themselves. Thus, the full geometric picture for our main finite groups is in fact as follows:
$$\mathbb Z_N\subset D_N\subset S_N\subset O_N$$

What is next? Many things, and as an obvious objective, we will definitely have to extend this list of examples. We will do this in the present Part II, with a look at all sorts of variations of $\mathbb Z_N,D_N,S_N$ that can be constructed, with geometric ideas in mind.

\bigskip

But probably too much talking. Let us start with something quite exciting, sort of a candidate for ``the most beautiful finite group ever'', namely:

\index{hyperoctahedral group}
\index{hypercube}

\begin{definition}
The hyperoctahedral group $H_N\subset O_N$ is the group formed by the symmetries of the unit cube in $\mathbb R^N$,
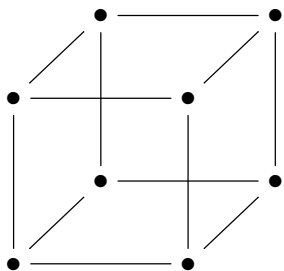
$$\xymatrix@R=19pt@C=20pt{
&\bullet\ar@{-}[rr]&&\bullet\\
\bullet\ar@{-}[rr]\ar@{-}[ur]&&\bullet\ar@{-}[ur]\\
&\bullet\ar@{-}[rr]\ar@{-}[uu]&&\bullet\ar@{-}[uu]\\
\bullet\ar@{-}[uu]\ar@{-}[ur]\ar@{-}[rr]&&\bullet\ar@{-}[uu]\ar@{-}[ur]
}$$
viewed as a graph, or equivalently, as a metric space.
\end{definition}

Here the equivalence at the end is clear from definitions, because any symmetry of the cube graph must preserve the lengths of the edges, and so we have:
$$G(\square_{graph})=G(\square_{metric})$$

The hyperoctahedral group is a quite interesting group, whose definition, as a symmetry group, reminds that of the dihedral group $D_N$. So, let us start our study in the same way as we did for $D_N$, with a discussion at small values of $N\in\mathbb N$:

\medskip

\underline{$N=1$}. Here the 1-cube is the segment, whose symmetries are the identity $id$, plus the symmetry $\tau$ with respect to the middle of the segment:
$$
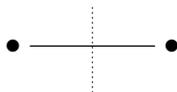
\xymatrix@R=10pt@C=20pt{
&\ar@{.}[dd]\\
\bullet\ar@{-}[rr]&&\bullet\\
&}$$

Thus, we obtain the group with 2 elements, which is a very familiar object: 
$$H_1=D_2=S_2=\mathbb Z_2$$

\underline{$N=2$}. Here the 2-cube is the square, whose symmetries are the 4 rotations, of angles $0^\circ,90^\circ,180^\circ,270^\circ$, and the 4 symmetries with respect to the 4 symmetry axes, which are the 2 diagonals, and the 2 segments joining the midpoints of opposite sides:
$$
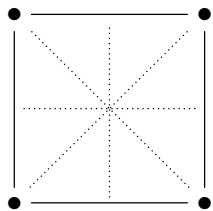
\xymatrix@R=26pt@C=26pt{
\bullet\ar@{-}[dd]\ar@{.}[ddrr]\ar@{-}[rr]&\ar@{.}[dd]&\bullet\ar@{-}[dd]\ar@{.}[ddll]\\
\ar@{.}[rr]&&\\
\bullet\ar@{-}[rr]&&\bullet
}$$

Thus, we obtain a group with 8 elements, which again is a very familiar object:
$$H_2=D_4=\mathbb Z_4\rtimes\mathbb Z_2$$

\underline{$N=3$}. Here the 3-cube is the usual cube in $\mathbb R^3$, pictured as follows:
$$
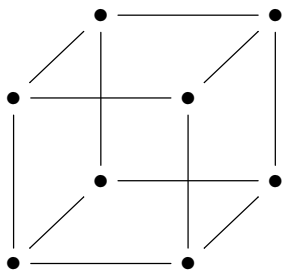
\xymatrix@R=19pt@C=20pt{
&\bullet\ar@{-}[rr]&&\bullet\\
\bullet\ar@{-}[rr]\ar@{-}[ur]&&\bullet\ar@{-}[ur]\\
&\bullet\ar@{-}[rr]\ar@{-}[uu]&&\bullet\ar@{-}[uu]\\
\bullet\ar@{-}[uu]\ar@{-}[ur]\ar@{-}[rr]&&\bullet\ar@{-}[uu]\ar@{-}[ur]
}$$

However, in relation with the symmetries, the situation now is considerably more complicated, because, thinking well, this cube has no less than 48 symmetries.

\medskip

All this looks quite complicated, but fortunately we can count $H_N$, at $N=3$, and at higher $N$ as well, by using a geometric trick, the result being as follows:

\begin{theorem}
We have the cardinality formula
$$|H_N|=2^NN!$$
coming from the fact that $H_N$ is the symmetry group of the coordinate axes of $\mathbb R^N$.
\end{theorem}

\begin{proof}
This follows from some geometric thinking, as follows:

\medskip

(1) Consider the standard cube in $\mathbb R^N$, centered at 0, and having as vertices the points having coordinates $\pm1$. With this picture in hand, it is clear that the symmetries of the cube coincide with the symmetries of the $N$ coordinate axes of $\mathbb R^N$.

\medskip

(2) In order to count now these latter symmetries, a bit as we did for the dihedral group, observe first that we have $N!$ permutations of these $N$ coordinate axes. 

\medskip

(3) But each of these permutations of the coordinate axes $\sigma\in S_N$ can be further ``decorated'' by a sign vector $e\in\{\pm1\}^N=\mathbb Z_2^N$, consisting of the possible $\pm1$ flips which can be applied to each coordinate axis, at the arrival. 

\medskip

(4) And the point is that, obviously, we obtain in this way all the elements of $H_N$. Thus, we have the following formula, for the cardinality of $H_N$:
$$|H_N|
=|S_N|\cdot|\mathbb Z_2^N|
=N!\cdot2^N$$

Thus, we are led to the conclusion in the statement. 
\end{proof}

The above result is quite interesting, in connection with the considerations from the beginning of this chapter, regarding our scheme $\mathbb Z_N\subset D_N\subset S_N\subset O_N$. Indeed, the hyperoctahedral group falls into this scheme, somehow as a new king, as follows:

\begin{comment}
We have embeddings of groups as follows,
$$\mathbb Z_N\subset D_N\subset S_N\subset H_N\subset O_N$$
with $H_N$ apparently replacing $S_N$, as the most important finite group.
\end{comment}

Of course, this is something a bit conjectural, and we will certainly have to work on this, in order to decide whether $H_N$ is indeed the new king. But this will be the case, trust me here, as the present chapter and Part II develops. Actually, the main opposition to the rule of $H_N$ will not come from below, from $S_N$, but rather from the above, from a certain complex version of it $K_N\subset U_N$ that can be constructed. More on this later.

\bigskip

Getting back now to what we have about $H_N$, as concrete results, namely Theorem 5.2, as in the dihedral group case, it is possible to go beyond this, with a crossed product decomposition, of quite special type, called wreath product decomposition:

\index{crossed product}
\index{wreath product}
\index{hyperoctahedral group}

\begin{theorem}
We have a wreath product decomposition $H_N=\mathbb Z_2\wr S_N$, which means by definition that we have a crossed product decomposition
$$H_N=\mathbb Z_2^N\rtimes S_N$$
with the permutations $\sigma\in S_N$ acting on the elements $e\in\mathbb Z_2^N$ as follows:
$$\sigma(e_1,\ldots,e_N)=(e_{\sigma(1)},\ldots,e_{\sigma(N)})$$
In particular we have, as found before, the cardinality formula $|H_N|=2^NN!$.
\end{theorem}

\begin{proof}
As explained in the proof of Theorem 5.2, the elements of $H_N$ can be identified with the pairs $g=(e,\sigma)$ consisting of a permutation $\sigma\in S_N$, and a sign vector $e\in\mathbb Z_2^N$, so that at the level of the cardinalities, we have the following formula:
$$|H_N|=|\mathbb Z_2^N\times S_N|$$

To be more precise, given an element $g\in H_N$, the element $\sigma\in S_N$ is the corresponding permutation of the $N$ coordinate axes, regarded as unoriented lines in $\mathbb R^N$, and $e\in\mathbb Z_2^N$ is the vector collecting the possible flips of these coordinate axes, at the arrival. Now observe that the product formula for two such pairs $g=(e,\sigma)$ is as follows, with the permutations $\sigma\in S_N$ acting on the elements $f\in\mathbb Z_2^N$ as in the statement:
$$(e,\sigma)(f,\tau)=(ef^\sigma,\sigma\tau)$$

Thus, we are precisely in the framework of the crossed products, as constructed in chapter 2, and we conclude that we have a crossed product decomposition, as follows:
$$H_N=\mathbb Z_2^N\rtimes S_N$$
 
Thus, we are led to the conclusion in the statement, with the formula $H_N=\mathbb Z_2\wr S_N$ being just a shorthand for the decomposition $H_N=\mathbb Z_2^N\rtimes S_N$ that we found.
\end{proof}

In view to further understand $H_N\subset O_N$, we can use the permutation matrices. We are led in this way to the following result, further extending what we know:

\begin{theorem}
The hyperoctahedral group $H_N$ appears also as
$$H_N=M_N(-1,0,1)\cap O_N$$
that is, as the group of permutation matrices with signed entries. We have
$$H_N=\mathbb Z_2\wr S_N\quad,\quad |H_N|=2^NN!$$
trivially, coming from this latter interpretation.
\end{theorem}

\begin{proof}
We have several assertions here, the idea being as follows:

\medskip

(1) To start with, according to Theorem 5.2, the hyperoctahedral group $H_N\subset O_N$ appears indeed as the group of permutation matrices with signed entries.

\medskip

(2) Next, this gives an embedding as follows, that we claim to be an equality:
$$H_N\subset M_N(-1,0,1)\cap O_N$$

Indeed, let us pick a matrix $U\in M_N(-1,0,1)\cap O_N$. Since by $U\in O_N$ the rows must have norm 1, there is exactly one $\pm1$ on each row, and the same goes for columms. Thus we must have $U\in H_N$, and our embedding above follows to be an equality.

\medskip

(3) Finally, the last assertions are imported from Theorem 5.4, with these being quite clear too in our new context, namely $H_N=M_N(-1,0,1)\cap O_N$.
\end{proof}

Moving forward now, and with Comment 5.3 in mind, let us formulate:

\begin{question}
What are the analogues of the signature map 
$$\varepsilon:S_N\to\mathbb Z_2$$
and of its kernel $\ker\varepsilon=A_N$, in the context of $H_N$?
\end{question}

And good question this is, hope you agree with me, because in order to truly dethrone $S_N$, the hyperoctahedral group $H_N$ would normally need such features, too.

\bigskip

In answer now, as we know well since chapter 3, the signature of the permutations $\sigma\in S_N$ is the determinant of the corresponding permutation matrices $\sigma\in S_N\subset O_N$, so the answer to the signature map question above is straightforward, as follows:

\begin{proposition}
We have a signature function on the hyperoctahedral group,
$$\varepsilon:H_N\to\mathbb Z_2$$
obtained via $H_N\subset O_N$, by taking the determinant.
\end{proposition}

\begin{proof}
This is indeed something self-explanatory, and with the value of the signature map $\varepsilon$ on a given group element $U\in M_N(-1,0,1)\cap O_N$ being the signature of the underlying permutation $\sigma\in S_N$, multiplied by the number of $-1$ entries of $U$.
\end{proof}

Getting now to our second question above, we know that the alternating group $A_N$ appears according to the formula $A_N=S_N\cap SO_N$. In view of this, let us formulate:

\begin{definition}
The special hyperoctahedral group $SH_N$ is:
$$SH_N=H_N\cap SO_N$$
That is, $SH_N$ is the group of orientation-preserving symmetries of the cube $\square_N$.
\end{definition}

As a first observation, $SH_N$ is indeed the kernel of the signature map $\varepsilon:H\to\mathbb Z_2$ from Proposition 5.7. Also, with respect to the various groups that we already know, namely $O_N$ and its subgroups $A_N,S_N,H_N$ and $SO_N$, we have a diagram as follows:
$$\xymatrix@R=56pt@C=50pt{
S_N\ar[r]&H_N\ar[r]&O_N\\
A_N\ar[u]\ar[r]&SH_N\ar[u]\ar[r]&SO_N\ar[u]}$$

To be more precise, here all arrows are inclusions, and this is an intersection diagram, in the sense that any square subdiagram $A\subset B,C\subset D$ satisfies $A=B\cap C$.

\bigskip

More concretely now, let us start now our study of $SH_N$ in the same way as we did for $H_N$, with a discussion at small values of $N\in\mathbb N$. The situation here is as follows:

\medskip

\underline{$N=1$}. Here the 1-cube is the segment, whose symmetries are the identity $id$, plus the symmetry $\tau$ with respect to the middle of the segment:
$$
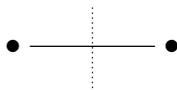
\xymatrix@R=10pt@C=20pt{
&\ar@{.}[dd]\\
\bullet\ar@{-}[rr]&&\bullet\\
&}$$

Thus, we obtain the group with 1 element, which is a very familiar object: 
$$SH_1=\{1\}$$

\underline{$N=2$}. Here the 2-cube is the square, whose symmetries are the 4 rotations, of angles $0^\circ,90^\circ,180^\circ,270^\circ$, and the 4 symmetries with respect to the 4 symmetry axes, which are the 2 diagonals, and the 2 segments joining the midpoints of opposite sides:
$$
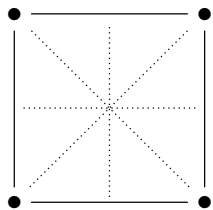
\xymatrix@R=26pt@C=26pt{
\bullet\ar@{-}[dd]\ar@{.}[ddrr]\ar@{-}[rr]&\ar@{.}[dd]&\bullet\ar@{-}[dd]\ar@{.}[ddll]\\
\ar@{.}[rr]&&\\
\bullet\ar@{-}[rr]&&\bullet
}$$

Thus, we obtain a group with 4 elements, which again is a very familiar object:
$$SH_2=\mathbb Z_4$$

\underline{$N=3$}. Here the 3-cube is the usual cube in $\mathbb R^3$, pictured as follows:
$$
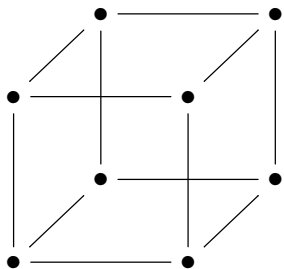
\xymatrix@R=19pt@C=20pt{
&\bullet\ar@{-}[rr]&&\bullet\\
\bullet\ar@{-}[rr]\ar@{-}[ur]&&\bullet\ar@{-}[ur]\\
&\bullet\ar@{-}[rr]\ar@{-}[uu]&&\bullet\ar@{-}[uu]\\
\bullet\ar@{-}[uu]\ar@{-}[ur]\ar@{-}[rr]&&\bullet\ar@{-}[uu]\ar@{-}[ur]
}$$

In relation with symmetries, we know that this cube has a total of $|H_3|=48$ symmetries, and it is easy to see that half of these are orientation-preserving, $|SH_3|=24$. In fact, there is one such orientation-preserving symmetry for each permutation of the 4 diagonals of the cube, so we reach once again to a very familiar object, namely:
$$SH_3=S_4$$

\underline{$N=4$}. Here things are considerably more complicated, because we have a total of $|H_4|=384$ symmetries, half of which are orientation-preserving, $|SH_4|=192$. Thus, we will stop our study here, without attempting to understand the structure of $SH_4$.

\bigskip

As a conclusion, we can see that in each case half of the elements of $H_N$ belong to $SH_N$, and the other half does not belong to $SH_N$. We have in fact the following result about this, along with a bit more, inspired by our previous study of $H_N$:

\begin{theorem}
The special hyperoctahedral group $SH_N$ appears also as
$$SH_N=M_N(-1,0,1)\cap SO_N$$
that is, as the group of signed permutation matrices of determinant $1$. We have
$$|SH_N|=2^{N-1}N!$$
showing that $SH_N\subset H_N$ has indeed index $2$, coming from this interpretation.
\end{theorem}

\begin{proof}
This is indeed something self-explanatory, with the first assertion coming from the following computation, based on Definition 5.8 and Theorem 5.5:
\begin{eqnarray*}
SH_N
&=&H_N\cap SO_N\\
&=&M_N(-1,0,1)\cap O_N\cap SO_N\\
&=&M_N(-1,0,1)\cap SO_N
\end{eqnarray*}

As for the element count, this is clear from this interpretation. Indeed, we know that $A_N\subset S_N$ accounts for half of the permutations. But when decorating with $\pm$ signs the proportions will not change, so $SH_N\subset H_N$ has indeed index 2, as claimed.
\end{proof}

Many other things can be said about $H_N$, $SH_N$, and we will be back to this on several occasions, directly in a more general setting, that of the complex reflection groups.

\section*{5b. Complex reflections} 

As a continuation of our main story in this chapter, namely battle of kings between $S_N$ and $H_N$, we can have a good look at this algebrically, according to:
$$S_N=M_N(0,1)\cap O_N\quad,\quad H_N=M_N(-1,0,1)\cap O_N$$

Indeed, we can see from this interpretation that $H_N$ is smarter than $S_N$, due to its use of $\mathbb Z_2=\{\pm1\}$ as nonzero scalars, instead of $\mathbb Z_1=\{1\}$. But with this idea in mind, we are naturally lead into using an arbitrary $\mathbb Z_s$ as nonzero scalars, as follows: 

\index{complex reflection group}

\begin{definition}
The complex reflection group $H_N^s\subset U_N$, depending on parameters
$$N\in\mathbb N\quad,\quad s\in\mathbb N\cup\{\infty\}$$
is the group of permutation-type matrices with $s$-th roots of unity as entries,
$$H_N^s=M_N(\mathbb Z_s\cup\{0\})\cap U_N$$
with the convention $\mathbb Z_\infty=\mathbb T$, at $s=\infty$.
\end{definition}

This construction is something quite tricky, that will keep us busy, for the remainder of this section. As a first observation, at $s=1,2$ we obtain the following groups:
$$H_N^1=S_N\quad,\quad 
H_N^2=H_N$$

Another important particular case of the above construction is $s=\infty$, where we obtain a group which is actually not finite, but is still compact, denoted as follows:
$$K_N\subset U_N$$

This latter group $K_N$ is called full complex reflection group, and will appear many times, in what follows. In view of this, let us highlight its definition, as follows:

\index{full reflection group}

\begin{definition}
The full complex reflection group is given by:
$$K_N=M_N(\mathbb T\cup\{0\})\cap U_N$$
That is, $K_N$ is the group of permutation-type matrices with entries from $\mathbb T$.
\end{definition}

In practice now, we seem to have a new king here, for the finite group theory, and with the aristocratic nature of $K_N$ being quite obvious, due to its non-finiteness.

\bigskip

Well, jokes left aside, $K_N$ seems to be a quite interesting object, with its potential remaining to be determined. So, let us first have a look at it at small values of $N$:

\medskip

\underline{$N=1$}.  What we have is the unit circle, $K_1=\mathbb T$.

\medskip 

\underline{$N=2$}. Here $K_2$ consists of the matrices as follows, with nonzero entries in $\mathbb T$:
$$\begin{pmatrix}x&0\\0&y\end{pmatrix}\qquad,\qquad
\begin{pmatrix}0&x\\y&0\end{pmatrix}$$

\underline{$N=3$}. Here $K_3$ consists of the matrices as follows, with nonzero entries in $\mathbb T$:
$$\begin{pmatrix}x&0&0\\0&y&0\\0&0&z\end{pmatrix}\quad,\quad
\begin{pmatrix}0&x&0\\y&0&0\\0&0&z\end{pmatrix}\quad,\quad
\begin{pmatrix}x&0&0\\0&0&y\\0&z&0\end{pmatrix}$$
$$\begin{pmatrix}0&0&x\\0&y&0\\z&0&0\end{pmatrix}\quad,\quad
\begin{pmatrix}0&0&x\\y&0&0\\0&z&0\end{pmatrix}\quad,\quad
\begin{pmatrix}0&x&0\\0&0&y\\z&0&0\end{pmatrix}$$

\underline{$N\geq4$}. And so on, you get the point, what we have is a bit like before for $H_N$, permutation matrices, but this time decorated by numbers in $\mathbb T$.

\medskip

The problem is now, is all this, passage to the complex numbers, really worth it? And in answer, certainly yes in relation with the basic groups that we know, as shown by:

\begin{fact}
The group $K_N$ is smarter than $S_N,H_N$, due to the following reasons:
\begin{enumerate}
\item We have $\mathbb Z_N\subset K_1$, which is better than $\mathbb Z_N\subset S_N\subset H_N$.

\item We have $D_N\subset K_2$, which better than $D_N\subset S_N\subset H_N$.
\end{enumerate}
\end{fact}

To be more precise here, in what regards (1), we certainly have $\mathbb Z_N\subset\mathbb T=K_1$, and this is something quite magic, making all of the sudden $\mathbb Z_N$ a 1-dimensional object, instead of the $N$-dimensional object that it used to be, according to $\mathbb Z_N\subset S_N\subset H_N$.

\bigskip

As for (2), which is even more impressive, this comes from the following result:

\begin{theorem}
We have an embedding $D_N\subset K_2$, coming as follows,
$$D_N=\left\{\begin{pmatrix}x&0\\0&y\end{pmatrix}\ ,\ 
\begin{pmatrix}0&x\\y&0\end{pmatrix}\Big|x=y^{-1}\in\mathbb Z_N\right\}\subset K_2$$
obtained by augmenting the standard copy $\mathbb Z_N\subset K_2$ with a twisted copy of it.
\end{theorem}

\begin{proof}
The matrices patterned as in the statement form indeed a group, and when adding the extra condition $xy=1$, this remains a group. In order now to establish the isomorphism with $D_N$, let us label our group elements as follows, with $xy=zt=1$:
$$R_x=\begin{pmatrix}x&0\\0&y\end{pmatrix}\quad,\quad 
S_z=\begin{pmatrix}0&z\\t&0\end{pmatrix}$$

We have then the following computations, for the products of these elements:
$$R_xR_z=\begin{pmatrix}x&0\\0&y\end{pmatrix}
\begin{pmatrix}z&0\\0&t\end{pmatrix}
=\begin{pmatrix}xz&0\\0&yt\end{pmatrix}
=R_{xz}$$
$$R_xS_z=\begin{pmatrix}x&0\\0&y\end{pmatrix}
\begin{pmatrix}0&z\\t&0\end{pmatrix}
=\begin{pmatrix}0&xz\\yt&0\end{pmatrix}
=S_{xz}$$
$$S_xR_z=
\begin{pmatrix}0&x\\y&0\end{pmatrix}
\begin{pmatrix}z&0\\0&t\end{pmatrix}\
=\begin{pmatrix}0&xt\\yz&0\end{pmatrix}
=S_{xz^{-1}}$$
$$S_xS_z=\begin{pmatrix}0&x\\y&0\end{pmatrix}
\begin{pmatrix}0&z\\t&0\end{pmatrix}
=\begin{pmatrix}xt&0\\0&yz\end{pmatrix}
=R_{xz^{-1}}$$

But, we recognize here the table of multiplication of $D_N$, as desired.
\end{proof}

Summarizing, Fact 5.12 understood, and with this being something remarkable. In addition to this, we have plenty of other good reasons for passing from the real to the complex numbers, at this point of this book, due to the following facts:

\bigskip

-- In chapter 6 we intend to talk about braid groups and knot invariants, which are best understood in terms of complex algebras. Then in chapter 7 we will talk about rotation groups in small dimensions, whose study requires the quotient map $SU_2\to SO_2$. And then in chapter 8 we will talk about Galois groups, needing $\mathbb C$ as ground field too.

\bigskip

-- As for Part III and Part IV, these will deal with representation theory and its applications, and with this being very related to linear algebra, which notoriously best lives over $\mathbb C$, due to the fact that most real matrices, as for instance the rotations of $\mathbb R^2$, have complex eigenvalues and eigenvectors, our study here will be over $\mathbb C$ too.

\bigskip

Summarizing, good time to pass to the complex numbers, and the complex reflection groups $H_N^s\subset U_N$ from Definition 5.10, with special attention to the group $H_N^\infty=K_N$ from Definition 5.11, which contains them all, will be our new objects of interest.

\bigskip

Let us start our study by summarizing some basic observations, as follows:

\begin{proposition}
The complex reflection groups $H_N^s\subset U_N$ are as follows:
\begin{enumerate}
\item At $s=1$ we have $H_N^1=S_N$, having cardinality $|S_N|=N!$.

\item At $s=2$ we have $H_N^2=H_N$, having cardinality $|H_N|=2^NN!$.

\item At $s=\infty$ we have $H_N^\infty=K_N$, having cardinality $|K_N|=\infty$.
\end{enumerate}
\end{proposition}

\begin{proof}
This is clear indeed from the discussion made after Definition 5.10, and with the cardinality results at $s=1$ and $s=2$ being something that we know well.
\end{proof}

Let us record as well the following result, which is something elementary too:

\begin{proposition}
We have inclusions as follows, for any $r,s\in\mathbb N\cup\{\infty\}$:
$$r|s\implies H_r\subset H_s$$
In particular, we have inclusions $S_N\subset H_N^s\subset K_N$, for any $s\in\mathbb N\cup\{\infty\}$.
\end{proposition}

\begin{proof}
With the cyclic group $\mathbb Z_s$ being viewed as group of the $s$-th roots of unity in the complex plane, as in Definition 5.10, we have inclusions as follows:
$$r|s\implies \mathbb Z_r\subset \mathbb Z_s$$

Thus, with the group $H_N^s$ constructed as in Definition 5.10, for $r|s$ we have:
\begin{eqnarray*}
H_N^r
&=&M_N(\mathbb Z_r\cup\{0\})\cap U_N\\
&\subset&M_N(\mathbb Z_s\cup\{0\})\cap U_N\\
&=&H_N^s
\end{eqnarray*} 

Finally, the last assertion is clear, and comes also from this, via $1|s|\infty$, for any $s$.
\end{proof}

Coming next, in analogy with what we know about $S_N,H_N$, we first have:

\begin{proposition}
The number of elements of $H_N^s$ with $s\in\mathbb N$ is:
$$|H_N^s|=s^NN!$$
At $s=\infty$, the group $K_N=H_N^\infty$ that we obtain is infinite.
\end{proposition}

\begin{proof}
This is indeed clear from our definition of $H_N^s$, as a matrix group, because there are $N!$ choices for a permutation-type matrix, and then $s^N$ choices for the corresponding $s$-roots of unity, which must decorate the $N$ nonzero entries.
\end{proof}

Once again in analogy with what we know at $s=1,2$, we have as well:

\index{wreath product}

\begin{theorem}
We have a wreath product decomposition
$$H_N^s=\mathbb Z_s^N\rtimes S_N=\mathbb Z_s\wr S_N$$
with the permutations $\sigma\in S_N$ acting on the elements $e\in\mathbb Z_s^N$ as follows:
$$\sigma(e_1,\ldots,e_N)=(e_{\sigma(1)},\ldots,e_{\sigma(N)})$$
In particular we have, as found before, the cardinality formula $|H_N^s|=s^NN!$.
\end{theorem}

\begin{proof}
As explained in the proof of Proposition 5.16, the elements of $H_N^s$ can be identified with the pairs $g=(e,\sigma)$ consisting of a permutation $\sigma\in S_N$, and a decorating vector $e\in\mathbb Z_s^N$, so that at the level of the cardinalities, we have:
$$|H_N|=|\mathbb Z_s^N\times S_N|$$

Now observe that the product formula for two such pairs $g=(e,\sigma)$ is as follows, with the permutations $\sigma\in S_N$ acting on the elements $f\in\mathbb Z_s^N$ as in the statement:
$$(e,\sigma)(f,\tau)=(ef^\sigma,\sigma\tau)$$

Thus, we are in the framework of the crossed products, and we obtain $H_N^s=\mathbb Z_s^N\rtimes S_N$. But this can be written, by definition, as $H_N^s=\mathbb Z_s\wr S_N$, and we are done.
\end{proof}

Finally, in relation with geometric aspects, the above groups appear as follows:

\begin{theorem}
The complex reflection group $H_N^s$ appears as a symmetry group,
$$H_N^s=G(C_s\ldots C_s)$$
with $C_s\ldots C_s$ consisting of $N$ disjoint copies of the oriented cycle $C_s$.
\end{theorem}

\begin{proof}
This is something elementary, the idea being as follows:

\medskip

(1) Consider first the oriented cycle $C_s$, which looks as follows:
$$
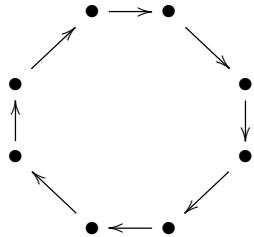
\xymatrix@R=15pt@C=16pt{
&\bullet\ar[r]&\bullet\ar[dr]\\
\bullet\ar[ur]&&&\bullet\ar[d]\\
\bullet\ar[u]&&&\bullet\ar[dl]\\
&\bullet\ar[ul]&\bullet\ar[l]}$$

It is then clear that the symmetry group of this graph is the cyclic group $\mathbb Z_s$.

\medskip

(2) In the general case now, where we have $N\in\mathbb N$ disjoint copies of the above cycle $C_s$, we must suitably combine the corresponding $N$ copies of the cyclic group $\mathbb Z_s$. But this leads to the wreath product group $H_N^s=\mathbb Z_s\wr S_N$, as stated.
\end{proof}

We will be back to this, graph symmetries, at the end of the present chapter.

\section*{5c. Abstract reflections}

As a next objective, in analogy with our discussion before of $A_N\subset S_N$ and of its analogue $SH_N\subset H_N$, we must discuss now the groups $SH_N^s\subset H_N^s$, obtained via the determinant 1 condition. Things here are quite straightforward, as follows:

\begin{theorem}
The determinant $\det:U_N\to\mathbb C$ restricts to a group morphism
$$\varepsilon:H_N^s\to\pm\mathbb Z_s$$
that we can call signature, whose kernel is the following subgroup $SH_N^s\subset H_N^s$:
$$SH_N^s=H_N^s\cap SO_N$$
When $s$ is even this subgroup has index $s$, and when $s$ is odd, the index is $2s$.
\end{theorem}

\begin{proof}
This is indeed quite clear, save for the parity considerations, as follows:

\medskip

(1) To start with, the determinant map $\det:U_N\to\mathbb C$ restricts indeed to the subgroup $H_N^s\subset U_N$ into a certain group morphism $\varepsilon:H_N^s\to\pm\mathbb Z_s$, with the usual convention that $\mathbb Z_s\subset\mathbb T$ is the group formed by the $s$-th roots of unity, so that the bigger group $\pm\mathbb Z_s\subset\mathbb T$ is this group itself, when $s$ is even, and is the group $\mathbb Z_{2s}$, when $s$ is odd.

\medskip

(2) Next, in analogy with what we already know at $s=2$, regarding the group $SH_N$, we can say that the signature $\varepsilon:H_N^s\to\pm\mathbb Z_s$ computes the signature of the underlying permutation $\sigma\in S_N$, multiplied by the product of the various $\mathbb Z_s$ decorations.

\medskip

(3) Then, we can talk about the group $SH_N^s=H_N^s\cap SO_N$ in the statement, which is indeed the kernel of $\varepsilon:H_N^s\to\pm\mathbb Z_s$, and which has index $s$ or $2s$, depending on the parity of $s$, as indicated, and with this coming from the parity discussion in (1).

\medskip

(4) Finally, observe that at $s=1$ we obtain the alternating group, $SH_N^1=A_N$, and that at $s=2$ we obtain the special hyperoctahedral group, $SH_N^2=SH_N$.
\end{proof}

Quite remarkably, the above result is not the end of the story, because we can use powers of the determinant, in order to construct new subgroups of $H_N^s$:

\begin{theorem}
We have subgroups of the basic complex reflection groups,
$$H_N^{s,d}\subset H_N^s$$
with $N\in\mathbb N$ and $s\in\mathbb N\cup\{\infty\}$, constructed via the following condition, with $d|[2,s]$,
$$(\det U)^d=1$$
which cover all the complex reflection groups that we have so far.
\end{theorem}

\begin{proof}
This is something quite self-explanatory, save for the discussion of the examples, and more specifically of the group $D_N$, with the details being as follows:

\medskip

(1) To start with, let us recall from Theorem 5.19 that the signature function on the group $H_N^s$ appears as follows, with $[2,s]$ standing for a least common multiple:
$$\varepsilon:H_N^s\to\mathbb Z_{[2,s]}$$

Thus, for any $d|[2,s]$, we can define a subgroup $H_N^{s,d}\subset H_N^s$ as follows:
$$H_N^{s,d}=\left\{U\in H_N^s\Big|(\det U)^d=1\right\}$$

(2) Regarding now the basic examples of such groups, we have the following equalities, and with these groups covering pretty much everything that we know:
$$H_N^s=H_N^{s,[2,s]}\quad,\quad SH_N^s=H_N^{s,1}$$

(3) In fact, the only group which needs some discussion is $D_N$. For this purpose, let us recall from Theorem 5.13 that this group can be viewed as follows:
$$D_N=\left\{\begin{pmatrix}x&0\\0&y\end{pmatrix}\ ,\ 
\begin{pmatrix}0&x\\y&0\end{pmatrix}\Big|x=y^{-1}\in\mathbb Z_N\right\}\subset K_2$$

We can see that we have $D_N\subset H_2^N$, and in fact even more, we have an embedding as follows, coming from the fact that the determinant of the above matrices is $\pm1$:
$$D_N\subset H_2^{N,2}$$

(4) The problem is now, is this embedding an isomorphism, or not. We have:
$$H_2^{N,2}=\left\{U\in H_2^N\Big|(\det U)^2=1\right\}$$

Now when $N$ is odd we have $-1\notin\mathbb Z_N$, so the determinant condition above can be interpreted as follows, with $\square$ standing for the product of nonzero entries:
$$(\det U)^2=1\iff\det U=\pm1\iff\square\, U=1$$

We conclude that when $N$ is odd, we have an equality as follows, as desired:
$$D_N=H_2^{N,2}$$

(5) In the remaining case, $N$ even, we have $-1\in\mathbb Z_N$, and the above reasoning will not work. In fact, in this case we have $H_2^{N,2}=D_N\sqcup E_N$, with $E_N$ being given by:
$$E_N=\left\{\begin{pmatrix}x&0\\0&-y\end{pmatrix}\ ,\ 
\begin{pmatrix}0&x\\-y&0\end{pmatrix}\Big|x=y^{-1}\in\mathbb Z_N\right\}$$

Thus, in this case, wrong cardinality, and we must find something else.

\medskip

(6) As an idea for a fix, let us look at the subgroup $H_2^{N,1}\subset H_2^{N,2}$, which has the correct cardinality. This subgroup appears via the $\det U=1$ condition, so is given by:
$$H_2^{N,1}=\left\{\begin{pmatrix}x&0\\0&y\end{pmatrix}\ ,\ 
\begin{pmatrix}0&x\\-y&0\end{pmatrix}\Big|x=y^{-1}\in\mathbb Z_N\right\}$$

So, let us try to redo the proof of Theorem 5.13, in this setting, with that $-$ sign added. We can label our group elements as follows, with $xy=zt=1$:
$$R_x=\begin{pmatrix}x&0\\0&y\end{pmatrix}\quad,\quad 
S_z=\begin{pmatrix}0&z\\-t&0\end{pmatrix}$$

We have then the following computations, for the products of these elements:
$$R_xR_z=\begin{pmatrix}x&0\\0&y\end{pmatrix}
\begin{pmatrix}z&0\\0&t\end{pmatrix}
=\begin{pmatrix}xz&0\\0&yt\end{pmatrix}
=R_{xz}$$
$$R_xS_z=\begin{pmatrix}x&0\\0&y\end{pmatrix}
\begin{pmatrix}0&z\\-t&0\end{pmatrix}
=\begin{pmatrix}0&xz\\-yt&0\end{pmatrix}
=S_{xz}$$
$$S_xR_z=
\begin{pmatrix}0&x\\-y&0\end{pmatrix}
\begin{pmatrix}z&0\\0&t\end{pmatrix}\
=\begin{pmatrix}0&xt\\-yz&0\end{pmatrix}
=S_{xz^{-1}}$$
$$S_xS_z=\begin{pmatrix}0&x\\-y&0\end{pmatrix}
\begin{pmatrix}0&z\\-t&0\end{pmatrix}
=\begin{pmatrix}-xt&0\\0&-yz\end{pmatrix}
=R_{-xz^{-1}}$$

Which is very frustrating, everything fine, save for that $-$ sign at the very end.

\medskip

(7) Nevermind, and as a second idea for a fix, by taking advantage of the fact that $N$ is even, let us look at the following group, which again has the correct cardinality:
$$H_2^{N/2,2}=\left\{U\in H_2^{N/2}\Big|(\det U)^2=1\right\}$$

This subgroup appears via the $\det U=\pm1$ condition, so is given by:
$$H_2^{N/2,2}=\left\{\begin{pmatrix}x&0\\0&\pm y\end{pmatrix}\ ,\ 
\begin{pmatrix}0&x\\\pm y&0\end{pmatrix}\Big|x=y^{-1}\in\mathbb Z_{N/2}\right\}$$

However, the bad news is that the matrices on the left have order $[2,N/2]$, instead of $N$ as wished. Thus, most likely wrong way too, and we should probably stop here.

\medskip

(8) This being said, let us see what we get, when $N/2$ is odd, so that:
$$[2,N/2]=N$$

We can label our group elements as follows, with $xy=zt=1$, and $a,b=\pm1$:
$$R_{ax}=\begin{pmatrix}ax&0\\0&y\end{pmatrix}\quad,\quad 
S_{bz}=\begin{pmatrix}0&bz\\t&0\end{pmatrix}$$

We have then the following computations, for the products of these elements:
$$R_{ax}R_{bz}=\begin{pmatrix}ax&0\\0&y\end{pmatrix}
\begin{pmatrix}bz&0\\0&t\end{pmatrix}
=\begin{pmatrix}axbz&0\\0&yt\end{pmatrix}
=R_{axbz}$$
$$R_{ax}S_{bz}=\begin{pmatrix}ax&0\\0&y\end{pmatrix}
\begin{pmatrix}0&bz\\t&0\end{pmatrix}
=\begin{pmatrix}0&axbz\\yt&0\end{pmatrix}
=S_{axbz}$$
$$S_{ax}R_{bz}=
\begin{pmatrix}0&ax\\y&0\end{pmatrix}
\begin{pmatrix}bz&0\\0&t\end{pmatrix}\
=\begin{pmatrix}0&axt\\ybz&0\end{pmatrix}
=S_{axbz^{-1}}$$
$$S_{ax}S_{bz}=\begin{pmatrix}0&ax\\y&0\end{pmatrix}
\begin{pmatrix}0&bz\\t&0\end{pmatrix}
=\begin{pmatrix}axt&0\\0&ybz\end{pmatrix}
=R_{axbz^{-1}}$$

Which is a win, providing us with an isomorphism as follows, for $N/2$ odd:
$$D_N=H_2^{N/2,2}$$

(9) As for the continuation of this, this is straightforward, with the next case being $N/4$ odd, then $N/8$ odd, and so on, discussed in a similar way. Thus, the dihedral groups $D_N$ are indeed covered by our construction, as claimed in the statement.
\end{proof}

Still with me I hope, after all these computations. Good work that we did, with all our previous examples of reflection groups being covered by our $H_N^{s,d}$ construction.

\bigskip

This being said, I don't know about you, but personally I find quite frustrating that the extra parameter is subject to $d|[2,s]$ instead of $d|s$, and also that the original construction of $D_N\subset K_2$, from Theorem 5.13, which was something quite simple and bright, had to get twisted as above, into something quite technical, for our covering purposes here.

\bigskip

As always in such difficult situations, time to ask the cat. And cat says:

\begin{cat}
Among the correct and the beautiful, I always chose the beautiful.
\end{cat}

Which sounds quite wise, and no wonder about this, cats are usually pretty good at groups, mathematics and physics, guess mine must be a reincarnation of Hermann Weyl. So, going now with the beauty, as cat suggests, and too bad if we scratch a bit our mathematical classification result, coming from the above, let us formulate:

\begin{theorem}
We have subgroups of the basic complex reflection groups,
$$H_N^{sd}\subset H_N^s$$
constructed via the following condition, $\square$ being the product of nonzero entries
$$(\square\, U)^d=1$$
which cover most of the complex reflection groups that we have so far.
\end{theorem}

\begin{proof}
As mentioned, this will be a remake of Theorem 5.21, with the groups $H_N^{s,d}$ with $d|[2,s]$ from there replaced by the groups $H_N^{sd}$ with $d|s$ from here, and with everything or almost meant to get simpler. The details for all this are as follows:

\medskip

(1) To start with, with $\square$ being as above, we have a group morphism as follows:
$$\square:H_N^s\to\mathbb Z_s$$

Thus, for any $d|s$, we can define a subgroup $H_N^{sd}\subset H_N^s$ as follows:
$$H_N^{sd}=\left\{U\in H_N^s\Big|(\square\, U)^d=1\right\}$$

(2) Regarding now the basic examples of such groups, these are as follows:
$$H_N^s=H_N^{ss}$$

Observe that the groups $SH_N^s$ are no longer covered, or at least are not covered in an obvious way. Nevermind, and we will leave some thinking here as an exercise.

\medskip

(3) Next, and getting now to the key part, recall from Theorem 5.13 that:
$$D_N=\left\{\begin{pmatrix}x&0\\0&y\end{pmatrix}\ ,\ 
\begin{pmatrix}0&x\\y&0\end{pmatrix}\Big|x=y^{-1}\in\mathbb Z_N\right\}\subset K_2$$

But this translates into $D_N=H_2^{N1}$, which is obviously something very nice.
\end{proof}

As a conclusion to all this, good work that we did, and we will stop here with our construction of complex reflection groups, due to a famous classification result of Shephard and Todd, that we would like to explain now. To start with, we can talk about complex reflections and about complex reflection groups abstractly, as follows:

\begin{definition}
We can talk about reflections and reflection groups, as follows:
\begin{enumerate}
\item A reflection is a symmetry $S\in U_N$ with respect to a hyperplane $P\subset\mathbb C^N$.

\item A reflection group is a group $G\subset U_N$ generated by reflections, $G=<S_i>$.

\item Such a reflection group is called irreducible when it has no invariant subspaces.
\end{enumerate}
\end{definition}

Observe that we have not assumed $G$ to be finite, in the above, and with this making the above formalism quite broad, for instance with many continuous groups $G\subset U_N$ being reflection groups, in the above sense. Still in this setting, with no finiteness assumption on $G$, these reflection groups are best investigated by writing them as follows:
$$G=\left<S_1,\ldots,S_n\Big|(S_iS_j)^{m_{ij}}=1\right>$$ 

And there has been a lot of work here, by Coxeter and others. Getting now to the finite group case, any reflection group appears as product of irreducible reflection groups, and in what regards these latter groups, we have the following classification result:

\index{complex reflection group}
\index{Shephard-Todd}

\begin{theorem}
The irreducible complex reflection groups are
$$H_N^{sd}=\left\{U\in H_N^s\Big|(\square\, U)^d=1\right\}$$
along with $34$ exceptional examples.
\end{theorem}

\begin{proof}
This is something quite advanced, that we will not attempt to prove here, or even explain in detail, with the list of 34 exceptional cases, and we refer here to the paper of Shephard and Todd \cite{sto}, and to the subsequent literature on the subject. Be said in passing, observe that cat was right, with $H_N^{sd}$ being indeed the correct objects.
\end{proof}

As a conclusion to this, certainly done with the construction of basic reflection groups, with $K_N$ being the master group, containing everyone else, then with the groups $H_N^s$ being the basic objects, and then with $H_N^{sd}$ being more specialized objects. And with the extra remark that, by Shephard-Todd, there are 34 exceptional examples as well.

\bigskip

Going now beyond the reflection group case, we can use the above technology in order to say things about the arbitrary subgroups $G\subset U_N$, finite or not. Indeed, we can talk about the reflection subgroup of any subgroup $G\subset U_N$, as follows: 

\begin{definition}
Given $G\subset U_N$, we can define its reflection subgroup to be
$$K=G\cap K_N$$
with the intersection taken inside $U_N$.
\end{definition}

Many things can be said in relation with this, but let us not stop here. Indeed, given an intermediate subgroup $H_N\subset G\subset U_N$, we can view it as follows:
$$
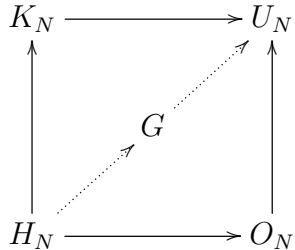
\xymatrix@R=25pt@C=25pt{
K_N\ar[rr]&&U_N\\
&G\ar@.[ur]&\\
H_N\ar[rr]\ar[uu]\ar@.[ur]&&O_N\ar[uu]}$$

Thus, we have some sort of 2D orientation for the subgroups $H_N\subset G\subset U_N$, and this suggests extending the construction in Definition 5.25, in the following way:

\begin{definition}
Associated to any intermediate compact group $H_N\subset G\subset U_N$ are its discrete, real, complex and smooth versions, given by the formulae
$$G^d=G\cap K_N\quad,\quad 
G^r=G\cap O_N$$
$$G^c=<G,K_N>\quad,\quad 
G^s=<G,O_N>$$
with $<\,,>$ being the topological generation operation, involving taking a closure.
\end{definition}

But with this in hand, it is natural now to formulate the following definition:

\index{oriented group}

\begin{definition}
A compact group $H_N\subset G\subset U_N$ is called oriented if
$$\xymatrix@R=41pt@C=45pt{
K_N\ar[r]&G^c\ar[r]&U_N\\
G^d\ar[u]\ar[r]&G\ar[r]\ar[u]&G^s\ar[u]\\
H_N\ar[r]\ar[u]&G^r\ar[u]\ar[r]&O_N\ar[u]}$$
is an intersection and generation diagram, in the sense that any of its square subdiagrams $A\subset B,C\subset D$ satisfies $A=B\cap C$ and $D=<B,C>$.
\end{definition}

And this notion is quite interesting, because most of the basic examples of closed subgroups $G\subset U_N$, finite or continuous, are oriented. In fact, we have:

\begin{question}
What are the oriented groups $H_N\subset G\subset U_N$? What about the oriented groups coming in families, $G=(G_N)$, with $N\in\mathbb N$?
\end{question}

And we will stop here our discussion, sometimes a good question is better as hunting trophy than a final theorem, or at least that's what my cats say. We will be back to this in chapters 9-12 below, under a number of supplementary assumptions on the groups $G$ that we consider, which will allow us to derive a number of classification results.

\section*{5d. Graph products}

Getting back now to the basics, namely the groups $H_N^s$, which remain the core objects of complex reflection group theory, we recall from Theorem 5.18 that these appear as follows, with $C_s\ldots C_s$ consisting of $N$ disjoint copies of the oriented cycle $C_s$:
$$H_N^s=G(C_s\ldots C_s)$$

This fact, which is quite interesting, is very easy to generalize, as follows:

\begin{theorem}
Given a connected graph $X$, we have the formula
$$G(X\ldots X)=G(X)\wr S_N$$
with $X\ldots X$ consisting of $N$ disjoint copies of $X$.
\end{theorem}

\begin{proof}
This is something quite clear, as before for $X=C_s$, because the graph $X\ldots X$ has an obvious action of $G(X)^N$ on it, along with an action of $S_N$ too. And these copies of $G(X)^N$ and $S_N$ do not exactly commute, but are rather related as follows:
$$G(X\ldots X)=G(X)^N\rtimes S_N$$

Thus, we are led to the wreath product formula in the statement.
\end{proof}

As a continuation of this, we can have a look at more general products of graphs. Let us start our study with the following standard graph theory definition:

\index{direct product}
\index{Cartesian product}
\index{lexicographic product}

\begin{definition}
Given two finite graphs $X,Y$, we can construct:
\begin{enumerate}
\item The direct product $X\times Y$ has vertex set $X\times Y$, and edges:
$$(i,\alpha)-(j,\beta)\Longleftrightarrow i-j,\, \alpha-\beta$$

\item The Cartesian product $X\,\square\,Y$ has vertex set $X\times Y$, and edges:
$$(i,\alpha)-(j,\beta)\Longleftrightarrow i=j,\, \alpha-\beta\mbox{ \rm{or} }i-j,\alpha=\beta$$

\item The lexicographic product $X\circ Y$ has vertex set $X\times Y$, and edges:
$$(i,\alpha)-(j,\beta)\Longleftrightarrow\alpha-\beta\mbox{ \rm{or} }\alpha=\beta,\,
 i-j$$
\end{enumerate}
We call these three products the standard products of graphs.
\end{definition}

Several comments can be made, in relation with this. First, the direct product $X\times Y$ is the usual one in a categorical sense, and we will leave clarifying this observation as an exercise. As for the Cartesian product $X\,\square\,Y$, this is quite natural too from a geometric perspective, for instance because a product by a segment gives a prism. 

\bigskip

Regarding now the lexicographic product $X\circ Y$, this is something interesting too, obtained by definition by putting a copy of $X$ at each vertex of $Y$. And with the basic example here, making the link with Theorem 5.29, being as follows:
$$X\ldots X=X\circ[\ldots]$$

Summarizing, we would like to have an extension of Theorem 5.29, in the general framework of Definition 5.30. Let us start with something elementary, as follows:

\index{wreath product}

\begin{theorem}
We have group embeddings as follows, for any graphs $X,Y$,
$$G(X)\times G(Y)\subset G(X \times Y)$$
$$G(X)\times G(Y)\subset G(X\,\square\,Y)$$
$$G(X)\wr G(Y)\subset G(X\circ Y)$$
but these embeddings are not always isomorphisms.
\end{theorem}

\begin{proof}
The fact that we have indeed embeddings as above is clear from definitions. As for the counterexamples, in each case, these are easy to construct as well, coming from a quick study of small graphs, at $N=2,3,4,\ldots\,$, that we will leave as an exercise.
\end{proof}

The problem now is that of deciding when the embeddings in Theorem 5.31 are isomorphisms. And in answer to this, we first have the following result:

\index{connected graph}
\index{regular graph}

\begin{theorem}
Let $X$ and $Y$ be finite connected regular graphs. If their spectra $\{\lambda\}$ and $\{\mu\}$ do not contain $0$ and satisfy
$$\big\{\lambda_i/\lambda_j\big\}\cap\big\{\mu_k/\mu_l\big\}=\{1\}$$
then $G(X\times Y)=G(X)\times G(Y)$. Also, if their spectra satisfy
$$\big\{\lambda_i-\lambda_j\big\}\cap\big\{\mu_k-\mu_l\big\}=\{0\}$$
then $G(X\,\square\,Y)=G(X)\times G(Y)$. 
\end{theorem}

\begin{proof}
This is something quite standard, the idea being as follows:

\medskip

(1) To start with, we know from Theorem 5.31 that we have embeddings as follows, valid for any two graphs $X,Y$, and coming from definitions:
$$G(X)\times G(Y)\subset G(X \times Y)$$
$$G(X)\times G(Y)\subset G(X\,\square\,Y)$$

(2) Now let $\lambda_1$ be the valence of $X$. Since $X$ is regular we have $\lambda_1\in Sp(X)$, with $1$ as eigenvector, and since $X$ is connected $\lambda_1$ has multiplicity 1. Thus if $P_1$ is the orthogonal projection onto $\mathbb C1$, the spectral decomposition of $d_X$ is of the following form:
$$d_X=\lambda_1P_1+\sum_{i\neq1}\lambda_iP_i$$

We have a similar formula for the adjacency matrix $d_Y$, namely:
$$d_Y=\mu_1Q_1+\sum_{j\neq1}\mu_jQ_j$$

(3) But this gives the following formulae for the first two graph products:
$$d_{X\times Y}=\sum_{ij}(\lambda_i\mu_j)P_{i}\otimes Q_{j}$$
$$d_{X\,\square\,Y}=\sum_{ij}(\lambda_i+\mu_i)P_i\otimes Q_j$$

Here the projections form partitions of unity, and the scalar are distinct, so these are spectral decompositions. Now since group actions commute with the spectral projections, these actions will commute with both $P_1\otimes1$, $1\otimes Q_1$, and this gives the results.
\end{proof}

Regarding the lexicographic products, we have a similar result here, as follows:

\begin{theorem}
Let $X,Y$ be regular graphs, with $X$ connected. If their spectra $\{\lambda_i\}$ and $\{\mu_j\}$ satisfy the condition
$$\big\{\lambda_1-\lambda_i\big|i\neq 1\big\}\cap\big\{-n\mu_j\big\}=\emptyset$$
where $n$ and $\lambda_1$ are the order and valence of $X$, then $G(X\circ Y)=G(X)\wr G(Y)$.   
\end{theorem}

\begin{proof}
This is something more tricky, the idea being as follows:

\medskip

(1) To start with, we know from Theorem 5.31 that we have an embedding as follows, valid for any two graphs $X,Y$, and coming from definitions:
$$G(X)\wr G(Y)\subset G(X\circ Y)$$

(2) But then, we can use the same method as in the proof of Theorem 5.32, in order to show that we have equality. However, this is something a bit more technical, and for the story and details here, you can have a look at my graph book \cite{ba1}.

\medskip

(3) Finally, let us mention that this generalizes indeed Theorem 5.29, with the spectral assumption in the statement being trivially verified, in the present case.
\end{proof}

As a conclusion to all this, Theorems 5.32 + 5.33 generalize Theorem 5.29, which in turn generalizes our original computation, $H_N^s=G(C_s\ldots C_s)$, from Theorem 5.18. Talk about pure mathematics, and generalizing things, whether needed or not.

\bigskip

More seriously now, the above results come in fact from work that I did with Julien Bichon in the mid 2000s, motivated by the understanding of the quantum symmetry groups $G^+(X)\supset G(X)$ of the transitive graphs $X$ having small number of vertices. We will not talk about such things here, but instead, we can have a look at:

\begin{question}
What can we say about the symmetry groups $G(X)$ of the transitive graphs $X$ having small number of vertices? Do these always decompose as products?
\end{question}

In answer now to this, certainly yes for graphs having very small number of vertices, namely $N=2,\ldots,9$. However, at $N=10$ we have the Petersen graph, as follows:
$$
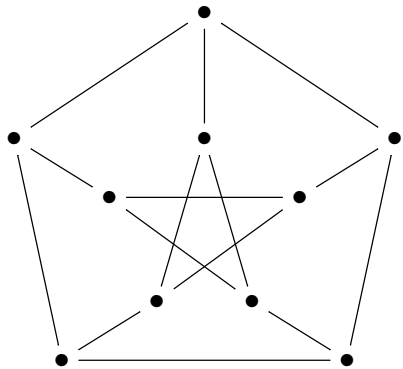
\xymatrix@R=1.5pt@C=5pt{
&&&&\bullet\ar@{-}[dddddrrrr]\ar@{-}[dddddllll]\\
\\
\\
\\
\\
\bullet\ar@{-}[ddddddddr]&&&&\bullet\ar@{-}[uuuuu]\ar@{-}[ddddddl]\ar@{-}[ddddddr]&&&&\bullet\ar@{-}[ddddddddl]\\
\\
&&
\bullet\ar@{-}[uull]\ar@{-}[ddddrrr]\ar@{-}[rrrr]&&&&\bullet\ar@{-}[uurr]\ar@{-}[ddddlll]\\
\\
\\
\\
&&&\bullet&&\bullet\\
\\
&\bullet\ar@{-}[rrrrrr]\ar@{-}[uurr]&&&&&&\bullet\ar@{-}[uull]}
$$

Intuition suggests that this graph should appear as a product of the 5-cycle with the segment. However, there are many bugs with this idea, which does not work. In order to view the Petersen graph as part of a larger family, and have some general theory going for it, we have to proceed in a quite unexpected way, as follows:

\begin{proposition}
The Petersen graph is part of the Kneser graph family,
$$P_{10}=K(5,2)$$
with $K(n,s)$ having as vertices the $s$-element subsets of $\{1,\ldots,n\}$, and with the edges being drawn between distinct subsets.
\end{proposition}

\begin{proof}
Consider indeed the Kneser graph $K(n,s)$, as constructed above. At $n=5$, $s=2$ the vertices are the $\binom{5}{2}=10$ subsets of $\{1,\ldots,5\}$ having 2 elements, and since each such subset $\{p,q\}$ is disjoint from exactly 3 other such subsets $\{u,v\}$, our graph is trivalent, and when drawing the picture, we obtain indeed the Petersen graph.
\end{proof}

Many things can be said about the Kneser graphs, and among others, we have:

\index{Petersen graph}
\index{Kneser graph}

\begin{theorem}
The Kneser graphs $K(n,s)$ have the following properties:
\begin{enumerate}
\item $K(n,1)$ is the complete graph $K_n$.

\item $K(n,2)$ is the complement of the line graph $L(K_n)$.

\item $K(2m-1,m-1)$ is the so-called odd graph $O_m$.

\item The symmetry group of $K(n,s)$ is the symmetric group $S_n$.
\end{enumerate}
In particular, $P_{10}=K(5,2)=L(K_5)^c=O_3$, having symmetry group $S_5$.
\end{theorem}

\begin{proof}
All this is quite self-explanatory, with (1) being trivial, (2) being clear too, with the line graph $L(X)$ of a given graph $X$ being by definition the incidence graph of the edges of $X$, then (3) standing for the definition of $O_m$, as being the graph $K(2m-1,m-1)$, and finally (4) being clear too, coming from the definition of the Kneser graphs.
\end{proof}

Many other things can be said, as a continuation of this. For more, including the story with the corresponding quantum symmetry groups $G^+(X)\supset G(X)$, have a look at \cite{ba1}.

\section*{5e. Exercises}

This was a quite exciting algebraic chapter, and as exercises on this, we have:

\begin{exercise}
Establish, with bare hands, the formula $|H_3|=48$.
\end{exercise}

\begin{exercise}
Meditate as well on the isomorphism $SH_3=S_4$.
\end{exercise}

\begin{exercise}
Compute the groups $H_N^s$, at small values of $N,s$.
\end{exercise}

\begin{exercise}
Compute the groups $H_N^{sd}$, at small values of $N,s,d$.
\end{exercise}

\begin{exercise}
Clarify the relation between the series $\{H_N^{s,d}\}$ and $\{H_N^{sd}\}$.
\end{exercise}

\begin{exercise}
Learn more about the classification result of Shephard and Todd.
\end{exercise}

\begin{exercise}
Compute the reflection groups of the various groups $G\subset U_N$.
\end{exercise}

\begin{exercise}
Compute the symmetries of graphs with small number of vertices.
\end{exercise}

As bonus exercise, learn as well about regular polyhedra, and their symmetries.

\chapter{Braid groups}

\section*{6a. Knots and links}

Time for a new topic, topology, and knots and links. However, as we will soon discover, when looking at this, groups are just around the corner. The relevant group $B_k$, called braid group, will be actually infinite. But, a bit like the reflection groups appear as ``versions'' of the symmetric group $S_k$, so will appear $B_k$, making it somehow belong to the finite group philosophy, and have its place here, at this point in this book.

\bigskip

Getting started, we already talked a bit about topology in chapter 1, with a discussion about homotopy groups, in relation with the free groups $F_N$. However, such things are rather advanced, and at the very basic level, topology fundamentals, we have:

\index{knot}

\begin{definition}
A knot is a smooth closed curve in $\mathbb R^3$,
$$\gamma:\mathbb T\to\mathbb R^3$$
regarded modulo smooth transformations of $\mathbb R^3$.
\end{definition}

Observe that our knots are by definition oriented. The reverse knot $z\to\gamma(z^{-1})$ can be isomorphic or not to the original knot $z\to\gamma(z)$, and we will discuss this in a moment. At the level of examples, we first have the unknot, represented as follows:
$$
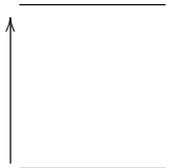
\xymatrix@R=55pt@C=55pt{
\ar@{-}[r]&\ar@{-}[d]\\
\ar[u]&\ar@{-}[l]}$$

The unknot is already a quite interesting mathematical object, suggesting a lot of exciting mathematical questions, for the most quite difficult, as follows:

\index{random knot}
\index{tied knot}

\begin{questions}
In relation with the unknot:
\begin{enumerate}
\item Given a closed curve in $\mathbb R^3$, via equations, can we decide if it is tied or not? 

\item Given the 2D picture of a knot, can we decide if the knot is tied or not?

\item A random cable or rope is usually tied. But, can we really prove this?
\end{enumerate}
\end{questions}

Obviously, difficult questions, and as you can see, knot theory is not an easy thing. But do not worry, we will manage to find our way through this jungle, and even come up with some mathematics for it. Going ahead now with examples, as the simplest possible true knot, meaning tied knot, we have the trefoil knot, which looks as follows:
$$
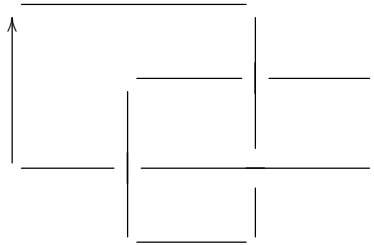
\xymatrix@R=15pt@C=35pt{
\ar@{-}[rr]&&\ar@{-}[dd]\\
&\ar@{-}[r]&|\ar@{-}[r]&\\
\ar@{-}[r]\ar[uu]&|\ar@{-}[rr]&-&\ar@{-}[u]\\
&\ar@{-}[r]\ar@{-}[uu]&\ar@{-}[u]
}$$

We also have the opposite trefoil knot, obtained by reversing the orientation, whose picture is identical to that of the trefoil knot, save for the orientation of the arrow:
$$
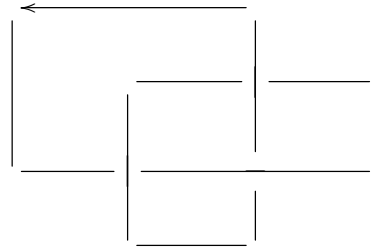
\xymatrix@R=15pt@C=35pt{
&&\ar[ll]\ar@{-}[dd]\\
&\ar@{-}[r]&|\ar@{-}[r]&\\
\ar@{-}[r]\ar@{-}[uu]&|\ar@{-}[rr]&-&\ar@{-}[u]\\
&\ar@{-}[r]\ar@{-}[uu]&\ar@{-}[u]
}$$

As before with the unknot, while the trefoil knot might look quite trivial, when it comes to formal mathematics regarding it, we are quickly led into delicate questions. Let us formulate a few intuitive observations about it, as follows:

\begin{fact}
In relation with the trefoil knot:
\begin{enumerate}
\item This knot is indeed tied, that is, not isomorphic to the unknot.

\item The trefoil knot and its opposite knot are not isomorphic.
\end{enumerate}
\end{fact}

To be more precise, here (1) is something which definitely holds, as we know it from real life, but if looking for a formal proof for this, based on Definition 6.1, we will certainly run into troubles. As for (2), here again we are looking for troubles, because when playing with two trefoil knots, made from rope, with opposite arrows marked on them, we can certainly see that our two beasts are not identical, but go find a formal proof for that.

\bigskip

In short, as before with the unknot, modesty. For the moment, let us keep exploring the subject, by recording as Questions and Facts things that we see and feel, but cannot prove yet, mathematically, based on Definition 6.1 alone, due to a lack of tools.

\bigskip

Getting back now to Definition 6.1, as stated, it is convenient to allow, in relation with certain mathematical questions, links in our discussion:

\index{link}

\begin{definition}
A link is a collection of disjoint knots in $\mathbb R^3$, taken as usual oriented, and regarded as usual up to isotopy.
\end{definition}

As before with the knots, which can be truly knotted or not, there is a discussion here with respect to the links, which can be truly linked or not, and with orientation involved too. Drawing some pictures here, with some basic examples, is very instructive, the idea being that two or several basic unknots can be linked in many possible ways. For instance, as simplest non-trivial link, made of two unknots, which are indeed linked, we have: 
$$
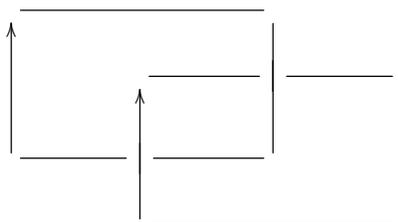
\xymatrix@R=12pt@C=40pt{
\ar@{-}[rr]&&\ar@{-}[dd]\\
&\ar@{-}[r]&|\ar@{-}[r]&\\
\ar@{-}[r]\ar[uu]&|\ar@{-}[r]&&\\
&\ar@{-}[rr]\ar[uu]&&\ar@{-}[uu]
}$$

By reversing the orientation of one unknot, we have as well the following link:
$$
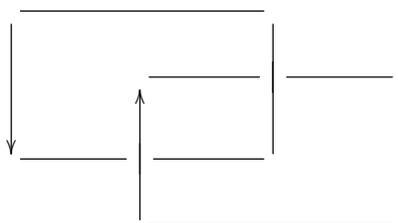
\xymatrix@R=12pt@C=40pt{
\ar@{-}[rr]\ar[dd]&&\ar@{-}[dd]\\
&\ar@{-}[r]&|\ar@{-}[r]&\\
\ar@{-}[r]&|\ar@{-}[r]&&\\
&\ar@{-}[rr]\ar[uu]&&\ar@{-}[uu]
}$$

This was for the story of two linked unknots, which is easily seen to stop here, with the above two links, but when trying to link $N$ unknots, with $N=3,4,5,\ldots\,$, many things can happen. Which leads us into the following philosophical question:

\begin{question}
Mathematically speaking, which are simpler to enumerate,
\begin{enumerate}
\item Usual knots, that is, links with one component,

\item Or links with several components, all being unknots,
\end{enumerate}
and this, in order to have some business started, for the links?
\end{question}

And with this being probably enough, as preliminary experimental work, time to draw some conclusions. Obviously, what we have so far, namely Questions 6.2, Fact 6.3 and Question 6.5, is very interesting, at the core of everything that can be called ``geometry and topology''. And by further thinking a bit, at how knots and links can be tied, in so many fascinating ways, we are led to the following philosophical conclusion:

\begin{conclusion}
Knots and links are to geometry and topology what prime numbers are to number theory.
\end{conclusion}

Which is very nice, but the question is, how to get started? Fortunately, graphs come to the rescue, via to the following simple fact, which will be our starting point:

\index{projection of knot}
\index{tetravalent graph}
\index{crossings of knot}

\begin{fact}
The plane projection of a knot or link is similar to an oriented graph with $4$-valent vertices, except for the fact that we have some extra data at the vertices,
$$
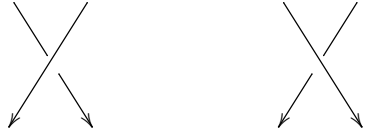
\xymatrix@C=10pt@R=20pt{
\ar@{-}[dr]&&\ar[ddll]&&&&\ar[ddrr]&&\ar@{-}[dl]\\
&\ar[dr]&&&&&&\ar[dl]&\\
&&&&&&&&}$$
telling us, about the $2$ strings crossing there, which goes up and which goes down.
\end{fact}

Based on this, let us try now to construct some knot invariants. A natural idea is that of defining the invariant on the 2D picture of the knot, that is, on a plane projection of the knot, and then proving that the invariant is indeed independent on the chosen plane. This method rests on the following technical result, which is well-known:

\index{Reidemeister moves}

\begin{theorem}
Two pictures correspond to plane projections of the same knot or link precisely when they differ by a sequence of Reidemeister moves, namely:
\begin{enumerate}
\item Moves of type I, given by $\propto\ \leftrightarrow\,|$.

\item Moves of type II, given by $)\hskip-1.7mm(\ \leftrightarrow\ )($.

\item Moves of type III, given by $\vartriangle\ \leftrightarrow\,\triangledown$.
\end{enumerate}
\end{theorem}

\begin{proof}
This is something very standard, as follows:

\medskip

(1) To start with, the Reidemeister moves of type I are by definition as follows:
$$
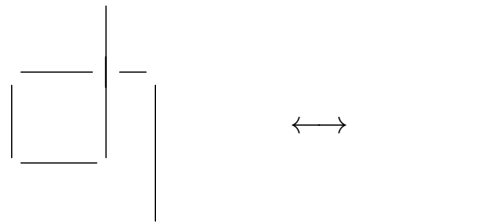
\xymatrix@R=5pt@C=10pt{
&&\ar@{-}[dddd]&&&&&&&\ar@{-}[dddddd]\\
\\
\ar@{-}[rr]\ar@{-}[dd]&&|\ar@{-}[r]&\ar@{-}[dddd]\\
&&&&&&\longleftrightarrow\\
\ar@{-}[rr]&&\\
\\
&&&&&&&&&
}$$

(2) Regarding the Reidemeister moves of type II, these are by definition as follows:
$$
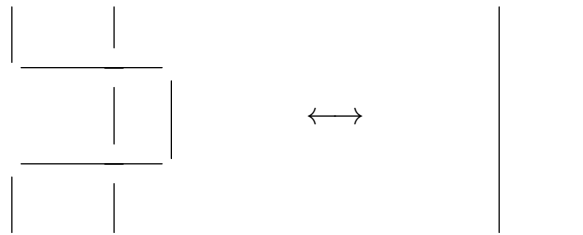
\xymatrix@R=5pt@C=10pt{
\ar@{-}[dd]&&\ar@{-}[dd]&&&&&&&\ar@{-}[dddddd]&&\ar@{-}[dddddd]\\
\\
\ar@{-}[rrr]&&-\ar@{-}[dd]&\ar@{-}[dd]\\
&&&&&&\longleftrightarrow\\
\ar@{-}[dd]\ar@{-}[rrr]&&-\ar@{-}[dd]&\\
\\
&&&&&&&&&&&
}$$

(3) As for the Reidemeister moves of type III, these are by definition as follows:
$$
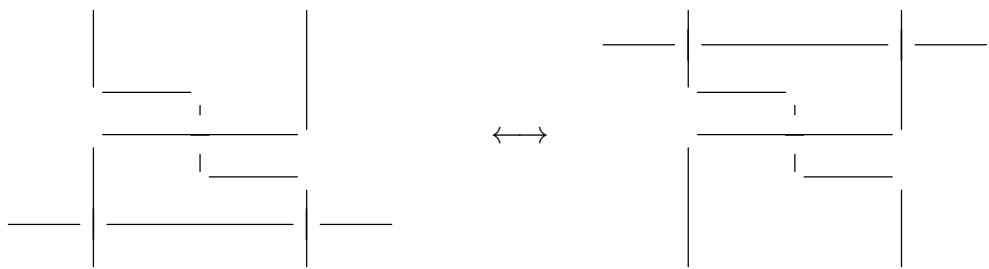
\xymatrix@R=5pt@C=10pt{
&&\ar@{-}[dd]&&&&\ar@{-}[ddd]&&&&&&&\ar@{-}[dd]&&&&\ar@{-}[ddd]\\
&&&&&&&&&&&\ar@{-}[rr]&&|\ar@{-}[rrrr]&&&&|\ar@{-}[rr]&&\\
&&\ar@{-}[rr]&&\ar@{-}[d]&&&&&&&&&\ar@{-}[rr]&&\ar@{-}[d]\\
&&\ar@{-}[rrrr]\ar@{-}[ddd]&&-\ar@{-}[d]&&&&&&\longleftrightarrow&&&\ar@{-}[rrrr]\ar@{-}[ddd]&&-\ar@{-}[d]&&\\
&&&&\ar@{-}[rr]&&\ar@{-}[dd]&&&&&&&&&\ar@{-}[rr]&&\ar@{-}[dd]\\
\ar@{-}[rr]&&|\ar@{-}[rrrr]&&&&|\ar@{-}[rr]&&\\
&&&&&&&&&&&&&&&&&
}$$

(4) This was for the precise statement of the theorem, and in what regards now the proof, this is somewhat clear from definitions, and in practice, this can be done by some sort of cut and paste procedure, or recurrence if you prefer, easy exercise for you.
\end{proof}

Summarizing, we have some mathematics started for the knots and links, and the idea in what follows will be that of defining invariants on 2D pictures, using Theorem 6.8.

\section*{6b. Braid groups}

At a more advanced level now, we will need the following key observation, making the connection with group theory, and algebra in general, due to Alexander:

\index{braid group}

\begin{theorem}
Any knot or link can be thought of as being the closure of a braid,
$$
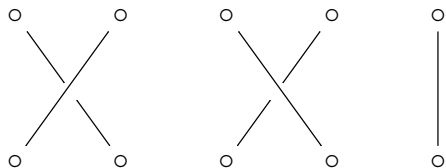
\xymatrix@C=10pt@R=18pt{
\circ\ar@{-}[dr]&&\circ\ar@{-}[ddll]&&\circ\ar@{-}[ddrr]&&\circ\ar@{-}[dl]&&\circ\ar@{-}[dd]\\
&\ar@{-}[dr]&&&&\ar@{-}[dl]&&&\\
\circ&&\circ&&\circ&&\circ&&\circ}
$$
with the braids forming a group $B_k$, called braid group.
\end{theorem}

\begin{proof}
Again, this is something quite self-explanatory, as follows:

\medskip

(1) Consider indeed the braids with $k$ strings, with the convention that things go from up to down. For instance the braid in the statement should be thought of as being:
$$
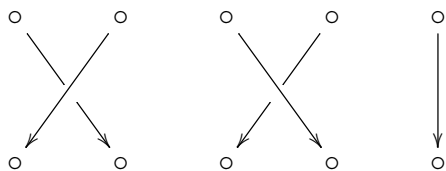
\xymatrix@C=10pt@R=18pt{
\circ\ar@{-}[dr]&&\circ\ar[ddll]&&\circ\ar[ddrr]&&\circ\ar@{-}[dl]&&\circ\ar[dd]\\
&\ar[dr]&&&&\ar[dl]&&&\\
\circ&&\circ&&\circ&&\circ&&\circ}
$$

But, with this convention, braids become some sort of permutations of $\{1,\ldots,k\}$, which are decorated at the level of crossings, with for instance the above braid corresponding to the following permutation of $\{1,2,3,4,5\}$, with due decorations: 
$$
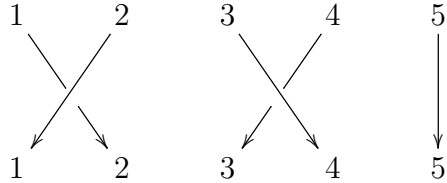
\xymatrix@C=10pt@R=18pt{
1\ar@{-}[dr]&&2\ar[ddll]&&3\ar[ddrr]&&4\ar@{-}[dl]&&5\ar[dd]\\
&\ar[dr]&&&&\ar[dl]&&&\\
1&&2&&3&&4&&5}
$$

In any case, we can see in this picture that $B_k$ is indeed a group, with composition law similar to that of the permutations in $S_k$, that is, going from up to down. 

\medskip

(2) Moreover, we can also see in this picture that we have a surjective group morphism $B_k\to S_k$, obtained by forgetting the decorations, at the level of crossings. For instance the braid pictured above is mapped in this way to the following permutation in $S_5$:
$$
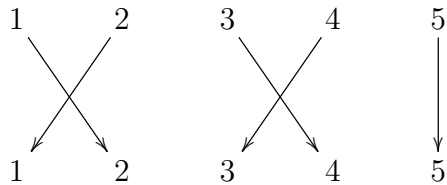
\xymatrix@C=10pt@R=18pt{
1\ar[ddrr]&&2\ar[ddll]&&3\ar[ddrr]&&4\ar[ddll]&&5\ar[dd]\\
&&&&&&&&\\
1&&2&&3&&4&&5}
$$

It is possible to do some more algebra here, in relation with the morphism $B_k\to S_k$, but we will not need this, right now. We will keep in mind, from the above, the fact that ``braids are not exactly permutations, but they compose like permutations''.

\medskip

(3) Regarding now the closure operation in the statement, this consists by definition in adding semicircles at right, which makes our braid into a certain oriented link. As an illustration, the closure of the braid pictured above is the following link: 
$$
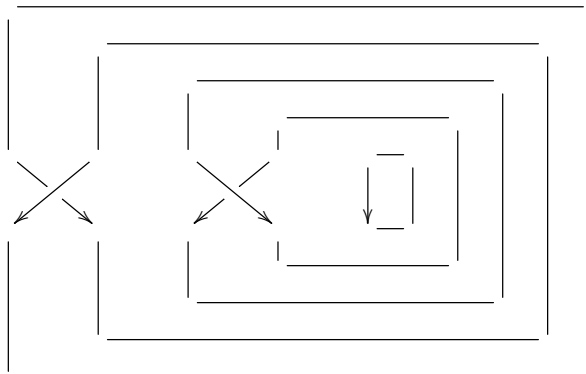
\xymatrix@C=10pt@R=7pt{
\ar@{-}[dddd]\ar@{-}[rrrrrrrrrrrrr]&&&&&&&&&&&&&\ar@{-}[dddddddddd]\\
&&\ar@{-}[ddd]\ar@{-}[rrrrrrrrrr]&&&&&&&&&&\ar@{-}[dddddddd]\\
&&&&\ar@{-}[dd]\ar@{-}[rrrrrrr]&&&&&&&\ar@{-}[dddddd]\\
&&&&&&\ar@{-}[d]\ar@{-}[rrrr]&&&&\ar@{-}[dddd]\\
\ar@{-}[dr]&&\ar[ddll]&&\ar[ddrr]&&\ar@{-}[dl]&&\ar[dd]\ar@{-}[r]&\ar@{-}[dd]\\
&\ar[dr]&&&&\ar[dl]&&&\\
&&&&&&&&\ar@{-}[r]&\\
&&&&&&\ar@{-}[u]\ar@{-}[rrrr]&&&&\\
&&&&\ar@{-}[uu]\ar@{-}[rrrrrrr]&&&&&&&\\
&&\ar@{-}[uuu]\ar@{-}[rrrrrrrrrr]&&&&&&&&&&\\
\ar@{-}[uuuu]\ar@{-}[rrrrrrrrrrrrr]&&&&&&&&&&&&&}
$$

(4) This was for the precise statement of the theorem, and in what regards now the proof, this can be done by some sort of cut and paste procedure, or recurrence if you prefer. As before with Theorem 6.8, we will leave this as an easy exercise for you.
\end{proof}

The question is now, what can we say about the braid group $B_k$, from a group-theoretic perspective? As a first observation, as mentioned above, we have a group morphism as follows, making $B_k$ a bit similar to the reflection groups from chapter 5:
$$B_k\to S_k$$

However, this similarity remains quite vague, because $B_k$ is infinite at $k\geq2$. Indeed, if we look at $B_2$, we can say that we have a vertical string on the right, and then a string on the left, which circles around it. But this circling can be $N$ times clockwise, say when looking from up, or $N$ times counterclockwise, which shows that we have:
$$B_2=\mathbb Z$$

Next, if we look at $B_3$, with some inspiration from $S_3=<XI,IX>$, we can say that this group is generated by two elements too, namely the preimages of $XI$, $IX$: 
$$\xymatrix@C=3pt@R=7pt{
&\circ\ar@{-}[dr]&&\circ\ar@{-}[ddll]&\circ\ar@{-}[dd]\\
g=&&\ar@{-}[dr]\\
&\circ&&\circ&\circ}\qquad\quad \qquad
\xymatrix@C=3pt@R=7pt{
&\circ\ar@{-}[dd]&\circ\ar@{-}[dr]&&\circ\ar@{-}[ddll]\\
h=&&&\ar@{-}[dr]&&&\\
&\circ&\circ&&\circ}$$

Thus we have a quotient map $F_2\to B_3$, and while the generators $g,h$ are not exactly free, being subject for instance to $ghg=hgh$, they are not very far from being free either, and it is possible in fact to come up with a group embedding, as follows:
$$F_2\subset B_3$$

At $k=4$ and higher, things can only become more complicated. So, let us summarize our conclusions so far in the form of an informal statement, as follows:

\begin{conclusion}
We can think of $B_k$ as being a decorated version of $S_k$, due to
$$B_k\to S_k$$
but with the warning that $B_k$ is infinite, and there is in fact a lot of freeness inside it.
\end{conclusion}

Moving on, with our various techniques from Part I in mind, you might perhaps wonder if we can say something about the order of the group elements $g\in B_k$. And here, things are quite tricky, because go compose a given $g\in B_k$ with itself several times, and see what you get, this does not look as an easy task. However, a few experiments here can of course be done, and based on this, let us record the following further conclusion:

\begin{conclusion}
The more you look at $B_k$, the more of infinite nature this appears, for instance with the order of the group elements $g\neq1$ being $ord(g)=\infty$.
\end{conclusion}

Summarizing, studying the braid group $B_k$ with bare hands is not an easy task, and what we have so far is rather a collection of conjectures. Fortunately, there is an answer to these difficulties, coming from the following key result, due to Artin:

\begin{theorem}
The braid group $B_k$ is generated by variables $g_1,\ldots,g_{k-1}$, with
$$g_ig_{i+1}g_i=g_{i+1}g_ig_{i+1}\quad,\quad 
g_ig_j=g_jg_i\ {\rm for}\ |i-j|\geq2$$
called Artin relations.
\end{theorem}

\begin{proof}
In order to prove this result, consider the following braids:
$$\xymatrix@C=3pt@R=6pt{
&\circ\ar@{-}[dr]&&\circ\ar@{-}[ddll]&\circ\ar@{-}[dd]&\circ\ar@{-}[dd]&&\circ\ar@{-}[dd]&\circ\ar@{-}[dd]\\
g_1=&&\ar@{-}[dr]&&&&\ldots&\\
&\circ&&\circ&\circ&\circ&&\circ&\circ}$$
$$\xymatrix@C=3pt@R=6pt{
&\circ\ar@{-}[dd]&\circ\ar@{-}[dr]&&\circ\ar@{-}[ddll]&\circ\ar@{-}[dd]&&\circ\ar@{-}[dd]&\circ\ar@{-}[dd]\\
g_2=&&&\ar@{-}[dr]&&&\ldots&\\
&\circ&\circ&&\circ&\circ&&\circ&\circ}$$
$$\ \ \ \ \vdots$$
$$\xymatrix@C=3pt@R=6pt{
&\circ\ar@{-}[dd]&\circ\ar@{-}[dd]&&\circ\ar@{-}[dd]&\circ\ar@{-}[dd]&\circ\ar@{-}[dr]&&\circ\ar@{-}[ddll]\\
\ \ \ \ \ \ \ g_{k-1}=&&&\ldots&&&&\ar@{-}[dr]\\
&\circ&\circ&&\circ&\circ&\circ&&\circ&&&&}$$

We have then $g_ig_j=g_jg_i$, for $|i-j|\geq2$. As for the relation $g_ig_{i+1}g_i=g_{i+1}g_ig_{i+1}$, by translation it is enough to check this at $i=1$. And here, we first have:
$$\xymatrix@C=3pt@R=7pt{
&\circ\ar@{-}[dr]&&\circ\ar@{-}[ddll]&&\circ\ar@{-}[dd]&\circ\ar@{-}[dd]&&\circ\ar@{-}[dd]&\circ\ar@{-}[dd]\\
&&\ar@{-}[dr]&&&&&\ldots&\\
&\circ\ar@{-}[dd]&&\circ\ar@{-}[dr]&&\circ\ar@{-}[ddll]&\circ\ar@{-}[dd]&&\circ\ar@{-}[dd]&\circ\ar@{-}[dd]\\
g_1g_2g_1=&&&&\ar@{-}[dr]&&&\ldots\\
&\circ\ar@{-}[dr]&&\circ\ar@{-}[ddll]&&\circ\ar@{-}[dd]&\circ\ar@{-}[dd]&&\circ\ar@{-}[dd]&\circ\ar@{-}[dd]\\
&&\ar@{-}[dr]&&&&&\ldots&\\
&\circ&&\circ&&\circ&\circ&&\circ&\circ}
$$

On the other hand, we have as well the following computation:
$$\xymatrix@C=3pt@R=7pt{
&\circ\ar@{-}[dd]&&\circ\ar@{-}[dr]&&\circ\ar@{-}[ddll]&\circ\ar@{-}[dd]&&\circ\ar@{-}[dd]&\circ\ar@{-}[dd]\\
&&&&\ar@{-}[dr]&&&\ldots\\
&\circ\ar@{-}[dr]&&\circ\ar@{-}[ddll]&&\circ\ar@{-}[dd]&\circ\ar@{-}[dd]&&\circ\ar@{-}[dd]&\circ\ar@{-}[dd]\\
g_2g_1g_2=&&\ar@{-}[dr]&&&&&\ldots&\\
&\circ\ar@{-}[dd]&&\circ\ar@{-}[dr]&&\circ\ar@{-}[ddll]&\circ\ar@{-}[dd]&&\circ\ar@{-}[dd]&\circ\ar@{-}[dd]\\
&&&&\ar@{-}[dr]&&&\ldots\\
&\circ&&\circ&&\circ&\circ&&\circ&\circ}$$

Now since the above two pictures are identical, up to isotopy, we have, as desired: 
$$g_1g_2g_1=g_2g_1g_2$$

Thus, the braid group $B_k$ is indeed generated by elements $g_1,\ldots,g_{k-1}$ with the relations in the statement, and in what regards now the proof of universality, this can only be something quite routine, and we will leave this as an instructive exercise.
\end{proof}

As further results now on the braid group, of various levels of difficulty, and that we will not really need in what follows, all coming via the Artin theorem, we have:

\begin{theorem}
The braid group $B_k$ has the following properties:
\begin{enumerate}
\item We have a group morphism $B_k\to\mathbb Z$ given by $g_i\to 1$ for any $i$, which is precisely the abelianization map, $B_k\to B_k/[B_k,B_k]$.

\item This morphism shows in particular that we have $ord(g_i)=\infty$, for any $i$. In fact, $B_k$ has no torsion, $g\neq1\implies ord(g)=\infty$.

\item The standard generators $g_1,g_2\in B_3$ are not free, but by suitably modifying them we can get an embedding $F_2\subset B_3$. Thus, $F_2\subset B_k$ for any $k\geq3$.

\item $B_k$ is the homotopy group of $X=(\mathbb C^k-\Delta)/S_k$, with $\Delta\subset\mathbb C^k$ standing for the set of points $z$ satisfying $z_i=z_j$ for some $i\neq j$.

\item $B_k$ is linear, in the sense that it admits an embedding $B_k\subset GL_n(\mathbb C)$.
\end{enumerate}
\end{theorem}

\begin{proof}
Tough theorem that we have here, among others solving our various conjectures before, and with most assertions being non-trivial. As mentioned above, we will not really need this in what follows, so here are a few details on this, for what this is worth, and in each case, exercise of course for you, to learn more about all this:

\medskip

(1) We have indeed a morphism $B_k\to\mathbb Z$ given by $g_i\to 1$, because $1_i\in\mathbb Z$ satisfy the relations in Theorem 6.12. As for the second assertion, this is something more tricky.

\medskip

(2) The first assertion is clear, because $g_i^n=1$ reads $n=0$ via the morphism $B_k\to\mathbb Z$ from (1). As for the second assertion, this is something more tricky.

\medskip

(3) This is again something tricky, the idea being that by suitably modifying $g_1,g_2$ we can get rid of the sole relation satisfied by them, namely $g_1g_2g_1=g_2g_1g_2$.

\medskip

(4) This is the only easy assertion, in all this, and we will leave some study here, coming as a continuation of our learning from chapter 1, as an exercise for you.

\medskip

(5) This is something truly tough, and quite recent, from the early 2000s, due to Bigelow and Krammer. As an exercise for you however, related to this, try proving first that the free group $F_2$ is linear, as a matter of getting used to the linear groups. And if you have difficulties, look up ``table tennis lemma'' on the internet, no joke here.
\end{proof}

\section*{6c. Markov theorem}

Getting back now to knots and links, a quick comparison between our main results so far, namely Theorem 6.8 due to Reidemeister, and then Theorem 6.9 due to Alexander, suggests the following question, whose answer will certainly advance us:

\begin{question}
What is the analogue of the Reidemeister theorem, in the context of braids? That is, when do two braids produce, via closing, the same link?
\end{question} 

And this is, and we insist, a very good question, because assuming that we have an answer to it, no need afterwards to bother with plane projections, decorated graphs, Reidemeister moves, and amateurish topology in general, it will be all about groups and algebra. Which groups and algebra questions, you guessed right, we will eat them raw. In answer now, we have the following theorem, due to Markov:

\begin{theorem}
Two elements of the full braid group, obtained as the increasing union of the various braid groups, with embeddings given by $\beta\to\beta\,|$,
$$B_\infty=\bigsqcup_{k=1}^\infty B_k$$
produce the same link, via closing, when one can pass from one to another via: 
\begin{enumerate}
\item Conjugation: $\beta\to\alpha\beta\alpha^{-1}$.

\item Markov move: $\beta\to g_k^{\pm1}\beta$.
\end{enumerate}
\end{theorem}

\begin{proof}
This is a version of the Reidemeister theorem, the idea being as follows:

\medskip

(1) To start with, it is clear that conjugating a braid, $\beta\to\alpha\beta\alpha^{-1}$, will produce the same link after closing, because we can pull the $\alpha,\alpha^{-1}$ to the right, in the obvious way, and there on the right, these $\alpha,\alpha^{-1}$ will annihilate, according to $\alpha\alpha^{-1}=1$.

\medskip

(2) Regarding now the Markov move from the statement, with $\beta\in B_k\subset B_{k+1}$ and with $g_1,\ldots,g_k\in B_{k+1}$ being the standard Artin generators, from Theorem 6.12 and its proof, this is the tricky move, which is worth a proof. Taking $k=3$ for an illustration, and representing $\beta\in B_3$ by a box, the link obtained by closing $g_4\beta$ is as follows, which is obviously the same link as the one obtained by closing $\beta$, and the same goes for $g_4^{-1}\beta$:
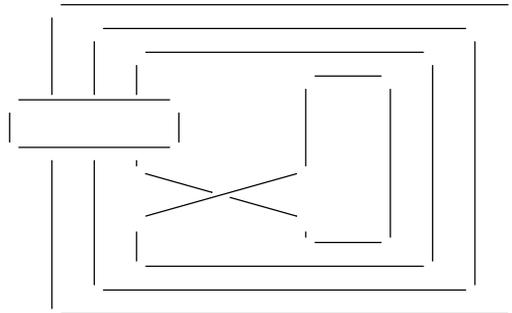
$$\xymatrix@C=9pt@R=2pt{
&\ar@{-}[dddd]\ar@{-}[rrrrrrrrrrr]&&&&&&&&&&&\ar@{-}[ddddddddddddd]\\
&&\ar@{-}[ddd]\ar@{-}[rrrrrrrrr]&&&&&&&&&\ar@{-}[ddddddddddd]\\
&&&\ar@{-}[dd]\ar@{-}[rrrrrrr]&&&&&&&\ar@{-}[ddddddddd]\\
&&&&&&&\ar@{-}[rr]\ar@{-}[dddd]&&\ar@{-}[ddddddd]\\
\ar@{-}[rrrr]\ar@{-}[dd]&&&&\ar@{-}[dd]\\
\\
\ar@{-}[rrrr]&\ar@{-}[ddddddd]&\ar@{-}[dddddd]&\ar@{-}[d]&&&\\
&&&\ar@{-}[drr]&&&&\ar@{-}[ddllll]\\
&&&&&\ar@{-}[drr]&\\
&&&\ar@{-}[dd]&&&&\ar@{-}[d]&\\
&&&&&&&\ar@{-}[rr]&&\\
&&&\ar@{-}[rrrrrrr]&&&&&&&\\
&&\ar@{-}[rrrrrrrrr]&&&&&&&&&\\
&\ar@{-}[rrrrrrrrrrr]&&&&&&&&&&&
}$$

(3) Thus, the links produced by braids are indeed invariant under the two moves in the statement. As for the proof of the converse, this comes from the Reidemeister theorem, applied in the context of the Alexander theorem, or perhaps simpler, by reasoning directly, a bit as in the proof of the Reidemeister theorem. We will leave this as an exercise.
\end{proof}

As explained before, the above kind of theorem is exactly what we need, in order to reformulate everything in terms of groups and algebra. To be more precise, looking now more in detail at what Theorem 6.15 exactly says, we are led to the following strategy:

\begin{strategy}
In order to construct numeric invariants for knots and links:
\begin{enumerate}
\item We must map $B_\infty$ somewhere, and then apply the trace.

\item And if the trace is preserved by Markov moves, it's a win.
\end{enumerate}
\end{strategy}

You get the point with all this, if we are do (1) then, by using the trace property $tr(ab)=tr(ba)$ of the trace, we will have $tr(\beta)=tr(\alpha\beta\alpha^{-1})$, in agreement with what Theorem 6.15 (1) requires. And if we do (2) too, whatever that condition exactly means, and more on this in a moment, we will have as well $tr(\beta)=tr(g_k^{\pm1}\beta)$, in agreement with what Theorem 6.15 (2) requires, so we will have our invariant for knots and links.

\bigskip

This sound very good, but before getting into details, let us be a bit megalomaniac, and add two more ambitious points to our war plan, as follows:

\begin{addendum}
Our victory will be total, with a highly reliable invariant, if:
\begin{enumerate}
\item The representation and trace are faithful as possible.

\item And they depend, if possible, on several parameters.
\end{enumerate}
\end{addendum}

Here (1) and (2) are obviously related, because the more parameters we have in (2), the more chances for our constructions in (1) to be faithful will be. In short, what we are wishing here for is an invariant which distinguishes well between various knots and links, and this can only come via a mixture of faithfulness, and parameters involved.

\bigskip

So long for the plan, and in practice now, getting back to what Strategy 6.16 says, we are faced right away with a problem, coming from the fact that $B_\infty$ is not that easy to represent. You might actually already know this, if you have struggled a bit with the exercises that I left for you, at the end of the previous section. So, we are led to:

\begin{question}
How to represent the braid group $B_\infty$?
\end{question}

So, this was the question that Reidemeister, Alexander, Markov, Artin and the others were fighting with, a long time ago, in the first half of the 20th century. Quite surprisingly, the answer to it came very late, in the 80s, from Jones \cite{jo2}, with inspiration from operator algebras, and more specifically, from his previous paper \cite{jo1} about subfactors.

\bigskip

Retrospectively looking at all this, what really matters in Jones' answer to Question 6.18 is the algebra constructed by Temperley-Lieb \cite{tli}, in the context of questions from statistical mechanics. But then, by looking even more retrospectively at all this, we can even say that the answer to Question 6.18 comes from nothing at all, meaning basic category theory. So, this will be our approach in what follows, with our answer being:

\begin{answer}
Thinking well, $B_\infty$ is self-represented, without help from the outside.
\end{answer}

So, ready for some basics, category theory? We first need objects, and our set of objects will be the good old $\mathbb N$. As for the arrows, somehow in relation with topology and braids, we will choose something very simple too, our definition being as follows:

\index{meander}
\index{partition}
\index{noncrossing partition}
\index{Temperley-Lieb category}

\begin{definition}
The Temperley-Lieb category $TL_N$ has the positive integers $\mathbb N$ as objects, with the space of arrows $k\to l$ being the formal span
$$TL_N(k,l)=span(NC_2(k,l))$$
of noncrossing pairings between an upper row of $k$ points, and a lower row of $l$ points
$$
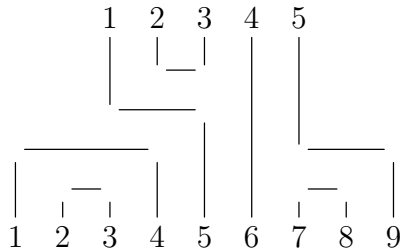
\xymatrix@R=8pt@C=5pt{
&&1\ar@{-}[dd]&2\ar@{-}[d]&3\ar@{-}[d]&4&5\ar@{-}[ddd]\\
&&&\ar@{-}[r]&&&&&\\
&&\ar@{-}[rr]&&&&&&\\
\ar@{-}[rrr]&&&&&&\ar@{-}[rr]&&\\
&\ar@{-}[r]&&&&&\ar@{-}[r]&&\\
1\ar@{-}[uu]&2\ar@{-}[u]&3\ar@{-}[u]&4\ar@{-}[uu]&5\ar@{-}[uuu]&6\ar@{-}[uuuuu]&7\ar@{-}[u]&8\ar@{-}[u]&9\ar@{-}[uu]
}$$
and with the composition of arrows appearing by composing the pairings, in the obvious way, with the rule $\bigcirc=N$, for the closed loops that might appear.
\end{definition}

This definition is something quite subtle, hiding several non-trivial things, and is worth a detailed discussion, our comments about it being as follows:

\bigskip

(1) First of all, our scalars in this chapter will be complex numbers, $\lambda\in\mathbb C$, and the ``formal span'' in the above must be understood in this sense, namely abstract complex vector space spanned by the elements of $NC_2(k,l)$. Of course it is possible to use an arbitrary field, at least at this stage of things, but remember that we are interested in quantum mechanics, and related mathematics, where the field of scalars is $\mathbb C$.

\bigskip

(2) Regarding the composition of arrows, this is by vertical concatenation, with our usual convention that things go ``from up to down''. And with this coming from care for our planet, and for entropy at the galactic level, I mean why pushing things from left to right, when we can have gravity work for us, pulling them from up to down:
$$\xymatrix@R=35pt@C=5pt{
up\ar[d]\\
down
}$$

(3) Alternatively, this comes from our previous convention for the braids, as hanging from up to down, which, when thinking well, is something quite natural too.

\bigskip

(4) Less poetically now, the ``from up to down'' convention for everything diagrams is also something useful for purely mathematical purposes, because the left-right direction can be in this way reserved for the intervention of sums $\Sigma$ and scalars $\lambda\in\mathbb C$, while the up-down direction can be reserved for ``action'', and with this being something extremely useful, trust me here, when doing complicated mathematics or physics.

\bigskip

(5) Let us discuss now what happens with the closed circles, when concatenating. As an example, let us consider a full capping of noncrossing pairings, also called meander:
$$
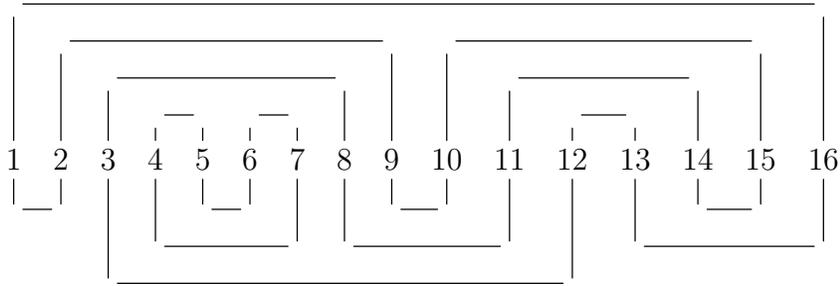
\xymatrix@R=7pt@C=5pt{
\ar@{-}[rrrrrrrrrrrrrrr]&&&&&&&&&&&&&&&\\
&\ar@{-}[rrrrrrr]&&&&&&&&\ar@{-}[rrrrr]&&&&&&\\
&&\ar@{-}[rrrrr]&&&&&&&&\ar@{-}[rrr]&&&&&\\
&&&\ar@{-}[r]&&\ar@{-}[r]&&&&&&\ar@{-}[r]&&&&\\
1\ar@{-}[uuuu]&2\ar@{-}[uuu]&3\ar@{-}[uu]&4\ar@{-}[u]&5\ar@{-}[u]&6\ar@{-}[u]&7\ar@{-}[u]&8\ar@{-}[uu]&9\ar@{-}[uuu]&10\ar@{-}[uuu]&11\ar@{-}[uu]&12\ar@{-}[u]&13\ar@{-}[u]&14\ar@{-}[uu]&15\ar@{-}[uuu]&16\ar@{-}[uuuu]\\
\ar@{-}[u]\ar@{-}[r]&\ar@{-}[u]&&&\ar@{-}[u]\ar@{-}[r]&\ar@{-}[u]&&&\ar@{-}[u]\ar@{-}[r]&\ar@{-}[u]&&&&\ar@{-}[u]\ar@{-}[r]&\ar@{-}[u]&\\
&&&\ar@{-}[uu]\ar@{-}[rrr]&&&\ar@{-}[uu]&\ar@{-}[uu]\ar@{-}[rrr]&&&\ar@{-}[uu]&&\ar@{-}[uu]\ar@{-}[rrr]&&&\ar@{-}[uu]\\
&&\ar@{-}[uuu]\ar@{-}[rrrrrrrrr]&&&&&&&&&\ar@{-}[uuu]&&&&
}$$

According to our conventions, this meander appears as the product $\pi\sigma\in NC_2(0,0)$ between the upper pairing $\sigma\in NC_2(0,16)$ and the lower pairing $\pi\in NC_2(16,0)$. But, what is the value of this product? We have two loops appearing, namely:
$$1-2-9-10-15-14-11-8-3-12-13-16$$
$$4-5-6-7$$

Thus, according to Definition 6.20, the value of this meander is $N^2$, with one $N$ for each  of the above loops, and with these two values of $N$ multiplying each other.

\bigskip

(6) The same discussion applies to an arbitrary composition $\pi\sigma\in NC_2(k,m)$ between an upper pairing $\sigma\in NC_2(k,l)$ and a lower pairing $\pi\in NC_2(l,m)$, with a certain number of loops appearing in this way, each contributing with a multiplicative factor $N$.

\bigskip

(7) Finally, in Definition 6.20 the value of the circle $N=\bigcirc$ can be pretty much anything, but due to some positivity reasons to become clear later, we will assume in what follows $N\in[1,\infty)$. Also, we will call this parameter $N$ the ``index'', with the precise reasons for calling this index to become clear later, too, as this book develops.

\bigskip

(8) As a last comment on this, we can understand now why we called the braid group $B_k$ instead of $B_N$, and this, since the beginning of this chapter, the idea being that, when doing topology and knot invariants, the number $N$ is reserved for something else.

\bigskip

With all this discussed, what is next? More category theory I guess, and matter of having a theorem formulated too, instead of definitions only, let us formulate:

\index{tensor category}

\begin{theorem}
The Temperley-Lieb category $TL_N$ is a tensor $*$-category, with:
\begin{enumerate}
\item Composition of arrows: by vertical concatenation.

\item Tensoring of arrows: by horizontal concatenation.

\item Star operation: by turning the arrows upside-down.
\end{enumerate}
\end{theorem}

\begin{proof}
This is more of a definition, disguised as a theorem. To be more precise, we already know about (1), from Definition 6.20, and we can talk as well about (2) and (3), constructed as above, with (2) using of course multiplicativity with respect to the scalars, and with (3) using antimultiplicativity with respect to the scalars:
$$\left(\sum_i\lambda_i\pi_i\right)\otimes\left(\sum_j\mu_j\sigma_j\right)=\sum_{ij}\lambda_i\mu_j\pi_i\otimes\sigma_j$$
$$\left(\sum_i\lambda_i\pi_i\right)^*=\sum_i\bar{\lambda}_i\pi_i^*$$

And the point now is that our three operations are compatible with each other via all sorts of compatibility formulae, which are all clear from definitions, with the conclusion being that what we have a tensor $*$-category, as stated. We will leave the details here, basically amounting in figuring out what a tensor $*$-category exactly is, as an exercise.
\end{proof}

In order to further understand the category $TL_N$, let us focus on its diagonal part, formed by the End spaces of various objects. With the convention that these End spaces embed into each other by adding bars at right, this is a graded algebra, as follows:
$$\Delta TL_N=\bigcup_{k\geq0}TL_N(k,k)$$

Moreover, for further fine-tuning our study, let us actually focus on the individual components of this graded algebra. These components will play a key role in what follows, and they are worth a dedicated definition, and new notation and name, as follows:

\index{Temperley-Lieb algebra}

\begin{definition}
The Temperley-Lieb algebra $TL_N(k)$ is the formal span
$$TL_N(k)=span(NC_2(k,k))$$
with multiplication coming by concatenating, with the rule $\bigcirc=N$.
\end{definition}

In other words, $TL_N(k)$ appears as the formal span of the noncrossing pairings between an upper row of $k$ points, and a lower row of $k$ points, with multiplication coming by concatenating, with our usual rule for the circles that might appear, namely:
$$\bigcirc=N$$

Getting back now to what we know about $TL_N$, from Theorem 6.21, the tensor product operation makes sense in the context of the diagonal algebra $\Delta TL_N$, but does not apply to its individual components $TL_N(k)$. However, the involution is useful, and we have:

\begin{proposition}
The Temperley-Lieb algebra $TL_N(k)$ is a $*$-algebra, with involution coming by turning the diagrams upside-down.
\end{proposition}

\begin{proof}
This is something trivial, which follows from Theorem 6.21, and can be verified as well directly, and we will leave this as an instructive exercise.
\end{proof}

Getting back now to knots and links, we first have to make the connection between braids and Temperley-Lieb diagrams. But this can be done as follows:

\begin{theorem}
The following happen:
\begin{enumerate}
\item We have a braid group representation $B_k\to TL_N(k)$, mapping standard generators to standard generators.

\item We have a trace $tr:TL_N(k)\to\mathbb C$, obtained by closing the diagrams, which is positive, and has a suitable Markov invariance property.
\end{enumerate}
\end{theorem}

\begin{proof}
Again, this is something quite intuitive, with the generators in (1) being by definition the standard ones, on both sides, and with the closing operation in (2) being similar to the one for braids, from Theorem 6.9. To be more precise:

\medskip

(1) The idea here is to map the Artin generators of the braid group to suitable modifications of the following Temperley-Lieb diagrams, called Jones projections:
$$
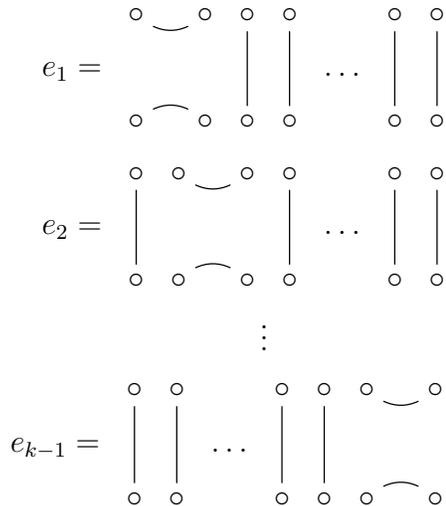
\xymatrix@C=3pt@R=7pt{
&\circ\ar@/_/@{-}[rr]&&\circ&\circ\ar@{-}[dd]&\circ\ar@{-}[dd]&&\circ\ar@{-}[dd]&\circ\ar@{-}[dd]\\
e_1=&&&&&&\ldots&\\
&\circ\ar@/^/@{-}[rr]&&\circ&\circ&\circ&&\circ&\circ}
$$
$$\xymatrix@C=3pt@R=7pt{
&\circ\ar@{-}[dd]&\circ\ar@/_/@{-}[rr]&&\circ&\circ\ar@{-}[dd]&&\circ\ar@{-}[dd]&\circ\ar@{-}[dd]\\
e_2=&&&&&&\ldots&\\
&\circ&\circ\ar@/^/@{-}[rr]&&\circ&\circ&&\circ&\circ}
$$
$$\ \ \ \ \vdots$$
$$\xymatrix@C=3pt@R=7pt{
&\circ\ar@{-}[dd]&\circ\ar@{-}[dd]&&\circ\ar@{-}[dd]&\circ\ar@{-}[dd]&\circ\ar@/_/@{-}[rr]&&\circ\\
\ \ \ \ \ \ \ e_{k-1}=&&&\ldots&&&&\\
&\circ&\circ&&\circ&\circ&\circ\ar@/^/@{-}[rr]&&\circ&&&&}
$$

As a first observation, these diagrams satisfy $e_i^2=Ne_i$, with $N=\bigcirc$ being as usual the value of the circle, so it is rather the rescaled versions $f_i=e_i/N$ which are projections, but we will not bother with this, and use our terminology above. Next, our Jones projections certainly satisfy the Artin relations $e_ie_j=e_je_i$, for $|i-j|\geq2$. Our claim now is that is that we have as well the formula $e_ie_{i\pm1}e_i=e_i$. Indeed, by translation it is enough to check $e_ie_{i+1}e_i=e_i$ at $i=1$, and this follows from the following computation:
$$\xymatrix@C=3pt@R=7pt{
&\circ\ar@/_/@{-}[rr]&&\circ&&\circ\ar@{-}[dd]&\circ\ar@{-}[dd]&&\circ\ar@{-}[dd]&\circ\ar@{-}[dd]\\
&&&&&&&\ldots&\\
&\circ\ar@{-}[dd]\ar@/^/@{-}[rr]&&\circ\ar@/_/@{-}[rr]&&\circ&\circ\ar@{-}[dd]&&\circ\ar@{-}[dd]&\circ\ar@{-}[dd]\\
e_1e_2e_1=&&&&&&&\ldots&&&=e_1\\
&\circ\ar@/_/@{-}[rr]&&\circ\ar@/^/@{-}[rr]&&\circ\ar@{-}[dd]&\circ\ar@{-}[dd]&&\circ\ar@{-}[dd]&\circ\ar@{-}[dd]\\
&&&&&&&\ldots&\\
&\circ\ar@/^/@{-}[rr]&&\circ&&\circ&\circ&&\circ&\circ}
$$

As for the verification of the relation $e_2e_1e_2=e_2$, this is similar, as follows:
$$\xymatrix@C=3pt@R=7pt{
&\circ\ar@{-}[dd]&&\circ\ar@/_/@{-}[rr]&&\circ&\circ\ar@{-}[dd]&&\circ\ar@{-}[dd]&\circ\ar@{-}[dd]\\
&&&&&&&\ldots\\
&\circ\ar@/_/@{-}[rr]&&\circ\ar@/^/@{-}[rr]&&\circ\ar@{-}[dd]&\circ\ar@{-}[dd]&&\circ\ar@{-}[dd]&\circ\ar@{-}[dd]\\
e_2e_1e_2=&&&&&&&\ldots&&&=e_2\\
&\circ\ar@/^/@{-}[rr]&&\circ\ar@/_/@{-}[rr]&&&\circ\ar@{-}[dd]&&\circ\ar@{-}[dd]&\circ\ar@{-}[dd]\\
&&&&&&&\ldots\\
&\circ\ar@{-}[uu]&&\circ\ar@/^/@{-}[rr]&&\circ&\circ&&\circ&\circ}$$

Now with the relations $e_ie_{i\pm1}e_i=e_i$ in hand, let us try to reach to the Artin relations $g_ig_{i+1}g_i=g_{i+1}g_ig_{i+1}$. For this purpose, let us set $g_i=te_i-1$. We have then:
\begin{eqnarray*}
g_ig_{i+1}g_i
&=&(te_i-1)(te_{i+1}-1)(te_i-1)\\
&=&t^3e_i-t^2(Ne_i+e_ie_{i+1}+e_{i+1}e_i)+t(2e_i+e_{i+1})-1\\
&=&t(t^2-Nt+2)e_i+te_{i+1}-t^2(e_ie_{i+1}+e_{i+1}e_i)
\end{eqnarray*}

On the other hand, we have as well the following computation:
\begin{eqnarray*}
g_{i+1}g_ig_{i+1}
&=&(te_{i+1}-1)(te_i-1)(te_{i+1}-1)\\
&=&t^3e_{i+1}-t^2(Ne_{i+1}+e_ie_{i+1}+e_{i+1}e_i)+t(2e_{i+1}+e_i)-1\\
&=&t(t^2-Nt+2)e_{i+1}+te_i-t^2(e_ie_{i+1}+e_{i+1}e_i)
\end{eqnarray*}

Thus with $t^2-Nt+1=0$ we have a representation $B_k\to TL_N(k)$, as desired.

\medskip

(2) This is something more subtle, especially in what regards the positivity properties of the trace $tr:TL_N(k)\to\mathbb C$, which requires a bit more mathematics. So, no hurry with this, and we will discuss all this, and applications, in the remainder of this chapter.
\end{proof}

\section*{6d. Knot invariants}

Let us discuss now the positivity property of the Temperley-Lieb trace, constructed as indicated in Theorem 6.24. This is something quite subtle, which in the operator algebra context, that of the original paper of Jones \cite{jo1}, comes for free, or almost, and more on this later in this book. In the meantime, we will present a more pedestrian approach to the question, based on pure combinatorics, due to Di Francesco \cite{dif}.

\bigskip

The positivity will come from a systematic study of the partitions. Let us start with:

\index{lattice of partitions}
\index{order of partitions}
\index{supremum of partitions}

\begin{definition}
Let $P(k)$ be the set of partitions of $\{1,\ldots,k\}$, and $\pi,\sigma\in P(k)$.
\begin{enumerate}
\item We write $\pi\leq\sigma$ if each block of $\pi$ is contained in a block of $\sigma$.

\item We let $\pi\vee\sigma\in P(k)$ be the partition obtained by superposing $\pi,\sigma$.
\end{enumerate}
Also, we denote by $|.|$ the number of blocks of the partitions $\pi\in P(k)$.
\end{definition}

As an illustration here, at $k=2$ we have $P(2)=\{||,\sqcap\}$, and we have:
$$||\leq\sqcap$$

Also, at $k=3$ we have $P(3)=\{|||,\sqcap|,\sqcap\hskip-3.2mm{\ }_|\,,|\sqcap,\sqcap\hskip-0.7mm\sqcap\}$, and the order relation is as follows:
$$|||\ \leq\ \sqcap|\ ,\ \sqcap\hskip-3.2mm{\ }_|\ ,\ |\sqcap\ \leq\ \sqcap\hskip-0.7mm\sqcap$$

In order to study the Gram matrix $G_k(\pi,\sigma)=N^{|\pi\vee\sigma|}$, and more specifically to compute its determinant, we will use several standard facts about the partitions. We have:

\index{Gram matrix}
\index{M\"obius function}

\begin{definition}
The M\"obius function of any lattice, and so of $P$, is given by
$$\mu(\pi,\sigma)=\begin{cases}
1&{\rm if}\ \pi=\sigma\\
-\sum_{\pi\leq\tau<\sigma}\mu(\pi,\tau)&{\rm if}\ \pi<\sigma\\
0&{\rm if}\ \pi\not\leq\sigma
\end{cases}$$
with the construction being performed by recurrence.
\end{definition}

As an illustration here, for $P(2)=\{||,\sqcap\}$, we have by definition:
$$\mu(||,||)=\mu(\sqcap,\sqcap)=1$$

Also, $||<\sqcap$, with no intermediate partition in between, so we obtain:
$$\mu(||,\sqcap)=-\mu(||,||)=-1$$

Finally, we have $\sqcap\not\leq||$, and so we have as well the following formula:
$$\mu(\sqcap,||)=0$$

Thus, as a conclusion to this, we have computed the M\"obius matrix $M_2(\pi,\sigma)=\mu(\pi,\sigma)$ of the lattice $P(2)=\{||,\sqcap\}$, the formula being as follows:
$$M_2=\begin{pmatrix}1&-1\\ 0&1\end{pmatrix}$$

Back to the general case now, the main interest in the M\"obius function comes from the M\"obius inversion formula, which states that the following happens:
$$f(\sigma)=\sum_{\pi\leq\sigma}g(\pi)\quad
\implies\quad g(\sigma)=\sum_{\pi\leq\sigma}\mu(\pi,\sigma)f(\pi)$$

In linear algebra terms, the statement and proof of this formula are as follows:

\index{M\"obius inversion}

\begin{theorem}
The inverse of the adjacency matrix of $P(k)$, given by
$$A_k(\pi,\sigma)=\begin{cases}
1&{\rm if}\ \pi\leq\sigma\\
0&{\rm if}\ \pi\not\leq\sigma
\end{cases}$$
is the M\"obius matrix of $P$, given by $M_k(\pi,\sigma)=\mu(\pi,\sigma)$.
\end{theorem}

\begin{proof}
This is well-known, coming for instance from the fact that $A_k$ is upper triangular. Indeed, when inverting, we are led into the recurrence from Definition 6.26.
\end{proof}

Now back to our Gram matrix considerations, we have the following key result:

\begin{proposition}
The Gram matrix $G_k(\pi,\sigma)=N^{|\pi\vee\sigma|}$ decomposes as a product of upper/lower triangular matrices, $G_k=A_kL_k$, where
$$L_k(\pi,\sigma)=
\begin{cases}
N(N-1)\ldots(N-|\pi|+1)&{\rm if}\ \sigma\leq\pi\\
0&{\rm otherwise}
\end{cases}$$
and where $A_k$ is the adjacency matrix of $P(k)$.
\end{proposition}

\begin{proof}
We have indeed the following computation:
\begin{eqnarray*}
G_k(\pi,\sigma)
&=&N^{|\pi\vee\sigma|}\\
&=&\#\left\{i_1,\ldots,i_k\in\{1,\ldots,N\}\Big|\ker i\geq\pi\vee\sigma\right\}\\
&=&\sum_{\tau\geq\pi\vee\sigma}\#\left\{i_1,\ldots,i_k\in\{1,\ldots,N\}\Big|\ker i=\tau\right\}\\
&=&\sum_{\tau\geq\pi\vee\sigma}N(N-1)\ldots(N-|\tau|+1)
\end{eqnarray*}

According now to the definition of $A_k,L_k$, this formula reads:
\begin{eqnarray*}
G_k(\pi,\sigma)
&=&\sum_{\tau\geq\pi}L_k(\tau,\sigma)\\
&=&\sum_\tau A_k(\pi,\tau)L_k(\tau,\sigma)\\
&=&(A_kL_k)(\pi,\sigma)
\end{eqnarray*}

Thus, we are led to the formula in the statement.
\end{proof}

We are led in this way to the following formula, due to Lindst\"om \cite{lin}:

\index{Gram matrix}
\index{Gram determinant}
\index{Lindst\"om formula}
\index{linear independence}

\begin{theorem}
The determinant of the Gram matrix $G_k$ is given by
$$\det(G_k)=\prod_{\pi\in P(k)}\frac{N!}{(N-|\pi|)!}$$
with the convention that in the case $N<k$ we obtain $0$.
\end{theorem}

\begin{proof}
If we order $P(k)$ as usual, with respect to the number of blocks, and then lexicographically, $A_k$ is upper triangular, and $L_k$ is lower triangular. Thus, we have:
\begin{eqnarray*}
\det(G_k)
&=&\det(A_k)\det(L_k)\\
&=&\det(L_k)\\
&=&\prod_\pi L_k(\pi,\pi)\\
&=&\prod_\pi N(N-1)\ldots(N-|\pi|+1)
\end{eqnarray*}

Thus, we are led to the formula in the statement.
\end{proof}

Getting now to what we wanted to do, Gram determinants, we will need:

\index{fattening of partitions}
\index{shrinking of partitions}

\begin{proposition}
We have a bijection $NC(k)\simeq NC_2(2k)$, as follows:
\begin{enumerate}
\item The application $NC(k)\to NC_2(2k)$ is the ``fattening'' one, obtained by doubling all the legs, and doubling all the strings too.

\item Its inverse $NC_2(2k)\to NC(k)$ is the ``shrinking'' application, obtained by collapsing pairs of consecutive neighbors.
\end{enumerate}
\end{proposition}

\begin{proof}
This is something self-explanatory, and in order to see how this works, let us discuss an example. Consider a noncrossing partition, say the following one:
$$
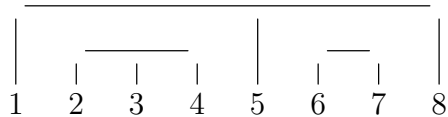
\xymatrix@R=10pt@C=10pt{
\ar@{-}[rrrrrrr]&&&&&&&\\
&\ar@{-}[rr]&&&&\ar@{-}[r]&&\\
1\ar@{-}[uu]&2\ar@{-}[u]&3\ar@{-}[u]&4\ar@{-}[u]&5\ar@{-}[uu]&6\ar@{-}[u]&7\ar@{-}[u]&8\ar@{-}[uu]
}$$

Now let us ``fatten'' this partition, by doubling everything, as follows:
$$
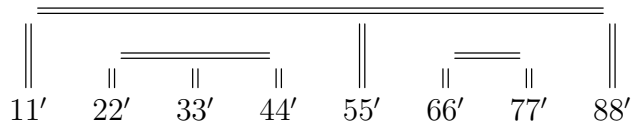
\xymatrix@R=10pt@C=10pt{
\ar@{=}[rrrrrrr]&&&&&&&\\
&\ar@{=}[rr]&&&&\ar@{=}[r]&&\\
11'\ar@{=}[uu]&22'\ar@{=}[u]&33'\ar@{=}[u]&44'\ar@{=}[u]&55'\ar@{=}[uu]&66'\ar@{=}[u]&77'\ar@{=}[u]&88'\ar@{=}[uu]
}$$

Now by relabeling the points $1,\ldots,16$, what we have is indeed a noncrossing pairing. As for the reverse operation, that is obviously obtained by ``shrinking'' our pairing, by collapsing pairs of consecutive neighbors, that is, by identifying $1=2$, then $3=4$, then $5=6$, and so on, up to $15=16$. Thus, we are led to the conclusion in the statement.
\end{proof}

At the level of the associated Gram matrices, the result is as follows:

\index{Gram matrix}

\begin{proposition}
The Gram matrices of $NC_2(2k)\simeq NC(k)$ are related by
$$G_{2k,n}(\pi,\sigma)=n^k(\Delta_{kn}^{-1}G_{k,n^2}\Delta_{kn}^{-1})(\pi',\sigma')$$
where $\pi\to\pi'$ is the shrinking operation, and $\Delta_{kn}$ is the diagonal of $G_{kn}$.
\end{proposition}

\begin{proof}
In the context of the bijection from Proposition 6.30, we have:
$$|\pi\vee\sigma|=k+2|\pi'\vee\sigma'|-|\pi'|-|\sigma'|$$

We therefore have the following formula, valid for any $n\in\mathbb N$:
$$n^{|\pi\vee\sigma|}=n^{k+2|\pi'\vee\sigma'|-|\pi'|-|\sigma'|}$$

Thus, we are led to the formula in the statement.
\end{proof}

Getting back now to our business, Gram determinants, we first have:

\begin{proposition}
The first Gram matrices and determinants for $NC_2$ are
$$\det\begin{pmatrix}N^2&N\\N&N^2\end{pmatrix}=N^2(N^2-1)$$
$$\det\begin{pmatrix}
N^3&N^2&N^2&N^2&N\\
N^2&N^3&N&N&N^2\\
N^2&N&N^3&N&N^2\\
N^2&N&N&N^3&N^2\\
N&N^2&N^2&N^2&N^3
\end{pmatrix}=N^5(N^2-1)^4(N^2-2)$$
with the matrices being written by using the lexicographic order on $NC_2(2k)$.
\end{proposition}

\begin{proof}
The formula at $k=2$, where $NC_2(4)=\{\sqcap\sqcap,\bigcap\hskip-4.9mm{\ }_\cap\,\}$, is clear. At $k=3$ however, things are tricky. We have $NC(3)=\{|||,\sqcap|,\sqcap\hskip-3.2mm{\ }_|\,,|\sqcap,\sqcap\hskip-0.7mm\sqcap\}$, and the corresponding Gram matrix and its determinant are, according to Theorem 6.29:
$$\det\begin{pmatrix}
N^3&N^2&N^2&N^2&N\\
N^2&N^2&N&N&N\\
N^2&N&N^2&N&N\\
N^2&N&N&N^2&N\\
N&N&N&N&N
\end{pmatrix}=N^5(N-1)^4(N-2)$$

By using Proposition 6.31, the Gram determinant of $NC_2(6)$ is given by:
\begin{eqnarray*}
\det(G_{6N})
&=&\frac{1}{N^2\sqrt{N}}\times N^{10}(N^2-1)^4(N^2-2)\times\frac{1}{N^2\sqrt{N}}\\
&=&N^5(N^2-1)^4(N^2-2)
\end{eqnarray*}

Thus, we have obtained the formula in the statement.
\end{proof}

In general, the above tricks won't work, because $NC(k)$ is strictly smaller than $P(k)$ at $k\geq4$. However, following Di Francesco \cite{dif}, we have the following result:

\index{meander determinant}
\index{Gram determinant}

\begin{theorem}
The determinant of the Gram matrix for $NC_2$ is given by
$$\det(G_{kN})=\prod_{r=1}^{[k/2]}P_r(N)^{d_{k/2,r}}$$
where $P_r$ are the Chebycheff polynomials, given by
$$P_0=1\quad,\quad 
P_1=X\quad,\quad 
P_{r+1}=XP_r-P_{r-1}$$
and $d_{kr}=f_{kr}-f_{k,r+1}$, with $f_{kr}$ being the following numbers, depending on $k,r\in\mathbb Z$,
$$f_{kr}=\binom{2k}{k-r}-\binom{2k}{k-r-1}$$
with the convention $f_{kr}=0$ for $k\notin\mathbb Z$. 
\end{theorem}

\begin{proof}
This is something quite technical, obtained by using a decomposition as follows of the Gram matrix $G_{kN}$, with the matrix $T_{kN}$ being lower triangular:
$$G_{kN}=T_{kN}T_{kN}^t$$

Thus, a bit as in the proof of the Lindst\"om formula, we obtain the result, but the problem lies however in the construction of $T_{kN}$, which is non-trivial. See \cite{dif}.
\end{proof}

Let us record as well the following result, also from Di Francesco \cite{dif}:

\index{meander determinant}
\index{Gram determinant}

\begin{theorem}
The determinant of the Gram matrix for $NC$ is given by
$$\det(G_{kN})=(\sqrt{N})^{a_k}\prod_{r=1}^kP_r(\sqrt{N})^{d_{kr}}$$
where $P_r$ are the Chebycheff polynomials, given by
$$P_0=1\quad,\quad 
P_1=X\quad,\quad 
P_{r+1}=XP_r-P_{r-1}$$
and $d_{kr}=f_{kr}-f_{k,r+1}$, with $f_{kr}$ being the following numbers, depending on $k,r\in\mathbb Z$,
$$f_{kr}=\binom{2k}{k-r}-\binom{2k}{k-r-1}$$
with the convention $f_{kr}=0$ for $k\notin\mathbb Z$, and where $a_k=\sum_{\pi\in P(k)}(2|\pi|-k)$.
\end{theorem}

\begin{proof}
This follows indeed from Theorem 6.33, by using Proposition 6.31.
\end{proof}

Now back to the knots and links, we have all the needed ingredients. Indeed, we can put everything together, and we obtain, following Jones \cite{jo2}:

\index{Jones polynomial}
\index{Alexander polynomial}
\index{skein relations}
\index{Temperley-Lieb trace}
\index{Kauffmann polynomial}

\begin{theorem}
We can define the Jones polynomial of an oriented knot or link as being  the image of the corresponding braid producing it via the map
$$tr:B_k\to TL_N(k)\to\mathbb C$$
with the following change of variables: 
$$N=q^{1/2}+q^{-1/2}$$
We obtain a Laurent polynomial in $q^{1/2}$, which is an invariant, up to planar isotopy.
\end{theorem}

\begin{proof}
There is a long story here, the idea being as follows:

\medskip

(1) To start with, the result follows indeed by combining the above ingredients, the idea being that the various algebraic properties of $tr:TL_N(k)\to\mathbb C$ are exactly what is needed for the above composition, up to a normalization, to be invariant under the Reidemeister moves of type I, II, III, and so to produce indeed a knot invariant. 

\medskip

(2) More specifically, the result follows from Theorem 6.15, combined with what we have in Theorem 6.24, which is now fully proved, with the positivity part coming from Theorem 6.33, and with the change of variables $N=q^{1/2}+q^{-1/2}$ in the statement coming from the equation $t^2-Nt+1=0$ that we found in the proof of Theorem 6.24.

\medskip

(3) As an illustration for how this works, consider first the unknot:
$$
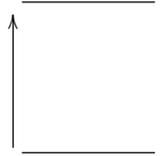
\xymatrix@R=50pt@C=50pt{
\ar@{-}[r]&\ar@{-}[d]\\
\ar[u]&\ar@{-}[l]}$$

For this knot, or rather unknot, the corresponding Jones polynomial is:
$$V=1$$

(4) Next, let us look at the link formed by two unlinked unknots:
$$
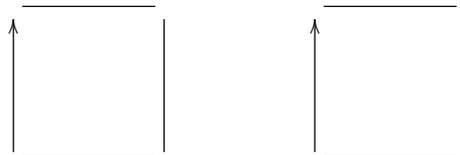
\xymatrix@R=50pt@C=50pt{
\ar@{-}[r]&\ar@{-}[d]&\ar@{-}[r]&\ar@{-}[d]\\
\ar[u]&\ar@{-}[l]&\ar[u]&\ar@{-}[l]
}$$

For this link, or rather unlink, the corresponding Jones polynomial is:
$$V=-q^{-1/2}-q^{1/2}$$

(5) Next, let us look at the link formed by two linked unknots, namely:
$$
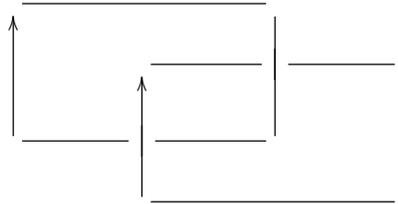
\xymatrix@R=10pt@C=40pt{
\ar@{-}[rr]&&\ar@{-}[dd]\\
&\ar@{-}[r]&|\ar@{-}[r]&\\
\ar@{-}[r]\ar[uu]&|\ar@{-}[r]&&\\
&\ar@{-}[rr]\ar[uu]&&\ar@{-}[uu]
}$$

For this link, the corresponding Jones polynomial is given by:
$$V=q^{1/2}+q^{5/2}$$

(6) Finally, let us look at the trefoil knot, which is as follows:
$$
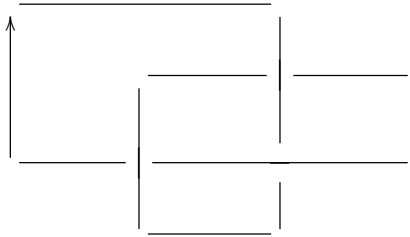
\xymatrix@R=14pt@C=40pt{
\ar@{-}[rr]&&\ar@{-}[dd]\\
&\ar@{-}[r]&|\ar@{-}[r]&\\
\ar@{-}[r]\ar[uu]&|\ar@{-}[rr]&-&\ar@{-}[u]\\
&\ar@{-}[r]\ar@{-}[uu]&\ar@{-}[u]
}$$

For this knot, the corresponding Jones polynomial is as follows:
$$V=q+q^3-q^4$$

Observe that, as previously for the unknot, this is a Laurent polynomial in $q$. This is part of a more general phenomenon, the point being that for knots, or more generally for links having an odd number of components, we get a Laurent polynomial in $q$.

\medskip

(7) In practice now, far more things can be said, about this. For instance the change of variables $N=q^{1/2}+q^{-1/2}$ is something well-known in planar algebras, and with all this being related to operator algebras and subfactor theory. More about this later.

\medskip

(8) From a purely topological perspective, however, nothing beats the skein relation interpretation of the Jones polynomial $V_L(q)$, which is as follows, with $L_+,L_-,L_0$ being knots, or rather links, differing at exactly 1 crossing, in the 3 possible ways:
$$q^{-1}V_{L_+}-qV_{L_-}=(q^{1/2}+q^{-1/2})V_{L_0}$$

To be more precise, here are the conventions for $L_+,L_-,L_0$, that you need to know, in order to play with the above formula, and compute Jones polynomials at wish:
$$
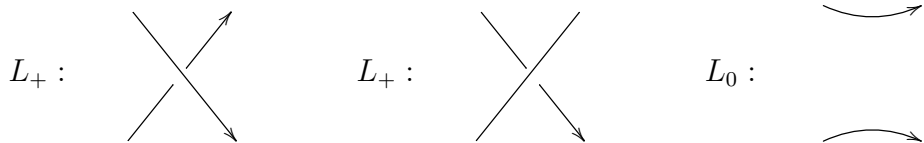
\xymatrix@R=15pt@C=15pt{
&\ar[ddrr]&&&&&\ar@{-}[dr]&&&&&\ar@/_/[rr]&&\\
L_+:&&\ar[ur]&&&L_+:&&\ar[dr]&&&L_0:\\
&\ar@{-}[ur]&&&&&\ar@{-}[uurr]&&&&&\ar@/^/[rr]&&}$$

As for the proof of the above formula, this comes from our definition of the Jones polynomial, because thinking well, ``unclosing'' links as to get braids, and then closing Temperley-Lieb diagrams as to get scalars, as required by the construction of $V_L(q)$, seemingly is some sort of identity operation, but the whole point comes from the fact that the Artin braids $g_1,\ldots,g_{k-1}$ and the Jones projections $e_1,\ldots,e_{k-1}$ differ precisely by a crossing being replaced by a non-crossing. Exercise for you, to figure out all this.

\medskip

(9) In short, up to you to learn all this, in detail, and its generalizations too, with link polynomials defined more generally via relations of the following type:
$$xP_{L_+}+yP_{L_-}+zP_{L_0}=0$$

Equivalently, we can define these more general invariants by using various versions of the Temperley-Lieb algebra. As usual, check here the papers of Jones \cite{jo1}, \cite{jo2}, \cite{jo3}.

\medskip

(10) With the comment here that, among all these polynomial invariants, the first came, historically, the Alexander polynomial, a long time ago. However, from a modern point of view, the Jones polynomial is the central invariant of knots and links.

\medskip 

(11) As another comment, with all this pure mathematics digested, physics strikes back, via a very interesting relation with statistical mechanics, happening in 2D as well, the idea being that ``interactions happen at crossings'', and it is these interactions that produce the knot invariant, as a kind of partition function. See Jones \cite{jo4}.

\medskip

(12) Quite remarkably, the above invariants can be directly understood in 3D as well, in a purely geometric way, with elegance, and no need for 2D projection. But this is a more complicated story, involving ideas from quantum field theory. See Witten \cite{wit}.
\end{proof}

\section*{6e. Exercises}

This was a quite exciting chapter, mixing groups, algebra, graphs, geometry, topology, and perhaps some physics too, who knows, and as exercises on this, we have:

\begin{exercise}
Try to enumerate the links made of linked unknots.
\end{exercise}

\begin{exercise}
Clarify all details in relation with the Reidemeister moves.
\end{exercise}

\begin{exercise}
Clarify all the details in relation with the closure of braids.
\end{exercise}

\begin{exercise}
Learn more about the braid group, and its representations.
\end{exercise}

\begin{exercise}
Learn more about the Temperley-Lieb algebra, and its applications.
\end{exercise}

\begin{exercise}
Have some fun with meanders, and the Di Francesco formula.
\end{exercise}

\begin{exercise}
Clarify all details for the construction of the Jones polynomial.
\end{exercise}

\begin{exercise}
Compute the Jones polynomial, for knots and links of your choice.
\end{exercise}

As bonus exercise, and no surprise here, learn more about knots and links, the Jones polynomial, and other invariants, first from Jones, and then from Witten.

\chapter{Low dimensions}

\section*{7a. Rotation groups}

Getting back to the finite groups, we know from chapter 5 that many of these appear as real rotation groups $G\subset O_N$, or as complex rotation groups $G\subset U_N$. In this chapter we study such finite rotation groups, $G\subset O_N$ or $G\subset U_N$, for the most in the low dimensional situation, $N=1,2,3$, where many interesting things can be said. 

\bigskip

Observe that $G\subset O_N$ implies of course $G\subset U_N$, but it is better to keep things like this, separated between real and complex, because when coming with a concrete dimension $N=1,2,3,\ldots\,$, things will be certainly easier for $O_N$ than for $U_N$. In fact, along the same lines, fine-tuning our problem, it is convenient to add right away the determinant to our discussion, and formulate the question that we want to solve as follows:

\begin{question}
What are the finite subgroups of the $4$ main rotation groups,
$$\xymatrix@R=50pt@C=50pt{
SU_N\ar[r]&U_N\\
SO_N\ar[u]\ar[r]&O_N\ar[u]}$$
in low dimensions, $N=1,2,3,\ldots$?
\end{question}

Observe that we can further fine-tune our problem by using the various secondary continuous groups that we know from chapter 1, namely the bistochastic groups $B_N\subset O_N$ and $C_N\subset U_N$, or the symplectic group $Sp_N\subset U_N$, with $N\in2\mathbb N$. And there is in fact more, because we can also consider special versions of these, obtained by intersecting with $SU_N$, and also talk about projective versions, $O_N\to PO_N$ and $U_N\to PU_N$. We will keep such things for later, typically for situations where the answer to Question 7.1 appears to be fairly complicated, and further fine-tunings of our problem are welcome.

\bigskip 

In answer now, let us start with a few observations regarding the general case, where $N\in\mathbb N$ is arbitrary. According to our general theory from chapter 5, we have:

\begin{theorem}
The reflection subgroups of the $4$ main rotation groups are
$$\xymatrix@R=50pt@C=50pt{
SK_N\ar[r]&K_N\\
SH_N\ar[u]\ar[r]&H_N\ar[u]}$$
with $K_N=\mathbb T\wr S_N$ and $H_N=\mathbb Z_2\wr S_N$ being the master reflection groups.
\end{theorem}

\begin{proof}
This is something a bit informal, with the first assertion being trivial, and with the whole point coming from the second assertion, which summarizes our knowledge of real and complex reflection groups from chapter 5, the idea being that the real ones appear as subgroups of $H_N$, and the complex ones appear as subgroups of $K_N$. For more on all this, details and comments, we refer to the material from chapter 5.
\end{proof}

Summarizing, we have many interesting finite subgroups of the 4 main rotation groups, coming from our work in chapter 5, and we will see later how these particularize at the low dimensional situation, typically $N=1,2,3$, that we will be interested in here.

\bigskip

Moving on, but still talking generalities, in the case where $N\in\mathbb N$ is arbitrary, let us record as well the following result regarding the 4 main rotation groups themselves, which is something very useful, that we will use many times, in what follows:

\begin{theorem}
The following happen, regarding the main rotation groups:
\begin{enumerate}
\item $U\in O_N\implies\det U=\pm1$.

\item $O_N=SO_N\sqcup(-SO_N)$, when $N$ is odd.

\item $U\in U_N\implies|\det U|=1$.

\item $U_N=\bigcup_{w\in\mathbb T}wSU_N$, for any $N$.
\end{enumerate}
\end{theorem}

\begin{proof}
This is something elementary, coming from definitions, as follows:

\medskip

(1) This comes indeed from the following computation:
\begin{eqnarray*}
U\in O_N
&\implies&U^t=U^{-1}\\
&\implies&\det(U^t)=\det(U^{-1})\\
&\implies&\det U=(\det U)^{-1}\\
&\implies&\det U=\pm1
\end{eqnarray*}

(2) According to (1) we have the following decomposition formula, with $\overline{SO}_N\subset O_N$ standing for the set of orthogonal matrices having determinant $-1$:
$$O_N=SO_N\sqcup\overline{SO}_N$$

Now the point is that when $N$ is odd we have $\det(-U)=-\det U$, for any matrix $U\in M_N(\mathbb R)$, and by using this, we can see right away that we have:
$$\overline{SO}_N=-SO_N$$

Thus, we are led to the decomposition formula in the statement, namely:
$$O_N=SO_N\sqcup(-SO_N)$$

By the way, observe that this fails when $N$ is even, and in a quite drastic way, for instance because at $N=2$ the group $SO_2$ consists of the rotations of the plane, while the other component $\overline{SO}_2$ consists of the symmetries of the plane. More on this later.

\medskip

(3) This follows from the following computation, similar to that in (1):
\begin{eqnarray*}
U\in U_N
&\implies&U^*=U^{-1}\\
&\implies&\det(U^*)=\det(U^{-1})\\
&\implies&\overline{\det U}=(\det U)^{-1}\\
&\implies&|\det U|=1
\end{eqnarray*}

(4) According to (3) we have the following decomposition formula, with $SU_N^{(z)}\subset U_N$ standing for the set of unitary matrices having determinant $z\in\mathbb T$, and coming with the warning that, contrary to the decomposition in (2), this is not a decomposition into connected components, due to the continuous nature of the parameter $z\in\mathbb T$:
$$U_N=\bigsqcup_{z\in\mathbb T}SU_N^{(z)}$$

Still following (2), let us try now to relate the components $SU_N^{(z)}$ to the main component, $SU_N=SU_N^{(1)}$. But this is an easy task in the present complex case, because we can extract $N$-th roots of any complex number. Indeed, let $w\in\mathbb T$ be such that:
$$w^N=z$$

Now given an arbitrary matrix $U\in SU_N^{(z)}$, the rescaled matrix $V=U/w$ is unitary, $V\in U_N$. As for the determinant of this latter matrix, this is given by:
\begin{eqnarray*}
\det(V)
&=&\det(U/w)\\
&=&\det U/w^N\\
&=&z/z\\
&=&1
\end{eqnarray*}

Thus we have $V\in SU_N$, and so $U\in wSU_N$, and with this in hand, our previous decomposition of $U_N$ takes the following form, which is the one in the statement:
$$U_N=\bigcup_{w\in\mathbb T}wSU_N$$

(5) Finally, observe that this latter decomposition is no longer a disjoint union, due to the choice needed in the above, when solving $w^N=z$. As yet another remark, getting back now to (2), all this suggests some complex number trickery, based on $i^2=-1$, in order to deal with $O_N$ when $N$ is even. We will leave some exploration here as an interesting exercise, and with the remark however that the $N=2$ case, discussed in (2), shows that we cannot really expect very concrete things to arise, in this way.
\end{proof}

As an interesting observation now, based on the technology from Theorem 7.3, we can fine-tune if we want the diagrams in Question 7.1 and Theorem 7.2, as follows: 

\begin{comment}
We can fine-tune the main diagram of rotation groups, into
$$\xymatrix@R=50pt@C=50pt{
SU_N\ar[r]&U_N^d\ar[r]&U_N\\
SO_N\ar[r]\ar[u]&O_N^d\ar[u]\ar[r]&O_N\ar[u]}$$
with the middle objects coming via $(\det U)^d=1$. The corresponding reflection groups are
$$\xymatrix@R=50pt@C=50pt{
SK_N\ar[r]&K_N^{\infty,d}\ar[r]&K_N\\
SH_N\ar[u]\ar[r]&H_N^{2,d}\ar[u]\ar[r]&H_N\ar[u]}$$
with this making a better link with the series $\{H_N^{sd}\}$, or rather $\{H_N^{s,d}\}$, from chapter 5.
\end{comment}

Which sounds interesting, but as already mentioned after Question 7.1, we will use such tricks only when needed, and stay for the moment with Question 7.1, as stated. 

\bigskip

With this discussed, time for some classification work, at small values of $N$? To start with, I am pretty much sure that we can deal with $N=1$. Indeed, in this case all our matrices are just numbers, and the main rotation groups themselves are as follows:
$$\xymatrix@R=17pt@C=18pt{
SU_1\ar[rr]&&U_1&&\{1\}\ar[rr]&&\mathbb T\\
&&&=\\
SO_1\ar[uu]\ar[rr]&&O_1\ar[uu]&&\{1\}\ar[uu]\ar[rr]&&\{\pm1\}\ar[uu]}$$

Equivalently, with $\mathbb Z_s$ standing as usual for the group of $s$-th roots of unity, and with the extra convention $\mathbb Z_\infty=\mathbb T$, that we already used in chapter 5, the diagram is:
$$\xymatrix@R=17pt@C=18pt{
\mathbb Z_1\ar[rr]&&\mathbb Z_\infty\\
\\
\mathbb Z_1\ar[uu]\ar[rr]&&\mathbb Z_2\ar[uu]}$$

Now the point is that, with the finite subgroups of the cyclic groups being cyclic, we are led to the following result, answering Question 7.1 at $N=1$:

\begin{theorem}
The finite subgroups of the basic continuous groups at $N=1$ are:
$$\xymatrix@R=17pt@C=18pt{
SU_1\ar[rr]&&U_1&&\mathbb Z_1\ar[rr]&&\{\mathbb Z_n|n\in\mathbb N\}\\
&&&:\\
SO_1\ar[uu]\ar[rr]&&O_1\ar[uu]&&\mathbb Z_1\ar[uu]\ar[rr]&&\{\mathbb Z_1,\mathbb Z_2\}\ar[uu]}$$
That is, all the finite rotation groups at $N=1$ are cyclic.
\end{theorem}

\begin{proof}
This is certainly something trivial, with only some explanations regarding the subgroups of $U_1=\mathbb T$ being needed, with the situation here being as follows:

\medskip

(1) To start with, the unit circle $\mathbb T$ has many subgroups, as you can see by picking some random numbers $\{z_i\}\subset\mathbb T$, finitely many, or countably many, or even uncountably many, and looking at the group $G=<z_i>$ that they generate, which can vary a lot. 

\medskip

(2) However, when looking at the finite subgroups $G\subset\mathbb T$, things are easy, due to:
\begin{eqnarray*}
|G|=m
&\implies&g^m=1, \forall g\in G\\
&\implies&g\in\mathbb Z_m,\ \forall g\in G\\
&\implies&G\subset\mathbb Z_m\\
&\implies&G=\mathbb Z_n,\ n|m
\end{eqnarray*}

Thus, end of the story, and we are led to the conclusion in the statement. 

\medskip

(3) Finally, let us mention that in what regards the infinite subgroups $G\subset\mathbb T$, when restricting the attention to those which are closed, we only have one solution, namely $G=\mathbb T$ itself. Thus, as a generalization of the present result, we can say that all closed rotation groups at $N=1$, finite or not, are cyclic, with our usual convention $\mathbb Z_\infty=\mathbb T$.
\end{proof}
 
At $N=2$ now, our experience with the real and complex reflections from chapter 5 shows that things are considerably more subtle over $\mathbb C$ than over $\mathbb R$. So, it is perhaps wiser to focus first on the real case, and leave the complex case for later.

\bigskip

So, let us study $SO_2$, $O_2$ are their subgroups. In what regards the groups $SO_2$, $O_2$ themselves, these are groups that we know well, and this since chapter 1, but always good to talk about them again. Their basic theory can be summarized as follows:

\begin{theorem}
We have the following results:
\begin{enumerate}
\item $SO_2$ is the group of usual rotations in the plane, which are given by:
$$R_t=\begin{pmatrix}\cos t&-\sin t\\ \sin t&\cos t\end{pmatrix}$$

\item $O_2$ consists in addition of the usual symmetries in the plane, given by:
$$S_t=\begin{pmatrix}\cos t&\sin t\\ \sin t&-\cos t\end{pmatrix}$$

\item Abstractly speaking, we have group isomorphisms as follows,
$$SO_2\simeq\mathbb T\quad,\quad O_2=\mathbb T\rtimes\mathbb Z_2$$
with the second one coming from $O_2=SO_2\rtimes<S_t>$, for any symmetry $S_t$.
\end{enumerate}
\end{theorem}

\begin{proof}
These are basically things that we know, as follows:

\medskip

(1) This is clear, because the only isometries of the plane which preserve the orientation are the usual rotations. As for the formula of $R_t$, rotation of angle $t$, this is something that we know well from chapter 1, obtained by computing $R_t\binom{1}{0}$ and $R_t\binom{0}{1}$.

\medskip

(2) This is clear too, because rotations left aside, we are left with the symmetries of the plane, in the usual sense. As for formula of $S_t$, symmetry with respect to $Ox$ rotated by $t/2$, this is something that we know too, obtained by computing $S_t\binom{1}{0}$ and $S_t\binom{0}{1}$.

\medskip

(3) The first assertion is clear, because the angles $t\in\mathbb R$, taken as usual modulo $2\pi$, form the group $\mathbb T$. As for the second assertion, the proof here is similar to the proof of the crossed product decomposition $D_n=\mathbb Z_n\rtimes\mathbb Z_2$ for the dihedral groups.
\end{proof}

Getting now to the subgroups of $SO_2,O_2$, we have the following result:

\begin{theorem}
The finite subgroups of $SO_2,O_2$ are as follows:
\begin{enumerate}
\item The finite subgroups of $SO_2$ are the cyclic groups $\mathbb Z_n$.

\item For $O_2$, we obtain in addition the dihedral groups $D_n$.
\end{enumerate}
\end{theorem}

\begin{proof}
This is again something elementary, as follows:

\medskip

(1) This is indeed something clear, geometrically, which formally comes from $SO_2\simeq\mathbb T$, via the discussion from Theorem 7.5, regarding the same group there, $U_1\simeq\mathbb T$. 

\medskip

(2) In order to prove this, consider a finite subgroup as follows:
$$G\subset O_2\quad,\quad G\not\subset SO_2$$

According to (1), we have a formula as follows, for a certain $n\in\mathbb N$:
$$G\cap SO_2=\mathbb Z_n$$

Now let us pick $S\in G-SO_2$. Since products of symmetries are rotations, any other element $T\in G-SO_2$ must satisfy $ST\in\mathbb Z_n$, and so $T\in S\mathbb Z_n$. We conclude that our group $G$ must appear as follows, coming from a subgroup $\mathbb Z_n\subset\mathbb T$, and a symmetry $S\in O_2$:
$$G=\mathbb Z_n\sqcup S\mathbb Z_n$$

But this latter group must have the same multiplication table as the dihedral group $D_n$, and  conclude that we have an isomorphism $G\simeq D_n$, as desired.
\end{proof}

Quite nice the above, and in fact we can do better, as follows:

\begin{theorem}
The finite rotation groups in $2$ dimensions appear as the symmetry groups of the regular polygons,
$$
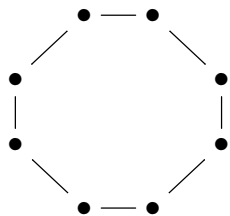
\xymatrix@R=12pt@C=13pt{
&\bullet\ar@{-}[r]\ar@{-}[dl]&\bullet\ar@{-}[dr]\\
\bullet\ar@{-}[d]&&&\bullet\ar@{-}[d]\\
\bullet\ar@{-}[dr]&&&\bullet\ar@{-}[dl]\\
&\bullet\ar@{-}[r]&\bullet}$$
with these polygons being taken unoriented as above, or oriented.
\end{theorem}

\begin{proof}
This is indeed self-explanatory, based on Theorem 7.7, and with the remark that in what regards $SO_2$, looking at the symmetries of an oriented polygon, or at the orientation-preserving symmetries of an unoriented polygon, is the same thing.
\end{proof}

The above result looks quite exciting, and it is tempting at this point to forget our next task, namely understanding what happens in 2 complex dimensions, and move instead to 3 real dimensions, with the following enigmatic question in mind:

\begin{question}
Can we have a $3D$ analogue of Theorem 7.8 going, with regular polygons replaced by regular polyhedra?
\end{question}

And good question this is. We will see in the next section, following Plato, and then Euler, and Klein and others, that the answer to this question is remarkably ``yes'', and with this solving our group theory problems, in 3 real dimensions. 

\bigskip

As for the 2 complex dimensions, these will be not forgotten either, and we will see later, following again Euler, Klein and others, including this time Rodrigues, Hamilton, and also Pauli, Dirac and other physicists, that things are quite interesting here too.

\bigskip

In short, very good news, we have discovered here a goldmine for mathematics.

\section*{7b. Platonic solids}

With Question 7.9 in mind, let us get now into 3D geometry, and symmetries. At the start of everything, we have the following remarkable result, going back to Plato:

\index{octahedron}
\index{icosahedron}
\index{dodecahedron}
\index{Platonic solid}
\index{regular polyhedron}

\begin{theorem}
There are $5$ regular polyhedra, called Platonic solids, namely:
\begin{enumerate}
\item Tetrahedron, having $4$ vertices and $4$ faces.

\item Octahedron, having $6$ vertices and $8$ faces.

\item Cube, having $8$ vertices and $6$ faces.

\item Icosahedron, having $12$ vertices and $20$ faces.

\item Dodecahedron, having $20$ vertices and $12$ faces.
\end{enumerate}
\end{theorem} 

\begin{proof}
Many things can be said here, the idea being as follows:

\medskip

(1) Let us try to figure out how a regular polyhedron looks like. There are a number of faces meeting at each vertex, $\geq3$ faces to be more precise, and when flattening the polyhedron there, we can see appear an angle $t$, called angle defect at that vertex:
$$
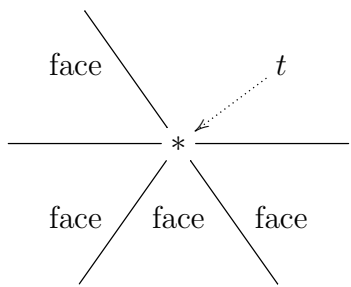
\xymatrix@R=15pt@C=12pt{
&\ar@{-}[ddr]&&\\
&{\rm face}&&t\ar@.[dl]\\
\ar@{-}[rr]&&\ast\ar@{-}[rr]&&\\
&{\rm face}&{\rm face}&{\rm face}\\
&\ar@{-}[uur]&&\ar@{-}[uul]
}$$

Now since hexagons and higher have angles $\geq120^\circ$, these cannot be used for constructing polyhedra, due to $t>0$. In fact, still due to $t>0$, we are left with 5 cases:

\medskip

-- Polyhedron made of triangles, with 3 or 4 or 5 faces meeting at each vertex.

\medskip

-- Polyhedron made of squares, with 3 faces meeting at each vertex.

\medskip

-- Polyhedron made of pengatons, with 3 faces meeting at each vertex.

\medskip

(2) Now let us try to construct the solutions. In the first case, polyhedron made of triangles, with 3 faces meeting at each vertex, we obtain the tetrahedron:
$$
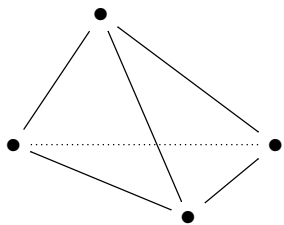
\xymatrix@R=15pt@C=20pt{
&\bullet\ar@{-}[ddl]\ar@{-}[ddrr]\ar@{-}[dddr]\\
\\
\bullet\ar@{-}[drr]\ar@{.}[rrr]&&&\bullet\ar@{-}[dl]\\
&&\bullet}$$

(3) Two other obvious solutions, corresponding to the second and fourth cases above, triangles meeting $\times4$, and squares meeting $\times3$, are the octahedron and the cube:
$$
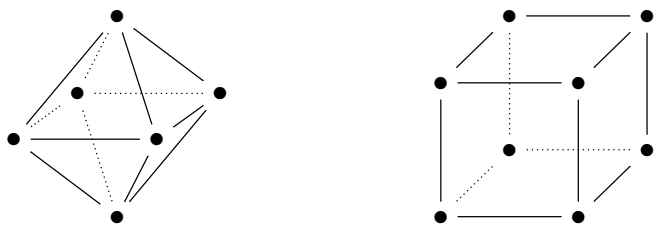
\xymatrix@R=5pt@C=2pt{
&&&\bullet\ar@{.}[ddl]\ar@{-}[dddlll]\ar@{-}[ddrrr]\ar@{-}[dddr]\\
\\
&&\bullet\ar@{.}[rrrr]&&&&\bullet\\
\bullet\ar@{-}[rrrr]\ar@{.}[urr]&&&&\bullet\ar@{-}[urr]\\
\\
&&&\bullet\ar@{-}[uur]\ar@{-}[uuurrr]\ar@{-}[uulll]\ar@{.}[uuul]}
\qquad\qquad\qquad
\xymatrix@R=13pt@C=13pt{
&\bullet\ar@{-}[rr]&&\bullet\\
\bullet\ar@{-}[rr]\ar@{-}[ur]&&\bullet\ar@{-}[ur]\\
&\bullet\ar@{.}[rr]\ar@{.}[uu]&&\bullet\ar@{-}[uu]\\
\bullet\ar@{-}[uu]\ar@{.}[ur]\ar@{-}[rr]&&\bullet\ar@{-}[uu]\ar@{-}[ur]
}$$

Before going further, observe that there is a relation between these two polyhedra, with the vertices of the octahedron appearing at the middle of the faces of the cube, and vice versa. Due to this, we say that the octahedron and cube are dual, and with this explaining why their number of vertices and faces are interchanged, as follows:
$$(6,8)\leftrightarrow(8,6)$$

By the way, observe that the tetrahedron is self-dual, $(4,4)\leftrightarrow(4,4)$. These dualities will be quite important to us later, when looking at the symmetries of our polyhedra.

\medskip

(4) Back to constructing solutions, we are left with studying the third and fifth cases in (1), namely triangles meeting $\times5$, and pentagons meeting $\times3$. And here, by some kind of miracle, we have indeed solutions, namely the icosahedron and dodecahedron, which look as follows, with in each case half of the faces, those facing us, represented:
$$
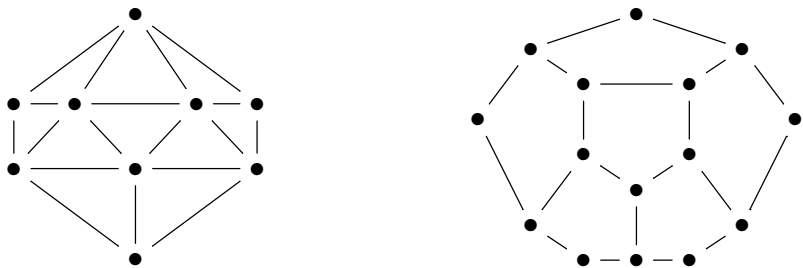
\xymatrix@R=2.6pt@C=10pt{
&&\bullet\ar@{-}[dddl]\ar@{-}[dddll]\ar@{-}[dddr]\ar@{-}[dddrr]\\
\\
\\
\bullet\ar@{-}[r]\ar@{-}[dd]&\bullet\ar@{-}[rr]\ar@{-}[ddl]\ar@{-}[ddr]&&\bullet\ar@{-}[r]\ar@{-}[ddl]\ar@{-}[ddr]&\bullet\ar@{-}[dd]\\
\\
\bullet\ar@{-}[rr]\ar@{-}[dddrr]&&\bullet\ar@{-}[rr]\ar@{-}[ddd]&&\bullet\ar@{-}[dddll]\\
\\
\\
&&\bullet}\qquad\qquad\qquad
\xymatrix@R=1pt@C=7pt{
&&&\bullet\ar@{-}[dll]\ar@{-}[drr]\\
&\bullet\ar@{-}[dr]\ar@{-}[ddl]&&&&\bullet\ar@{-}[dl]\ar@{-}[ddr]\\
&&\bullet\ar@{-}[rr]\ar@{-}[dd]&&\bullet\ar@{-}[dd]\\
\bullet\ar@{-}[dddr]&&&&&&\bullet\ar@{-}[dddl]\\
&&\bullet\ar@{-}[dr]\ar@{-}[ddl]&&\bullet\ar@{-}[dl]\ar@{-}[ddr]\\
&&&\bullet\ar@{-}[dd]\\
&\bullet\ar@{-}[dr]&&&&\bullet\ar@{-}[dl]\\
&&\bullet\ar@{-}[r]&\bullet&\bullet\ar@{-}[l]}$$

As before with the octahedron and cube, these two latter polyhedra are dual, with this intechanging their number of vertices and faces, $(12,20)\leftrightarrow(20,12)$.
\end{proof}

Many other things can be said about the Platonic solids, and their symmetry groups, and more on this in a moment. In the meantime, let us have a look at the various numbers appearing in Theorem 7.10, namely the number of vertices $v$, and the number of faces $f$. We can equally talk about the number of edges $e$, and the data is as follows:
$$\begin{array}{cccccccccccc}
&T&O&C&I&D\\
v&4&6&8&12&20\\
f&4&8&6&20&12\\
e&6&12&12&30&30
\end{array}$$

Now if you are a number theory nerd, and hope so are you, I am pretty much sure that, after contemplating a bit the above table, you will notice that:
$$v+f=e+2$$

Quite remarkably, this is something which is true for any connected planar graph, with the precise result here, due to Euler, being as follows:

\index{Euler formula}

\begin{theorem}
For a connected planar graph we have the Euler formula
$$v+f=e+2$$
with $v,e,f$ being the number of vertices, edges and faces. 
\end{theorem}

\begin{proof}
This is something very standard, the idea being as follows:

\medskip

(1) Regarding the precise statement, given a connected planar graph, that is, drawn in a planar way, without crossings, we can certainly talk about the numbers $v$ and $e$, as for any graph, and also about $f$, as being the number of faces that our graph has, in our picture, with these including by definition the outer face too, the one going to $\infty$. With these conventions, the claim is that the Euler formula $v+f=e+2$ holds indeed.

\medskip

(2) As a first illustration for how this formula works, consider a triangle:
$$
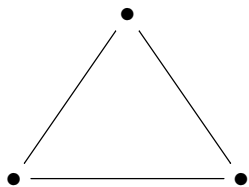
\xymatrix@R=50pt@C=30pt{
&\bullet\ar@{-}[dl]\ar@{-}[dr]\\
\bullet\ar@{-}[rr]&&\bullet}$$

Here we have $v=e=3$, and $f=2$, with this accounting for the interior and exterior, and we conclude that the Euler formula holds indeed in this case, as follows:
$$3+2=3+2$$

(3) More generally now, let us look at an arbitrary $N$-gon graph:
$$
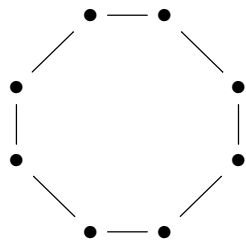
\xymatrix@R=15pt@C=15pt{
&\bullet\ar@{-}[r]\ar@{-}[dl]&\bullet\ar@{-}[dr]\\
\bullet\ar@{-}[d]&&&\bullet\ar@{-}[d]\\
\bullet\ar@{-}[dr]&&&\bullet\ar@{-}[dl]\\
&\bullet\ar@{-}[r]&\bullet}$$

Then, for this graph, the Euler formula holds indeed, as follows:
$$N+2=N+2$$

(4) With these examples discussed, let us look now for a proof. The idea will be to proceed by recurrence on the number of faces $f$. And here, as a first observation, the result holds at $f=1$, where our graph must be planar and without cycles, and so must be a tree. Indeed, with $N$ being the number of vertices, the Euler formula holds, as:
$$N+1=(N-1)+2$$

(5) At $f=2$ now, our graph must be an $N$-gon as above, but with some trees allowed to grow from the vertices, with an illustrating example here being as follows:
$$
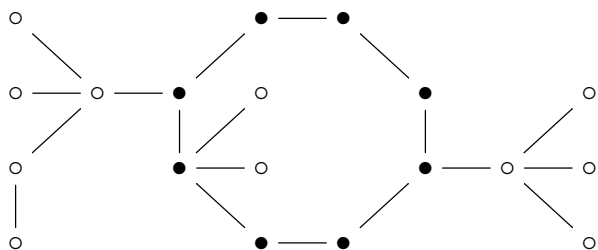
\xymatrix@R=16pt@C=18pt{
\circ\ar@{-}[dr]&&&\bullet\ar@{-}[r]\ar@{-}[dl]&\bullet\ar@{-}[dr]\\
\circ\ar@{-}[r]&\circ\ar@{-}[r]&\bullet\ar@{-}[d]&\circ&&\bullet\ar@{-}[d]&&\circ\\
\circ\ar@{-}[ur]\ar@{-}[d]&&\bullet\ar@{-}[dr]\ar@{-}[ur]\ar@{-}[r]&\circ&&\bullet\ar@{-}[dl]\ar@{-}[r]&\circ\ar@{-}[r]\ar@{-}[ur]\ar@{-}[dr]&\circ\\
\circ&&&\bullet\ar@{-}[r]&\bullet&&&\circ}$$

But here we can argue, again based on the fact that for a rooted tree, the non-root vertices are in obvious bijection with the edges, that removing all these trees won't change the problem. So, we are left with the problem for the $N$-gon, already solved in (3).

\medskip

(6) And so on, the idea being that we can first remove all the trees, by using the argument in (5), and then we are left with some sort of agglomeration of $N$-gons, for which we can check the Euler formula directly, a bit as in (3), or by recurrence.

\medskip

(7) To be more precise, let us try to do the recurrence on the number of faces $f$. For this purpose, consider one of the faces of our graph, which looks as follows, with $v_i$ denoting the number of vertices on each side, with the endpoints excluded:
$$
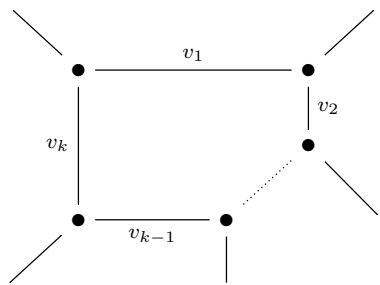
\xymatrix@R=16pt@C=18pt{
&&&&&\\
&\bullet\ar@{-}[dd]_{v_k}\ar@{-}[rrr]^{v_1}\ar@{-}[ul]&&&\bullet\ar@{-}[d]^{v_2}\ar@{-}[ur]\\
&&&&\bullet\ar@{.}[dl]\ar@{-}[dr]\\
&\bullet\ar@{-}[rr]_{v_{k-1}}\ar@{-}[dl]&&\bullet\ar@{-}[d]&&&\\
&&&&}$$

(8) Now let us collapse this chosen face to a single point, in the obvious way. In this process, the total number of vertices $v$, faces $f$ and edges $e$ evolves as follows:
$$v\to v-k+1-\sum v_i$$
$$f\to f-1$$
$$e\to e-\sum(v_i+1)$$

Thus, in this process, the Euler quantity $v+f-e$ evolves as follows:
\begin{eqnarray*}
v+f-e
&\to&v-k+1-\sum v_i+f-1-e+\sum(v_i+1)\\
&=&v-k+1-\sum v_i+f-1-e+\sum v_i+k\\
&=&v+f-e
\end{eqnarray*}

So, done with the recurrence, and the Euler formula $v+f=e+2$ is proved.
\end{proof}

As a famous application, or rather version, of the Euler formula, let us record:

\index{convex polyhedron}

\begin{proposition}
For a convex polyhedron we have the Euler formula
$$v+f=e+2$$
with $v,e,f$ being the number of vertices, edges and faces. 
\end{proposition}

\begin{proof}
This is more or less the same thing as Theorem 7.11, as follows:

\medskip

(1) In one sense, consider a convex polyhedron $P$. We can then enlarge one face, as much as needed, and then smash our polyhedron with a big hammer, as to get a planar graph $X$. As an illustration, here is how this method works, for a cube:
$$
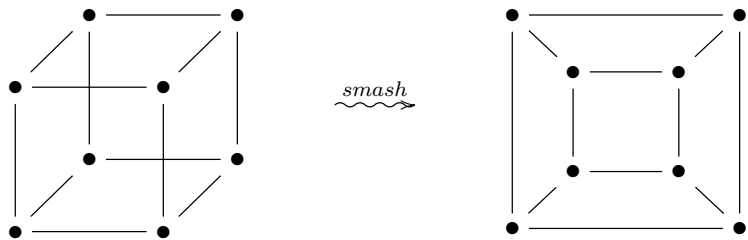
\xymatrix@R=15pt@C=15pt{
&\bullet\ar@{-}[rr]&&\bullet\\
\bullet\ar@{-}[rr]\ar@{-}[ur]&&\bullet\ar@{-}[ur]\\
&\bullet\ar@{-}[rr]\ar@{-}[uu]&&\bullet\ar@{-}[uu]\\
\bullet\ar@{-}[uu]\ar@{-}[ur]\ar@{-}[rr]&&\bullet\ar@{-}[uu]\ar@{-}[ur]
}
\qquad
\xymatrix@R=10pt@C=30pt{\\ \\ \ar@{~>}[r]^{smash}&\\ \\ }
\qquad
\xymatrix@R=9.2pt@C=10pt{
\bullet\ar@{-}[rrrr]&&&&\bullet\\
&\bullet\ar@{-}[rr]\ar@{-}[ul]&&\bullet\ar@{-}[ur]\\
\\
&\bullet\ar@{-}[uu]\ar@{-}[dl]\ar@{-}[rr]&&\bullet\ar@{-}[uu]\ar@{-}[dr]\\
\bullet\ar@{-}[rrrr]\ar@{-}[uuuu]&&&&\bullet\ar@{-}[uuuu]
}$$

But, in this process, each of the numbers $v,f,e$ stays the same, so we get the Euler formula for $P$, as a consequence of the Euler formula for $X$, from Theorem 7.11.

\medskip

(2) Conversely, consider a connected planar graph $X$. Then, save for getting rid of the internal trees, as explained in the proof of Theorem 7.11, we can assume that we are dealing with an agglomeration of $N$-gons, again as explained in the proof of Theorem 7.11. But now, we can inflate our graph as to obtain a convex polyhedron $P$:
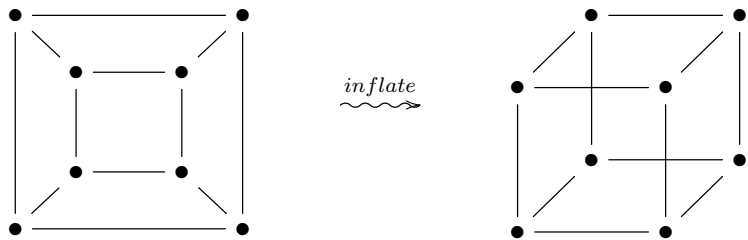
$$\xymatrix@R=9.2pt@C=10pt{
\bullet\ar@{-}[rrrr]&&&&\bullet\\
&\bullet\ar@{-}[rr]\ar@{-}[ul]&&\bullet\ar@{-}[ur]\\
\\
&\bullet\ar@{-}[uu]\ar@{-}[dl]\ar@{-}[rr]&&\bullet\ar@{-}[uu]\ar@{-}[dr]\\
\bullet\ar@{-}[rrrr]\ar@{-}[uuuu]&&&&\bullet\ar@{-}[uuuu]}
\qquad
\xymatrix@R=10pt@C=30pt{\\ \\ \ar@{~>}[r]^{inflate}&\\ \\ }
\qquad
\xymatrix@R=15pt@C=15pt{
&\bullet\ar@{-}[rr]&&\bullet\\
\bullet\ar@{-}[rr]\ar@{-}[ur]&&\bullet\ar@{-}[ur]\\
&\bullet\ar@{-}[rr]\ar@{-}[uu]&&\bullet\ar@{-}[uu]\\
\bullet\ar@{-}[uu]\ar@{-}[ur]\ar@{-}[rr]&&\bullet\ar@{-}[uu]\ar@{-}[ur]
}$$

Again, in this process, each of the numbers $v,f,e$ will stay the same, and so we get the Euler formula for $X$, as a consequence of the Euler formula for $P$.
\end{proof}

Summarizing, Euler formula understood. Getting back now to the regular polyhedra, we have the following version of Theorem 7.10, obtained by using this technology:

\begin{theorem}
The Euler formula $v+f=e+2$ allows for only $5$ regular polyhedra, which are the $5$ Platonic solids, as follows:
\begin{enumerate}
\item Tetrahedron: $4+4=6+2$.

\item Octahedron: $6+8=12+2$.

\item Cube: $8+6=12+2$.

\item Icosahedron: $12+20=30+2$.

\item Dodecahedron: $20+12=30+2$.
\end{enumerate}
\end{theorem} 

\begin{proof}
As mentioned, this is Euler's remake of the Plato theorem, using his formula $v+f=e+2$. Indeed, given a regular polyhedron, let us denote by $m$ the number of edges meeting at any vertex, and by $n$ the size of the polygons used. We have then:
$$v=\frac{2e}{m}\quad,\quad f=\frac{2e}{n}$$

Thus, for a regular polyhedron, the Euler formula $v+f=e+2$ becomes:
$$\frac{2e}{m}+\frac{2e}{n}=e+2$$

Now observe that this latter formula gives the following estimate:
$$\frac{1}{m}+\frac{1}{n}=\frac{1}{2}+\frac{1}{e}>\frac{1}{2}$$

Thus, $m,n\geq3$ must be small, and more specifically, the allowed values are:
$$(m,n)=(3,3),(4,3),(5,3),(3,4),(3,5)$$

But these are exactly the 5 cases that we met in the proof of Plato, and the continuation is the same as there, further helped by the fact that, now that we have $m,n$, we can recapture $e$, and then $v,f$ as well, with the corresponding 5 cases being as follows:
$$(v,f,e)=(4,4,6),(6,8,12),(12,20,30),(8,6,12),(20,12,30)$$

To be more precise, as before in the proof of Theorem 7.10, we obtain in this way the tetrahedron, the octahedron, the icosahedron, the cube and the dodecahedron.
\end{proof}

\section*{7c. Klein subgroups}

Getting back now to Question 7.9, we would like to compute the symmetry groups $G\subset O_3$ and $SG\subset SO_3$ of the various Platonic solids that we found, and then try to prove that these are basically all the finite subgroups of $O_3$ and $SO_3$.

\bigskip

In order to do so, let us begin with some generalities regarding $O_3,SO_3$ and their subgroups. We have here the following elementary result, further building on what we know from Theorem 7.3, regarding the groups $O_N,SO_N$ with $N$ odd:

\index{central symmetry}

\begin{proposition}
The following happen, regarding $O_3,SO_3$ and their subgroups:
\begin{enumerate}
\item The central symmetry $-1\in O_3$ is not orientation-preserving, $-1\notin SO_3$.

\item We have a disjoint union decomposition $O_3=SO_3\sqcup(-SO_3)$.

\item This decomposition gives an identification $O_3=SO_3\times\mathbb Z_2$.

\item More generally, assuming $G\subset O_3$, $-1\in G$, we have $G=SG\times\mathbb Z_2$.
\end{enumerate}
\end{proposition}

\begin{proof}
This is something elementary, as follows:

\medskip

(1) This is best viewed by using the determinant, $\det(-1)=-1$.

\medskip

(2) This follows indeed from $\det U=\pm1$ for $U\in O_3$, and from (1).

\medskip

(3) This is the group-theoretical reformulation of the decomposition in (2).

\medskip

(4) This is similar, based on $G=SG\sqcup(-SG)$, coming from $-1\in G$.
\end{proof}

Getting now to the symmetry groups that we are interested in, those of the Platonic solids found in the previous section, we have the following result, about them:

\begin{theorem}
The symmetry groups $G\subset O_3$ and the orientation-preserving symmetry groups $SG\subset SO_3$ of the Platonic solids are as follows:
\begin{enumerate}
\item Tetrahedron: $G=S_4$, $SG=A_4$.

\item Octahedron and cube: $G=S_4\times\mathbb Z_2$, $SG=S_4$.

\item Icosahedron and dodecahedron: $G=A_5\times\mathbb Z_2$, $SG=A_5$.
\end{enumerate}
\end{theorem} 

\begin{proof}
This basically comes from our experience from chapter 5, with some extra work needed for the  icosahedron and dodecahedron, the idea being as follows:

\medskip

(1) In what regards the tetrahedron, we certainly have $G=S_4$, and then $SG=A_4$, and with this latter group being usually called tetrahedral group. Observe that, contrary to what happens for the other polyhedra, the central symmetry $-1\in O_3$ is not a symmetry of the tetrahedron, so Proposition 7.14 (4) does not apply. In fact, we have $S_4\neq A_4\times\mathbb Z_2$ as abstract groups, because none of the transpositions $\tau\in S_4$ is central.

\medskip

(2) Regarding now the cube, here we have $G=H_3$, and $SG=SH_3=S_4$, as we know well since chapter 5, and with this latter $S_4$ being best understood as acting on the diagonals of the cube. Then, due to $-1\in G$, Proposition 7.14 (4) applies, and gives:
$$G=S_4\times\mathbb Z_2$$

(3) As for the ocathedron, this being dual to the cube, the symmetry groups are the same. Let us mention also that $SG=S_4$ is called octahedral group, and with this explaining why $G=H_3$, which is twice as big, is called hyperoctahedral group.

\medskip

(4) In what regards now the icosahedron and dodecahedron, these are dual too, so they have the same symmetry groups. In order to compute these common symmetry groups, let us look at the dodecahedron, whose picture, facing us, was as follows:
$$
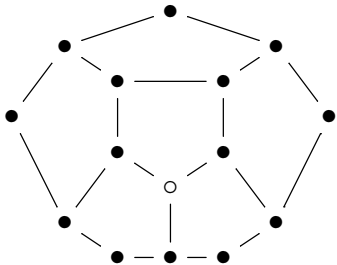
\xymatrix@R=1pt@C=7pt{
&&&\bullet\ar@{-}[dll]\ar@{-}[drr]\\
&\bullet\ar@{-}[dr]\ar@{-}[ddl]&&&&\bullet\ar@{-}[dl]\ar@{-}[ddr]\\
&&\bullet\ar@{-}[rr]\ar@{-}[dd]&&\bullet\ar@{-}[dd]\\
\bullet\ar@{-}[dddr]&&&&&&\bullet\ar@{-}[dddl]\\
&&\bullet\ar@{-}[dr]\ar@{-}[ddl]&&\bullet\ar@{-}[dl]\ar@{-}[ddr]\\
&&&\circ\ar@{-}[dd]\\
&\bullet\ar@{-}[dr]&&&&\bullet\ar@{-}[dl]\\
&&\bullet\ar@{-}[r]&\bullet&\bullet\ar@{-}[l]}$$

Now let us pick a vertex, say the one marked $\circ$, and look at the 3 faces meeting at this vertex. A symmetry $g\in SG$ must then send this vertex $\circ$ to one of the 20 available vertices $*$ of the dodecahedron, and then there is an extra $\times3$ choice, coming from the permutation of the 3 faces, at the arrival, around $*$. Thus, we conclude that we have:
$$|SG|=20\times3=60$$

(5) Before further commenting on the dodecahedron, it is worth noticing that our method above applies to any regular polyhedron $P$. Indeed, if we denote by $v$ the number of vertices, and by $m$ the number of faces meeting at any vertex, we obtain:
$$|SG|=vm=2e$$

In addition, $|SG|=2e$ can be seen as well directly, because any symmetry $g\in SG$ is uniquely determined by its action on a given edge, up to a $\times2$ choice at the arrival. Needless to say, all this fits with the data for our various polyhedra, as follows:
$$e=6,12,30\implies|SG|=12,24,60$$

(6) Getting back now to the dodecahedron, as a conclusion to the above discussion, we have two ways at looking at the corresponding group $SG$, coming from:
$$|SG|=20\times3=30\times2$$

Observe also that we have an embedding $\mathbb Z_5\subset SG$, obtained by rotating any given face of the dodecahedron. Now by putting everything together, this shows, via some routine abstract algebra that we will leave as an exercise, that we have, as claimed:
$$SG=A_5$$

(7) But you might wonder if there is a simpler proof for this, using a clever embedding $SG\subset S_5$, say a bit as before with $SH_3=S_4$ acting on the diagonals of the cube. In answer, yes, but with this being a bit neuron-burning, the idea being that we have exactly 5 cubes having vertices among the 20 vertices of the dodecahedron, and the symmetries $g\in SG$ come from permutations of these 5 cubes, which must be alternating.

\medskip

(8) Finally, still talking dodecahedron and isocahedron, these have central symmetry $-1\in G$, so by Proposition 7.14 (4) we obtain $G=A_5\times\mathbb Z_2$, as claimed.
\end{proof}

Good work that we did, and time now to answer Question 7.9, regarding the classification of finite groups of 3D rotations. In order to deal with this, we will need:

\index{rotation axis}

\begin{theorem}[Euler]
Any usual rotation in $3D$ space
$$U\in SO_3$$
has a rotation axis. 
\end{theorem}

\begin{proof}
We have the following computation, using some linear algebra magic:
\begin{eqnarray*}
\det(U-1)
&=&\det(U^t-1)\\
&=&\det(U^t(1-U))\\
&=&\det(U^t)\det(1-U)\\
&=&\det(1-U)
\end{eqnarray*}

Thus $\det(U-1)=0$, which tells us that $U$ must have a $1$-eigenvector:
$$U\xi=\xi$$

Thus, we got our rotation axis for our abstract rotation $U\in SO_3$, as desired. 
\end{proof}

We can now answer Question 7.9 positively, as follows:

\index{tetrahedral group}
\index{octahedral group}
\index{icosahedral group}
\index{Burnside formula}
\index{rotation group}
\index{Klein subgroups}

\begin{theorem}[Klein]
The finite subgroups of $SO_3$ are as follows,
\begin{enumerate}
\item Cyclic, $\mathbb Z_n$.

\item Dihedral, $D_n$.

\item Tetrahedral, $A_4$.

\item Octahedral, $S_4$.

\item Icosahedral, $A_5$.
\end{enumerate}
all appearing as symmetry groups of regular polygons and polyhedra.
\end{theorem}

\begin{proof}
This is something truly remarbable, the idea being as follows:

\medskip

(1) To start with, we certainly have as examples the groups in the statement. Indeed, those in (1,2) come from Theorem 7.7, via the following standard embedding:
$$O_2\subset SO_3\quad,\quad U\to 
\begin{pmatrix}U&0\\0&\det U\end{pmatrix}$$ 

As for those in (3,4,5), these are the groups that we found in Theorem 7.15. 

\medskip

(2) Regarding now the converse, assume that $G\subset SO_3$ is finite. Given $g\in G-\{1\}$, consider its rotation axis coming from Theorem 7.16, and then the two points $\pm x$ where this axis intersects the unit sphere $S^2\subset\mathbb R^3$, called poles of $g$. We can consider then the set $X\subset S^2$ of all poles of all elements $g\in G-\{1\}$, and we have an action as follows:
$$G\curvearrowright X$$

(3) In order to exploit this latter action, we can use the following counting trick, due to Burnside, which is valid for any finite group action on a finite set, $G\curvearrowright X$:
\begin{eqnarray*}
\sum_{g\in G}|X^g|
&=&\sum_{x\in X}|G_x|\\
&=&|G|\sum_{x\in X}\frac{1}{|Gx|}\\
&=&|G|\sum_{O\in X/G}|O|\cdot\frac{1}{|O|}\\
&=&|G|\cdot|X/G|
\end{eqnarray*}

To be more precise, here $X^g\subset X$ is the set of fixed points by $g\in G$, and $G_x\subset G$ is the stabilizer of $x\in X$, and we have used the general theory from chapter 4.

\medskip

(4) Now let us see what the Burnside formula gives, for the action in (2). If we denote by $N$ the number of orbits of our action $G\curvearrowright X$, this formula reads:
$$|X|+2(|G|-1)=N|G|$$

Now observe that this latter formula can be further processed in the following way, with $\{x_1,\ldots,x_N\}\subset X$ being a set of representatives for the orbits of $G\curvearrowright X$:
\begin{eqnarray*}
2\left(1-\frac{1}{|G|}\right)
&=&N-\frac{|X|}{|G|}\\
&=&N-\frac{1}{|G|}\sum_{i=1}^N[G:G_{x_i}]\\
&=&\sum_{i=1}^N1-\frac{1}{|G_{x_i}|}
\end{eqnarray*}

(5) And the point is that this latter formula is exactly what we need. Indeed, observe that the left term and the right components are subject to the following estimates:
$$2\left(1-\frac{1}{|G|}\right)<2\quad,\quad1-\frac{1}{|G_{x_i}|}\geq\frac{1}{2}$$

We conclude that we must have $N=2,3$, which is a big win, we are almost there.

\medskip

(6) In practice now, in the case $N=2$, the formula that we found in (4) reads:
$$\frac{2}{|G|}=\frac{1}{|G_x|}+\frac{1}{|G_y|}$$

But a quick study shows that the solution here is $G=\mathbb Z_n$, corresponding to:
$$\frac{2}{n}=\frac{1}{n}+\frac{1}{n}$$

(7) Regarding now the case $N=3$, here the formula found in (4) reads:
$$1+\frac{2}{|G|}=\frac{1}{|G_x|}+\frac{1}{|G_y|}+\frac{1}{|G_z|}$$

But here we have 4 possible cases, corresponding to the following solutions of this:
$$1+\frac{2}{2n}=\frac{1}{2}+\frac{1}{2}+\frac{1}{n}\qquad,\qquad
1+\frac{2}{12}=\frac{1}{2}+\frac{1}{3}+\frac{1}{3}$$
$$1+\frac{2}{24}=\frac{1}{2}+\frac{1}{3}+\frac{1}{4}\qquad,\qquad 
1+\frac{2}{60}=\frac{1}{2}+\frac{1}{3}+\frac{1}{5}$$

And a study of these cases, that we will leave as an instructive exercise, leads to the other solutions in the statement, namely $G=D_n$, $G=A_4$, $G=S_4$, $G=A_5$.
\end{proof}

Very nice all this. We should mention that, with a bit more work, based on the above, the finite subgroups of $O_3$ can be classified too, using Proposition 7.14, and with this being something quite straightforward. We will leave this, again, as an instructive exercise.

\section*{7d. Complex rotations}

Moving forward, let us go back now to $N=2$ dimensions, but with a study in the complex case. We first have here the following result, that we know since chapter 1:

\index{unit sphere}

\begin{proposition}
We have the following formula,
$$SU_2=\left\{\begin{pmatrix}a&b\\ -\bar{b}&\bar{a}\end{pmatrix}\ \Big|\ |a|^2+|b|^2=1\right\}$$
which makes $SU_2$ isomorphic to the unit complex sphere $S^1_\mathbb C\subset\mathbb C^2$.
\end{proposition}

\begin{proof}
We know this from chapter 1, the idea being that for a matrix $U=\binom{a\ b}{c\ d}$ having determinant 1, the complex isometry condition $U^*=U^{-1}$ reads:
$$\begin{pmatrix}\bar{a}&\bar{c}\\\bar{b}&\bar{d}\end{pmatrix}
=\begin{pmatrix}d&-b\\-c&a\end{pmatrix}$$

Thus $U$ is as in the statement, and with $|a|^2+|b|^2=1$ coming from $\det U=1$.
\end{proof}

Here is a useful reformulation of the above result, using real numbers:

\begin{proposition}
We have the formula
$$SU_2=\left\{\begin{pmatrix}x+iy&z+it\\ -z+it&x-iy\end{pmatrix}\ \Big|\ x^2+y^2+z^2+t^2=1\right\}$$
which makes $SU_2$ isomorphic to the unit real sphere $S^3_\mathbb R\subset\mathbb R^3$.
\end{proposition}

\begin{proof}
This is indeed something self-explanatory, coming from Proposition 7.18.
\end{proof}

At a more advanced level now, here is yet another reformulation of what we have:

\index{Pauli matrices}
\index{quaternions}
\index{spin matrices}

\begin{theorem}
We have the following formula,
$$SU_2=\left\{xc_1+yc_2+zc_3+tc_4\ \Big|\ x^2+y^2+z^2+t^2=1\right\}$$
where $c_1,c_2,c_3,c_4$ are matrices given by
$$c_1=\begin{pmatrix}1&0\\ 0&1\end{pmatrix}\quad,\quad
c_2=\begin{pmatrix}i&0\\ 0&-i\end{pmatrix}\quad,\quad 
c_3=\begin{pmatrix}0&1\\ -1&0\end{pmatrix}\quad,\quad 
c_4=\begin{pmatrix}0&i\\ i&0\end{pmatrix}$$
called Pauli spin matrices.
\end{theorem}

\begin{proof}
According to Proposition 7.19 the elements $U\in SU_2$ are the matrices as follows, depending on parameters $x,y,z,t\in\mathbb R$ satisfying $x^2+y^2+z^2+t^2=1$:
$$U=x\begin{pmatrix}1&0\\ 0&1\end{pmatrix}
+y\begin{pmatrix}i&0\\ 0&-i\end{pmatrix}
+z\begin{pmatrix}0&1\\ -1&0\end{pmatrix}+
t\begin{pmatrix}0&i\\ i&0\end{pmatrix}$$

Thus, we are led to the conclusions in the statement.
\end{proof}

The above result is often the most convenient one, when dealing with $SU_2$. This is because the Pauli matrices have a number of remarkable properties, which are very useful when doing computations. These properties can be summarized as follows:

\begin{proposition}
The Pauli matrices multiply according to the formulae
$$c_2^2=c_3^2=c_4^2=-1$$
$$c_2c_3=-c_3c_2=c_4$$
$$c_3c_4=-c_4c_3=c_2$$
$$c_4c_2=-c_2c_4=c_3$$
they conjugate according to the following rules,
$$c_1^*=c_1\ ,\ c_2^*=-c_2\ ,\ c_3^*=-c_3\ ,\ c_4^*=-c_4$$
and they form an orthonormal basis of $M_2(\mathbb C)$, with respect to the scalar product
$$<a,b>=tr(ab^*)$$
with $tr:M_2(\mathbb C)\to\mathbb C$ being the normalized trace of $2\times 2$ matrices, $tr=Tr/2$.
\end{proposition}

\begin{proof}
The first two assertions, regarding the multiplication and conjugation rules for the Pauli matrices, follow from some elementary computations. As for the last assertion, this follows by using these rules. Indeed, the fact that the Pauli matrices are pairwise orthogonal follows from computations of the following type, for $i\neq j$:
$$<c_i,c_j>
=tr(c_ic_j^*)
=tr(\pm c_ic_j)
=tr(\pm c_k)
=0$$

As for the fact that the Pauli matrices have norm 1, this follows from:
$$<c_i,c_i>
=tr(c_ic_i^*)
=tr(\pm c_i^2)
=tr(c_1)
=1$$

Thus, we are led to the conclusion in the statement.
\end{proof}

Moving on, we would like to discuss now a key relation between $SU_2$ and $SO_3$. Let us start with the following construction, whose goal will become clear in a moment:

\index{adjoint action}

\begin{proposition}
The adjoint action $SU_2\curvearrowright M_2(\mathbb C)$, given by
$$T_U(M)=UMU^*$$
leaves invariant the following real vector subspace of $M_2(\mathbb C)$,
$$E=span_\mathbb R(c_1,c_2,c_3,c_4)$$
and we obtain in this way a group morphism $SU_2\to GL_4(\mathbb R)$.
\end{proposition}

\begin{proof}
We have two assertions to be proved, as follows:

\medskip

(1) We must first prove that, with $E\subset M_2(\mathbb C)$ being the real vector space in the statement, we have the following implication:
$$U\in SU_2,M\in E\implies UMU^*\in E$$

But this is clear from the multiplication rules for the Pauli matrices, from Proposition 7.21. Indeed, let us write our matrices $U,M$ as follows:
$$U=xc_1+yc_2+zc_3+tc_4$$
$$M=ac_1+bc_2+cc_3+dc_4$$

We know that the coefficients $x,y,z,t$ and $a,b,c,d$ are real, due to $U\in SU_2$ and $M\in E$. The point now is that when computing $UMU^*$, by using the various rules from Proposition 7.21, we obtain a matrix of the same type, namely a combination of $c_1,c_2,c_3,c_4$, with real coefficients. Thus, we have $UMU^*\in E$, as desired.

\medskip

(2) In order to conclude, let us identify $E\simeq\mathbb R^4$, by using the basis $c_1,c_2,c_3,c_4$. The result found in (1) shows that we have a correspondence as follows:
$$SU_2\to M_4(\mathbb R)\quad,\quad U\to (T_U)_{|E}$$

Now observe that for any $U\in SU_2$ and any $M\in M_2(\mathbb C)$ we have:
$$T_{U^*}T_U(M)=U^*UMU^*U=M$$

Thus $T_{U^*}=T_U^{-1}$, and so the correspondence that we found can be written as:
$$SU_2\to GL_4(\mathbb R)\quad,\quad U\to (T_U)_{|E}$$

But this a group morphism, due to the following computation:
$$T_UT_V(M)=UVMV^*U^*=T_{UV}(M)$$

Thus, we are led to the conclusion in the statement.
\end{proof}

The point now, which makes the link with $SO_3$, and which will ultimately elucidate the structure of $SO_3$, is that Proposition 7.22 can be improved as follows:

\begin{theorem}
The adjoint action $SU_2\curvearrowright M_2(\mathbb C)$ leaves invariant the space
$$F=span_\mathbb R(c_2,c_3,c_4)$$
and we obtain in this way a group morphism $SU_2\to SO_3$.
\end{theorem}

\begin{proof}
We can do this in several steps, as follows:

\medskip

(1) Our first claim is that the group morphism $SU_2\to GL_4(\mathbb R)$ constructed in Proposition 7.22 is in fact a morphism $SU_2\to O_4$. In order to prove this, recall the following formula, valid for any $U\in SU_2$, from the proof of Proposition 7.22:
$$T_{U^*}=T_U^{-1}$$

We want to prove that the matrices $T_U\in GL_4(\mathbb R)$ are orthogonal, and in view of the above formula, it is enough to prove that we have:
$$T_U^*=(T_U)^t$$

So, let us prove this. For any two matrices $M,N\in E$, we have:
\begin{eqnarray*}
<T_{U^*}(M),N>
&=&<U^*MU,N>\\
&=&tr(U^*MUN)\\
&=&tr(MUNU^*)
\end{eqnarray*}

On the other hand, we have as well the following formula:
\begin{eqnarray*}
<(T_U)^t(M),N>
&=&<M,T_U(N)>\\
&=&<M,UNU^*>\\
&=&tr(MUNU^*)
\end{eqnarray*}

Thus we have indeed $T_U^*=(T_U)^t$, which proves our $SU_2\to O_4$ claim.

\medskip

(2) In order now to finish, recall that we have by definition $c_1=1$, as a matrix. Thus, the action of $SU_2$ on the vector $c_1\in E$ is given by:
$$T_U(c_1)=Uc_1U^*=UU^*=1=c_1$$ 

We conclude that $c_1\in E$ is invariant under $SU_2$, and by orthogonality the following subspace of $E$ must be invariant as well under the action of $SU_2$:
$$e_1^\perp=span_\mathbb R(c_2,c_3,c_4)$$

Now if we call this subspace $F$, and we identify $F\simeq\mathbb R^3$ by using the basis $c_2,c_3,c_4$, we obtain by restriction to $F$ a morphism of groups as follows:
$$SU_2\to O_3$$

But since this morphism is continuous and $SU_2$ is connected, its image must be connected too. Now since the target group decomposes as $O_3=SO_3\sqcup(-SO_3)$, and $1\in SU_2$ gets mapped to $1\in SO_3$, the whole image must lie inside $SO_3$, and we are done.
\end{proof}

We can now formulate a key result, due to Euler-Rodrigues, as follows:

\index{double cover map}
\index{Euler-Rodrigues formula}
\index{rotation}

\begin{theorem}
We have a double cover map, obtained via the adjoint representation,
$$SU_2\to SO_3$$
and this map produces the Euler-Rodrigues formula
$$U=\begin{pmatrix}
x^2+y^2-z^2-t^2&2(yz-xt)&2(xz+yt)\\
2(xt+yz)&x^2+z^2-y^2-t^2&2(zt-xy)\\
2(yt-xz)&2(xy+zt)&x^2+t^2-y^2-z^2
\end{pmatrix}$$
for the generic elements of $SO_3$.
\end{theorem}

\begin{proof}
We have several things to be proved here, the idea being as follows:

\medskip

(1) Our first claim is that, with respect to the standard basis $c_1,c_2,c_3,c_4$ of the vector space $\mathbb R^4=span(c_1,c_2,c_3,c_4)$, the morphism $T:SU_2\to GL_4(\mathbb R)$ is given by:
$$T_U=\begin{pmatrix}
1&0&0&0\\
0&x^2+y^2-z^2-t^2&2(yz-xt)&2(xz+yt)\\
0&2(xt+yz)&x^2+z^2-y^2-t^2&2(zt-xy)\\
0&2(yt-xz)&2(xy+zt)&x^2+t^2-y^2-z^2
\end{pmatrix}$$

(2) Indeed, with notations from Proposition 7.22 and its proof, let us first look at the action $L:SU_2\curvearrowright\mathbb R^4$ by left multiplication, which is by definition given by:
$$L_U(M)=UM$$

In order to compute the matrix of this action, let us write, as usual:
$$U=xc_1+yc_2+zc_3+tc_4$$
$$M=ac_1+bc_2+cc_3+dc_4$$

By using the multiplication formulae in Proposition 7.21, we obtain:
\begin{eqnarray*}
UM
&=&(xc_1+yc_2+zc_3+tc_4)(ac_1+bc_2+cc_3+dc_4)\\
&=&(xa-yb-zc-td)c_1\\
&+&(xb+ya+zd-tc)c_2\\
&+&(xc-yd+za+tb)c_3\\
&+&(xd+yc-zb+ta)c_4
\end{eqnarray*}

We conclude that the matrix of the left action considered above is:
$$L_U=\begin{pmatrix}
x&-y&-z&-t\\
y&x&-t&z\\
z&t&x&-y\\
t&-z&y&x
\end{pmatrix}$$

(3) Similarly, let us look now at the action $R:SU_2\curvearrowright\mathbb R^4$ by right multiplication, which is by definition given by the following formula:
$$R_U(M)=MU^*$$

In order to compute the matrix of this action, let us write, as before:
$$U=xc_1+yc_2+zc_3+tc_4$$
$$M=ac_1+bc_2+cc_3+dc_4$$

By using the multiplication formulae in Proposition 7.21, we obtain:
\begin{eqnarray*}
MU^*
&=&(ac_1+bc_2+cc_3+dc_4)(xc_1-yc_2-zc_3-tc_4)\\
&=&(ax+by+cz+dt)c_1\\
&+&(-ay+bx-ct+dz)c_2\\
&+&(-az+bt+cx-dy)c_3\\
&+&(-at-bz+cy+dx)c_4
\end{eqnarray*}

We conclude that the matrix of the right action considered above is:
$$R_U=\begin{pmatrix}
x&y&z&t\\
-y&x&-t&z\\
-z&t&x&-y\\
-t&-z&y&x
\end{pmatrix}$$

(4) Now by composing, the matrix of the adjoint matrix in the statement is:
\begin{eqnarray*}
T_U
&=&R_UL_U\\
&=&\begin{pmatrix}
x&y&z&t\\
-y&x&-t&z\\
-z&t&x&-y\\
-t&-z&y&x
\end{pmatrix}
\begin{pmatrix}
x&-y&-z&-t\\
y&x&-t&z\\
z&t&x&-y\\
t&-z&y&x
\end{pmatrix}\\
&=&\begin{pmatrix}
1&0&0&0\\
0&x^2+y^2-z^2-t^2&2(yz-xt)&2(xz+yt)\\
0&2(xt+yz)&x^2+z^2-y^2-t^2&2(zt-xy)\\
0&2(yt-xz)&2(xy+zt)&x^2+t^2-y^2-z^2
\end{pmatrix}
\end{eqnarray*}

(5) Summarizing, we have proved our claim in (1). We conclude that, when looking at $T:SU_2\to GL_4(\mathbb R)$ as a group morphism $SU_2\to O_4$, what we have in fact is a group morphism $SU_2\to O_3$, and even $SU_2\to SO_3$, given by the Euler-Rodrigues formula.

\medskip

(6) Next, the kernel of this morphism is elementary to compute, as follows:
\begin{eqnarray*}
\ker(SU_2\to SO_3)
&=&\left\{U\in SU_2\Big|T_U(M)=M,\forall M\in E\right\}\\
&=&\left\{U\in SU_2\Big|Uc_i=c_iU,\forall i\right\}\\
&=&\{\pm1\}
\end{eqnarray*}

(7) Finally, in what regards the surjectivity, we can argue here for instance that since each rotation $U\in SO_3$ is uniquely determined by its rotation axis, plus its rotation angle $t\in[0,2\pi)$, we are led to the conclusion that $U$ is uniquely determined by an element of $SU_2/\{\pm 1\}$, and so appears indeed via the Euler-Rodrigues formula, as stated.
\end{proof}

Getting back now to our finite subgroup questions, we have:

\index{Klein subgroups}
\index{dicyclic group}
\index{binary tetrahedral group}
\index{binary octahedral group}
\index{binary icosahedral group}

\begin{theorem}[Klein]
The subgroups of $SU_2$ are as follows:
\begin{enumerate}
\item Cyclic, $\mathbb Z_n$.

\item Dicyclic, $DC_n$.

\item Binary tetrahedral, lifting $A_4$.

\item Binary octahedral, lifting $S_4$.

\item Binary icosahedral, lifting $A_5$.
\end{enumerate}
\end{theorem}

\begin{proof}
This is indeed something quite standard, from what we have, the idea being that the various groups in (2-5) appear as lifts via $SU_2\to SO_3$ of the groups in Theorem 7.17 (2-5). We will leave some further learning here as an instructive exercise.
\end{proof}

Good work that we did, but the story is not over with this, because we can talk about $SU_3$ as well. As usual, exercise for you, to learn more about all this.

\section*{7e. Exercises}

This was an exciting chapter, of physics flavor, and as exercises here, we have:

\begin{exercise}
Learn more, from Pauli and physicists, about the Pauli matrices.
\end{exercise}

\begin{exercise}
Learn also, from Dirac and the others, about the Dirac matrices.
\end{exercise}

\begin{exercise}
Learn too about Weyl matrices, Clifford algebras, and related topics.
\end{exercise}

\begin{exercise}
And do not forget to have a look at quaternions and octonions too.
\end{exercise}

\begin{exercise}
Nor should you forget to learn about the Gell-Mann matrices.
\end{exercise}

\begin{exercise}
Learn about Clebsch-Gordan rules, and ADE classification.
\end{exercise}

\begin{exercise}
Learn more about the Klein subgroups, of both $SU_2$ and $SO_3$.
\end{exercise}

\begin{exercise}
Have a look as well at the subgroups of $U_2$, $O_3$, $SU_3$, $U_3$.
\end{exercise}

As bonus exercise, read about the various ADE classifications in mathematics.

\chapter{Galois groups}

\section*{8a. Equations, roots}

Time to end the present Part II of this book, on basic applications of groups to geometry, and to mathematics in general, and no discussion here would be complete without a word on Galois theory. So, we will talk about this in the present chapter.

\bigskip

Before anything, we should mention that Galois theory and its applications are something quite vast, requiring a lot of algebraic knowledge, including a bit of advanced theory of polynomials, field theory of all types, also rings, modules and ideals, not to forget the algebras, and with all this going well beyond the group theory that we are familiar with. We will certainly not have time to discuss all this in detail, and our aim here will be that of explaining what Galois theory is about, rather than what Galois theory is.

\bigskip

This being said, since you just asked me, what is Galois theory? In answer:

\begin{principle}
Galois theory is the art of studying the polynomials $P$ by permuting their roots, and with the corresponding permutation groups $Gal(P)$ being called Galois groups. And with this having countless applications, across mathematics.
\end{principle}

But probably too much talking, let us get to work. We must first talk about polynomials, and as a starting point here we have the following result, that you surely know, and which by the way is called fundamental theorem of algebra, no less than that:

\index{roots of polynomial}
\index{complex roots}

\begin{theorem}
Any polynomial $P\in\mathbb C[X]$ decomposes as
$$P=c(X-a_1)\ldots (X-a_N)$$
with $c\in\mathbb C$ and with $a_1,\ldots,a_N\in\mathbb C$.
\end{theorem}

\begin{proof}
This is something very classical, the idea being as follows:

\medskip

(1) The problem is that of proving that our polynomial has at least one root, because afterwards we can proceed by recurrence. We prove this by contradiction. So, assume that $P$ has no roots, and pick a number $z\in\mathbb C$ where $|P|$ attains its minimum:
$$|P(z)|=\min_{x\in\mathbb C}|P(x)|>0$$ 

(2) Since $Q(t)=P(z+t)-P(z)$ is a polynomial which vanishes at $t=0$, this polynomial must be of the form $ct^k$ + higher terms, with $c\neq0$, and with $k\geq1$ being an integer. We obtain from this that, with $t\in\mathbb C$ small, we have the following estimate:
$$P(z+t)\simeq P(z)+ct^k$$

(3) If we write $t=rw$, with $r>0$ small, and with $|w|=1$, our estimate becomes:
$$P(z+rw)\simeq P(z)+cr^kw^k$$

(4) Now recall that we assumed $P(z)\neq0$. We can therefore choose $w\in\mathbb T$ such that $cw^k$ points in the opposite direction to that of $P(z)$, and we obtain in this way:
\begin{eqnarray*}
|P(z+rw)|
&\simeq&|P(z)+cr^kw^k|\\
&=&|P(z)|(1-|c|r^k)
\end{eqnarray*}

(5) Now by choosing $r>0$ small enough, as for the error in the first estimate to be small, and overcame by the negative quantity $-|c|r^k$, we obtain from this:
$$|P(z+rw)|<|P(z)|$$

(6) But this contradicts our definition of $z\in\mathbb C$, as a point where $|P|$ attains its minimum. Thus $P$ has a root, and by recurrence it has $N$ roots, as stated.
\end{proof}

As a comment now, the above proof being by contradiction, it cannot be used in order to find the roots of a given polynomial $P\in\mathbb C[X]$. And in relation with this latter question, that we would like to investigate now, many things can be said.

\bigskip

Let us start with a discussion in degree 2. The main result here is well-known:

\begin{theorem}
The solutions of $ax^2+bx+c=0$ with $a,b,c\in\mathbb C$ are
$$x_{1,2}=\frac{-b\pm\sqrt{b^2-4ac}}{2a}$$
with the square root of complex numbers being defined as $\sqrt{re^{it}}=\sqrt{r}e^{it/2}$.
\end{theorem}

\begin{proof}
This is indeed something that you surely know well, coming from:
\begin{eqnarray*}
ax^2+bx+c=0
&\iff&x^2+\frac{b}{a}x+\frac{c}{a}=0\\
&\iff&\left(x+\frac{b}{2a}\right)^2=\frac{b^2-4ac}{4a^2}\\
&\iff&x+\frac{b}{2a}=\pm\frac{\sqrt{b^2-4ac}}{2a}
\end{eqnarray*}

Thus, we are led to the conclusion in the statement.
\end{proof}

In practice now, Theorem 8.3 is not exactly the end of the story with degree 2, because that formula relies on some trigonometry for extracting the square roots. In practice, it is often better to use the following result, involving real square roots only:

\index{square root}

\begin{theorem}
Any complex number $x=p+iq$ has two square roots, given by
$$\sqrt{x}=\pm\sqrt{\frac{p+\sqrt{p^2+q^2}}{2}}\pm i\sqrt{\frac{-p+\sqrt{p^2+q^2}}{2}}$$
with the signs being identical when $q>0$, and opposite when $q<0$.
\end{theorem}

\begin{proof}
This is something quite routine, the idea being as follows:

\medskip

(1) With $x=p+iq$ as in the statement, and $\sqrt{x}=r+is$, our equation is:
$$(r+is)^2=p+iq$$

In terms of the real and imaginary parts, we have two equations, as follows:
$$r^2-s^2=p\quad,\quad 2rs=q$$

(2) Let us first compute the number $u=r^2$. The equation for it is as follows:
$$u-\frac{p^2}{4u}=q$$

Thus, the number $u=r^2$ satisfies the following degree 2 equation:
$$u^2-pu-\frac{q^2}{4}=0$$

But this latter equation has a unique positive solution, given by:
$$u=\frac{p+\sqrt{p^2+q^2}}{2}$$

Thus, we are led to the formula of $r=\pm\sqrt{u}$ in the statement.

\medskip

(3) Similarly, let us compute now $v=s^2$. The equation for it is as follows:
$$\frac{q^2}{4v}-v=p$$

Thus, the number $v=s^2$ satisfies the following degree 2 equation:
$$v^2+pv-\frac{q^2}{4}=0$$

But this latter equation has a unique positive solution, given by:
$$v=\frac{-p+\sqrt{p^2+q^2}}{2}$$

Thus, we are led to the formula of $s=\pm\sqrt{v}$ in the statement, and this gives the result, with the last assertion regarding signs being clear, coming from $2rs=q$.
\end{proof}

Moving now to degree 3 and higher, we would first like to understand what the analogue of $\Delta=b^2-4ac$ is. In order to discuss this question, let us start with:

\begin{proposition}
Given a monic polynomial $P\in\mathbb C[X]$, factorized as
$$P=(X-a_1)\ldots(X-a_k)$$
the following happen:
\begin{enumerate}
\item The coefficients of $P$ are symmetric functions in $a_1,\ldots,a_k$.

\item The symmetric functions in $a_1,\ldots,a_k$ are polynomials in the coefficients of $P$.
\end{enumerate}
\end{proposition}

\begin{proof}
This is something very standard, the idea being as follows:

\medskip

(1) By expanding our polynomial, we have the following formula:
$$P=\sum_{r=0}^k(-1)^r\sum_{i_1<\ldots<i_r}a_{i_1}\ldots a_{i_r}\cdot X^{k-r}$$

Thus the coefficients of $P$ are, up to some signs, the following functions:
$$f_r=\sum_{i_1<\ldots<i_r}a_{i_1}\ldots a_{i_r}$$

But these are indeed symmetric functions in $a_1,\ldots,a_k$, as claimed. 

\medskip

(2) Conversely now, let us look at the symmetric functions in the roots $a_1,\ldots,a_k$. These appear as linear combinations of the basic symmetric functions, given by:
$$S_r=\sum_ia_i^r$$

Moreover, when allowing polynomials instead of linear combinations, we need in fact only the first $k$ such sums, namely $S_1,\ldots,S_k$. That is, the symmetric functions $\mathcal F$ in our variables $a_1,\ldots,a_k$, with integer coefficients, appear as follows:
$$\mathcal F=\mathbb Z[S_1,\ldots,S_k]$$ 

(3) The point now is that, alternatively, the symmetric functions in our variables $a_1,\ldots,a_k$ appear as well as linear combinations of the functions $f_r$ that we found in (1), and that when allowing polynomials instead of linear combinations, we need in fact only the first $k$ functions, namely $f_1,\ldots,f_k$. That is, we have as well:
$$\mathcal F=\mathbb Z[f_1,\ldots,f_k]$$

But this gives the result, because we can pass from $\{S_r\}$ to $\{f_r\}$, and vice versa.
\end{proof}

Getting back now to our original question, regarding finding the analogue of the quantity $\Delta=b^2-4ac$, we have the following key result, putting us on the right track:

\begin{theorem}
Given two polynomials $P,Q\in\mathbb C[X]$, written as
$$P=c(X-a_1)\ldots(X-a_k)\quad,\quad 
Q=d(X-b_1)\ldots(X-b_l)$$
the following quantity, which is called resultant of $P,Q$,
$$R(P,Q)=c^ld^k\prod_{ij}(a_i-b_j)$$
is a certain polynomial in the coefficients of $P,Q$, with integer coefficients, and we have $R(P,Q)=0$ precisely when $P,Q$ have a common root.
\end{theorem}

\begin{proof}
This is something quite tricky, the idea being as follows:

\medskip

(1) Given two polynomials $P,Q\in\mathbb C[X]$, we can certainly construct the quantity $R(P,Q)$ in the statement, with the role of the normalization factor $c^ld^k$ to become clear later on, and then we have $R(P,Q)=0$ precisely when $P,Q$ have a common root:
$$R(P,Q)=0\iff \exists i,j,a_i=b_j$$

(2) As bad news, however, this quantity $R(P,Q)$, defined in this way, is a priori not very useful in practice, because it depends on the roots $a_i,b_j$ of our polynomials $P,Q$, that we cannot compute in general. However, and here comes our point, as we will prove below, it turns out that $R(P,Q)$ is in fact a polynomial in the coefficients of $P,Q$, with integer coefficients, and this is where the power of $R(P,Q)$ comes from.

\medskip

(3) Getting started now, let us expand the formula of $R(P,Q)$, by making all the multiplications there, abstractly, in our head. Everything being symmetric in $a_1,\ldots,a_k$, we obtain in this way certain symmetric functions in these variables, which will be therefore certain polynomials in the coefficients of $P$. Moreover, due to our normalization factor $c^l$, these polynomials in the coefficients of $P$ will have integer coefficients.

\medskip

(4) With this done, let us look now what happens with respect to the remaining variables $b_1,\ldots,b_l$, which are the roots of $Q$. Once again what we have here are certain symmetric functions in these variables $b_1,\ldots,b_l$, and these symmetric functions must be certain polynomials in the coefficients of $Q$. Moreover, due to our normalization factor $d^k$, these polynomials in the coefficients of $Q$ will have integer coefficients.

\medskip

(5) Thus, we are led to the conclusion in the statement, that $R(P,Q)$ is a polynomial in the coefficients of $P,Q$, with integer coefficients, and with the remark that the $c^ld^k$ factor is there for these latter coefficients to be indeed integers, instead of rationals.
\end{proof}

All the above might seem a bit complicated, so as an illustration, let us work out an example. Consider the case of a polynomial of degree 2, and a polynomial of degree 1:
$$P=ax^2+bx+c\quad,\quad 
Q=dx+e$$

In order to compute the resultant, let us factorize our polynomials:
$$P=a(x-p)(x-q)\quad,\quad 
Q=d(x-r)$$

The resultant can be then computed as follows, by using the method above:
\begin{eqnarray*}
R(P,Q)
&=&ad^2(p-r)(q-r)\\
&=&ad^2(pq-(p+q)r+r^2)\\
&=&cd^2+bd^2r+ad^2r^2\\
&=&cd^2-bde+ae^2
\end{eqnarray*}

Finally, observe that $R(P,Q)=0$ corresponds indeed to the fact that $P,Q$ have a common root. Indeed, the root of $Q$ is $r=-e/d$, and we have:
$$P(r)
=\frac{ae^2}{d^2}-\frac{be}{d}+c
=\frac{R(P,Q)}{d^2}$$

Regarding now the explicit formula of the resultant $R(P,Q)$, this is something quite complicated, and there are several methods for dealing with this problem. We have:

\index{resultant}

\begin{theorem}
The resultant of two polynomials, written as 
$$P=p_kX^k+\ldots+p_1X+p_0\quad,\quad 
Q=q_lX^l+\ldots+q_1X+q_0$$
appears as the determinant of an associated matrix, as follows,
$$R(P,Q)=
\begin{vmatrix}
p_k&&&q_l\\
\vdots&\ddots&&\vdots&\ddots\\
p_0&&p_k&q_0&&q_l\\
&\ddots&\vdots&&\ddots&\vdots\\
&&p_0&&&q_0
\end{vmatrix}
$$
with the matrix having size $k+l$, and having $0$ coefficients at the blank spaces.
\end{theorem}

\begin{proof}
This is something quite clever, due to Sylvester, as follows:

\medskip

(1) Consider the vector space $\mathbb C_k[X]$ formed by the polynomials of degree $<k$:
$$\mathbb C_k[X]=\left\{P\in\mathbb C[X]\Big|\deg P<k\right\}$$

This is a vector space of dimension $k$, having as basis the monomials $1,X,\ldots,X^{k-1}$. Now given polynomials $P,Q$ as in the statement, consider the following linear map:
$$\Phi:\mathbb C_l[X]\times\mathbb C_k[X]\to\mathbb C_{k+l}[X]\quad,\quad (A,B)\to AP+BQ$$

\medskip

(2) Our first claim is that with respect to the standard bases for all the vector spaces involved, namely those consisting of the monomials $1,X,X^2,\ldots$, the matrix of $\Phi$ is the matrix in the statement. But this is something which is clear from definitions.

\medskip

(3) Our second claim is that $\det\Phi=0$ happens precisely when $P,Q$ have a common root. Indeed, our polynomials $P,Q$ having a common root means that we can find $A,B$ such that $AP+BQ=0$, and so that $(A,B)\in\ker\Phi$, which reads $\det\Phi=0$.

\medskip

(4) Finally, our claim is that we have $\det\Phi=R(P,Q)$. But this follows from the uniqueness of the resultant, up to a scalar, and with this uniqueness being elementary to establish, along the lines of the proofs of Proposition 8.5 and Theorem 8.6.
\end{proof}

As an illustration for this, consider our favorite polynomials, as before:
$$P=ax^2+bx+c\quad,\quad 
Q=dx+e$$

According to the above result, the resultant should be then, as it should:
$$R(P,Q)
=\begin{vmatrix}
a&d&0\\
b&e&d\\
c&0&e
\end{vmatrix}
=ae^2-bde+cd^2$$

We can go back now to our original question, regarding $\Delta$, and we have:

\index{discriminant}

\begin{theorem}
Given a polynomial $P\in\mathbb C[X]$, written as
$$P(X)=aX^N+bX^{N-1}+cX^{N-2}+\ldots$$
its discriminant, defined as being the following quantity,
$$\Delta(P)=\frac{(-1)^{\binom{N}{2}}}{a}\,R(P,P')$$
is a polynomial in the coefficients of $P$, with integer coefficients, and $\Delta(P)=0$ happens precisely when $P$ has a double root.
\end{theorem}

\begin{proof}
The fact that the discriminant $\Delta(P)$ is a polynomial in the coefficients of $P$, with integer coefficients, comes from Theorem 8.6, coupled with the fact that the division by the leading coefficient $a$ is indeed possible, under $\mathbb Z$, as being shown by the following formula, which is written a bit informally, coming from Theorem 8.7:
$$R(P,P')=
\begin{vmatrix}
a&&&Na\\
\vdots&\ddots&&\vdots&\ddots\\
z&&a&y&&Na\\
&\ddots&\vdots&&\ddots&\vdots\\
&&z&&&y
\end{vmatrix}
$$

Also, the fact that we have $\Delta(P)=0$ precisely when $P$ has a double root is clear from Theorem 8.6. Finally, let us mention that the sign $(-1)^{\binom{N}{2}}$ is there for various reasons, including the compatibility with $\Delta=b^2-4ac$ in degree 2, that we will discuss next.
\end{proof}

As a first illustration, let us see what happens in degree 2. We have:
$$P=aX^2+bX+c\quad,\quad 
P'=2aX+b$$

Thus, the discriminant is given by the following formula, as it should:
$$\Delta(P)
=-\frac{1}{a}\begin{vmatrix}
a&2a&\\
b&b&2a\\
c&&b
\end{vmatrix}=b^2-4ac$$

At the theoretical level now, we have the following result, which is non-trivial:

\begin{theorem}
The discriminant of a polynomial $P$ is given by the formula
$$\Delta(P)=a^{2N-2}\prod_{i<j}(r_i-r_j)^2$$
where $a$ is the leading coefficient, and $r_1,\ldots,r_N$ are the roots.
\end{theorem}

\begin{proof}
This is something quite tricky, the idea being as follows:

\medskip

(1) The first thought here goes to the formula in Theorem 8.6, so let us see what that formula teaches us, in the case $Q=P'$. Let us write $P,P'$ as follows:
$$P=a(x-r_1)\ldots(x-r_N)$$ 
$$P'=Na(x-p_1)\ldots(x-p_{N-1})$$

According to Theorem 8.6, the resultant of $P,P'$ is then given by:
$$R(P,P')=a^{N-1}(Na)^N\prod_{ij}(r_i-p_j)$$

And bad news, this is not exactly what we wished for, namely the formula in the statement. That is, we are on the good way, but certainly have to work some more.

\medskip

(2) Obviously, we must get rid of the roots $p_1,\ldots,p_{N-1}$ of the polynomial $P'$. In order to do this, let us rewrite the formula that we found in (1) in the following way:
\begin{eqnarray*}
R(P,P')
&=&N^Na^{2N-1}\prod_i\left(\prod_j(r_i-p_j)\right)\\
&=&N^Na^{2N-1}\prod_i\frac{P'(r_i)}{Na}\\
&=&a^{N-1}\prod_iP'(r_i)
\end{eqnarray*}

(3) In order to compute now $P'$, and more specifically the values $P'(r_i)$ that we are interested in, we can use the Leibnitz rule. So, consider our polynomial:
$$P(x)=a(x-r_1)\ldots(x-r_N)$$ 

The Leibnitz rule for derivatives tells us that $(fg)'=f'g+fg'$, but then also that $(fgh)'=f'gh+fg'h+fgh'$, and so on. Thus, for our polynomial, we obtain:
$$P'(x)=a\sum_i(x-r_1)\ldots\underbrace{(x-r_i)}_{missing}\ldots(x-r_N)$$ 

Now when applying this formula to one of the roots $r_i$, we obtain:
$$P'(r_i)=a(r_i-r_1)\ldots\underbrace{(r_i-r_i)}_{missing}\ldots(r_i-r_N)$$ 

By making now the product over all indices $i$, this gives the following formula:
$$\prod_iP'(r_i)=a^N\prod_{i\neq j}(r_i-r_j)$$

(4) Time now to put everything together. By taking the formula in (2), making the normalizations in Theorem 8.8, and then using the formula found in (3), we obtain:
\begin{eqnarray*}
\Delta(P)
&=&(-1)^{\binom{N}{2}}a^{N-2}\prod_iP'(r_i)\\
&=&(-1)^{\binom{N}{2}}a^{2N-2}\prod_{i\neq j}(r_i-r_j)\\
&=&a^{2N-2}\prod_{i<j}(r_i-r_j)^2
\end{eqnarray*}

Thus, we are led to the conclusion in the statement.
\end{proof}

\section*{8b. Cardano formulae}

Getting now to the point where we wanted to get, solving equations of small degree $N=2,3,4,\ldots\,$, let us start with some applications of Theorem 8.9. We have:

\index{real roots}
\index{complex roots}

\begin{proposition}
Consider a polynomial with real coefficients, $P\in\mathbb R[X]$, assumed for simplicity to have nonzero discriminant, $\Delta\neq0$.
\begin{enumerate}
\item In degree $2$, the roots are real when $\Delta>0$, and complex when $\Delta<0$.

\item In degree $3$, all roots are real precisely when $\Delta>0$.
\end{enumerate}
\end{proposition}

\begin{proof}
This is indeed something self-explanatory, coming from Theorem 8.9.
\end{proof}

Let us work out now explicitly what happens in degree 3. We must first compute the discriminant of arbitrary degree 3 polynomials, and the result here is as follows:

\begin{theorem}
The discriminant of a degree $3$ polynomial,
$$P=aX^3+bX^2+cX+d$$
is given by $\Delta(P)=b^2c^2-4ac^3-4b^3d-27a^2d^2+18abcd$.
\end{theorem}

\begin{proof}
Consider indeed a polynomial of degree 3, and its derivative:
$$P=aX^3+bX^2+cX+d$$
$$P'=3aX^2+2bX+c$$

By using Theorem 8.7 and computing the determinant, we obtain:
\begin{eqnarray*}
R(P,P')
&=&\begin{vmatrix}
a&&3a\\
b&a&2b&3a\\
c&b&c&2b&3a\\
d&c&&c&2b\\
&d&&&c
\end{vmatrix}\\
&=&a(4b^3d+27a^2d^2-18abcd+4ac^3-b^2c^2)
\end{eqnarray*}

Now according to Theorem 8.8, the discriminant of our polynomial is given by:
\begin{eqnarray*}
\Delta(P)
&=&-\frac{R(P,P')}{a}\\
&=&-4b^3d-27a^2d^2+18abcd-4ac^3+b^2c^2\\
&=&b^2c^2-4ac^3-4b^3d-27a^2d^2+18abcd
\end{eqnarray*}

Thus, we have obtained the formula in the statement.
\end{proof}

Still talking degree 3 equations, let us try now to solve such an equation $P=0$, with $P=aX^3+bX^2+cX+d$ as above. By linear transformations we can assume $a=1,b=0$, and then it is convenient to write $c=3p,d=2q$. Thus, our equation becomes:
$$x^3+3px+2q=0$$

Regarding such equations, many things can be said, and we notably have the following famous result, due to Cardano, solving the problem in the case $\Delta<0$:

\index{degree 3 equation}
\index{Cardano formula}
\index{depressed cubic}

\begin{theorem}
For a normalized degree $3$ equation, namely 
$$x^3+3px+2q=0$$
the discriminant is $\Delta=-108(p^3+q^2)$. Assuming $p,q\in\mathbb R$ and $\Delta<0$, the numbers
$$x=w\sqrt[3]{-q+\sqrt{p^3+q^2}}+w^2\sqrt[3]{-q-\sqrt{p^3+q^2}}$$
with $w=1,e^{2\pi i/3},e^{4\pi i/3}$ are the solutions of our equation.
\end{theorem}

\begin{proof}
The formula of $\Delta$ is clear from the above, and with $108=4\times 27$. Now with $x$ as in the statement, by using $(a+b)^3=a^3+b^3+3ab(a+b)$, we have:
\begin{eqnarray*}
x^3
&=&\left(w\sqrt[3]{-q+\sqrt{p^3+q^2}}+w^2\sqrt[3]{-q-\sqrt{p^3+q^2}}\right)^3\\
&=&-2q+3\sqrt[3]{-q+\sqrt{p^3+q^2}}\cdot\sqrt[3]{-q-\sqrt{p^3+q^2}}\cdot x\\
&=&-2q+3\sqrt[3]{q^2-p^3-q^2}\cdot x\\
&=&-2q-3px
\end{eqnarray*}

Thus, we are led to the conclusion in the statement.
\end{proof}

As a comment here, the formula in Theorem 8.12 holds of course in the case $\Delta>0$ too, and also when the coefficients are complex numbers, $p,q\in\mathbb C$. However, these extensions are not useful, because when it comes to extract the above square and cubic roots, for complex numbers, you end up with the initial question, solving a degree 3 equation. 

\bigskip

In higher degree things become quite complicated. In degree 4, to start with, we first have the following result, dealing with the discriminant and its applications:

\index{degree 4 polynomial}

\begin{theorem}
The discriminant of $P=ax^4+bx^3+cx^2+dx+e$ is given by:
\begin{eqnarray*}
\Delta
&=&256a^3e^3-192a^2bde^2-128a^2c^2e^2+144a^2cd^2e-27a^2d^4\\
&&+144ab^2ce^2-6ab^2d^2e-80abc^2de+18abcd^3+16ac^4e\\
&&-4ac^3d^2-27b^4e^2+18b^3cde-4b^3d^3-4b^2c^3e+b^2c^2d^2
\end{eqnarray*}
In the case $\Delta<0$ we  have $2$ real roots and $2$ complex conjugate roots, and in the case $\Delta>0$ the roots are either all real or all complex.
\end{theorem}

\begin{proof}
The formula of $\Delta$ follows from the definition of the discriminant, from Theorem 8.8, with the resultant computed via Theorem 8.7, as follows:
$$\Delta
=\frac{1}{a}\begin{vmatrix}
a&&&4a\\
b&a&&3b&4a\\
c&b&a&2c&3b&4a\\
d&c&b&d&2c&3b&4a\\
e&d&c&&d&2c&3b\\
&e&d&&&d&2c\\
&&e&&&&d
\end{vmatrix}$$

As for the last assertion, the study here is routine, a bit as in degree $N=2,3$.
\end{proof}

In practice, as in degree 3, we can do first some manipulations on our polynomials, as to have them in simpler form, and we have the following version of Theorem 8.13:

\index{depressed quartic}

\begin{theorem}
The discriminant of a normalized degree $4$ polynomial, written as
$$P=x^4+6px^2+4qx+3r$$ 
is given by the following formula:
$$\Delta=256\times27\times\big(9p^4r-2p^3q^2-6p^2r^2+6pq^2r-q^4+r^3\big)$$
In the case $\Delta<0$ we have $2$ real roots and $2$ complex conjugate roots, and in the case $\Delta>0$ the roots are either all real or all complex.
\end{theorem}

\begin{proof}
This follows from Theorem 8.13, with $a,b,c,d,e$ as above, but we can deduce this as well directly. Indeed, the formula of $\Delta$ follows, quite easily, from:
$$\Delta
=\begin{vmatrix}
1&&&4\\
&1&&&4\\
6p&&1&12p&&4\\
4q&6p&&4q&12p&&4\\
3r&4q&6p&&4q&12p&\\
&3r&4q&&&4q&12p\\
&&3r&&&&4q
\end{vmatrix}$$

As for the last assertion, this is something that we know from Theorem 8.13.
\end{proof}

Time now to get to the real thing, solving the equation. We have here:

\index{degree 4 equation}
\index{Cardano formula}

\begin{theorem}
The roots of a normalized degree $4$ equation, written as
$$x^4+6px^2+4qx+3r=0$$ 
are as follows, with $y$ satisfying the equation $(y^2-3r)(y-3p)=2q^2$,
$$x_1=\frac{1}{\sqrt{2}}\left(-\sqrt{y-3p}+\sqrt{-y-3p+\frac{4q}{\sqrt{2y-6p}}}\right)$$
$$x_2=\frac{1}{\sqrt{2}}\left(-\sqrt{y-3p}-\sqrt{-y-3p+\frac{4q}{\sqrt{2y-6p}}}\right)$$
$$x_3=\frac{1}{\sqrt{2}}\left(\sqrt{y-3p}+\sqrt{-y-3p-\frac{4q}{\sqrt{2y-6p}}}\right)$$
$$x_4=\frac{1}{\sqrt{2}}\left(\sqrt{y-3p}-\sqrt{-y-3p-\frac{4q}{\sqrt{2y-6p}}}\right)$$
and with $y$ being computable via the Cardano formula.
\end{theorem}

\begin{proof}
This is something quite tricky, the idea being as follows:

\medskip

(1) To start with, let us write our equation in the following form:
$$x^4=-6px^2-4qx-3r$$

The idea will be that of adding a suitable common term, to both sides, as to make square on both sides, as to eventually end with a sort of double quadratic equation. For this purpose, our claim is that what we need is a number $y$ satisfying:
$$(y^2-3r)(y-3p)=2q^2$$

Indeed, assuming that we have this number $y$, our equation becomes:
\begin{eqnarray*}
(x^2+y)^2
&=&x^4+2x^2y+y^2\\
&=&-6px^2-4qx-3r+2x^2y+y^2\\
&=&(2y-6p)x^2-4qx+y^2-3r\\
&=&(2y-6p)x^2-4qx+\frac{2q^2}{y-3p}\\
&=&\left(\sqrt{2y-6p}\cdot x-\frac{2q}{\sqrt{2y-6p}}\right)^2
\end{eqnarray*}

(2) Which looks very good, leading us to the following degree 2 equations:
$$x^2+y+\sqrt{2y-6p}\cdot x-\frac{2q}{\sqrt{2y-6p}}=0$$
$$x^2+y-\sqrt{2y-6p}\cdot x+\frac{2q}{\sqrt{2y-6p}}=0$$

Now let us write these two degree 2 equations in standard form, as follows:
$$x^2+\sqrt{2y-6p}\cdot x+\left(y-\frac{2q}{\sqrt{2y-6p}}\right)=0$$
$$x^2-\sqrt{2y-6p}\cdot x+\left(y+\frac{2q}{\sqrt{2y-6p}}\right)=0$$

(3) Regarding the first equation, the solutions there are as follows:
$$x_1=\frac{1}{2}\left(-\sqrt{2y-6p}+\sqrt{-2y-6p+\frac{8q}{\sqrt{2y-6p}}}\right)$$
$$x_2=\frac{1}{2}\left(-\sqrt{2y-6p}-\sqrt{-2y-6p+\frac{8q}{\sqrt{2y-6p}}}\right)$$

As for the second equation, the solutions there are as follows:
$$x_3=\frac{1}{2}\left(\sqrt{2y-6p}+\sqrt{-2y-6p-\frac{8q}{\sqrt{2y-6p}}}\right)$$
$$x_4=\frac{1}{2}\left(\sqrt{2y-6p}-\sqrt{-2y-6p-\frac{8q}{\sqrt{2y-6p}}}\right)$$

(4) Now by cutting a $\sqrt{2}$ factor from everything, this gives the formulae in the statement. As for the last claim, regarding the nature of $y$, this comes from Cardano.
\end{proof}

We still have to compute the number $y$ appearing in the above via Cardano, and the result here, adding to what we already have in Theorem 8.15, is as follows:

\begin{theorem}[continuation]
The value of $y$ in the previous theorem is
$$y=t+p+\frac{a}{t}$$
where the number $t$ is given by the formula
$$t=\sqrt[3]{b+\sqrt{b^2-a^3}}$$
with $a=p^2+r$ and $b=2p^2-3pr+q^2$.
\end{theorem}

\begin{proof}
The legend has it that this is what comes from Cardano, but depressing and normalizing and solving $(y^2-3r)(y-3p)=2q^2$ makes it for too many operations, so the most pragmatic is to simply check this equation. With $y$ as above, we have:
\begin{eqnarray*}
y^2-3r
&=&t^2+2pt+(p^2+2a)+\frac{2pa}{t}+\frac{a^2}{t^2}-3r\\
&=&t^2+2pt+(3p^2-r)+\frac{2pa}{t}+\frac{a^2}{t^2}
\end{eqnarray*}

With this in hand, we have the following computation:
\begin{eqnarray*}
(y^2-3r)(y-3p)
&=&\left(t^2+2pt+(3p^2-r)+\frac{2pa}{t}+\frac{a^2}{t^2}\right)\left(t-2p+\frac{a}{t}\right)\\
&=&t^3+(a-4p^2+3p^2-r)t+(2pa-6p^3+2pr+2pa)\\
&&+(3p^2a-ra-4p^2a+a^2)\frac{1}{t}+\frac{a^3}{t^3}\\
&=&t^3+(a-p^2-r)t+2p(2a-3p^2+r)+a(a-p^2-r)\frac{1}{t}+\frac{a^3}{t^3}\\
&=&t^3+2p(-p^2+3r)+\frac{a^3}{t^3}
\end{eqnarray*}

Now by using the formula of $t$ in the statement, this gives:
\begin{eqnarray*}
(y^2-3r)(y-3p)
&=&b+\sqrt{b^2-a^3}-4p^2+6pr+\frac{a^3}{b+\sqrt{b^2-a^3}}\\
&=&b+\sqrt{b^2-a^3}-4p^2+6pr+b-\sqrt{b^2-a^3}\\
&=&2b-4p^2+6pr\\
&=&2(2p^2-3pr+q^2)-4p^2+6pr\\
&=&2q^2
\end{eqnarray*}

Thus, we are led to the conclusion in the statement.
\end{proof}

\section*{8c. Abstract algebra}

In degree 5 and higher things get considerably more complicated, with no formula for the solutions being available, or even possible. In order to discuss this phenomenon, following Abel, Galois and others, we will need some algebra. Let us start with:

\begin{definition}
We have notions of rings, modules and ideals, as follows:
\begin{enumerate}
\item A ring $R$ is a set with operations $+$ and $\times$, satisfying the usual conditions for such operations, except for $ab=ba$, and for $a\neq0\implies\exists a^{-1}$.

\item A module $V$ over a ring $R$ is a vector space, but we will call it ring, and keep the name vector spaces for the modules over fields, $R=F$.

\item An ideal $I\subset R$ is a subgroup with the left ideal property $i\in I,r\in R\implies ir\in I$, or the right ideal property $i\in I,r\in R\implies ri\in I$, or both.
\end{enumerate}
\end{definition}

This was a quite crowded statement, but you get the point, with (1) and (2) we are sort of trying to do field and vector space mathematics, over things which are not necessarily fields and vector spaces over them, and (3) is something technical, non-field specific. At the level of examples, these abound, and we have three important ones, as follows:

\bigskip

(1) The integers form a ring, $R=\mathbb Z$, which in addition is commutative, $ab=ba$. As obvious module over $\mathbb Z$, we have the lattice $V=\mathbb Z^N$. Finally, since $R=\mathbb Z$ is commutative, the 3 notions of ideals coincide, and these are the subsets $I=a\mathbb Z$, with $a\in\mathbb Z$.

\bigskip

(2) The matrices over the integers form a ring, $R=M_k(\mathbb Z)$, which is noncommutative at $k>1$. As obvious module over $M_k(\mathbb Z)$, we have the lattice $V=\mathbb Z^k$. As for the ideals, things here are a bit more complicated, but since at $k=2$ the matrices of type $(^a_0{\ }^b_0)$ form a left ideal which is not a right ideal, and the matrices of type $(^a_b{\ }^0_0)$ form a right ideal which is not a left ideal, at least we know that our 3 types of ideals make sense.

\bigskip

(3) In relation with arithmetic, the integers modulo $N\in\mathbb N$ form a ring $R=\mathbb Z_N$, which is commutative. Also, the matrices over these integers modulo $N$ form a ring $R'=M_k(\mathbb Z_N)$, which is noncommutative at $k>1$. As for a preliminary study of the modules and ideals, for these two rings, I will leave it to you, and enjoy.

\bigskip

The question that you surely have in mind is, what are ideals good for? Answer:

\begin{proposition}
For a subgroup $I\subset R$, the following are equivalent:
\begin{enumerate}
\item $I$ is a two-sided ideal.

\item $R/I$ is a ring.
\end{enumerate}
\end{proposition}

\begin{proof}
This is something which requires some thinking, as follows:

\medskip

(1) Since the additive group $(R,+)$ is abelian, given an additive subgroup $I\subset R$ we can form the quotient group $R/I$, which is abelian too, with addition as follows:
$$(a+I)+(b+I)=(a+b+I)$$

(2) The question is now, can we turn this abelian group $R/I$ into a ring? Normally the multiplication can only be as follows, and with this clarifying our statement:
$$(a+I)(b+I)=(ab+I)$$

(3) But, will this work. As a first observation, there is a bit of analogy here with group theory, where $H\subset G$ must be normal in order for $G/H$ to be a group. Thus, our claim is that the ideal condition is somehow the ``analogue of normality, in the ring setting''.

\medskip

(4) In practice now, it is quite clear, exactly as in the group theory setting, that everything will be fine, provided that our multiplication is well-defined. And for this multiplication to be well-defined, the following condition must be satisfied:
$$(a+I)=(a'+I)\ ,\ (b+I)=(b'+I)\quad\implies\quad (ab+I)=(a'b'+I)$$

But this amounts in the following condition to be satisfied:
$$a-a'\in I\ ,\ b-b'\in I\quad\implies\quad ab-a'b'\in I$$

(5) Now comes the math. We have the following identity, which shows that if $I\subset R$ is a two-sided ideal, then the above condition is satisfied, and so done:
$$ab-a'b'=a(b-b')+(a-a')b'$$

(6) Conversely now, if the condition in (4) is satisfied, we have in particular:
$$i-0\in I\ ,\ r-r\in I\quad\implies\quad ir-0r\in I$$
$$r-r\in I\ ,\ i-0\in I\quad\implies\quad ri-r0\in I$$

Thus $I\subset R$ must be a two-sided ideal, and this finishes the proof.
\end{proof}

Many other things can be said about rings, modules and ideals, especially in the commutative case. For formulating a theorem on the subject, we have:

\begin{theorem}
Assuming that $R$ is commutative and $I\subset R$ is a maximal ideal, in the sense that it is a proper ideal, $I\neq R$, and there is no bigger proper ideal 
$$I\subset J\subset R$$
the quotient ring $F=R/I$ is a field.
\end{theorem}

\begin{proof}
Here is the proof, and with this being guaranteed to be useful learning:

\medskip

(1) Before starting, a quick example. We know that over $R=\mathbb Z$, the ideals are the subsets $I=p\mathbb Z$ with $p\in\mathbb N$. But such an ideal is maximal precisely when $p$ is prime, and this is the same as asking for the quotient ring $R/I=\mathbb Z_p$ to be a field. 

\medskip

(2) In general now, assume first that $R/I$ is a field. This means that any nonzero element of $R/I$ is invertible, and with our usual conventions for $R/I$, this reads:
$$\forall a\notin I\ ,\ \exists b\in R\ ,\ (ab+I)=(1+I)$$

Now assume by contradiction that $I\subset R$ is not maximal, so that we have a bigger ideal $I\subset J\subset R$. If we pick $a\in J-I$, we obtain, by the above, the following:
$$a\in J-I\ ,\ b\in R\ ,\ ab=1+i\ ,\ i\in I$$

But this is contradictory, because since $J$ is an ideal, containing $I$, we must have $ab,i\in J$, so we conclude that we have $1\in J$, and so $J=R$, contradiction.

\medskip

(3) Conversely, assume now that $I$ is maximal, and assume too, by contradiction, that $R/I$ is not a field. Then we can find a zero divisor in $R/I$, which reads:
$$(a+I)(b+I)=(I)\ ,\ a,b\notin I$$

In other words, we can find $ab\in I$ with $a,b\notin I$. But then, let us look at:
$$I\subset I+aR\subset R$$

(4) What we have in the middle is an ideal, and it is also clear, from $a\notin I$, that the inclusion on the left is proper. As for the inclusion on the right, our claim is that this is proper too. Indeed, assuming otherwise, we would have a formula as follows:
$$i+ac=1\ ,\ i\in I$$

Now by multiplying everything by $b$, we obtain from this:
$$ib+acb=b\ ,\ i\in I$$

But this is contradictory, because on the left we have $ib\in I$ and $acb=(ab)c\in I$, which gives $b\in I$, contradicting the condition $b\notin I$. Thus, our claim is proved.

\medskip

(5) But this is the end of the story, because what we just proved is that what we have in (3) is indeed a proper ideal, contradicting the maximality of $I$, as desired.
\end{proof}

Still with me, I hope, after all these abstractions, and please believe me, Theorem 8.19 is something of key importance, be that for algebraic geometry, or for arithmetic, as we will soon discover, or even for analysis, in the context of the Banach algebras.

\bigskip

Going ahead now with our general abstract algebra program, let us formulate:

\begin{definition}
We have notions of fields, vector spaces and algebras, as follows:
\begin{enumerate}
\item A field $F$ is a field $F$ as we know them, with in algebra parlance these being the commutative rings $R$ with each nonzero element being invertible.

\item A vector space $V$ over a field $F$ is a vector space as we know them, in algebra parlance these being the modules $V$ over a field $F$.

\item An algebra $A$ over a field $F$ is a vector space over $F$, with a ring multiplication operation $\times$, compatible with the vector space structure.
\end{enumerate}
\end{definition}

In relation with this, fields are something quite familiar, and we notably have here $\mathbb Q,\mathbb R,\mathbb C$, or the field $\mathbb F_p$ of remainders modulo a prime number $p$. In what regards the vector spaces, we know about them since ever, and finally, regarding algebras, we know many algebras of functions from analysis. But, thinking well, from a purely algebraic perspective, all these objects have many operations, and this is why they came at last.

\bigskip

As basic examples now, passed the fields $F$ and the vector spaces $V$ that we know well, we are left with finding interesting examples of algebras $A$. And here the examples abound,  with all being related to geometry or analysis of some sort, as follows: 

\bigskip

(1)  First we have algebra of polynomials $A=F[X]$. This is a very basic algebra, important for geometry, and with the extra feature that it is commutative, $PQ=QP$.

\bigskip

(2) More generally, we have the algebra of polynomials $A=F[X_1,\ldots,X_N]$. Again, this algebra is important for algebraic geometry, and is commutative, $PQ=QP$.

\bigskip

(3) Still talking commutative algebras, we have many of them coming from analysis, the general principle being that ``functions form algebras''. More on this in a moment.

\bigskip

(4) We have as well the algebra of matrices $A=M_N(F)$. Again this is a very basic example, that we know well, which this time is not commutative, $ST\neq TS$.

\bigskip

More in detail now, getting to the various algebras of functions, we have here the following key result, bringing among others some further light on Theorem 8.19 too:

\begin{theorem}
Given a compact space $X$, the following happen:
\begin{enumerate}
\item The continuous functions $f:X\to\mathbb C$ form a complex algebra $C(X)$.

\item Given $x\in X$, the functions satisfying $f(x)=0$, form an ideal $I\subset C(X)$. 

\item This ideal is maximal, and any maximal ideal $I\subset C(X)$ appears in this way. 

\item In this picture, the fact that the quotient is a field, $C(X)/I=\mathbb C$, is clear.
\end{enumerate}
\end{theorem}

\begin{proof}
All this is self-explanatory, the idea being as follows:

\medskip

(1) This is clear. Observe that our algebra is commutative, $fg=gf$.

\medskip

(2) This is again clear, because $f(x)=0$ implies $(fg)(x)=0$.

\medskip

(3) This follows from basic topology, via a suitable open cover for $X$.

\medskip

(4) This is clear, because $C(X)\to C(X)/I$ maps $f\to f(x)\in\mathbb C$.
\end{proof}

And good news, that is all. Done with abstract algebra, and good learning this was, that we can use later when needed, and we can now turn to more concrete things. 

\section*{8d. Galois theory}

We are now ready for Galois theory. Let us start with a brief discussion about fields, and their characteristic. The basic results here can be summarized as follows:

\index{field}
\index{arbitrary field}
\index{characteristic of field}
\index{prime field}
\index{Fermat theorem}
\index{Fermat polynomial}
\index{finite field}
\index{multiplicative group}

\begin{theorem}
Given a field $F$, define its characteristic $p=char(F)$ as being the smallest $p\in\mathbb N$ such that the following happens, and as $p=0$, if this never happens:
$$\underbrace{1+\ldots+1}_{p\ times}=0$$
Then, assuming $p>0$, this number $p$ must be prime, we have a field embedding $\mathbb F_p\subset F$, and $q=|F|$ must be of the form $q=p^k$, with $k\in\mathbb N$. Also, we have the formulae
$$(a+b)^p=a^p+b^p\quad,\quad a^q=a$$
valid for any $a,b\in F$, and the Fermat poynomial $X^q-X$ factorizes as:
$$X^q-X=\prod_{a\in F}(X-a)$$
Also, regardless of $p$, any finite multiplicative subgroup $G\subset F-\{0\}$ must be cyclic.
\end{theorem}

\begin{proof}
This is a very crowded statement, containing all sorts of things that we know well from basic algebra, the idea with all this being as follows:

\medskip

(1) The fact that $p>0$ must be prime comes by contradiction, by using:
$$(\underbrace{1+\ldots+1}_{a\ times})\times(\underbrace{1+\ldots+1}_{b\ times})=\underbrace{1+\ldots+1}_{ab\ times}$$

Indeed, assuming that we have $p=ab$ with $a,b>1$, the above formula corresponds to an equality of type $AB=0$ with $A,B\neq0$ inside $F$, which is impossible.

\medskip

(2) Back to the general case, $F$ has a smallest subfield $E\subset F$, called prime field, consisting of the various sums $1+\ldots+1$, and their quotients. In the case $p=0$ we obviously have $E=\mathbb Q$. In the case $p>0$ now, the multiplication formula in (1) shows that the set $S=\{1+\ldots+1\}$ is stable under taking quotients, and so $E=S$.

\medskip

(3) Now with $E=S$ in hand, we obviously have $(E,+)=\mathbb Z_p$, and since the multiplication is given by the formula in (1), we conclude that we have $E=\mathbb F_p$, as a field. Thus, in the case $p>0$, we have constructed an embedding $\mathbb F_p\subset F$, as claimed.

\medskip

(4) In the context of the above embedding $\mathbb F_p\subset F$, we can say that $F$ is a vector space over $\mathbb F_p$, and so we have $|F|=p^k$, with $k\in\mathbb N$ being the dimension of this space.

\medskip

(5) The baby Fermat formula $(a+b)^p=a^p+b^p$, which reminds the Fermat little theorem, $a^p=a(p)$ over $\mathbb Z$, follows in the same way, namely from the binomial formula, because all the non-trivial binomial coefficients $\binom{p}{s}$ are multiples of $p$:
$$(a+b)^p=\sum_{k=0}^p\binom{p}{k}a^kb^{p-k}=a^p+b^p$$

(6) As for the Fermat formula $a^q=a$ itself, which implies the assertion about $X^q-X$, this follows from the last assertion, which can be proved via some basic arithmetic inside $F$, and which for $G=F-\{0\}$ itself, with $|F|=q$, gives $a^{q-1}=1$, for any $a\neq0$.

\medskip

(7) Let us pick indeed an element $g\in G$ of highest order, $n=ord(g)$. Our claim, which will prove the results, is that the order $m=ord(h)$ of any $h\in G$ satisfies $m|n$.

\medskip

(8) In order to prove this claim, let $d=(m,n)$, write $d=am+bn$ with $a,b\in\mathbb Z$, and set $k=g^ah^b$. We have then the following computation:
$$k^m
=g^{am}h^{bm}
=g^{am}
=g^{d-bn}
=g^d$$

Similarly, $k^n=h^d$. By using $k^m=g^d$, we obtain the following formula:
$$k^{[m,n]}
=k^{mn/d}
=(k^m)^{n/d}
=(g^d)^{n/d}
=g^n
=1$$

Thus $ord(k)|[m,n]$, and our claim is that we have in fact $ord(k)=[m,n]$.

\medskip

(9) In order to prove this latter claim, assume first that we are in the case $d=1$. But here the result is clear, because the formulae $k^m=g^d$ and $k^n=h^d$ in (8) read:
$$g=k^m\quad,\quad h=g^n$$

Indeed, since the numbers $n=ord(g)$ and $m=ord(g)$ are prime to each other, we conclude from this that we have $ord(k)=mn$, as desired. Thus, result proved for $d=1$, and in what regards the general case, where $d$ is arbitrary, this follows from this.

\medskip

(10) Summarizing, we have proved our claim in (8). Now since the order $n=ord(g)$ was assumed to be maximal, we must have $[m,n]|n$, and so $m|n$. Thus, we have proved our claim in (7), namely that the order $m=ord(h)$ of any $h\in G$ satisfies $m|n$.

\medskip

(11) But with this claim in hand, the result follows. Indeed, since the polynomial $x^n-1$ has all the elements $h\in G$ as roots, its degree must satisfy $n\geq|G|$. On the other hand, from $n=ord(g)$ with $g\in G$, we have $n||G|$. We therefore conclude that we have $n=|G|$, which shows that $G$ is indeed cyclic, generated by the element $g\in G$.

\medskip

(12) Finally, assuming $|F|=q<\infty$, we know that the multiplicative group $F-\{0\}$ is cyclic, of order $q-1$. Thus, the following formula is satisfied, for any $a\in F-\{0\}$:
$$a^{q-1}=1$$

Now by multiplying by $a$, this gives the Fermat formula $a^q=a$, with of course the remark that this formula trivially holds as well for $a=0$.
\end{proof}

The above result raises many questions. Since most of these questions seem to have something to do with field extensions, let us start by discussing this. We first have:

\index{field extension}
\index{separable extension}
\index{Galois theorem}
\index{root of polynomial}

\begin{theorem} 
Given a field extension $E\subset F$, we can talk about its Galois group $G$, as the group of automorphisms of $F$ fixing $E$. The intermediate fields
$$E\subset K\subset F$$
are then in correspondence with the subgroups $H\subset G$, with such a field $K$ corresponding to the subgroup $H$ consisting of automorphisms $g\in G$ fixing $K$.
\end{theorem}

\begin{proof}
This is something quite self-explanatory, and follows indeed from some algebra, under suitable assumptions, in order for that algebra to properly apply. There are many good books, that you can learn this technology from, such as Lang \cite{la1}.
\end{proof}

Getting now towards polynomials and their roots, we have here:

\index{splitting field}
\index{algebraic closure}
\index{algebrically closed}

\begin{theorem}
Given a field $F$ and a polynomial $P\in F[X]$, we can talk about the abstract splitting field of $P$, where this polynomial decomposes as:
$$P(X)=c\prod_i(X-a_i)$$
In particular, any field $F$ has a certain algebraic closure $\bar{F}$, where all the polynomials $P\in F[X]$, and in fact all polynomials $P\in\bar{F}[X]$ too, have roots.
\end{theorem}

\begin{proof}
This is again something self-explanatory, which follows from Theorem 8.23 and from some extra algebra, under suitable assumptions, in order for that extra algebra to properly apply. Regarding the construction at the end, as main example here we have $\bar{\mathbb R}=\mathbb C$. However, as an interesting fact, $\bar{\mathbb Q}\subset\mathbb C$ is a proper subfield. See \cite{la1}.
\end{proof}

Good news, with this in hand, we can now elucidate the structure of finite fields:

\index{finite field}
\index{uniqueness of finite fields}

\begin{theorem}
For any prime power $q=p^k$ there is a unique field $\mathbb F_q$ having $q$ elements. At $k=1$ this is the usual $\mathbb F_p$. In general, this is the splitting field of:
$$P=X^q-X$$
Moreover, we can construct an explicit model for $\mathbb F_q$, at $q=p^2$ or higher, as
$$\mathbb F_q=\mathbb F_p[X]/(Q)$$
with $Q\in\mathbb F_p[X]$ being a suitable irreducible polynomial, of degree $k$. 
\end{theorem}

\begin{proof}
There are several assertions here, the idea being as follows:

\medskip

(1) The first assertion, regarding the existence and uniqueness of $\mathbb F_q$, follows from Theorem 8.22 and Theorem 8.24. Indeed, we know from Theorem 8.30 that given a finite field, $|F|=q$ with $k\in\mathbb N$, the Fermat polynomial $P=X^q-X$ factorizes as:
$$X^q-X=\prod_{a\in F}(X-a)$$

But this shows, via the general theory from Theorem 8.24, that our field $F$ must be the splitting field of $P$, and so is unique. As for the existence, this follows again from Theorem 8.22, telling us that the splitting field always exists.

\medskip

(2) In what regards now the modeling of $\mathbb F_q$, at $q=p$ there is nothing to do, because we have our usual $\mathbb F_p$ here. At $q=p^2$ and higher, we know from commutative algebra that we have an isomorphism as follows, whenever $Q\in\mathbb F_p[X]$ is taken irreducible:
$$\mathbb F_q=\mathbb F_p[X]/(Q)$$

(3) Regarding now the best choice of the irreducible polynomial $Q\in\mathbb F_p[X]$, providing us with a good model for the finite field $\mathbb F_q$, that we can use in practice, this question depends on the value of $q=p^k$, and many things can be said here. All in all, our models are quite similar to $\mathbb C=\mathbb R[i]$, with $i$ being a formal number satisfying $i^2=-1$.

\medskip

(4) To be more precise, at the simplest exponent, $q=4$, to start with, we can use $Q=X^2+X+1$, with this being actually the unique possible choice of a degree 2 irreducible polynomial $Q\in\mathbb F_2[X]$, and this leads to a model as follows:
$$\mathbb F_4=\left\{0,1,a,a+1\,\Big|\,a^2=a+1\right\}$$

To be more precise here, we assume of course that the characteristic of our model is $p=2$, which reads $x+x=0$ for any $x$, and so determines the addition table. As for the multiplication table, this is uniquely determined by $a^2=-a-1=a+1$.

\medskip

(5) Next, at exponents of type $q=p^2$ with $p\geq3$ prime, we can use $Q=X^2-r$, with $r$ being a non-square modulo $p$, and with $(p-1)/2$ choices here. We are led to:
$$\mathbb F_{p^2}=\left\{a+b\gamma\,\Big|\,\gamma^2=r\right\}$$

Here, as before with $\mathbb F_4$, our formula is something self-explanatory. Observe the analogy with $\mathbb C=\mathbb R[i]$, with $i$ being a formal number satisfying $i^2=-1$.

\medskip

(6) Finally, at $q=p^k$ with $k\geq3$ things become more complicated, but the main idea remains the same. We have for instance models for $\mathbb F_8$, $\mathbb F_{27}$ using $Q=X^3-X-1$, and a model for $\mathbb F_{16}$ using $Q=X^4+X+1$. Many other things can be said here.
\end{proof}

As another application of the above, which motivated Galois, we have:

\index{degree 5 polynomial}
\index{Galois theory}
\index{solvable group}
\index{tower of extensions}
\index{field extension}
\index{separable extension}
\index{Abel-Ruffini}

\begin{theorem}
Unlike in degree $N\leq4$, there is no formula for the roots of polynomials of degree $N=5$ and higher, with the reason for this, coming from Galois theory, being that $S_5$ is not solvable. The simplest numeric example is $P=X^5-X-1$.
\end{theorem}

\begin{proof}
This is something quite tricky, the idea being as follows:

\medskip

(1) The first assertion, for generic polynomials, is due to Abel-Ruffini, but Galois theory helps in better understanding this, and comes with a number of bonus points too, namely the possibility of formulating a finer result, with Abel-Ruffini's original ``generic'', which was something algebraic, being now replaced by an analytic ``generic'', and also with the possibility of dealing with concrete polynomials, such as $P=X^5-X-1$.

\medskip

(2) Regarding now the details of the Galois proof of the Abel-Ruffini theorem, assume that the roots of a polynomial $P\in F[X]$ can be computed by using iterated roots, a bit as for the degree 2 equation, or for the degree 3 and 4 equations, via Cardano. Then, algebrically speaking, this gives rise to a tower of fields as folows, with $F_0=F$, and each $F_{i+1}$ being obtained from $F_i$ by adding a root, $F_{i+1}=F_i(x_i)$, with $x_i^{n_i}\in F_i$:
$$F_0\subset F_1\subset\ldots\subset F_k$$

(3) In order for Galois theory to apply well to this situation, we must make all the extensions normal, which amounts in replacing each $F_{i+1}=F_i(x_i)$ by its extension $K_i(x_i)$, with $K_i$ extending $F_i$ by adding a $n_i$-th root of unity. Thus, with this replacement, we can assume that the tower in (2) in normal, meaning that all Galois groups are cyclic.

\medskip

(4) Now by Galois theory, at the level of the corresponding Galois groups we obtain a tower of groups as follows as follows, which is a resolution of the last group $G_k$, the Galois group of $P$, in the sense of group theory, in the sense that all quotients are cyclic:
$$G_1\subset G_2\subset\ldots\subset G_k$$

As a conclusion, Galois theory tells us that if the roots of a polynomial $P\in F[X]$ can be computed by using iterated roots, then its Galois group $G=G_k$ must be solvable.

\medskip

(5) In the generic case, the conclusion is that Galois theory tells us that, in order for all polynomials of degree 5 to be solvable, via square roots, the group $S_5$, which appears there as Galois group, must be solvable, in the sense of group theory. But this is wrong, because the alternating subgroup $A_5\subset S_5$ is simple, and therefore not solvable.

\medskip

(6) Finally, regarding the polynomial $P=X^5-X-1$, some elementary computations here, based on arithmetic over $\mathbb F_2,\mathbb F_3$, and involving various cycles of length $2,3,5$, show that its Galois group is $S_5$. Thus, we have our counterexample.

\medskip

(7) To be more precise, our polynomial factorizes over $\mathbb F_2$ as follows:
$$X^5-X-1=(X^2+X+1)(X^3+X^2+1)$$

We deduce from this the existence of an element $\tau\sigma\in G\subset S_5$, with $\tau\in S_5$ being a transposition, and with $\sigma\in S_5$ being a 3-cycle, disjoint from it. Thus, we have:
$$\tau=(\tau\sigma)^3\in G$$

(8) On the other hand since $P=X^5-X-1$ is irreducible over $\mathbb F_5$, we have as well available a certain 5-cycle $\rho\in G$. Now since $<\tau,\rho>=S_5$, we conclude that the Galois group of $P$ is full, $G=S_5$, and by (4) and (5) we have our counterexample.

\medskip

(9) Finally, as mentioned in (1), all this shows as well that a random polynomial of degree 5 or higher is not solvable by square roots, and with this being an elementary consequence of the main result from (5), via some standard analysis arguments.
\end{proof}

So long for Galois theory, and its basic applications, quicky explained. For more on all this, we recommend any solid abstract algebra book, such as Lang \cite{la1}.

\section*{8e. Exercises}

This was a quite technical algebra chapter, and as exercises here, we have:

\begin{exercise}
Learn more about fields having characteristic $0$.
\end{exercise}

\begin{exercise}
Learn as well about fields which are algebrically closed.
\end{exercise}

\begin{exercise}
Learn more about quadratic extensions.
\end{exercise}

\begin{exercise}
Learn also about cyclotomic fields.
\end{exercise}

\begin{exercise}
Learn more about rings, modules and ideals.
\end{exercise}

\begin{exercise}
Learn more about algebras, and their various properties.
\end{exercise}

\begin{exercise}
Learn more about finite fields, and their modeling.
\end{exercise}

\begin{exercise}
Learn more about higher degree polynomials, and their roots.
\end{exercise}

As bonus exercise, and no surprise here, systematically learn Galois theory.

\part{Representations}

\ \vskip50mm

\begin{center}
{\em It was dark all around, there was frost in the ground

When the tigers broke free

And no one survived

From the Royal Fusiliers Company Z}
\end{center}

\chapter{Representations}

\section*{9a. Representations}

We have seen so far many interesting examples of finite groups, often appearing as permutation or rotation groups, $G\subset S_N$ or $G\subset O_N$ or $G\subset U_N$, and techniques for dealing with them. However, while in what regards the examples our discussion was reasonably complete, at least at an introductory level, concerning the techniques, there are still many things to be said, of quite fundamental nature. In fact, and you might probably find this a bit odd, coming at this place of this book, well into our study, page 200, the basic modern techniques for studying the finite groups are yet to come.

\bigskip

We will discuss these techniques here, in the present Part III of this book, following Hermann Weyl and others. The idea is something very simple, as follows:

\begin{principle}
A finite group $G$, coming either abstractly, or as a concrete group, $G\subset S_N$ or $G\subset O_N$ or $G\subset U_N$, is best studied via its representations
$$u:G\to U_N$$ 
and more specifically via its irreducible representations, with the functions $f:G\to\mathbb C$ decomposing as sums of coefficients of these. Also, in the concrete group situation,
$$G\subset [S_N\subset O_N\subset U_N]$$
it all comes down to splitting the tensor powers of $G\subset U_N$ into irreducible representations.
\end{principle}

Well, so this is the general principle, basically a visionary new idea, due to Weyl, and if you don't understand, no worries, nor did I when first seeing such things, back in the days, as a student. Actually, in case the above principle looks totally opaque to you, here is a softer version of it, making some links with things that we already know:

\begin{principle}[soft version]
In order to study an abstract finite group $G$, we can make it concrete, by looking at its representations
$$u:G\to U_N$$
with each such representation containing a piece of useful information about $G$. And with a bit of luck, all these pieces will provide us with a good picture of $G$.
\end{principle}

But probably too much talking, let us get to work, and we will understand later what exactly we are doing. To start with, let us formulate the following definition:

\index{representation}
\index{character}
\index{finite group}
\index{character of representation}
\index{trace of representation}

\begin{definition}
A representation of a finite group $G$ is a morphism as follows:
$$u:G\to U_N$$
The character of such a representation is the function $\chi:G\to\mathbb C$ given by
$$g\to Tr(u_g)$$
where $Tr$ is the usual trace of the $N\times N$ matrices, $Tr(M)=\sum_iM_{ii}$.
\end{definition}

As a basic example here, for any finite group we always have available the trivial 1-dimensional representation, or character, which is by definition as follows:
$$u:G\to U_1\quad,\quad 
g\to(1)$$

In fact, talking 1-dimensional representations, we already know about these, from chapter 4, with the summary of our results there being as follows:

\begin{proposition}
The $1$-dimensional representations of $G$ are the morphisms
$$u:G\to\mathbb T$$
and we have $u=\chi$ in this case. These morphisms, or characters, must come via
$$u:G\to G_{ab}\to\mathbb T$$
from the characters $G_{ab}\to\mathbb T$, which themselves form a group, isomorphic to $G_{ab}$.
\end{proposition}

\begin{proof}
This is indeed something self-explanatory, coming from our discussion from chapter 4, and with the non-trivial part being the last assertion, which comes from the structure theorem for the finite abelian groups, such as $G_{ab}$, as explained there.
\end{proof}

Moving now to higher dimensions, as another class of basic examples, when our finite group $G$ appears by definition, or by theorem, a group of unitary matrices, $G\subset U_N$, the embedding $G\subset U_N$ itself is a representation, called fundamental one:
$$u:G\subset U_N\quad,\quad 
g\to g$$

And here, again, we are in familiar territory, because we can say the following, coming from our knowledge of the Cayley theorem and related topics, from chapter 3:

\begin{proposition}
Available for any finite group $G$ is its regular representation
$$u:G\subset S_N\subset O_N\subset U_N$$
with $N=|G|$, obtained via Cayley and permutation matrices, the formula being
$$u_g(e_h)=e_{gh}$$
with $\{e_h|h\in G\}$ being the standard basis of $\mathbb C^N$. Its character is $\chi(g)=N\delta_{g1}$.
\end{proposition}

\begin{proof}
This is again something self-explanatory, coming from our discussion from chapter 3, on the Cayley theorem, permutation matrices and related topics, and with the character computation being something elementary too, as follows:
\begin{eqnarray*}
\chi(g)
&=&Tr(u_g)\\
&=&\sum_{h\in G}<u_g(e_h),e_h>\\
&=&\sum_{h\in G}<e_{gh},e_h>\\
&=&N\delta_{g1}
\end{eqnarray*}

Thus, we are led to the conclusions in the statement.
\end{proof}

Summarizing, we definitely have interesting illustrations for Definition 9.3, and even some beginning of theory on the way, based on our material from chapters 3-4.

\bigskip

What is next? You guessed it right, more examples. Inspired by the above, let us formulate the following question, which looks like something quite interesting:

\begin{question}
Given a subgroup $G\subset U_N$, besides its fundamental representation
$$u:G\subset U_N\quad,\quad g\to g$$
we can equally talk about its conjugate fundamental representation
$$\bar{u}:G\subset U_N\quad,\quad g\to\bar{g}$$
and probably about many more, coming via other operations. What exactly are these?
\end{question} 

To be more precise here, consider the usual conjugation of the unitary matrices, $(\bar{U})_{ij}=\bar{U}_{ij}$. This can be viewed as a group isomorphism, as follows:
$$U_N\simeq U_N\quad,\quad U\to\bar{U}$$

Now given an embedding $u:G\subset U_N$, we can compose it with this isomorphism $U_N\simeq U_N$, and we obtain another embedding $\bar{u}:G\subset U_N$. And with $\bar{u}$ being in general different from $u$ itself, as the 1D examples, in the context of Proposition 9.4, show.

\bigskip

Moreover, and here comes the point, this construction $u\to\bar{u}$ is probably just the tip of the iceberg, because assuming that we have constructed other operations $*$ on the representations, we can associate to $u$ a whole family of representations, as follows:
$$u_{a_1b_1\ldots a_kb_k}=u^{*a_1}*\bar{u}^{*b_1}*\ldots*u^{*a_k}*\bar{u}^{*b_k}$$

In order to clarify all this, and understand which representations are available for a subgroup $G\subset U_N$, or more generally for an already represented group, $G\to U_N$, let us first discuss the various operations on the representations, in the general context of Definition 9.3. We have here the following result, which is something elementary:

\index{sum of representations}
\index{product of representations}
\index{conjugate representation}
\index{spinned representation}

\begin{theorem}
The representations of a finite group $G$ are subject to:
\begin{enumerate}
\item Making sums. Given representations $u,v$, having dimensions $N,M$, their sum is the $N+M$-dimensional representation $u+v=diag(u,v)$.

\item Making products. Given representations $u,v$, having dimensions $N,M$, their tensor product is the $NM$-dimensional representation $(u\otimes v)_{ia,jb}=u_{ij}v_{ab}$.

\item Taking conjugates. Given a representation $u$, having dimension $N$, its complex conjugate is the $N$-dimensional representation $(\bar{u})_{ij}=\bar{u}_{ij}$.

\item Spinning by unitaries. Given a representation $u$, having dimension $N$, and a unitary $V\in U_N$, we can spin $u$ by this unitary, $u\to VuV^*$.
\end{enumerate}
\end{theorem}

\begin{proof}
The fact that the operations in the statement are indeed well-defined, among maps from $G$ to unitary groups, can be checked as follows:

\medskip

(1) This follows from the trivial fact that if $g\in U_N$ and $h\in U_M$ are two unitaries, then their diagonal sum is a unitary too, as follows:
$$\begin{pmatrix}g&0\\ 0&h\end{pmatrix}\in U_{N+M}$$

(2) This follows from the fact that if $g\in U_N$ and $h\in U_M$ are two unitaries, then $g\otimes h\in U_{NM}$ is a unitary too. Given unitaries $g,h$, let us set indeed:
$$(g\otimes h)_{ia,jb}=g_{ij}h_{ab}$$

This matrix is then a unitary too, as shown by the following computation:
\begin{eqnarray*}
[(g\otimes h)(g\otimes h)^*]_{ia,jb}
&=&\sum_{kc}(g\otimes h)_{ia,kc}((g\otimes h)^*)_{kc,jb}\\
&=&\sum_{kc}(g\otimes h)_{ia,kc}\overline{(g\otimes h)_{jb,kc}}\\
&=&\sum_{kc}g_{ik}h_{ac}\bar{g}_{jk}\bar{h}_{bc}\\
&=&\sum_kg_{ik}\bar{g}_{jk}\sum_ch_{ac}\bar{h}_{bc}\\
&=&\delta_{ij}\delta_{ab}
\end{eqnarray*}

(3) This simply follows from the fact that if $g\in U_N$ is unitary, then so is its complex conjugate, $\bar{g}\in U_N$, and this due to the following formula, obtained by conjugating:
$$g^*=g^{-1}\implies g^t=\bar{g}^{-1}$$

(4) This is clear as well, because if $g\in U_N$ is unitary, and $V\in U_N$ is another unitary, then we can spin $g$ by this unitary, and we obtain a unitary as follows:
$$VgV^*\in U_N$$

Thus, our operations are well-defined, and this leads to the above conclusions.
\end{proof}

Before going further, let us record as well a result about the corresponding characters. The formulae here, which are all straightforward, are as follows:

\index{character}

\begin{proposition}
We have the following formulae, regarding characters,
$$\chi_{u+v}=\chi_u+\chi_v\quad,\quad 
\chi_{u\otimes v}=\chi_u\chi_v$$ 
$$\chi_{\bar{u}}=\bar{\chi}_u\quad,\quad 
\chi_{VuV^*}=\chi_u$$
in relation with the basic operations for the representations.
\end{proposition}

\begin{proof}
All these assertions are elementary, by using the following well-known trace formulae, valid for any two square matrices $g,h$, and any unitary $V$:
$$Tr(diag(g,h))=Tr(g)+Tr(h)\quad,\quad 
Tr(g\otimes h)=Tr(g)Tr(h)$$
$$Tr(\bar{g})=\overline{Tr(g)}\quad,\quad 
Tr(VgV^*)=Tr(g)$$

To be more precise, the first formula is clear from definitions. Regarding now the second formula, the computation here is immediate too, as follows:
\begin{eqnarray*}
Tr(g\otimes h)
&=&\sum_{ia}(g\otimes h)_{ia,ia}\\
&=&\sum_{ia}g_{ii}h_{aa}\\
&=&\sum_ig_{ii}\sum_ah_{aa}\\
&=&Tr(g)Tr(h)
\end{eqnarray*}

Regarding now the third formula, this is clear from definitions, by conjugating. Finally, regarding the fourth formula, this can be established as follows:
$$Tr(VgV^*)=Tr(gV^*V)=Tr(g)$$

Thus, we are led to the conclusions in the statement.
\end{proof}

Assume now that we are given a finite group $G\subset U_N$. By using the above operations, we can construct a whole family of representations of $G$, as follows:

\index{Peter-Weyl representations}
\index{colored integer}

\begin{definition}
Given a finite group $G\subset U_N$, its Peter-Weyl representations are the tensor products between the fundamental representation and its conjugate:
$$u:G\subset U_N\quad,\quad 
\bar{u}:G\subset U_N$$ 
We denote these tensor products $u^{\otimes k}$, with $k=\circ\bullet\bullet\circ\ldots$ being a colored integer, with the colored tensor powers being defined according to the rules 
$$u^{\otimes\circ}=u\quad,\quad 
u^{\otimes\bullet}=\bar{u}\quad,\quad
u^{\otimes kl}=u^{\otimes k}\otimes u^{\otimes l}$$
and with the convention that $u^{\otimes\emptyset}$ is the trivial representation $1:G\to U_1$.
\end{definition}

Here are a few examples of such Peter-Weyl representations, namely those coming from the colored integers of length 2, to be often used in what follows:
$$u^{\otimes\circ\circ}=u\otimes u\quad,\quad 
u^{\otimes\circ\bullet}=u\otimes\bar{u}$$
$$u^{\otimes\bullet\circ}=\bar{u}\otimes u\quad,\quad 
u^{\otimes\bullet\bullet}=\bar{u}\otimes\bar{u}$$

In relation now with characters, we have the following result:

\index{colored powers}

\begin{theorem}
The characters of Peter-Weyl representations are given by
$$\chi_{u^{\otimes k}}=(\chi_u)^k$$
with the colored powers of a variable $\chi$ being by definition given by
$$\chi^\circ=\chi\quad,\quad
\chi^\bullet=\bar{\chi}\quad,\quad
\chi^{kl}=\chi^k\chi^l$$
and with the convention that $\chi^\emptyset$ equals by definition $1$.
\end{theorem}

\begin{proof}
This follows indeed from the additivity, multiplicativity and conjugation formulae established in Proposition 9.8, via the conventions in Definition 9.9.
\end{proof}

The above result, although being something rather trivial, is quite interesting, because it makes a connection with our analytic work for the group $S_N$ from chapter 3. To be more precise, with that work in mind, we can formulate the following question:

\begin{question}
With the character of $u:S_N\subset O_N$ counting the fixed points,
$$\chi(\sigma)=\#\left\{i\in\{1,\ldots,N\}\Big|\sigma(i)=i\right\}$$
our previous counting results tell us that this character, regarded as random variable
$$\chi:S_N\to\mathbb N$$
follows the Poisson law $p_1$, with $N\to\infty$. What about other subgroups $G\subset U_N$?
\end{question}

And good question this is, providing us with an excellent motivation for the study of the variables in Theorem 9.10, the point being that, for an arbitrary subgroup $G\subset U_N$, the moments of the main character $\chi_u:G\to\mathbb C$, and hopefully the law of this main character too, can be computed by averaging over $G$ that variables. 

\bigskip

But more on this later in this book, on several occasions. In fact, for answering the above question for various other subgroups $G\subset U_N$, the quite ad-hoc counting techniques that we used in chapter 3 for the symmetric group $S_N$ will obviously not work, and we are precisely in need, as a main tool, of more advanced representation theory.

\bigskip

In any case, as a conclusion to all this, and to our previous abstract considerations too, given a closed subgroup $G\subset U_N$, we would like to understand its Peter-Weyl representations. In order to do so, let us formulate the following key definition:

\index{Hom space}
\index{End space}
\index{Fix space}

\begin{definition}
Given a finite group $G$, and two of its representations,
$$u:G\to U_N\quad,\quad 
v:G\to U_M$$
we define the linear space of intertwiners between these representations as being 
$$Hom(u,v)=\left\{T\in M_{M\times N}(\mathbb C)\Big|Tu_g=v_gT,\forall g\in G\right\}$$
and we use the following conventions:
\begin{enumerate}
\item We use the notations $Fix(u)=Hom(1,u)$, and $End(u)=Hom(u,u)$.

\item We write $u\sim v$ when $Hom(u,v)$ contains an invertible element.

\item We say that $u$ is irreducible, and write $u\in Irr(G)$, when $End(u)=\mathbb C1$.
\end{enumerate}
\end{definition}

In the above the terminology is very standard, with Hom and End standing respectively for ``homomorphisms'' and ``endomorphisms'', and with Fix standing for ``fixed points''. In practice, it is useful to think of the representations of $G$ as being the objects of some kind of abstract combinatorial structure associated to $G$, and of the intertwiners between these representations as being the ``arrows'' between these objects. 

\bigskip

We have in fact the following result about this, that you will certainly appreciate if you are a bit familiar with category theory, which clarifies the situation:

\index{tensor category}

\begin{theorem}
The following happen, regarding the above notions:
\begin{enumerate}
\item The intertwiners are stable under composition:
$$T\in Hom(u,v)\ ,\ 
S\in Hom(v,w)
\implies ST\in Hom(u,w)$$

\item The intertwiners are stable under taking tensor products:
$$S\in Hom(u,v)\ ,\ 
T\in Hom(w,t)\\
\implies S\otimes T\in Hom(u\otimes w,v\otimes t)$$

\item The intertwiners are stable under taking adjoints:
$$T\in Hom(u,v)
\implies T^*\in Hom(v,u)$$

\item Thus, the Hom spaces form a tensor $*$-category.
\end{enumerate}
\end{theorem}

\begin{proof}
All this is clear from definitions, the verifications being as follows:

\medskip

(1) This follows indeed from the following computation, valid for any $g\in G$:
$$STu_g=Sv_gT=w_gST$$

(2) Again, this is clear, because we have the following computation:
\begin{eqnarray*}
(S\otimes T)(u_g\otimes w_g)
&=&Su_g\otimes Tw_g\\
&=&v_gS\otimes t_gT\\
&=&(v_g\otimes t_g)(S\otimes T)
\end{eqnarray*}

(3) This follows from the following computation, valid for any $g\in G$:
\begin{eqnarray*}
Tu_g=v_gT
&\implies&u_g^*T^*=T^*v_g^*\\
&\implies&T^*v_g=u_gT^*
\end{eqnarray*}

(4) This is just an abstract conclusion of (1,2,3), with a tensor $*$-category being by definition an abstract beast satisfying these conditions (1,2,3). We will be back to tensor categories later on in this book, with more details on all this.
\end{proof}

And with this, end of our preliminary discussion regarding the representations. With respect to the guiding principles from the beginning of this chapter, we can say that Principle 9.2 is reasonably understood. So, time now to get to the more specialized Principle 9.1, and see how that splitting technology for representations and their coefficients works. We will explain this, called Peter-Weyl theory, in the remainder of this chapter.

\section*{9b. Peter-Weyl}

Our claim now is that Theorem 9.13 gives us everything that we need, in order to have some advanced representation theory started, for our finite groups $G$. Indeed, as a main consequence of Theorem 9.13, we have the following key result:

\begin{theorem}
Given a representation $u:G\to U_N$, the linear space
$$End(u)\subset M_N(\mathbb C)$$
is a $*$-algebra, with respect to the usual involution of the matrices.
\end{theorem}

\begin{proof}
We know from Theorem 9.13 (1) that $End(u)$ is a subalgebra of $M_N(\mathbb C)$, and we know as well from Theorem 9.13 (3) that this subalgebra is stable under the involution $*$. Thus, what we have here is a $*$-subalgebra of $M_N(\mathbb C)$, as claimed.
\end{proof}

The point now is that we can combine the above result with the following standard fact, that you might know or not, from advanced linear algebra:

\index{operator algebra}
\index{finite dimensional algebra}

\begin{theorem}
Let $A\subset M_N(\mathbb C)$ be a $*$-algebra.
\begin{enumerate}
\item We can write $1=p_1+\ldots+p_k$, with $p_i\in A$ being central minimal projections.

\item The linear spaces $A_i=p_iAp_i$ are non-unital $*$-subalgebras of $A$.

\item We have a non-unital $*$-algebra sum decomposition $A=A_1\oplus\ldots\oplus A_k$.

\item We have unital $*$-algebra isomorphisms $A_i\simeq M_{n_i}(\mathbb C)$, with $n_i=rank(p_i)$.

\item Thus, we have a $*$-algebra isomorphism $A\simeq M_{n_1}(\mathbb C)\oplus\ldots\oplus M_{n_k}(\mathbb C)$.
\end{enumerate}
\end{theorem}

\begin{proof}
This is something very standard. Consider indeed an arbitrary $*$-algebra of the $N\times N$ matrices, $A\subset M_N(\mathbb C)$. Let us first look at the center of this algebra, $Z(A)=A\cap A'$. This center, viewed as an algebra, is then of the following form:
$$Z(A)\simeq\mathbb C^k$$

Consider now the standard basis $e_1,\ldots,e_k\in\mathbb C^k$, and let  $p_1,\ldots,p_k\in Z(A)$ be the images of these vectors via the above identification. In other words, these elements $p_1,\ldots,p_k\in A$ are central minimal projections, summing up to 1:
$$p_1+\ldots+p_k=1$$

The idea is then that this partition of the unity will eventually lead to the block decomposition of $A$, as in the statement. We prove this in 4 steps, as follows:

\medskip

\underline{Step 1}. We first construct the matrix blocks, our claim here being that each of the following linear subspaces of $A$ are non-unital $*$-subalgebras of $A$:
$$A_i=p_iAp_i$$

But this is clear, with the fact that each $A_i$ is closed under the various non-unital $*$-subalgebra operations coming from the projection equations $p_i^2=p_i^*=p_i$.

\medskip

\underline{Step 2}. We prove now that the above algebras $A_i\subset A$ are in a direct sum position, in the sense that we have a non-unital $*$-algebra sum decomposition, as follows:
$$A=A_1\oplus\ldots\oplus A_k$$

As with any direct sum question, we have two things to be proved here. First, by using the formula $p_1+\ldots+p_k=1$ and the projection equations $p_i^2=p_i^*=p_i$, we conclude that we have the needed generation property, namely:
$$A_1+\ldots+ A_k=A$$

As for the fact that the sum is indeed direct, this follows as well from the formula $p_1+\ldots+p_k=1$, and from the projection equations $p_i^2=p_i^*=p_i$.

\medskip

\underline{Step 3}. Our claim now, which will finish the proof, is that each of the $*$-subalgebras $A_i=p_iAp_i$ constructed above is in fact a full matrix algebra. To be more precise, with $n_i=rank(p_i)$, our claim is that we have isomorphisms, as follows:
$$A_i\simeq M_{n_i}(\mathbb C)$$

In order to prove this claim, recall that the projections $p_i\in A$ were chosen central and minimal. Thus, the center of each of the algebras $A_i$ reduces to the scalars:
$$Z(A_i)=\mathbb C$$

But this shows, either via a direct computation, or via the bicommutant theorem, that the each of the algebras $A_i$ is a full matrix algebra, as claimed.

\medskip

\underline{Step 4}. We can now obtain the result, by putting together what we have. Indeed, by using the results from Step 2 and Step 3, we obtain an isomorphism as follows:
$$A\simeq M_{n_1}(\mathbb C)\oplus\ldots\oplus M_{n_k}(\mathbb C)$$

In addition to this, a careful look at the isomorphisms established in Step 3 shows that at the global level, of the algebra $A$ itself, the above isomorphism simply comes by twisting the following standard multimatrix embedding, discussed in the beginning of the proof, (1) above, by a certain unitary matrix $U\in U_N$:
$$M_{n_1}(\mathbb C)\oplus\ldots\oplus M_{n_k}(\mathbb C)\subset M_N(\mathbb C)$$

Now by putting everything together, we obtain the result.
\end{proof}

Good news, we can now formulate our first Peter-Weyl theorem, as follows:

\index{Peter-Weyl}

\begin{theorem}[PW1]
Let $u:G\to U_N$ be a representation, consider the algebra $A=End(u)$, and write its unit as above, with $p_i$ being central minimal projections:
$$1=p_1+\ldots+p_k$$
The representation $u$ decomposes then as a direct sum, as follows,
$$u=u_1+\ldots+u_k$$
with each $u_i$ being an irreducible representation, obtained by restricting $u$ to $Im(p_i)$.
\end{theorem}

\begin{proof}
This basically follows from Theorem 9.14 and Theorem 9.15, as follows:

\medskip

(1) As a first observation, by replacing $G$ with its image $u(G)\subset U_N$, we can assume if we want that our representation $u$ is faithful, $G\subset_uU_N$. However, this replacement will not be really needed, and we will keep using $u:G\to U_N$, as above.

\medskip

(2) In order to prove the result, we will need some preliminaries. We first associate to our representation $u:G\to U_N$ the corresponding action map on $\mathbb C^N$. If a linear subspace $V\subset\mathbb C^N$ is invariant, the restriction of the action map to $V$ is an action map too, which must come from a subrepresentation $v\subset u$. This is clear indeed from definitions, and with the remark that the unitaries, being isometries, restrict indeed into unitaries.

\medskip

(3) Consider now a projection $p\in End(u)$. From $pu=up$ we obtain that the linear space $V=Im(p)$ is invariant under $u$, and so this space must come from a subrepresentation $v\subset u$. It is routine to check that the operation $p\to v$ maps subprojections to subrepresentations, and minimal projections to irreducible representations.

\medskip

(4) To be more precise here, the condition $p\in End(u)$ reformulates as follows:
$$pu_g=u_gp\quad,\quad\forall g\in G$$

As for the condition that $V=Im(p)$ is invariant, this reformulates as follows:
$$pu_gp=u_gp\quad,\quad\forall g\in G$$

Thus, we are in need of a technical linear algebra result, stating that for a projection $P\in M_N(\mathbb C)$ and a unitary $U\in U_N$, the following happens:
$$PUP=UP\implies PU=UP$$

(5) But this can be established with some $*$-algebra know-how, as follows:
\begin{eqnarray*}
tr[(PU-UP)(PU-UP)^*]
&=&tr[(PU-UP)(U^*P-PU^*)]\\
&=&tr[P-PUPU^*-UPU^*P+UPU^*]\\
&=&tr[P-UPU^*-UPU^*+UPU^*]\\
&=&tr[P-UPU^*]\\
&=&0
\end{eqnarray*}

Indeed, by positivity this gives $PU-UP=0$, as desired.

\medskip

(6) With these preliminaries in hand, let us decompose the algebra $End(u)$ as in Theorem 9.15, by using the decomposition $1=p_1+\ldots+p_k$ into minimal projections. If we denote by $u_i\subset u$ the subrepresentation coming from the vector space $V_i=Im(p_i)$, then we obtain in this way a decomposition $u=u_1+\ldots+u_k$, as in the statement.
\end{proof}

In order to formulate our second Peter-Weyl theorem, let us formulate:

\index{space of coefficients}

\begin{definition}
Given a finite subgroup $G\subset U_N$, and a unitary representation $v:G\to U_M$, the space of coefficients of this representation is:
$$C_v=\left\{f\circ v\Big|f\in M_M(\mathbb C)^*\right\}$$
In other words, by delinearizing, $C_\nu\subset C(G)$ is the following linear space,
$$C_v=span\Big[g\to(v_g)_{ij}\Big]$$
with $g\to(v_g)_{ij}$ being the standard matrix coefficients of $v:G\to U_M$.
\end{definition}

As a basic example of coefficient we have, besides the matrix coefficients $g\to(v_g)_{ij}$, the character, which appears as the diagonal sum of these coefficients:
$$\chi_v(g)=\sum_i(v_g)_{ii}$$

Here is now our second Peter-Weyl theorem, complementing Theorem 9.16:

\index{Peter-Weyl}

\begin{theorem}[PW2]
Given a subgroup $G\subset_uU_N$, any irreducible representation 
$$v:G\to U_M$$
appears inside a tensor product of the fundamental representation $u$ and its adjoint $\bar{u}$.
\end{theorem}

\begin{proof}
In order to prove the result, we will use the following three elementary facts, regarding the spaces of coefficients introduced above:

\medskip

(1) The construction $v\to C_v$ is functorial, in the sense that it maps subrepresentations into linear subspaces. This is indeed something which is routine to check.

\medskip

(2) By the Stone-Weierstrass theorem, which tells us that we have $<g_{ij}>=C(G)$, we conclude that have an inclusion of linear spaces as follows:
$$C_v\subset<g_{ij}>$$

(3) By definition of the Peter-Weyl representations, as arbitrary tensor products between the fundamental representation $u$ and its conjugate $\bar{u}$, we have:
$$<g_{ij}>=\sum_kC_{u^{\otimes k}}$$

(4) Now by putting together the observations (2,3) we conclude that we must have an inclusion of linear spaces as follows, for certain exponents $k_1,\ldots,k_p$:
$$C_v\subset C_{u^{\otimes k_1}\oplus\ldots\oplus\pi^{\otimes k_p}}$$

By using now the functoriality result from (1), we deduce from this inclusion of coefficient spaces that we have an inclusion of representations, as follows:
$$v\subset u^{\otimes k_1}\oplus\ldots\oplus u^{\otimes k_p}$$

Together with Theorem 9.16, this leads to the conclusion in the statement.
\end{proof}

As a conclusion to what we have so far, the problem to be solved is that of splitting the Peter-Weyl representations into sums of irreducible representations.

\section*{9c. Averaging results}

In order to further advance, and complete the Peter-Weyl theory, we need to talk about averaging over $G$. In the present finite group case the situation is trivial, as follows:

\begin{proposition}
Any finite group $G$ has a unique probability measure which is invariant under left and right translations,
$$\mu(E)=\mu(gE)=\mu(Eg)$$
and this is the normalized counting measure on $G$, given by $\mu(E)=|E|/|G|$.
\end{proposition}

\begin{proof}
The uniformity condition in the statement gives, with $E=\{h\}$:
$$\mu\{h\}=\mu\{gh\}=\mu\{hg\}$$

Thus $\mu$ must be the usual counting measure, normalized as to have mass 1.
\end{proof}

For our purposes here, we need to know more about averaging over $G$, in relation with the coefficients of representations. Also, in view of some further generalizations later, to the case of the compact groups, it is convenient to use a functional analytic approach, to such questions. More specifically, it is convenient to work with the integration functionals with respect to the various measures on $G$, instead of the measures themselves. 

\bigskip

With these conventions, which are quite standard, we first have:

\begin{proposition}
Given a unital positive linear form $\varphi:C(G)\to\mathbb C$, the limit
$$\int_\varphi f=\lim_{n\to\infty}\frac{1}{n}\sum_{k=1}^n\varphi^{*k}(f)$$
exists, and for a coefficient of a representation $f=(\tau\otimes id)v$ we have
$$\int_\varphi f=\tau(P)$$
where $P$ is the orthogonal projection onto the $1$-eigenspace of $(id\otimes\varphi)v$.
\end{proposition}

\begin{proof}
Given a finite subgroup $u:G\subset U_N$, the coordinate functions $u_{ij}:G\to\mathbb C$ obviously separate the points $g\in G$, so by Stone-Weierstrass we have:
$$C(G)=<u_{ij}>$$

We conclude, by linearity, that it is enough to prove the first assertion for functions of the following type, where $v$ is a Peter-Weyl representation, and $\tau$ is a linear form:
$$f=(\tau\otimes id)v$$

Thus we are led into the second assertion, and more precisely, we can have the whole result proved if we can establish the following formula, with $f=(\tau\otimes id)v$:
$$\lim_{n\to\infty}\frac{1}{n}\sum_{k=1}^n\varphi^{*k}(f)=\tau(P)$$

In order to prove this latter formula, observe that we have:
$$\varphi^{*k}(f)
=(\tau\otimes\varphi^{*k})v
=\tau((id\otimes\varphi^{*k})v)$$

Now let us set $M=(id\otimes\varphi)v$. In terms of this matrix, we have:
\begin{eqnarray*}
((id\otimes\varphi^{*k})v)_{i_0i_{k+1}}
&=&\sum_{i_1\ldots i_k}((id\otimes\varphi)v)_{i_0i_1}\ldots((id\otimes\varphi)v)_{i_ki_{k+1}}\\
&=&\sum_{i_1\ldots i_k}M_{i_0i_1}\ldots M_{i_ki_{k+1}}\\
&=&(M^k)_{i_0i_{k+1}}
\end{eqnarray*}

We conclude that we have the following formula, valid for any $k\in\mathbb N$:
$$(id\otimes\varphi^{*k})v=M^k$$

It follows that our Ces\`aro limit is given by the following formula:
\begin{eqnarray*}
\lim_{n\to\infty}\frac{1}{n}\sum_{k=1}^n\varphi^{*k}(f)
&=&\lim_{n\to\infty}\frac{1}{n}\sum_{k=1}^n\tau(M^k)\\
&=&\tau\left(\lim_{n\to\infty}\frac{1}{n}\sum_{k=1}^nM^k\right)
\end{eqnarray*}

Now since $v$ is unitary we have $||v||=1$, and from $M=(id\otimes\varphi)v$ we obtain:
$$||M||\leq1$$

We conclude that the last Ces\`aro limit appearing above converges, and equals the orthogonal projection onto the $1$-eigenspace of $M$. That is, we conclude that we have:
$$\lim_{n\to\infty}\frac{1}{n}\sum_{k=1}^nM^k=P$$

Thus our initial Ces\`aro limit converges as well, to the quantity $\tau(P)$, as desired. 
\end{proof}

The point now is that when the linear form $\varphi\in C(G)^*$ from the above result is chosen faithful, we obtain the following finer result, which is something very useful:

\begin{proposition}
Given a faithful unital linear form $\varphi\in C(G)^*$, the limit
$$\int_\varphi f=\lim_{n\to\infty}\frac{1}{n}\sum_{k=1}^n\varphi^{*k}(f)$$
exists, and is independent of $\varphi$, given on coefficients of representations by
$$\left(id\otimes\int_\varphi\right)v=P$$
where $P$ is the orthogonal projection onto the space $Fix(v)=\left\{\xi\in\mathbb C^n\big|v\xi=\xi\right\}$.
\end{proposition}

\begin{proof}
In view of Proposition 9.20, it remains to prove that when $\varphi$ is faithful, the $1$-eigenspace of the matrix $M=(id\otimes\varphi)v$ equals the space $Fix(v)$.

\medskip

``$\supset$'' This is clear, and for any $\varphi$, because we have the following implication:
$$v\xi=\xi\implies M\xi=\xi$$

``$\subset$'' Here we must prove that, when $\varphi$ is faithful, we have:
$$M\xi=\xi\implies v\xi=\xi$$

For this purpose, assume that we have $M\xi=\xi$, and consider the following function:
$$f=\sum_i\left(\sum_jv_{ij}\xi_j-\xi_i\right)\left(\sum_kv_{ik}\xi_k-\xi_i\right)^*$$

We must prove that we have $f=0$. Since $v$ is unitary, we have:
\begin{eqnarray*}
f
&=&\sum_{ijk}v_{ij}v_{ik}^*\xi_j\bar{\xi}_k-\frac{1}{N}v_{ij}\xi_j\bar{\xi}_i-\frac{1}{N}v_{ik}^*\xi_i\bar{\xi}_k+\frac{1}{N^2}\xi_i\bar{\xi}_i\\
&=&\sum_j|\xi_j|^2-\sum_{ij}v_{ij}\xi_j\bar{\xi}_i-\sum_{ik}v_{ik}^*\xi_i\bar{\xi}_k+\sum_i|\xi_i|^2\\
&=&||\xi||^2-<v\xi,\xi>-\overline{<v\xi,\xi>}+||\xi||^2\\
&=&2(||\xi||^2-Re(<v\xi,\xi>))
\end{eqnarray*}

By using now our assumption $M\xi=\xi$, we obtain from this:
\begin{eqnarray*}
\varphi(f)
&=&2\varphi(||\xi||^2-Re(<v\xi,\xi>))\\
&=&2(||\xi||^2-Re(<M\xi,\xi>))\\
&=&2(||\xi||^2-||\xi||^2)\\
&=&0
\end{eqnarray*}

Now since $\varphi$ is faithful, this gives $f=0$, and so $v\xi=\xi$, as claimed.
\end{proof}

We can now formulate a main result about averaging, as follows:

\index{Haar measure}
\index{Haar integration}

\begin{theorem}
The averaging over a finite group $G\subset U_N$ can be constructed by starting with any faithful positive unital linear form $\varphi\in C(G)^*$, and setting:
$$\int_G=\lim_{n\to\infty}\frac{1}{n}\sum_{k=1}^n\varphi^{*k}$$
Moreover, for any representation $v$ we have the following formula,
$$\left(id\otimes\int_G\right)v=P$$
where $P$ is the orthogonal projection onto $Fix(v)=\left\{\xi\in\mathbb C^n\big|v\xi=\xi\right\}$.
\end{theorem}

\begin{proof}
We can prove this from what we have, in several steps, as follows:

\medskip

(1) Let us first go back to the general context of Proposition 9.20. Since convolving one more time with $\varphi$ will not change the Ces\`aro limit appearing there, the functional $\int_\varphi\in C(G)^*$ constructed there has the following invariance property:
$$\int_\varphi*\,\varphi=\varphi*\int_\varphi=\int_\varphi$$

In the case where $\varphi$ is assumed to be faithful, as in Proposition 9.21, our claim is that we have the following formula, valid this time for any $\psi\in C(G)^*$:
$$\int_\varphi*\,\psi=\psi*\int_\varphi=\psi(1)\int_\varphi$$

Moreover, it is enough to prove this formula on a coefficient of a representation:
$$f=(\tau\otimes id)v$$

(2) In order to do so, consider the following two matrices:
$$P=\left(id\otimes\int_\varphi\right)v\quad,\quad 
Q=(id\otimes\psi)v$$

We have then the following two computations, involving these matrices:
$$\left(\int_\varphi*\,\psi\right)f
=\left(\tau\otimes\int_\varphi\otimes\,\psi\right)(v_{12}v_{13})
=\tau(PQ)$$
$$\left(\psi*\int_\varphi\right)f
=\left(\tau\otimes\psi\otimes\int_\varphi\right)(v_{12}v_{13})
=\tau(QP)$$

Also, regarding the term on the right in our formula in (1), this is given by:
$$\psi(1)\int_\varphi f=\psi(1)\tau(P)$$

We conclude from all this that our claim is equivalent to the following equality:
$$PQ=QP=\psi(1)P$$

(3) But this latter equality holds indeed, coming from the fact, that we know from Proposition 9.21, that $P=(id\otimes\int_\varphi)v$ equals the orthogonal projection onto $Fix(v)$. Thus, we have proved our claim in (1), namely that the following formula holds: 
$$\int_\varphi*\,\psi=\psi*\int_\varphi=\psi(1)\int_\varphi$$

(4) In order to finish now, it is convenient to introduce the following abstract operation, on the continuous functions $f,f':C(G)\to\mathbb C$ on our group:
$$\Delta(f\otimes f')(g\otimes h)=f(g)f'(h)$$

With this convention, the formula that we established above can be written as:
$$\psi\left(\int_\varphi\otimes\, id\right)\Delta
=\psi\left(id\otimes\int_\varphi\right)\Delta
=\psi\int_\varphi(.)1$$

This formula being true for any $\psi\in C(G)^*$, we can simply delete $\psi$. We conclude that the following invariance formula holds indeed, with $\int_G=\int_\varphi$:
$$\left(\int_G\otimes\, id\right)\Delta=\left(id\otimes\int_G\right)\Delta=\int_G(.)1$$

But this is exactly the left and right invariance formula we were looking for.

\medskip

(5) Finally, observe that we can recover the uniqueness as well, because assuming that we have two invariant integrals $\int_G,\int_G'$, we have, according to the invariance formula:
$$\left(\int_G\otimes\int_G'\right)\Delta
=\left(\int_G'\otimes\int_G\right)\Delta
=\int_G(.)1
=\int_G'(.)1$$

Thus we have $\int_G=\int_G'$, and this finishes the proof.
\end{proof}

Summarizing, we can now integrate over $G$. As a first application, we have:

\index{main character}
\index{moments}

\begin{theorem}
Given a finite group $G\subset_uU_N$, we have the following formula, valid for any unitary group representation $v:G\to U_M$:
$$\int_G\chi_v=\dim(Fix(v))$$
In particular, the moments of the main character $\chi=\chi_u$ are given by the formula
$$\int_G\chi^k=\dim(Fix(u^{\otimes k}))$$
so knowing the law of $\chi$ is the same as knowing the dimensions on the right.
\end{theorem}

\begin{proof}
We have three assertions here, the idea being as follows:

\medskip

(1) Given a unitary representation $v:G\to U_M$ as in the statement, its character $\chi_v$ is a coefficient, so we can use the integration formula for coefficients in Theorem 9.22. If we denote by $P$ the projection onto $Fix(v)$, that formula gives, as desired:
\begin{eqnarray*}
\int_G\chi_v
&=&Tr(P)\\
&=&\dim(Im(P))\\
&=&dim(Fix(v))
\end{eqnarray*}

(2) This follows from (1), applied to the Peter-Weyl representations, as follows:
\begin{eqnarray*}
\int_G\chi^k
&=&\int_G\chi_u^k\\
&=&\int_G\chi_{u^{\otimes k}}\\
&=&\dim(Fix(u^{\otimes k}))
\end{eqnarray*}

(3) This follows from (2), and from the standard fact, which follows from definitions, that a probability measure is uniquely determined by its moments.
\end{proof}

Observe that the above result makes an interesting link with Question 9.11, raising the possibility of answering that question by using representation theory techniques. We will come back to this later in this book, with some concrete applications.

\bigskip

Along the same lines, but at a more advanced level, let us record as well the following key result, dealing with the computation of the arbitrary integrals over $G$:

\index{Weingarten formula}
\index{Gram matrix}
\index{Weingarten matrix}

\begin{theorem}
The Haar integration over a finite group $G\subset_uU_N$ is given by the Weingarten type formula
$$\int_Gg_{i_1j_1}^{e_1}\ldots g_{i_kj_k}^{e_k}\,dg=\sum_{\pi,\sigma\in D_k}\delta_\pi(i)\delta_\sigma(j)W_k(\pi,\sigma)$$
valid for any colored integer $k=e_1\ldots e_k$ and any multi-indices $i,j$, where $D_k$ is a linear basis of $Fix(u^{\otimes k})$, the associated generalized Kronecker symbols are given by
$$\delta_\pi(i)=<\pi,e_{i_1}\otimes\ldots\otimes e_{i_k}>$$
and $W_k=G_k^{-1}$ is the inverse of the Gram matrix, $G_k(\pi,\sigma)=<\pi,\sigma>$.
\end{theorem}

\begin{proof}
As a first observation, this result computes indeed all the integrals over $G$, and this due to the Stone-Weierstrass theorem, which shows that we have:
$$C(G)=<g_{ij}>$$

Regarding now the proof, we know from Peter-Weyl theory that the integrals in the statement form altogether the orthogonal projection $P^k$ onto the following space:
$$Fix(u^{\otimes k})=span(D_k)$$

Consider now the following linear map, with $D_k=\{\xi_k\}$ being as in the statement:
$$E(x)=\sum_{\pi\in D_k}<x,\xi_\pi>\xi_\pi$$

By a standard linear algebra computation, it follows that we have $P=WE$, where $W$ is the inverse of the restriction of $E$ to the following space:
$$K=span\left(T_\pi\Big|\pi\in D_k\right)$$

But this restriction is precisely the linear map given by the matrix $G_k$, and so $W$ itself is the linear map given by the matrix $W_k$, and this gives the result.
\end{proof}

\section*{9d. More Peter-Weyl}

Getting back now to abstract algebra, in order to further develop the Peter-Weyl theory that we have, we will need, besides the above averaging techniques, an extra ingredient, which is something having its own interest, as follows:

\index{Frobenius isomorphism}

\begin{theorem}
We have a Frobenius type isomorphism
$$Hom(v,w)\simeq Fix(v\otimes\bar{w})$$
valid for any two representations $v,w$.
\end{theorem}

\begin{proof}
According to the definitions, we have the following equivalences:
\begin{eqnarray*}
T\in Hom(v,w)
&\iff&Tv=wT\\
&\iff&\sum_jT_{aj}v_{ji}=\sum_bw_{ab}T_{bi},\forall a,i
\end{eqnarray*}

On the other hand, we have as well the following equivalences:
\begin{eqnarray*}
T\in Fix(v\otimes\bar{w})
&\iff&(v\otimes\bar{w})T=\xi\\
&\iff&\sum_{jb}v_{ij}w_{ab}^*T_{bj}=T_{ai}\forall a,i
\end{eqnarray*}

With these formulae in hand, both inclusions follow from the unitarity of $v,w$.
\end{proof}

Good news, we can now formulate our third Peter-Weyl theorem, as follows:

\index{Frobenius isomorphism}
\index{Peter-Weyl}

\begin{theorem}[PW3]
We have a direct sum decomposition of linear spaces
$$C(G)=\bigoplus_{v\in Irr(G)}M_{\dim(v)}(\mathbb C)$$
with the summands being pairwise orthogonal with respect to the scalar product
$$<f,g>=\int_Gf\bar{g}$$
where $\int_G$ is the averaging over $G$.
\end{theorem}

\begin{proof}
This is something very standard, the idea being as follows:

\medskip

(1) By combining our previous two Peter-Weyl results, from Theorem 9.16 and Theorem 9.18, we deduce that we have a linear space decomposition as follows:
$$C(G)
=\sum_{v\in Irr(G)}C_v
=\sum_{v\in Irr(G)}M_{\dim(v)}(\mathbb C)$$

Thus, in order to conclude, it is enough to prove that for any two irreducible representations $v,w\in Irr(A)$, the corresponding spaces of coefficients are orthogonal:
$$v\not\sim w\implies C_v\perp C_w$$ 

(2) For this purpose, we will use Theorem 9.22, telling us that for any representation $v$ we have the following formula, where $P$ is the orthogonal projection on $Fix(v)$:
$$\left(id\otimes\int_G\right)v=P$$

(3) We will also use the fact, that we just established in Theorem 9.25, that for any two representations $v,w$ we have available the Frobenius isomorphism, as follows:
$$Hom(v,w)\simeq Fix(v\otimes\bar{w})$$

(4) Now back to our orthogonality question from (1), let us set indeed:
$$P_{ia,jb}=\int_Gv_{ij}w_{ab}^*$$

Then $P$ is the orthogonal projection onto the following vector space:
$$Fix(v\otimes\bar{w})
\simeq Hom(v,w)
=\{0\}$$

Thus we have $P=0$, and this gives the result.
\end{proof}

The story is not over with this, and as a final result on the subject, we have the following theorem regarding the characters, completing the Peter-Weyl theory:

\index{Peter-Weyl}
\index{central function}
\index{algebra of characters}

\begin{theorem}[PW4]
The characters of irreducible representations belong to
$$C(G)_{central}=\left\{f\in C(G)\Big|f(gh)=f(hg),\forall g,h\in G\right\}$$
called algebra of central functions on $G$, and form an orthonormal basis of it.
\end{theorem}

\begin{proof}
We have several things to be proved, the idea being as follows:

\medskip

(1) Observe first that $C(G)_{central}$ is indeed an algebra, which contains all the characters. Conversely, consider a function $f\in C(G)$, written as follows:
$$f=\sum_{v\in Irr(G)}f_v$$

The condition $f\in C(G)_{central}$ states then that for any $v\in Irr(G)$, we must have:
$$f_v\in C(G)_{central}$$

But this means precisely that the coefficient $f_v$ must be a scalar multiple of $\chi_v$, and so the characters form a basis of $C(G)_{central}$, as stated. 

\medskip

(2) The fact that we have an orthogonal basis follows from Theorem 9.26. 

\medskip

(3) As for the fact that the characters have norm 1, this follows from:
\begin{eqnarray*}
\int_G\chi_v\chi_v^*
&=&\sum_{ij}\int_Gv_{ii}v_{jj}^*\\
&=&\sum_i\frac{1}{N}\\
&=&1
\end{eqnarray*}

To be more precise, we have used here the fact that the above integrals $\int_Gv_{ij}v_{kl}^*$ form the orthogonal projection onto the following vector space: 
$$Fix(v\otimes\bar{v})
\simeq End(v)
=\mathbb C1$$

Thus, the proof of our theorem is now complete.
\end{proof}

As a key observation now, complementing Theorem 9.27, observe that a function $f:G\to\mathbb C$ is central, in the sense that it satisfies $f(gh)=f(hg)$, precisely when it satisfies the following condition, saying that it must be constant on conjugacy classes:
$$f(ghg^{-1})=f(h),\forall g,h\in G$$

Now the point is that this makes the algebra of central functions something quite easy to compute, via standard algebra, and this puts us on the right track for computing $Irr(G)$. Or at least, this is how the theory goes, because there are many tricks too. 

\bigskip

As a basic illustration for this method, which clarifies some previous considerations from chapter 4, in relation with our study there of the finite abelian groups, we have:

\begin{theorem}
For a finite abelian group $G$ the irreducible representations are all $1$-dimensional, equal to their own characters,
$$\chi:G\to\mathbb T$$
and these characters form the dual discrete abelian group $\widehat{G}$.
\end{theorem}

\begin{proof}
This comes indeed from the Peter-Weyl theory, as follows:

\medskip

(1) Since our group $G$ was assumed to be abelian, any function $f:G\to\mathbb C$ is obviously central, so the algebra of central functions is $C(G)$ itself:
$$C(G)_{central}=C(G)$$

(2) Thus the decomposition of $C(G)$ from Theorem 9.26 reduces in this case to the decomposition of $C(G)_{central}$ from Theorem 9.27, and in particular, the irreducible representations $u\in Irr(G)$ must be all 1-dimensional, equal to their own characters $\chi_u$.

\medskip

(3) Finally, the last assertion is something that we know well from chapter 4, and with the extra comment that we have in fact an isomorphism $\widehat{G}\simeq G$, coming from the structure theorem for the finite abelian groups, as explained there. 

\medskip

(4) As a final comment on this, observe that $\widehat{G}\simeq G$, or the structure theorem for the finite abelian groups, do not come from Peter-Weyl for the abelian groups, whose conclusions reduce to what is said in the statement. Thus, although Peter-Weyl for the finite abelian groups does part of the job that we did in chapter 4, this is not everything, and our arithmetic work there remains something needed, going beyond Peter-Weyl.
\end{proof}

Getting now to the non-abelian case, things here can be quite complicated. For the simplest non-abelian group that we know, namely $S_3=D_3$, the result is as follows:

\begin{theorem}
The group $S_3=D_3$ has $3$ irreducible representations, namely:
\begin{enumerate}
\item The trivial representation, $g\to 1$.

\item The signature representation, $g\to\varepsilon(g)$.

\item The $2D$ representation $u-1$, with $u$ being the standard $3D$ representation.
\end{enumerate}
\end{theorem}

\begin{proof}
We certainly have the representations in (1) and (2), which are obviously irreducible, and non-equivalent. Now let us look at the 3D representation:
$$u:[S_3=D_3]\subset O_3\subset U_3$$

Since this representation appears via the permutation matrices, which sum up to 1 on each row, we conclude that the all-one vector is fixed by this representation:
$$\begin{pmatrix}1\\ 1\\1\end{pmatrix}\in Fix(u)$$

Thus, we can consider the following representation, which is 2-dimensional:
$$v=u-1$$

And we can stop here, because our group being non-abelian, and of order 6, a quick look at Theorem 9.26 shows that the decomposition there must come from:
$$6=1+1+4$$

Thus, $1,\varepsilon,v$ are indeed the irreducible representations, as stated.
\end{proof}

Regarding now more complicated groups, with a bit more work the ideas in the above proof extend to all dihedral groups $D_N$, and we will leave the study here as an exercise. As for the symmetric groups $S_N$, the situation here is more complicated, but with these being our favorite finite groups, we will keep them on our radar. More later. 

\section*{9e. Exercises}

This was a quite standard chapter, on the basics of representation theory, and as exercises on this, which are all basic as well, we have:

\begin{exercise}
Review if needed the basics of the $\otimes$ operation.
\end{exercise}

\begin{exercise}
Learn a bit about categories, and about tensor categories too.
\end{exercise}

\begin{exercise}
Review if needed the Stone-Weierstrass theorem, and its applications.
\end{exercise}

\begin{exercise}
Clarify the linear algebra details, in our Ces\`aro limiting proofs.
\end{exercise}

\begin{exercise}
Compute fixed points, for some representations of your choice.
\end{exercise}

\begin{exercise}
Compute the polynomial integrals, for some $G\subset U_N$ of your choice.
\end{exercise}

\begin{exercise}
Learn more about characters, and about character tables too.
\end{exercise}

\begin{exercise}
Learn more about central functions, and their applications.
\end{exercise}

As bonus exercise, try understanding what happens to the Peter-Weyl theory developed above, when the group in question $G\subset U_N$ is only assumed to be compact.

\chapter{Brauer theorems}

\section*{10a. Tensor categories}

You might probably think, with the Peter-Weyl theory learned, and with this being certainly a great thing, end of our discussion regarding general representation theory, and time to get into applications, and all sorts of explicit computations.

\bigskip

In answer, welcome to modern algebra, where theories come in layers. In what concerns group representations, the Peter-Weyl theory that we learned is just the first layer, foundations. Several other layers to come, substantially improving our knowledge.

\bigskip

As a starting point, we have seen in the previous chapter that, no matter what we want to do with our group $G\subset U_N$, of representation theory type, we must compute the spaces $End(u^{\otimes k})$. In view of this, it is convenient to introduce the following notion:

\index{Lie group}
\index{Lie algebra}
\index{Tannakian category}
\index{tensor category}

\begin{definition}
The Tannakian category of $G\subset_uU_N$ is the collection of spaces
$$C_{kl}=Hom(u^{\otimes k},u^{\otimes l})$$
of intertwiners between its various Peter-Weyl representations.
\end{definition}

This notion in something quite subtle, that will keep us busy for this whole chapter, and for some time afterwards too. Here are a few basic observations, about this:

\bigskip

(1) To start with, as already mentioned, we are mainly interested in computing the spaces $C_{kk}=End(u^{\otimes k})$. However, passing to the whole collection of spaces $C_{kl}$ is not a big deal, and this due to the Frobenius isomorphisms $Hom(v,w)\simeq Fix(v\otimes\bar{w})$, which show that computing $\{C_{kk}\}$ is more or less the same thing as computing $\{C_{kl}\}$.

\bigskip

(2) This being said, which one to choose, $\{C_{kk}\}$ or $\{C_{kl}\}$? In favor of $\{C_{kk}\}$, we have the fact that these spaces embed in each other, and form an algebra $\cup_kC_{kk}$. But then, we can equally say that in favor of $\{C_{kl}\}$, these form a category, which is nice too.

\bigskip

(3) Adding to our dillema, the algebras of type $\cup_kC_{kk}$ are something very particular, and powerful, namely planar algebras in the sense of Jones \cite{jo6}. However, due to the Frobenius duality mentioned in (1), somewhat the same goes for the categories $\{C_{kl}\}$.

\bigskip

(4) In short, we seem to have a tie here, related to many deep things. In order to decide what to do, let us make the observation that the groups $\{1\}$ and $\{\pm1\}$ have the same $\{C_{kk}\}$ spaces, but are distinguished by the $\{C_{kl}\}$ spaces. So, we will choose for now the categories $\{C_{kl}\}$, as in Definition 10.1, and leave the algebras $\cup_kC_{kk}$ for later.

\bigskip

This was our first series of comments on Definition 10.1, of rather philosophical nature. Getting now to more concrete things, let us make the following observations:

\begin{proposition}
The Tannakian category associated to a subgroup $G\subset U_N$ is formed by the following vector spaces, with $H=\mathbb C^N$,
$$C_{kl}=\left\{T\in\mathcal L(H^{\otimes k},H^{\otimes l})\Big|Tg^{\otimes k}=g^{\otimes l}T,\forall g\in G\right\}$$
and with this making it clear that we have a category. In fact, what we have is a tensor $*$-category, appearing as subcategory of $E_{kl}=\mathcal L(H^{\otimes k},H^{\otimes l})$.
\end{proposition} 

\begin{proof}
These are things which follow from our general theory from chapter 9, but which are probably worth an independent look, now since everything is quite concrete:

\medskip

(1) To start with, due to $g\in G\subset U_N\subset\mathcal L(H)$, we have $g^{\otimes k}\in\mathcal L(H^{\otimes k})$ for any $k$, so the equality $Tg^{\otimes k}=g^{\otimes l}T$ makes indeed sense, as an equality of maps as follows: 
$$Tg^{\otimes k},g^{\otimes l}T\in \mathcal L(H^{\otimes k},H^{\otimes l})$$

(2) It is also clear by definition that each $C_{kl}$ is a complex vector space. Moreover, it is also clear from definitions that $C=(C_{kl})$ is indeed a category, in the sense that:
$$T\in C_{kl}\ ,\ S\in C_{lm}\implies ST\in C_{km}$$

(3) Next, we have as well a tensor structure on our category, coming from:
$$T\in C_{kl}\ ,\ S\in C_{mn}\implies S\otimes T\in C_{km,ln}$$

And finally, we have an involution too on our category, coming from:
$$T\in C_{kl}\implies T^*\in C_{lk}$$

(4) Thus, we are led to the conclusions in the statement, and with everything being of course something that we already knew, coming from our theory from chapter 9.
\end{proof}

Quite remarkably, the subgroup $G\subset U_N$ can be reconstructed from its Tannakian category $C=(C_{kl})$, and in a very simple and efficient way, as follows:

\index{Tannakian duality}
\index{Peter-Weyl theory}

\begin{claim}
Given a finite subgroup $G\subset U_N$, we have
$$G=\left\{g\in U_N\Big|Tg^{\otimes k}=g^{\otimes l}T,\forall k,l,\forall T\in C_{kl}\right\}$$
where $C=(C_{kl})$ is the associated Tannakian category.
\end{claim}

As a first observation, in case you are a bit familiar with operator algebras, this reminds the bicommutant theorem of von Neumann, stating that an operator algebra $A\subset B(H)$ must equal its bicommutant, $A=A''$. To be more precise, as in the case of the bicommutant theorem, where we have an inclusion of algebras $A\subset A''$, that we must prove to be an equality, what we have in Claim 10.3 is an inclusion of groups, appearing by a double commutation procedure, that we must prove to be an equality:

\begin{proposition}
Given a finite subgroup $G\subset U_N$, the set
$$\widetilde{G}=\left\{g\in U_N\Big|Tg^{\otimes k}=g^{\otimes l}T,\forall k,l,\forall T\in C_{kl}\right\}$$
is closed subgroup of $U_N$, and we have inclusions of groups as follows:
$$G\subset\widetilde{G}\subset U_N$$
Claim 10.3 states that the inclusion $G\subset\widetilde{G}$ must be an equality.
\end{proposition}

\begin{proof}
Let us first prove that $\widetilde{G}$ is a group. Assuming $g,h\in\widetilde{G}$, we have $gh\in\widetilde{G}$, due to the following computation, valid for any $k,l$ and any $T\in C_{kl}$:
\begin{eqnarray*}
T(gh)^{\otimes k}
&=&Tg^{\otimes k}h^{\otimes k}\\
&=&g^{\otimes l}Th^{\otimes k}\\
&=&g^{\otimes l}h^{\otimes l}T\\
&=&(gh)^{\otimes l}T
\end{eqnarray*}

Also, we have $1\in\widetilde{G}$, trivially. Finally, assuming $g\in\widetilde{G}$, we have:
\begin{eqnarray*}
T(g^{-1})^{\otimes k}
&=&(g^{-1})^{\otimes l}[g^{\otimes l}T](g^{-1})^{\otimes k}\\
&=&(g^{-1})^{\otimes l}[Tg^{\otimes k}](g^{-1})^{\otimes k}\\
&=&(g^{-1})^{\otimes l}T
\end{eqnarray*}

Thus we have $g^{-1}\in\widetilde{G}$, and so $\widetilde{G}$ is a group, as claimed. Finally, the fact that we have an inclusion $G\subset\widetilde{G}$, and that $\widetilde{G}\subset U_N$ is closed, are both clear from definitions.
\end{proof}

In practice now, Claim 10.3 is something which comes from Peter-Weyl theory, or from Tannaka duality. However, one problem comes from the fact that the technical results that we need, Peter-Weyl or Tannaka, are in fact those in general form, for the compact subgroups $G\subset U_N$. Thus, we will not prove the above claim right away, and defer the discussion for later in this chapter, when we will talk about compact groups $G\subset U_N$. In fact, Claim 10.3 holds for such groups, and this is what we will prove, later.

\bigskip

This being said, we have certainly work to do for the finite subgroups $G\subset U_N$ that we know, consisting in computing the corresponding categories $\{C_{kl}\}$, which is, and we should remind this, what matters. Here is the computation for $S_N$, which is something quite fundamental, and which proves as well, as a bonus, Claim 10.3 in this case:

\index{symmetric group}
\index{partitions}
\index{Kronecker symbols}
\index{fitting indices}

\begin{theorem}
For the symmetric group $S_N\subset O_N$ we have the formula
$$C_{kl}=span\left(T_\pi\Big|\pi\in P(k,l)\right)$$
with $P(k,l)$ being the set of partitions of $k$ upper points, and $l$ lower points, and where
$$T_\pi(e_{i_1}\otimes\ldots\otimes e_{i_k})=\sum_{j_1\ldots j_l}\delta_\pi
\begin{pmatrix}i_1&\ldots&i_k\\ j_1&\ldots&j_l\end{pmatrix}
e_{j_1}\otimes\ldots\otimes e_{j_l}$$
where $\delta_\pi\in\{0,1\}$ is $1$ when the indices fit, and is $0$ otherwise. Also, Claim 10.3 holds.
\end{theorem}

\begin{proof}
This is something quite elementary, and importantly for us, fundamental for what is to follow in the remainder of this book, the idea being as follows:

\medskip

(1) First, at the level of notations, we denote as mentioned by $P(k,l)$ the set of partitions of $k$ upper points, and $l$ lower points, and with these partitions being represented as pictures, with the blocks being represented by strings. As an example, here is a partition in $P(3,3)$, with two blocks, represented by strings, in the obvious way:
$$
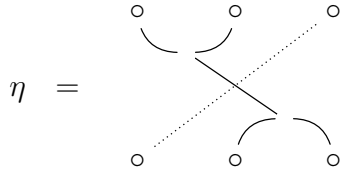
\xymatrix@R=2mm@C=3mm{\\ \\ \eta\ \ =\\ \\}\ \ \ 
\xymatrix@R=2mm@C=3mm{
\circ\ar@/_/@{-}[dr]&&\circ&&\circ\ar@{.}[ddddllll]\\
&\ar@/_/@{-}[ur]\ar@{-}[ddrr]\\
\\
&&&\ar@/^/@{-}[dr]\\
\circ&&\circ\ar@/^/@{-}[ur]&&\circ}$$

(2) Now given $\pi\in P(k,l)$ and multi-indices $i=(i_1,\ldots,i_k)$ and $j=(j_1,\ldots,j_l)$, we can put $i,j$ on the legs of $\pi$, in the obvious way. Then, if the ``indices fit'', meaning that all the strings of $\pi$ join equal indices of $i,j$, we set $\delta_\pi\binom{i}{j}=1$. Otherwise, we set $\delta_\pi\binom{i}{j}=0$. As an example, for the above partition $\eta\in P(3,3)$, we have the following formula:
$$\delta_\eta\begin{pmatrix}a&b&c\\d&e&f\end{pmatrix}=\delta_{abef}\delta_{cd}$$

(3) In order to prove now the result, let us first work out the case $k=0$, that we will regularly need in what follows. It is traditional here, and convenient, to change a bit notations. So, let us associate to any $\pi\in P(k)$ a vector, as follows:
$$\xi_\pi=\sum_{i_1\ldots i_k}\delta_\pi(i_1 \ldots i_k)e_{i_1}\otimes\ldots\otimes e_{i_k}$$

With this notation, we must prove that we have the following equality:
$$C_{0k}=span\left(\xi_\pi\Big|\pi\in P(k)\right)$$

(4) Let us first prove that we have $\supset$. Given $\sigma\in S_N$, we have indeed:
\begin{eqnarray*}
\sigma^{\otimes k}\xi_\pi
&=&\sum_{i_1\ldots i_k}\delta_\pi(i_1 \ldots i_k)\sigma(e_{i_1})\otimes\ldots\otimes\sigma(e_{i_k})\\
&=&\sum_{i_1\ldots i_k}\delta_\pi(i_1 \ldots i_k)e_{\sigma(i_1)}\otimes\ldots\otimes e_{\sigma(i_k)}\\
&=&\sum_{i_1\ldots i_k}\delta_\pi(\sigma^{-1}(j_1)\ldots\sigma^{-1}(j_k))e_{j_1}\otimes\ldots\otimes e_{j_k}\\
&=&\sum_{i_1\ldots i_k}\delta_\pi(j_1\ldots j_k)e_{j_1}\otimes\ldots\otimes e_{j_k}\\
&=&\xi_\pi
\end{eqnarray*}

(5) In order to prove now $\subset$, consider an arbitrary vector of $\mathbb C^N$, as follows:
$$\xi=\sum_{i_1\ldots i_k}\lambda_{i_1\ldots i_k}e_{i_1}\otimes\ldots\otimes e_{i_k}$$

Given $\sigma\in S_N$, by reasoning as before, we have the following formula:
$$\sigma^{\otimes k}\xi=\sum_{i_1\ldots i_k}\lambda_{\sigma^{-1}(i_1)\ldots\sigma^{-1}(i_k)}e_{i_1}\otimes\ldots\otimes e_{i_k}$$

Thus the condition $\sigma^{\otimes k}\xi=\xi$ for any $\sigma\in S_N$ is equivalent to:
$$\lambda_{i_1\ldots i_k}=\lambda_{\sigma(i_1)\ldots\sigma(i_k)}\quad,\quad\forall i,\sigma$$

But this latter condition is equivalent to the following condition:
$$\ker i=\ker j\implies\lambda_i=\lambda_j$$

Thus, we are led to the conclusion that $\lambda:\{1,\ldots,N\}^k\to\mathbb C$ must come from a function $\varphi:P(k)\to\mathbb C$, via a formula of type $\lambda_i=\varphi(\ker i)$, and it follows that the inclusion $\supset$ that we established in (4) is indeed an equality, as desired.

\medskip

(6) Summarizing, and getting back now to our theorem as stated, we have proved the formula of $C_{kl}$ there, in the case $k=0$. In order to pass now to the general case, two methods are available. We can either fine-tune the above computation, and we will leave this as an instructive exercise, or we can argue that the result at $k=0$ gives the result in general, via Frobenius duality, and we will leave this as an instructive exercise too.

\medskip

(7) Regarding now Claim 10.3, consider the intermediate subgroup $S_N\subset\widetilde{S}_N\subset O_N$, constructed as in Proposition 10.4, that we must prove to be equal to $S_N$. In order to prove this equality, we use the following one-block ``fork'' partition:
$$\xymatrix@R=1mm@C=2mm{\\ \\ \mu\ \ =\\ \\ }\ \ \ 
\xymatrix@R=2mm@C=3mm{
\circ\ar@/_/@{-}[dr]&&\circ\\
&\ar@/_/@{-}[ur]\ar@{-}[dd]\\
&&&\\
&\circ}$$

The linear map associated to $\mu$ is then given by the following formula:
$$T_\mu(e_i\otimes e_j)=\delta_{ij}e_i$$

We therefore have the following formula, valid for any $g\in O_N$:
$$(T_\mu g^{\otimes 2})_{i,jk}
=\sum_{lm}(T_\mu)_{i,lm}(g^{\otimes 2})_{lm,jk}
=g_{ij}g_{ik}$$

On the other hand, we have as well the following formula:
$$(gT_\mu)_{i,jk}
=\sum_lg_{il}(T_\mu)_{l,jk}
=\delta_{jk}g_{ij}$$

Thus, we have the following equivalence, valid for any $g\in O_N$:
$$T_\mu g^{\otimes 2}=gT_\mu\iff g_{ij}g_{ik}=\delta_{jk}g_{ij},\forall i,j,k$$

Now by assuming $g\in\widetilde{S}_N$, the formula on the left holds, so the formula on the right must hold too. But this shows that we must have $g_{ij}\in\{0,1\}$, with exactly one 1 entry on each of the rows of $g$. Thus we must have $g\in S_N$, which finishes the proof. 
\end{proof}

The above result is quite encouraging, and suggests looking into other reflection groups. For $H_N$, the group of symmetries of the unit cube in $\mathbb R^N$, we have:

\index{hyperoctahedral group}
\index{hypercube}
\index{blocks of even size}

\begin{theorem}
For the hyperoctadedral group $H_N=\mathbb Z_2\wr S_N\subset O_N$ we have
$$C_{kl}=span\left(T_\pi\Big|\pi\in P_{even}(k,l)\right)$$
where $P_{even}$ are the partitions all whose blocks have even size. Also, Claim 10.3 holds.
\end{theorem}

\begin{proof}
This follows a bit as for $S_N$. Consider indeed a vector of $\mathbb C^N$, as follows:
$$\xi=\sum_{i_1\ldots i_k}\lambda_{i_1\ldots i_k}e_{i_1}\otimes\ldots\otimes e_{i_k}$$

Then the condition $g^{\otimes k}\xi=\xi$ for any $g=\sigma^w\in H_N$ is equivalent to:
$$\lambda_{i_1\ldots i_k}=w_{i_1}\ldots w_{i_k}\lambda_{\sigma(i_1)\ldots\sigma(i_k)}\quad,\quad\forall i,\sigma$$

But this latter condition is equivalent to the following condition, along with the fact that we must have $w_{i_1}\ldots w_{i_k}=1$, which amounts in saying that $\ker i\in P_{even}$:
$$\ker i=\ker j\implies\lambda_i=\lambda_j$$

Thus, we are led to the conclusion in the statement. Finally, regarding the verification of Claim 10.3, this follows again as for $S_N$, by using the following partition:
$$\xymatrix@R=0.5mm@C=2mm{\\ \\ \\ \chi\ \ =\\ \\ }\ \ \ 
\xymatrix@R=2mm@C=3mm{
\circ\ar@/_/@{-}[dr]&&\circ\\
&\ar@/_/@{-}[ur]\ar@{-}[dd]\\
&&&\\
&\ar@/^/@{-}[dr]\ar@/_/@{-}[dl]\\
\circ&&\circ}$$

And we will leave the verifications here, which are elementary, as an exercice.
\end{proof}

Quite remarkably, the results for $S_N,H_N$ can be generalized as follows:

\index{wreath product}
\index{complex reflection group}

\begin{theorem}
For the complex reflection group $H_N^s=\mathbb Z_s\wr S_N\subset U_N$ we have
$$C_{kl}=span\left(T_\pi\Big|\pi\in P^s(k,l)\right)$$
where $P^s(k,l)$ is the set of partitions of $k$ upper points and $l$ lower points, satisfying
$$\#\circ=\#\bullet(s)$$
as a weighted equality, in each block. Also, Claim 10.3 holds.
\end{theorem}

\begin{proof}
This follows a bit as for $S_N,H_N$. Consider indeed a vector, as follows:
$$\xi=\sum_{i_1\ldots i_k}\lambda_{i_1\ldots i_k}e_{i_1}\otimes\ldots\otimes e_{i_k}$$

Then the condition $g^{\otimes k}\xi=\xi$ for any $g=\sigma^w\in H_N^s$ is equivalent to:
$$\lambda_{i_1\ldots i_k}=w_{i_1}\ldots w_{i_k}\lambda_{\sigma(i_1)\ldots\sigma(i_k)}\quad,\quad\forall i,\sigma$$

But this latter condition is equivalent to the following condition, along with the fact that we must have $w_{i_1}\ldots w_{i_k}=1$, which amounts in saying that $\ker i\in P^s(k)$:
$$\ker i=\ker j\implies\lambda_i=\lambda_j$$

Thus, we are led to the conclusion in the statement.
\end{proof}

Finally, let us record the $s=\infty$ particular case of Theorem 10.7, as follows:

\index{full reflection group}

\begin{theorem}
For the full complex reflection group $K_N=\mathbb T\wr S_N\subset U_N$ we have
$$C_{kl}=span\left(T_\pi\Big|\pi\in\mathcal P_{even}(k,l)\right)$$
where $\mathcal P_{even}(k,l)$ is the set of partitions of $k$ upper points and $l$ lower points, satisfying
$$\#\circ=\#\bullet$$
as a weighted equality, in each block. Also, Claim 10.3 holds.
\end{theorem}

\begin{proof}
This appears indeed as the $s=\infty$ particular case of Theorem 10.7.
\end{proof}

\section*{10b. Compact groups}

Moving on, let us upgrade now our whole representation theory discussion, started at the beginning of chapter 9, to the level of arbitrary compact groups $G\subset U_N$. This can help indeed in our finite group quest, with our objectives being multiple, as follows:

\bigskip

(1) First, we would like to prove Claim 10.3, for the subgroups $G\subset U_N$ which are finite, or more generally compact, by using the Peter-Weyl theory for the groups $G\subset\widetilde{G}$ from Proposition 1.4, the point being that the group $\widetilde{G}$ constructed there is compact.

\bigskip

(2) We would like as well to compute the Tannakian categories of $O_N,U_N$, in a combinatorial way, in analogy with what we have for $H_N,K_N$. Finally, we would like to investigate $SU_2,SO_3$ too, with this being related to the material from chapter 7.

\bigskip

Getting to work now, there is in fact not much work to be done, because the Peter-Weyl theory for the finite groups, as developed in chapter 9, generalizes to the compact groups, in a quite straightforward way. In fact, things in chapter 9 were written precisely as for their generalization now, to the compact case, to be straightforward.

\bigskip

To start with, we can certainly talk about representations, as follows:

\begin{definition}
A representation of a compact group $G$ is a morphism
$$u:G\to U_N$$
assumed to be continuous. The character of $u$ is the function $\chi:G\to\mathbb C$ given by
$$g\to Tr(u_g)$$
where $Tr$ is the usual trace of the $N\times N$ matrices, $Tr(M)=\sum_iM_{ii}$.
\end{definition}

As in the finite group case, we can talk about operations on representations, and about intertwining spaces too, and we have a first Peter-Weyl result, as follows:

\begin{theorem}[PW1]
Let $u:G\to U_N$ be a representation, consider the algebra $A=End(u)$, and write its unit as follows, with $p_i$ being central minimal projections:
$$1=p_1+\ldots+p_k$$
The representation $u$ decomposes then as a direct sum, as follows,
$$u=u_1+\ldots+u_k$$
with each $u_i$ being an irreducible representation, obtained by restricting $u$ to $Im(p_i)$.
\end{theorem}

\begin{proof}
This is something that we know from chapter 9 for the finite groups, and the proof in the present setting, that of the compact groups, is identical.
\end{proof}

Regarding now the second Peter-Weyl result, this extends too, as follows:

\begin{theorem}[PW2]
Given $G\subset_uU_N$, any smooth irreducible representation 
$$v:G\to U_M$$
appears inside a tensor product of the fundamental representation $u$ and its adjoint $\bar{u}$.
\end{theorem}

\begin{proof}
This is something a bit more tricky, the idea being that any closed subgroup $G\subset U_N$ can be shown to be a smooth, or Lie group, and with this in hand, a representation $v:G\to U_M$ is smooth precisely when its matrix coefficients $g\to(v_g)_{ij}$ appear as polynomials in the standard matrix coordinates $g\to g_{ij}$, and their conjugates $g\to\overline{g}_{ij}$. But with this, we obtain the result, exactly as in chapter 9 for the finite groups.
\end{proof}

Next, we have an extension of the third Peter-Weyl result, as follows:

\begin{theorem}[PW3]
Given a closed subgroup $G\subset U_N$, we have
$$C(G)=\overline{\bigoplus}_{v\in Irr(G)}M_{\dim(v)}(\mathbb C)$$
with the summands being pairwise orthogonal with respect to the scalar product
$$<f,g>=\int_Gf\bar{g}$$
where $\int_G$ is the averaging over $G$.
\end{theorem}

\begin{proof}
Consider indeed the dense subalgebra $\mathcal C(G)\subset C(G)$ generated by the coefficients of the fundamental representation. By following the proof from chapter 9, notably with the remark that everything there regarding integration extends to the compact case, this algebra decomposes as follows, with the summands being pairwise orthogonal:
$$\mathcal C(G)=\bigoplus_{v\in Irr(G)}M_{\dim(v)}(\mathbb C)$$

Now by taking the closure, we are led to the conclusion in the statement.
\end{proof}

Finally, the fourth Peter-Weyl result extends as well, as follows:

\begin{theorem}[PW4]
The characters of irreducible representations belong to
$$\mathcal C(G)_{central}=\left\{f\in\mathcal C(G)\Big|f(gh)=f(hg),\forall g,h\in G\right\}$$
called algebra of smooth central functions on $G$, and form an orthonormal basis of it.
\end{theorem}

\begin{proof}
This is again something straightforward, by following the proof from chapter 9, for the finite groups, and with the algebra $\mathcal C(G)\subset C(G)$ being, as before, the dense subalgebra generated by the coefficients of the fundamental representation.
\end{proof}

Next, as in the finite group case, we can talk about Tannakian categories:

\index{Tannakian category}
\index{Peter-Weyl representations}

\begin{definition}
The Tannakian category associated to a closed subgroup $G\subset_uU_N$ is the collection $C=(C_{kl})$ of vector spaces
$$C_{kl}=Hom(u^{\otimes k},u^{\otimes l})$$
where the representations $u^{\otimes k}$ with $k=\circ\bullet\bullet\circ\ldots$ colored integer, defined by
$$u^{\otimes\emptyset}=1\quad,\quad 
u^{\otimes\circ}=u\quad,\quad
u^{\otimes\bullet}=\bar{u}$$
and multiplicativity, $u^{\otimes kl}=u^{\otimes k}\otimes u^{\otimes l}$, are the Peter-Weyl representations.
\end{definition}

As a first observation, the knowledge of the Tannakian category is more or less the same thing as the knowledge of the fixed point spaces, which appear as:
$$Fix(u^{\otimes k})=C_{\emptyset k}$$

Indeed, these latter spaces fully determine all the spaces $C_{kl}$, because of the Frobenius isomorphisms, which for the Peter-Weyl representations read:
\begin{eqnarray*}
C_{kl}
&=&Hom(u^{\otimes k},u^{\otimes l})\\
&\simeq&Hom(1,\bar{u}^{\otimes k}\otimes u^{\otimes l})\\
&=&Hom(1,u^{\otimes\bar{k}l})\\
&=&Fix(u^{\otimes\bar{k}l})
\end{eqnarray*}

However, in what regards the relation with the algebras $C_{kk}=End(u^{\otimes k})$, things here are more tricky, because the groups $\{1\}$ and $\{\pm 1\}$ are distinguished by their Tannakian categories, but not by their algebras $C_{kk}=End(u^{\otimes k})$. More on this later.

\bigskip

Good news, as a first application of our technology, we can now prove Claim 10.3, or rather the extension of Claim 10.3 to the case of the compact groups, as follows:

\index{Tannakian duality}
\index{Peter-Weyl theory}

\begin{theorem}
Given a compact subgroup $G\subset U_N$, we have
$$G=\left\{g\in U_N\Big|Tg^{\otimes k}=g^{\otimes l}T,\forall k,l,\forall T\in C_{kl}\right\}$$
where $C=(C_{kl})$ is the associated Tannakian category.
\end{theorem}

\begin{proof}
This is something coming from Peter-Weyl theory, as follows:

\medskip

(1) Consider the set constructed in the statement, namely:
$$\widetilde{G}=\left\{g\in U_N\Big|Tg^{\otimes k}=g^{\otimes l}T,\forall k,l,\forall T\in C_{kl}\right\}$$

As before in the finite group case, what we have here is an intermediate compact group, as follows, that we want to prove to be equal to $G$:
$$G\subset\widetilde{G}\subset U_N$$

(2) In order to prove this, consider the Tannakian category of $\widetilde{G}$, namely:
$$\widetilde{C}_{kl}=\left\{T\in\mathcal L(H^{\otimes k},H^{\otimes l})\Big|Tg^{\otimes k}=g^{\otimes l}T,\forall g\in\widetilde{G}\right\}$$

By functoriality, from $G\subset\widetilde{G}$ we obtain $\widetilde{C}\subset C$. On the other hand, according to the definition of $\widetilde{G}$, we have $C\subset\widetilde{C}$. Thus, we have the following equality:
$$C=\widetilde{C}$$

(3) Assume now by contradiction that $G\subset\widetilde{G}$ is not an equality. Then, at the level of algebras of functions, the following quotient map is not an isomorphism either:
$$C(\widetilde{G})\to C(G)$$ 

On the other hand, we know from Peter-Weyl that we have decompositions as follows, with the sums being over all irreducible unitary representations:
$$C(\widetilde{G})=\overline{\bigoplus}_{v\in Irr(\widetilde{G})}M_{\dim v}(\mathbb C)
\quad,\quad 
C(G)=\overline{\bigoplus}_{w\in Irr(G)}M_{\dim w}(\mathbb C)$$

Now observe that each unitary representation $v:\widetilde{G}\to U_K$ restricts into a certain representation $v':G\to U_K$. Since the quotient map $C(\widetilde{G})\to C(G)$ is not an isomorphism, we conclude that there is at least one representation $v$ satisfying:
$$v\in Irr(\widetilde{G})\quad,\quad v'\notin Irr(G)$$

(4) We are now in position to conclude. By using Peter-Weyl theory again, the above representation $v\in Irr(\widetilde{G})$ appears in a certain tensor power of the fundamental representation $u:\widetilde{G}\subset U_N$. Thus, we have inclusions of representations, as follows:
$$v\in u^{\otimes k}\quad,\quad v'\in u'^{\otimes k}$$

Now since we know that $v$ is irreducible, and that $v'$ is not, by using one more time Peter-Weyl theory, we conclude that we have a strict inequality, as follows:
\begin{eqnarray*}
\dim(\widetilde{C}_{kk})
&=&dim(End(u^{\otimes k}))\\
&<&dim(End(u'^{\otimes k})\\
&=&\dim(C_{kk})
\end{eqnarray*}

But this contradicts the equality $C=\widetilde{C}$ found in (2), which finishes the proof.

\medskip

(5) Alternatively, and we should mention this too, we can use Tannakian duality, in its general, standard form. This duality, that we will not prove here, states that any compact group $G$ appears as the group of endomorphisms of the canonical inclusion functor $Rep(G)\subset\mathcal H$, where $Rep(G)$ is the category of final dimensional continuous unitary representations of $G$, and $\mathcal H$ is the category of finite dimensional Hilbert spaces. 

\medskip

(6) Now in the case of a closed subgroup $G\subset_uU_N$, we know from Peter-Weyl theory that any $r\in Rep(G)$ appears as a subrepresentation $r\in u^{\otimes k}$. In categorical terms, this means that, with suitable definitions, $Rep(G)$ appears as a ``completion'' of the category $C=(C_{kl})$. Thus $C$ uniquely determines $G$, and we obtain the result.

\medskip

(7) As yet another proof, we can use the Tannakian duality in its soft, reduced form, that we will prove later in this chapter, with this giving right away the formula in the statement. More on this, when discussing this latter version of Tannakian duality.
\end{proof}

As a conclusion to all this, Claim 10.3 fully proved, both in its initial formulation, for the finite groups, and also in its more general formulation, for the compact groups, as above. All the above was of course quite brief, but we will be back to this topic, and to Tannakian duality in general, on numerous occasions, in what follows.

\section*{10c. Brauer theorems}

We have so far a beginning of theory based on Claim 10.3 and its philosophy, with the main examples being the reflection groups $G=\mathbb Z_s\wr S_N$. We would like now to look into the continuous groups $G\subset U_N$. Let us start with a basic investigation of the simplest such group, namely $O_N\subset U_N$. We can say the following, about this group:

\begin{proposition}
For the orthogonal group $O_N$ we have inclusions as follows,
$$span\left(T_\pi\Big|\pi\in P_2(k,l)\right)\subset C_{kl}\subset span\left(T_\pi\Big|\pi\in P(k,l)\right)$$
where $P_2(k,l)$ is the set of pairings of $k$ upper points, and $l$ lower points.
\end{proposition}

\begin{proof}
Since we have $S_N\subset O_N$, by functoriality and Theorem 10.5 we have:
$$C_{kl}\subset span\left(T_\pi\Big|\pi\in P(k,l)\right)$$

For the other inclusion, the one on the left, let us first work out, as usual, the case $k=0$, and with the change $k\leftrightarrow l$. For a pairing $\pi\in P_2(k)$ we set, as before:
$$\xi_\pi=\sum_{i_1\ldots i_k}\delta_\pi(i_1 \ldots i_k)e_{i_1}\otimes\ldots\otimes e_{i_k}$$

We must prove $\xi_\pi\in C_{0k}$. For this purpose, let us pick $g\in O_N$, and write:
$$g(e_i)=\sum_jg_{ji}e_j$$

We have then the following computation:
\begin{eqnarray*}
g^{\otimes k}\xi_\pi
&=&\sum_{i_1\ldots i_k}\delta_\pi(i_1 \ldots i_k)ge_{i_1}\otimes\ldots\otimes ge_{i_k}\\
&=&\sum_{i_1\ldots i_k}\sum_{j_1\ldots j_k}\delta_\pi(i_1 \ldots i_k)g_{j_1i_1}\ldots g_{j_ki_k}e_{j_1}\otimes\ldots\otimes e_{j_k}
\end{eqnarray*}

As an illustration now, let us see what happens for a simple pairing, such as $\pi=\cap\,\cap$. Here the above computation can be continued as follows:
\begin{eqnarray*}
g^{\otimes k}\xi_{\cap\,\cap}
&=&\sum_{i_1\ldots i_4}\sum_{j_1\ldots j_4}\delta_{\cap\,\cap}(i_1 \ldots i_4)g_{j_1i_1}\ldots g_{j_4i_4}e_{j_1}\otimes\ldots\otimes e_{j_4}\\
&=&\sum_{ab}\sum_{j_1\ldots j_4}g_{j_1a}g_{j_2a}g_{j_3b}g_{j_4b}e_{j_1}\otimes\ldots\otimes e_{j_4}\\
&=&\sum_{j_1\ldots j_4}\delta_{j_1j_2}\delta_{j_3j_4}e_{j_1}\otimes\ldots\otimes e_{j_4}\\
&=&\xi_{\cap\,\cap}
\end{eqnarray*}

The same computation works in general, and by using $gg^t=1$, we obtain $g^{\otimes k}\xi_\pi=\xi_\pi$, for any $\pi\in P_2(k)$. Thus, we have indeed inclusions as in the statement.
\end{proof}

The above result, perhaps coupled with a few more computations, that you can do by yourself, in order to evaluate the situation, suggests that we should have:
$$C_{kl}=span\left(T_\pi\Big|\pi\in P_2(k,l)\right)$$

However, this is hard to prove with bare hands, and we will have to trick. Our trick will be something quite natural, the idea being first to prove that the above spaces $span(T_\pi)$ ``qualify'' for what is expected from a Tannakian category, and then, conclude that we have equality, because $span(T_\pi)$ can only correspond to $O_N$. Let us start with:

\index{tensor category}
\index{flip operator}
\index{colored integer}

\begin{definition}
A tensor category over $H=\mathbb C^N$ is a collection $C=(C_{kl})$ of linear spaces $C_{kl}\subset\mathcal L(H^{\otimes k},H^{\otimes l})$ satisfying the following conditions:
\begin{enumerate}
\item $S,T\in C$ implies $S\otimes T\in C$.

\item If $S,T\in C$ are composable, then $ST\in C$.

\item $T\in C$ implies $T^*\in C$.

\item Each $C_{kk}$ contains the identity operator.

\item $C_{\emptyset k}$ with $k=\circ\bullet,\bullet\circ$ contain the operator $R:1\to\sum_ie_i\otimes e_i$.

\item $C_{kl,lk}$ with $k,l=\circ,\bullet$ contain the flip operator $\Sigma:a\otimes b\to b\otimes a$.
\end{enumerate}
\end{definition}

Here, as usual, the tensor powers $H^{\otimes k}$, which are Hilbert spaces depending on a colored integer $k=\circ\bullet\bullet\circ\ldots\,$, are defined by the following formulae, and multiplicativity:
$$H^{\otimes\emptyset}=\mathbb C\quad,\quad 
H^{\otimes\circ}=H\quad,\quad
H^{\otimes\bullet}=\bar{H}\simeq H$$

We have already met such categories, when dealing with the Tannakian categories of the closed subgroups $G\subset U_N$, and our knowledge can be summarized as follows:

\begin{theorem}
For a closed subgroup $G\subset_uU_N$, the associated Tannakian category
$$C_{kl}=Hom(u^{\otimes k},u^{\otimes l})$$
is a tensor category over the Hilbert space $H=\mathbb C^N$.
\end{theorem}

\begin{proof}
We know that the fundamental representation $u$ acts on the Hilbert space $H=\mathbb C^N$, and that its conjugate $\bar{u}$ acts on the Hilbert space $\bar{H}=\mathbb C^N$. Now by multiplicativity we conclude that any Peter-Weyl representation $u^{\otimes k}$ acts on the Hilbert space $H^{\otimes k}$, so that we have embeddings as in Definition 10.17, as follows:
$$C_{kl}\subset\mathcal L(H^{\otimes k},H^{\otimes l})$$

Regarding now the fact that the axioms (1-6) in Definition 10.17 are satisfied:

\medskip

(1-4) These axioms all follow from definitions, as we already know.

\medskip

(5) This follows from the fact that each element $g\in G$ is a unitary, which can be reformulated as follows, with $R:1\to\sum_ie_i\otimes e_i$ being the map in Definition 10.17:
$$R\in Hom(1,g\otimes\bar{g})\quad,\quad 
R\in Hom(1,\bar{g}\otimes g)$$

Indeed, given an arbitrary matrix $g\in M_N(\mathbb C)$, we have the following computation:
\begin{eqnarray*}
(g\otimes\bar{g})(R(1)\otimes1)
&=&\left(\sum_{ijkl}e_{ij}\otimes e_{kl}\otimes g_{ij}\bar{g}_{kl}\right)\left(\sum_ae_a\otimes e_a\otimes 1\right)\\
&=&\sum_{ika}e_i\otimes e_k\otimes g_{ia}\bar{g}_{ka}^*\\
&=&\sum_{ik}e_i\otimes e_k\otimes(gg^*)_{ik}
\end{eqnarray*}

We conclude from this that we have the following equivalence:
\begin{eqnarray*}
R\in Hom(1,g\otimes\bar{g})
&\iff&gg^*=1
\end{eqnarray*}

By replacing $g$ with its conjugate matrix $\bar{g}$, we have as well:
$$R\in Hom(1,\bar{g}\otimes g)\iff\bar{g}g^t=1$$

Thus, the two intertwining conditions in Definition 10.17 (5) are both equivalent to the fact that $g$ is unitary, and so these conditions are indeed satisfied, as desired.

\medskip

(6) This is again something elementary, coming from the fact that the various matrix coefficients $g\to g_{ij}$ and their complex conjugates $g\to\bar{g}_{ij}$ commute with each other. To be more precise, with $\Sigma:a\otimes b\to b\otimes a$ being the flip operator, we have:
\begin{eqnarray*}
(g\otimes h)(\Sigma\otimes id)(e_a\otimes e_b\otimes 1)
&=&\left(\sum_{ijkl}e_{ij}\otimes e_{kl}\otimes g_{ij}h_{kl}\right)(e_b\otimes e_a\otimes1)\\
&=&\sum_{ik}e_i\otimes e_k\otimes g_{ib}h_{ka}
\end{eqnarray*}

On the other hand, we have as well the following computation:
\begin{eqnarray*}
(\Sigma\otimes id)(h\otimes g)(e_a\otimes e_b\otimes 1)
&=&(\Sigma\otimes id)\left(\sum_{ijkl}e_{ij}\otimes e_{kl}\otimes h_{ij}g_{kl}\right)(e_a\otimes e_b\otimes1)\\
&=&(\Sigma\otimes id)\left(\sum_{ik}e_i\otimes e_k\otimes h_{ia}g_{kb}\right)\\
&=&\sum_{ik}e_k\otimes e_i\otimes h_{ia}g_{kb}\\
&=&\sum_{ik}e_i\otimes e_k\otimes h_{ka}g_{ib}
\end{eqnarray*}

Now since functions commute, $g_{ib}h_{ka}=h_{ka}g_{ib}$, this gives the result.
\end{proof}

With the above in hand, our purpose now will be that of showing that any closed subgroup $G\subset U_N$ is uniquely determined by its Tannakian category $C=(C_{kl})$:
$$G\leftrightarrow C$$

Getting started, we want to construct a correspondence $G\leftrightarrow C$, and we already know from Theorem 10.18 how the correspondence $G\to C$ appears, namely via:
$$C_{kl}=Hom(u^{\otimes k},u^{\otimes l})$$

Regarding now the construction in the other sense, $C\to G$, this is something very simple as well, coming from the following elementary result: 

\begin{theorem}
Given a tensor category $C=(C_{kl})$ over the space $H\simeq\mathbb C^N$,
$$G=\left\{g\in U_N\Big|Tg^{\otimes k}=g^{\otimes l}T\ ,\ \forall k,l,\forall T\in C_{kl}\right\}$$
is a closed subgroup $G\subset U_N$.
\end{theorem}

\begin{proof}
Consider indeed the closed subset $G\subset U_N$ constructed in the statement. We want to prove that $G$ is indeed a group, and the verifications here go as follows:

\medskip

(1) Given two matrices $g,h\in G$, their product satisfies $gh\in G$, due to the following computation, valid for any $k,l$ and any $T\in C_{kl}$:
\begin{eqnarray*}
T(gh)^{\otimes k}
&=&Tg^{\otimes k}h^{\otimes k}\\
&=&g^{\otimes l}Th^{\otimes k}\\
&=&g^{\otimes l}h^{\otimes l}T\\
&=&(gh)^{\otimes l}T
\end{eqnarray*}

(2) Also, we have $1\in G$, trivially. Finally, for $g\in G$ and $T\in C_{kl}$, we have:
\begin{eqnarray*}
T(g^{-1})^{\otimes k}
&=&(g^{-1})^{\otimes l}[g^{\otimes l}T](g^{-1})^{\otimes k}\\
&=&(g^{-1})^{\otimes l}[Tg^{\otimes k}](g^{-1})^{\otimes k}\\
&=&(g^{-1})^{\otimes l}T
\end{eqnarray*}

Thus we have $g^{-1}\in G$, and so $G$ is a group, as claimed.
\end{proof}

Summarizing, we have so far precise axioms for the tensor categories $C=(C_{kl})$, given in Definition 10.17, as well as correspondences as follows:
$$G\to C_G\quad,\quad 
C\to G_C$$

We will prove in what follows that these correspondences are inverse to each other. In order to get started, we first have the following technical result:

\begin{proposition}
Consider the following conditions:
\begin{enumerate}
\item $C=C_{G_C}$, for any tensor category $C$.

\item $G=G_{C_G}$, for any closed subgroup $G\subset U_N$.
\end{enumerate}
We have then $(1)\implies(2)$. Also, $C\subset C_{G_C}$ is automatic.
\end{proposition}

\begin{proof}
Given $G\subset U_N$, we have $G\subset G_{C_G}$. On the other hand, by using (1) we have $C_G=C_{G_{C_G}}$. Thus, we have an inclusion of closed subgroups of $U_N$, which becomes an isomorphism at the level of the associated Tannakian categories, so $G=G_{C_G}$. Finally, the fact that we have an inclusion $C\subset C_{G_C}$ is clear from definitions.
\end{proof}

In order to establish Tannakian duality, we will need some abstract constructions. Following Malacarne \cite{mal}, let us start with the following elementary fact:

\begin{proposition}
Given a tensor category $C=(C_{kl})$ over a Hilbert space $H$,
$$E_C
=\bigoplus_{k,l}C_{kl}
\subset\bigoplus_{k,l}B(H^{\otimes k},H^{\otimes l})
\subset B\left(\bigoplus_kH^{\otimes k}\right)$$
is a closed $*$-subalgebra. Also, inside this algebra,
$$E_C^{(s)}
=\bigoplus_{|k|,|l|\leq s}C_{kl}
\subset\bigoplus_{|k|,|l|\leq s}B(H^{\otimes k},H^{\otimes l})
=B\left(\bigoplus_{|k|\leq s}H^{\otimes k}\right)$$
is a finite dimensional $*$-subalgebra.
\end{proposition}

\begin{proof}
This is clear indeed from the categorical axioms from Definition 10.17.
\end{proof}

Now back to our reconstruction question, according to Proposition 10.20 we want to prove $C=C_{G_C}$, which is the same as proving the following equality:
$$E_C=E_{C_{G_C}}$$

But for this purpose, we can use a standard commutant trick, as follows:

\index{bicommutant theorem}

\begin{theorem}
For any $*$-algebra $A\subset M_N(\mathbb C)$ we have the equality
$$A=A''$$
where prime denotes the commutant, $X'=\left\{T\in M_N(\mathbb C)\big|Tx=xT,\forall x\in X\right\}$.
\end{theorem}

\begin{proof}
This is a particular case of von Neumann's bicommutant theorem, which follows from the explicit description of $A$ worked out in chapter 9, namely:
$$A=M_{n_1}(\mathbb C)\oplus\ldots\oplus M_{n_k}(\mathbb C)$$

Indeed, the center of each matrix algebra being reduced to the scalars, the commutant of this algebra is as follows, with each copy of $\mathbb C$ corresponding to a matrix block:
$$A'=\mathbb C\oplus\ldots\oplus\mathbb C$$

Now when taking once again the commutant, the computation is trivial, and we obtain in this way $A$ itself, and this leads to the conclusion in the statement.
\end{proof}

By using now the above bicommutant theorem, we are led to the following result:

\begin{proposition}
Given a Tannakian category $C$, the following are equivalent:
\begin{enumerate}
\item $C=C_{G_C}$.

\item $E_C=E_{C_{G_C}}$.

\item $E_C^{(s)}=E_{C_{G_C}}^{(s)}$, for any $s\in\mathbb N$.

\item $E_C^{(s)'}=E_{C_{G_C}}^{(s)'}$, for any $s\in\mathbb N$.
\end{enumerate}
In addition, the inclusions $\subset$, $\subset$, $\subset$, $\supset$ are automatically satisfied.
\end{proposition}

\begin{proof}
The equivalences are clear from definitions, and from the bicommutant theorem. As for the last assertion, we have $C\subset C_{G_C}$, which shows that we have as well: 
$$E_C\subset E_{C_{G_C}}$$

We therefore obtain by truncating $E_C^{(s)}\subset E_{C_{G_C}}^{(s)}$, and by taking the commutants, this gives $E_C^{(s)}\supset E_{C_{G_C}}^{(s)}$. Thus, we are led to the conclusion in the statement.
\end{proof}

Summarizing, we would like to prove that we have $E_C^{(s)'}\subset E_{C_{G_C}}^{(s)'}$. But this can be done by doing some algebra, and we are led to the following conclusion:

\index{Tannakian duality}

\begin{theorem}
The Tannakian duality constructions 
$$C\to G_C\quad,\quad 
G\to C_G$$
are inverse to each other.
\end{theorem}

\begin{proof}
This follows by doing some algebra, in order to prove that we have indeed $E_C^{(s)'}\subset E_{C_{G_C}}^{(s)'}$, as mentioned above, and we refer here to the paper of Malacarne \cite{mal}. Alternatively, this can be proved via standard Tannakian methods, and we refer here to the literature on the subject. For more on all this, you have as well my book \cite{ba2}.
\end{proof}

With this piece of general theory in hand, let us go back now to the orthogonal group $O_N$, and to partitions and pairings, as in Proposition 10.16. In order to construct a Tannakian category out of the pairings, via the operation $\pi\to T_\pi$, we will need:

\begin{proposition}
The assignement $\pi\to T_\pi$ is categorical, in the sense that
$$T_\pi\otimes T_\sigma=T_{[\pi\sigma]}\quad,\quad 
T_\pi T_\sigma=N^{c(\pi,\sigma)}T_{[^\sigma_\pi]}\quad,\quad 
T_\pi^*=T_{\pi^*}$$
where $c(\pi,\sigma)$ are certain integers, coming from the erased components in the middle.
\end{proposition}

\begin{proof}
This is something elementary, the computations being as follows:

\medskip

(1) The concatenation axiom follows from the following computation:
\begin{eqnarray*}
&&(T_\pi\otimes T_\sigma)(e_{i_1}\otimes\ldots\otimes e_{i_p}\otimes e_{k_1}\otimes\ldots\otimes e_{k_r})\\
&=&\sum_{j_1\ldots j_q}\sum_{l_1\ldots l_s}\delta_\pi\begin{pmatrix}i_1&\ldots&i_p\\j_1&\ldots&j_q\end{pmatrix}\delta_\sigma\begin{pmatrix}k_1&\ldots&k_r\\l_1&\ldots&l_s\end{pmatrix}e_{j_1}\otimes\ldots\otimes e_{j_q}\otimes e_{l_1}\otimes\ldots\otimes e_{l_s}\\
&=&\sum_{j_1\ldots j_q}\sum_{l_1\ldots l_s}\delta_{[\pi\sigma]}\begin{pmatrix}i_1&\ldots&i_p&k_1&\ldots&k_r\\j_1&\ldots&j_q&l_1&\ldots&l_s\end{pmatrix}e_{j_1}\otimes\ldots\otimes e_{j_q}\otimes e_{l_1}\otimes\ldots\otimes e_{l_s}\\
&=&T_{[\pi\sigma]}(e_{i_1}\otimes\ldots\otimes e_{i_p}\otimes e_{k_1}\otimes\ldots\otimes e_{k_r})
\end{eqnarray*}

(2) The composition axiom follows from the following computation:
\begin{eqnarray*}
&&T_\pi T_\sigma(e_{i_1}\otimes\ldots\otimes e_{i_p})\\
&=&\sum_{j_1\ldots j_q}\delta_\sigma\begin{pmatrix}i_1&\ldots&i_p\\j_1&\ldots&j_q\end{pmatrix}
\sum_{k_1\ldots k_r}\delta_\pi\begin{pmatrix}j_1&\ldots&j_q\\k_1&\ldots&k_r\end{pmatrix}e_{k_1}\otimes\ldots\otimes e_{k_r}\\
&=&\sum_{k_1\ldots k_r}N^{c(\pi,\sigma)}\delta_{[^\sigma_\pi]}\begin{pmatrix}i_1&\ldots&i_p\\k_1&\ldots&k_r\end{pmatrix}e_{k_1}\otimes\ldots\otimes e_{k_r}\\
&=&N^{c(\pi,\sigma)}T_{[^\sigma_\pi]}(e_{i_1}\otimes\ldots\otimes e_{i_p})
\end{eqnarray*}

(3) Finally, the involution axiom follows from the following computation:
\begin{eqnarray*}
&&T_\pi^*(e_{j_1}\otimes\ldots\otimes e_{j_q})\\
&=&\sum_{i_1\ldots i_p}<T_\pi^*(e_{j_1}\otimes\ldots\otimes e_{j_q}),e_{i_1}\otimes\ldots\otimes e_{i_p}>e_{i_1}\otimes\ldots\otimes e_{i_p}\\
&=&\sum_{i_1\ldots i_p}\delta_\pi\begin{pmatrix}i_1&\ldots&i_p\\ j_1&\ldots& j_q\end{pmatrix}e_{i_1}\otimes\ldots\otimes e_{i_p}\\
&=&T_{\pi^*}(e_{j_1}\otimes\ldots\otimes e_{j_q})
\end{eqnarray*}

Summarizing, our correspondence is indeed categorical.
\end{proof}

Good news, we can now prove the Brauer theorem for $O_N$, as follows:

\index{pairings}
\index{orthogonal group}
\index{Brauer theorem}

\begin{theorem}
For the orthogonal group $O_N$ we have
$$C_{kl}=span\left(T_\pi\Big|\pi\in P_2(k,l)\right)$$
where $P_2(k,l)$ is the set of pairings of $k$ upper points, and $l$ lower points.
\end{theorem}

\begin{proof}
We know from Proposition 10.16 that we have inclusions as follows:
$$C_{kl}\supset span\left(T_\pi\Big|\pi\in P_2(k,l)\right)$$

On the other hand, Proposition 10.25 shows that the spaces on the right form a Tannakian category, in the sense of Definition 10.17. Thus the Tannakian duality result from Theorem 10.24 applies, and provides us with a subgroup $G\subset U_N$, such that:
$$G=\left\{g\in U_N\Big|T_\pi g^{\otimes k}=g^{\otimes l}T_\pi,\forall k,l,\forall \pi\in P_2(k,l)\right\}$$

Moreover, by functoriality of Tannakian duality, we have $O_N\subset G$. But the relation $g^{\otimes 2}T_\cap=T_\cap$ with $g\in U_N$ implies $g\in O_N$, as explained in the proof of Theorem 10.18, so we have as well $G\subset O_N$. Thus, we have $G=O_N$, which gives the result.
\end{proof}

For the unitary group $U_N$ now, the result is similar, as follows:

\index{unitary group}
\index{matching pairings}
\index{Brauer theorem}

\begin{theorem}
For the unitary group $U_N$ we have
$$C_{kl}=span\left(T_\pi\Big|\pi\in\mathcal P_2(k,l)\right)$$
where $\mathcal P_2(k,l)$ is the set of matching pairings of $k$ upper points, and $l$ lower points.
\end{theorem}

\begin{proof}
The proof here is very similar to the proof for $O_N$, and in fact, even a bit simpler, after a close examination, and with the convention, in the statement, that matching means joining $\circ-\circ$ or $\bullet-\bullet$ on the horizontal, and $\circ-\bullet$ on the vertical.
\end{proof}

Quite nice all this, and more in this topic in the next chapter, and afterwards.

\section*{10d. Clebsch-Gordan}

As a last piece of Lie group theory, we are now in position of dealing with the representations of $SU_2$ and $SO_3$. The idea will be that of using the following isomorphism:
$$SU_2\simeq S^3_\mathbb R$$

To start with, we can integrate over the spheres, as follows:

\begin{theorem}
The polynomial integrals over the unit sphere $S^{N-1}_\mathbb R\subset\mathbb R^N$, with respect to the normalized, mass $1$ measure, are given by the following formula,
$$\int_{S^{N-1}_\mathbb R}x_1^{k_1}\ldots x_N^{k_N}\,dx=\frac{(N-1)!!k_1!!\ldots k_N!!}{(N+\Sigma k_i-1)!!}$$
valid when all exponents $k_i$ are even. If an exponent is odd, the integral vanishes.
\end{theorem}

\begin{proof}
This is something very standard, the idea being as follows:

\medskip

(1) Let us first discuss the case $N=2$. Here the sphere is the unit circle $\mathbb T$, with $z=e^{it}$ the standard coordinates are $\cos t,\sin t$, and we can use here the well-known Wallis formula, coming by double recurrence, which is as follows:
$$\int_0^{\pi/2}\cos^pt\sin^qt\,dt=\left(\frac{\pi}{2}\right)^{\varepsilon(p)\varepsilon(q)}\frac{p!!q!!}{(p+q+1)!!}$$

To be more precise, here $\varepsilon(p)=1$ if $p$ is even, and $\varepsilon(p)=0$ if $p$ is odd, and the double factorials are as follows, with the product ending at $2$ or $1$, depending on $m$:
$$m!!=(m-1)(m-3)(m-5)\ldots$$

Now observe that for even exponents $p,q$, the Wallis formula reads:
$$\int_0^{\pi/2}\cos^pt\sin^qt\,dt
=\frac{p!!q!!}{(p+q+1)!!}$$

Thus, we obtain the formula in the statement, so done with the case $N=2$. 

\medskip

(2) In the general case now, where the dimension $N\in\mathbb N$ is arbitrary, the integral in the statement can be written in spherical coordinates, as follows:
$$I=\frac{2^N}{A}\int_0^{\pi/2}\ldots\int_0^{\pi/2}x_1^{k_1}\ldots x_N^{k_N}J\,dt_1\ldots dt_{N-1}$$

Here $A$ is the area of the sphere, $J$ is the Jacobian, and the $2^N$ factor comes from the restriction to the $1/2^N$ part of the sphere where all the coordinates are positive. By using the Wallis formula, the normalization constant in front of the integral is:
$$\frac{2^N}{A}=\left(\frac{2}{\pi}\right)^{[N/2]}(N-1)!!$$

As for the unnormalized integral, this is given by:
\begin{eqnarray*}
I'=\int_0^{\pi/2}\ldots\int_0^{\pi/2}
&&(\cos t_1)^{k_1}
(\sin t_1\cos t_2)^{k_2}\\
&&\vdots\\
&&(\sin t_1\sin t_2\ldots\sin t_{N-2}\cos t_{N-1})^{k_{N-1}}\\
&&(\sin t_1\sin t_2\ldots\sin t_{N-2}\sin t_{N-1})^{k_N}\\
&&\sin^{N-2}t_1\sin^{N-3}t_2\ldots\sin^2t_{N-3}\sin t_{N-2}\\
&&dt_1\ldots dt_{N-1}
\end{eqnarray*}

(3) By rearranging the terms, the unnormalized integral is given by:
\begin{eqnarray*}
I'
&=&\int_0^{\pi/2}\cos^{k_1}t_1\sin^{k_2+\ldots+k_N+N-2}t_1\,dt_1\\
&&\int_0^{\pi/2}\cos^{k_2}t_2\sin^{k_3+\ldots+k_N+N-3}t_2\,dt_2\\
&&\vdots\\
&&\int_0^{\pi/2}\cos^{k_{N-2}}t_{N-2}\sin^{k_{N-1}+k_N+1}t_{N-2}\,dt_{N-2}\\
&&\int_0^{\pi/2}\cos^{k_{N-1}}t_{N-1}\sin^{k_N}t_{N-1}\,dt_{N-1}
\end{eqnarray*}

By using now the Wallis formula at $N=2$, we obtain the following formula:
\begin{eqnarray*}
I'
&=&\frac{k_1!!(k_2+\ldots+k_N+N-2)!!}{(k_1+\ldots+k_N+N-1)!!}\left(\frac{\pi}{2}\right)^{\varepsilon(N-2)}\\
&&\frac{k_2!!(k_3+\ldots+k_N+N-3)!!}{(k_2+\ldots+k_N+N-2)!!}\left(\frac{\pi}{2}\right)^{\varepsilon(N-3)}\\
&&\vdots\\
&&\frac{k_{N-2}!!(k_{N-1}+k_N+1)!!}{(k_{N-2}+k_{N-1}+l_N+2)!!}\left(\frac{\pi}{2}\right)^{\varepsilon(1)}\\
&&\frac{k_{N-1}!!k_N!!}{(k_{N-1}+k_N+1)!!}\left(\frac{\pi}{2}\right)^{\varepsilon(0)}
\end{eqnarray*}

(4) Now let $F$ be the part involving the double factorials, and $P$ be the part involving the powers of $\pi/2$, so that $I'=F\cdot P$. Regarding $F$, by cancelling terms we have:
$$F=\frac{k_1!!\ldots k_N!!}{(\Sigma k_i+N-1)!!}$$

As in what regards $P$, by summing the exponents, we obtain $P=\left(\frac{\pi}{2}\right)^{[N/2]}$. Now by putting everything together, we obtain the formula in the statement.
\end{proof}

Good news, we can now come back to $SU_2$, and we have the following result:

\begin{theorem}
The irreducible representations of $SU_2$ are all self-adjoint, and can be labelled by positive integers, with their fusion rules being as follows,
$$r_k\otimes r_l=r_{|k-l|}+r_{|k-l|+2}+\ldots+r_{k+l}$$
called Clebsch-Gordan rules. The corresponding dimensions are $\dim r_k=k+1$.
\end{theorem}

\begin{proof}
There are several proofs for this fact, the simplest one, with the knowledge that we have, being via purely algebraic methods, as follows:

\medskip

(1) Our first claim is that we have the following estimate, telling us that the even moments of the main character are smaller than the Catalan numbers:
$$\int_{SU_2}\chi^{2k}\leq C_k$$

(2) But this can be proved, and even in a strong form, with $=$ sign, by using $SU_2\simeq S^3_\mathbb R$. Indeed, in this picture, the main character is the following variable:
$$\chi\begin{pmatrix}x+iy&z+it\\ -z+it&x-iy\end{pmatrix}=2x$$

Now by using the formula in Theorem 10.28 we obtain, as desired:
\begin{eqnarray*}
\int_{S^3_\mathbb R}(2x)^{2k}
&=&4^k\cdot\frac{3!!(2k)!!}{(2k+3)!!}\\
&=&2^{2k+1}\cdot\frac{(2k)!}{2^kk!2^{k+1}(k+1)!}\\
&=&\frac{1}{k+1}\binom{2k}{k}\\
&=&C_k
\end{eqnarray*}

(3) Alternatively, the estimate in (1) can be deduced with purely algebraic methods, by using a Brauer type argument for $SU_2$, as follows:
\begin{eqnarray*}
\int_{SU_2}\chi^{2k}
&=&\dim(Fix(u^{\otimes 2k}))\\
&=&\dim\left(span\left(T_\pi'\Big|\pi\in NC_2(2k)\right)\right)\\
&\leq&|NC_2(2k)|\\
&=&C_k
\end{eqnarray*}

To be more precise, we have a Brauer theorem for $SU_2$, involving pairings, but coming from a different implementation $\pi\to T_\pi'$ of the pairings, involving some signs. And with this being proved as the Brauer theorem for $O_N$, with modifications where needed.

\medskip

(4) Long story short, we have our estimate in (1), and this is all that we need. Our claim is that we can construct, by recurrence on $k\in\mathbb N$, a sequence $r_k$ of irreducible, self-adjoint and distinct representations of $SU_2$, satisfying:
$$r_0=1\quad,\quad
r_1=u\quad,\quad 
r_k+r_{k-2}=r_{k-1}\otimes r_1$$

Indeed, assume that $r_0,\ldots,r_{k-1}$ are constructed, and let us construct $r_k$. We have:
$$r_{k-1}+r_{k-3}=r_{k-2}\otimes r_1$$

Thus $r_{k-1}\subset r_{k-2}\otimes r_1$, and since $r_{k-2}$ is irreducible, by Frobenius we have:
$$r_{k-2}\subset r_{k-1}\otimes r_1$$

We conclude there exists a certain representation $r_k$ such that:
$$r_k+r_{k-2}=r_{k-1}\otimes r_1$$

(5) By recurrence, $r_k$ is self-adjoint. Now observe that according to our recurrence formula, we can split $u^{\otimes k}$ as a sum of the following type, with positive coefficients:  
$$u^{\otimes k}=c_kr_k+c_{k-2}r_{k-2}+\ldots$$

We conclude by Peter-Weyl that we have an inequality as follows, with equality precisely when $r_k$ is irreducible, and non-equivalent to the other summands $r_i$:
$$\sum_ic_i^2\leq\dim(End(u^{\otimes k}))$$

(6) But by (1) the number on the right is $\leq C_k$, and some straightforward combinatorics, based on the fusion rules, shows that the number on the left is $C_k$ as well:
$$C_k=\sum_ic_i^2\leq\dim(End(u^{\otimes k}))=\int_{SU_2}\chi^{2k}\leq C_k$$

Thus we have equality in our estimate, so our representation $r_k$ is irreducible, and non-equivalent to $r_{k-2},r_{k-4},\ldots$ Moreover, this representation $r_k$ is not equivalent to $r_{k-1},r_{k-3},\ldots$ either, with this coming from $r_p\subset u^{\otimes p}$ for any $p$, and from:
$$\dim(Fix(u^{\otimes 2s+1}))=\int_{SU_2}\chi^{2s+1}=0$$

(7) Thus, we proved our claim. Now since each irreducible representation of $SU_2$ appears into some $u^{\otimes k}$, and we know how to decompose each $u^{\otimes k}$ into sums of representations $r_k$, these representations $r_k$ are all the irreducible representations of $SU_2$, and we are done with the main assertion. As for the dimension formula, this is clear.
\end{proof}

Regarding now $SO_3$, we have here a similar result, as follows:

\begin{theorem}
The irreducible representations of $SO_3$ are all self-adjoint, and can be labelled by positive integers, with their fusion rules being as follows,
$$r_k\otimes r_l=r_{|k-l|}+r_{|k-l|+1}+\ldots+r_{k+l}$$ 
also called Clebsch-Gordan rules. The corresponding dimensions are $\dim r_k=2k+1$.
\end{theorem}

\begin{proof}
As before with $SU_2$, there are many possible proofs here, which are all instructive. Here is our take on the subject, in the spirit of our proof for $SU_2$:

\medskip

(1) Our first claim is that we have the following formula, telling us that the moments of the main character equal the Catalan numbers:
$$\int_{SO_3}\chi^k=C_k$$

But this is something that we know from $SU_2$, coming via Euler-Rodrigues. Alternatively, this can be deduced as well from Tannakian duality, a bit as for $SU_2$.

\medskip

(2) Our claim now is that we can construct, by recurrence on $k\in\mathbb N$, a sequence $r_k$ of irreducible, self-adjoint and distinct representations of $SO_3$, satisfying:
$$r_0=1\quad,\quad
r_1=u-1\quad,\quad 
r_k+r_{k-1}+r_{k-2}=r_{k-1}\otimes r_1$$

Indeed, assume that $r_0,\ldots,r_{k-1}$ are constructed, and let us construct $r_k$. The Frobenius trick from the proof for $SU_2$ will no longer work, due to some technical reasons, so we have to invoke (1). To be more precise, by integrating characters we obtain:
$$r_{k-1},r_{k-2}\subset r_{k-1}\otimes r_1$$

We conclude that there exists a representation $r_k$ such that:
$$r_{k-1}\otimes r_1=r_k+r_{k-1}+r_{k-2}$$

(3) Once again by integrating characters, we conclude that $r_k$ is irreducible, and non-equivalent to $r_1,\ldots,r_{k-1}$, and this proves our claim. Also, since any irreducible representation of $SO_3$ must appear in some tensor power of $u$, and we can decompose each $u^{\otimes k}$ into sums of representations $r_p$, we conclude that these representations $r_p$ are all the irreducible representations of $SO_3$. Finally, the dimension formula is clear.
\end{proof}

\section*{10e. Exercises}

This was a quite technical algebra chapter, and as exercises on this, we have:

\begin{exercise}
Learn more about categories in general, in mathematics.
\end{exercise}

\begin{exercise}
Do some computations for the diagonal groups, $G\subset\mathbb Z_2^N$.
\end{exercise}

\begin{exercise}
Explore a bit the structure of the diagonal algebras $(C_{kk})$.
\end{exercise}

\begin{exercise}
Learn more about von Neumann's bicommutant theorem.
\end{exercise}

\begin{exercise}
Learn about the other available proofs of Tannakian duality.
\end{exercise}

\begin{exercise}
Learn also about the results of Doplicher-Roberts, and Deligne.
\end{exercise}

\begin{exercise}
Try proving the various Brauer theorems, with bare hands.
\end{exercise}

\begin{exercise}
Explore as well the Brauer theorems for $SU_2$ and for $SO_3$.
\end{exercise}

As bonus exercise, actually good for myself too, read a category theory book.

\chapter{Easy groups}

\section*{11a. Easy groups}

We have seen in the previous chapter that the Tannakian duals of the basic finite groups $S_N,H_N,K_N$, and of the basic continuous groups $O_N,U_N$ too, are very simple objects. To be more precise, the Brauer theorems for these various groups state that we have equalities as follows, with $D=P,P_{even},{\mathcal P}_{even}$ and $D=P_2,\mathcal P_2$ respectively:
$$Hom(u^{\otimes k},u^{\otimes l})=span\left(T_\pi\Big|\pi\in D(k,l)\right)$$

Our goal here will be that of axiomatizing and studying the closed subgroups $G\subset U_N$, finite or not, which are of this type, with $D$ being allowed to be a general category of partitions. We will call such groups ``easy'', and our results will be as follows:

\bigskip 

-- At the level of continuous examples, we will see that besides $O_N,U_N$, we have for instance the bistochastic groups $B_N,C_N$. Which is something quite interesting. 

\medskip

-- At the level of discrete examples, we already know that the reflection groups $H_N^s$ are all easy, but we will provide here a more conceptual proof for this fact.

\medskip

-- Some key Lie groups such as $SU_2,SO_3$, or the symplectic group $Sp_N$, are not easy, but these are however covered by a suitable ``super-easiness'' formalism.

\medskip

-- At the level of the general theory, we will develop some algebraic theory in this chapter, for the most in relation with various product operations.

\medskip

-- Also at the level of the general theory, we will develop as well some analytic theory, later in Part IV, based on the same idea, computations with partitions.

\bigskip

All this sounds quite exciting, good theory that we will be developing here, hope you agree with me. Let us summarize this discussion as a general principle, as follows:

\begin{principle}
Most compact groups $G\subset U_N$, finite or not, are easy, or not far from being easy. And partitions and easiness are the key to their advanced study.
\end{principle}

Getting started now, with the notion of easiness being basically a reorganization of various things that we know from chapter 10, we intend to go quite slow, with the aim of turning this reorganization into a new start, as self-contained as possible. 

\bigskip

We will first need the following key definition, that we know well from chapter 10, which was shown there to appear in relation with the representation theory of the basic compact groups $G\subset U_N$, both finite and infinite, such as $S_N,H_N,K_N$, or $O_N,U_N$:

\index{Kronecker symbols}
\index{maps associated to partitions}

\begin{definition}
Given a partition $\pi\in P(k,l)$ and an integer $N\in\mathbb N$, we define
$$T_\pi:(\mathbb C^N)^{\otimes k}\to(\mathbb C^N)^{\otimes l}$$
by the following formula, with $e_1,\ldots,e_N$ being the standard basis of $\mathbb C^N$,
$$T_\pi(e_{i_1}\otimes\ldots\otimes e_{i_k})=\sum_{j_1\ldots j_l}\delta_\pi\begin{pmatrix}i_1&\ldots&i_k\\ j_1&\ldots&j_l\end{pmatrix}e_{j_1}\otimes\ldots\otimes e_{j_l}$$
and with the coefficients on the right being Kronecker type symbols.
\end{definition}

To be more precise here, in order to compute the Kronecker type symbols $\delta_\pi(^i_j)\in\{0,1\}$, we proceed exactly as in the pairing case, namely by putting the multi-indices $i=(i_1,\ldots,i_k)$ and $j=(j_1,\ldots,j_l)$ on the legs of $\pi$, in the obvious way. In case all the blocks of $\pi$ contain equal indices of $i,j$, we set $\delta_\pi(^i_j)=1$. Otherwise, we set $\delta_\pi(^i_j)=0$.

\bigskip

With the above notion in hand, we can now formulate the following key definition, motivated by the Brauer theorems for $S_N,H_N,K_N$ and $O_N,U_N$, as indicated before:

\index{easy group}

\begin{definition}
A closed subgroup $G\subset U_N$ is called easy when
$$Hom(u^{\otimes k},u^{\otimes l})=span\left(T_\pi\Big|\pi\in D(k,l)\right)$$
for any two colored integers $k,l=\circ\bullet\circ\bullet\ldots\,$, for certain sets of partitions
$$D(k,l)\subset P(k,l)$$
where $\pi\to T_\pi$ is the standard implementation of the partitions, as linear maps.
\end{definition}

In other words, we call a group $G$ easy when its Tannakian category appears in the simplest possible way, from the linear maps associated to partitions. The terminology is quite natural, because Tannakian duality is basically our only serious tool. See \cite{bsp}.

\bigskip

As basic examples, comimg from the above-mentioned Brauer theorems, we have:

\begin{fact}
The following subgroups $G\subset U_N$ are easy:
\begin{enumerate}
\item The groups $G=S_N,H_N,K_N$, coming from $D=P,P_{even},{\mathcal P}_{even}$.

\item The continuous groups $G=O_N,U_N$, coming from $D=P_2,\mathcal P_2$.
\end{enumerate}
\end{fact}

We will actually come back to this later in this chapter, with new proofs for these results, more conceptual than the proofs from chapter 10, which were a bit ad-hoc.

\bigskip

In the general case now, as an important theoretical remark, in the context of Definition 11.3, consider the following collection of sets of partitions:
$$D=\bigsqcup_{k,l}D(k,l)$$

This collection of sets $D$ obviously determines $G$, but the converse is not true. Indeed, at $N=1$ for instance, both the choices $D=P_2,\mathcal P_2$ produce the same easy group, namely $G=\{1\}$. We will be back to this issue on several occasions, with results about it. 

\bigskip

In order to advance, our first goal will be that of establishing a duality between easy groups and certain special classes of collections of sets as above. Let us start with:

\index{category of partitions}

\begin{definition}
Let $P(k,l)$ be the set of partitions between an upper colored integer $k$, and a lower colored integer $l$. A collection of subsets 
$$D=\bigsqcup_{k,l}D(k,l)$$
with $D(k,l)\subset P(k,l)$ is called a category of partitions when it has the following properties:
\begin{enumerate}
\item Stability under the horizontal concatenation, $(\pi,\sigma)\to[\pi\sigma]$.

\item Stability under vertical concatenation $(\pi,\sigma)\to[^\sigma_\pi]$, with matching middle symbols.

\item Stability under the upside-down turning $*$, with switching of colors, $\circ\leftrightarrow\bullet$.

\item Each set $P(k,k)$ contains the identity partition $||\ldots||$.

\item The sets $P(\emptyset,\circ\bullet)$ and $P(\emptyset,\bullet\circ)$ both contain the semicircle $\cap$.

\item The sets $P(k,\bar{k})$ with $|k|=2$ contain the crossing partition $\slash\hskip-2.0mm\backslash$.
\end{enumerate}
\end{definition} 

Observe the similarity with the axioms for Tannakian categories, from chapter 10. We will see in a moment that this similarity can be turned into something very precise, the idea being that such a category produces a family of easy quantum groups $(G_N)_{N\in\mathbb N}$, one for each $N\in\mathbb N$, via the formula in Definition 11.2, and Tannakian duality.

\bigskip

As basic examples of such categories, that we already met in Fact 11.4, we have the categories of partitions $P,P_{even},\mathcal P_{even}$, related to the representation theory of $S_N,H_N,K_N$, as well as the categories of pairings $P_2,\mathcal P_2$, related to the representation theory of $O_N,U_N$. We will see in what follows that there are many other examples.
 
\bigskip

The relation with the Tannakian categories comes from the following result:

\begin{proposition}
The assignement $\pi\to T_\pi$ is categorical, in the sense that
$$T_\pi\otimes T_\sigma=T_{[\pi\sigma]}\quad,\quad 
T_\pi T_\sigma=N^{c(\pi,\sigma)}T_{[^\sigma_\pi]}\quad,\quad 
T_\pi^*=T_{\pi^*}$$
where $c(\pi,\sigma)$ are certain integers, coming from the erased components in the middle.
\end{proposition}

\begin{proof}
This is something that we know well from chapter 10, with the concatenation axiom for instance coming from the following computation:
\begin{eqnarray*}
&&(T_\pi\otimes T_\sigma)(e_{i_1}\otimes\ldots\otimes e_{i_p}\otimes e_{k_1}\otimes\ldots\otimes e_{k_r})\\
&=&\sum_{j_1\ldots j_q}\sum_{l_1\ldots l_s}\delta_\pi\begin{pmatrix}i_1&\ldots&i_p\\j_1&\ldots&j_q\end{pmatrix}\delta_\sigma\begin{pmatrix}k_1&\ldots&k_r\\l_1&\ldots&l_s\end{pmatrix}e_{j_1}\otimes\ldots\otimes e_{j_q}\otimes e_{l_1}\otimes\ldots\otimes e_{l_s}\\
&=&\sum_{j_1\ldots j_q}\sum_{l_1\ldots l_s}\delta_{[\pi\sigma]}\begin{pmatrix}i_1&\ldots&i_p&k_1&\ldots&k_r\\j_1&\ldots&j_q&l_1&\ldots&l_s\end{pmatrix}e_{j_1}\otimes\ldots\otimes e_{j_q}\otimes e_{l_1}\otimes\ldots\otimes e_{l_s}\\
&=&T_{[\pi\sigma]}(e_{i_1}\otimes\ldots\otimes e_{i_p}\otimes e_{k_1}\otimes\ldots\otimes e_{k_r})
\end{eqnarray*}

As for the proofs of the composition axiom and of the involution axiom, these come via similar computations, and we refer here to chapter 10 for the details.
\end{proof}

Time now to put everyting together. We have the following result:

\begin{theorem}
Each category of partitions $D=(D(k,l))$ produces a family of compact groups $G=(G_N)$, one for each $N\in\mathbb N$, via the formula
$$Hom(u^{\otimes k},u^{\otimes l})=span\left(T_\pi\Big|\pi\in D(k,l)\right)$$
and the Tannakian duality correspondence.
\end{theorem}

\begin{proof}
Given an integer $N\in\mathbb N$, consider the correspondence $\pi\to T_\pi$ constructed in Definition 11.2, and then the collection of linear spaces in the statement, namely:
$$C_{kl}=span\left(T_\pi\Big|\pi\in D(k,l)\right)$$

According to the formulae in Proposition 11.6, and to our axioms for the categories of partitions, from Definition 11.5, this collection of spaces $C=(C_{kl})$ satisfies the axioms for the Tannakian categories, from chapter 10. Thus the Tannakian duality result there applies, and provides us with a closed subgroup $G_N\subset U_N$ such that:
$$C_{kl}=Hom(u^{\otimes k},u^{\otimes l})$$

Thus, we are led to the conclusion in the statement.
\end{proof}

In relation with the easiness property, we can now formulate a key result, which can serve as an alternative definition for the easy groups, as follows:

\index{easy group}

\begin{theorem}
A closed subgroup $G\subset U_N$ is easy precisely when
$$Hom(u^{\otimes k},u^{\otimes l})=span\left(T_\pi\Big|\pi\in D(k,l)\right)$$
for any colored integers $k,l$, for a certain category of partitions $D\subset P$.
\end{theorem}

\begin{proof}
This basically follows from Theorem 11.7, as follows:

\medskip

(1) In one sense, we know from Theorem 11.7 that any category of partitions $D\subset P$ produces a family of closed groups $G\subset U_N$, one for each $N\in\mathbb N$, according to Tannakian duality and to the Hom space formula there, namely:
$$Hom(u^{\otimes k},u^{\otimes l})=span\left(T_\pi\Big|\pi\in D(k,l)\right)$$

But these groups $G\subset U_N$ are indeed easy, in the sense of Definition 11.3.

\medskip

(2) In the other sense now, assume that $G\subset U_N$ is easy, in the sense of Definition 11.3, coming via the above Hom space formula, from a collection of sets as follows:
$$D=\bigsqcup_{k,l}D(k,l)$$

Consider now the category of partitions $\widetilde{D}=<D>$ generated by this family. This is by definition the smallest category of partitions containing $D$, whose existence follows by starting with $D$, and performing the various categorical operations, namely horizontal and vertical concatenation, and upside-down turning. It follows then, via another application of Tannakian duality, that we have the following formula, for any $k,l$:
$$Hom(u^{\otimes k},u^{\otimes l})=span\left(T_\pi\Big|\pi\in\widetilde{D}(k,l)\right)$$

Thus, our group $G\subset U_N$ can be viewed as well as coming from $\widetilde{D}$, and so appearing as particular case of the construction in Theorem 11.7, and this gives the result. 
\end{proof}

As already mentioned above, Theorem 11.8 can be regarded as an alternative definition for easiness, with the assumption that $D\subset P$ must be a category of partitions being added. In what follows we will rather use this new definition, which is more precise.

\bigskip

We will see in what follows that many interesting examples of compact quantum groups are easy.  Moreover, most of the known series of ``basic'' compact quantum groups,  $G=(G_N)$ with $N\in\mathbb N$, can be in principle made fit into some suitable extensions of the easy quantum group formalism. We will discuss this too, in what follows.

\bigskip

Getting now to the basic examples of easy groups, these are the groups $S_N,H_N,K_N$ and $O_N,U_N$ from Fact 11.4. However, we have promised to come here with new proofs for all this, a bit more conceptual. In what regards $O_N,U_N$, the result is as follows:

\index{Brauer theorem}
\index{easy group}
\index{orthogonal group}
\index{unitary group}

\begin{theorem}
We have the following results:
\begin{enumerate}
\item The unitary group $U_N$ is easy, coming from the category $\mathcal P_2$.

\item The orthogonal group $O_N$ is easy as well, coming from the category $P_2$.
\end{enumerate}
\end{theorem}

\begin{proof}
This comes indeed from Tannakian duality, coupled with some knowledge of the maps $T_\pi$, for the basic partitions $\pi$, that we have since chapter 10, as follows:

\medskip

(1) The group $U_N$ being defined via the relations $u^*=u^{-1}$, $u^t=\bar{u}^{-1}$, we conclude that the associated Tannakian category must be $C=span(T_\pi|\pi\in D)$, with:
$$D
=<{\ }^{\,\cap}_{\circ\bullet}\,\,,{\ }^{\,\cap}_{\bullet\circ}>
=\mathcal P_2$$

(2) As for the subgroup $O_N\subset U_N$, this being defined by imposing the relations $u_{ij}=\bar{u}_{ij}$, the associated Tannakian category must be $C=span(T_\pi|\pi\in D)$, with:
$$D
=<\mathcal P_2,|^{\hskip-1.32mm\circ}_{\hskip-1.32mm\bullet},|_{\hskip-1.32mm\circ}^{\hskip-1.32mm\bullet}>
=P_2$$
  
Thus, we are led to the conclusions in the statement.
\end{proof}

Still in the continuous case, and coming as something truly new, we have:

\index{Brauer theorem}
\index{easy group}
\index{bistochastic group}

\begin{theorem}
We have the following results:
\begin{enumerate}
\item The unitary bistochastic group $C_N$ is easy, coming from the category $\mathcal P_{12}$ of matching singletons and pairings.

\item The orthogonal bistochastic group $B_N$ is easy too, coming from the category $P_{12}$ of singletons and pairings.
\end{enumerate}
\end{theorem}

\begin{proof}
The proof here is quite similar to the proof of Theorem 11.9. To be more precise, we can use the results there, and the proof goes as follows:

\medskip

(1) The group $C_N\subset U_N$ is defined by imposing the following relations, with $\xi$ being the all-one vector, which correspond to the bistochasticity condition:
$$u\xi=\xi\quad,\quad 
\bar{u}\xi=\xi$$

But these relations tell us precisely that the following two operators, with the partitions on the right being singletons, must be in the associated Tannakian category $C$:
$$T_\pi\quad:\quad\pi=|_{\hskip-1.32mm\circ}\ ,\ |_{\hskip-1.32mm\bullet}$$

Thus the associated Tannakian category is $C=span(T_\pi|\pi\in D)$, with:
$$D
=<\mathcal P_2,|_{\hskip-1.32mm\circ},|_{\hskip-1.32mm\bullet}>
=\mathcal P_{12}$$

(2) In order to deal now with the real bistochastic group $B_N$, we can either use a similar argument, or simply use the following intersection formula: 
$$B_N=C_N\cap O_N$$

Indeed, at the categorical level, this intersection formula tells us that the associated Tannakian category is given by $C=span(T_\pi|\pi\in D)$, with:
$$D
=<\mathcal P_{12},P_2>
=P_{12}$$

Thus, we are led to the conclusions in the statement.
\end{proof}

As a comment here, we have used in the above the fact, which is something quite trivial, that the category of partitions associated to an intersection of easy quantum groups is generated by the corresponding categories of partitions. We will be back to this, and to some other product operations as well, with similar results, later on.

\bigskip

We can put now the results that we have together, as follows:

\begin{theorem}
The basic unitary and bistochastic groups,
$$\xymatrix@R=50pt@C=50pt{
C_N\ar[r]&U_N\\
B_N\ar[u]\ar[r]&O_N\ar[u]}$$
are all easy, coming from the various categories of singletons and pairings.
\end{theorem}

\begin{proof}
We know from the above that the groups in the statement are indeed easy, the corresponding diagram of categories of partitions being as follows:
$$\xymatrix@R=16mm@C=18mm{
\mathcal P_{12}\ar[d]&\mathcal P_2\ar[l]\ar[d]\\
P_{12}&P_2\ar[l]}$$

Thus, we are led to the conclusion in the statement.
\end{proof}

Summarizing, we have a notion of easiness based on the Brauer theorems for $O_N,U_N$, and their versions for $B_N,C_N$, with everything coming up together quite nicely.

\section*{11b. Reflection groups}

Getting now to finite group case, we know from chapter 10 that the complex reflection groups $H_N^s=\mathbb Z_s\wr S_N$ are all easy. However, the proofs there were a bit ad-hoc, and all this is worth a slow revisit. In what regards the symmetric group $S_N$, we have:

\index{Brauer theorem}
\index{symmetric group}

\begin{theorem}
The symmetric group $S_N$, regarded as group of unitary matrices,
$$S_N\subset O_N\subset U_N$$
via the permutation matrices, is easy, coming from the category of all partitions $P$.
\end{theorem}

\begin{proof}
Consider indeed the group $S_N$, regarded as a group of unitary matrices, with each permutation $\sigma\in S_N$ corresponding to the associated permutation matrix:
$$\sigma(e_i)=e_{\sigma(i)}$$

In order to prove the result, consider the one-block ``fork'' partition, namely:
$$\xymatrix@R=1mm@C=2mm{\\ \\ \mu\ \ =\\ \\ }\ \ \ 
\xymatrix@R=2mm@C=3mm{
\circ\ar@/_/@{-}[dr]&&\circ\\
&\ar@/_/@{-}[ur]\ar@{-}[dd]\\
&&&\\
&\circ}$$

The linear map associated to $\mu$ is then given by the following formula:
$$T_\mu(e_i\otimes e_j)=\delta_{ij}e_i$$

In order to do the computations, we use the following formulae:
$$u=(u_{ij})_{ij}\quad,\quad 
u^{\otimes 2}=(u_{ij}u_{kl})_{ik,jl}\quad,\quad 
T_\mu=(\delta_{ijk})_{i,jk}$$

By using these formulae, we obtain the following equality:
$$(T_\mu u^{\otimes 2})_{i,jk}
=\sum_{lm}(T_\mu)_{i,lm}(u^{\otimes 2})_{lm,jk}
=u_{ij}u_{ik}$$

On the other hand, we have as well the following equality:
$$(uT_\mu)_{i,jk}
=\sum_lu_{il}(T_\mu)_{l,jk}
=\delta_{jk}u_{ij}$$

We therefore conclude that we have an equivalence, as follows:
$$T_\mu\in Hom(u^{\otimes 2},u)\iff u_{ij}u_{ik}=\delta_{jk}u_{ij},\forall i,j,k$$

In other words, the elements $u_{ij}$ must be projections, which must be pairwise orthogonal on the rows of $u=(u_{ij})$. But this reformulates into the following equality:
$$C(S_N)=C(O_N)\Big/\Big<T_\mu\in Hom(u^{\otimes 2},u)\Big>$$

According now to our general conventions for easiness, this means that the symmetric group $S_N$ is easy, coming from the following category of partitions:
$$D=<\mu>=P$$

Thus, we are led to the conclusion in the statement.
\end{proof}

Next, regarding the hyperoctahedral group $H_N$, we have the following result:

\index{Brauer theorem}
\index{hyperoctahedral group}

\begin{theorem}
The hyperoctahedral group $H_N$, regarded as group of matrices,
$$H_N\subset O_N\subset U_N$$
is easy, coming from the category of partitions with even blocks $P_{even}$.
\end{theorem}

\begin{proof}
This follows as usual from Tannakian duality. To be more precise, consider the following one-block partition $\chi\in P(2,2)$, which looks like a $\chi$ letter:
$$\xymatrix@R=0.5mm@C=2mm{\\ \\ \\ \chi\ \ =\\ \\ }\ \ \ 
\xymatrix@R=2mm@C=3mm{
\circ\ar@/_/@{-}[dr]&&\circ\\
&\ar@/_/@{-}[ur]\ar@{-}[dd]\\
&&&\\
&\ar@/^/@{-}[dr]\ar@/_/@{-}[dl]\\
\circ&&\circ}$$

The linear map associated to this partition is then given by:
$$T_\chi(e_i\otimes e_j)=\delta_{ij}e_i\otimes e_i$$

By using this formula, we have the following computation:
\begin{eqnarray*}
(T_\chi\otimes id)u^{\otimes 2}(e_a\otimes e_b)
&=&(T_\chi\otimes id)\left(\sum_{ijkl}e_{ij}\otimes e_{kl}\otimes u_{ij}u_{kl}\right)(e_a\otimes e_b)\\
&=&(T_\chi\otimes id)\left(\sum_{ik}e_i\otimes e_k\otimes u_{ia}u_{kb}\right)\\
&=&\sum_ie_i\otimes e_i\otimes u_{ia}u_{ib}
\end{eqnarray*}

On the other hand, we have as well the following computation:
\begin{eqnarray*}
u^{\otimes 2}(T_\chi\otimes id)(e_a\otimes e_b)
&=&\delta_{ab}\left(\sum_{ijkl}e_{ij}\otimes e_{kl}\otimes u_{ij}u_{kl}\right)(e_a\otimes e_a)\\
&=&\delta_{ab}\sum_{ij}e_i\otimes e_k\otimes u_{ia}u_{ka}
\end{eqnarray*}

We conclude from this that we have the following equivalence:
$$T_\chi\in End(u^{\otimes 2})\iff \delta_{ik}u_{ia}u_{ib}=\delta_{ab}u_{ia}u_{ka},\forall i,k,a,b$$

But the relations on the right tell us that the entries of $u=(u_{ij})$ must satisfy $\alpha\beta=0$ on each row and column of $u$, and so that the corresponding closed subgroup $G\subset O_N$ consists of the matrices $g\in O_N$ which are permutation-like, with $\pm1$ nonzero entries. Thus, the corresponding group is $G=H_N$, and as a conclusion to this, we have:
$$C(H_N)=C(O_N)\Big/\Big<T_\chi\in End(u^{\otimes 2})\Big>$$

According now to our conventions for easiness, this means that the hyperoctahedral group $H_N$ is easy, coming from the following category of partitions:
$$D=<\chi>=P_{even}$$

Thus, we are led to the conclusion in the statement.
\end{proof}

Next, regarding the full reflection group $K_N$, we have the following result:

\begin{theorem}
The full reflection group $K_N=\mathbb T\wr S_N$, regarded as subgroup
$$K_N\subset U_N$$
comes from $\mathcal P_{even}$, the partitions satisfying $\#\circ=\#\bullet$, weighted equality, in each block. 
\end{theorem}

\begin{proof}
We are now dealing with unitary matrices, so we must use colored partitions. Consider the following partition $\chi\in P(\circ\bullet\,,\bullet\circ)$, that we already met above, uncolored:
$$
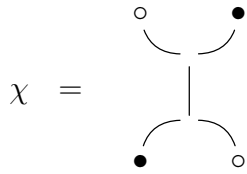
\xymatrix@R=0.5mm@C=2mm{\\ \\ \\ \chi\ \ =\\ \\ }\ \ \ 
\xymatrix@R=2mm@C=3mm{
\circ\ar@/_/@{-}[dr]&&\bullet\\
&\ar@/_/@{-}[ur]\ar@{-}[dd]\\
&&&\\
&\ar@/^/@{-}[dr]\ar@/_/@{-}[dl]\\
\bullet&&\circ}$$

Our computations from the previous proof, for the group $H_N$, modify into:
$$(T_\chi\otimes id)(u\otimes\bar{u})(e_a\otimes e_b)=\sum_ie_i\otimes e_i\otimes u_{ia}\bar{u}_{ib}$$
$$(\bar{u}\otimes u)(T_\chi\otimes id)(e_a\otimes e_b)=\delta_{ab}\sum_{ij}e_i\otimes e_k\otimes\bar{u}_{ia}u_{ka}$$

We conclude from this that we have the following equivalence:
$$T_\chi\in Hom(u\otimes\bar{u},\bar{u}\otimes u)\iff \delta_{ik}u_{ia}\bar{u}_{ib}=\delta_{ab}\bar{u}_{ia}u_{ka},\forall i,k,a,b$$

But the relations on the right tell us that the entries of $u=(u_{ij})$ must satisfy $\alpha\beta=0$ on each row and column of $u$, and as a conclusion to this, we have:
$$C(K_N)=C(U_N)\Big/\Big<T_\chi\in Hom(u\otimes\bar{u},\bar{u}\otimes u)\Big>$$

Thus the group $K_N$ is easy, coming from the following category of partitions:
$$D=<\chi>=\mathcal P_{even}$$

We are therefore led to the conclusion in the statement.
\end{proof}

More generally now, we have in fact the following grand result:

\index{Brauer theorem}
\index{complex reflection group}

\begin{theorem}
The complex reflection group $H_N^s=\mathbb Z_s\wr S_N$ is easy, the corresponding category $P^s$ consisting of the partitions satisfying the condition
$$\#\circ=\#\bullet(s)$$
as a weighted sum, in each block. In particular, we have the following results:
\begin{enumerate}
\item $S_N$ is easy, coming from the category $P$.

\item $H_N=\mathbb Z_2\wr S_N$ is easy, coming from the category $P_{even}$.

\item $K_N=\mathbb T\wr S_N$ is easy, coming from the category $\mathcal P_{even}$.
\end{enumerate}
\end{theorem}

\begin{proof}
This is something coming at $s=1,2,\infty$ from Theorems 11.12, 11.13 and 11.14, as indicated in (1,2,3), with this to be discussed in a moment, and in general, the proof is similar. Consider indeed the following partition, with $s+2$ legs:
$$
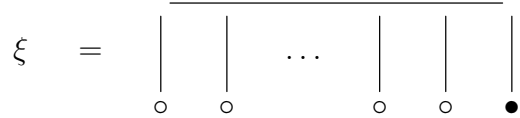
\xymatrix@R=6pt@C=12pt{
&&\ar@{-}[dd]\ar@{-}[rrrrr]&\ar@{-}[dd]&&\ar@{-}[dd]&\ar@{-}[dd]&\ar@{-}[dd]\\
\xi&=&&&\ldots\\
&&\circ&\circ&&\circ&\circ&\bullet}$$

Observe that, up to rotation and some discussion regarding the colors, this coincides with the partitions $\mu,\chi$ that we used before at $s=1,2$. In general now, we have:
$$T_\xi=\sum_je_j^{\otimes s+2}$$

Our claim, which will prove the result, is that we have the following formula:
$$C(H_N^s)=C(K_N)\Big/\Big<T_\xi\in Fix(u^{\otimes s+1}\otimes\bar{u})\Big>$$

Indeed, by using the above formula of $T_\xi$, we have the following computation:
\begin{eqnarray*}
(u^{\otimes s+1}\otimes\bar{u})(T_\xi\otimes 1)
&=&\sum_{ij}e_{i_1}\otimes\ldots\otimes e_{i_{s+2}}\otimes u_{i_1j}\ldots u_{i_{s+1}j}\bar{u}_{i_{s+2}j}\\
&=&\sum_{ij}e_i\otimes\ldots\otimes e_i\otimes u_{ij}^{s+1}\bar{u}_{ij}\\
&=&\sum_ie_i^{\otimes s+2}\otimes\left(\sum_ju_{ij}^{s+1}\bar{u}_{ij}\right)
\end{eqnarray*}

We conclude that, for a subgroup of $K_N$, we have the following equivalence:
$$T_\xi\in Fix(u^{\otimes s+1}\otimes\bar{u})\iff\sum_ju_{ij}^{s+1}\bar{u}_{ij}=1$$

Now the conditions on the right being those defining the subgroup $H_N^s\subset K_N$, we conclude that we have the equality announced above, namely:
$$C(H_N^s)=C(K_N)\Big/\Big<T_\xi\in Fix(u^{\otimes s+1}\otimes\bar{u})\Big>$$

But with this, we can finish the proof of the main assertion. Indeed, it follows that the group $H_N^s$ is easy, coming from the following category of partitions:
$$D=<\mathcal P_{even},\xi>=P^s$$

Summarizing, theorem proved, and in what regards the particular cases, which generalize what we knew from Theorems 11.12, 11.13 and 11.14, these are as follows:

\medskip

(1) At $s=1$ we know that we have $H_N^1=S_N$. Regarding now the corresponding category, here the condition $\#\circ=\#\bullet(1)$ is automatic, and so $P^1=P$.

\medskip

(2) At $s=2$ we know that we have $H_N^2=H_N$. Regarding now the corresponding category, here the condition $\#\circ=\#\bullet(2)$ reformulates as follows:
$$\#\circ+\,\#\bullet=0(2)$$

Thus each block must have even size, and we obtain, as claimed, $P^2=P_{even}$.

\medskip

(3) At $s=\infty$ we know that we have $H_N^\infty=K_N$. Regarding now the corresponding category, here the condition $\#\circ=\#\bullet(\infty)$ reads:
$$\#\circ=\#\bullet$$

But this is the condition defining $\mathcal P_{even}$, and so $P^\infty=\mathcal P_{even}$, as claimed.
\end{proof}

And with this, end of our discussion regarding various Brauer theorems, and subsequent examples of easy groups. Actually, our list of easy groups has currently become quite big, and here is a selection of the main results that we have so far: 

\begin{theorem}
We have a diagram of compact groups as follows,
$$\xymatrix@R=50pt@C=50pt{
K_N\ar[r]&U_N\\
H_N\ar[u]\ar[r]&O_N\ar[u]}$$
where $H_N=\mathbb Z_2\wr S_N$ and $K_N=\mathbb T\wr S_N$, and all these groups are easy.
\end{theorem}

\begin{proof}
This follows from the above results. To be more precise, we know that the above groups are all easy, the corresponding categories of partitions being as follows:
$$\xymatrix@R=16mm@C=18mm{
\mathcal P_{even}\ar[d]&\mathcal P_2\ar[l]\ar[d]\\
P_{even}&P_2\ar[l]}$$

Thus, we are led to the conclusion in the statement.
\end{proof}

Observed that we have removed from the above result our beloved symmetric group $S_N$, and this for a quite subtle reason, coming from the fact that, once we agree that $K_N$ should appear on top left, $H_N$ should appear on bottom left, according to:
$$H_N=K_N\cap O_N$$

In practice, the diagram in Theorem 11.16 is something very useful, providing us with some 2D orientation into the world of easy groups, and of subgroups $G\subset U_N$ in general. In fact, we have already met this philosophy in chapter 5, and more about this later.

\section*{11c. Basic operations} 

Let us discuss now some basic operations for the closed subgroups $G\subset U_N$, in general, and in the easy case. We will be mainly interested in the following operations:

\begin{definition}
The closed subgroups of $U_N$ are subject to intersection and generation operations, constructed as follows:
\begin{enumerate}
\item Intersection: $H\cap K$ is the usual intersection of $H,K$.

\item Generation: $<H,K>$ is the closed subgroup generated by $H,K$.
\end{enumerate}
\end{definition}

Alternatively, we can define these operations at the function algebra level, by performing certain operations on the associated ideals, as follows:

\begin{proposition}
Assuming that we have presentation results as follows,
$$C(H)=C(U_N)/I\quad,\quad 
C(K)=C(U_N)/J$$
the groups $H\cap K$ and $<H,K>$ are given by the following formulae,
$$C(H\cap K)=C(U_N)/<I,J>$$
$$C(<H,K>)=C(U_N)/(I\cap J)$$
at the level of the associated algebras of functions.
\end{proposition}

\begin{proof}
This is indeed clear from the definition of the operations $\cap$ and $<\,,\,>$, as formulated above, and from the Stone-Weierstrass theorem.
\end{proof}

In what follows we will need Tannakian formulations of the above two operations. The result here, that we have already used a couple of times in the above, is as follows:

\begin{theorem}
The intersection and generation operations $\cap$ and $<\,,>$ can be constructed via the Tannakian correspondence $G\to C_G$, as follows:
\begin{enumerate}
\item Intersection: defined via $C_{G\cap H}=<C_G,C_H>$.

\item Generation: defined via $C_{<G,H>}=C_G\cap C_H$.
\end{enumerate}
\end{theorem}

\begin{proof}
This follows from Proposition 11.18, and from Tannakian duality. Indeed, it follows from Tannakian duality that given a closed subgroup $G\subset U_N$, with fundamental representation $v$, the algebra of functions $C(G)$ has the following presentation:
$$C(G)=C(U_N)\Big/\left<T\in Hom(u^{\otimes k},u^{\otimes l})\Big|\forall k,\forall l,\forall T\in Hom(v^{\otimes k},v^{\otimes l})\right>$$

In other words, given a closed subgroup $G\subset U_N$, we have a presentation of the following type, with $I_G$ being the ideal coming from the Tannakian category of $G$:
$$C(G)=C(U_N)/I_G$$

But this leads to the conclusion in the statement.
\end{proof}

In relation now with our easiness questions, we first have the following result:

\begin{proposition}
Assuming that $H,K$ are easy, then so is $H\cap K$, and we have
$$D_{H\cap K}=<D_H,D_K>$$
at the level of the corresponding categories of partitions.
\end{proposition}

\begin{proof}
We have indeed the following elementary computation:
\begin{eqnarray*}
C_{H\cap K}
&=&<C_H,C_K>\\
&=&<span(D_H),span(D_K)>\\
&=&span(<D_H,D_K>)
\end{eqnarray*}

Thus, by Tannakian duality we obtain the result.
\end{proof}

Regarding now the generation operation, the situation here is more complicated, due to a number of technical reasons, and we only have the following statement:

\begin{proposition}
Assuming that $H,K$ are easy, we have an inclusion 
$$<H,K>\subset\{H,K\}$$
coming from an inclusion of Tannakian categories as follows,
$$C_H\cap C_K\supset span(D_H\cap D_K)$$
where $\{H,K\}$ is the easy group having as category of partitions $D_H\cap D_K$.
\end{proposition}

\begin{proof}
This follows from the definition and properties of the generation operation, explained above, and from the following computation:
\begin{eqnarray*}
C_{<H,K>}
&=&C_H\cap C_K\\
&=&span(D_H)\cap span(D_K)\\
&\supset&span(D_H\cap D_K)
\end{eqnarray*}

Indeed, by Tannakian duality we obtain from this all the assertions.
\end{proof}

Summarizing, we have some problems here, and we must proceed as follows:

\index{easy generation}

\begin{theorem}
The intersection and easy generation operations $\cap$ and $\{\,,\}$ can be constructed via the Tannakian correspondence $G\to D_G$, as follows:
\begin{enumerate}
\item Intersection: defined via $D_{G\cap H}=<D_G,D_H>$.

\item Easy generation: defined via $D_{\{G,H\}}=D_G\cap D_H$.
\end{enumerate}
\end{theorem}

\begin{proof}
Here the situation is as follows:

\medskip

(1) This is a true and honest result, coming from Proposition 11.20.

\medskip

(2) This is more of an empty statement, coming from Proposition 11.21.
\end{proof}

We will be back to these questions later, with some further details and results. With the above notions in hand, however, even if not fully satisfactory, we can formulate a nice result, which improves our main result so far, namely Theorem 11.16, as follows:

\begin{theorem}
The basic unitary and reflection groups, namely
$$\xymatrix@R=50pt@C=50pt{
K_N\ar[r]&U_N\\
H_N\ar[u]\ar[r]&O_N\ar[u]}$$
are all easy, and they form an intersection and easy generation diagram, in the sense that the above square diagram satisfies $H_N=K_N\cap O_N$ and $U_N=\{K_N,O_N\}$.
\end{theorem}

\begin{proof}
We know from Theorem 11.16 that the groups in the statement are all easy, the corresponding categories of partitions being as follows:
$$\xymatrix@R=16mm@C=18mm{
\mathcal P_{even}\ar[d]&\mathcal P_2\ar[l]\ar[d]\\
P_{even}&P_2\ar[l]}$$

Now observe that this latter diagram is an intersection and generation diagram. By using Theorem 11.22, this reformulates into the fact that the diagram of groups in the statement is an intersection and easy generation diagram, as claimed.
\end{proof}

It is possible to further improve the above result, by proving that the diagram there is actually a plain generation diagram. However, this is something more technical, and for a discussion here, you can check for instance my quantum group book \cite{ba2}.

\bigskip

Moving forward, as a continuation of the above, it is possible to develop some more general theory, along the above lines. Let us start with the following notion:

\index{homogeneous group}

\begin{definition}
A closed subgroup $G\subset U_N$ is called homogeneous when
$$S_N\subset G\subset U_N$$
with $S_N\subset U_N$ being the standard embedding, via permutation matrices.
\end{definition}

We will be interested in such groups, which cover for instance all the easy groups, and many more. At the Tannakian level, we have the following result:

\begin{theorem}
The homogeneous groups $S_N\subset G\subset U_N$ are in one-to-one correspondence with the intermediate tensor categories
$$span\left(T_\pi\Big|\pi\in\mathcal P_2\right)\subset C\subset span\left(T_\pi\Big|\pi\in P\right)$$
where $P$ is the category of all partitions, $\mathcal P_2$ is the category of the matching pairings, and $\pi\to T_\pi$ is the standard implementation of partitions, as linear maps.
\end{theorem}

\begin{proof}
This follows from Tannakian duality, and from the Brauer type results for $S_N,U_N$. To be more precise, we know from Tannakian duality that each closed subgroup $G\subset U_N$ can be reconstructed from its Tannakian category $C=(C_{kl})$, as follows:
$$C(G)=C(U_N)\Big/\left<T\in Hom(u^{\otimes k},u^{\otimes l})\Big|\forall k,l,\forall T\in C_{kl}\right>$$

Thus we have a one-to-one correspondence $G\leftrightarrow C$, given by Tannakian duality, and since the endpoints $G=S_N,U_N$ are both easy, corresponding to the categories $C=span(T_\pi|\pi\in D)$ with $D=P,\mathcal P_2$, this gives the result.
\end{proof}

Coming next, and making the link with easiness, we have:

\begin{proposition}
Given a homogeneous group $S_N\subset G\subset U_N$, with associated Tannakian category $C=(C_{kl})$, the sets
$$D^1(k,l)=\left\{\pi\in P(k,l)\Big|T_\pi\in C_{kl}\right\}$$ 
form a category of partitions, in the sense of Definition 11.5.
\end{proposition}

\begin{proof}
We use the basic categorical properties of the correspondence $\pi\to T_\pi$ between partitions and linear maps, that we established in the above, namely:
$$T_{[\pi\sigma]}=T_\pi\otimes T_\sigma\quad,\quad 
T_{[^\sigma_\pi]}\sim T_\pi T_\sigma\quad,\quad 
T_{\pi^*}=T_\pi^*$$

Together with the fact that $C$ is a tensor category, we deduce from these formulae that we have the following implication:
\begin{eqnarray*}
\pi,\sigma\in D^1
&\implies&T_\pi,T_\sigma\in C\\
&\implies&T_\pi\otimes T_\sigma\in C\\
&\implies&T_{[\pi\sigma]}\in C\\
&\implies&[\pi\sigma]\in D^1
\end{eqnarray*}

On the other hand, we have as well the following implication:
\begin{eqnarray*}
\pi,\sigma\in D^1
&\implies&T_\pi,T_\sigma\in C\\
&\implies&T_\pi T_\sigma\in C\\
&\implies&T_{[^\sigma_\pi]}\in C\\
&\implies&[^\sigma_\pi]\in D^1
\end{eqnarray*}

Finally, we have as well the following implication:
\begin{eqnarray*}
\pi\in D^1
&\implies&T_\pi\in C\\
&\implies&T_\pi^*\in C\\
&\implies&T_{\pi^*}\in C\\
&\implies&\pi^*\in D^1
\end{eqnarray*}

Thus $D^1$ is indeed a category of partitions, as claimed.
\end{proof}

We can further refine the above observation, in the following way:

\begin{proposition}
Given a compact group $S_N\subset G\subset U_N$, construct $D^1\subset P$ as above, and let $S_N\subset G^1\subset U_N$ be the easy group associated to $D^1$. Then:
\begin{enumerate}

\item We have $G\subset G^1$, as subgroups of $U_N$.

\item $G^1$ is the smallest easy group containing $G$.

\item $G$ is easy precisely when $G\subset G^1$ is an isomorphism.
\end{enumerate}
\end{proposition}

\begin{proof}
All this is elementary, the proofs being as follows:

\medskip

(1) We know that the Tannakian category of $G^1$ is given by:
$$C_{kl}^1=span\left(T_\pi\Big|\pi\in D^1(k,l)\right)$$

Thus we have $C^1\subset C$, and so $G\subset G^1$, as subgroups of $U_N$.

\medskip

(2) Assuming that we have $G\subset G'$, with $G'$ easy, coming from a Tannakian category $C'=span(D')$, we must have $C'\subset C$, and so $D'\subset D^1$. Thus, $G^1\subset G'$, as desired.

\medskip

(3) This is a trivial consequence of (2).
\end{proof}

Summarizing, we have now a notion of ``easy envelope'', as follows:

\index{easy envelope}

\begin{definition}
The easy envelope of a homogeneous group $S_N\subset G\subset U_N$ is the easy group $S_N\subset G^1\subset U_N$ associated to the category of partitions
$$D^1(k,l)=\left\{\pi\in P(k,l)\Big|T_\pi\in C_{kl}\right\}$$ 
where $C=(C_{kl})$ is the Tannakian category of $G$.
\end{definition}

At the level of examples, most of the known homogeneous groups $S_N\subset G\subset U_N$ are in fact easy. However, there are non-easy interesting examples as well, as follows:

\bigskip

-- First we have the subgroup $U_N^d\subset U_N$ constructed via the relation $(\det U)^d=1$, that we already met in chapter 7, whose easy envelope is $U_N$.

\bigskip

-- Then, we have the subgroup $H_N^{s,d}\subset H_N^s$ constructed via exactly the same relation, $(\det U)^d=1$, that we know since chapter 5, whose easy envelope is $H_N^s$.

\bigskip

-- With the remark however that, in relation with this latter fact, there is an exception, namely the group $H_2^{4,2}$ which is easy, and so equal to its easy envelope.

\bigskip

-- Summarizing, we have examples. As for the proofs, all this follows either directly, or by using standard classification results for the easy groups, from \cite{twe}.

\bigskip

As a technical observation now, we can in fact generalize the above construction to any closed subgroup $G\subset U_N$, and we have the following result:

\begin{proposition}
Given a closed subgroup $G\subset U_N$, construct $D^1\subset P$ as above, and let $S_N\subset G^1\subset U_N$ be the easy group associated to $D^1$. We have then
$$G^1=(<G,S_N>)^1$$
where $<G,S_N>\subset U_N$ is the smallest closed subgroup containing $G,S_N$.
\end{proposition}

\begin{proof}
According to our Tannakian results, the subgroup $<G,S_N>\subset U_N$ in the statement exists indeed, and can be obtained by intersecting categories, as follows:
$$C_{<G,S_N>}=C_G\cap C_{S_N}$$

We conclude from this that for any $\pi\in P(k,l)$ we have:
$$T_\pi\in C_{<G,S_N>}(k,l)\iff T_\pi\in C_G(k,l)$$

It follows that the $D^1$ categories for the groups $<G,S_N>$ and $G$ coincide, and so the easy envelopes $(<G,S_N>)^1$ and $G^1$ coincide as well, as stated.
\end{proof}

In order now to fine-tune all this, by using an arbitrary parameter $p\in\mathbb N$, which can be thought of as being an ``easiness level'', we can proceed as follows:

\index{easiness level}

\begin{definition}
Given a compact group $S_N\subset G\subset U_N$, and an integer $p\in\mathbb N$, we construct the family of linear spaces
$$E^p(k,l)=\left\{\alpha_1T_{\pi_1}+\ldots+\alpha_pT_{\pi_p}\in C(k,l)\Big|\alpha_i\in\mathbb C,\pi_i\in P(k,l)\right\}$$
and we denote by $C^p$ the smallest tensor category containing $E^p=(E^p(k,l))$, and by $S_N\subset G^p\subset U_N$ the compact group corresponding to this category $C^p$.
\end{definition}

As a first observation, at $p=1$ we have $C^1=E^1=span(D^1)$, where $D^1$ is the category of partitions constructed in Proposition 11.27. Thus the group $G^1$ constructed above coincides with the ``easy envelope'' of $G$, from Definition 11.28.

\bigskip

In the general case, $p\in\mathbb N$, the family $E^p=(E^p(k,l))$ constructed above is not necessarily a tensor category, but we can of course consider the tensor category $C^p$ generated by it, as indicated. Finally, in the above definition we have used of course the Tannakian duality results, in order to perform the operation $C^p\to G^p$.

\bigskip

In practice, the construction in Definition 11.30 is often something quite complicated, and it is convenient to use the following observation:

\begin{proposition}
The category $C^p$ constructed above is generated by the spaces
$$E^p(l)=\left\{\alpha_1T_{\pi_1}+\ldots+\alpha_pT_{\pi_p}\in C(l)\Big|\alpha_i\in\mathbb C,\pi_i\in P(l)\right\}$$
where $C(l)=C(0,l),P(l)=P(0,l)$, with $l$ ranging over the colored integers.
\end{proposition}

\begin{proof}
We use the well-known fact, that we know from chapter 9, that given a closed subgroup $G\subset U_N$, we have a Frobenius type isomorphism, as follows:
$$Hom(u^{\otimes k},u^{\otimes l})\simeq Fix(u^{\otimes\bar{k}l})$$

If we apply this to the group $G^p$, we obtain an isomorphism as follows:
$$C(k,l)\simeq C(\bar{k}l)$$

On the other hand, we have as well an isomorphism $P(k,l)\simeq P(\bar{k}l)$, obtained by performing a counterclockwise rotation to the partitions $\pi\in P(k,l)$. According to the above definition of the spaces $E^p(k,l)$, this induces an isomorphism as follows:
$$E^p(k,l)\simeq E^p(\bar{k}l)$$

We deduce from this that for any partitions $\pi_1,\ldots,\pi_p\in C(k,l)$, having rotated versions $\rho_1,\ldots,\rho_p\in C(\bar{k}l)$, and for any scalars $\alpha_1,\ldots,\alpha_p\in\mathbb C$, we have:
$$\alpha_1T_{\pi_1}+\ldots+\alpha_pT_{\pi_p}\in C(k,l)\iff\alpha_1T_{\rho_1}+\ldots+\alpha_pT_{\rho_p}\in C(\bar{k}l)$$

But this gives the conclusion in the statement, and we are done.
\end{proof}

The main properties of the construction $G\to G^p$ can be summarized as follows:

\begin{theorem}
Given a compact group $S_N\subset G\subset U_N$, the compact groups $G^p$ constructed above form a decreasing family, whose intersection is $G$:
$$G=\bigcap_{p\in\mathbb N}G^p$$
Moreover, $G$ is easy when this decreasing limit is stationary, $G=G^1$.
\end{theorem}

\begin{proof}
By definition of $E^p(k,l)$, and by using Proposition 11.31, these linear spaces form an increasing filtration of $C(k,l)$. The same remains true when completing into tensor categories, and so we have an increasing filtration, as follows:
$$C=\bigcup_{p\in\mathbb N}C^p$$

At the compact group level now, we obtain the decreasing intersection in the statement. Finally, the last assertion is clear from Proposition 11.31.
\end{proof}

As a main consequence of the above results, we can now formulate: 

\begin{definition}
We say that a homogeneous compact group 
$$S_N\subset G\subset U_N$$
is easy at order $p$ when $G=G^p$, with $p$ being chosen minimal with this property.
\end{definition}

Observe that the order 1 notion corresponds to the usual easiness. In general, all this is quite abstract, but there are several explicit examples, that can be worked out. For more on all this, you can check my quantum group book \cite{ba2}.

\section*{11d. Classification results}

Let us go back now to plain easiness, and discuss some classification results, following \cite{bsp}, and then the more recent paper of Tarrago-Weber \cite{twe}. In order to cut from the complexity, we must impose an extra axiom, and we will use here:

\index{uniform group}

\begin{theorem}
For an easy group $G=(G_N)$, coming from a category of partitions $D\subset P$, the following conditions are equivalent:
\begin{enumerate}
\item $G_{N-1}=G_N\cap U_{N-1}$, via the embedding $U_{N-1}\subset U_N$ given by $u\to diag(u,1)$.

\item $G_{N-1}=G_N\cap U_{N-1}$, via the $N$ possible diagonal embeddings $U_{N-1}\subset U_N$.

\item $D$ is stable under the operation which consists in removing blocks.
\end{enumerate}
If these conditions are satisfied, we say that $G=(G_N)$ is uniform.
\end{theorem}

\begin{proof}
We use the general easiness theory explained above, as follows:

\medskip

$(1)\iff(2)$ This is something standard, coming from the inclusion $S_N\subset G_N$, which makes everything $S_N$-invariant. The result follows as well from the proof of $(1)\iff(3)$ below, which can be converted into a proof of $(2)\iff(3)$, in the obvious way.

\medskip

$(1)\iff(3)$ Given a subgroup $K\subset U_{N-1}$, with fundamental representation $u$, consider the $N\times N$ matrix $v=diag(u,1)$. Our claim is that for any $\pi\in P(k)$ we have:
$$\xi_\pi\in Fix(v^{\otimes k})\iff\xi_{\pi'}\in Fix(v^{\otimes k'}),\,\forall\pi'\in P(k'),\pi'\subset\pi$$

In order to prove this, we must study the condition on the left. We have:
\begin{eqnarray*}
\xi_\pi\in Fix(v^{\otimes k})
&\iff&(v^{\otimes k}\xi_\pi)_{i_1\ldots i_k}=(\xi_\pi)_{i_1\ldots i_k},\forall i\\
&\iff&\sum_j(v^{\otimes k})_{i_1\ldots i_k,j_1\ldots j_k}(\xi_\pi)_{j_1\ldots j_k}=(\xi_\pi)_{i_1\ldots i_k},\forall i\\
&\iff&\sum_j\delta_\pi(j_1,\ldots,j_k)v_{i_1j_1}\ldots v_{i_kj_k}=\delta_\pi(i_1,\ldots,i_k),\forall i
\end{eqnarray*}

Now let us recall that our representation has the special form $v=diag(u,1)$. We conclude from this that for any index $a\in\{1,\ldots,k\}$, we must have:
$$i_a=N\implies j_a=N$$

With this observation in hand, if we denote by $i',j'$ the multi-indices obtained from $i,j$ obtained by erasing all the above $i_a=j_a=N$ values, and by $k'\leq k$ the common length of these new multi-indices, our condition becomes:
$$\sum_{j'}\delta_\pi(j_1,\ldots,j_k)(v^{\otimes k'})_{i'j'}=\delta_\pi(i_1,\ldots,i_k),\forall i$$

Here the index $j$ is by definition obtained from $j'$ by filling with $N$ values. In order to finish now, we have two cases, depending on $i$, as follows:

\medskip

\underline{Case 1}. Assume that the index set $\{a|i_a=N\}$ corresponds to a certain subpartition $\pi'\subset\pi$. In this case, the $N$ values will not matter, and our formula becomes:
$$\sum_{j'}\delta_\pi(j'_1,\ldots,j'_{k'})(v^{\otimes k'})_{i'j'}=\delta_\pi(i'_1,\ldots,i'_{k'})$$

\underline{Case 2}. Assume now the opposite, namely that the set $\{a|i_a=N\}$ does not correspond to a subpartition $\pi'\subset\pi$. In this case the indices mix, and our formula reads:
$$0=0$$

Thus, we are led to $\xi_{\pi'}\in Fix(v^{\otimes k'})$, for any subpartition $\pi'\subset\pi$, as claimed.

\medskip

Now with this claim in hand, the result follows from Tannakian duality.
\end{proof}

We can now formulate a first classification result, as follows:

\begin{theorem}
The uniform orthogonal easy groups are as follows,
$$\xymatrix@R=50pt@C=50pt{
B_N\ar[r]&O_N\\
S_N\ar[u]\ar[r]&H_N\ar[u]}$$
and this diagram is an intersection and easy generation diagram.
\end{theorem}

\begin{proof}
We know that the quantum groups in the statement are indeed easy and uniform, the corresponding categories of partitions being as follows:
$$\xymatrix@R=50pt@C50pt{
P_{12}\ar[d]&P_2\ar[d]\ar[l]\\
P&P_{even}\ar[l]}$$

Since this latter diagram is an intersection and generation diagram, we conclude that we have an intersection and easy generation diagram of quantum groups, as stated. Regarding now the classification, consider an arbitrary easy group, as follows:
$$S_N\subset G_N\subset O_N$$

This group must then come from a category of partitions, as follows:
$$P_2\subset D\subset P$$

Now if we assume $G=(G_N)$ to be uniform, this category of partitions $D$ is uniquely determined by the subset $L\subset\mathbb N$ consisting of the sizes of the blocks of the partitions in $D$. Our claim now is that the admissible sets are as follows:

\medskip

\begin{enumerate}
\item $L=\{2\}$, producing $O_N$.

\medskip

\item $L=\{1,2\}$, producing $B_N$.

\medskip

\item $L=\{2,4,6,\ldots\}$, producing $H_N$.

\medskip

\item $L=\{1,2,3,\ldots\}$, producing $S_N$.
\end{enumerate}

\medskip

Indeed, in one sense, this follows from our easiness results for $O_N,B_N,H_N,S_N$. In the other sense now, assume that $L\subset\mathbb N$ is such that the set $P_L$ consisting of partitions whose sizes of the blocks belong to $L$ is a category of partitions. We know from the axioms of the categories of partitions that the semicircle $\cap$ must be in the category, so we have $2\in L$. Our claim is that the following conditions must be satisfied as well:
$$k,l\in L,\,k>l\implies k-l\in L$$
$$k\in L,\,k\geq 2\implies 2k-2\in L$$

Indeed, we will prove that both conditions follow from the axioms of the categories of
partitions. Let us denote by $b_k\in P(0,k)$ the one-block partition, as follows:
$$b_k=\left\{\begin{matrix}\sqcap\hskip-0.7mm \sqcap&\ldots&\sqcap\\
1\hskip2mm 2&\ldots&k\end{matrix} \right\}$$

For $k>l$, we can write $b_{k-l}$ in the following way:
$$b_{k-l}=\left\{\begin{matrix}\sqcap\hskip-0.7mm
\sqcap&\ldots&\ldots&\ldots&\ldots&\sqcap\\ 1\hskip2mm 2&\ldots&l&l+1&\ldots&k\\
\sqcup\hskip-0.7mm \sqcup&\ldots&\sqcup&|&\ldots&|\\ &&&1&\ldots&k-l\end{matrix}\right\}$$

In other words, we have the following formula:
$$b_{k-l}=(b_l^*\otimes |^{\otimes k-l})b_k$$

Since all the terms of this composition are in $P_L$, we have $b_{k-l}\in P_L$, and this proves our first formula. As for the second formula, this can be proved in a similar way, by capping two adjacent $k$-blocks with a $2$-block, in the middle.

\medskip

With the above two formulae in hand, we can conclude in the following way:

\medskip

\underline{Case 1}. Assume $1\in L$. By using the first formula with $l=1$ we get:
$$k\in L\implies k-1\in L$$

This condition shows that we must have $L=\{1,2,\ldots,m\}$, for a certain number $m\in\{1,2,\ldots,\infty\}$. On the other hand, by using the second formula we get:
\begin{eqnarray*}
m\in L
&\implies&2m-2\in L\\
&\implies&2m-2\leq m\\
&\implies&m\in\{1,2,\infty\}
\end{eqnarray*}

The case $m=1$ being excluded by the condition $2\in L$, we reach to one of the two sets producing the groups $S_N,B_N$.

\medskip

\underline{Case 2}. Assume $1\notin L$. By using the first formula with $l=2$ we get:
$$k\in L\implies k-2\in L$$

This condition shows that we must have $L=\{2,4,\ldots,2p\}$, for a certain number $p\in\{1,2,\ldots,\infty\}$. On the other hand, by using the second formula we get:
\begin{eqnarray*}
2p\in L
&\implies&4p-2\in L\\
&\implies&4p-2\leq 2p\\
&\implies&p\in\{1,\infty\}
\end{eqnarray*}

Thus $L$ must be one of the two sets producing $O_N,H_N$, and we are done.
\end{proof}

All the above is very nice, but the continuation of the story is more complicated. When lifting the uniformity assumption, the final classification results become more technical, due to the presence of various copies of $\mathbb Z_2$, that can be added, while keeping the easiness property still true. To be more precise, in the real case it is known that we have exactly 6 solutions, which are as follows, with the convention $G_N'=G_N\times\mathbb Z_2$:
$$\xymatrix@R=50pt@C=50pt{
B_N\ar[r]&B_N'\ar[r]&O_N\\
S_N\ar[u]\ar[r]&S_N'\ar[u]\ar[r]&H_N\ar[u]}$$

In the unitary case now, the classification is quite similar, but more complicated, as explained in the paper of Tarrago-Weber \cite{twe}. In particular we have:

\begin{theorem}
The uniform easy groups which are purely unitary, in the sense that they appear as complexifications of real easy groups, are as follows,
$$\xymatrix@R=50pt@C=50pt{
C_N\ar[r]&U_N\\
S_N\ar[u]\ar[r]&K_N\ar[u]}$$
and this diagram is an intersection and easy generation diagram.
\end{theorem}

\begin{proof}
We know from the above that the groups in the statement are indeed easy and uniform, the corresponding categories of partitions being as follows:
$$\xymatrix@R=50pt@C50pt{
\mathcal P_{12}\ar[d]&\mathcal P_2\ar[d]\ar[l]\\
P&\mathcal P_{even}\ar[l]}$$

Since this latter diagram is an intersection and generation diagram, we conclude that we have an intersection and easy generation diagram of groups, as stated. As for the uniqueness result, the proof here is similar to the proof from the real case, from Theorem 11.35, by examining the possible sizes of the blocks of the partitions in the category, and doing some direct combinatorics. For details here, we refer to Tarrago-Weber \cite{twe}.
\end{proof}

Finally, let us mention that the easy quantum group formalism can be extended into a ``super-easy'' group formalism, covering as well the symplectic group $Sp_N$. This is something a bit technical, and we refer here to the paper of Collins-\'Sniady \cite{csn}.

\section*{11e. Exercises}

Welcome to easiness, and as exercises on this, all easy I guess, we have:

\begin{exercise}
Get to know more about set partitions, by whatever means.
\end{exercise}

\begin{exercise}
Explicitly compute some further linear maps of type $T_\pi$.
\end{exercise}

\begin{exercise}
Try improving our easy generation formalism, and results.
\end{exercise}

\begin{exercise}
Prove that the easy envelope of $U_N^d$ is the group $U_N$.
\end{exercise}

\begin{exercise}
Prove that the easy envelope of $H_N^{s,d}$ is generically $H_N^s$.
\end{exercise}

\begin{exercise}
Learn the full classification results in the orthogonal case.
\end{exercise}

\begin{exercise}
Have a look at the classification results in the unitary case.
\end{exercise}

\begin{exercise}
Have a look as well at $Sp_N$, and at $SU_2,SO_3$ too.
\end{exercise}

As bonus exercise, experiment with computing integrals, over the easy groups.

\chapter{Quantum groups}

\section*{12a. Quantum groups}

With easiness understood, what is next? You would probably say, lots of computations for the easy groups that we know, based on the following principle:

\begin{principle}
Everything advanced in relation with the easy groups that we know, such as $S_N,H_N,K_N$ or $O_N,U_N$, should come from the easiness formula for them,
$$Hom(u^{\otimes k},u^{\otimes l})=span\left(T_\pi\Big|\pi\in D(k,l)\right)$$
with $D$ being $P,P_{even},{\mathcal P}_{even}$ and $P_2,\mathcal P_2$ for our main examples. For instance the representations can be found via $End(u^{\otimes k})$, and the main character law, via $dim(Fix(u^{\otimes k}))$.
\end{principle}

This being said, such computations can be quite complex, and before getting into them, we should perhaps have one more look at generalities, I mean representation theory in general, and see if we have not missed something, which can help us in our quest.

\bigskip

In answer now, assuming some knowledge in quantum algebra, that you probably don't have, but I do, and so will you soon, after reading this chapter, here is something very interesting, comimg as a sort of big surprise, in relation with all this:

\begin{fact}
The free counterparts of the basic groups that we know, which are the quantum groups $S_N^+,H_N^+,K_N^+$ and $O_N^+,U_N^+$, are easy too, in a similar sense,
$$Hom(u^{\otimes k},u^{\otimes l})=span\left(T_\pi\Big|\pi\in D(k,l)\right)$$
with $D$ being $NC,NC_{even},{\mathcal NC}_{even}$ and $NC_2,{\mathcal NC}_2$, categories of noncrossing partitions. And, their combinatorics is much simpler than that of $S_N,H_N,K_N$ or $O_N,U_N$.
\end{fact}

Which sounds quite amazing, doesn't it, looks like at the level of general theory, we have missed something. So, believe me here, and we are therefore led to the following principle, which will replace Principle 12.1, and will be our guiding light, next:

\begin{principle}
Quantum comes first, being simpler than classical.
\end{principle}

Getting to work now, what is a quantum group? Or, perhaps at an even more foundational level, what is a quantum space? In answer, quantum spaces should certainly have something to do with quantum physics, and browsing through the early work on the subject, due to Heisenberg and Schr\"odinger, then Dirac, von Neumann and others, here is something due to von Neumann which, while physically speaking is a bit too vague to be true, mathematically is just something amazingly simple and beautiful:

\begin{principle}[Von Neumann]
A quantum space is the dual of an operator algebra.
\end{principle}

So, we will go with this. Our plan will be to talk about operator algebras, and the related notion of quantum space, then about quantum groups, appearing via suitable classes of operator algebras, and finally about $S_N^+,H_N^+,K_N^+$ and $O_N^+,U_N^+$, easiness, and Fact 12.2. As for Principle 12.1, in view of Principle 12.3, that can remain for later.

\bigskip

We must first talk about operator algebras. For our purposes here, with defining quantum groups being our main goal, it is actually convenient to deviate a bit from what von Neumann was exactly saying, and have as starting point the following related definition, due to Gelfand, chosen as usual for its amazing simplicity and beauty:

\begin{definition}
A $C^*$-algebra is a complex algebra $A$, having a norm $||.||$ making it a Banach algebra, and an involution $*$, related to the norm by the formula 
$$||aa^*||=||a||^2$$
which must hold for any $a\in A$.
\end{definition}

As a basic example here, making the link with von Neumann and physics, given a complex Hilbert space $H$, such as the Schr\"odinger space $L^2(\mathbb R^3)$ of wave functions of the electron, or its version $L^2(\mathbb R^3)\otimes\mathbb C^2$ by Pauli, accounting for the electron spin, the bounded linear operators $T:H\to H$ form a $C^*$-algebra, with norm and involution as follows:
$$||T||=\sup_{||x||=1}||Tx||\quad,\quad<Tx,y>=<x,T^*y>$$ 

This algebra is denoted $B(H)$, and more generally, any norm closed $*$-subalgebra $A\subset B(H)$ is a $C^*$-algebra. It is possible to prove that a converse of this holds, in the sense that any $C^*$-algebra appears as an operator algebra, $A\subset B(H)$. More on this in a moment, but in any case, we have our definition for the operator algebras, for sure.

\bigskip

Getting now towards quantum spaces, the original approach by von Neumann was by using the spectral theorem for normal operators $T\in B(H)$. In the context of Definition 12.5, no such operators are available, and this forces us to do some algebra instead. And as good news, we reach to the following remarkable result, due to Gelfand:

\index{Gelfand theorem}
\index{commutative algebra}
\index{compact space}

\begin{theorem}
The commutative $C^*$-algebras are precisely the algebras
$$A=C(X)$$
with $X$ being a compact space.
\end{theorem}

\begin{proof}
In one sense this is clear, because given a compact space $X$, the algebra $C(X)$ is a $C^*$-algebra, with the usual norm and involution of functions, namely:
$$||f||=\sup_{x\in X}|f(x)|\quad,\quad f^*(x)=\overline{f(x)}$$

Conversely, assuming that $A$ is commutative, consider the compact space $X$ of Banach algebra characters $\chi:A\to\mathbb C$. We have then a morphism of algebras as follows:
$$ev:A\to C(X)\quad,\quad a\to ev_a=[\chi\to\chi(a)]$$

In order to prove that $ev$ is involutive, we can argue that it is enough to prove that we have $ev_{a^*}=ev_a^*$ for the self-adjoint elements $a$. But this follows from the following computation, with the inclusion at the end coming via some complex analysis:
$$ev_a(\chi)=\chi(a)\in\sigma(a)\subset\mathbb R$$

Next, since $A$ is commutative, each element is normal, so $ev$ is isometric, with this coming from the following computation, with $\rho$ denoting the spectral radius:
$$||ev_a||=\rho(a)=||a||$$

It remains to prove that $ev$ is surjective. But this follows from the Stone-Weierstrass theorem, because $ev(A)$ is a closed subalgebra of $C(X)$, which separates the points.
\end{proof}

In view of Theorem 12.6, we can now formulate the following definition:

\begin{definition}
Given an arbitrary $C^*$-algebra $A$, we can write 
$$A=C(X)$$
and call the abstract space $X$ a compact quantum space.
\end{definition}

In other words, we can define the category of compact quantum spaces $X$ as being the category of the $C^*$-algebras $A$, with the arrows reversed. A morphism $f:X\to Y$ corresponds by definition to a morphism $\Phi:C(Y)\to C(X)$, a product of spaces $X\times Y$ corresponds by definition to a product of algebras $C(X)\otimes C(Y)$, and so on.

\bigskip

All this is of course, while certainly categorical and rigorous, a bit speculative. As a first true result, making appear the simplest non-classical quantum spaces, we have:

\begin{theorem}
The finite quantum spaces are exactly the disjoint unions of type
$$X=M_{n_1}\sqcup\ldots\sqcup M_{n_k}$$
where $M_n$ is the finite quantum space given by $C(M_n)=M_n(\mathbb C)$.
\end{theorem}

\begin{proof}
For a compact quantum space $X$, coming from a $C^*$-algebra $A$ via the formula $A=C(X)$, being finite can only mean that the following number is finite:
$$|X|=\dim_\mathbb CA<\infty$$

But now, we can use a result from chapter 9, or rather a straightforward extension of that result, stating that our finite dimensional $C^*$-algebra must decompose as follows:
$$A=M_{n_1}(\mathbb C)\oplus\ldots\oplus M_{n_k}(\mathbb C)$$

Now since direct sums of algebras $A$ correspond to disjoint unions of quantum spaces $X$, via the correspondence $A=C(X)$, this leads to the conclusion in the statement.
\end{proof}

At the general level, we have as well the following key result, announced before, due to Gelfand, Naimark and Segal, and called GNS representation theorem:

\begin{theorem}
Any $C^*$-algebra appears as an operator algebra:
$$A\subset B(H)$$
Moreover, when $A$ is separable, which is usually the case, $H$ can be taken separable.
\end{theorem}

\begin{proof}
Let us first prove that the result holds in the commutative case, $A=C(X)$. Here, we can pick a positive measure on $X$, and construct our embedding as follows:
$$C(X)\subset B(L^2(X))\quad,\quad f\to[g\to fg]$$

In general the proof is similar, the idea being that given a $C^*$-algebra $A$ we can construct a Hilbert space $H=L^2(A)$, and then an embedding as follows:
$$A\subset B(L^2(A))\quad,\quad a\to[b\to ab]$$

Finally, the last assertion is clear, because when $A$ is separable, meaning that it has a countable algebraic basis, so does the associated Hilbert space $H=L^2(A)$.
\end{proof}

And with this, end of our learning of operator algebras and quantum spaces, with Definitions 12.5, 12.7 and Theorems 12.6, 12.8, 12.9 being the starting package.

\bigskip

We are now ready to introduce the quantum groups. The axioms here, due to Woronowicz \cite{wor}, and slightly modified for our present purposes, are as follows:

\begin{definition}
A Woronowicz algebra is a $C^*$-algebra $A$, given with a unitary matrix $u\in M_N(A)$ whose coefficients generate $A$, such that the formulae
$$\Delta(u_{ij})=\sum_ku_{ik}\otimes u_{kj}\quad,\quad
\varepsilon(u_{ij})=\delta_{ij}\quad,\quad 
S(u_{ij})=u_{ji}^*$$
define morphisms of $C^*$-algebras $\Delta:A\to A\otimes A$, $\varepsilon:A\to\mathbb C$ and $S:A\to A^{opp}$, called comultiplication, counit and antipode. 
\end{definition}

Here the tensor product needed for $\Delta$ can be any $C^*$-algebra tensor product, and more on this later. In order to get rid of redundancies, coming from this and from amenability issues, we will divide everything by an equivalence relation, as follows:

\begin{definition}
We agree to identify two Woronowicz algebras, $(A,u)=(B,v)$, when we have an isomorphism of $*$-algebras
$$<u_{ij}>\simeq<v_{ij}>$$
mapping standard coordinates to standard coordinates, $u_{ij}\to v_{ij}$.
\end{definition}

We say that $A$ is cocommutative when $\Sigma\Delta=\Delta$, where $\Sigma(a\otimes b)=b\otimes a$ is the flip. We have then the following key result, from \cite{wor}, providing us with examples:

\begin{theorem}
The following are Woronowicz algebras:
\begin{enumerate}
\item $C(G)$, with $G\subset U_N$ compact Lie group, with structural maps as follows:
$$\Delta(\varphi)=\big[(g,h)\to \varphi(gh)\big]\quad,\quad 
\varepsilon(\varphi)=\varphi(1)\quad,\quad
S(\varphi)=\big[g\to\varphi(g^{-1})\big]$$

\item $C^*(\Gamma)$, with $F_N\to\Gamma$ finitely generated group, with structural maps as follows:
$$\Delta(g)=g\otimes g\quad,\quad 
\varepsilon(g)=1\quad,\quad
S(g)=g^{-1}$$
\end{enumerate}
Moreover, we obtain in this way all the commutative/cocommutative algebras.
\end{theorem}

\begin{proof}
For the first assertion we can use the matrix $u=(u_{ij})$ formed by the standard matrix coordinates of $G$, which is by definition given by:
$$g=\begin{pmatrix}
u_{11}(g)&\ldots&u_{1N}(g)\\
\vdots&&\vdots\\
u_{N1}(g)&\ldots&u_{NN}(g)
\end{pmatrix}$$

For the second assertion, where $C^*(\Gamma)$ stands for the formal $C^*$-algebra completion of the usual group algebra $\mathbb C[\Gamma]$, we can use the diagonal matrix formed by generators:
$$u=\begin{pmatrix}
g_1&&0\\
&\ddots&\\
0&&g_N
\end{pmatrix}$$

Finally, regarding the last assertion, in the commutative case this follows from the Gelfand theorem, and in the cocommutative case, we will be back to this.
\end{proof}

In order to get now to quantum groups, we will need as well:

\begin{proposition}
Assuming that $G\subset U_N$ is abelian, we have an identification of Woronowicz algebras $C(G)=C^*(\Gamma)$, with $\Gamma$ being the Pontrjagin dual of $G$:
$$\Gamma=\big\{\chi:G\to\mathbb T\big\}$$
Conversely, assuming that $F_N\to\Gamma$ is abelian, we have an identification of Woronowicz algebras $C^*(\Gamma)=C(G)$, with $G$ being the Pontrjagin dual of $\Gamma$:
$$G=\big\{\chi:\Gamma\to\mathbb T\big\}$$
Thus, the Woronowicz algebras which are both commutative and cocommutative are exactly those of type $A=C(G)=C^*(\Gamma)$, with $G,\Gamma$ being abelian, in Pontrjagin duality.
\end{proposition}

\begin{proof}
This follows from the Gelfand theorem applied to $C^*(\Gamma)$, and from the fact that the characters of a group algebra come from the characters of the group.
\end{proof}

In view of this result, and of the findings from Theorem 12.12 too, we have the following definition, complementing Definition 12.10 and Definition 12.11:

\begin{definition}
Given a Woronowicz algebra, we write it as follows, and call $G$ a compact quantum Lie group, and $\Gamma$ a finitely generated discrete quantum group:
$$A=C(G)=C^*(\Gamma)$$
Also, we say that $G,\Gamma$ are dual to each other, and write $G=\widehat{\Gamma},\Gamma=\widehat{G}$.
\end{definition}

Let us discuss now some tools for studying the Woronowicz algebras, and the underlying quantum groups. First, we have the following result:

\begin{proposition}
Let $(A,u)$ be a Woronowicz algebra.
\begin{enumerate} 
\item $\Delta,\varepsilon$ satisfy the usual axioms for a comultiplication and a counit, namely:
$$(\Delta\otimes id)\Delta=(id\otimes \Delta)\Delta$$
$$(\varepsilon\otimes id)\Delta=(id\otimes\varepsilon)\Delta=id$$

\item $S$ satisfies the antipode axiom, on the $*$-algebra generated by entries of $u$: 
$$m(S\otimes id)\Delta=m(id\otimes S)\Delta=\varepsilon(.)1$$

\item In addition, the square of the antipode is the identity, $S^2=id$.
\end{enumerate}
\end{proposition}

\begin{proof}
As a first observation, the result holds in the commutative case, $A=C(G)$ with $G\subset U_N$. Indeed, here we know from Theorem 12.12 that $\Delta,\varepsilon,S$ appear as functional analytic transposes of the multiplication, unit and inverse maps $m,u,i$:
$$\Delta=m^t\quad,\quad 
\varepsilon=u^t\quad,\quad 
S=i^t$$

Thus, the various conditions in the statement on $\Delta,\varepsilon,S$ simply come from the group axioms satisfied by $m,u,i$. In general now, the first axiom follows from:
$$(\Delta\otimes id)\Delta(u_{ij})=(id\otimes \Delta)\Delta(u_{ij})=\sum_{kl}u_{ik}\otimes u_{kl}\otimes u_{lj}$$

As for the other axioms, the verifications here are similar.
\end{proof}

In order to reach to more advanced results, the idea will be that of doing representation theory. Following Woronowicz \cite{wor}, let us start with the following definition:

\begin{definition}
Given $(A,u)$, we call corepresentation of it any unitary matrix $v\in M_n(\mathcal A)$, with $\mathcal A=<u_{ij}>$, satisfying the same conditions as $u$, namely:
$$\Delta(v_{ij})=\sum_kv_{ik}\otimes v_{kj}\quad,\quad
\varepsilon(v_{ij})=\delta_{ij}\quad,\quad
S(v_{ij})=v_{ji}^*$$
We also say that $v$ is a representation of the underlying compact quantum group $G$.
\end{definition}

In the commutative case, $A=C(G)$ with $G\subset U_N$, we obtain in this way the finite dimensional unitary smooth representations $v:G\to U_n$, via the following formula:
$$v(g)=\begin{pmatrix}
v_{11}(g)&\ldots&v_{1n}(g)\\
\vdots&&\vdots\\
v_{n1}(g)&\ldots&v_{nn}(g)
\end{pmatrix}$$

With these conventions, we have the following fundamental result, from \cite{wor}:

\begin{theorem}
Any Woronowicz algebra has a unique Haar integration functional, 
$$\left(\int_G\otimes\, id\right)\Delta=\left(id\otimes\int_G\right)\Delta=\int_G(.)1$$
which can be constructed by starting with any faithful positive form $\varphi\in A^*$, and setting
$$\int_G=\lim_{n\to\infty}\frac{1}{n}\sum_{k=1}^n\varphi^{*k}$$
where $\phi*\psi=(\phi\otimes\psi)\Delta$. Moreover, for any corepresentation $v\in M_n(\mathbb C)\otimes A$ we have
$$\left(id\otimes\int_G\right)v=P$$
where $P$ is the orthogonal projection onto $Fix(v)=\{\xi\in\mathbb C^n|v\xi=\xi\}$.
\end{theorem}

\begin{proof}
This is something that we know for $A=C(G)$ with $G\subset U_N$, and the proof in general is similar. In case you doubt, I should tell you that the Haar proof from chapter 9 was more or less copy/pasted from my quantum book group \cite{ba2}. So, this works.
\end{proof}

We can now develop, again following \cite{wor}, the Peter-Weyl theory for the corepresentations of $A$. Consider the dense subalgebra $\mathcal A\subset A$ generated by the coefficients of the fundamental corepresentation $u$, and endow it with the following scalar product: 
$$<a,b>=\int_Gab^*$$

With this convention, we have the following result, also from \cite{wor}:

\begin{theorem}
We have the following Peter-Weyl type results:
\begin{enumerate}
\item Any corepresentation decomposes as a sum of irreducible corepresentations.

\item Each irreducible corepresentation appears inside a certain $u^{\otimes k}$.

\item $\mathcal A=\bigoplus_{v\in Irr(A)}M_{\dim(v)}(\mathbb C)$, the summands being pairwise orthogonal.

\item The characters of irreducible corepresentations form an orthonormal system.
\end{enumerate}
\end{theorem}

\begin{proof}
This is something that we know well from chapter 10, when $G\subset U_N$ is a compact group. In general, the proof is quite similar, by using Theorem 12.17.
\end{proof}

Finally, no discussion about compact and discrete quantum groups would be complete without a word on amenability. The result here, again from \cite{wor}, is as follows:

\begin{theorem}
Let $A_{full}$ be the enveloping $C^*$-algebra of $\mathcal A$, and $A_{red}$ be the quotient of $A$ by the null ideal of the Haar integration. The following are then equivalent:
\begin{enumerate}
\item The Haar functional of $A_{full}$ is faithful.

\item The projection map $A_{full}\to A_{red}$ is an isomorphism.

\item The counit map $\varepsilon:A_{full}\to\mathbb C$ factorizes through $A_{red}$.

\item We have $N\in\sigma(Re(\chi_u))$, the spectrum being taken inside $A_{red}$.
\end{enumerate}
If this is the case, we say that the underlying discrete quantum group $\Gamma$ is amenable.
\end{theorem}

\begin{proof}
This is well-known in the group dual case, $A=C^*(\Gamma)$, with $\Gamma$ being a usual discrete group. In general, the result follows by adapting the group dual case proof:

\medskip

$(1)\iff(2)$ This simply follows from the fact that the GNS construction for the algebra $A_{full}$ with respect to the Haar functional produces the algebra $A_{red}$.

\medskip

$(2)\iff(3)$ Here $\implies$ is trivial, and conversely, a counit $\varepsilon:A_{red}\to\mathbb C$ produces an isomorphism $\Phi:A_{red}\to A_{full}$, by slicing the map $\widetilde{\Delta}:A_{red}\to A_{red}\otimes A_{full}$.

\medskip

$(3)\iff(4)$ Here $\implies$ is clear, coming from $\varepsilon(N-Re(\chi (u)))=0$, and the converse can be proved by doing some functional analysis. See \cite{wor}.
\end{proof}

This was for the basic theory of the quantum groups in the sense of Woronowicz, quickly explained. For more on all this, you have for instance my book \cite{ba2}.

\section*{12b. Quantum permutations}

The compact quantum groups include the compact Lie groups, $G\subset U_N$, and the abstract duals $G=\widehat{\Gamma}$ of the finitely generated groups $F_N\to\Gamma$. Following Wang \cite{wan}, let us discuss now a number of truly ``new'' quantum groups. We first have:

\index{free orthogonal group}
\index{free unitary group}

\begin{theorem}
The following universal algebras are Woronowicz algebras,
$$C(O_N^+)=C^*\left((u_{ij})_{i,j=1,\ldots,N}\Big|u=\bar{u},u^t=u^{-1}\right)$$
$$C(U_N^+)=C^*\left((u_{ij})_{i,j=1,\ldots,N}\Big|u^*=u^{-1},u^t=\bar{u}^{-1}\right)$$
so the underlying quantum spaces $O_N^+,U_N^+$ are compact quantum groups.
\end{theorem}

\begin{proof}
This comes from the elementary fact that if a matrix $u=(u_{ij})$ is orthogonal or biunitary, then so must be the following matrices:
$$(u^\Delta)_{ij}=\sum_ku_{ik}\otimes u_{kj}\quad,\quad 
(u^\varepsilon)_{ij}=\delta_{ij}\quad,\quad 
(u^S)_{ij}=u_{ji}^*$$

Thus we can define $\Delta,\varepsilon,S$ by using the universal property of $C(O_N^+)$, $C(U_N^+)$.
\end{proof}

The basic properties of these beasts $O_N^+,U_N^+$ can be summarized as follows:

\begin{theorem}
The quantum groups $O_N^+,U_N^+$ have the following properties:
\begin{enumerate}
\item The closed subgroups $G\subset U_N^+$ are exactly the $N\times N$ compact quantum groups. As for the closed subgroups $G\subset O_N^+$, these are those satisfying $u=\bar{u}$.

\item We have liberation embeddings $O_N\subset O_N^+$ and $U_N\subset U_N^+$, obtained by dividing the algebras $C(O_N^+),C(U_N^+)$ by their respective commutator ideals.

\item We have as well embeddings $\widehat{L}_N\subset O_N^+$ and $\widehat{F}_N\subset U_N^+$, where $L_N$ is the free product of $N$ copies of $\mathbb Z_2$, and where $F_N$ is the free group on $N$ generators.
\end{enumerate}
\end{theorem}

\begin{proof}
All these assertions are elementary, as follows:

\medskip

(1) This is clear from definitions, with the remark that, in the context of Definition 12.10, the formula $S(u_{ij})=u_{ji}^*$ shows that the matrix $\bar{u}$ must be unitary too. 

\medskip

(2) This follows from the Gelfand theorem. To be more precise, this shows that we have presentation results for $C(O_N),C(U_N)$, similar to those in Theorem 12.20, but with the commutativity between the standard coordinates and their adjoints added:
\begin{eqnarray*}
C(O_N)&=&C^*_{comm}\left((u_{ij})_{i,j=1,\ldots,N}\Big|u=\bar{u},u^t=u^{-1}\right)\\
C(U_N)&=&C^*_{comm}\left((u_{ij})_{i,j=1,\ldots,N}\Big|u^*=u^{-1},u^t=\bar{u}^{-1}\right)
\end{eqnarray*}

Thus, we are led to the conclusion in the statement.

\medskip

(3) This follows indeed from (1) and from Theorem 12.12, with the remark that with $u=diag(g_1,\ldots,g_N)$, the condition $u=\bar{u}$ is equivalent to $g_i^2=1$, for any $i$.
\end{proof}

This was for the basic theory of $O_N^+,U_N^+$, and I can hear you thinking that all this is not that easy to digest, and no match for the simplicity of $O_N,U_N$. In answer, wait for it. We will do later representation theory for $O_N^+,U_N^+$, and you will be amazed.

\bigskip

Still following Wang \cite{wan}, let us discuss now the construction and basic properties of the quantum permutation group $S_N^+$. Let us first look at $S_N$. We have here:

\index{magic matrix}
\index{magic unitary}
\index{symmetric group}

\begin{proposition}
Consider the symmetric group $S_N$, viewed as permutation group of the $N$ coordinate axes of $\mathbb R^N$. The coordinate functions on $S_N\subset O_N$ are given by
$$u_{ij}=\chi\left(\sigma\in G\Big|\sigma(j)=i\right)$$
and the matrix $u=(u_{ij})$ that these functions form is magic, in the sense that its entries are projections $(p^2=p^*=p)$, summing up to $1$ on each row and each column.
\end{proposition}

\begin{proof}
The action of $S_N$ on the standard basis $e_1,\ldots,e_N\in\mathbb R^N$ being given by $\sigma:e_j\to e_{\sigma(j)}$, this gives the formula of $u_{ij}$ in the statement. As for the fact that the matrix $u=(u_{ij})$ that these functions form is magic, this is clear.
\end{proof}

With a bit more effort, we obtain the following nice characterization of $S_N$:

\index{Gelfand theorem}
\index{commutative algebra}

\begin{theorem}
The algebra of functions on $S_N$ has the following presentation,
$$C(S_N)=C^*_{comm}\left((u_{ij})_{i,j=1,\ldots,N}\Big|u={\rm magic}\right)$$
and the multiplication, unit and inversion map of $S_N$ appear from the maps
$$\Delta(u_{ij})=\sum_ku_{ik}\otimes u_{kj}\quad,\quad 
\varepsilon(u_{ij})=\delta_{ij}\quad,\quad 
S(u_{ij})=u_{ji}$$
defined at the algebraic level, of functions on $S_N$, by transposing.
\end{theorem}

\begin{proof}
The universal algebra $A$ in the statement being commutative, by the Gelfand theorem it must be of the form $A=C(X)$, with $X$ being a certain compact space. Now since we have coordinates $u_{ij}:X\to\mathbb R$, we have an embedding $X\subset M_N(\mathbb R)$. Also, since we know that these coordinates form a magic matrix, the elements $g\in X$ must be 0-1 matrices, having exactly one 1 entry on each row and each column, and so $X=S_N$. Thus we have proved the first assertion, and the second assertion is clear as well.
\end{proof}

We can now liberate $S_N$, into a certain quantum group $S_N^+$, as follows:

\begin{theorem}
The following universal $C^*$-algebra, with magic meaning as usual formed by projections $(p^2=p^*=p)$, summing up to $1$ on each row and each column,
$$C(S_N^+)=C^*\left((u_{ij})_{i,j=1,\ldots,N}\Big|u={\rm magic}\right)$$
is a Woronowicz algebra, with comultiplication, counit and antipode given by:
$$\Delta(u_{ij})=\sum_ku_{ik}\otimes u_{kj}\quad,\quad 
\varepsilon(u_{ij})=\delta_{ij}\quad,\quad 
S(u_{ij})=u_{ji}$$
Thus the space $S_N^+$ is a compact quantum group, called quantum permutation group.
\end{theorem}

\begin{proof}
As a first observation, the universal $C^*$-algebra in the statement is indeed well-defined, because the conditions $p^2=p^*=p$ satisfied by the coordinates give:
$$||u_{ij}||\leq1$$

In order to prove now that we have a Woronowicz algebra, we must construct maps $\Delta,\varepsilon,S$ given by the formulae in the statement. Consider the following matrices:
$$u^\Delta_{ij}=\sum_ku_{ik}\otimes u_{kj}\quad,\quad 
u^\varepsilon_{ij}=\delta_{ij}\quad,\quad 
u^S_{ij}=u_{ji}$$

Our claim is that, since $u$ is magic, so are these three matrices. Indeed, regarding $u^\Delta$, its entries are idempotents, as shown by the following computation:
$$(u_{ij}^\Delta)^2
=\sum_{kl}u_{ik}u_{il}\otimes u_{kj}u_{lj}
=\sum_{kl}\delta_{kl}u_{ik}\otimes\delta_{kl}u_{kj}
=u_{ij}^\Delta$$

These elements are self-adjoint as well, as shown by the following computation:
$$(u_{ij}^\Delta)^*
=\sum_ku_{ik}^*\otimes u_{kj}^*
=\sum_ku_{ik}\otimes u_{kj}
=u_{ij}^\Delta$$

The row and column sums for the matrix $u^\Delta$ can be computed as follows:
$$\sum_ju_{ij}^\Delta
=\sum_{jk}u_{ik}\otimes u_{kj}
=\sum_ku_{ik}\otimes 1
=1$$
$$\sum_iu_{ij}^\Delta
=\sum_{ik}u_{ik}\otimes u_{kj}
=\sum_k1\otimes u_{kj}
=1$$

Thus, $u^\Delta$ is magic. Regarding now $u^\varepsilon,u^S$, these matrices are magic too, and this for obvious reasons. Thus, all our three matrices $u^\Delta,u^\varepsilon,u^S$ are magic, so we can define $\Delta,\varepsilon,S$ by the formulae in the statement, by using the universality property of $C(S_N^+)$.
\end{proof}

Our first task now is to make sure that Theorem 12.24 produces indeed a new quantum group, which does not collapse to $S_N$. Following Wang, we have:

\begin{theorem}
We have an embedding $S_N\subset S_N^+$, given at the algebra level by: 
$$u_{ij}\to\chi\left(\sigma\in S_N\Big|\sigma(j)=i\right)$$
This is an isomorphism at $N\leq3$, but not at $N\geq4$, where $S_N^+$ is not classical, nor finite.
\end{theorem} 

\begin{proof}
The fact that we have indeed an embedding as above follows from Theorem 12.23. Observe that in fact more is true, because Theorems 12.23 and 12.24 give:
$$C(S_N)=C(S_N^+)\Big/\Big<ab=ba\Big>$$

Regarding now the second assertion, we can prove this in four steps, as follows:

\medskip

\underline{Case $N=2$}. The fact that $S_2^+$ is indeed classical, and hence collapses to $S_2$, is trivial, because the $2\times2$ magic matrices are as follows, with $p$ being a projection:
$$U=\begin{pmatrix}p&1-p\\1-p&p\end{pmatrix}$$

Thus $C(S_2^+)$ is commutative, and equals its biggest commutative quotient, $C(S_2)$.

\medskip

\underline{Case $N=3$}. It is enough to check that $u_{11},u_{22}$ commute. But this follows from:
\begin{eqnarray*}
u_{11}u_{22}
&=&u_{11}u_{22}(u_{11}+u_{12}+u_{13})\\
&=&u_{11}u_{22}u_{11}+u_{11}u_{22}u_{13}\\
&=&u_{11}u_{22}u_{11}+u_{11}(1-u_{21}-u_{23})u_{13}\\
&=&u_{11}u_{22}u_{11}
\end{eqnarray*}

Indeed, by conjugating, $u_{22}u_{11}=u_{11}u_{22}u_{11}$, so $u_{11}u_{22}=u_{22}u_{11}$, as desired.

\medskip

\underline{Case $N=4$}. Consider the following matrix, with $p,q$ being projections:
$$U=\begin{pmatrix}
p&1-p&0&0\\
1-p&p&0&0\\
0&0&q&1-q\\
0&0&1-q&q
\end{pmatrix}$$ 

This matrix is magic, and we can choose $p,q\in B(H)$ as for the algebra $<p,q>$ to be noncommutative and infinite dimensional. We conclude that $C(S_4^+)$ is noncommutative and infinite dimensional as well, and so $S_4^+$ is non-classical and infinite, as claimed.

\medskip

\underline{Case $N\geq5$}. Here we can use the standard embedding $S_4^+\subset S_N^+$, obtained at the level of the corresponding magic matrices in the following way:
$$u\to\begin{pmatrix}u&0\\ 0&1_{N-4}\end{pmatrix}$$

Indeed, with this in hand, the fact that $S_4^+$ is a non-classical, infinite compact quantum group implies that $S_N^+$ with $N\geq5$ has these two properties as well.
\end{proof}

The above result might seem quite surprising. How on Earth can the set $\{1,2,3,4\}$ have an infinity of quantum permutations, and will us be able to fully understand this, one day. But do not worry, the remainder of this chapter will be here for that. As a first observation, as a matter of doublechecking our findings, we are not wrong with our formalism, because as discovered once again by Wang \cite{wan}, we have as well:

\begin{theorem}
The quantum permutation group $S_N^+$ acts on the set $X=\{1,\ldots,N\}$, the corresponding coaction map $\Phi:C(X)\to C(X)\otimes C(S_N^+)$ being given by:
$$\Phi(e_i)=\sum_je_j\otimes u_{ji}$$
In fact, $S_N^+$ is the biggest compact quantum group acting on $X$, by leaving the counting measure invariant, in the sense that $(tr\otimes id)\Phi=tr(.)1$, where $tr(e_i)=\frac{1}{N},\forall i$.
\end{theorem}

\begin{proof}
Our claim is that given a compact matrix quantum group $G$, the following formula defines a morphism of algebras, which is a coaction map, leaving the trace invariant, precisely when the matrix $u=(u_{ij})$ is a magic corepresentation of $C(G)$: 
$$\Phi(e_i)=\sum_je_j\otimes u_{ji}$$

Indeed, let us first determine when $\Phi$ is multiplicative. We have:
$$\Phi(e_i)\Phi(e_k)
=\sum_{jl}e_je_l\otimes u_{ji}u_{lk}
=\sum_je_j\otimes u_{ji}u_{jk}$$
$$\Phi(e_ie_k)
=\delta_{ik}\Phi(e_i)
=\delta_{ik}\sum_je_j\otimes u_{ji}$$

We conclude that the multiplicativity of $\Phi$ is equivalent to the following conditions:
$$u_{ji}u_{jk}=\delta_{ik}u_{ji}\quad,\quad\forall i,j,k$$

Similarly, $\Phi$ is unital when $\sum_iu_{ji}=1$, $\forall j$. Finally, the fact that $\Phi$ is a $*$-morphism translates into $u_{ij}=u_{ij}^*$, $\forall i,j$. Summing up, in order for $\Phi(e_i)=\sum_je_j\otimes u_{ji}$ to be a morphism of $C^*$-algebras, the elements $u_{ij}$ must be projections, summing up to 1 on each row of $u$. Regarding now the preservation of the trace, observe that we have:
$$(tr\otimes id)\Phi(e_i)=\frac{1}{N}\sum_ju_{ji}$$

Thus the trace is preserved precisely when the elements $u_{ij}$ sum up to 1 on each of the columns of $u$. We conclude from this that $\Phi(e_i)=\sum_je_j\otimes u_{ji}$ is a morphism of $C^*$-algebras preserving the trace precisely when $u$ is magic, and this gives the result.
\end{proof}

In order to further study $S_N^+$, we will need the notion of dual free product of compact quantum groups $\hat{*}$, which is defined according to the following formula:
$$C(G\,\hat{*}\,H)=C(G)*C(H)\quad,\quad w=diag(u,v)$$

With this convention, we have the following result, improving Theorem 12.25:

\begin{theorem}
The quantum groups $S_N^+$ have the following properties:
\begin{enumerate}
\item We have $S_N^+\,\hat{*}\,S_M^+\subset S_{N+M}^+$, for any $N,M$.

\item In particular, we have an embedding $\widehat{D_\infty}\subset S_4^+$. 

\item If $\mathbb Z_{N_1}*\ldots*\mathbb Z_{N_k}\to\Gamma$, with $N=\sum N_i$, then $\widehat{\Gamma}\subset S_N^+$.

\item The quantum groups $S_N^+$ with $N\geq5$ are not coamenable.
\end{enumerate}
\end{theorem}

\begin{proof}
This is something quite standard, the idea being as follows:

\medskip

(1) If we denote by $u,v$ the fundamental corepresentations of $C(S_N^+),C(S_M^+)$, the fundamental corepresentation of $C(S_N^+\,\hat{*}\,S_M^+)$ is by definition:
$$w=\begin{pmatrix}u&0\\0&v\end{pmatrix}$$

But this matrix is magic, because both $u,v$ are magic, and this gives the result.

\medskip

(2) This result, which refines our $N=4$ trick from the proof of Theorem 12.25, follows from (1) with $N=M=2$. Indeed, we have the following computation:
\begin{eqnarray*}
S_2^+\,\hat{*}\,S_2^+
&=&S_2\,\hat{*}\, S_2
=\mathbb Z_2\,\hat{*}\, \mathbb Z_2\\
&\simeq&\widehat{\mathbb Z_2}\,\hat{*}\, \widehat{\mathbb Z_2}
=\widehat{\mathbb Z_2*\mathbb Z_2}\\
&=&\widehat{D_\infty}
\end{eqnarray*}

(3) This result, which generalizes (2), can be deduced as follows:
\begin{eqnarray*}
\widehat{\Gamma}
&\subset&\widehat{\mathbb Z_{N_1}*\ldots*\mathbb Z_{N_k}}
=\widehat{\mathbb Z_{N_1}}\,\hat{*}\,\ldots\,\hat{*}\,\widehat{\mathbb Z_{N_k}}\\
&\simeq&\mathbb Z_{N_1}\,\hat{*}\,\ldots\,\hat{*}\,\mathbb Z_{N_k}
\subset S_{N_1}\,\hat{*}\,\ldots\,\hat{*}\,S_{N_k}\\
&\subset&S_{N_1}^+\,\hat{*}\,\ldots\,\hat{*}\,S_{N_k}^+
\subset S_N^+
\end{eqnarray*}

(4) This follows from (3), because at $N=5$ the dual of the group $\Gamma=\mathbb Z_2*\mathbb Z_3$, which is well-known not to be amenable, embeds into $S_5^+$. As for the general case, that of $S_N^+$ with $N\geq5$, here the result follows by using the embedding $S_5^+\subset S_N^+$.
\end{proof}

As a conclusion to this, the quantum group $S_N^+$ is classical at $N=1,2,3$, non-classical but probably quite reasonable at $N=4$, and wild at $N\geq5$. Good to know.

\section*{12c. Diagrams, easiness}

In order to further study the quantum groups $O_N^+,U_N^+,S_N^+$ we can use representation theory. We have indeed the following version of Tannakian duality:

\begin{theorem}
The following operations are inverse to each other:
\begin{enumerate}
\item The construction $A\to C$, which associates to any Woronowicz algebra $A$ the tensor category formed by the intertwiner spaces $C_{kl}=Hom(u^{\otimes k},u^{\otimes l})$.

\item The construction $C\to A$, which associates to a tensor category $C$ the Woronowicz algebra $A$ presented by the relations $T\in Hom(u^{\otimes k},u^{\otimes l})$, with $T\in C_{kl}$.
\end{enumerate}
\end{theorem}

\begin{proof}
As before with such general results, this is something that we already know for the classical groups, and the extension to quantum groups is straightforward. And for further making my case, the Tannakian discussion from chapter 10 was in fact more or less copy/pasted from by quantum group book \cite{ba2}. So, this definitely works.
\end{proof}

Next comes the notion of easiness. Let us start with the following definition:

\begin{definition}
Let $P(k,l)$ be the set of partitions between an upper row of $k$ points, and a lower row of $l$ points. A set $D=\bigsqcup_{k,l}D(k,l)$ with $D(k,l)\subset P(k,l)$ is called a category of partitions when it has the following properties:
\begin{enumerate}
\item Stability under the horizontal concatenation, $(\pi,\sigma)\to[\pi\sigma]$.

\item Stability under the vertical concatenation, $(\pi,\sigma)\to[^\sigma_\pi]$.

\item Stability under the upside-down turning, $\pi\to\pi^*$.

\item Each set $P(k,k)$ contains the identity partition $||\ldots||$.

\item The sets $P(\emptyset,\circ\bullet)$, $P(\emptyset,\bullet\circ)$ contain the semicircle partition $\cap$.
\end{enumerate}
\end{definition} 

Observe that this is precisely the definition that we used in chapter 11, with the condition there on the basic crossing $\slash\hskip-2.1mm\backslash$, which produces commutativity via Tannakian duality, removed. Next, in relation with quantum groups, we have the following notion:

\begin{definition}
A compact quantum matrix group $G$ is called easy when
$$Hom(u^{\otimes k},u^{\otimes l})=span\left(T_\pi\Big|\pi\in D(k,l)\right)$$
for any colored integers $k,l$, for a certain category of partitions $D\subset P$, where
$$T_\pi(e_{i_1}\otimes\ldots\otimes e_{i_k})=\sum_{j_1\ldots j_l}\delta_\pi\begin{pmatrix}i_1&\ldots&i_k\\ j_1&\ldots&j_l\end{pmatrix}e_{j_1}\otimes\ldots\otimes e_{j_l}$$
with the Kronecker type symbols $\delta_\pi\in\{0,1\}$ depending on whether the indices fit or not. 
\end{definition}

Again, this is something coming as a continuation of the material from chapter 11, and for more on this definition, and its meaning, we refer to the material there. Getting now to the main examples, adding to those from chapter 11, these are as follows:

\begin{theorem}
The following quantum groups are easy:
\begin{enumerate}
\item $U_N^+$, coming from the category of noncrossing matching pairings ${\mathcal NC}_2$.

\item $O_N^+$, coming from the category of noncrossing pairings $NC_2$.

\item $S_N^+$, coming from the category of all noncrossing partitions $NC$.
\end{enumerate}
\end{theorem}

\begin{proof}
This comes by using the same arguments as for $U_N,O_N,S_N$, as follows:

\medskip

(1) The quantum group $U_N^+$ is defined via the following relations:
$$u^*=u^{-1}\quad,\quad 
u^t=\bar{u}^{-1}$$ 

Thus, the associated Tannakian category is $C=span(T_\pi|\pi\in D)$, with:
$$D
=<{\ }^{\,\cap}_{\circ\bullet}\,\,,{\ }^{\,\cap}_{\bullet\circ}>
={\mathcal NC}_2$$

(2) The subgroup $O_N^+\subset U_N^+$ is defined by imposing the following relations:
$$u_{ij}=\bar{u}_{ij}$$

Thus, the associated Tannakian category is $C=span(T_\pi|\pi\in D)$, with:
$$D
=<\mathcal{NC}_2,|^{\hskip-1.32mm\circ}_{\hskip-1.32mm\bullet},|_{\hskip-1.32mm\circ}^{\hskip-1.32mm\bullet}>
=NC_2$$

(3) The subgroup $S_N^+\subset O_N^+$ is defined by imposing the following relations:
$$u={\rm magic}$$

But, as explained in chapter 11, the magic condition reads $T_\mu\in Hom(u^{\otimes 2},u)$, with $\mu\in P(2,1)$ being the fork partition. We conclude that $S_N^+$ is indeed easy, coming from:
$$D=<NC_2,\mu>=NC$$

Thus, we are led to the conclusions in the statement.
\end{proof}

As a first consequence of the above result, of theoretical nature, we have:

\begin{theorem}
The various liberation operations that we have, namely
$$U_N\to U_N^+\quad,\quad O_N\to O_N^+\quad,\quad S_N\to S_N^+$$
all appear, at the level of the corresponding categories of partitions, via:
$$D\to D\cap NC$$
As for the reverse operations $G^+\to G$, these appear via $D\to<D,\slash\hskip-2.1mm\backslash>$.
\end{theorem}

\begin{proof}
This is indeed something quite self-explanatory, coming from Theorem 12.31, and for some further comments on the meaning of all this, we refer to \cite{bsp}.
\end{proof}

Time perhaps for some applications, in relation with Fact 12.2 and Principle 12.3? Hang on, many spectacular things to come. As a first piece of good news, we have:

\begin{theorem}
Unlike the passage $O_N\to U_N$, which is something non-trivial, coming from $\mathfrak o_N\to\mathfrak u_N$ at Lie algebra level, the passage $O_N^+\to U_N^+$ is a free complexification
$$U_N^+=\widetilde{O_N^+}$$
with this meaning that we have an embedding of algebras as follows,
$$C(U_N^+)\subset C(\mathbb T)*C(O_N^+)\quad,\quad u=zv$$
up to identifying the full and reduced versions of these algebras.
\end{theorem}

\begin{proof}
Many things going on here, the idea being as follows:

\medskip

(1) The orthogonal and unitary groups $O_N,N_N$ are both Lie groups, in the sense that they are smooth manifolds, and the corresponding Lie algebras $\mathfrak o_N,\mathfrak u_N$, which are by definition the respective tangent spaces at 1, can be computed by differentiating the equations defining $O_N,U_N$, with the conclusion being as follows:
$$\mathfrak o_N=\left\{ A\in M_N(\mathbb R)\Big|A^t=-A\right\}\quad,\quad 
\mathfrak u_N=\left\{ B\in M_N(\mathbb C)\Big|B^*=-B\right\}$$

By doing the linear algebra, we conclude that we have a direct sum decomposition of real linear spaces as follows, with $\Delta\subset M_N(\mathbb R)$ being the diagonal matrices:
$$\mathfrak u_N\simeq\mathfrak o_N\oplus\Delta\oplus\mathfrak o_N$$

Now based on this, we can say that $\mathfrak o_N\to\mathfrak u_N$ is some sort of ``complexification'', and so that the passage $O_N\to U_N$ appears too as a complexification, as follows:
$$O_N\to\mathfrak o_N\to\mathfrak u_n\to U_N$$

(2) Getting now to quantum groups, given an arbitrary pair $(G,v)$ we can associate to it its free complexification $(\widetilde{G},\tilde{v})$ via the formula in the statement, namely:
$$C(\widetilde{G})=<\tilde{v}_{ij}>\subset C(\mathbb T)*C(G)\quad,\quad \tilde{v}=zv$$

And our claim is that we have the equality in the statement, $\widetilde{O_N^+}=U_N^+$.

\medskip

(3) In order to prove this, we can use the easiness property of $O_N^+,U_N^+$. Indeed, let us denote by $v,zv,u$ the fundamental representations of the following quantum groups:
$$O_N^+\subset\widetilde{O_N^+}\subset U_N^+$$

At the level of the associated Hom spaces we obtain reverse inclusions, as follows:
$$Hom(v^{\otimes k},v^{\otimes l})
\supset Hom((zv)^{\otimes k},(zv)^{\otimes l})
\supset Hom(u^{\otimes k},u^{\otimes l})$$

But the spaces on the left and on the right are known from Theorem 12.31, with the result there stating that these spaces are as follows:
$$span\left(T_\pi\Big|\pi\in NC_2(k,l)\right)\supset span\left(T_\pi\Big|\pi\in\mathcal{NC}_2(k,l)\right)$$

Regarding now the spaces in the middle, these are obtained from those on the left by some sort of ``coloring'', so we obtain the same spaces as those on the right. Thus, our embedding $\widetilde{O_N^+}\subset U_N^+$ must be an isomorphism, modulo the usual equivalence relation.
\end{proof}

As an interesting consequence of the above result, we have:

\begin{theorem}
We have an identification of quantum groups, as follows:
$$PO_N^+=PU_N^+$$
Thus, free projective geometry is scalarless, at the same time real and complex.
\end{theorem}

\begin{proof}
Again, we have several things going on here, as follows:

\medskip

(1) To start with, in the general context of the Woronowicz quantum groups, given a pair $(G,v)$, we can talk about its projective version $(PG,u)$, constructed as follows:
$$C(PG)=<u_{ij}>\subset C(G)\quad,\quad u=v\otimes\bar{v}$$

(2) Now in the case of $G=U_N^+$, by using Theorem 12.33, we obtain right away:
$$PU_N^+=P\widetilde{O_N^+}=PO_N^+$$

(3) We can deduce this as well directly. Indeed, with notations as before, we have:
$$Hom\left((v\otimes v)^k,(v\otimes v)^l\right)=span\left(T_\pi\Big|\pi\in NC_2((\circ\bullet)^k,(\circ\bullet)^l)\right)$$
$$Hom\left((u\otimes\bar{u})^k,(u\otimes\bar{u})^l\right)=span\left(T_\pi\Big|\pi\in \mathcal{NC}_2((\circ\bullet)^k,(\circ\bullet)^l)\right)$$

Now the sets on the right being equal, $PO_N^+\subset PU_N^+$ must be an isomorphism.

\medskip

(4) Finally, the last assertion is something a bit subjective, coming from $PO_N^+=PU_N^+$. For more on this, we refer to the free projective geometry literature.
\end{proof}

Getting now to the classification of irreducible representations, and their fusion rules, good news here too, with everything simplifying in the quantum case. We first have:

\begin{theorem}
The irreducible representations of $O_N^+$  with $N\geq2$ can be labeled by positive integers, $r_k$ with $k\in\mathbb N$, the fusion rules for these representations are
$$r_k\otimes r_l=r_{|k-l|}+r_{|k-l|+2}+\ldots+r_{k+l}$$ 
exactly as the Clebsch-Gordan rules for $SU_2$, and the corresponding dimensions are 
$$\dim r_k=(q^{k+1}-q^{-k-1})/(q-q^{-1})$$
with $q+q^{-1}=N$. And with this being something far simpler than for $O_N$.
\end{theorem}

\begin{proof}
The proof is nearly identical to the proof for $SU_2$ from chapter 10, with the needed estimate there, involving the Catalan numbers, coming for $O_N^+$ as follows:
\begin{eqnarray*}
\int_{O_N^+}\chi^{2k}
&=&\dim(Fix(u^{\otimes 2k}))\\
&=&\dim\left(span\left(T_\pi|\pi\in NC_2(2k)\right)\right)\\
&\leq&|NC_2(2k)|\\
&=&C_k
\end{eqnarray*}

Thus, by reasoning as in chapter 10, we can construct the irreducible representations by recurrence, with the conclusions in the statement. As for the last assertion, this is something subjective, but try proving me wrong, by quicky dealing with $O_N$.
\end{proof}

Regarding now the quantum group $S_N^+$, we have a similar result for it, as follows:

\begin{theorem}
The irreducible representations of $S_N^+$ with $N\geq4$ can be labeled by positive integers, $r_k$ with $k\in\mathbb N$, the fusion rules for these representations are
$$r_k\otimes r_l=r_{|k-l|}+r_{|k-l|+1}+\ldots+r_{k+l}$$ 
exactly as the Clebsch-Gordan rules for $SO_3$, and the corresponding dimensions are 
$$\dim r_k=(q^{k+1}-q^{-k})/(q-1)$$
with $q+q^{-1}=N-2$. And with this being something far simpler than for $S_N$.
\end{theorem}

\begin{proof}
Similar story here, with however a subtlety coming at the technical level, due to the fact, explained in the case of $SO_3$ in chapter 10, that we need here an equality involving the Catalan numbers, instead of just an estimate. But this equality can be deduced as follows, by using easiness, and the Di Francesco results from chapter 6:
\begin{eqnarray*}
\int_{S_N^+}\chi^k
&=&\dim(Fix(u^{\otimes k}))\\
&=&\dim\left(span\left(T_\pi|\pi\in NC(k)\right)\right)\\
&=&|NC(k)|\\
&=&C_k
\end{eqnarray*}

As for the last assertion, as before with $O_N$ vs $O_N^+$, prove me wrong here.
\end{proof}

Finally, let us record as well the result for $U_N^+$, which is as follows:

\begin{theorem}
The irreducible representations of the quantum group $U_N^+$ with $N\geq2$ can be indexed by $\mathbb N*\mathbb N$, with the fusion rules being as follows,
$$r_k\otimes r_l=\sum_{k=xy,l=\bar{y}z}r_{xz}$$
and the corresponding dimensions $\dim r_k$ can be computed by recurrence. And with this being something far simpler than for $U_N$.
\end{theorem}

\begin{proof}
This can be proved either directly, as for $O_N^+$, with the needed estimate coming from easiness, or via Theorem 12.33, with the fusion rules in the statement being some sort of ``free complexification'' of the Clebsch-Gordan rules for $O_N^+$. See \cite{ba2}.
\end{proof}

As a conclusion now, we have strong evidence for Fact 12.2 and Principle 12.3. Want it or not, that statements definitely hold, and we will have to live with this.

\bigskip

In practice, all this goes back to work of Wang and mine from the 1990s, and many things have happened, since. To start with, Theorems 12.35 and 12.36 suggest stronger ties of type $O_N^+\sim SU_2$ and $S_N^+\sim SO_3$, and there is nowadays a whole variety of quite satisfactory results in this sense, of both algebraic and analytic nature.

\bigskip

As in what regards $U_N^+$, here again countless things coming on top of Theorem 12.37 are currently known. Finally, the same goes for $PO_N^+=PU_N^+$, which has same fusion rules as $SO_3$, and with this coming from something of type $PO_N^+=PU_N^+\sim SO_3$. For more on all this, you can have a look at my quantum group book \cite{ba2}.

\section*{12d. Quantum reflections}

Getting now to quantum reflections, things here are more tricky. Following \cite{bbc}, the free analogues of the reflection groups $H_N^s$ can be constructed as follows:

\begin{definition}
The algebra $C(H_N^{s+})$ is the universal $C^*$-algebra generated by $N^2$ normal elements $u_{ij}$, subject to the following relations,
\begin{enumerate}
\item $u=(u_{ij})$ is unitary,

\item $u^t=(u_{ji})$ is unitary,

\item $p_{ij}=u_{ij}u_{ij}^*$ is a projection,

\item $u_{ij}^s=p_{ij}$,
\end{enumerate}
with Woronowicz algebra maps $\Delta,\varepsilon,S$ constructed by universality.
\end{definition}

Here we allow the value $s=\infty$, with the convention that the last axiom simply disappears in this case. Observe also that at $s<\infty$ the normality condition is actually redundant, because a partial isometry $a$ subject to the relation $aa^*=a^s$ is normal. 

\bigskip

As a first result now, making the connection with $H_N^s$, we have:

\begin{theorem}
We have an inclusion of quantum groups
$$H_N^s\subset H_N^{s+}$$
which is a liberation, in the sense that the classical version of $H_N^{s+}$, obtained by dividing by the commutator ideal, is the group $H_N^s$.
\end{theorem}

\begin{proof}
This follows as before for $O_N\subset O_N^+$ or for $S_N\subset S_N^+$, by using the Gelfand theorem, applied to the quotient of $C(H_N^{s+})$ by its commutator ideal.
\end{proof}

Let us prove now that $H_N^{s+}$ with $s<\infty$ is a quantum permutation group. For this purpose, we must change the fundamental representation. Let us start with:

\begin{definition}
A $(s,N)$-sudoku matrix is a magic unitary of size $sN$, of the form
$$m=\begin{pmatrix}
a^0&a^1&\ldots&a^{s-1}\\
a^{s-1}&a^0&\ldots&a^{s-2}\\
\vdots&\vdots&&\vdots\\
a^1&a^2&\ldots&a^0
\end{pmatrix}$$
where $a^0,\ldots,a^{s-1}$ are $N\times N$ matrices.
\end{definition}

The basic examples of such matrices come from the group $H_n^s$. Indeed, with $w=e^{2\pi i/s}$, each of the $N^2$ matrix coordinates $u_{ij}:H_N^s\to\mathbb C$ takes values in the following set:
$$S=\{0\}\cup\{1,w,\ldots,w^{s-1}\}$$

Thus, this coordinate function $u_{ij}:H_N^s\to\mathbb C$ decomposes as follows:
$$u_{ij}=\sum_{r=0}^{s-1}w^ra^r_{ij}$$

Here each function $a^r_{ij}$ is a function taking values in $\{0,1\}$, and so is a projection in the $C^*$-algebra sense, and it follows from definitions that these projections form indeed a sudoku matrix. Now with this notion in hand, we have the following result:

\begin{theorem}
The following happen:
\begin{enumerate}
\item The algebra $C(H_N^s)$ is isomorphic to the universal commutative $C^*$-algebra generated by the entries of a $(s,N)$-sudoku matrix.

\item The algebra $C(H_N^{s+})$ is isomorphic to the universal $C^*$-algebra generated by the entries of a $(s,N)$-sudoku matrix.
\end{enumerate}
\end{theorem}

\begin{proof}
The first assertion follows from the second one, via Theorem 12.39. In order to prove the second assertion, consider the universal algebra in the statement, namely:
$$A=C^*\left(a_{ij}^p\ \Big\vert \left(a^{q-p}_{ij}\right)_{pi,qj}=(s,N)-\mbox{sudoku }\right)$$

Consider also the algebra $C(H_N^{s+})$. According to Definition 12.38, this is presented by certain relations $R$, that we will call here level $s$ cubic conditions:
$$C(H_N^{s+})=C^*\left(u_{ij}\ \Big\vert\  u=N\times N\mbox{ level $s$ cubic }\right)$$

We will construct a pair of inverse morphisms between these algebras.

\medskip

(1) Our first claim is that $U_{ij}=\sum_pw^{-p}a^p_{ij}$ is a level $s$ cubic unitary. Indeed, by using the sudoku condition, the verification of (1-4) in Definition 12.38 is routine.

\medskip

(2) Our second claim is that the elements $A^p_{ij}=\frac{1}{s}\sum_rw^{rp}u^r_{ij}$, with the convention $u_{ij}^0=p_{ij}$, form a level $s$ sudoku unitary. Once again, the proof here is routine.

\medskip

(3) According to the above, we can define a morphism $\Phi:C(H_N^{s+})\to A$ by the formula $\Phi(u_{ij})=U_{ij}$, and a morphism $\Psi:A\to C(H_N^{s+})$ by the formula $\Psi(a^p_{ij})=A^p_{ij}$.

\medskip

(4) We check now the fact that $\Phi,\Psi$ are indeed inverse morphisms:
\begin{eqnarray*}
\Psi\Phi(u_{ij})
&=&\sum_pw^{-p}A^p_{ij}\\
&=&\frac{1}{s}\sum_pw^{-p}\sum_rw^{rp}u_{ij}^r\\
&=&\frac{1}{s}\sum_{pr}w^{(r-1)p}u_{ij}^r\\
&=&u_{ij}
\end{eqnarray*}

As for the other composition, we have the following computation:
\begin{eqnarray*}
\Phi\Psi(a^p_{ij})
&=&\frac{1}{s}\sum_rw^{rp}U_{ij}^r\\
&=&\frac{1}{s}\sum_rw^{rp}\sum_qw^{-rq}a_{ij}^q\\
&=&\frac{1}{s}\sum_qa_{ij}^q\sum_rw^{r(p-q)}\\
&=&a^p_{ij}
\end{eqnarray*}

Thus we have an isomorphism $C(H_N^{s+})=A$, as claimed.
\end{proof}

In fact, we can say more about the quantum permutation group meaning of $H_N^{s+}$, in analogy with what we know from before regarding the group $H_N^s$, as follows:

\begin{theorem}
The disjoint union of $N$ copies of the oriented $s$-cycle
$$
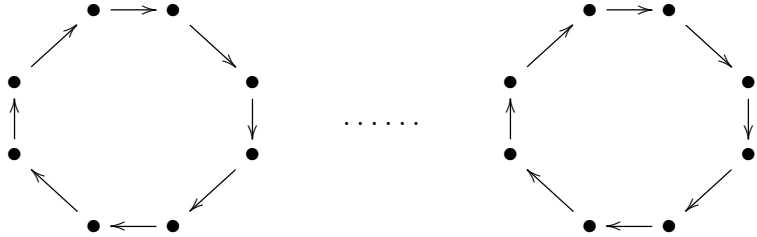
\xymatrix@R=15pt@C=17pt{
&\bullet\ar[r]&\bullet\ar[dr]\\
\bullet\ar[ur]&&&\bullet\ar[d]\\
\bullet\ar[u]&&&\bullet\ar[dl]\\
&\bullet\ar[ul]&\bullet\ar[l]}\qquad
\xymatrix@R=13pt@C=17pt{\\ \\ \ldots\ldots}
\qquad
\xymatrix@R=15pt@C=17pt{
&\bullet\ar[r]&\bullet\ar[dr]\\
\bullet\ar[ur]&&&\bullet\ar[d]\\
\bullet\ar[u]&&&\bullet\ar[dl]\\
&\bullet\ar[ul]&\bullet\ar[l]}$$
has $H_N^s$ as symmetry group, and $H_N^{s+}$ as quantum symmetry group.
\end{theorem}

\begin{proof}
In what regards the $H_N^s$ assertion, this is something that we know well, and which is obvious. Regarding now the $H_N^{s+}$ assertion, our claim here is that we have the following formula, with $d$ being the adjacency matrix of the graph in the statement:
$$C(H_N^{s+})=C(S_{Ns}^+)\big/\big<du=ud\Big>$$

In order to prove this, observe that the adjacency matrix $d$ is given by:
$$d=
\begin{pmatrix}
0&1_N&0&\ldots&0\\
0&0&1_N&\ldots&0\\
\vdots&\vdots&&\ddots&\\
0&0&0&\ldots&1_N\\
1_N&0&0&\ldots&0
\end{pmatrix}$$

Now given a $sN\times sN$ magic unitary, commutation with $d$ means that the matrix is circulant. Thus, we obtain the sudoku relations from Definition 12.40, as desired.
\end{proof}

Next in line, we can talk about wreath products, the result here being as follows:

\begin{theorem}
We have the following results:
\begin{enumerate}
\item $H_N^s=\mathbb Z_s\wr S_N$, with $\wr$ being a wreath product.

\item $H_N^{s+}=\mathbb Z_s\wr_*S_N^+$, with $\wr_*$ being a free wreath product.
\end{enumerate}
\end{theorem}

\begin{proof}
In what regards the first result, this is something that we know well. Regarding the second result, the free wreath product of two quantum permutation groups $G\subset S_s^+$, $H\subset S_N^+$, with fundamental corepresentations $u,v$, can be constructed via:
$$C(G\wr_*H)=(C(G)^{*N}*C(H))\big/\big<[u_{ij}^{(a)},v_{ab}]=0\big>$$

But in our case, it is routine to show that $G=\mathbb Z_s$ and $H=S_N^+$ produce in this way the quantum group $G\wr_*H=H_N^{s+}$. Alternatively, all this can come as well from Theorem 12.42, via the following formulae, valid for any connected graph $X$:
$$G(X\ldots X)=G(X)\wr S_N\quad,\quad 
G^+(X\ldots X)=G^+(X)\wr_*S_N^+$$

Thus, one way or another, we are led to the conclusions in the statement.
\end{proof}

Regarding now the easiness property of $H_N^s,H_N^{s+}$, we have here:

\begin{theorem}
The quantum groups $H_N^s,H_N^{s+}$ are easy, the corresponding categories
$$P^s\subset P\quad,\quad 
NC^s\subset NC$$
consisting of partitions satisfying $\#\circ=\#\bullet(s)$, as a weighted sum, in each block.
\end{theorem}

\begin{proof}
Here the first result is something that we know well, from chapter 11, and the proof of the second result is nearly identical, based on the same computations.
\end{proof}

Next in line, remember the fundamental square of easy groups, from chapter 11? We can extend now that square into a cube, even more fundamental, as follows:

\begin{theorem}
We have quantum rotation and reflection groups, as follows,
$$
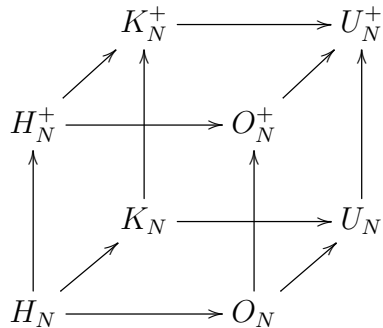
\xymatrix@R=18pt@C=17pt{
&K_N^+\ar[rr]&&U_N^+\\
H_N^+\ar[rr]\ar[ur]&&O_N^+\ar[ur]\\
&K_N\ar[rr]\ar[uu]&&U_N\ar[uu]\\
H_N\ar[uu]\ar[ur]\ar[rr]&&O_N\ar[uu]\ar[ur]
}$$
which are all easy, the corresponding categories of partitions being as follows,
$$
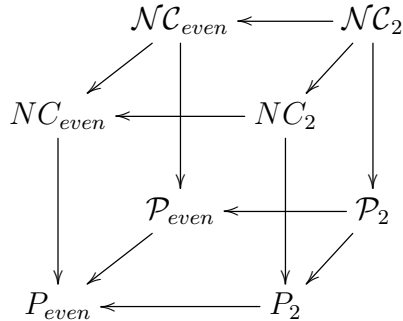
\xymatrix@R=19pt@C=2.5pt{
&\mathcal{NC}_{even}\ar[dl]\ar[dd]&&\mathcal {NC}_2\ar[dl]\ar[ll]\ar[dd]\\
NC_{even}\ar[dd]&&NC_2\ar[dd]\ar[ll]\\
&\mathcal P_{even}\ar[dl]&&\mathcal P_2\ar[dl]\ar[ll]\\
P_{even}&&P_2\ar[ll]
}$$
with on top, the symbol $NC$ standing everywhere for noncrossing partitions.
\end{theorem}

\begin{proof}
This follows indeed by putting together the various results above.
\end{proof}

The same comments as those in the classical case apply, the idea being that the above cubes provide us with some very useful 3D orientation inside quantum algebra. In fact, looking at what has been done recently, say in the last 20 years, of rather fundamental nature, either algebraic or analytic, it is hard to find anything not inspired by the above cubes. Moral of the story I guess, 3D orientation matters, even for math nerds.

\bigskip

Now back to $H_N^{s+}$, regarding the fusion rules, the result here is as follows:

\begin{theorem}
The irreducible representations of $H_N^{s+}$ with $N\geq4$ can be labeled $r_x$ with $x\in<\mathbb Z_s>$, the fusion rules being as follows,
$$r_x\otimes r_y=\sum_{x=vz,y=\bar{z}w}r_{vw}+r_{v\cdot w}$$
with the dot product on the right being given by 
$$(i_1\ldots i_k)\cdot (j_1\ldots j_l)=i_1\ldots i_{k-1}(i_k+j_1)j_2\ldots j_l$$
and with this being far simpler than for the group $H_N^s$ itself.
\end{theorem}

\begin{proof}
As before with other results of this type, this is something of Clebsch-Gordan nature, that we already know, from Theorem 12.36, to hold at $s=1$. Observe also the similarity with Theorem 12.37. As for the proof, this is via our standard methods, namely recurrence, with the needed estimates coming from easiness. For the story with all this, technical details and more, you can have a look at my quantum group book \cite{ba2}.
\end{proof}

And with this, good news, end of our quantum group excursion. We have learned many interesting things, and in particular, we have reached to a good understanding of Fact 12.2 and Principle 12.3, which themselves are something quite remarkable. As for our original Principle 12.1, well, guess that will be left as an exercise for you.

\section*{12e. Exercises}

As just mentioned, now that we understood the simple side of the story, the quantum one, left is the understanding of the complicated part, the classical one:

\begin{exercise}
Compute the fusion rules for $S_N$.
\end{exercise}

\begin{exercise}
Compute the fusion rules for $H_N$.
\end{exercise}

\begin{exercise}
Compute the fusion rules for $K_N$.
\end{exercise}

\begin{exercise}
Compute the fusion rules for $H_N^s$.
\end{exercise}

\begin{exercise}
Compute the fusion rules for $O_N$.
\end{exercise}

\begin{exercise}
Compute the fusion rules for $U_N$.
\end{exercise}

\begin{exercise}
Compute the fusion rules for $B_N$.
\end{exercise}

\begin{exercise}
Compute the fusion rules for $C_N$.
\end{exercise}

As bonus exercise, and no surprise here, learn of course more about quantum groups too. We will be actually back to this, on several occasions, in what follows.

\part{Analytic aspects}

\ \vskip50mm

\begin{center}
{\em Never mind

I'll find someone like you

I wish nothing but the best

For you too}
\end{center}

\chapter{Poisson laws}

\section*{13a. Character laws}

Time for a new topic, analysis. You would probably say, not much analysis to do on the finite groups. But this is wrong, remember the following result, that we learned in chapter 3, which looks quite interesting, and in look for a better understanding:

\index{fixed points}
\index{Poisson law}

\begin{fact}
For the symmetric group $S_N$, the number of fixed points
$$\chi:S_N\to\mathbb N$$
follows with $N\to\infty$ limit the Poisson law $p_1$.
\end{fact}

So, what to do with this? Many things. Let us start with some generalities. We can regard the symmetric group $S_N$ as being the group of permutations of the $N$ coordinate axes of $\mathbb R^N$, which gives an embedding via the permutation matrices, as follows:
$$S_N\subset O_N\quad,\quad \sigma(e_i)=e_{\sigma(i)}$$ 

Thus, we have a representation $S_N\subset O_N$, and the character of this representation $\chi:S_N\to\mathbb R$ is precisely the above variable  counting the fixed points, due to:
\begin{eqnarray*}
\chi(\sigma)
&=&\sum_i<\sigma(e_i),e_i>\\
&=&\sum_i<e_{\sigma(i)},e_i>\\
&=&\#\left\{i\in\{1,\ldots,N\}\Big|\sigma(i)=i\right\}
\end{eqnarray*}

Summarizing, Fact 13.1 is something about characters, as we know them since Part III. But in view of this, we can ask, more generally, the following question:

\begin{question}
Given a subgroup $G\subset U_N$, what is the law of the main character
$$\chi(g)=Tr(g)$$
regarded as random variable $\chi:G\to\mathbb C$? What about asymptotics, when $G=G_N$?
\end{question}

To be more precise here, in what regards the first question, we can certainly view the main character as a variable $\chi:G\to\mathbb C$, and ask for the computation of its law, with respect to the uniform measure on $G$. As for the second question, we have seen all across this book that the main groups usually come in series, $G_N\subset U_N$, and we can ask, a bit as in Fact 13.1, about the $N\to\infty$ asymptotics of our computations.

\bigskip

The above question sounds quite good, and among others, is obviously related to the representation theory that we have been doing in Part III. In fact, we have the following representation theoretic version of the above question, coming from Peter-Weyl:

\begin{question}[version]
Given a subgroup $G\subset_uU_N$, what are the moments
$$\int_G\chi^k=\dim(Fix(u^{\otimes k}))$$
of the main character, $\chi=Tr$? What is $law(\chi)$? What about asymptotics?
\end{question}

To be more precise, this is based on the standard fact, that we know since chapter 9, that integrating a character gives the number of fixed points, which gives:
\begin{eqnarray*}
\int_G\chi^k
&=&\int_G(\chi_u)^k\\
&=&\int_G\chi_{u^{\otimes k}}\\
&=&\dim(Fix(u^{\otimes k}))
\end{eqnarray*}

As a further comment, we can equally view things the other way around, say with solving Question 13.2 first, via some ad-hoc counting techniques, a bit as we did in chapter 3 for the symmetric group $S_N$, and then asking for the representation theory consequences of this, via the formula in Question 13.3. Let us record this question too, as follows: 

\begin{question}[twisted]
What does the knowledge of $law(\chi)$, and its moments
$$\int_G\chi^k=\dim(Fix(u^{\otimes k}))$$
tell us about the representation theory of $G$? For instance, using Frobenius
$$End(u^{\otimes k})\simeq Fix(u^{\otimes k\bar{k}})$$
can we split $u^{\otimes k}$ into irreducibles, via the knowledge of the numbers $\dim(End(u^{\otimes k}))$?
\end{question}

And good question this is too, because we have seen in Part III that many interesting groups $G\subset U_N$, such as the ubiquitous $SU_2$ and $SO_3$, are asking precisely for this type of game, namely computing or estimating the above dimensions, by whatever techniques available, and then construcing the irreducible representations, by recurrence.

\bigskip

Finally, no discussion about all this would be complete without a word on easiness. Indeed, this is something which is certainly present in the context of Fact 13.1, with the symmetric group $S_N$ considered being easy. Also, this can make the last part of Question 13.2 more precise, with the easy groups coming by definition in series $G_N\subset U_N$, as required there. As for the relation with Questions 13.3 and 13.4, certainly no worries here, with easiness being something of deep representation theory nature. So, let us formulae as well an ``easy'' version of our various questions, as follows:

\begin{question}[easy]
Given an easy group $G_N\subset U_N$, coming from a category of partitions $D\subset P$, can we compute the moments of the main character,
$$\int_G\chi^k=\dim\left(span\left(T_\pi\big|\pi\in D(k)\right)\right)$$
then the law of $\chi$, and its $N\to\infty$ asymptotics? Also, how is all this related to the representation theory of $G_N$, when the number $N$ is generic, or fixed?
\end{question}

And again, good question this is, reminding us all sorts of interesting things that we previously did in Part III, in relation with group easiness, and related topics.

\bigskip

As a last twist to the story, getting back now to our starting point, Fact 13.1, we have seen in fact in chapter 3 that, more generally, given $t\in(0,1]$, the number of fixed points of permutations $\sigma\in S_N$ among $\{1,\ldots,[tN]\}$ follows with $N\to\infty$ the Poisson law $p_t$. Which sounds quite interesting, and suggests the following update of our questions:

\begin{question}[parametric]
In the context of our various questions above, what about looking at the truncations of the main character, given by
$$\chi_t=\sum_{i=1}^{[tN]}u_{ii}$$
with $t\in(0,1]$, can we have some theory going for these, both computations and applications? What about other general variables, related to the main character $\chi$?
\end{question}

To be more precise, this is indeed something coming from the above-mentioned computation from chapter 3 for the symmetric group $S_N\subset O_N$, coming via:
\begin{eqnarray*}
\chi_t(\sigma)
&=&\sum_{i=1}^{[tN]}<\sigma(e_i),e_i>\\
&=&\sum_{i=1}^{[tN]}<e_{\sigma(i)},e_i>\\
&=&\#\left\{i\in\{1,\ldots,[tN]\}\Big|\sigma(i)=i\right\}
\end{eqnarray*}

And with this, end of our motivational discussion. In the hope that I convinced you that the questions to be investigated in the present Part IV, which are basically the above ones, are quite interesting. And in case you are still not convinced, and hesitating with keep reading this book, now with the basic algebra aspects learned, here is for you: 

\begin{comment}
On the bottom line, what is to follow is probability theory explained by a fellow algebraist. Which sounds quite nice, unique occasion of learning this.
\end{comment}

So, this was for the story, with probability theory in the finite group context, and in practice now, we better hurry up, because with all this philosophical discussion, the bell is soon ringing, and we have not learned today any new formula or anything. 
\smallskip
$$***$$

\smallskip

Getting started now, we have about 100 pages for discussing the above, and we intend to go quite slowly. Let us start with something simple and elementary, as follows:

\begin{proposition}
For a permutation group $G\subset S_N\subset O_N$, the main character
$$\chi(g)=\sum_ig_{ii}$$
counts the number of fixed points, and the law of this character is given by
$$P\big(\chi=k\big)=\frac{1}{|G|}\cdot\#\left\{g\in G\Big|g\ {\rm has}\ k\ {\rm fixed\ points}\right\}$$
for any $k\in\mathbb N$. Moreover, $k$ ranges in the set $\{0,1,\ldots,N-2,N\}$.
\end{proposition}

\begin{proof}
Here the first assertion is something that we know well, from our discussion before regarding the symmetric group $S_N$, with the computation in general being:
\begin{eqnarray*}
\chi(g)
&=&\sum_i<g(e_i),e_i>\\
&=&\sum_i<e_{g(i)},e_i>\\
&=&\#\left\{i\in\{1,\ldots,N\}\Big|g(i)=i\right\}
\end{eqnarray*}

Regarding the second assertion, this comes from this, trivially, as follows:
\begin{eqnarray*}
P\big(\chi=k\big)
&=&\frac{1}{|G|}\cdot\#\left\{g\in G\Big|\chi(g)=k\right\}\\
&=&\frac{1}{|G|}\cdot\#\left\{g\in G\Big|g\ {\rm has}\ k\ {\rm fixed\ points}\right\}
\end{eqnarray*}

As for the last assertion, the number of fixed points must be certainly an integer $k\in\{0,\ldots,N\}$, but with the value $k=N-1$ being excluded, for obvious reasons. 
\end{proof}

Now let us do some computations, for the various permutation groups $G\subset S_N$ that we know. For the simplest such group, which is the cyclic group $\mathbb Z_N$, we have:

\begin{proposition}
The main character of $\mathbb Z_N\subset O_N$ is given by:
$$\chi(g)=\begin{cases}
0&{\rm if}\ g\neq1\\
N&{\rm if}\ g=1
\end{cases}$$
Thus, at the probabilistic level, we have the following formula,
$$law(\chi)=\left(1-\frac{1}{N}\right)\delta_0+\frac{1}{N}\delta_N$$
telling us that the main character $\chi$ follows a Bernoulli law.
\end{proposition}

\begin{proof}
The first formula is clear, by counting fixed points as in Proposition 3.8, and in the matrix picture too, because the cyclic permutation matrices have 0 on the diagonal, and so 0 as trace, unless the matrix is the identity, having trace $N$. As for the second formula, this is a probabilistic reformulation of the first one.
\end{proof}

For the dihedral group now, which is the next one in our hierarchy, the computation is more interesting, and the final answer is no longer uniform in $N$, as follows:

\index{dihedral group}

\begin{proposition}
For the dihedral group $D_N\subset S_N$ we have
$$law(\chi)=\begin{cases}
\left(\frac{3}{4}-\frac{1}{2N}\right)\delta_0+\frac{1}{4}\delta_2+\frac{1}{2N}\delta_N&(N\ even)\\
&\\
\left(\frac{1}{2}-\frac{1}{2N}\right)\delta_0+\frac{1}{2}\delta_1+\frac{1}{2N}\delta_N&(N\ odd)
\end{cases}$$
and with this law being no longer uniform in $N$.
\end{proposition}

\begin{proof}
The dihedral group $D_N$ consists indeed of:

\medskip

-- $N$ symmetries, having each $1$ fixed point when $N$ is odd, and having 0 or 2 fixed points, distributed $50-50$, when $N$ is even.

\medskip

-- $N$ rotations, each having $0$ fixed points, except for the identity, which is technically a rotation too, and which has $N$ fixed points.

\medskip

Thus, we are led via Proposition 3.8 to the formulae in the statement.
\end{proof}

Quite nice all this, but can we conclude something, out of these computations? In answer, obviously, nothing much, but there is a reason to this, as follows:

\begin{conclusion}
For the cyclic and dihedral groups $\mathbb Z_N\subset D_N\subset S_N$ the character computations are not very interesting, notably due to the lack of interesting $N\to\infty$ asymptotics. But, we can put here at blame the non-easiness of $\mathbb Z_N,D_N$.
\end{conclusion}

To be more precise here, regarding the $N\to\infty$ asymptotics, for the cyclic group $\mathbb Z_N$ we obtain the Dirac mass $\delta_0$, which does not look very exciting, while for the dihedral group $D_N$ we obtain two limiting measures, namely two Bernoulli laws, as follows, which do not look very exciting either, and which, importantly, do not coincide:
$$\mu_{even}=\frac{3}{4}\,\delta_0+\frac{1}{4}\,\delta_2\quad,\quad 
\mu_{odd}=\frac{1}{2}\,\delta_0+\frac{1}{2}\,\delta_1$$

So, not very good all this, and what is to be put at blame? In answer, in what regards $D_N$, what we have obviously comes the fact that $D=(D_N)$ is not a true series, uniform in $N\in\mathbb N$, and thinking a bit at this, from the perspective of representation theory, leads to the conclusion that the problem with $D_N$ comes from its non-easiness. As for $\mathbb Z_N$, here we do have a limit, but with this limit being something not very interesting, we can declare ourselves unhappy too, and put this on the account of non-easiness too.

\bigskip

Moving on, let us investigate now the symmetric group $S_N$ itself. Here we already know the answer, which is interesting, since chapter 3, but since this was some time ago, it is perhaps better to talk again about all this, more in detail. We have here:

\index{random permutation}
\index{derangement}
\index{fixed points}

\begin{theorem}
For the symmetric group $S_N\subset O_N$ we have the formula
$$P(\chi=k)=\frac{1}{k!}\left(1-\frac{1}{1!}+\frac{1}{2!}-\ldots+(-1)^{N-1}\frac{1}{(N-1)!}+(-1)^N\frac{1}{N!}\right)$$
which in the $N\to\infty$ limit simplifies into the formula
$$P(\chi=k)=\frac{1}{ek!}$$
which tells us that the asymptotic law of $\chi$ is the Poisson law $p_1$.
\end{theorem}

\begin{proof}
This follows from a direct fixed point count, as follows:

\medskip

(1) Let us first prove the result at $k=0$. Here we must count the permutations $\sigma\in S_N$ having no fixed points at all, called derangements, and for this purpose, let us set:
$$S_N^i=\left\{\sigma\in S_N\Big|\sigma(i)=i\right\}$$

A derangement is then a permutation $\sigma\in S_N$ not belonging to any of these sets $S_N^i$, and so, according to Proposition 13.8, the probability that we are interested in is:
$$P(\chi=0)=\frac{1}{N!}\Big|\left(S_N^1\cup\ldots\cup S_N^N\right)^c\Big|$$

But we can compute the cardinality of the set on the right by using the inclusion-exclusion principle, and we are led in this way to the following formula, as desired:
\begin{eqnarray*}
P(\chi=0)
&=&\frac{1}{N!}\left(|S_N|-\sum_i|S_N^i|+\sum_{i<j}|S_N^i\cap S_N^j|-\sum_{i<j<k}|S_N^i\cap S_N^j\cap S_N^k|+\ldots\right)\\
&=&\frac{1}{N!}\sum_{k=0}^N(-1)^k\sum_{i_1<\ldots<i_k}|S_N^{i_1}\cap\ldots\cap S_N^{i_k}|\\
&=&\frac{1}{N!}\sum_{k=0}^N(-1)^k\sum_{i_1<\ldots<i_k}(N-k)!\\
&=&\frac{1}{N!}\sum_{k=0}^N(-1)^k\binom{N}{k}(N-k)!\\
&=&\sum_{k=0}^N\frac{(-1)^k}{k!}
\end{eqnarray*}

(2) In the general case now, we have to count the permutations $\sigma\in S_N$ having exactly $k$ points. Since having such a permutation amounts in choosing $k$ points among $1,\ldots,N$, and then permuting the $N-k$ points left, without fixed points allowed, we have:
\begin{eqnarray*}
\#\left\{\sigma\in S_N\Big|\chi(\sigma)=k\right\}
&=&\binom{N}{k}\#\left\{\sigma\in S_{N-k}\Big|\chi(\sigma)=0\right\}\\
&=&\frac{N!}{k!(N-k)!}\#\left\{\sigma\in S_{N-k}\Big|\chi(\sigma)=0\right\}\\
&=&N!\times\frac{1}{k!}\times\frac{\#\left\{\sigma\in S_{N-k}\Big|\chi(\sigma)=0\right\}}{(N-k)!}
\end{eqnarray*}

Now by dividing everything by $N!$, we obtain from this the following formula:
\begin{eqnarray*}
P(\chi=k)
&=&\frac{\#\left\{\sigma\in S_N\Big|\chi(\sigma)=k\right\}}{N!}\\
&=&\frac{1}{k!}\times\frac{\#\left\{\sigma\in S_{N-k}\Big|\chi(\sigma)=0\right\}}{(N-k)!}\\
&=&\frac{1}{k!}\left(1-\frac{1}{1!}+\frac{1}{2!}-\ldots+(-1)^{N-1}\frac{1}{(N-1)!}+(-1)^N\frac{1}{N!}\right)
\end{eqnarray*}

(3) Finally, the asymptotics are clear, coming from $1/e=\sum_r(-1)^r/r!$, and the last assertion follows from the definition of the Poisson law, namely $P(k)=1/ek!$. 
\end{proof}

More generally, we have the following result, that we also know from chapter 3:

\begin{theorem}
For the symmetric group $S_N\subset O_N$, the truncated character
$$\chi_t(g)=\sum_{i=1}^{[tN]}g_{ii}$$
with $t\in(0,1]$ counts the corresponding truncated number of fixed points,
$$\chi_t(g)=\#\left\{i\in\{1,\ldots,[tN]\Big|g(i)=i\right\}$$
and follows the Poisson law $p_t$, in the $N\to\infty$ limit.
\end{theorem}

\begin{proof}
As before in the proof of Theorem 13.12, we get by inclusion-exclusion:
\begin{eqnarray*}
P(\chi_t=0)
&=&\frac{1}{N!}\sum_{r=0}^{[tN]}(-1)^r\sum_{k_1<\ldots<k_r<[tN]}|S_N^{k_1}\cap\ldots\cap S_N^{k_r}|\\
&=&\frac{1}{N!}\sum_{r=0}^{[tN]}(-1)^r\binom{[tN]}{r}(N-r)!\\
&=&\sum_{r=0}^{[tN]}\frac{(-1)^r}{r!}\cdot\frac{[tN]!(N-r)!}{N!([tN]-r)!}
\end{eqnarray*}

Now with $N\to\infty$, we obtain from this the following estimate:
$$P(\chi_t=0)
\simeq\sum_{r=0}^{[tN]}\frac{(-1)^r}{r!}\cdot t^r
\simeq\frac{1}{e^t}$$

More generally, by counting the permutations $\sigma\in S_N$ having exactly $k$ fixed points among $1,\ldots,[tN]$, again as in the proof of Theorem 13.12, we obtain:
$$P(\chi_t=k)\simeq\frac{t^k}{e^tk!}$$

Thus, we obtain in the limit a Poisson law of parameter $t$, as stated.
\end{proof}

\section*{13b. Poisson variables}

Done with the symmetric group? Not really, because the results that we have, namely Theorem 13.12, generalized into Theorem 13.13, suggest many interesting questions. We intend to discuss these questions, or at least some of them, in the present section.

\bigskip

As a first natural question, is upgrading from Theorem 13.12, involving characters which are something quite conceptual, to Theorem 13.13, involving that truncated characters which look like something quite specialized, something really worth it?

\bigskip

In answer, a bit of general probability theory knowledge can help us, as follows:

\begin{principle}
The main probability measures, as for instance the Poisson or normal laws, do not come alone, but rather in one-parameter families 
$$\mu=\big\{\mu_t\big|t>0\big\}$$
which usually form convolution semigroups, and with this being something essential to their understanding. So, definitely green light for the truncated characters.
\end{principle}

Which sounds quite interesting, what we have here is probability theory teaching us that, despite our escalation of representation theory methods in chapters 9, 10, 11, 12, we might have missed something there. Never too old, for learning new things.

\bigskip

In view of this, it is perhaps wiser to temporarily forget about groups, and have a look at theoretical probability, who knows what new things can we learn here, that can be of help. Regarding the Poisson laws $p_t$, their basic theory can be summarized as follows:

\begin{theorem}
The Poisson laws $p_t$ with $t>0$, given by
$$p_t=e^{-t}\sum_{k\geq0}\frac{t^k}{k!}\,\delta_k$$
have the following properties:
\begin{enumerate}
\item Their mean is $E=t$, and the variance is $V=t$.

\item We have the convolution semigroup formula $p_s*p_t=p_{s+t}$.

\item They appear as formal exponentials, $p_t=\sum_kt^k(\delta_1-\delta_0)^{*k}/k!$.

\item The Fourier transform is given by $F_{p_t}(y)=\exp((e^{iy}-1)t)$.

\item They appear via Poisson limits, $((1-\frac{t}{n})\delta_0+\frac{t}{n}\delta_1)^{*n}\to p_t$.

\item The moments of $p_1$ are the Bell numbers, $M_k(p_1)=|P(k)|$.

\item In fact, $M_k(p_t)=\sum_{\pi\in P(k)}t^{|\pi|}$, with $|.|$ being the number of blocks.
\end{enumerate}
\end{theorem}

\begin{proof}
Many things going on here, with the computations being as follows:

\medskip

(1) Regarding the mean, this can be computed as follows:
\begin{eqnarray*}
E
&=&e^{-t}\sum_{k\geq1}\frac{t^k}{k!}\cdot k\\
&=&te^{-t}\sum_{k\geq1}\frac{t^{k-1}}{(k-1)!}\\
&=&t
\end{eqnarray*}

For the variance, we first compute the second moment, as follows:
\begin{eqnarray*}
M_2
&=&e^{-t}\sum_{k\geq1}\frac{t^k}{k!}\cdot k^2\\
&=&te^{-t}\sum_{k\geq1}\frac{t^{k-1}k}{(k-1)!}\\
&=&te^{-t}\sum_{l\geq0}\frac{t^ll}{l!}+te^{-t}\sum_{l\geq0}\frac{t^l}{l!}\\
&=&t^2e^{-t}\sum_{m\geq0}\frac{t^m}{m!}+t\\
&=&t^2+t
\end{eqnarray*}

Thus, we obtain $V=M_2-E^2=(t^2+t)-t^2=t$, as claimed.

\medskip

(2) Regarding the convolution semigroup property, this comes as follows:
\begin{eqnarray*}
p_s*p_t
&=&e^{-s}\sum_k\frac{s^k}{k!}\,\delta_k*e^{-t}\sum_l\frac{t^l}{l!}\,\delta_l\\
&=&e^{-s-t}\sum_n\delta_n\sum_{k+l=n}\frac{s^kt^l}{k!l!}\\
&=&e^{-s-t}\sum_n\frac{\delta_n}{n!}\sum_{k+l=n}\frac{n!}{k!l!}s^kt^l\\\
&=&e^{-s-t}\sum_n\frac{(s+t)^n}{n!}\,\delta_n\\
&=&p_{s+t}
\end{eqnarray*}

(3) Next, regarding the formal exponential property, this comes from:
\begin{eqnarray*}
\sum_k\frac{t^k(\delta_1-\delta_0)^{*k}}{k!}
&=&\sum_k\frac{t^k}{k!}\sum_{r+s=k}(-1)^s\frac{k!}{r!s!}\delta_r\\
&=&\sum_kt^k\sum_{r+s=k}(-1)^s\frac{\delta_r}{r!s!}\\
&=&\sum_r\frac{t^r\delta_r}{r!}\sum_s\frac{(-1)^st^s}{s!}\\
&=&\frac{1}{e^t}\sum_r\frac{t^r\delta_r}{r!}\\
&=&p_t
\end{eqnarray*}

(4) For the Fourier transform, we have indeed the following computation:
\begin{eqnarray*}
F_{p_t}(y)
&=&e^{-t}\sum_k\frac{t^k}{k!}F_{\delta_k}(y)\\
&=&e^{-t}\sum_k\frac{t^k}{k!}\,e^{iky}\\
&=&e^{-t}\sum_k\frac{(e^{iy}t)^k}{k!}\\
&=&\exp(-t)\exp(e^{iy}t)\\
&=&\exp\left((e^{iy}-1)t\right)
\end{eqnarray*}

Observe that this provides us with an alternative proof for (2), via the well-known fact that the logarithm of the Fourier transform linearizes the convolution.

\medskip

(5) Let us establish now the Poisson Limit Theorem, stating that we have:
$$\left(\left(1-\frac{t}{n}\right)\delta_0+\frac{t}{n}\,\delta_1\right)^{*n}\to p_t$$

For this purpose, let us denote by $\mu_n$ the measure under the convolution sign. We have the following computation, for the Fourier transform of the limit: 
\begin{eqnarray*}
F_{\delta_r}(y)=e^{iry}
&\implies&F_{\mu_n}(y)=\left(1-\frac{t}{n}\right)+\frac{t}{n}\,e^{iy}\\
&\implies&F_{\mu_n^{*n}}(y)=\left(\left(1-\frac{t}{n}\right)+\frac{t}{n}\,e^{iy}\right)^n\\
&\implies&F_{\mu_n^{*n}}(y)=\left(1+\frac{(e^{iy}-1)t}{n}\right)^n\\
&\implies&F(y)=\exp\left((e^{iy}-1)t\right)
\end{eqnarray*}

Thus, we obtain indeed the Fourier transform of $p_t$, as desired.

\medskip

(6) Getting now to the moments of $p_1$, we know from (1) that we have:
$$M_1=1\quad,\quad M_2=2$$

In general, the moments of $p_1$ are given by the following formula:
$$M_k=\frac{1}{e}\sum_r\frac{r^k}{r!}$$

We therefore have the following recurrence formula for these moments:
\begin{eqnarray*}
M_{k+1}
&=&\frac{1}{e}\sum_r\frac{(r+1)^{k+1}}{(r+1)!}\\
&=&\frac{1}{e}\sum_r\frac{r^k}{r!}\left(1+\frac{1}{r}\right)^k\\
&=&\frac{1}{e}\sum_r\frac{r^k}{r!}\sum_s\binom{k}{s}r^{-s}\\
&=&\sum_s\binom{k}{s}\cdot\frac{1}{e}\sum_r\frac{r^{k-s}}{r!}\\
&=&\sum_s\binom{k}{s}M_{k-s}
\end{eqnarray*}

With this done, let us try now to find a recurrence for the Bell numbers, $B_k=|P(k)|$. A partition of $\{1,\ldots,k+1\}$ appears by choosing $s$ neighbors for $1$, among the $k$ numbers available, and then partitioning the $k-s$ elements left. Thus, we have:
$$B_{k+1}=\sum_s\binom{k}{s}B_{k-s}$$

Also, we have $B_1=1$ and $B_2=2$. Thus we obtain $M_k=B_k$, as claimed.

\medskip

(7) The moments of the Poisson law $p_t$ with $t>0$ are given by:
$$M_k=e^{-t}\sum_r\frac{t^rr^k}{r!}$$

We therefore have the following recurrence formula for these moments:
\begin{eqnarray*}
M_{k+1}
&=&e^{-t}\sum_r\frac{t^{r+1}(r+1)^{k+1}}{(r+1)!}\\
&=&e^{-t}\sum_r\frac{t^{r+1}r^k}{r!}\left(1+\frac{1}{r}\right)^k\\
&=&e^{-t}\sum_r\frac{t^{r+1}r^k}{r!}\sum_s\binom{k}{s}r^{-s}\\
&=&\sum_s\binom{k}{s}\cdot e^{-t}\sum_r\frac{t^{r+1}r^{k-s}}{r!}\\
&=&t\sum_s\binom{k}{s}M_{k-s}
\end{eqnarray*}

Regarding now the initial values, we know from (1) that we have:
$$M_1=t\quad,\quad M_2=t+t^2$$

On the other hand, consider the numbers in the statement, namely:
$$S_k=\sum_{\pi\in P(k)}t^{|\pi|}$$

Since a partition of $\{1,\ldots,k+1\}$ appears by choosing $s$ neighbors for $1$, among the $k$ numbers available, and then partitioning the $k-s$ elements left, we have:
$$S_{k+1}=t\sum_s\binom{k}{s}S_{k-s}$$

As for the initial values of these numbers, these are $S_1=t$, $S_2=t+t^2$. Thus the initial values coincide, and so these numbers are the moments of $p_t$, as stated.
\end{proof}

And with this, good news, end of our basic probability learning. Theorem 13.15 contains a wealth of information, that help us, in relation with our algebra questions.

\section*{13c. Characters, revised} 

Back now to the symmetric group $S_N$, we would like to present some alternative proofs for Theorems 13.12 and 13.13, which are all instructive. To start with, we can replace the inclusion-exclusion principle with something more conceptual and analytic, namely:

\index{polynomial integrals}

\begin{theorem}
The integration over $S_N$ is given by the formula
$$\int_{S_N}g_{i_1j_1}\ldots g_{i_kj_k}=\begin{cases}
\frac{(N-|\ker i|)!}{N!}&{\rm if}\ \ker i=\ker j\\
0&{\rm otherwise}
\end{cases}$$
where $\ker i$ denotes the partition of $\{1,\ldots,k\}$ whose blocks collect the equal indices of the multi-index $i=(i_1\ldots i_k)$, and where $|.|$ is the number of blocks.
\end{theorem}

\begin{proof}
According to our usual conventions for the embedding $S_N\subset O_N$, coming by permuting the standard basis of $\mathbb R^N$, the coordinate functions on $S_N$ are given by:
$$g_{ij}(\sigma)=\delta_{\sigma(j)i}$$

Thus, the integrals in the statement are given by the following formula:
$$\int_{S_N}g_{i_1j_1}\ldots g_{i_kj_k}=\frac{1}{N!}\#\left\{\sigma\in S_N\Big|\sigma(j_1)=i_1,\ldots,\sigma(j_k)=i_k\right\}$$

Now observe that the existence of a permutation $\sigma\in S_N$ as above requires:
$$i_m=i_n\iff j_m=j_n$$

We conclude that the above integral must vanish in the following case:
$$\ker i\neq\ker j$$

In the remaining case, $\ker i=\ker j$, if we denote by $b\in\{1,\ldots,k\}$ the number of blocks of this partition $\ker i=\ker j$, we have $N-b$ points to be sent bijectively to $N-b$ points, and so $(N-b)!$ solutions, and so the integral is $\frac{(N-b)!}{N!}$, as claimed.
\end{proof}

We can now recover our previous character computations, as follows:

\begin{theorem}
For the symmetric group $S_N\subset O_N$, regarded as group of matrices, $S_N\subset O_N$, via the standard permutation matrices, the truncated character
$$\chi_t(g)=\sum_{i=1}^{[tN]}g_{ii}$$
counts the number of fixed points among $\{1,\ldots,[tN]\}$, and its law with respect to the counting measure becomes, with $N\to\infty$, a Poisson law of parameter $t$. 
\end{theorem}

\begin{proof}
With $S_{kb}$ being the Stirling numbers, counting the partitions of $\{1,\ldots,k\}$ having exactly $b$ blocks, we have indeed the following computation:
\begin{eqnarray*}
\int_{S_N}\chi_t^k
&=&\sum_{i_1\ldots i_k=1}^{[tN]}\int_{S_N}g_{i_1i_1}\ldots g_{i_ki_k}\\
&=&\sum_{\pi\in P(k)}\frac{[tN]!}{([tN]-|\pi|!)}\cdot\frac{(N-|\pi|!)}{N!}\\
&=&\sum_{b=1}^{[tN]}\frac{[tN]!}{([tN]-b)!}\cdot\frac{(N-b)!}{N!}\cdot S_{kb}
\end{eqnarray*}

In particular with $N\to\infty$ we obtain the following formula:
$$\lim_{N\to\infty}\int_{S_N}\chi_t^k=\sum_{b=1}^kS_{kb}t^b$$

But this is the $k$-th moment of the Poisson law $p_t$, as desired.
\end{proof}

As another result now regarding $S_N$, here is a useful related formula:

\begin{theorem}
We have the law formula
$${\rm law}(g_{11}+\ldots +g_{ss})=\frac{s!}{N!}\sum_{p=0}^s\frac{(N-p)!}{(s-p)!}
\cdot\frac{\left(\delta_1-\delta_0\right)^{*p}}{p!}$$ 
where $g_{ij}$ are the standard coordinates of $S_N\subset O_N$.
\end{theorem}

\begin{proof}
We have the following moment formula, where $m_f$ is the number of permutations of $\{1,\ldots ,N\}$ having exactly $f$ fixed points in the set $\{1,\ldots ,s\}$: 
$$\int_{S_N}(u_{11}+\ldots +u_{ss})^k=\frac{1}{N!}\sum_{f=0}^sm_ff^k$$

Thus the law in the statement, say $\nu_{sN}$, is the following average
of Dirac masses:
$$\nu_{sN}=\frac{1}{N!}\sum_{f=0}^s m_{f}\,\delta_f$$

Now observe that the permutations contributing to $m_f$ are obtained by choosing $f$
points in the set $\{1,\ldots ,s\}$, then by permuting the remaining $N-f$ points in $\{1,\ldots ,n\}$ in such a way that there is no fixed point in $\{1,\ldots,s\}$. But these latter permutations are counted as follows: we start with all permutations, we substract those having one fixed point, we add those having two fixed points, and so on. We obtain in this way:
\begin{eqnarray*}
\nu_{sN}
&=&\frac{1}{N!}\sum_{f=0}^s\begin{pmatrix}s\\
f\end{pmatrix}\left(\sum_{k=0}^{s-f}(-1)^k
\begin{pmatrix}s-f\\ k\end{pmatrix}(N-f-k)!\right)\,\delta_f\\
&=&\sum_{f=0}^s\sum_{k=0}^{s-f}(-1)^k\frac{1}{N!}\cdot
\frac{s!}{f!(s-f)!}\cdot\frac{(s-f)!(N-f-k)!}{k!(s-f-k)!}\,\delta_f\\
&=&\frac{s!}{N!}\sum_{f=0}^s\sum_{k=0}^{s-f}\frac{(-1)^k(N-f-k)!}{f!k!(s-f-k)!}\,\delta_f
\end{eqnarray*}

We can proceed as follows, by using the new index $p=f+k$:
\begin{eqnarray*}
\nu_{sN}
&=&\frac{s!}{N!}\sum_{p=0}^s\sum_{k=0}^{p}\frac{(-1)^k
(N-p)!}{(p-k)!k!(s-p)!}\,\delta_{p-k}\\
&=&\frac{s!}{N!}\sum_{p=0}^s\frac{(N-p)!}{(s-p)!p!}
\sum_{k=0}^{p}(-1)^k\begin{pmatrix}p\\
k\end{pmatrix}\,\delta_{p-k}\\
&=&\frac{s!}{N!}\sum_{p=0}^s\frac{(N-p)!}{(s-p)!}\cdot
\frac{\left(\delta_1-\delta_0\right)^{*p}}{p!}
\end{eqnarray*}

Here $*$ is convolution of real measures, and the assertion follows.
\end{proof}

As a consequence of Theorem 13.18, improving some previous results, we have:

\begin{theorem}
Let $g_{ij}$ be the standard coordinates of $C(S_N)$.
\begin{enumerate}
\item $u_{11}+\ldots +u_{ss}$ with $s=o(N)$ is a projection of trace $s/N$. 

\item $u_{11}+\ldots +u_{ss}$ with $s=tN+o(N)$ is Poisson of parameter $t$.
\end{enumerate}
\end{theorem}

\begin{proof}
We can use indeed the formula in Theorem 13.18, as follows:

\medskip

(1) With $s$ fixed and $N\to\infty$ we have the following estimate:
\begin{eqnarray*}
&&{\rm law}(u_{11}+\ldots +u_{ss})\\
&=&\sum_{p=0}^s\frac{(N-p)!}{N!}\cdot\frac{s!}{(s-p)!}
\cdot\frac{\left(\delta_1-\delta_0\right)^{*p}}{p!}\\
&=&\delta_0+\frac{s}{N}\,(\delta_1-\delta_0)+O(N^{-2})
\end{eqnarray*}

But the law on the right is that of a projection of trace $s/N$, as desired.

\medskip

(2) We have a law formula of the following type:
$${\rm law}(u_{11}+\ldots +u_{ss})=
\sum_{p=0}^sc_p\cdot\frac{(\delta_1-\delta_0)^{*p}}{p!}$$

The coefficients $c_p$ can be estimated by using the Stirling formula, as follows:
\begin{eqnarray*}
c_p
&=&\frac{(tN)!}{N!}\cdot\frac{(N-p)!}{(tN-p)!}\\
&\simeq&\frac{(tN)^{tN}}{N^N}\cdot\frac{(N-p)^{N-p}}{(tN-p)^{tN-p}}\\
&=&\left(\frac{tN}{tN-p}\right)^{tN-p} \left(
\frac{N-p}{N}\right)^{N-p}\left( \frac{tN}{N}\right)^p
\end{eqnarray*}

But the last expression can be estimated by using the definition of the exponentials, and we obtain in this way the following estimate:
$$c_p
\simeq e^pe^{-p}t^p
=t^p$$

We can now compute the Fourier transform with respect to a variable $y$:
\begin{eqnarray*}
{\mathcal F}\left( {\rm law}(u_{11}+\ldots +u_{ss})\right)
&\simeq&\sum_{p=0}^st^p\cdot\frac{(e^y-1)^p}{p!}\\
&=&e^{t(e^y-1)}
\end{eqnarray*}

But this is precisely the Fourier transform of the Poisson law $p_t$, as desired.
\end{proof}

Finally, as yet another variation of our techniques and results, we have:

\begin{theorem}
For the alternating group $A_N\subset O_N$, the truncated character
$$\chi_t(g)=\sum_{i=1}^{[tN]}g_{ii}$$
follows with $N\to\infty$ the Poisson law of parameter $t$, exactly as for $S_N$.
\end{theorem}

\begin{proof}
We have many methods available, and the simplest is to follow the proofs of Theorems 13.16 and 13.17. Observe first that the integrals over $A_N$ are given by:
$$\int_{A_N}g_{i_1j_1}\ldots g_{i_kj_k}=\frac{1}{N!/2}\#\left\{\sigma\in A_N\Big|\sigma(j_1)=i_1,\ldots,\sigma(j_k)=i_k\right\}$$

As before in the symmetric group case, this integral vanishes, unless $\ker i=\ker j$. Now in this latter case, if we denote by $b\in\{1,\ldots,k\}$ the number of blocks of the partition $\ker i=\ker j$, we have $N-b$ points to be sent bijectively to $N-b$ points. But when assuming $N>>0$, and more specifically $N>k$, half of these bijections will be alternating, and so we have  $(N-b)!/2$ solutions. Thus, the integral is:
\begin{eqnarray*}
\int_{A_N}g_{i_1j_1}\ldots g_{i_kj_k}
&=&\frac{1}{N!/2}\#\left\{\sigma\in A_N\Big|\sigma(j_1)=i_1,\ldots,\sigma(j_k)=i_k\right\}\\
&=&\frac{(N-b)!/2}{N!/2}\\
&=&\frac{(N-b)!}{N!}
\end{eqnarray*}

As a conclusion to this study, the asymptotic integrals over $A_N$ are given by the following formula, exactly as in Theorem 13.16, for the symmetric group $S_N$:
$$\int_{A_N}g_{i_1j_1}\ldots g_{i_kj_k}\simeq\begin{cases}
\frac{(N-|\ker i|)!}{N!}&{\rm if}\ \ker i=\ker j\\
0&{\rm otherwise}
\end{cases}$$

But with this integration formula in hand, the computations in the proof of Theorem 13.17 apply, without any changes needed, and give the result.
\end{proof}

\section*{13d. Reflection groups}

Let us discuss now the computation for $H_N$. The result here is as follows:

\begin{theorem}
For the hyperoctahedral group $H_N\subset O_N$ we have
$$\chi_t\sim e^{-t}\sum_{k=-\infty}^\infty\delta_k\sum_{p=0}^\infty \frac{(t/2)^{|k|+2p}}{(|k|+p)!p!}$$ 
with $N\to\infty$, with this measure being called Bessel law of parameter $t>0$.
\end{theorem}

\begin{proof}
Let us regard $H_N$ as being the symmetry group of the graph $\{I^1,\ldots ,I^N\}$ formed by $N$ segments. The diagonal coefficients are then given by:
$$g_{ii}(g)=\begin{cases}
\ 0\ \mbox{ if $g$ moves $I^i$}\\
\ 1\ \mbox{ if $g$ fixes $I^i$}\\
-1\mbox{ if $g$ returns $I^i$}
\end{cases}$$

We denote by $\uparrow g,\downarrow g$ the number of segments among $\{I^1,\ldots ,I^s\}$ which are fixed, respectively returned by an element $g\in H_N$. With this notation, we have:
$$g_{11}+\ldots+g_{ss}=\uparrow g-\downarrow g$$

Let us denote by $P_N$ probabilities computed over the group $H_N$. The density of the law of $g_{11}+\ldots+g_{ss}$ at a point $k\geq 0$ is then given by the following formula:
\begin{eqnarray*}
D(k)
&=&P_N(\uparrow g-\downarrow g=k)\\
&=&\sum_{p=0}^\infty P_N(\uparrow g=k+p, \downarrow g=p)
\end{eqnarray*}

Assume first that we are dealing with $t=1$. We have the following computation:
\begin{eqnarray*}
\lim_{N\to\infty}D(k)
&=&\lim_{N\to\infty}\sum_{p=0}^\infty(1/2)^{k+2p}\begin{pmatrix}k+2p\\ k+p\end{pmatrix} P_N(\uparrow g+\downarrow g=k+2p)\\ 
&=&\sum_{p=0}^\infty(1/2)^{k+2p}\begin{pmatrix}k+2p\\
k+p\end{pmatrix}\frac{1}{e(k+2p)!}\\
&=&\frac{1}{e}\sum_{p=0}^\infty \frac{(1/2)^{k+2p}}{(k+p)!p!}
\end{eqnarray*}

As for the general case $0<t\leq 1$, here the result follows by performing some modifications in the above computation. The asymptotic density is computed as follows:
\begin{eqnarray*}
\lim_{N\to\infty}D(k)
&=&\lim_{N\to\infty}\sum_{p=0}^\infty(1/2)^{k+2p}\begin{pmatrix}k+2p\\ k+p\end{pmatrix} P_N(\uparrow g+\downarrow g=k+2p)\\
&=&\sum_{p=0}^\infty(1/2)^{k+2p}\begin{pmatrix}k+2p\\
k+p\end{pmatrix}\frac{t^{k+2p}}{e^t(k+2p)!}\\
&=&e^{-t}\sum_{p=0}^\infty \frac{(t/2)^{k+2p}}{(k+p)!p!}
\end{eqnarray*}

Together with $D(-k)=D(k)$, this gives the formula in the statement.
\end{proof}

In the above result the terminology comes from the fact, up to $t\to t/2$, the density of the law is the following function, called Bessel function of the first kind:
$$f_k(t)=\sum_{p=0}^\infty \frac{t^{|k|+2p}}{(|k|+p)!p!}$$

Let us further study now these Bessel laws. A key result regarding the Poisson laws was the semigroup formula $p_s*p_t=p_{s+t}$, and in analogy with this, we have:

\index{convolution}

\begin{theorem}
The Bessel laws $b_t$ have the property
$$b_s*b_t=b_{s+t}$$
so they form a truncated one-parameter semigroup
with respect to convolution.
\end{theorem}

\begin{proof}
With notations as above, the Fourier transform of $b_t$ is given by:
$$Fb_t(y)=e^{-t}\sum_{k=-\infty}^\infty e^{ky}\,f_k(t/2)$$

We can compute now the derivative with respect to $t$, as follows:
$$Fb_t(y)'=-Fb_t(y)+\frac{e^{-t}}{2}\sum_{k=-\infty}^\infty e^{ky}\,f_k'(t/2)$$

On the other hand, the derivative of $f_k$ with $k\geq 1$ is given by:
\begin{eqnarray*}
f_k'(t)
&=&\sum_{p=0}^\infty \frac{(k+2p)t^{k+2p-1}}{(k+p)!p!}\\
&=&\sum_{p=0}^\infty \frac{(k+p)t^{k+2p-1}}{(k+p)!p!}+\sum_{p=0}^\infty\frac{p\,t^{k+2p-1}}{(k+p)!p!}\\
&=&\sum_{p=0}^\infty \frac{t^{k+2p-1}}{(k+p-1)!p!}+\sum_{p=1}^\infty\frac{t^{k+2p-1}}{(k+p)!(p-1)!}\\
&=&\sum_{p=0}^\infty \frac{t^{(k-1)+2p}}{((k-1)+p)!p!}+\sum_{p=1}^\infty\frac{t^{(k+1)+2(p-1)}}{((k+1)+(p-1))!(p-1)!}\\
&=&f_{k-1}(t)+f_{k+1}(t)
\end{eqnarray*}

This computation works in fact for any $k$, so we get:
\begin{eqnarray*}
Fb_t(y)'
&=&-Fb_t(y)+\frac{e^{-t}}{2}
\sum_{k=-\infty}^\infty e^{ky} (f_{k-1}(t/2)+f_{k+1}(t/2))\\
&=&-Fb_t(y)+\frac{e^{-t}}{2} \sum_{k=-\infty}^\infty
e^{(k+1)y}f_{k}(t/2)+e^{(k-1)y}f_{k}(t/2)\\
&=&-Fb_t(y)+\frac{e^{y}+e^{-y}}{2}\,Fb_t(y)\\
&=&\left(\frac{e^{y}+e^{-y}}{2}-1\right)Fb_t(y)
\end{eqnarray*}

Thus the log of the Fourier transform is linear in $t$, and we get
the assertion.
\end{proof}

In order to further interpret the above results, and discuss as well, more generally, the complex reflection groups $H_N^s=\mathbb Z_s\wr S_N$, let us begin with some probability. We have the following notion, inspired from the Poisson limit theory developed before:

\index{compound Poisson law}

\begin{definition}
Associated to any compactly supported positive measure $\nu$ on $\mathbb C$ is the probability measure
$$p_\nu=\lim_{n\to\infty}\left(\left(1-\frac{c}{n}\right)\delta_0+\frac{1}{n}\nu\right)^{*n}$$
where $c=mass(\nu)$, called compound Poisson law.
\end{definition}

In what follows we will be mainly interested in the case where $\nu$ is discrete, as is for instance the measure $\nu=t\delta_1$ with $t>0$, which produces via the above procedure the Poisson law $p_t$. To be more precise, we will be mainly interested in the case where $\nu$ is a multiple of the uniform measure on the $s$-th roots of unity. More on this later.

\bigskip

The following result allows us to detect the compound Poisson laws:

\index{Fourier transform}

\begin{proposition}
For $\nu=\sum_{i=1}^sc_i\delta_{z_i}$ with $c_i>0$ and $z_i\in\mathbb C$ we have
$$F_{p_\nu}(y)=\exp\left(\sum_{i=1}^sc_i(e^{iyz_i}-1)\right)$$
where $F$ denotes as usual the Fourier transform.
\end{proposition}

\begin{proof}
Let $\mu_n$ be the measure appearing in Definition 13.23, namely:
$$\mu_n=\left(1-\frac{c}{n}\right)\delta_0+\frac{1}{n}\nu$$

We have the following computation, in the context of Definition 13.23:
\begin{eqnarray*}
F_{\mu_n}(y)=\left(1-\frac{c}{n}\right)+\frac{1}{n}\sum_{i=1}^sc_ie^{iyz_i}
&\implies&F_{\mu_n^{*n}}(y)=\left(\left(1-\frac{c}{n}\right)+\frac{1}{n}\sum_{i=1}^sc_ie^{iyz_i}\right)^n\\
&\implies&F_{p_\nu}(y)=\exp\left(\sum_{i=1}^sc_i(e^{iyz_i}-1)\right)
\end{eqnarray*}

Thus, we have obtained the formula in the statement.
\end{proof}

We have the following result, providing an alternative to Definition 13.23, and which will be our formulation here of the Compound Poisson Limit Theorem (CPLT):

\index{compound Poisson Limit theorem}
\index{CPLT}

\begin{theorem}
For $\nu=\sum_{i=1}^sc_i\delta_{z_i}$ with $c_i>0$ and $z_i\in\mathbb C$ we have
$$p_\nu={\rm law}\left(\sum_{i=1}^sz_i\alpha_i\right)$$
where the variables $\alpha_i$ are Poisson $(c_i)$, independent.
\end{theorem}

\begin{proof}
Let $\alpha$ be the sum of Poisson variables in the statement. We have:
\begin{eqnarray*}
F_{\alpha_i}(y)=\exp(c_i(e^{iy}-1)
&\implies&F_{z_i\alpha_i}(y)=\exp(c_i(e^{iyz_i}-1))\\
&\implies&F_\alpha(y)=\exp\left(\sum_{i=1}^sc_i(e^{iyz_i}-1)\right)
\end{eqnarray*}

Thus we have the same formula as in Proposition 13.24, as desired.
\end{proof}

Getting back now to the Bessel laws, we have the following result:

\index{Bessel law}

\begin{theorem}
The Bessel laws $b_t$ are compound Poisson laws, given by
$$b_t=p_{t\varepsilon}$$
where $\varepsilon=\frac{1}{2}(\delta_{-1}+\delta_1)$ is the uniform measure on $\mathbb Z_2$.
\end{theorem}

\begin{proof}
This follows indeed by comparing the formula of the Fourier transform of $b_t$, from the proof of Theorem 13.22, with the formula in Proposition 13.24.
\end{proof}

Our next task will be that of unifying and generalizing the results that we have for $S_N,H_N$. For this purpose, following \cite{bbc}, let us introduce the following notion:

\begin{definition}
The Bessel law of level $s\in\mathbb N\cup\{\infty\}$ and parameter $t>0$ is
$$b_t^s=p_{t\varepsilon_s}$$
with $\varepsilon_s$ being the uniform measure on the $s$-th roots of unity.
\end{definition}

Observe that at $s=1,2$ we obtain the Poisson and real Bessel laws:
$$b^1_t=p_t\quad,\quad 
b^2_t=b_t$$

Another important particular case is $s=\infty$, where we obtain a measure which is actually not discrete, that we will denote here as follows:
$$b^\infty_t=B_t$$

As a basic result on these laws, generalizing those before about $p_t,b_t$, we have:

\begin{theorem}
The generalized Bessel laws $b^s_t$ have the property
$$b^s_t*b^s_{t'}=b^s_{t+t'}$$
so they form a truncated one-parameter semigroup with respect to convolution.
\end{theorem}

\begin{proof}
This follows indeed from the Fourier transform formula from Proposition 13.24, because for the Bessel laws, the log of this Fourier transform is linear in $t$.
\end{proof}

We can go back now to the reflection groups, and we have the following result:

\begin{theorem}
For the group $H_N^s=\mathbb Z_s\wr S_N$ we have, with $N\to\infty$,
$$\chi_t\sim b^s_t$$
where $b_t^s=p_{t\varepsilon_s}$ is the Bessel law of level $s\in\mathbb N\cup\{\infty\}$ and parameter $t>0$.
\end{theorem}

\begin{proof}
In the case $t=1$, by arguing as before at $s=2$, since the limit probability for a random permutation to have exactly $k$ fixed points is $e^{-1}/k!$, we obtain:
$$\lim_{N\to\infty}law(\chi_1)=e^{-1}\sum_{k=0}^\infty \frac{1}{k!}\,\varepsilon_s^{*k}$$

On the other hand, we get from the definition of the Bessel law $b^s_1$, as desired:
\begin{eqnarray*}
b^s_1
&=&\lim_{N\to\infty}\left(\left(1-\frac{1}{N}\right)\delta_0+\frac{1}{N}\,\varepsilon_s\right)^{*N}\\
&=&\lim_{N\to\infty}\sum_{k=0}^N\begin{pmatrix}N\\ k\end{pmatrix}\left(1-\frac{1}{N}\right)^{N-k}\frac{1}{N^k}\,\varepsilon_s^{*k}\\
&=&e^{-1}\sum_{k=0}^\infty\frac{1}{k!}\,\varepsilon_s^{*k}
\end{eqnarray*}

When $t>0$ is arbitrary, we can use the same method, with some modifications where needed, again by arguing as before at $s=2$, and we obtain the result.
\end{proof}

Many other interesting things can be said about the Bessel laws found above, and we will be back to them, in the next chapter. Finally, in what regards the general reflection groups $H_N^{sd}$ from chapter 5, the arguments in the proof of Theorem 13.20 apply to them, and the asymptotic laws follow to be the same as those computed above, for $H_N^s$.

\section*{13e. Exercises}

This was a pleasant and quite elementary chapter, and as exercises, we have:

\begin{exercise}
Learn more about the PLT and its meaning, from probabilists.
\end{exercise}

\begin{exercise}
Experiment a bit with the Bell numbers, and their evaluation.
\end{exercise}

\begin{exercise}
Learn more about the Bessel functions of the first kind.
\end{exercise}

\begin{exercise}
Learn more about compound Poisson laws, and their properties.
\end{exercise}

\begin{exercise}
Learn more about the Bessel laws, and their various properties.
\end{exercise}

\begin{exercise}
Fill in all the details, for the Bessel law computation for $H_N^s$.
\end{exercise}

\begin{exercise}
Clarify what we said regarding $A_N$, then investigate $H_N^{sd}$.
\end{exercise}

\begin{exercise}
What are the densities and moments of the Bessel laws?
\end{exercise}

As bonus exercise, and no surprise here, start reading some probability theory.

\chapter{Normal laws}

\section*{14a. Moment method}

We have learned in chapter 13 many interesting things regarding the main character of $S_N$, and more generally of $H_N^s=\mathbb Z_s\wr S_N$, but the story is far from being over, with this. We have indeed the following two pressing questions, to be investigated:

\bigskip

(1) What can we say, again regarding the main character of $S_N$ and $H_N^s$, in relation with representation theory? More specifically, how does the easiness property of $S_N$ and $H_N^s$ come into play? Can we go further, with finer computations, by using it?

\medskip

(2) What about continuous groups, such as $O_N$, $U_N$? It sounds reasonable to conjecture here that we obtain real and complex normal laws, but is this correct? Then, what about the free versions $O_N^+,U_N^+$? Also, what about the free versions $S_N^+,H_N^{s+}$?

\bigskip

We will discuss here these questions, with all this being technically not that complicated, in view of our current representation theory knowledge, acquired the hard way in Part III. The only thing that needs some thinking is a good plan, in order to properly deal with the above jungle of questions, and in answer, our plan will be as follows:

\bigskip

-- We will first talk a bit more about general probability theory, with this being something quite natural, as apprentice algebraists, analysis is our weak point, right.

\medskip

-- Then we will get into representation theory and easiness matters, by starting with the simplest classical objects which are around, namely the groups $O_N,U_N$.

\medskip

-- And then we will navigate through the remaining problems, with a discussion regarding $S_N,H_N^s$ at the end of this chapter, and with the free versions, in chapter 15.

\bigskip

Getting started now, we need a crash course in general probability. With the idea in mind of having things centered around the notion of moments, let us first formulate:

\index{probability space}
\index{random variable}
\index{moments}
\index{law}
\index{distribution}

\begin{definition}
Let $X$ be a probability space, that is, a space with a probability measure, and with the corresponding integration denoted $E$, and called expectation.
\begin{enumerate}
\item The random variables are the real functions $f\in L^\infty(X)$.

\item The moments of such a variable are the numbers $M_k(f)=E(f^k)$.

\item The law of such a variable is the measure given by $M_k(f)=\int_\mathbb Rx^kd\mu_f(x)$.
\end{enumerate}
\end{definition}

As a comment on this, the fact that $\mu_f$ as above exists indeed is not exactly trivial. But we can do this by looking at formulae of the following type:
$$E(\varphi(f))=\int_\mathbb R\varphi(x)d\mu_f(x)$$

Indeed, having this for monomials $\varphi(x)=x^n$, as above, is the same as having it for polynomials $\varphi\in\mathbb R[X]$, which in turn is the same as having it for the characteristic functions $\varphi=\chi_I$ of measurable sets $I\subset\mathbb R$. Thus, in the end, what we need is:
$$P(f\in I)=\mu_f(I)$$

But this formula can serve as a definition for $\mu_f$, and we are done. 

\bigskip

Regarding now independence, which is the key notion in probability, we have here:

\index{independence}

\begin{definition}
Two variables $f,g\in L^\infty(X)$ are called independent when
$$E(f^kg^l)=E(f^k)\,E(g^l)$$
happens, for any $k,l\in\mathbb N$.
\end{definition}

Again, this definition hides some non-trivial things, the idea being a bit as before, namely that of looking at formulae of the following type:
$$E[\varphi(f)\psi(g)]=E[\varphi(f)]\,E[\psi(g)]$$

To be more precise, passing as before from monomials to polynomials, then to characteristic functions, we are led to the usual definition of independence, namely:
$$P(f\in I,g\in J)=P(f\in I)\,P(g\in J)$$

As a first result now, in the study of independence, we have:

\index{convolution}

\begin{theorem}
Assuming that $f,g\in L^\infty(X)$ are independent, we have
$$\mu_{f+g}=\mu_f*\mu_g$$
where $*$ is the convolution of real probability measures.
\end{theorem}

\begin{proof}
We have the following computation, using the independence of $f,g$:
$$\int_\mathbb Rx^kd\mu_{f+g}(x)
=E((f+g)^k)
=\sum_r\binom{k}{r}M_r(f)M_{k-r}(g)$$

On the other hand, we have as well the following computation:
\begin{eqnarray*}
\int_\mathbb Rx^kd(\mu_f*\mu_g)(x)
&=&\int_{\mathbb R\times\mathbb R}(x+y)^kd\mu_f(x)d\mu_g(y)\\
&=&\sum_r\binom{k}{r}M_r(f)M_{k-r}(g)
\end{eqnarray*}

Thus $\mu_{f+g}$ and $\mu_f*\mu_g$ have the same moments, so they coincide, as claimed.
\end{proof}

As a second result on independence, which is more advanced, we have:

\index{independence}
\index{Fourier transform}

\begin{theorem}
Assuming that $f,g\in L^\infty(X)$ are independent, we have
$$F_{f+g}=F_fF_g$$
where $F_f(x)=E(e^{ixf})$ is the Fourier transform.
\end{theorem}

\begin{proof}
This is something which comes from Theorem 14.3, as follows:
\begin{eqnarray*}
F_{f+g}(x)
&=&\int_\mathbb Re^{ixz}d(\mu_f*\mu_g)(z)\\
&=&\int_{\mathbb R\times\mathbb R}e^{ix(z+t)}d\mu_f(z)d\mu_g(t)\\
&=&\int_\mathbb Re^{ixz}d\mu_f(z)\int_\mathbb Re^{ixt}d\mu_g(t)\\
&=&F_f(x)F_g(x)
\end{eqnarray*}

Thus, we are led to the conclusion in the statement.
\end{proof}

So long for laws, moments and independence, at the theoretical level. In what follows Definitions 14.1 and 14.2 and Theorems 14.3 and 14.4 are what we need, I mean we will certainly get away with that, as general theory, for our algebraic purposes here. 

\bigskip

Now remembering what we wanted to do, namely using representation theory and moments, we are led to the question of reconstructing a measure out of its moments. And the answer here is something quite tricky, involving complex analysis, as follows:

\begin{theorem}
The density of a real probability measure $\mu$ can be recaptured from the sequence of moments $\{M_k\}_{k\geq0}$ via the Stieltjes inversion formula
$$d\mu (x)=\lim_{t\searrow 0}-\frac{1}{\pi}\,Im\left(G(x+it)\right)\cdot dx$$
where the function on the right, given in terms of moments by
$$G(\xi)=\xi^{-1}+M_1\xi^{-2}+M_2\xi^{-3}+\ldots$$
is the Cauchy transform of the measure $\mu$.
\end{theorem}

\begin{proof}
The Cauchy transform of our measure $\mu$ is given by:
\begin{eqnarray*}
G(\xi)
&=&\xi^{-1}\sum_{k=0}^\infty M_k\xi^{-k}\\\
&=&\int_\mathbb R\frac{\xi^{-1}}{1-\xi^{-1}y}\,d\mu(y)\\
&=&\int_\mathbb R\frac{1}{\xi-y}\,d\mu(y)
\end{eqnarray*}

Now with $\xi=x+it$, we obtain the following formula:
\begin{eqnarray*}
Im(G(x+it))
&=&\int_\mathbb RIm\left(\frac{1}{x-y+it}\right)d\mu(y)\\
&=&\int_\mathbb R\frac{1}{2i}\left(\frac{1}{x-y+it}-\frac{1}{x-y-it}\right)d\mu(y)\\
&=&-\int_\mathbb R\frac{t}{(x-y)^2+t^2}\,d\mu(y)
\end{eqnarray*}

By integrating over $[a,b]$ we obtain, with the change of variables $x=y+tz$:
\begin{eqnarray*}
\int_a^bIm(G(x+it))dx
&=&-\int_\mathbb R\int_a^b\frac{t}{(x-y)^2+t^2}\,dx\,d\mu(y)\\
&=&-\int_\mathbb R\int_{(a-y)/t}^{(b-y)/t}\frac{t}{(tz)^2+t^2}\,t\,dz\,d\mu(y)\\
&=&-\int_\mathbb R\int_{(a-y)/t}^{(b-y)/t}\frac{1}{1+z^2}\,dz\,d\mu(y)\\
&=&-\int_\mathbb R\left(\arctan\frac{b-y}{t}-\arctan\frac{a-y}{t}\right)d\mu(y)
\end{eqnarray*}

Now observe that with $t\searrow0$ we have:
$$\lim_{t\searrow0}\left(\arctan\frac{b-y}{t}-\arctan\frac{a-y}{t}\right)
=\begin{cases}
\frac{\pi}{2}-\frac{\pi}{2}=0& (y<a)\\
\frac{\pi}{2}-0=\frac{\pi}{2}& (y=a)\\
\frac{\pi}{2}-(-\frac{\pi}{2})=\pi& (a<y<b)\\
0-(-\frac{\pi}{2})=\frac{\pi}{2}& (y=b)\\
-\frac{\pi}{2}-(-\frac{\pi}{2})=0& (y>b)
\end{cases}$$

We therefore obtain the following formula:
$$\lim_{t\searrow0}\int_a^bIm(G(x+it))dx=-\pi\left(\mu(a,b)+\frac{\mu(a)+\mu(b)}{2}\right)$$

Thus, we are led to the conclusion in the statement.
\end{proof}

Before getting further, let us mention that the above result does not fully solve the moment problem, because we still have the question of understanding when a sequence of numbers $M_1,M_2,M_3,\ldots$ can be the moments of a measure $\mu$.  We have here:

\begin{theorem}
A sequence of numbers $M_0,M_1,M_2,M_3,\ldots\in\mathbb R$, with $M_0=1$, is the series of moments of a real probability measure $\mu$ precisely when:
$$\begin{vmatrix}M_0\end{vmatrix}\geq0\quad,\quad 
\begin{vmatrix}
M_0&M_1\\
M_1&M_2
\end{vmatrix}\geq0\quad,\quad 
\begin{vmatrix}
M_0&M_1&M_2\\
M_1&M_2&M_3\\
M_2&M_3&M_4\\
\end{vmatrix}\geq0\quad,\quad 
\ldots$$
That is, the associated Hankel determinants must be all positive.
\end{theorem}

\begin{proof}
This is something a bit more advanced, the idea being as follows:

\medskip

(1) As a first observation, the positivity conditions in the statement tell us that the following associated linear forms must be positive:
$$\sum_{i,j=1}^nc_i\bar{c}_jM_{i+j}\geq0$$

(2) But this is something very classical, in one sense the result being elementary, coming from the following computation, which shows that we have positivity indeed:
$$\int_\mathbb R\left|\sum_{i=1}^nc_ix^i\right|^2d\mu(x)
=\int_\mathbb R\sum_{i,j=1}^nc_i\bar{c}_jx^{i+j}d\mu(x)
=\sum_{i,j=1}^nc_i\bar{c}_jM_{i+j}$$

(3) As for the other sense, here the result comes once again from the above formula, this time via some standard functional analysis.
\end{proof}

As a basic application of the Stieltjes formula, let us solve the moment problem for the simplest numbers on the planet, and perhaps universe, combinatorially speaking, which are the Catalan numbers $C_k$, and the central binomial coefficients $D_k$. We first have:

\begin{theorem}
The real measure having as even moments the Catalan numbers, $C_k=\frac{1}{k+1}\binom{2k}{k}$, and having all odd moments $0$ is the measure
$$\gamma_1=\frac{1}{2\pi}\sqrt{4-x^2}dx$$
called Wigner semicircle law on $[-2,2]$.
\end{theorem}

\begin{proof}
In order to apply the inversion formula, our starting point will be the standard calculus formula for the generating series of the Catalan numbers, namely:
$$\sum_{k=0}^\infty C_kz^k=\frac{1-\sqrt{1-4z}}{2z}$$

Observe that this formula holds indeed, by using the Taylor series expansion of $\sqrt{1-4z}$. Now by using this with $z=\xi^{-2}$, we obtain the following formula:
\begin{eqnarray*}
G(\xi)
&=&\xi^{-1}\sum_{k=0}^\infty C_k\xi^{-2k}\\
&=&\xi^{-1}\cdot\frac{1-\sqrt{1-4\xi^{-2}}}{2\xi^{-2}}\\
&=&\frac{\xi}{2}\left(1-\sqrt{1-4\xi^{-2}}\right)\\
&=&\frac{\xi}{2}-\frac{1}{2}\sqrt{\xi^2-4}
\end{eqnarray*}

Now let us apply Theorem 14.5. The study here goes as follows:

\medskip

(1) According to the general philosophy of the Stieltjes formula, the first term, namely $\xi/2$, which is ``trivial'', will not contribute to the density. 

\medskip

(2) As for the second term, which is something non-trivial, this will contribute to the density, the rule here being that the square root $\sqrt{\xi^2-4}$ will be replaced by the ``dual'' square root $\sqrt{4-x^2}\,dx$, and that we have to multiply everything by $-1/\pi$. 

\medskip

(3) As a conclusion, by Stieltjes inversion we obtain the following density:
$$d\mu(x)
=-\frac{1}{\pi}\cdot-\frac{1}{2}\sqrt{4-x^2}\,dx
=\frac{1}{2\pi}\sqrt{4-x^2}dx$$

Thus, we have obtained the mesure in the statement, and we are done.
\end{proof}

We have the following version of the above result:

\begin{theorem}
The real measure having as sequence of moments the Catalan numbers, $C_k=\frac{1}{k+1}\binom{2k}{k}$, is the measure
$$\pi_1=\frac{1}{2\pi}\sqrt{4x^{-1}-1}\,dx$$
called Marchenko-Pastur law on $[0,4]$.
\end{theorem}

\begin{proof}
As before, we use the standard formula for the generating series of the Catalan numbers. With $z=\xi^{-1}$ in that formula, we obtain the following formula:
\begin{eqnarray*}
G(\xi)
&=&\xi^{-1}\sum_{k=0}^\infty C_k\xi^{-k}\\
&=&\xi^{-1}\cdot\frac{1-\sqrt{1-4\xi^{-1}}}{2\xi^{-1}}\\
&=&\frac{1}{2}\left(1-\sqrt{1-4\xi^{-1}}\right)\\
&=&\frac{1}{2}-\frac{1}{2}\sqrt{1-4\xi^{-1}}
\end{eqnarray*}

With this in hand, let us apply now the Stieltjes inversion formula, from Theorem 14.5. We obtain, a bit as before in Theorem 14.7, the following density:
$$d\mu(x)
=-\frac{1}{\pi}\cdot-\frac{1}{2}\sqrt{4x^{-1}-1}\,dx
=\frac{1}{2\pi}\sqrt{4x^{-1}-1}\,dx$$

Thus, we are led to the conclusion in the statement.
\end{proof}

Regarding now the central binomial coefficients, we have here:

\begin{theorem}
The real probability measure having as moments the central binomial coefficients, $D_k=\binom{2k}{k}$, is the measure
$$\alpha_1=\frac{1}{\pi\sqrt{x(4-x)}}\,dx$$
called arcsine law on $[0,4]$.
\end{theorem}

\begin{proof}
We have the following computation, using some standard calculus formulae, and more specifically, the Taylor series expansion of $1/\sqrt{1-4/\xi}$:
\begin{eqnarray*}
G(\xi)
&=&\xi^{-1}\sum_{k=0}^\infty D_k\xi^{-k}\\
&=&\frac{1}{\xi}\sum_{k=0}^\infty D_k\left(-\frac{t}{4}\right)^k\\
&=&\frac{1}{\xi}\cdot\frac{1}{\sqrt{1-4/\xi}}\\
&=&\frac{1}{\sqrt{\xi(\xi-4)}} 
\end{eqnarray*}

But this gives the density in the statement, via Theorem 14.5. 
\end{proof}

Finally, we have the following version of the above result:

\begin{theorem}
The real probability measure having as moments the middle binomial coefficients, $E_k=\binom{k}{[k/2]}$, is the following law on $[-2,2]$,
$$\sigma_1=\frac{1}{2\pi}\sqrt{\frac{2+x}{2-x}}\,dx$$
called modified the arcsine law on $[-2,2]$.
\end{theorem}

\begin{proof}
In terms of the central binomial coefficients $D_k$, we have:
$$E_{2k}=\binom{2k}{k}=\frac{(2k)!}{k!k!}=D_k$$
$$E_{2k-1}=\binom{2k-1}{k}=\frac{(2k-1)!}{k!(k-1)!}=\frac{D_k}{2}$$

Standard calculus based on the Taylor formula for $(1+t)^{-1/2}$ gives:
$$\frac{1}{2x}\left(\sqrt{\frac{1+2x}{1-2x}}-1\right)=\sum_{k=0}^\infty E_kx^k$$

With $x=\xi^{-1}$ we obtain the following formula for the Cauchy transform:
\begin{eqnarray*}
G(\xi)
&=&\xi^{-1}\sum_{k=0}^\infty E_k\xi^{-k}\\
&=&\frac{1}{\xi}\left(\sqrt{\frac{1+2/\xi}{1-2/\xi}}-1\right)\\
&=&\frac{1}{\xi}\left(\sqrt{\frac{\xi+2}{\xi-2}}-1\right)
\end{eqnarray*}

By Stieltjes inversion we obtain the density in the statement.
\end{proof}

As a comment now, while the numbers $C_k,D_k$ were advertised as being the simplest on our planet, the corresponding measures rather seem to come from the outer space. However, this is not exactly true, and quite the opposite happens, because we have:

\begin{theorem}
The following happen:
\begin{enumerate}
\item The main character of $SU_2$ follows the Wigner law $\gamma_1$.

\item The main character of $SO_3$ follows the Marchenko-Pastur law $\pi_1$.
\end{enumerate}
\end{theorem}

\begin{proof}
This follows indeed from a close inspection of our computations for $SU_2,SO_3$ from chapter 10, which show that the main character moments are precisely those in Theorems 14.7 and 14.8. Thus, we are led to the conclusions in the statement.
\end{proof}

Quite nice all this, obviously the Stieltjes inversion formula is something very useful, and $\gamma_1,\pi_1$ are the ``simplest'' continuous measures. As for the arcsine laws $\alpha_1,\sigma_1$, these can be shown to appear as well in various representation theory contexts.

\section*{14b. Lindst\"om formula}

Getting now to what we wanted to do, as a main job for this chapter, namely the character laws for $O_N,U_N$, what in did in chapter 13 for $S_N,H_N^{s+}$, based on inclusion-exclusion, or direct integration, will certainly not apply to our situation. Also, what we just did for $SU_2,SO_3$, ultimately coming from 3D geometry, will not apply either. 

\bigskip

Nevermind. As explained in chapter 13, an alternative approach to the computation of character laws is via Peter-Weyl theory. Now since we know that $O_N,U_N$ are easy, problem basically solved, we just have to convert that partition combinatorics into probability.

\bigskip

Let us start our discussion with a general result, about easiness:

\begin{theorem}
Given an easy group $G\subset U_N$, coming from a category of partitions $D\subset P$, the moments of the main character $\chi=Tr(u)$ are the numbers
$$M_k=\dim\left(span\left(\xi_\pi\Big|\pi\in D(k)\right)\right)$$
with the conventions $D(k)=D(\emptyset,k)$, and $\xi_\pi=T_\pi$ for $\pi\in D(k)$. Thus we have
$$M_k\leq|D(k)|$$
with equality when the vectors $\{\xi_\pi|\pi\in D(k)\}$ are linearly independent.
\end{theorem}

\begin{proof}
We know from Peter-Weyl theory that integrating a character gives the number of fixed points, and in our situation, $G\subset_uU_N$, this gives:
\begin{eqnarray*}
\int_G\chi^k
&=&\int_G(\chi_u)^k\\
&=&\int_G\chi_{u^{\otimes k}}\\
&=&\dim(Fix(u^{\otimes k}))
\end{eqnarray*}

On the other hand, the easiness condition for $G\subset_uU_N$ tells us that we have:
$$Hom(u^{\otimes k},u^{\otimes l})=span\left(T_\pi\Big|\pi\in D(k,l)\right)$$

Now with $(k,l)\to(\emptyset,k)$, and with the conventions $D(k)=D(\emptyset,k)$, and $\xi_\pi=T_\pi$ for $\pi\in D(k)$ in the statement, which are both standard, we obtain:
$$Fix(u^{\otimes k})=span\left(\xi_\pi\Big|\pi\in D(k)\right)$$

Thus, we are led to the moment formula. As for the last assertion, that is clear.
\end{proof}

In view of the above result, we are led to the following question:

\begin{question}
Given a category of partitions $D\subset P$, when are the vectors
$$\xi_\pi=\sum_{\ker i\leq\pi}e_{i_1}\otimes\ldots\otimes e_{i_k}\in\mathbb R^N$$
with $\pi\in D(k)$ linearly independent?
\end{question}

And good question this is, because the answer is no in general, as you can easily check by doing some computations, at small values of $N,k$, and with $D=P$. In fact, there is not even need for computations, because $D=P$ is the category producing $S_N$, and a positive answer to Question 14.13 would lead, via Theorem 14.12, to the conclusion that the law of the main character of $S_N$ does not depend on $N\in\mathbb N$, which is wrong.

\bigskip

For other categories, such as $P_2,\mathcal P_2$, pretty much the same can be said, with counterexamples being easy to construct, at small values of $N,k$, and with this coming too from the fact that the main character laws of the corresponding easy groups, which are $O_N,U_N$, cannot be independent of $N\in\mathbb N$, for various geometric reasons.

\bigskip

This being said, and getting back now to $D=P$, a ray of light comes from our computations in chapter 13 for the main character of $S_N$, which lead us to:

\begin{answer}
For $D=P$, corresponding to the group $S_N$, the vectors
$$\left\{\xi_\pi\Big|\pi\in D(k)\right\}$$ 
are linearly independent with $N>>0$, and more precisely, when $N\geq k$.
\end{answer}

In practice now, although this can be deduced from our computations from chapter 13, all this is quite fundamental for our considerations to follow, and it is better to have a simple, independent proof of it. Let us begin with some standard combinatorics:

\index{lattice of partitions}
\index{order of partitions}
\index{supremum of partitions}

\begin{definition}
Let $P(k)$ be the set of partitions of $\{1,\ldots,k\}$, and $\pi,\sigma\in P(k)$.
\begin{enumerate}
\item We write $\pi\leq\sigma$ if each block of $\pi$ is contained in a block of $\sigma$.

\item We let $\pi\vee\sigma\in P(k)$ be the partition obtained by superposing $\pi,\sigma$.
\end{enumerate}
Also, we denote by $|.|$ the number of blocks of the partitions $\pi\in P(k)$.
\end{definition}

As an illustration here, at $k=2$ we have $P(2)=\{||,\sqcap\}$, and we have:
$$||\leq\sqcap$$

Also, at $k=3$ we have $P(3)=\{|||,\sqcap|,\sqcap\hskip-3.2mm{\ }_|\,,|\sqcap,\sqcap\hskip-0.7mm\sqcap\}$, and the order relation is as follows:
$$|||\ \leq\ \sqcap|\ ,\ \sqcap\hskip-3.2mm{\ }_|\ ,\ |\sqcap\ \leq\ \sqcap\hskip-0.7mm\sqcap$$

In relation with our linear independence questions, the idea will be that of using:

\begin{proposition}
The Gram matrix of the vectors $\xi_\pi$ is given by the formula
$$<\xi_\pi,\xi_\sigma>=N^{|\pi\vee\sigma|}$$
where $\vee$ is the superposition operation, and $|.|$ is the number of blocks.
\end{proposition}

\begin{proof}
According to the formula of the vectors $\xi_\pi$, we have:
\begin{eqnarray*}
<\xi_\pi,\xi_\sigma>
&=&\sum_{i_1\ldots i_k}\delta_\pi(i_1,\ldots,i_k)\delta_\sigma(i_1,\ldots,i_k)\\
&=&\sum_{i_1\ldots i_k}\delta_{\pi\vee\sigma}(i_1,\ldots,i_k)\\
&=&N^{|\pi\vee\sigma|}
\end{eqnarray*}

Thus, we have obtained the formula in the statement.
\end{proof}

In order to study the Gram matrix $G_k(\pi,\sigma)=N^{|\pi\vee\sigma|}$, and more specifically to compute its determinant, we will use several standard facts about the partitions. We have:

\index{Gram matrix}
\index{M\"obius function}

\begin{definition}
The M\"obius function of any lattice, and so of $P$, is given by
$$\mu(\pi,\sigma)=\begin{cases}
1&{\rm if}\ \pi=\sigma\\
-\sum_{\pi\leq\tau<\sigma}\mu(\pi,\tau)&{\rm if}\ \pi<\sigma\\
0&{\rm if}\ \pi\not\leq\sigma
\end{cases}$$
with the construction being performed by recurrence.
\end{definition}

As an illustration here, for $P(2)=\{||,\sqcap\}$, we have by definition:
$$\mu(||,||)=\mu(\sqcap,\sqcap)=1$$

Also, $||<\sqcap$, with no intermediate partition in between, so we obtain:
$$\mu(||,\sqcap)=-\mu(||,||)=-1$$

Finally, we have $\sqcap\not\leq||$, and so we have as well the following formula:
$$\mu(\sqcap,||)=0$$

Thus, as a conclusion, we have computed the M\"obius matrix $M_2(\pi,\sigma)=\mu(\pi,\sigma)$ of the lattice $P(2)=\{||,\sqcap\}$, the formula being as follows:
$$M_2=\begin{pmatrix}1&-1\\ 0&1\end{pmatrix}$$

Back to the general case now, the main interest in the M\"obius function comes from the M\"obius inversion formula, which states that the following happens:
$$f(\sigma)=\sum_{\pi\leq\sigma}g(\pi)\quad
\implies\quad g(\sigma)=\sum_{\pi\leq\sigma}\mu(\pi,\sigma)f(\pi)$$

In linear algebra terms, the statement and proof of this formula are as follows:

\index{M\"obius inversion}

\begin{theorem}
The inverse of the adjacency matrix of $P(k)$, given by
$$A_k(\pi,\sigma)=\begin{cases}
1&{\rm if}\ \pi\leq\sigma\\
0&{\rm if}\ \pi\not\leq\sigma
\end{cases}$$
is the M\"obius matrix of $P$, given by $M_k(\pi,\sigma)=\mu(\pi,\sigma)$.
\end{theorem}

\begin{proof}
This is well-known, coming for instance from the fact that $A_k$ is upper triangular. Indeed, when inverting, we are led into the recurrence from Definition 14.17.
\end{proof}

As an illustration, for $P(2)$ the formula $M_2=A_2^{-1}$ appears as follows:
$$\begin{pmatrix}1&-1\\ 0&1\end{pmatrix}=
\begin{pmatrix}1&1\\ 0&1\end{pmatrix}^{-1}$$

Now back to our Gram matrix considerations, we have the following key result:

\begin{proposition}
The Gram matrix of the vectors $\xi_\pi$ with $\pi\in P(k)$,
$$G_{\pi\sigma}=N^{|\pi\vee\sigma|}$$
decomposes as a product of upper/lower triangular matrices, $G_k=A_kL_k$, where
$$L_k(\pi,\sigma)=
\begin{cases}
N(N-1)\ldots(N-|\pi|+1)&{\rm if}\ \sigma\leq\pi\\
0&{\rm otherwise}
\end{cases}$$
and where $A_k$ is the adjacency matrix of $P(k)$.
\end{proposition}

\begin{proof}
We have the following computation, based on Proposition 14.16:
\begin{eqnarray*}
G_k(\pi,\sigma)
&=&N^{|\pi\vee\sigma|}\\
&=&\#\left\{i_1,\ldots,i_k\in\{1,\ldots,N\}\Big|\ker i\geq\pi\vee\sigma\right\}\\
&=&\sum_{\tau\geq\pi\vee\sigma}\#\left\{i_1,\ldots,i_k\in\{1,\ldots,N\}\Big|\ker i=\tau\right\}\\
&=&\sum_{\tau\geq\pi\vee\sigma}N(N-1)\ldots(N-|\tau|+1)
\end{eqnarray*}

According now to the definition of $A_k,L_k$, this formula reads:
\begin{eqnarray*}
G_k(\pi,\sigma)
&=&\sum_{\tau\geq\pi}L_k(\tau,\sigma)\\
&=&\sum_\tau A_k(\pi,\tau)L_k(\tau,\sigma)\\
&=&(A_kL_k)(\pi,\sigma)
\end{eqnarray*}

Thus, we are led to the formula in the statement.
\end{proof}

As an illustration for the above result, at $k=2$ we have $P(2)=\{||,\sqcap\}$, and the above decomposition $G_2=A_2L_2$ appears as follows:
$$\begin{pmatrix}N^2&N\\ N&N\end{pmatrix}
=\begin{pmatrix}1&1\\ 0&1\end{pmatrix}
\begin{pmatrix}N^2-N&0\\N&N\end{pmatrix}$$

We are led in this way to the following formula, due to Lindst\"om \cite{lin}:

\index{Gram matrix}
\index{Gram determinant}
\index{Lindst\"om formula}
\index{linear independence}

\begin{theorem}
The determinant of the Gram matrix $G_k$ is given by
$$\det(G_k)=\prod_{\pi\in P(k)}\frac{N!}{(N-|\pi|)!}$$
with the convention that in the case $N<k$ we obtain $0$.
\end{theorem}

\begin{proof}
If we order $P(k)$ as usual, with respect to the number of blocks, and then lexicographically, $A_k$ is upper triangular, and $L_k$ is lower triangular. Thus, we have:
\begin{eqnarray*}
\det(G_k)
&=&\det(A_k)\det(L_k)\\
&=&\det(L_k)\\
&=&\prod_\pi L_k(\pi,\pi)\\
&=&\prod_\pi N(N-1)\ldots(N-|\pi|+1)
\end{eqnarray*}

Thus, we are led to the formula in the statement.
\end{proof}

Now back to easiness and laws of characters, we can formulate:

\index{main character}
\index{stationary convergence}
\index{asymptotic moments}

\begin{theorem}
For an easy group $G=(G_N)$, coming from a category of partitions $D=(D(k,l))$, the asymptotic moments of the main character are given by
$$\lim_{N\to\infty}\int_{G_N}\chi^k=|D(k)|$$
where $D(k)=D(\emptyset,k)$, with the limiting sequence on the left consisting of certain integers, and being stationary at least starting from the $k$-th term.
\end{theorem}

\begin{proof}
We know from Theorem 14.12 that we have the following formula:
$$\int_{G_N}\chi^k
=\dim\left( span\left(\xi_\pi\Big|\pi\in D(k)\right)\right)$$

Now since by Theorem 14.20 the vectors $\xi_\pi$ are linearly independent with $N\geq k$, and in particular with $N\to\infty$, we obtain the formula in the statement.
\end{proof}

This is very nice, and as a first application of our technology, we can recover the main Poisson law result from chapter 13, at $t=1$, by using easiness, as follows:

\index{Bell numbers}
\index{Poisson law}

\begin{theorem}
For the symmetric group $S_N$, the main character becomes Poisson
$$\chi\sim p_1$$
in the $N\to\infty$ limit.
\end{theorem}

\begin{proof}
According to Theorem 14.21, we have the following formula:
$$\lim_{N\to\infty}\int_{S_N}\chi^k=|P(k)|$$

Now the Bell numbers on the right being the moments of $p_1$, this gives the result.
\end{proof}

Quite powerful all this, hope you agree with me. Gone all the inclusion-exclusion computations from chapter 13, or other tricks there, all we need is easiness.

\section*{14c. Normal variables} 

Welcome to $O_N,U_N$, take three. We have learned many things in this chapter, concerning various groups that we are interested in, but still no word about these. This being said, no worries, because Theorem 14.21 is what we need at $t=1$, and a modification of that, based on the Weingarten formula from chapter 9, will do the job at $t>0$.

\bigskip

In practice, we can safely conjecture that as limiting laws we will obtain the real and complex normal laws, so let us start with a crash course on these normal laws, coming as a continuation of our previous learning of probability, in general. We will need:

\begin{theorem}
We have the following formula,
$$\int_\mathbb Re^{-x^2}dx=\sqrt{\pi}$$
called Gauss integral formula.
\end{theorem}

\begin{proof}
Let $I$ be the above integral. By using polar coordinates, we obtain:
\begin{eqnarray*}
I^2
&=&\int_\mathbb R\int_\mathbb Re^{-x^2-y^2}dxdy\\
&=&\int_0^{2\pi}\int_0^\infty e^{-r^2}rdrdt\\
&=&2\pi\left[-\frac{e^{-r^2}}{2}\right]_0^\infty\\
&=&\pi
\end{eqnarray*}

Thus, we are led to the formula in the statement.
\end{proof}

As a main application of the Gauss formula, we can now formulate:

\begin{definition}
The normal law of parameter $1$ is the following measure:
$$g_1=\frac{1}{\sqrt{2\pi}}e^{-x^2/2}dx$$
More generally, the normal law of parameter $t>0$ is the following measure:
$$g_t=\frac{1}{\sqrt{2\pi t}}e^{-x^2/2t}dx$$
These are also called Gaussian distributions, with ``g'' standing for Gauss.
\end{definition}

Observe that the above laws have indeed mass 1, as they should. This follows indeed from the Gauss formula, which gives, with $x=\sqrt{2t}\,y$:
\begin{eqnarray*}
\int_\mathbb R e^{-x^2/2t}dx
&=&\int_\mathbb R e^{-y^2}\sqrt{2t}\,dy\\
&=&\sqrt{2t}\int_\mathbb R e^{-y^2}dy\\
&=&\sqrt{2t}\times\sqrt{\pi}\\
&=&\sqrt{2\pi t}
\end{eqnarray*}

Generally speaking, the normal laws appear as bit everywhere, in real life. The reasons behind this phenomenon come from the Central Limit Theorem (CLT), that we will explain in a moment, after developing some general theory. As a first result, we have:

\begin{proposition}
We have the variance formula
$$V(g_t)=t$$
valid for any $t>0$.
\end{proposition}

\begin{proof}
The first moment is 0, because our normal law $g_t$ is centered. As for the second moment, this can be computed as follows:
\begin{eqnarray*}
M_2
&=&\frac{1}{\sqrt{2\pi t}}\int_\mathbb Rx^2e^{-x^2/2t}dx\\
&=&\frac{1}{\sqrt{2\pi t}}\int_\mathbb R(tx)\left(-e^{-x^2/2t}\right)'dx\\
&=&\frac{1}{\sqrt{2\pi t}}\int_\mathbb Rte^{-x^2/2t}dx\\
&=&t
\end{eqnarray*}

We conclude from this that the variance is $V=M_2=t$.
\end{proof}

Here is another result, which is the key one for the study of the normal laws:

\begin{theorem}
We have the following formula, valid for any $t>0$:
$$F_{g_t}(x)=e^{-tx^2/2}$$
In particular, the normal laws satisfy $g_s*g_t=g_{s+t}$, for any $s,t>0$.
\end{theorem}

\begin{proof}
The Fourier transform formula can be established as follows:
\begin{eqnarray*}
F_{g_t}(x)
&=&\frac{1}{\sqrt{2\pi t}}\int_\mathbb Re^{-y^2/2t+ixy}dy\\
&=&\frac{1}{\sqrt{2\pi t}}\int_\mathbb Re^{-(y/\sqrt{2t}-\sqrt{t/2}ix)^2-tx^2/2}dy\\
&=&\frac{1}{\sqrt{2\pi t}}\int_\mathbb Re^{-z^2-tx^2/2}\sqrt{2t}dz\\
&=&\frac{1}{\sqrt{\pi}}e^{-tx^2/2}\cdot\sqrt{\pi}\\
&=&e^{-tx^2/2}
\end{eqnarray*}

As for the last assertion, this follows from the fact that $\log F_{g_t}$ is linear in $t$.
\end{proof}

We are now ready to state and prove the CLT, as follows:

\begin{theorem}[CLT]
Given random variables $f_1,f_2,f_3,\ldots\in L^\infty(X)$ which are i.i.d., centered, and with variance $t>0$, we have, with $n\to\infty$, in moments,
$$\frac{1}{\sqrt{n}}\sum_{i=1}^nf_i\sim g_t$$
where $g_t$ is the Gaussian law of parameter $t$, having as density $\frac{1}{\sqrt{2\pi t}}e^{-y^2/2t}dy$.
\end{theorem}

\begin{proof}
In terms of moments, the Fourier transform is given by:
$$F_f(x)
=E\left(\sum_{k=0}^\infty\frac{(ixf)^k}{k!}\right)
=\sum_{k=0}^\infty\frac{i^kM_k(f)}{k!}\,x^k$$

Thus, the Fourier transform of the variable in the statement is:
\begin{eqnarray*}
F(x)
&=&\left[F_f\left(\frac{x}{\sqrt{n}}\right)\right]^n\\
&=&\left[1-\frac{tx^2}{2n}+O(n^{-2})\right]^n\\
&\simeq&\left[1-\frac{tx^2}{2n}\right]^n\\
&\simeq&e^{-tx^2/2}
\end{eqnarray*}

But this latter function being the Fourier transform of $g_t$, we obtain the result.
\end{proof}

Let us discuss now some further properties of the normal law. We first have:

\begin{proposition}
The even moments of the normal law are the numbers
$$M_k(g_t)=t^{k/2}\times k!!$$
where $k!!=(k-1)(k-3)(k-5)\ldots\,$, and the odd moments vanish. 
\end{proposition}

\begin{proof}
We have the following computation, valid for any integer $k\in\mathbb N$:
\begin{eqnarray*}
M_k
&=&\frac{1}{\sqrt{2\pi t}}\int_\mathbb Ry^ke^{-y^2/2t}dy\\
&=&\frac{1}{\sqrt{2\pi t}}\int_\mathbb R(ty^{k-1})\left(-e^{-y^2/2t}\right)'dy\\
&=&\frac{1}{\sqrt{2\pi t}}\int_\mathbb Rt(k-1)y^{k-2}e^{-y^2/2t}dy\\
&=&t(k-1)\times\frac{1}{\sqrt{2\pi t}}\int_\mathbb Ry^{k-2}e^{-y^2/2t}dy\\
&=&t(k-1)M_{k-2}
\end{eqnarray*}

Now recall from the proof of Proposition 14.25 that we have $M_0=1$, $M_1=0$. Thus by recurrence, we are led to the formula in the statement.
\end{proof}

We have the following alternative formulation of the above result:

\begin{proposition}
The moments of the normal law are the numbers
$$M_k(g_t)=t^{k/2}|P_2(k)|$$
where $P_2(k)$ is the set of pairings of $\{1,\ldots,k\}$.
\end{proposition}

\begin{proof}
Let us count the pairings of $\{1,\ldots,k\}$. In order to have such a pairing, we must pair $1$ with one of the numbers $2,\ldots,k$, and then use a pairing of the remaining $k-2$ numbers. Thus, we have the following recurrence formula:
$$|P_2(k)|=(k-1)|P_2(k-2)|$$

As for the initial data, this is $P_1=0$, $P_2=1$. Thus, we are led to the result.
\end{proof}

We are not done yet, and here is one more improvement of the above:

\begin{theorem}
The moments of the normal law are the numbers
$$M_k(g_t)=\sum_{\pi\in P_2(k)}t^{|\pi|}$$
where $P_2(k)$ is the set of pairings of $\{1,\ldots,k\}$, and $|.|$ is the number of blocks.
\end{theorem}

\begin{proof}
This follows indeed from Proposition 14.29, because the number of blocks of a pairing of $\{1,\ldots,k\}$ is trivially $k/2$, independently of the pairing.
\end{proof}

Moving on, let us discuss now the complex analogues of all this, with a notion of complex normal, or Gaussian law. To start with, we have the following definition:

\begin{definition}
The complex normal, or Gaussian law of parameter $t>0$ is
$$G_t=law\left(\frac{1}{\sqrt{2}}(a+ib)\right)$$
where $a,b$ are independent, each following the law $g_t$.
\end{definition}

In short, the complex normal laws appear as natural complexifications of the real normal laws. As in the real case, these measures form convolution semigroups:

\begin{proposition}
The complex Gaussian laws have the property
$$G_s*G_t=G_{s+t}$$
for any $s,t>0$, and so they form a convolution semigroup.
\end{proposition}

\begin{proof}
This follows indeed from the real result, namely $g_s*g_t=g_{s+t}$, established in Theorem 14.26, simply by taking the real and imaginary parts.
\end{proof}

We have as well the following complex analogue of the CLT:

\index{CCLT}
\index{Complex CLT}

\begin{theorem}[CCLT]
Given complex variables $f_1,f_2,f_3,\ldots\in L^\infty(X)$ which are i.i.d., centered, and with common variance $t>0$, we have
$$\frac{1}{\sqrt{n}}\sum_{i=1}^nf_i\sim G_t$$
with $n\to\infty$, in moments.
\end{theorem}

\begin{proof}
This follows indeed from the real CLT, established in Theorem 14.27, simply by taking the real and imaginary parts of all the variables involved.
\end{proof}

Regarding now the moments, the situation here is more complicated than in the real case, because in order to have good results, we have to deal with both the complex variables, and their conjugates. Let us formulate the following definition:

\index{colored integers}
\index{colored moments}

\begin{definition}
The moments a complex variable $f\in L^\infty(X)$ are the numbers
$$M_k=E(f^k)$$
depending on colored integers $k=\circ\bullet\bullet\circ\ldots\,$, with the conventions
$$f^\emptyset=1\quad,\quad f^\circ=f\quad,\quad f^\bullet=\bar{f}$$
and multiplicativity, in order to define the colored powers $f^k$.
\end{definition}

Observe that, since $f,\bar{f}$ commute, we can permute terms, and restrict the attention to exponents of type $k=\ldots\circ\circ\circ\bullet\bullet\bullet\ldots\,$, if we want to. However, our results about the complex Gaussian laws, and other complex laws too, later on, will actually look better without doing is, so we will use Definition 14.34 as stated. We first have:

\begin{theorem}
The moments of the complex normal law are given by
$$M_k(G_t)=\begin{cases}
t^pp!&(k\ {\rm uniform, of\ length}\ 2p)\\
0&(k\ {\rm not\ uniform})
\end{cases}$$
where $k=\circ\bullet\bullet\circ\ldots$ is called uniform when it contains the same number of $\circ$ and $\bullet$.
\end{theorem}

\begin{proof}
We must compute the moments, with respect to colored integer exponents $k=\circ\bullet\bullet\circ\ldots$\,, of the variable from Definition 14.31, namely:
$$f=\frac{1}{\sqrt{2}}(a+ib)$$

We can assume that we are in the case $t=1$, and the proof here goes as follows:

\medskip

(1) As a first observation, in the case where our exponent $k=\circ\bullet\bullet\circ\ldots$ is not uniform, a standard rotation argument shows that the corresponding moment of $f$ vanishes. To be more precise, the variable $f'=wf$ is complex Gaussian too, for any complex number $w\in\mathbb T$, and from $M_k(f)=M_k(f')$ we obtain $M_k(f)=0$, in this case.

\medskip

(2) In the uniform case now, where the exponent $k=\circ\bullet\bullet\circ\ldots$ consists of $p$ copies of $\circ$ and $p$ copies of $\bullet$\,, the corresponding moment can be computed as follows:
\begin{eqnarray*}
M_k
&=&\int(f\bar{f})^p\\
&=&\frac{1}{2^p}\int(a^2+b^2)^p\\
&=&\frac{1}{2^p}\sum_r\binom{p}{r}\int a^{2r}\int b^{2p-2r}\\
&=&\frac{1}{2^p}\sum_r\binom{p}{r}(2r)!!(2p-2r)!!\\
&=&\frac{1}{2^p}\sum_r\frac{p!}{r!(p-r)!}\cdot\frac{(2r)!}{2^rr!}\cdot\frac{(2p-2r)!}{2^{p-r}(p-r)!}\\
&=&\frac{p!}{4^p}\sum_r\binom{2r}{r}\binom{2p-2r}{p-r}
\end{eqnarray*}

(3) In order to finish now the computation, let us recall that we have the following formula, coming from the generalized binomial formula, or from the Taylor formula:
$$\frac{1}{\sqrt{1+t}}=\sum_{q=0}^\infty\binom{2q}{q}\left(\frac{-t}{4}\right)^q$$

By taking the square of this series, we obtain the following formula:
$$\frac{1}{1+t}
=\sum_p\left(\frac{-t}{4}\right)^p\sum_r\binom{2r}{r}\binom{2p-2r}{p-r}$$

Now by looking at the coefficient of $t^p$ on both sides, we conclude that the sum on the right equals $4^p$. Thus, we can finish the moment computation in (2), as follows:
$$M_k=\frac{p!}{4^p}\times 4^p=p!$$

We are therefore led to the conclusion in the statement.
\end{proof}

As before with the real laws, a better-looking statement  is in terms of partitions:

\index{matching pairings}

\begin{theorem}
The moments of the complex normal law are the numbers
$$M_k(G_t)=\sum_{\pi\in\mathcal P_2(k)}t^{|\pi|}$$
where $\mathcal P_2(k)$ are the matching pairings of $\{1,\ldots,k\}$, and $|.|$ is the number of blocks.
\end{theorem}

\begin{proof}
This is a reformulation of Theorem 14.35. Indeed, we can assume that we are in the case $t=1$, and here we know from Theorem 14.35 that the moments are:
$$M_k=\begin{cases}
(|k|/2)!&(k\ {\rm uniform})\\
0&(k\ {\rm not\ uniform})
\end{cases}$$

On the other hand, the numbers $|\mathcal P_2(k)|$ are given by exactly the same formula. Indeed, in order to have a matching pairing of $k$, our exponent $k=\circ\bullet\bullet\circ\ldots$ must be uniform, consisting of $p$ copies of $\circ$ and $p$ copies of $\bullet$, with $p=|k|/2$. But then the matching pairings of $k$ correspond to the permutations of the $\bullet$ symbols, as to be matched with $\circ$ symbols, and so we have $p!$ such pairings. Thus, we have the same formula as for the moments of $f$, and we are led to the conclusion in the statement.
\end{proof}

\section*{14d. Rotation groups}

Welcome to $O_N,U_N$, take four. I think we can do now our computations. To start with, we have the following result, which solves the problem at $t=1$:

\begin{theorem}
The asymptotic laws of main characters are as follows:
\begin{enumerate}
\item For $O_N$, we obtain the real normal law $g_1$.

\item For $U_N$, we obtain the complex normal law $G_1$.
\end{enumerate}
\end{theorem}

\begin{proof}
This comes from Theorem 14.21, which for $O_N,U_N$ gives:
$$\int_{O_N}\chi^k\simeq|P_2(k)|\quad,\quad \int_{U_N}\chi^k\simeq|\mathcal P_2(k)|$$

Now these numbers being the moments of $g_1,G_1$, we obtain the result.
\end{proof}

In order to deal now with the general case $t>0$, we will need a modification of Theorem 14.21, based on the integration results from chapter 9. Let us start with:

\index{Gram matrix}
\index{Weingarten matrix}
\index{Weingarten formula}
\index{Haar functional}

\begin{theorem}
Assuming that a subgroup $G\subset U_N$ is easy, coming from a category of partitions $D\subset P$, we have the Weingarten formula
$$\int_Gu_{i_1j_1}^{e_1}\ldots u_{i_kj_k}^{e_k}=\sum_{\pi,\sigma\in D(k)}\delta_\pi(i)\delta_\sigma(j)W_{kN}(\pi,\sigma)$$
where $\delta\in\{0,1\}$ are the usual Kronecker type symbols, and where the Weingarten matrix $W_{kN}=G_{kN}^{-1}$ is the inverse of the Gram matrix $G_{kN}(\pi,\sigma)=N^{|\pi\vee\sigma|}$. 
\end{theorem}

\begin{proof}
We know from chapter 9 that the integrals in the statement, with $i,j$ varying, form altogether the orthogonal projection $P$ onto the following space:
$$Fix(u^{\otimes k})=span\left(\xi_\pi\Big|\pi\in D(k)\right)$$

In order to prove the result, consider the following linear map:
$$E(x)=\sum_{\pi\in D(k)}<x,\xi_\pi>\xi_\pi$$

By a standard linear algebra computation, it follows that we have $P=WE$, where $W$ is the inverse on $Fix(u^{\otimes k})$ of the restriction of $E$. But this restriction is the linear map coming from the Gram matrix of the vectors $\xi_\pi$, which is given by:
$$G_{kN}(\pi,\sigma)=<\xi_\pi,\xi_\sigma>=N^{|\pi\vee\sigma|}$$

Thus $W$ is the linear map coming from $W_{kN}=G_{kN}^{-1}$, and this gives the result.
\end{proof}

Observe that there is a bit of confusion in Theorem 14.38, because the partitions in $D(k)$ are known to be linearly independent only at $N\geq k$, and so their Gram matrix is invertible only when adding such a $N\geq k$ assumption. In what follows we will only use Theorem 14.38 in the $N>>0$ regime, where the above formulation will do. In general, the point is that $W_{kN}$ should be taken to be the quasi-inverse of $G_{kN}$. See \cite{csn}.

\bigskip

In relation now with truncated characters, we have the following formula:

\index{truncated character}

\begin{theorem}
The moments of truncated characters are given by the formula
$$\int_G(u_{11}+\ldots +u_{ss})^k=Tr(W_{kN}G_{ks})$$
where $G_{kN}$ and $W_{kN}=G_{kN}^{-1}$ are the associated Gram and Weingarten matrices.
\end{theorem}

\begin{proof}
We have the following computation, using Theorem 14.38:
\begin{eqnarray*}
\int_G(u_{11}+\ldots +u_{ss})^k
&=&\sum_{i_1=1}^{s}\ldots\sum_{i_k=1}^s\int u_{i_1i_1}\ldots u_{i_ki_k}\\
&=&\sum_{\pi,\sigma\in D(k)}W_{kN}(\pi,\sigma)\sum_{i_1=1}^{s}\ldots\sum_{i_k=1}^s\delta_\pi(i)\delta_\sigma(i)\\
&=&\sum_{\pi,\sigma\in D(k)}W_{kN}(\pi,\sigma)G_{ks}(\sigma,\pi)\\
&=&Tr(W_{kN}G_{ks})
\end{eqnarray*}

Thus, we have obtained the formula in the statement.
\end{proof}

In order to advance, we must impose an uniformity condition, coming from:

\index{uniformity}

\begin{proposition}
For an easy group $G=(G_N)$, coming from a category of partitions $D\subset P$, the following conditions are equivalent:
\begin{enumerate}
\item $G_{N-1}=G_N\cap U_{N-1}$, via the embedding $U_{N-1}\subset U_N$ given by $u\to diag(u,1)$.

\item $G_{N-1}=G_N\cap U_{N-1}$, via the $N$ possible diagonal embeddings $U_{N-1}\subset U_N$.

\item $D$ is stable under the operation which consists in removing blocks.
\end{enumerate}
If these conditions are satisfied, we say that $G=(G_N)$ is uniform.
\end{proposition}

\begin{proof}
This is something that we already met in chapter 11, in the context of our classification discussion there for the easy groups, the idea being as follows:

\medskip

$(1)\iff(2)$ This equivalence is elementary, coming from the inclusion $S_N\subset G_N$, which makes everything $S_N$-invariant. 

\medskip

$(1)\iff(3)$ Given a closed subgroup $K\subset U_{N-1}$, with fundamental representation $u$, consider the following $N\times N$ matrix:
$$v=\begin{pmatrix}u&0\\0&1\end{pmatrix}$$

Then for any $\pi\in P(k)$ a standard computation shows that we have:
$$\xi_\pi\in Fix(v^{\otimes k})\iff\xi_{\pi'}\in Fix(v^{\otimes k'}),\,\forall\pi'\in P(k'),\pi'\subset\pi$$

Now with this in hand, the result follows from Tannakian duality.
\end{proof}
 
Getting back now to the truncated characters, we have the following result:

\begin{theorem}
For a uniform easy group $G=(G_N)$, we have the formula
$$\lim_{N\to\infty}\int_{G_N}\chi_t^k=\sum_{\pi\in D(k)}t^{|\pi|}$$
with $D\subset P$ being the associated category of partitions.
\end{theorem}

\begin{proof}
In the uniform case the Gram matrix, and so the Weingarten matrix too, are asymptotically diagonal, so the asymptotic moments are given by:
\begin{eqnarray*}
\int_{G_N}\chi_t^k
&=&Tr(W_{kN}G_{k[tN]})\\
&\simeq&\sum_{\pi\in D(k)}N^{-|\pi|}[tN]^{|\pi|}\\
&\simeq&\sum_{\pi\in D(k)}t^{|\pi|}
\end{eqnarray*}

Thus, we are led to the conclusion in the statement.
\end{proof}

Good news, we can eventually finish our work on $O_N,U_N$, with:

\begin{theorem}
The asymptotic laws of truncated characters are as follows:
\begin{enumerate}
\item For $O_N$, we obtain the real normal law $g_t$.

\item For $U_N$, we obtain the complex normal law $G_t$.
\end{enumerate}
\end{theorem}

\begin{proof}
This comes from Theorem 14.41, which for $O_N,U_N$ gives:
$$\int_{O_N}\chi_t^k\simeq\sum_{\pi\in P_2(k)}t^{|\pi|}\quad,\quad 
\int_{U_N}\chi_t^k\simeq\sum_{\pi\in\mathcal P_2(k)}t^{|\pi|}$$

Now these numbers being the moments of $g_t,G_t$, we obtain the result.
\end{proof}

Next, it is quite clear that our technology will apply in fact to all the easy groups that we know, and we have the following grand result here, closing the discussion:

\begin{theorem}
The asymptotic laws of truncated characters are as follows:
\begin{enumerate}
\item For $S_N$, we obtain the Poisson law $p_t$.

\item For $H_N$, we obtain the real Bessel law $b_t$.

\item For $K_N$, we obtain the complex Bessel law $B_t$.

\item For $H_N^s$, we obtain the general Bessel law $b^s_t$.

\item For $O_N$, we obtain the real normal law $g_t$.

\item For $U_N$, we obtain the complex normal law $G_t$.

\item For $B_N$, we obtain the shifted real normal law $g_t'$.

\item For $C_N$, we obtain the shifted complex normal law $G_t'$.
\end{enumerate}
\end{theorem}

\begin{proof}
This collects the various results that we know, along with a bit more at the end, and with the claim that all this is best viewed as coming from easiness:

\medskip

(1-4) These are things that we know from chapter 13, obtained there via various counting methods, but as explained in Theorem 14.22, the central result, for $G=S_N$ at $t=1$, is best viewed via easiness. And the same goes for any $G=H_N^s$ and any $t>0$, by using Theorem 14.41, and a bit more combinatorics, that we will leave as an exercise.

\medskip

(5-8) Here the results for $O_N,U_N$ are those that we just discussed, and those for $B_N,C_N$ follow in a similar way, the idea being that the singletons produce a shift. Alternatively, we can get these latter results by using the isomorphisms $B_N\simeq O_{N-1}$ and $C_N\simeq U_{N-1}$ discussed in chapter 1, and we will leave the work here as an exercise too.
\end{proof}

As a comment now, the above result is not end of the story with the groups in question, because passed the characters, truncated or not, there are many other interesting variables that can be considered. Also, in what regards the laws that we found, while the normal and Poisson laws are certainly well-understood, and the same goes, to a certain extent, for the real Bessel laws $b_t$, in what regards the general Bessel laws $b^s_t$, many other interesting things can be said, for instance regarding the densities. For more here, we refer to \cite{bbc}.

\bigskip

As another comment, passed the easiness framework, we have already mentioned, in chapter 13, that for the general reflection groups $H_N^{sd}\subset H_N^s$ we obtain the same asymptotic laws as for $H_N^s$, and the same can be shown to hold for the groups $U_N^d\subset U_N$.

\bigskip

Finally, we have special groups like $SU_2,SO_3$, or the symplectic groups $Sp_N\subset U_N$, which are not easy, but are not far from being easy either. These groups are best investigated by using ``super-easiness'', that is, easiness involving some signs in the implementation of the partitions $\pi\to T_\pi$, in the spirit of our general ``super-space'' discussion from chapter 1, and for more on this, we refer to \cite{csn} and related papers.

\section*{14e. Exercises}

Good time for a remake of our exercise session from chapter 12, as follows:

\begin{exercise}
Further study $S_N$, both representations and characters.
\end{exercise}

\begin{exercise}
Further study $H_N$, both representations and characters.
\end{exercise}

\begin{exercise}
Further study $K_N$, both representations and characters.
\end{exercise}

\begin{exercise}
Further study $H_N^s$, both representations and characters.
\end{exercise}

\begin{exercise}
Further study $O_N$, both representations and characters.
\end{exercise}

\begin{exercise}
Further study $U_N$, both representations and characters.
\end{exercise}

\begin{exercise}
Further study $B_N$, both representations and characters.
\end{exercise}

\begin{exercise}
Further study $C_N$, both representations and characters.
\end{exercise}

As bonus exercise, do not forget about $SU_2,SO_3$, and about $Sp_N\subset U_N$ too.

\chapter{Liberation theory}

\section*{15a. Free probability}

With most of the groups $G\subset U_N$ discussed, end of the story, you would say. Well, not so quick, because remember from chapter 12 the main principle established there, namely that quantum is simpler than classical, and should therefore come first.

\bigskip

Normally this principle should be valid for analytic questions too, which after all are just analytic reformulations of the algebra that we know. So, in the context of what we have been doing so far in the present Part IV, the following question appears:

\begin{question}
What are the asymptotic character laws for $S_N^+,O_N^+,U_N^+$ and the other quantum groups that we know? Do these liberate the laws for $S_N,O_N,U_N$ and the other groups that we know, in which exact sense, and do things get indeed simpler in this way? What about other invariants, such as Gram determinants, do these liberate too?
\end{question}

And good question this is. In answer now, in view of what we have been doing in chapter 12, we can safely conjecture that the answer is yes. As for the proof, what we have been doing recently, in chapter 14, obviously does not rely on the fact that the groups there are classical, so with a bit of combinatorial work, we will have the results. 

\bigskip

This being said, without some further training in probability, of free type, we are at risk of not understanding much, what we are doing, and even worse, we are at risk of compromising the yes answer to the ``things get simpler'' part of the above question.

\bigskip

In short, welcome to free probability, to start with, and for Question 15.1, we will see later. Let us start our discussion with something coming from chapter 12, namely:

\begin{definition}
Let $A$ be a $C^*$-algebra, given with a trace $tr:A\to\mathbb C$.
\begin{enumerate}
\item The elements $a\in A$ are called random variables.

\item The moments of such a variable are the numbers $M_k(a)=tr(a^k)$.

\item The law of such a variable is the functional $\mu:P\to tr(P(a))$.
\end{enumerate}
\end{definition}

Here the exponent $k=\circ\bullet\bullet\circ\ldots$ is as before a colored integer, with the powers $a^k$ being defined by multiplicativity and the usual formulae, namely:
$$a^\emptyset=1\quad,\quad
a^\circ=a\quad,\quad
a^\bullet=a^*$$

As for the polynomial $P$, this is a noncommuting $*$-polynomial in one variable: 
$$P\in\mathbb C<X,X^*>$$

Generally speaking, the above definition is something quite abstract, but there is no other way of doing things, at least at this level of generality. However, in the special case where our variable $a\in A$ is self-adjoint, or more generally normal, we have:

\index{spectral measure}
\index{functional calculus}
\index{Riesz theorem}

\begin{proposition}
The law of a normal variable $a\in A$ can be identified with the corresponding spectral measure $\mu\in\mathcal P(\mathbb C)$, according to the following formula,
$$tr(f(a))=\int_{\sigma(a)}f(x)d\mu(x)$$
valid for any $f\in L^\infty(\sigma(a))$, coming from the Riesz and Gelfand theorems. In the self-adjoint case the spectral measure is real, $\mu\in\mathcal P(\mathbb R)$.
\end{proposition}

\begin{proof}
This is something very standard, coming from the Riesz theorem, via the Gelfand theorem for the $C^*$-algebras, that we know well since chapter 12.
\end{proof}

Let us discuss now independence, and its noncommutative versions. As a starting point, we have the following update of the classical notion of independence:

\index{abstract independence}
\index{independent algebras}

\begin{definition}
We call two subalgebras $B,C\subset A$ independent when the following condition is satisfied, for any $x\in B$ and $y\in C$: 
$$tr(xy)=tr(x)tr(y)$$
Equivalently, the following condition must be satisfied, for any $x\in B$ and $y\in C$: 
$$tr(x)=tr(y)=0\implies tr(xy)=0$$
Also, $b,c\in A$ are called independent when $B=<b>$ and $C=<c>$ are independent.
\end{definition}

Observe that the above two independence conditions are indeed equivalent, with this following from the following computation, with the convention $x'=x-tr(x)$:
\begin{eqnarray*}
tr(xy)
&=&tr[(x'+tr(x))(y'+tr(y))]\\
&=&tr(x'y')+t(x')tr(y)+tr(x)tr(y')+tr(x)tr(y)\\
&=&tr(x'y')+tr(x)tr(y)\\
&=&tr(x)tr(y)
\end{eqnarray*}

It is possible to develop some theory here, but this leads to the usual CLT. As a much more interesting notion now, we have Voiculescu's freeness \cite{voi}:

\index{freeness}
\index{free algebras}

\begin{definition}
Given a pair $(A,tr)$, we call two subalgebras $B,C\subset A$ free when the following condition is satisfied, for any $x_i\in B$ and $y_i\in C$:
$$tr(x_i)=tr(y_i)=0\implies tr(x_1y_1x_2y_2\ldots)=0$$
Also, $b,c\in A$ are called free when $B=<b>$ and $C=<c>$ are free.
\end{definition}

As a first observation on this notion, which is similar to independence, there is a certain lack of symmetry between Definition 15.4 and Definition 15.5, because the latter does not include an explicit formula for the following type of quantities:
$$tr(x_1y_1x_2y_2\ldots)$$

However, this is not an issue, and is simply due to the fact that the formula in the free case is something more complicated, the precise result being as follows:

\index{free trace}

\begin{proposition}
If $B,C\subset A$ are free, the restriction of $tr$ to $<B,C>$ can be computed in terms of the restrictions of $tr$ to $B,C$. To be more precise, we have
$$tr(x_1y_1x_2y_2\ldots)=P\Big(\{tr(x_{i_1}x_{i_2}\ldots)\}_i,\{tr(y_{j_1}y_{j_2}\ldots)\}_j\Big)$$
where $P$ is certain polynomial, depending on the length of $x_1y_1x_2y_2\ldots\,$, having as variables the traces of products $x_{i_1}x_{i_2}\ldots$ and $y_{j_1}y_{j_2}\ldots\,$, with $i_1<i_2<\ldots$ and $j_1<j_2<\ldots$
\end{proposition}

\begin{proof}
This is something a bit theoretical, so let us begin with an example. The computation after Definition 15.4 perfectly works when $x,y$ are free, and gives:
$$tr(xy)=tr(x)tr(y)$$

In general, the situation is of course more complicated than this, but the same trick applies. To be more precise, we can start our computation as follows:
\begin{eqnarray*}
tr(x_1y_1x_2y_2\ldots)
&=&tr\big[(x_1'+tr(x_1))(y_1'+tr(y_1))(x_2'+tr(x_2))\ldots\big]\\
&=&tr(x_1'y_1'x_2'y_2'\ldots)+{\rm other\ terms}\\
&=&{\rm other\ terms}
\end{eqnarray*}

Now regarding the ``other terms'', those which are left, each of them will consist of a product of traces of type $tr(x_i)$ and $tr(y_i)$, and then a trace of a product still remaining to be computed, which is of the following form, for some elements $\beta_i\in B$ and $\gamma_i\in C$:
$$tr(\beta_1\gamma_1\beta_2\gamma_2\ldots)$$

Now since the length of $\beta_1\gamma_1\beta_2\gamma_2\ldots$ is smaller than the length of the original product $x_1y_1x_2y_2\ldots$, we are led into of recurrence, and this gives the result.
\end{proof}

Let us discuss now some examples of independence and freeness. We first have:

\begin{theorem}
Given two algebras $(A,tr)$ and $(B,tr)$, the following hold:
\begin{enumerate}
\item $A,B$ are independent inside their tensor product $A\otimes B$, endowed with its canonical tensor product trace, given on basic tensors by $tr(x\otimes y)=tr(x)tr(y)$.

\item $A,B$ are free inside their free product $A*B$, endowed with its canonical free product trace, given by the formulae in Proposition 15.6.
\end{enumerate}
\end{theorem}

\begin{proof}
Both the above assertions are clear from definitions, as follows:

\medskip

(1) This is clear with either of the definitions of the independence, from Definition 15.4, because we have by construction of the product trace:
\begin{eqnarray*}
tr(xy)
&=&tr[(x\otimes1)(1\otimes y)]\\
&=&tr(x\otimes y)\\
&=&tr(x)tr(y)
\end{eqnarray*}

(2) This is clear too from definitions, the only point being that of showing that the notion of freeness, or the recurrence formulae in Proposition 15.6, can be used in order to construct a canonical free product trace, on the free product of the algebras involved:
$$tr:A*B\to\mathbb C$$

But this can be checked for instance by using a GNS construction. Indeed, consider the GNS constructions for the algebras $(A,tr)$ and $(B,tr)$:
$$A\to B(l^2(A))\quad,\quad
B\to B(l^2(B))$$

By taking the free product of these representations, we obtain a representation as follows, with the $*$ on the right being a free product of pointed Hilbert spaces:
$$A*B\to B(l^2(A)*l^2(B))$$

Now by composing with the linear form $T\to<T\xi,\xi>$, where $\xi=1_A=1_B$ is the common distinguished vector of $l^2(A)$, $l^2(B)$, we obtain a linear form, as follows:
$$tr:A*B\to\mathbb C$$

It is routine then to check that $tr$ is indeed a trace, and this is the ``canonical free product trace'' from the statement. Then, an elementary computation shows that $A,B$ are free inside $A*B$, with respect to this trace, and this finishes the proof.
\end{proof}

More concretely now, we have the following result, also from \cite{voi}:

\begin{theorem}
We have the following results, valid for group algebras:
\begin{enumerate}
\item $C^*(\Gamma),C^*(\Lambda)$ are independent inside $C^*(\Gamma\times\Lambda)$.

\item $C^*(\Gamma),C^*(\Lambda)$ are free inside $C^*(\Gamma*\Lambda)$.
\end{enumerate}
\end{theorem}

\begin{proof}
In order to prove these results, we can use the general results in Theorem 15.7, along with the following two isomorphisms, which are both standard:
$$C^*(\Gamma\times\Lambda)=C^*(\Lambda)\otimes C^*(\Gamma)\quad,\quad
C^*(\Gamma*\Lambda)=C^*(\Lambda)*C^*(\Gamma)$$

Alternatively, we can check the independence and freeness formulae on group elements, which is something trivial, and then conclude by linearity.
\end{proof}

We have seen before the basic limiting theorems in classical probability. In order to deal now with freeness, let us develop some tools. First, we have:

\index{free convolution}

\begin{theorem}
We have a free convolution operation $\boxplus$ for the distributions
$$\mu:\mathbb C<X,X^*>\to\mathbb C$$
which is well-defined by the following formula, with $a,b$ taken to be free:
$$\mu_a\boxplus\mu_b=\mu_{a+b}$$
This restricts to an operation, still denoted $\boxplus$, on the real probability measures.
\end{theorem}

\begin{proof}
We basically have to prove that if two variables $a,b$ are free, then the distribution of their sum $\mu_{a+b}$ depends only on the distributions $\mu_a,\mu_b$. But for this purpose, we can use the general formula from Proposition 15.6, namely:
$$tr(a_1b_1a_2b_2\ldots)=P\Big(\{tr(a_{i_1}a_{i_2}\ldots)\}_i,\{tr(b_{j_1}b_{j_2}\ldots)\}_j\Big)$$

Now by plugging in arbitrary powers of $a,b$ as variables $a_i,b_j$, we obtain a family of formulae of the following type, with $Q$ being certain polyomials:
$$tr(a^{k_1}b^{l_1}a^{k_2}b^{l_2}\ldots)=Q\Big(\{tr(a^k)\}_k,\{tr(b^l)\}_l\Big)$$

Thus the moments of $a+b$ depend only on the moments of $a,b$, with of course colored exponents in all this, according to our moment conventions, and this gives the result.
\end{proof}

Here is now the key technical result of free probability, due to Voiculescu \cite{voi}:

\index{Cauchy transform}
\index{R-transform}

\begin{theorem}
Given a probability measure $\mu$, define its $R$-transform as follows:
$$G_\mu(\xi)=\int_\mathbb R\frac{d\mu(t)}{\xi-t}\implies G_\mu\left(\ R_\mu(\xi)+\frac{1}{\xi}\right)=\xi$$
The free convolution operation is then linearized by the $R$-transform.
\end{theorem}

\begin{proof}
This is something quite tricky, the idea being as follows:

\medskip

(1) In order to model the free convolution, the best is to use creation operators on free Fock spaces, corresponding to the semigroup operator algebras $C^*(\mathbb N^{*k})$. Indeed, we have some freeness here, a bit in the same way as in the free group algebras $C^*(F_k)$.

\medskip

(2) The point now, motivating this choice, is that the variables of type $S^*+f(S)$, with $S\in L(\mathbb N)$ being the shift, and with $f\in\mathbb C[X]$ being an arbitrary  polynomial, are easily seen to model in moments all the possible distributions $\mu:\mathbb C[X]\to\mathbb C$.

\medskip

(3) Now let $f,g\in\mathbb C[X]$ and consider the variables $S^*+f(S)$ and $T^*+g(T)$, where $S,T\in L(\mathbb N*\mathbb N)$ are the shifts corresponding to the generators of $\mathbb N*\mathbb N$. These variables are free, and by using a $45^\circ$ argument, their sum has the same law as $S^*+(f+g)(S)$.

\medskip

(4) Thus the operation $\mu\to f$ linearizes the free convolution. We are therefore left with a computation inside $L(\mathbb N)$, which is elementary, and whose conclusion is that $R_\mu=f$ can be recaptured from $\mu$ via the Cauchy transform $G_\mu$, as in the statement.
\end{proof}

We have the following remarkable free analogue of the CLT, also from \cite{voi}: 

\index{free CLT}
\index{semicircle law}
\index{Wigner law}

\begin{theorem}[FCLT]
Given self-adjoint variables $x_1,x_2,x_3,\ldots,$ which are f.i.d., centered, with variance $t>0$, we have, with $n\to\infty$, in moments,
$$\frac{1}{\sqrt{n}}\sum_{i=1}^nx_i\sim\gamma_t$$\
where $\gamma_t=\frac{1}{2\pi t}\sqrt{4t-x^2}dx$ is the Wigner semicircle law of parameter $t$.
\end{theorem}

\begin{proof}
At $t=1$, the $R$-transform of the variable in the statement can be computed by using the linearization property from Theorem 15.10, and is given by:
$$R(\xi)
=nR_x\left(\frac{\xi}{\sqrt{n}}\right)
\simeq\xi$$

On the other hand, some standard computations, that we will leave as an exercise, show that the Cauchy transform of the Wigner law $\gamma_1$ satisfies the following equation:
$$G_{\gamma_1}\left(\xi+\frac{1}{\xi}\right)=\xi$$

Thus $R_{\gamma_1}(\xi)=\xi$, so the laws in the statement have the same $R$-transforms, and so they are equal. The passage to the general case, $t>0$, is routine, by dilation.
\end{proof}

In the complex case now, we have a similar result, also from \cite{voi}, as follows:

\begin{theorem}[FCCLT]
Given random variables $x_1,x_2,x_3,\ldots$ which are f.i.d., centered, with variance $t>0$, we have, with $n\to\infty$, in moments,
$$\frac{1}{\sqrt{n}}\sum_{i=1}^nx_i\sim\Gamma_t$$
where $\Gamma_t=law\big((a+ib)/\sqrt{2}\big)$, with $a,b$ being free, each following the Wigner semicircle law $\gamma_t$, is the Voiculescu circular law of parameter $t$.
\end{theorem}

\begin{proof}
This follows indeed from the FCLT, by taking real and imaginary parts.
\end{proof}

We have as well a free Poisson Limiting Theorem, as follows:

\begin{theorem}[FPLT]
The following Poisson limit converges, for any $t>0$,
$$\pi_t=\lim_{n\to\infty}\left(\left(1-\frac{t}{n}\right)\delta_0+\frac{t}{n}\delta_1\right)^{\boxplus n}$$
the limiting measure being the Marchenko-Pastur law $\pi_t$, given by
$$\pi_t=\max(1-t,0)\delta_0+\frac{\sqrt{4t-(x-1-t)^2}}{2\pi x}\,dx$$
with at $t=1$, this Marchenko-Pastur law being $\pi_1=\frac{1}{2\pi}\sqrt{4x^{-1}-1}\,dx$. 
\end{theorem}

\begin{proof}
Let $\mu$ be the measure in the statement, appearing under the convolution sign. The Cauchy transform of this measure is elementary to compute, given by:
$$G_{\mu}(\xi)=\left(1-\frac{t}{n}\right)\frac{1}{\xi}+\frac{t}{n}\cdot\frac{1}{\xi-1}$$

By using Theorem 15.10, we want to compute the following $R$-transform:
$$R
=R_{\mu^{\boxplus n}}(y)
=nR_\mu(y)$$

We know that the equation for this function $R$ is as follows:
$$\left(1-\frac{t}{n}\right)\frac{1}{y^{-1}+R/n}+\frac{t}{n}\cdot\frac{1}{y^{-1}+R/n-1}=y$$

Now with $n\to\infty$ we obtain from this the following formula:
$$R=\frac{t}{1-y}$$

But this being the $R$-transform of $\pi_t$, via some calculus, we are done.
\end{proof}

Finally, we have the following ``compound'' generalization of the FPLT:

\index{CPLT}
\index{Compound PLT}
\index{compound Poisson law}

\begin{theorem}[FCPLT]
Given a compactly supported positive measure $\rho$, of mass $c=mass(\rho)$, the following compound Poisson limit converges,
$$\pi_\rho=\lim_{n\to\infty}\left(\left(1-\frac{c}{n}\right)\delta_0+\frac{1}{n}\rho\right)^{\boxplus n}$$
and if we write $\rho=\sum_{i=1}^sc_i\delta_{z_i}$ with $c_i>0$ and $z_i\in\mathbb C$ we have the formula
$$\pi_\rho={\rm law}\left(\sum_{i=1}^sz_i\alpha_i\right)$$
where the variables $\alpha_i$ are free, following the Marchenko-Pastur laws $\pi_{c_i}$.
\end{theorem}

\begin{proof}
Let $\eta_n$ be the measure appearing in the above limit, under the free convolution sign. The Cauchy transform of $\eta_n$ is then given by the following formula:
$$G_{\eta_n}(\xi)=\left(1-\frac{c}{n}\right)\frac{1}{\xi}+\frac{1}{n}\sum_{i=1}^s\frac{c_i}{\xi-z_i}$$

By using Theorem 15.10, the equation of $R=R_{\eta_n^{\boxplus n}}$ is then as follows:
\begin{eqnarray*}
&&\left(1-\frac{c}{n}\right)\frac{1}{1/y+R/n}+\frac{1}{n}\sum_{i=1}^s\frac{c_i}{1/y+R/n-z_i}=y\\
&\implies&\left(1-\frac{c}{n}\right)\frac{1}{1+yR/n}+\frac{1}{n}\sum_{i=1}^s\frac{c_i}{1+yR/n-yz_i}=1
\end{eqnarray*}

Now multiplying by $n$, then rearranging the terms, and letting $n\to\infty$, we get:
\begin{eqnarray*}
\frac{c+yR}{1+yR/n}=\sum_{i=1}^s\frac{c_i}{1+yR/n-yz_i}
&\implies&c+yR_{\pi_\rho}(y)=\sum_{i=1}^s\frac{c_i}{1-yz_i}\\
&\implies&R_{\pi_\rho}(y)=\sum_{i=1}^s\frac{c_iz_i}{1-yz_i}
\end{eqnarray*}

On the other hand, let $\alpha=\sum_iz_i\alpha_i$ be the sum in the statement. We have:
\begin{eqnarray*}
R_{\alpha_i}(y)=\frac{c_i}{1-y}
&\implies&R_{z_i\alpha_i}(y)=\frac{c_iz_i}{1-yz_i}\\
&\implies&R_\alpha(y)=\sum_{i=1}^s\frac{c_iz_i}{1-yz_i}
\end{eqnarray*}

Thus we have the same formula as the one above, and this gives the result.
\end{proof}

So long for the basics of free probability theory, quickly explained. For more on all this there are many good books available, including the original one \cite{vdn}.

\section*{15b. Cumulants, bijection}

We would like to discuss now the correspondence between classical and free, which is something that we must better understand. Following Rota, we first have:

\index{cumulant}
\index{Fourier transform}

\begin{definition}
Associated to any real probability measure $\mu=\mu_f$ is the following modification of the logarithm of the Fourier transform $F_\mu(\xi)=E(e^{i\xi f})$,
$$K_\mu(\xi)=\log E(e^{\xi f})$$
called cumulant-generating function. The Taylor coefficients $k_n(\mu)$ of this series, given by
$$K_\mu(\xi)=\sum_{n=1}^\infty k_n(\mu)\,\frac{\xi^n}{n!}$$
are called cumulants of the measure $\mu$. We also use the notations $k_f,K_f$ for these cumulants and their generating series, where $f$ is a variable following the law $\mu$.
\end{definition}

In other words, the cumulants are more or less the coefficients of the logarithm of the Fourier transform $\log F_\mu$, up to some normalizations. To be more precise, we have $K_\mu(\xi)=\log F_\mu(-i\xi)$, so the formula relating $\log F_\mu$ to the cumulants $k_n(\mu)$ is:
$$\log F_\mu(-i\xi)=\sum_{n=1}^\infty k_n(\mu)\,\frac{\xi^n}{n!}$$

The interest in cumulants comes from the fact that $\log F_\mu$, and so the cumulants $k_n(\mu)$ too, linearize the convolution. To be more precise, we have the following result:

\index{independence}

\begin{theorem}
The cumulants have the following properties:
\begin{enumerate}
\item $k_n(cf)=c^nk_n(f)$.

\item $k_1(f+d)=k_1(f)+d$, and $k_n(f+d)=k_n(f)$ for $n>1$.

\item $k_n(f+g)=k_n(f)+k_n(g)$, if $f,g$ are independent.
\end{enumerate}
\end{theorem}

\begin{proof}
Here (1) and (2) are both clear from definitions, because we have the following computation, valid for any $c,d\in\mathbb R$, which gives the results:
\begin{eqnarray*}
K_{cf+d}(\xi)
&=&\log E(e^{\xi(cf+d)})\\
&=&\log[e^{\xi d}\cdot E(e^{\xi cf})]\\
&=&\xi d+K_f(c\xi)
\end{eqnarray*}

As for (3), this follows from the fact that the Fourier transform $F_f(\xi)=E(e^{i\xi f})$ satisfies the following formula, whenever $f,g$ are independent random variables:
$$F_{f+g}(\xi)=F_f(\xi)F_g(\xi)$$

Indeed, by applying the logarithm, we obtain the following formula:
$$\log F_{f+g}(\xi)=\log F_f(\xi)+\log F_g(\xi)$$

With the change of variables $\xi\to-i\xi$, we obtain the following formula:
$$K_{f+g}(\xi)=K_f(\xi)+K_g(\xi)$$

Thus, at the level of coefficients, we obtain $k_n(f+g)=k_n(f)+k_n(g)$, as claimed.
\end{proof}

At the level of the main examples now, we have the following result:

\index{Dirac mass}

\begin{theorem}
The sequence of cumulants $k_1,k_2,k_3,\ldots$ is as follows:
\begin{enumerate}
\item For $\mu=\delta_c$ the cumulants are $c,0,0,\ldots$

\item For $\mu=g_t$ the cumulants are $0,t,0,0,\ldots$

\item For $\mu=p_t$ the cumulants are $t,t,t,\ldots$

\item For $\mu=b_t$ the cumulants are $0,t,0,t,\ldots$

\item For $\mu=p_\nu$ the cumulants are the moments $M_n(\nu)$.
\end{enumerate}
\end{theorem}

\begin{proof}
We have several computations to be done, the idea being as follows:

\medskip

(1) For $\mu=\delta_c$ we have indeed the following computation:
$$K_\mu(\xi)
=\log E(e^{c\xi})
=\log(e^{c\xi})
=c\xi$$

But the plain coefficients of this series are the numbers $c,0,0,\ldots\,$, and so the Taylor coefficients of this series are these same numbers $c,0,0,\ldots\,$, as claimed.

\medskip

(2) For $\mu=g_t$ we have the following computation, which gives the result:
\begin{eqnarray*}
K_\mu(\xi)
&=&\log F_\mu(-i\xi)\\
&=&\log\exp\left[-t(-i\xi)^2/2\right]\\
&=&t\xi^2/2
\end{eqnarray*}

(3) For $\mu=p_t$ we have the following computation, which gives the result:
\begin{eqnarray*}
K_\mu(\xi)
&=&\log F_\mu(-i\xi)\\
&=&\log\exp\left[(e^{i(-i\xi)}-1)t\right]\\
&=&(e^\xi-1)t
\end{eqnarray*}

(4) For $\mu=b_t$ we have the following computation, which gives the result:
\begin{eqnarray*}
K_\mu(\xi)
&=&\log F_\mu(-i\xi)\\
&=&\log\exp\left[\left(\frac{e^\xi+e^{-\xi}}{2}-1\right)t\right]\\
&=&\left(\frac{e^\xi+e^{-\xi}}{2}-1\right)t
\end{eqnarray*}

(5) By using the Fourier transform formula for $p_\nu$ from chapter 13, we obtain:
\begin{eqnarray*}
K_{p_\nu}(\xi)
&=&\log\exp\left[\sum_it_i(e^{\xi z_i}-1)\right]\\
&=&\sum_{n\geq1}\frac{\xi^n}{n!}\sum_it_iz_i^n\\
&=&\sum_{n\geq1}\frac{\xi^n}{n!}\,M_n(\nu)
\end{eqnarray*}

Thus, we are led to the conclusion in the statement.
\end{proof}

Getting back now to theory, the sequence of cumulants $k_1,k_2,k_3,\ldots$ appears as a modification of the sequence of moments $M_1,M_2,M_3,\ldots\,$, and understanding the relation between moments and cumulants will be our next task. Let us start with:

\index{cumulant}
\index{generalized cumulant}

\begin{definition}
We define quantities $M_\pi(f),k_\pi(f)$, depending on partitions 
$$\pi\in P(k)$$
by starting with $M_n(f),k_n(f)$, and using multiplicativity over the blocks. 
\end{definition}

To be more precise, the convention here is that for the one-block partition $1_n\in P(n)$, the corresponding moment and cumulant are the usual ones, namely:
$$M_{1_n}(f)=M_n(f)\quad,\quad k_{1_n}(f)=k_n(f)$$

Then, for an arbitrary partition $\pi\in P(k)$, we decompose this partition into blocks, having sizes $b_1,\ldots,b_s$, and we set, by multiplicativity over blocks:
$$M_\pi(f)=M_{b_1}(f)\ldots M_{b_s}(f)\quad,\quad k_\pi(f)=k_{b_1}(f)\ldots k_{b_s}(f)$$

With this convention, following Rota, we can now formulate a key result, fully clarifying the relation between moments and cumulants, as follows:

\index{moment-cumulant formula}

\begin{theorem}
We have the moment-cumulant formulae
$$M_n(f)=\sum_{\nu\in P(n)}k_\nu(f)\quad,\quad 
k_n(f)=\sum_{\nu\in P(n)}\mu(\nu,1_n)M_\nu(f)$$
or, equivalently, we have the moment-cumulant formulae
$$M_\pi(f)=\sum_{\nu\leq\pi}k_\nu(f)\quad,\quad 
k_\pi(f)=\sum_{\nu\leq\pi}\mu(\nu,\pi)M_\nu(f)$$
where $\mu$ is the M\"obius function of $P(n)$.
\end{theorem}

\begin{proof}
There are several things going on here, the idea being as follows:

\medskip

(1) The first set of formulae is equivalent to the second set of formulae, by multiplicativity over blocks. Also, due to the M\"obius inversion formula, in the second set of formulae, the two formulae there are in fact equivalent. Thus, the 4 formulae in the statement are all equivalent. We will prove in what follows the first formula.

\medskip

(2) In order to do this, we use the very definition of the cumulants, namely:
$$\log E(e^{\xi f})=\sum_{s=1}^\infty k_s(f)\,\frac{\xi^s}{s!}$$

By exponentiating, we obtain from this the following formula:
$$E(e^{\xi f})=\exp\left(\sum_{s=1}^\infty k_s(f)\,\frac{\xi^s}{s!}\right)$$

(3) Let us first compute the function on the left. This is easily done, as follows:
$$E(e^{\xi f})
=E\left(\sum_{n=0}^\infty\frac{(\xi f)^n}{n!}\right)
=\sum_{n=0}^\infty M_n(f)\,\frac{\xi^n}{n!}$$

(4) Regarding now the function on the right, this is given by:
\begin{eqnarray*}
\exp\left(\sum_{s=1}^\infty k_s(f)\,\frac{\xi^s}{s!}\right)
&=&\sum_{p=0}^\infty\frac{\left(\sum_{s=1}^\infty k_s(f)\,\frac{\xi^s}{s!}\right)^p}{p!}\\
&=&\sum_{p=0}^\infty\frac{1}{p!}\sum_{s_1=1}^\infty k_{s_1}(f)\,\frac{\xi^{s_1}}{s_1!}\ldots\ldots\sum_{s_p=1}^\infty k_{s_p}(f)\,\frac{\xi^{s_p}}{s_p!}\\
&=&\sum_{p=0}^\infty\frac{1}{p!}\sum_{s_1=1}^\infty\ldots\sum_{s_p=1}^\infty k_{s_1}(f)\ldots k_{s_p}(f)\,\frac{\xi^{s_1+\ldots+s_p}}{s_1!\ldots s_p!}
\end{eqnarray*}

(5) The point now is that all this leads us into partitions. Indeed, we are summing over indices $s_1,\ldots,s_p\in\mathbb N$, which can be thought of as corresponding to a partition of $n=s_1+\ldots+s_p$. So, let us rewrite our sum, as a sum over partitions. For this purpose, recall that the number of partitions $\nu\in P(n)$ having blocks of sizes $s_1,\ldots,s_p$ is:
$$\binom{n}{s_1,\ldots,s_p}=\frac{n!}{p_1!\ldots p_s!}$$

Also, when resumming over partitions, there will be a $p!$ factor as well, coming from the permutations of $s_1,\ldots,s_p$. Thus, our sum can be rewritten as follows:
\begin{eqnarray*}
\exp\left(\sum_{s=1}^\infty k_s(f)\,\frac{\xi^s}{s!}\right)
&=&\sum_{n=0}^\infty\sum_{p=0}^\infty\frac{1}{p!}\sum_{s_1+\ldots+s_p=n}k_{s_1}(f)\ldots k_{s_p}(f)\,\frac{\xi^n}{s_1!\ldots s_p!}\\
&=&\sum_{n=0}^\infty\frac{\xi^n}{n!}\sum_{p=0}^\infty\frac{1}{p!}\sum_{s_1+\ldots+s_p=n}\binom{n}{s_1,\ldots,s_p}k_{s_1}(f)\ldots k_{s_p}(f)\\
&=&\sum_{n=0}^\infty\frac{\xi^n}{n!}\sum_{\nu\in P(n)}k_\nu(f)
\end{eqnarray*}

(6) We are now in position to conclude. According to (2,3,5), we have:
$$\sum_{n=0}^\infty M_n(f)\,\frac{\xi^n}{n!}=\sum_{n=0}^\infty\frac{\xi^n}{n!}\sum_{\nu\in P(n)}k_\nu(f)$$

Thus, we have the following formula, valid for any $n\in\mathbb N$:
$$M_n(f)=\sum_{\nu\in P(n)}k_\nu(f)$$

We are therefore led to the conclusions in the statement.
\end{proof}

Let us discuss now the free analogues of the above, following Speicher. We have:

\index{cumulant}
\index{free cumulant}
\index{Cauchy transform}

\begin{definition}
The free cumulants $\kappa_n(a)$ of a variable $a\in A$ are defined by:
$$R_a(\xi)=\sum_{n=1}^\infty\kappa_n(a)\xi^{n-1}$$
That is, the free cumulants are the coefficients of the $R$-transform, shifted by $1$.
\end{definition}

As before with the classical cumulants, we have a number of basic examples and illustrations, and a number of basic general results. Let us start with:

\begin{theorem}
The free cumulants have the following properties:
\begin{enumerate}
\item $\kappa_n(\lambda a)=\lambda^n\kappa_n(a)$.

\item $\kappa_n(a+b)=\kappa_n(a)+\kappa_n(b)$, if $a,b$ are free.
\end{enumerate}
\end{theorem}

\begin{proof}
In what regards (1), we have here the following computation:
\begin{eqnarray*}
G_{\lambda a}\left(\lambda R_a(\lambda\xi)+\frac{1}{\xi}\right)
&=&\frac{1}{\lambda}\,G_a\left(R_a(\lambda\xi)+\frac{1}{\lambda\xi}\right)\\
&=&\frac{1}{\lambda}\cdot\lambda\xi\\
&=&\xi
\end{eqnarray*}

Thus $R_{\lambda a}(\xi)=\lambda R_a(\lambda\xi)$, which gives (1). As for (2), this follows from the fact, that we know well, that the $R$-transform linearizes the free convolution operation. 
\end{proof}

Again in analogy with the classical case, at the level of examples, we have:

\index{Dirac mass}

\begin{theorem}
The sequence of free cumulants $\kappa_1,\kappa_2,\kappa_3,\ldots$ is as follows:
\begin{enumerate}
\item For $\mu=\delta_c$ the free cumulants are $c,0,0,\ldots$

\item For $\mu=\gamma_t$ the free cumulants are $0,t,0,0,\ldots$

\item For $\mu=\pi_t$ the free cumulants are $t,t,t,\ldots$

\item For $\mu=\beta_t$ the free cumulants are $0,t,0,t,\ldots$

\item For $\mu=\pi_\nu$ the free cumulants are the moments $M_n(\nu)$.
\end{enumerate}
\end{theorem}

\begin{proof}
The proofs are analogous to those from the classical case, as follows:

\medskip

(1) For $\mu=\delta_c$ we have $G_\mu(\xi)=1/(\xi-c)$, and so $R_\mu(\xi)=c$, as desired.

\medskip

(2) For $\mu=\gamma_t$ we have, as computed before, $R_\mu(\xi)=t\xi$, as desired.

\medskip

(3) For $\mu=\pi_t$ we have, also from before, $R_\mu(\xi)=t/(1-\xi)$, as desired.

\medskip

(4) For $\mu=\beta_t$ we have $R_\mu(\xi)=t/(1-\xi^2)$, which gives the result. 

\medskip

(5) With $\nu=\sum_ic_i\delta_{z_i}$, we have the following computation:
\begin{eqnarray*}
R_{\pi_\nu}(\xi)
&=&\sum_i\frac{c_iz_i}{1-\xi z_i}\\
&=&\sum_{n\geq1}\xi^{n-1}\sum_ic_iz_i^n\\
&=&\sum_{n\geq 1}\xi^{n-1}\,M_n(\nu)
\end{eqnarray*}

Thus, we are led to the conclusion in the statement.
\end{proof}

Also as before in the classical case, we can define generalized free cumulants $\kappa_\pi(a)$ with $\pi\in P(k)$ by starting with the numeric free cumulants $\kappa_n(a)$, as follows:

\begin{definition}
We define free cumulants $\kappa_\pi(a)$, depending on partitions 
$$\pi\in P(k)$$
by starting with $\kappa_n(a)$, and using multiplicativity over the blocks. 
\end{definition}

To be more precise, the convention here is that for the one-block partition $1_n\in P(n)$, the corresponding free cumulant is the usual one, namely:
$$\kappa_{1_n}(a)=\kappa_n(a)$$

Then, for an arbitrary partition $\pi\in P(k)$, we decompose this partition into blocks, having sizes $b_1,\ldots,b_s$, and we set, by multiplicativity over blocks:
$$\kappa_\pi(a)=\kappa_{b_1}(a)\ldots\kappa_{b_s}(a)$$

With this convention, we have the following result, due to Speicher:

\index{moment-cumulant formula}

\begin{theorem}
We have the moment-cumulant formulae
$$M_n(a)=\sum_{\nu\in NC(n)}\kappa_\nu(a)\quad,\quad 
\kappa_n(a)=\sum_{\nu\in NC(n)}\mu(\nu,1_n)M_\nu(a)$$
or, equivalently, we have the moment-cumulant formulae
$$M_\pi(a)=\sum_{\nu\leq\pi}\kappa_\nu(a)\quad,\quad 
\kappa_\pi(a)=\sum_{\nu\leq\pi}\mu(\nu,\pi)M_\nu(a)$$
where $\mu$ is the M\"obius function of $NC(n)$.
\end{theorem}

\begin{proof}
As before in the classical case, the 4 formulae in the statement are equivalent, via M\"obius inversion. Thus, it is enough to prove one of them. We have:
$$G_a(\xi)=\sum_{n=0}^\infty\frac{M_n(a)}{\xi^{n+1}}$$

Consider the inverse of this Cauchy transform $G_a$, with respect to composition:
$$G_a(K_a(\xi))=K_a(G_a(\xi))=\xi$$

According to Definition 15.20, the free cumulants $\kappa_n(a)$ appear then as follows:
$$K_a(\xi)=\frac{1}{\xi}+\sum_{n=1}^\infty\kappa_n(a)\xi^{n-1}$$

Thus, we can compute moments in terms of free cumulants by using either of the inversion formulae $G_a(K_a(\xi))=\xi$ and $K_a(G_a(\xi))=\xi$. In practice, the best is to use the second inversion formula, $K_a(G_a(\xi))=\xi$, which after some manipulations reads:
\begin{eqnarray*}
&&\kappa_1z+\kappa_2z^2(1+M_1z+M_2z^2+\ldots)+\kappa_3z^3(1+M_1z+M_2z^2+\ldots)^2+\ldots\\
&=&(M_1z+M_2z^2+\ldots)-(M_1z+M_2z^2+\ldots)^2+(M_1z+M_2z^2+\ldots)^3-\ldots
\end{eqnarray*}

But this gives, and after some computations, the first formula, namely:
$$M_n(a)=\sum_{\nu\in NC(n)}\kappa_\nu(a)$$

For details here, we refer for instance to the book of Nica-Speicher \cite{nsp}.
\end{proof}

With the above theory in hand, we can now formulate, following Bercovici-Pata \cite{bpa}, the following simple definition, making the connection between classical and free:

\index{convolution semigroup}
\index{Bercovici-Pata bijection}

\begin{definition}
A convolution semigroup of measures
$$\{m_t\}_{t>0}\quad:\quad m_s*m_t=m_{s+t}$$
is in Bercovici-Pata bijection with a free convolution semigroup of measures
$$\{\mu_t\}_{t>0}\quad:\quad \mu_s\boxplus\mu_t=\mu_{s+t}$$
when the classical cumulants of $m_t$ coincide with the free cumulants of $\mu_t$.
\end{definition}

As basic examples for this correspondence, coming from our various explicit computations of classical and free cumulants in this chapter, we have:

\begin{theorem}
The semigroups formed by the normal, Poisson and Bessel laws
$$\{g_t\}_{t>0}\quad,\quad\{G_t\}_{t>0}\quad,\quad\{p_t\}_{t>0}\quad,\quad\{b_t\}_{t>0}$$
are in Bercovici-Pata bijection with the free convolution semigroups
$$\{\gamma_t\}_{t>0}\quad,\quad\{\Gamma_t\}_{t>0}\quad,\quad\{\pi_t\}_{t>0}\quad,\quad\{\beta_t\}_{t>0}$$
formed by the Wigner, Voiculescu, Marchenko-Pastur and free Bessel laws.
\end{theorem}

\begin{proof}
This follows indeed from our various cumulant computations above, from Theorem 15.17 and Theorem 15.22, leading to the same sequences of cumulants.
\end{proof}

Many other things can be said about the Bercovici-Pata bijection, notably with an extension to the case of complex measures, and with a number of analytic features too. For more on all this, you can have a look at \cite{bpa}, \cite{nsp} and related literature.

\section*{15c. Rotations, reflections}

Getting back now to our usual business, namely groups and characters, we first have the following result, regarding the usual characters of the main free quantum groups:

\index{Marchenko-Pastur law}
\index{Bessel law}
\index{real Bessel law}
\index{complex Bessel law}
\index{Wigner law}
\index{Voiculescu law}
\index{semicircle law}
\index{circular law}

\begin{theorem}
The laws of main characters are as follows:
\begin{enumerate}
\item For $S_N^+$ with $N\geq4$, we obtain the Marchenko-Pastur law $\pi_1$.

\item For $H_N^+$ with $N\geq4$, we obtain the free real Bessel law $\beta_1$.

\item For $K_N^+$ with $N\geq4$, we obtain the free complex Bessel law $\mathfrak B_1$.

\item For $H_N^{s+}$ with $N\geq4$, we obtain the general free Bessel law $\beta^s_1$.

\item For $O_N^+$ with $N\geq2$, we obtain the Wigner semicircle law $\gamma_1$.

\item For $U_N^+$ with $N\geq2$, we obtain the Voiculescu circular law $\Gamma_1$.

\item For $B_N^+$ with $N\geq3$, we obtain the shifted Wigner law $\gamma_1'$.

\item For $C_N^+$ with $N\geq3$, we obtain the shifted Voiculescu law $\Gamma_1'$.
\end{enumerate}
\end{theorem}

\begin{proof}
As a first observation, this is something quite remarkable, much easier than in the classical case, due to the fact that the results are all valid for $N$ big enough, not just in the $N\to\infty$ limit. As for the proof, the idea is as follows:

\medskip

(1) To start with, as before in the classical case, for an easy quantum group $G_N\subset U_N$, coming from a category of partitions $D\subset P$, we have:
$$\lim_{N\to\infty}\int_{G_N}\chi^k=|D(k)|$$

But this gives all the results in the statement, in the $N\to\infty$ limit, via the various moment formulae established throughout this chapter, for the above measures.

\medskip

(2) However, the result says much more than this, being not asymptotic. In order to understand this, a first way is by looking at the basic moment formula for the easy quantum groups, which is similar to the one for the classical easy groups, namely:
$$\int_{G_N}\chi^k=\dim\left(span\left(\xi_\pi\Big|\pi\in D(k)\right)\right)$$

Thus, in order to prove the result as stated, we must show that for each of the categories of partitions $D\subset P$ involved, the above families of vectors $\{\xi_\pi|\pi\in D(k)\}$ are linearly independent, provided that $N\in\mathbb N$ is subject to the inequalities in the statement.

\medskip

(3) But here, the point is that we have in all cases $D\subset NC$, our quantum groups being free, and we can therefore use the meander determinant techniques of Di Francesco, explained in chapter 6. To be more precise, the Di Francesco results from chapter 6 solve the linear independence problem for $O_N^+$ with $N\geq2$, and for $S_N^+$ with $N\geq 4$ too, and the extension to the other quantum groups is straightforward, by using:
$$H_N^{s+}=\mathbb Z_s\wr_*S_N^+\quad,\quad U_N^+=\widetilde{O_N^+}\quad,\quad 
B_N^+\simeq O_{N-1}^+\quad,\quad C_N^+\simeq U_{N-1}^+$$

(4) However, all this remains excessively heavy, and in practice, it is more convenient to do a case-by-case study, based on the basic Peter-Weyl formula, namely:
$$\int_{G_N}\chi^k=\dim(Fix(u^{\otimes k}))$$

(5) Indeed, regarding $S_N^+$, we know from chapter 12 that the fusion rules are the Clebsch-Gordan ones for $SO_3$, independently of $N\geq4$. Thus the above moments are the same at any $N\geq4$, and so they coincide with the asymptotic moments, as desired.

\medskip

(6) Regarding now the quantum reflection groups $H_N^s=\mathbb Z_s\wr_*S_N^+$, these contain the quantum permutation groups $S_N^+$, so by functoriality the spaces of invariants are smaller, and we therefore have linear independence at $N\geq4$ too, as desired.

\medskip

(7) In what regards now $O_N^+$, we know from chapter 12 that the fusion rules are the Clebsch-Gordan ones for $SU_2$, independently of $N\geq2$. Thus the moments are the same at any $N\geq2$, and so they coincide with the asymptotic moments, as desired.

\medskip

(8) Regarding now $U_N^+$, we can get away here by using functoriality and the result for $O_N^+$, or via free complexification, or via the fusion rules from chapter 12.

\medskip

(9) Finally, regarding the bistochastic quantum groups, we can use here the standard isomorphisms $B_N^+\simeq O_{N-1}^+$ and $C_N^+\simeq U_{N-1}^+$, whose proof is similar to the corresponding isomorphisms in the classical case, and we obtain the results, via $O_N^+$ and $U_N^+$.
\end{proof}

Regarding now the corresponding truncated characters, the result here is more modest, namely asymptotic, similar to the one in the classical case, as follows:

\begin{theorem}
The asymptotic laws of truncated characters are as follows:
\begin{enumerate}
\item For $S_N^+$, we obtain the Marchenko-Pastur law $\pi_t$.

\item For $H_N^+$, we obtain the free real Bessel law $\beta_t$.

\item For $K_N^+$, we obtain the free complex Bessel law $\mathfrak B_t$.

\item For $H_N^{s+}$, we obtain the general free Bessel law $\beta^s_t$.

\item For $O_N^+$, we obtain the Wigner semicircle law $\gamma_t$.

\item For $U_N^+$, we obtain the Voiculescu circular law $\Gamma_t$.

\item For $B_N^+$, we obtain the shifted Wigner law $\gamma_t'$.

\item For $C_N^+$, we obtain the shifted Voiculescu law $\Gamma_t'$.
\end{enumerate}
\end{theorem}

\begin{proof}
Given an easy quantum group $G\subset_uU_N^+$, coming from a category of partitions $D\subset P$, we have the following Weingarten formula, with $W_{kN}=G_{kN}^{-1}$:
$$\int_Gu_{i_1j_1}^{e_1}\ldots u_{i_kj_k}^{e_k}=\sum_{\pi,\sigma\in D(k)}\delta_\pi(i)\delta_\sigma(j)W_{kN}(\pi,\sigma)$$

Then, we have the following computation, again as in the classical case:
\begin{eqnarray*}
\int_G\chi_t^k
&=&\sum_{i_1=1}^{[tN]}\ldots\sum_{i_k=1}^{[tN]}\int u_{i_1i_1}\ldots u_{i_ki_k}\\
&=&\sum_{\pi,\sigma\in D(k)}W_{kN}(\pi,\sigma)\sum_{i_1=1}^{[tN]}\ldots\sum_{i_k=1}^{[tN]}\delta_\pi(i)\delta_\sigma(i)\\
&=&\sum_{\pi,\sigma\in D(k)}W_{kN}(\pi,\sigma)G_{k[tN]}(\sigma,\pi)\\
&=&Tr(W_{kN}G_{k[tN]})
\end{eqnarray*}

In the uniform case $G_{kN}$ and $W_{kN}$ are asymptotically diagonal, which gives:
$$\int_{G_N}\chi_t^k
\simeq\sum_{\pi\in D(k)}N^{-|\pi|}[tN]^{|\pi|}
\simeq\sum_{\pi\in D(k)}t^{|\pi|}$$

But this gives the results, via our various easiness and moment formulae.
\end{proof}

As a consequence of the above results, we can draw some nice cubes, as follows:

\index{standard cube}
\index{truncated character}

\begin{theorem}
The main quantum rotation and reflection groups
$$
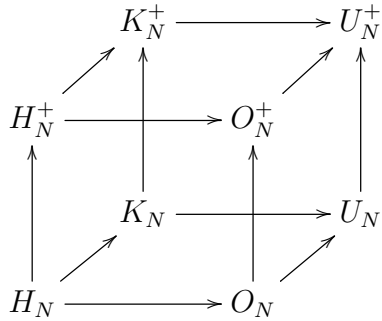
\xymatrix@R=17pt@C=17pt{
&K_N^+\ar[rr]&&U_N^+\\
H_N^+\ar[rr]\ar[ur]&&O_N^+\ar[ur]\\
&K_N\ar[rr]\ar[uu]&&U_N\ar[uu]\\
H_N\ar[uu]\ar[ur]\ar[rr]&&O_N\ar[uu]\ar[ur]
}$$
are all easy, coming from the following categories of partitions,
$$
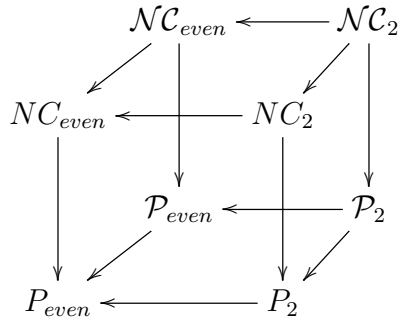
\xymatrix@R=18.5pt@C=2pt{
&\mathcal{NC}_{even}\ar[dl]\ar[dd]&&\mathcal {NC}_2\ar[dl]\ar[ll]\ar[dd]\\
NC_{even}\ar[dd]&&NC_2\ar[dd]\ar[ll]\\
&\mathcal P_{even}\ar[dl]&&\mathcal P_2\ar[dl]\ar[ll]\\
P_{even}&&P_2\ar[ll]
}$$
and the corresponding truncated characters follow with $N\to\infty$ the laws
$$
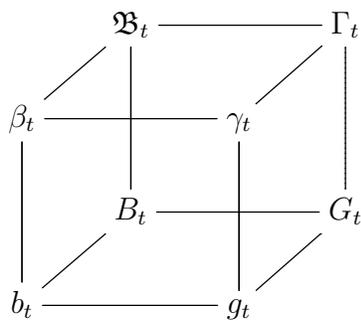
\xymatrix@R=18pt@C=22pt{
&\mathfrak B_t\ar@{-}[rr]\ar@{-}[dd]&&\Gamma_t\ar@{-}[dd]\\
\beta_t\ar@{-}[rr]\ar@{-}[dd]\ar@{-}[ur]&&\gamma_t\ar@{-}[dd]\ar@{-}[ur]\\
&B_t\ar@{-}[rr]\ar@{-}[uu]&&G_t\ar@{.}[uu]\\
b_t\ar@{-}[uu]\ar@{-}[ur]\ar@{-}[rr]&&g_t\ar@{-}[uu]\ar@{-}[ur]
}$$
which are the main limiting laws in classical and free probability.
\end{theorem}

\begin{proof}
This is indeed something self-explanatory, coming from the various results above, regarding classical and free probability, and classical and free quantum groups.
\end{proof}

The above cubes are something quite interesting, providing us with some 3D orientation, inside quantum algebra and probability, at large. Indeed, the $Ox$ axis corresponds to the making continuous/discretizing operation, the $Oy$ axis corresponds to complexification/taking the real version, and the $Oz$ axis corresponds to liberating/taking the classical version. Which is something remarkable, and very useful in practice.

\bigskip

As a technical comment now, all this comes with a price, because we are restricting in this way the attention to the intermediate quantum groups $H_N\subset G\subset U_N^+$, which can fit into the cube, therefore ignoring the quantum groups $S_N,S_N^+$, and more generally $H_N^s,H_N^{s+}$ with $s$ odd, as well as the bistochastic groups $B_N,B_N^+$ and $C_N,C_N^+$.

\bigskip

Regarding these forgotten quantum groups, you might probably think that these can be included too, by adding a 4th dimension. However, this is not really possible:

\begin{fact}
The classical and free uniform orthogonal easy quantum groups
$$
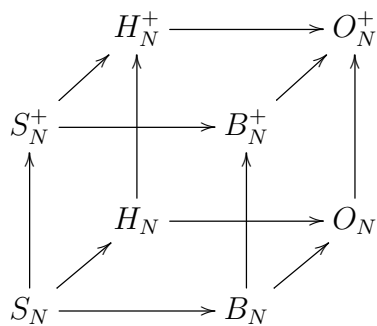
\xymatrix@R=17pt@C=17pt{
&H_N^+\ar[rr]&&O_N^+\\
S_N^+\ar[rr]\ar[ur]&&B_N^+\ar[ur]\\
&H_N\ar[rr]\ar[uu]&&O_N\ar[uu]\\
S_N\ar[uu]\ar[ur]\ar[rr]&&B_N\ar[uu]\ar[ur]
}$$
and the unitary analogues of these quantum groups, namely
$$
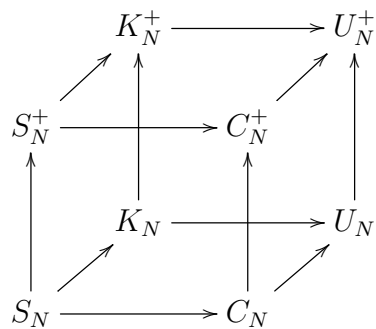
\xymatrix@R=17.5pt@C=17pt{
&K_N^+\ar[rr]&&U_N^+\\
S_N^+\ar[rr]\ar[ur]&&C_N^+\ar[ur]\\
&K_N\ar[rr]\ar[uu]&&U_N\ar[uu]\\
S_N\ar[uu]\ar[ur]\ar[rr]&&C_N\ar[uu]\ar[ur]
}$$
cannot be merged, without degeneration, into a $4$-dimensional cubic diagram.
\end{fact}

And we will leave some thinking here as an instructive exercise. So, as a conclusion, the 3D cubes in Theorem 15.29 remain the main tool, but for certain more technical questions, in relation with bistochasticity, some further 3D cubes are available too.

\section*{15d. Gram determinants}

Getting back now to our original Question 15.1, we have seen that in what regards the characters, truncated or not, the liberation operation simplifies things, and generally speaking, confirms our general principle, that quantum is simpler than classical.

\bigskip

Passed characters, many other probabilistic things can be done, such as investigating Diaconis-Shahshahani variables, invariance questions with De Finetti type theorems, all sorts of random walk questions, idempotent states following It\^o, classical and free hypespherical laws, classical and free hypergeometric laws, classical and free Meixner laws, and many more. For an introduction to this, you can check my book \cite{ba2}.

\bigskip

Along these lines, we would like to discuss now, in the remainder of this chapter, the Gram determinant question, which is perhaps central, in all this. By restricting the attention to the orthogonal case, the question that we would like to solve is as follows:

\begin{question}
In the context of the main orthogonal easy quantum groups
$$
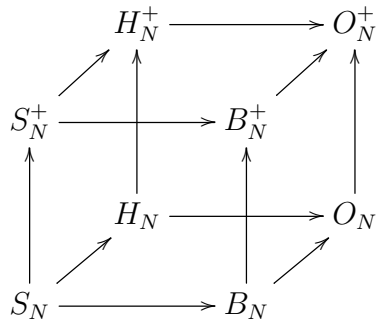
\xymatrix@R=17pt@C=17pt{
&H_N^+\ar[rr]&&O_N^+\\
S_N^+\ar[rr]\ar[ur]&&B_N^+\ar[ur]\\
&H_N\ar[rr]\ar[uu]&&O_N\ar[uu]\\
S_N\ar[uu]\ar[ur]\ar[rr]&&B_N\ar[uu]\ar[ur]
}$$
how do the Gram determinants behave, on the vertical?
\end{question}

And good question this is, because passed all sorts of algebraic manipulations, that we have been doing since long, involving representations, characters, laws, moments, cumulants and so on, it all comes down to the linear independence of the vectors $\xi_\pi$, and so, essentially, to the computation of the Gram determinant of these vectors.

\bigskip 

In answer now, we first need to know how to compute these determinants. In the classical discrete case, the answer is something quite simple, as follows:

\begin{theorem}
For the groups $G=S_N,H_N$ we have the Lindst\"om formula
$$\det(G_{kN})=\prod_{\pi\in D(k)}\frac{N!}{(N-|\pi|)!}$$
with $D=P,P_{even}$, and with $|.|$ being the number of blocks.
\end{theorem}

\begin{proof}
This is something that we already know for $G=S_N$, from chapter 14, and the proof for $G=H_N$ is similar, based in the fact that the corresponding categories of partitions $D=P,P_{even}$, have the property of forming semilattices under $\vee$. Consider indeed the following matrix, obtained by making determinant-preserving operations:
$$G_{kN}'(\pi,\sigma)=\sum_{\pi\leq\tau}\mu(\pi,\tau)N^{|\tau\vee\sigma|}$$

It follows then from the M\"obius inversion formula that we have:
$$G_{kN}'(\pi,\sigma)=
\begin{cases}
N(N-1)\ldots(N-|\sigma|+1)&{\rm if}\ \pi\leq\sigma\\
0&{\rm otherwise}
\end{cases}$$

Thus the matrix is upper triangular, and by computing the product on the diagonal we obtain the formula in the statement, exactly as in the $G=S_N$ case before.
\end{proof}

Next, let us discuss the case of the orthogonal group $O_N$. Here things are more complicated, the combinatorics being that of the Young diagrams. In relation with these, let us denote by $|.|$ the number of boxes, and use quantity $f^\lambda$, which gives the number of standard Young tableaux of shape $\lambda$. With these conventions, we have:

\index{Young tableaux}

\begin{theorem}
The determinant of the Gram matrix of $O_N$ is given by
$$\det(G_{kN})=\prod_{|\lambda|=k/2}f_N(\lambda)^{f^{2\lambda}}$$
where the quantities on the right are $f_N(\lambda)=\prod_{(i,j)\in\lambda}(N+2j-i-1)$.
\end{theorem}

\begin{proof}
For the group $O_N$ the Gram matrix is diagonalizable, as follows:
$$G_{kN}=\sum_{|\lambda|=k/2}f_N(\lambda)P_{2\lambda}$$

To be more precise, here $1=\sum P_{2\lambda}$ is the standard partition of unity associated to the Young diagrams having $k/2$ boxes, and the coefficients $f_N(\lambda)$ are those in the statement. Now since we have $Tr(P_{2\lambda})=f^{2\lambda}$, this gives the formula in the statement.
\end{proof}

As a variation of the above result, for the bistochastic group $B_N$, we have:

\index{bistochastic group}

\begin{theorem}
For the bistochastic group $B_N$ we have
$$\det(G_{kN})=N^{a_k}\prod_{|\lambda|\leq k/2}f_N(\lambda)^{\binom{k}{2|\lambda|}f^{2\lambda}}$$
where $a_k=\sum_{\pi\in P_{12}(k)}(2|\pi|-k)$, and $f_N(\lambda)=\prod_{(i,j)\in\lambda}(N+2j-i-2)$.
\end{theorem}

\begin{proof}
We recall that we have an isomorphism $B_N\simeq O_{N-1}$, given by $u=v+1$, where $u,v$ are the fundamental representations of $B_N,O_{N-1}$. But this gives:
\begin{eqnarray*}
Fix(u^{\otimes k})
&=&Fix\left((v+1)^{\otimes k}\right)\\
&=&Fix\left(\sum_{r=0}^k\binom{k}{r}v^{\otimes r}\right)
\end{eqnarray*}

Now if we denote by $\det',f'$ the objects in Theorem 15.33, we obtain:
\begin{eqnarray*}
\det(G_{kN})
&=&N^{a_k}\prod_{r=1}^k{\rm det}'(G_{r,N-1})^{\binom{k}{r}}\\
&=&N^{a_k}\prod_{r=1}^k\left(\prod_{|\lambda|=r/2}f'_{N-1}(\lambda)^{f^{2\lambda}}\right)^{\binom{k}{r}}
\end{eqnarray*}

This gives the formula in the statement.
\end{proof}

Summarizing, we have computed the Gram determinants for the main easy groups $G\subset O_N$. We should mention that the extension of these results to the case $G\subset U_N$ is quite straightforward, with of course colored integers needed for indexing the various Gram determinants, and we will leave some work here as an instructive exercise.

\bigskip

Getting now to the free case, the simplest object here is the free orthogonal group $O_N^+$, coming from the category of pairings $NC_2$. The associated Gram determinant, that we already met in chapter 6 in relation with questions in topology, is known as ``meander determinant'', and following Di Francesco \cite{dif}, we have the following result:

\index{meander determinant}
\index{Gram determinant}
\index{Di Francesco formula}

\begin{theorem}
The determinant of the Gram matrix for $O_N^+$ is given by
$$\det(G_{kN})=\prod_{r=1}^{[k/2]}P_r(N)^{d_{k/2,r}}$$
where $P_r$ are the Chebycheff polynomials, given by
$$P_0=1\quad,\quad 
P_1=X\quad,\quad 
P_{r+1}=XP_r-P_{r-1}$$
and $d_{kr}=f_{kr}-f_{k,r+1}$, with $f_{kr}$ being the following numbers, depending on $k,r\in\mathbb Z$,
$$f_{kr}=\binom{2k}{k-r}-\binom{2k}{k-r-1}$$
with the convention $f_{kr}=0$ for $k\notin\mathbb Z$. 
\end{theorem}

\begin{proof}
This is indeed something quite standard, but long and technical, that we explained a bit in chapter 6, in relation with our considerations there. For more on this, we refer to \cite{dif}, and to the more recent paper \cite{bcu} too, containing a short proof of this formula, based however on the advanced planar technology of Jones \cite{jo6}.
\end{proof}

Regarding now the quantum group $S_N^+$, coming from the category of partitions $NC$, we have here the following formula, also established by Di Francesco in \cite{dif}:

\index{meander determinant}
\index{Gram determinant}
\index{Chebycheff polynomials}
\index{Di Francesco formula}

\begin{theorem}
The determinant of the Gram matrix for $S_N^+$ is given by
$$\det(G_{kN})=(\sqrt{N})^{a_k}\prod_{r=1}^kP_r(\sqrt{N})^{d_{kr}}$$
where $P_r$ are as before the Chebycheff polynomials, given by
$$P_0=1\quad,\quad 
P_1=X\quad,\quad 
P_{r+1}=XP_r-P_{r-1}$$
then $d_{kr}=f_{kr}-f_{k,r+1}$, with $f_{kr}$ being the following numbers, depending on $k,r\in\mathbb Z$,
$$f_{kr}=\binom{2k}{k-r}-\binom{2k}{k-r-1}$$
with the convention $f_{kr}=0$ for $k\notin\mathbb Z$, and where $a_k=\sum_{\pi\in P(k)}(2|\pi|-k)$.
\end{theorem}

\begin{proof}
As explained in chapter 6, this comes from Theorem 15.35, via the following standard bijection, obtained via fattening and shrinking:
$$NC(k)\simeq NC_2(2k)$$

Indeed, if we denote by $G'$ the Gram matrix for $O_N^+$, we have the following formula, coming from the above bijection, with $D_{kN}=diag(N^{|\widetilde\pi|/2-k/4})$:
$$G_{kN}=D_{kN}G'_{2k,\sqrt{N}}D_{kN}$$

But with this formula in hand, the result follows from Theorem 15.35.
\end{proof}

As yet another version of Theorem 15.35, we can formulate as well:

\begin{theorem}
For the quantum group $B_N^+$ we have
$$\det(G_{kN})=N^{a_k}\prod_{r=1}^{[k/2]}P_r(N-1)^{\sum_{l=1}^{[k/2]}\binom{k}{2l}d_{lr}}$$
with $P_r$ and $d_{kr}$ being as before, and with $a_k=\sum_{\pi\in NC_{12}(k)}(2|\pi|-k)$.
\end{theorem}

\begin{proof}
The passage $O_N^+\to B_N^+$ is quite similar to the passage $O_N\to B_N$, from the proof of Theorem 15.34. To be more precise, we have an isomorphism $B_N^+\simeq O_{N-1}^+$, and by using prime exponents for the various $O_N^+$-related objects, we obtain:
\begin{eqnarray*}
\det(G_{kN})
&=&N^{a_k}\prod_{l=1}^{[k/2]}{\rm det}'(G_{2l,N-1})^{\binom{k}{2l}}\\
&=&N^{a_k}\prod_{l=1}^{[k/2]}\left(\prod_{r=1}^lP_r(N-1)^{d_{lr}}\right)^{\binom{k}{2l}}
\end{eqnarray*}

Together with Theorem 15.35, this gives the formula in the statement.
\end{proof}

Finally, in what regards the quantum group $H_N^+$, the result here is as follows:

\begin{theorem}
For the quantum group $H_N^+$ we have the formula
$$\det(G_{kN})=(\sqrt{N})^{a_k}\prod_{r=1}^{[k/2]}P_r(\sqrt{N})^{2d_{k/2,r}'}$$
with $d_{sr}'=f_{sr}'-f_{s,r+1}'$, where 
$f_{sr}'=\binom{3s}{s-r}-\binom{3s}{s-r-1}$ for $s\in\mathbb Z$, $f_{sr}'=0$ for $s\notin\mathbb Z$.
\end{theorem}

\begin{proof}
This can be viewed as yet another variation of Theorem 15.35, involving this time some colors on the strings, and for details here we refer to \cite{dif} and \cite{bcu}.
\end{proof}

Getting now to our original Question 15.31, we have formulae for all determinants there, and the remaining puzzle is that of putting everything together, with some conceptual results and formulae, regarding all this. For a discussion here, we refer to  \cite{bcu}.

\section*{15e. Exercises}

This was an exciting combinatorial chapter, and as exercises, we have:

\begin{exercise}
Clarify if needed our various results regarding characters.
\end{exercise}

\begin{exercise}
Learn about classical and free Diaconis-Shahshahani variables.
\end{exercise}

\begin{exercise}
Learn about classical and free De Finetti theorems.
\end{exercise}

\begin{exercise}
Learn about classical and free random walks.
\end{exercise}

\begin{exercise}
Learn about classical and free idempotent states.
\end{exercise}

\begin{exercise}
Learn about classical and free hyperspherical laws.
\end{exercise}

\begin{exercise}
Learn about classical and free hypergeometric laws.
\end{exercise}

\begin{exercise}
Learn also about classical and free Meixner laws.
\end{exercise}

As bonus exercise, try finding a global formula, for the various Gram determinants.

\chapter{Planar algebras}

\section*{16a. Operator algebras}

Time to end this book, and with the finite groups $G$ discussed, and their liberations $G^+$ discussed too, and with this quite often getting us into some other natural generalizations of the finite groups, namely the compact groups $G$, which can be classical and quantum, and the discrete groups too, again classical and quantum, $\Gamma=\widehat{G}$, you might wonder what can we talk about now, on top of this seemingly complete discussion.

\bigskip

In answer, subfactors and planar algebras in the sense of Jones \cite{jo1}, \cite{jo6}, the point being that this generalizes everything that we know, and much more, according to: 

\begin{principle}
With a suitable definition for the operator algebras $A\subset B(H)$, these decompose into factors, $A=\int_XA_x$, and an inclusion of such factors, called subfactor
$$A\subset B$$
can be thought of as coming from an action of a quantum group type object $G$, which widely generalizes all groups and quantum groups that we know, and that we don't.
\end{principle}

To be more precise here, past some operator algebra work, to be explained in a moment, the first part holds indeed, and in what regards the examples, these generalize indeed all sorts of quantum algebra beasts that can be constructed, including Hopf algebras, tensor categories, and various type of quantum groups that are available, provided that these have built-in a minimal amount of positivity, which is usually the case.

\bigskip

Importantly, this technology is not just for the sake of generalizing and unifying, as to have a very general theory, ready to welcome all sorts of beasts appearing from advanced theoretical physics, because there are concrete applications of this to the finite groups $G$ themselves. To be more precise, as a complement to the above principle, we have:

\begin{principle}[continuation]
In practice, the combinatorics of $A\subset B$ is described by the associated planar algebra $P=(P_k)$, with the simplest such algebras being
$$P_k=End(u^{\otimes k})$$
with $u$ being the fundamental representation of a group $G$, finite or compact or quantum. And with this making the link with Peter-Weyl, Tannaka, Brauer and easiness.
\end{principle}

But probably enough advertisement, let us get to work. We will need some basic von Neumann algebra theory, coming as a complement, or rather refinement, of the basic $C^*$-algebra theory developed in chapter 12. Let us start with:

\begin{theorem}
For a $*$-algebra of operators $A\subset B(H)$ the following conditions are equivalent, and if satisfied, we say that $A$ is a von Neumann algebra:
\begin{enumerate}
\item $A$ is closed under the weak topology, making each $T\to Tx$ continuous.

\item $A$ equals its bicommutant, $A=A''$, computed inside $B(H)$.
\end{enumerate}
\end{theorem}

\begin{proof}
This is von Neumann's bicommutant theorem, that we actually invoked a few times already, when talking Tannakian duality, in its finite dimensional particular case, which is elementary. In general the discussion, and then the proof, are as follows:

\medskip

(1) As a first comment, the weak topology on $B(H)$, making each $T\to Tx$ with $x\in H$ continuous, is indeed weaker than the norm topology, in the sense that we have:
$$||T_n-T||\to0\implies ||T_nx-Tx||\to0,\forall x\in H$$

In particular, we see that a von Neumann algebra in the sense of (1), that is, closed under the weak topology, must be a $C^*$-algebra, that is, closed under the norm.

\medskip

(2) Before getting further, let us see if the converse of this fact is true. This is certainly true in finite dimensions, $H=\mathbb C^N$, where we have $B(H)=M_N(\mathbb C)$, and where the operator $*$-algebras $A\subset B(H)$ are as follows, automatically closed both for the norm topology, and the weak topology, and with these two topologies actually coinciding:
$$A=M_{n_1}(\mathbb C)\oplus\ldots\oplus M_{n_k}(\mathbb C)$$

(3) In infinite dimensions, however, things change. Indeed, let us first take a look at the most basic examples of commutative $C^*$-algebras that we know, the commutative ones. These naturally appear from compact measured spaces $X$, as follows:
$$C(X)\subset B(L^2(X))\quad,\quad f\to[g\to fg]$$

(4) But, it is pretty much clear that such an algebra will not be weakly closed, unless $X$ is discrete, with the details here being left to you. So, in infinite dimensions, there are far less von Neumann algebras than $C^*$-algebras, with this being good to know.

\medskip

(5) Still talking about this, the following natural question appears, what happens if we take the weak closure of the algebra $C(X)\subset B(L^2(X))$ considered above? And the answer here, obtained via some basic measure theory and functional analysis, that we will leave as an exercise, is that we obtain the following algebra:
$$L^\infty(X)\subset B(L^2(X))\quad,\quad f\to[g\to fg]$$

(6) But this is quite interesting, because forgetting now about $C^*$-algebras, what we have here is a nice method of producing von Neumann algebras, in the weakly closed sense, and with the measured space $X$ being no longer required to be compact.

\medskip

(7) As a conclusion to all this, ``von Neumann algebras have something to do with measure theory, in the same way as $C^*$-algebras have something to do with topology''. Which sounds quite deep, so good, and time to stop here. More on this later.

\medskip

(8) Hang on, we are not done yet with the preliminaries, because all the above was in relation with the condition (1) in the statement, and we still have the condition (2) in the statement to comment on. So, here we go again, with a basic exploration, of that condition. To start with, given a subalgebra $A\subset B(H)$, or even a subset $A\subset B(H)$, we can talk about its commutant inside $B(H)$, constructed as follows:
$$A'=\left\{T\in B(H)\Big|TS=ST,\forall S\in A\right\}$$

Now if we take the commutant $A''$ of this commutant $A'$, it is obvious that the elements of the original algebra or set $A$ will be inside. Thus, we have an inclusion as follows:
$$A\subset A''$$

(9) The question is now, why $A=A''$ should be equivalent to $A$ being weakly closed, and why should we care about this? These are both good questions, so let us start with the first one. As a first observation, in finite dimensions the bicommutant condition is automatic, because with $A\subset M_N(\mathbb C)$ being as in (2) above, its commutant is:
$$A'=\mathbb C\oplus\ldots\oplus \mathbb C$$

But now, by taking again the commutant, we obtain the original algebra $A$:
$$A''=M_{n_1}(\mathbb C)\oplus\ldots\oplus M_{n_k}(\mathbb C)$$

(10) Moving now to infinite dimensions, the first thought goes into taking the commutant of the basic examples of $C^*$-algebras, $C(X)\subset B(L^2(X))$. But here, up to some mesure theory and functional analysis work, that we will leave as an exercise, we are led to the following conclusion, which proves the bicommutant theorem in this case:
$$C(X)''=L^\infty(X)$$

(11) Summarizing, we have some intuition on the condition $A=A''$ from the statement, and we can also say, based on the above, that the method for proving the bicommutant theorem would be that of establishing the following equality, for any $*$-subalgebra $A\subset B(H)$, with on the right being the closure with respect to the weak topology:
$$A''=\overline{A}^{\,w}$$

(12) Before getting to work, however, we still have a question to be answered, namely, why should we care about all this? I mean, the condition (1) in the statement, weak closedness, looks very nice and mathematical, that would be a good axiom for the von Neumann algebras, so why bothering with commutants, and with the condition (2).

\medskip

(13) In answer, at the elementary level, and with my apologies for calling these damn things ``elementary'', we have seen in chapter 10, when struggling with Tannakian duality, that the bicommutant operation and theorem can be something very useful.

\medskip

(14) In answer too, at the advanced level now, in abstract quantum mechanics the vectors of the Hilbert space $x\in H$ are the states of the system, and the linear self-adjoint operators $T:H\to H$ are the observables, and taking the commutant of a set or algebra of observables is something extremely natural. And this is how von Neumann came upon such things, back in the 1930s, and looking now retrospectively, we can even say that his bicommutant theorem is not only important in the context of quantum mechanics, but even ``makes abstract quantum mechanics properly work''. So, in short, trust me, with the present bicommutant theorem we are into first-class mathematics and physics.

\medskip

(15) Time perhaps for the proof? We recall from (11) that we would like to prove the following equality, for any $*$-algebra of operators $A\subset B(H)$:
$$A''=\overline{A}^{\,w}$$

(16) Let us first prove $\supset$. Since we have $A\subset A''$, we just have to prove that $A''$ is weakly closed. But, assuming $T_i\to T$ weakly, we have indeed:
\begin{eqnarray*}
T_i\in A''
&\implies&ST_i=T_iS,\ \forall S\in A'\\
&\implies&ST=TS,\ \forall S\in A'\\
&\implies&T\in A
\end{eqnarray*}

(17) Let us prove now $\subset$. Here we must establish the following implication:
$$T\in A''\implies T\in\overline{A}^{\,w}$$

For this purpose, we use an amplification trick. Consider indeed the Hilbert space $K$ obtained by summing $n$ times $H$ with itself:
$$K=H\oplus\ldots\oplus H$$

The operators over $K$ can be regarded as being square matrices with entries in $B(H)$, and in particular, we have a representation $\pi:B(H)\to B(K)$, as follows:
$$\pi(T)=\begin{pmatrix}
T\\
&\ddots\\
&&T
\end{pmatrix}$$

(18) The idea will be that of doing the computations in this representation. First, in this representation, the image of our algebra $A\subset B(H)$ is given by:
$$\pi(A)=\left\{\begin{pmatrix}
T\\
&\ddots\\
&&T
\end{pmatrix}\Big|T\in A\right\}$$

We can compute the commutant of this image, exactly as in the usual scalar matrix case, and we obtain the following formula:
$$\pi(A)'=\left\{\begin{pmatrix}
S_{11}&\ldots&S_{1n}\\
\vdots&&\vdots\\
S_{n1}&\ldots&S_{nn}
\end{pmatrix}\Big|S_{ij}\in A'\right\}$$

(19) We conclude from this that, given an operator $T\in A''$ as above, we have:
$$\begin{pmatrix}
T\\
&\ddots\\
&&T
\end{pmatrix}\in\pi(A)''$$

In other words, the conclusion of all this is that we have:
$$T\in A''\implies \pi(T)\in\pi(A)''$$

(20) Now given a vector $x\in K$, consider the orthogonal projection $P\in B(K)$ on the norm closure of the vector space $\pi(A)x\subset K$. Since the subspace $\pi(A)x\subset K$ is invariant under the action of $\pi(A)$, so is its norm closure inside $K$, and we obtain from this:
$$P\in\pi(A)'$$

By combining this with what we found above, we conclude that we have:
$$T\in A''\implies \pi(T)P=P\pi(T)$$

Now since this holds for any vector $x\in K$, we conclude that any operator $T\in A''$ belongs to the weak closure of $A$. Thus, we have $A''\subset\overline{A}^{\,w}$, as desired.  
\end{proof}

Very nice all this, but as you can see, the von Neumann algebras are far more subtle objects than the $C^*$-algebras, and their proper understanding, even at the very basic level, is a far more complicated business than what we quickly did in chapter 12, for the $C^*$-algebras. Welcome to the real quantum, the quantum mechanics one.

\bigskip

Moving ahead, the continuation of the story involves an accumulation of non-trivial results, due to Murray and von Neumann, from the 1930s and 1940s, and then due to Connes, much later, in the 1970s, the conclusions being as follows:

\index{von Neumann algebra}
\index{factor}

\begin{theorem}
The von Neumann algebras are as follows:
\begin{enumerate}
\item In the commutative case, these are the algebras $A=L^\infty(X)$, with $X$ measured space, represented on $H=L^2(X)$, up to a multiplicity.

\item If we write the center as $Z(A)=L^\infty(X)$, then we have a decomposition of type $A=\int_XA_x\,dx$, with the fibers $A_x$ having trivial center,  $Z(A_x)=\mathbb C$.

\item The factors, $Z(A)=\mathbb C$, can be fully classified in terms of ${\rm II}_1$ factors, which are those satisfying $\dim A=\infty$, and having a faithful trace $tr:A\to\mathbb C$.
\end{enumerate}
\end{theorem}

\begin{proof}
This is something quite heavy, the idea being as follows:

\medskip

(1) As already discussed above, it is clear that $L^\infty(X)$ is indeed a von Neumann algebra on $H=L^2(X)$. The converse can be proved as well, by using spectral theory, one way of viewing this being by saying that, given a commutative von Neumann algebra $A\subset B(H)$, its elements $T\in A$ are commuting normal operators, so the Spectral Theorem for such operators applies, and gives $A=L^\infty(X)$, for some measured space $X$.

\medskip

(2) This is von Neumann's reduction theory main result, whose statement is already quite hard to understand, and whose proof uses advanced functional analysis. To be more precise, in finite dimensions this is something that we know well, with the formula $A=\int_XA_x\,dx$ corresponding to our usual direct sum decomposition, namely:
$$A=M_{n_1}(\mathbb C)\oplus\ldots\oplus M_{n_k}(\mathbb C)$$

In infinite dimensions, things are more complicated, but the idea remains the same, namely using (1) for the commutative von Neumann algebra $Z(A)$, as to get a measured space $X$, and then making your way towards a decomposition of type $A=\int_XA_x\,dx$.

\medskip

(3) This is something fairly heavy, due to Murray-von Neumann and Connes, the idea being that the other factors can be basically obtained via crossed product constructions. To be more precise, the various type of factors can be classified as follows:

\medskip

-- Type I. These are the matrix algebras $M_N(\mathbb C)$, called of type ${\rm I}_N$, and their infinite generalization, $B(H)$ with $H$ infinite dimensional, called of type ${\rm I}_\infty$. Although these factors are very interesting and difficult mathematical objects, from the perspective of the general von Neumann algebra classification work, they are dismissed as ``trivial''.

\medskip

-- Type II. These are the infinite dimensional factors having a trace, which is a usual trace $tr:A\to\mathbb C$ in the type ${\rm II}_1$ case, and is something more technical, possibly infinite, in the remaining case, the type ${\rm II}_\infty$ one, with these latter factors being of the form $B(H)\otimes A$, with $A$ being a ${\rm II}_1$ factor, and with $H$ being an infinite dimensional Hilbert space.

\medskip

-- Type III. These are the factors which are infinite dimensional, and do not have a trace $tr:A\to\mathbb C$. Murray and von Neumann struggled a lot with such beasts, with even giving an example being a non-trivial task, but later Connes came and classified them, basically showing that they appear from ${\rm II}_1$ factors, via crossed product constructions.
\end{proof}

So long for basic, or rather advanced but foundational, von Neumann algebra theory. In what follows we will focus on the ${\rm II}_1$ factors, according to the following principle:

\index{factor}

\begin{principle}
The building blocks of the von Neumann algebra theory are the ${\rm II}_1$ factors, which are the von Neumann algebras having the following properties:
\begin{enumerate}
\item They are infinite dimensional, $\dim A=\infty$.

\item They are factors, their center being $Z(A)=\mathbb C$.

\item They have a faithful trace $tr:A\to\mathbb C$.
\end{enumerate}
\end{principle}

But you might perhaps ask, is it even clear that such beasts exist? Good point, and in answer, given a discrete group $\Gamma$, you can talk about its von Neumann algebra, obtained by talking the weak closure of the usual group algebra, or group $C^*$-algebra:
$$L(\Gamma)\subset B(l^2(\Gamma))$$

This algebra is then infinite dimensional when $\Gamma$ is infinite, and also has a trace, given on group elements by $tr(g)=\delta_{g1}$. As for the center, this consists of the functions on $\Gamma$ which are constant on the conjugacy classes, so when $\Gamma$ has infinite conjugacy classes, called ICC property, what we have is a factor. Thus, as a conclusion, when $\Gamma$ is infinite and has the ICC property, its von Neumann algebra $L(\Gamma)$ is a ${\rm II}_1$ factor.

\bigskip

Let us summarize this finding, along with a bit more, as follows:

\index{hyperfinite factor}

\begin{theorem}
We have the following examples of ${\rm II}_1$ factors:
\begin{enumerate}
\item The group von Neumann algebras $L(\Gamma)$, with $\Gamma$ being an infinite group, having the infinite conjugacy class (ICC) property.

\item The Murray-von Neumann hyperfinite factor $R=\overline{\cup_kM_{n_k}(\mathbb C)}^{\,w}$, with the limit being independent on the summands, and on the inclusions between them.

\item With the remark that when $\Gamma$ as above is assumed to be amenable, its associated ${\rm II}_1$ factor $L(\Gamma)$ is the Murray-von Neumann hyperfinite factor $R$.
\end{enumerate}
\end{theorem}

\begin{proof}
Here the first assertion comes from the above discussion, and the rest, regarding the factor $R$, is due to Murray and von Neumann, using standard functional analysis. With the remark however that the notion of hyperfiniteness can be plugged into the general considerations from Theorem 16.4, and with the resulting questions, which are of remarkable difficulty, having been solved only relatively recently, basically by Connes in the 1970s, and with a last contribution by Haagerup in the 1980s, the general idea being that, in the end, everything hyperfinite can be reconstructed from $R$.
\end{proof}

Many other things can be said, along these lines, and if truly interested in theoretical physics, be that quantum mechanics, or statistical mechanics, or other, have a look at all this, von Neumann algebras, this is first-class mathematical technology.

\section*{16b. Subfactor theory}

In view of Principle 16.5, and its quantum mechanics ramifications, it looks reasonable to forget about the Hilbert space $H$, about operators $T\in B(H)$, about other von Neumann algebras and factors $A\subset B(H)$ that might appear, about other mathematics and physics too, why not about your friends, spouse and hobbies too, but please keep teaching some calculus, that is first class mathematics, and focus on the ${\rm II}_1$ factors.

\bigskip

With this idea in mind, we have our objects, the ${\rm II}_1$ factors, but what about morphisms. And here, a natural idea is that of looking at the inclusions of such factors:

\begin{definition}
A subfactor is an inclusion of ${\rm II}_1$ factors $A_0\subset A_1$.
\end{definition}

So, these will be the objects that we will be interested in, in what follows. With the comment that, while quantum mechanics and von Neumann algebras have been around for a while, since the 1920s, and Definition 16.7 is something very natural emerging from this, it took mathematics and physics a lot of time to realize this, with Definition 16.7 basically dating back to the late 1970s, with the beginning of the work of Jones, on it. Moral of the story, sometimes it takes a lot of skill, to come up with simple things. 

\bigskip

Now given a subfactor $A_0\subset A_1$, a first question is that of defining its index, measuring how big $A_1$ is, when compared to $A_0$. But this can be done as follows:

\begin{theorem}
Given a subfactor $A_0\subset A_1$, the number
$$N=\frac{\dim_{A_0}H}{\dim_{A_1}H}$$
is independent of the ambient Hilbert space $H$, and is called index.
\end{theorem}

\begin{proof}
This is something quite standard, the idea being as follows:

\medskip

(1) To start with, given a representation of a ${\rm II}_1$ factor $A\subset B(H)$, we can talk about the corresponding coupling constant, as being a number as follows:
$$\dim_AH\in(0,\infty]$$

To be more precise, we can construct this coupling constant in the following way, with $u:H\to L^2(A)\otimes l^2(\mathbb N)$ being an isometry satisfying $ux=(x\otimes1)u$:
$$\dim_AH=tr(uu^*)$$

(2) Alternatively, we can use the following formula, after proving first that the number on the right is indeed independent of the choice on a nonzero vector $x\in H$:
$$\dim_AH=\frac{tr_A(P_{A'x})}{tr_{A'}(P_{Ax})}$$

This latter formula was in fact the original definition of the coupling constant, by Murray and von Neumann. However, technically speaking, it is better to use (1).

\medskip

(3) Now with this in hand, given a subfactor $A_0\subset A_1$, the fact that the index as defined above is indeed independent of the ambient Hilbert space $H$ comes from the various basic properties of the coupling constant, established by Murray and von Neumann.
\end{proof}

There are many examples of subfactors coming from groups, and every time we obtain the intuitive index. In general now, following Jones \cite{jo1}, let us start with:

\begin{proposition}
Given a subfactor $A_0\subset A_1$, there is a unique linear map
$$E:A_1\to A_0$$
which is positive, unital, trace-preserving and which is such that, for any $a_1,a_2\in A_0$:
$$E(a_1ba_2)=a_1E(b)a_2$$
This map is called conditional expectation from $A_1$ onto $A_0$.
\end{proposition}

\begin{proof}
We make use of the standard representation of the ${\rm II}_1$ factor $A_1$, with respect to its unique trace $tr:A_1\to\mathbb C$, namely:
$$A_1\subset L^2(A_1)$$

If we denote by $\Omega$ the standard cyclic and separating vector of $L^2(A_1)$, we have an identification of vector spaces $A_0\Omega=L^2(A_0)$. Consider now the following projection:
$$e:L^2(A_1)\to L^2(A_0)$$

It follows from definitions that we have an inclusion $e(A_1\Omega)\subset A_0\Omega$. Thus the above projection $e$ induces by restriction a certain linear map, as follows:
$$E:A_1\to A_0$$

This linear map $E$ and the orthogonal projection $e$ are related by:
$$exe=E(x)e$$

But this shows that the linear map $E$ satisfies the various conditions in the statement, namely positivity, unitality, trace preservation and bimodule property. As for the uniqueness assertion, this follows by using the same argument, applied backwards, the idea being that a map $E$ as in the statement must come from a projection $e$.
\end{proof}

We will be interested in what follows in the orthogonal projection $e:L^2(A_1)\to L^2(A_0)$ producing the expectation $E:A_1\to A_0$, rather than in $E$ itself:

\index{Jones projection}

\begin{definition}
Associated to any subfactor $A_0\subset A_1$ is the orthogonal projection
$$e:L^2(A_1)\to L^2(A_0)$$
producing the conditional expectation $E:A_1\to A_0$ via the following formula:
$$exe=E(x)e$$
This projection is called Jones projection for the subfactor $A_0\subset A_1$.
\end{definition}

Quite remarkably, the subfactor $A_0\subset A_1$, as well as its commutant, can be recovered from the knowledge of this projection, in the following way:

\begin{proposition}
Given a subfactor $A_0\subset A_1$, with Jones projection $e$, we have
$$A_0=A_1\cap\{e\}'\quad,\quad 
A_0'=(A_1'\cap\{e\})''$$
as equalities of von Neumann algebras, acting on the space $L^2(A_1)$.
\end{proposition}

\begin{proof}
The above two formulae both follow from $exe=E(x)e$, via some elementary computations, and for details here, we refer to Jones' paper \cite{jo1}.
\end{proof}

We are now ready to formulate a key definition, as follows:

\index{basic construction}

\begin{definition}
Associated to any subfactor $A_0\subset A_1$ is the basic construction
$$A_0\subset_eA_1\subset A_2$$
with $A_2=<A_1,e>$ being the algebra generated by $A_1$ and by the Jones projection
$$e:L^2(A_1)\to L^2(A_0)$$
acting on the Hilbert space $L^2(A_1)$.
\end{definition}

The idea now, following as before Jones \cite{jo1}, will be that the inclusion $A_1\subset A_2$ appears as a kind of ``reflection'' of the original inclusion $A_0\subset A_1$, and also that the basic construction can be iterated, with all this leading to non-trivial results. We first have:

\begin{proposition}
Given a subfactor $A_0\subset A_1$ having finite index, 
$$[A_1:A_0]<\infty$$
the basic construction $A_0\subset_eA_1\subset A_2$ has the following properties:
\begin{enumerate}
\item $A_2=JA_0'J$.

\item $A_2=\overline{A_1+A_1eb}$.

\item $A_2$ is a ${\rm II}_1$ factor.

\item $[A_2:A_1]=[A_1:A_0]$.

\item $eA_2e=A_0e$.

\item $tr(e)=[A_1:A_0]^{-1}$.

\item $tr(xe)=tr(x)[A_1:A_0]^{-1}$, for any $x\in A_1$.
\end{enumerate}
\end{proposition}

\begin{proof}
All this is standard, by using the same type of mathematics as in the proof of Proposition 16.9, and for details here, we refer to Jones' paper \cite{jo1}.
\end{proof}

Let us perform now twice the basic construction, and see what we get. The result here, which is something more technical, at least at the first glance, is as follows:

\begin{proposition}
Associated to $A_0\subset A_1$ is the double basic construction
$$A_0\subset_eA_1\subset_fA_2\subset A_3$$
with $e:L^2(A_1)\to L^2(A_0)$ and $f:L^2(A_2)\to L^2(A_1)$ having the following properties:
$$fef=[A_1:A_0]^{-1}f\quad,\quad 
efe=[A_1:A_0]^{-1}e$$
\end{proposition}

\begin{proof}
We have two formulae to be proved, the idea being as follows:

\medskip

(1) The first formula in the statement is clear, because we have:
$$fef
=E(e)f
=tr(e)f
=[A_1:A_0]^{-1}f$$

(2) Regarding now the second formula, it is enough to check this on the dense subset $(A_1+A_1eA_1)\Omega$. Thus, we must show that for any $x,y,z\in A_1$, we have:
$$efe(x+yez)\Omega=[A_1:A_0]^{-1}e(x+yez)\Omega$$

But this is something which is routine as well. See Jones \cite{jo1}.
\end{proof}

We can in fact perform the basic construction by recurrence, and we obtain:

\index{Jones tower}

\begin{theorem}
Associated to any subfactor $A_0\subset A_1$ is the Jones tower
$$A_0\subset_{e_1}A_1\subset_{e_2}A_2\subset_{e_3}A_3\subset\ldots\ldots$$
with the Jones projections having the following properties:
\begin{enumerate}
\item $e_i^2=e_i^*=e_i$.

\item $e_ie_j=e_je_i$ for $|i-j|\geq2$.

\item $e_ie_{i\pm1}e_i=[A_1:A_0]^{-1}e_i$.

\item $tr(we_{n+1})=[A_1:A_0]^{-1}tr(w)$, for any word $w\in<e_1,\ldots,e_n>$.
\end{enumerate}
\end{theorem}

\begin{proof}
This follows from Proposition 16.13 and Proposition 16.14, because the triple basic construction does not need in fact any further study. See \cite{jo1}.
\end{proof}

The relations found in Theorem 16.15 are in fact well-known, from the standard theory of the Temperley-Lieb algebra. This algebra, discovered by Temperley and Lieb in the context of statistical mechanics \cite{tli}, has a very simple definition, as follows:

\index{Temperley-Lieb algebra}
\index{noncrossing pairings}

\begin{definition}
The Temperley-Lieb algebra of index $N\in[1,\infty)$ is defined as
$$TL_N(k)=span(NC_2(k,k))$$
with product given by vertical concatenation, with the rule
$$\bigcirc=N$$
for the closed circles that might appear when concatenating.
\end{definition}

In other words, the algebra $TL_N(k)$, depending on parameters $k\in\mathbb N$ and $N\in[1,\infty)$, is the linear span of the pairings $\pi\in NC_2(k,k)$. The product operation is obtained by linearity, for the pairings which span $TL_N(k)$ this being the usual vertical concatenation, with the conventions that things go ``from top to bottom'', and that each circle that might appear when concatenating is replaced by a scalar factor, equal to $N$.

\bigskip

In what concerns us, we will just need some elementary results. First, we have:

\index{Jones projection}

\begin{proposition}
The Temperley-Lieb algebra $TL_N(k)$ is generated by the diagrams
$$\varepsilon_1={\ }^\cup_\cap\quad,\quad 
\varepsilon_2=|\!{\ }^\cup_\cap\quad,\quad
\varepsilon_3=||\!{\ }^\cup_\cap\quad,\quad 
\ldots$$
which are all multiples of projections, in the sense that their rescaled versions
$$e_i=N^{-1}\varepsilon_i$$
satisfy the abstract projection relations $e_i^2=e_i^*=e_i$.
\end{proposition}

\begin{proof}
We have two assertions here, the idea being as follows:

\medskip

(1) The fact that the Temperley-Lieb algebra $TL_N(k)$ is indeed generated by the sequence $\varepsilon_1,\varepsilon_2,\ldots$ follows by drawing pictures, and more specifically by decomposing each basis element $\pi\in NC_2(k,k)$ as a product of such elements $\varepsilon_i$.

\medskip

(2) Regarding now the projection assertion, when composing $\varepsilon_i$ with itself we obtain $\varepsilon_i$ itself, times a circle. Thus, according to our multiplication convention, we have:
$$\varepsilon_i^2=N\varepsilon_i$$

Also, when turning upside-down $\varepsilon_i$, we obtain $\varepsilon_i$ itself. Thus, according to our involution convention for the Temperley-Lieb algebra, we have the following formula:
$$\varepsilon_i^*=\varepsilon_i$$

We conclude that the rescalings $e_i=N^{-1}\varepsilon_i$ satisfy $e_i^2=e_i^*=e_i$, as desired.
\end{proof}

As a second result now, making the link with Theorem 16.15, we have:

\begin{proposition}
The standard generators $e_i=N^{-1}\varepsilon_i$ of the Temperley-Lieb algebra $TL_N(k)$ have the following properties, where $tr$ is the trace obtained by closing:
\begin{enumerate}
\item $e_ie_j=e_je_i$ for $|i-j|\geq2$.

\item $e_ie_{i\pm1}e_i=N^{-1}e_i$.

\item $tr(we_{n+1})=N^{-1}tr(w)$, for any word $w\in<e_1,\ldots,e_n>$.
\end{enumerate}
\end{proposition}

\begin{proof}
This follows indeed by doing some elementary computations with diagrams, in the spirit of those performed in the proof of Proposition 16.17.
\end{proof}

With the above results in hand, and still following Jones' paper \cite{jo1}, we can now reformulate Theorem 16.15 into something more conceptual, as follows:

\index{subfactor}
\index{Temperley-Lieb algebra}

\begin{theorem}
Given a subfactor $A_0\subset A_1$, construct its the Jones tower:
$$A_0\subset_{e_1}A_1\subset_{e_2}A_2\subset_{e_3}A_3\subset\ldots\ldots$$
The rescaled sequence of projections $e_1,e_2,e_3,\ldots\in B(H)$ produces then a representation 
$$TL_N\subset B(H)$$
of the Temperley-Lieb algebra of index $N=[A_1:A_0]$.
\end{theorem}

\begin{proof}
We know from Theorem 16.15 that the rescaled sequence of Jones projections $e_1,e_2,e_3,\ldots\in B(H)$ behaves algebrically exactly as the following $TL_N$ diagrams:
$$\varepsilon_1={\ }^\cup_\cap\quad,\quad 
\varepsilon_2=|\!{\ }^\cup_\cap\quad,\quad
\varepsilon_3=||\!{\ }^\cup_\cap\quad,\quad 
\ldots$$

But these diagrams generate $TL_N$, and so we have an embedding $TL_N\subset B(H)$, where $H$ is the Hilbert space where our subfactor $A_0\subset A_1$ lives, as claimed.
\end{proof}

Let us make now the following key observation, also from \cite{jo1}:

\index{relative commutant}
\index{higher commutant}
\index{planar algebra}

\begin{theorem} 
Given a finite index subfactor $A_0\subset A_1$, the graded algebra $P=(P_k)$ formed by the sequence of higher relative commutants
$$P_k=A_0'\cap A_k$$
contains the copy of the Temperley-Lieb algebra constructed above, $TL_N\subset P$. This graded algebra $P=(P_k)$ is called ``planar algebra'' of the subfactor.
\end{theorem}

\begin{proof}
As a first observation, since the Jones projection $e_1:A_1\to A_0$ commutes with $A_0$, we have $e_1\in P_2$. By translation we obtain, for any $k\in\mathbb N$:
$$e_1,\ldots,e_{k-1}\in P_k$$

Thus we have indeed an inclusion of graded algebras $TL_N\subset P$, as claimed.
\end{proof}

As an interesting consequence of the above results, also from \cite{jo1}, we have:

\index{index of subfactors}
\index{small index}
\index{Jones index}
\index{quantized index}

\begin{theorem}
The index of subfactors $A\subset B$ is ``quantized'' in the $[1,4]$ range,
$$N\in\left\{4\cos^2\left(\frac{\pi}{n}\right)\Big|n\geq3\right\}\cup[4,\infty]$$
with the obstruction coming from the existence of the representation $TL_N\subset B(H)$.
\end{theorem}

\begin{proof}
This comes from the basic construction, and more specifically from the combinatorics of the Jones projections $e_1,e_2,e_3,\ldots$, the idea being as folows:

\medskip

(1) In order to best comment on what happens, when iterating the basic construction, let us record the first few values of the numbers in the statement:
$$4\cos^2\left(\frac{\pi}{3}\right)=1\quad,\quad 
4\cos^2\left(\frac{\pi}{4}\right)=2$$
$$4\cos^2\left(\frac{\pi}{5}\right)=\frac{3+\sqrt{5}}{2}\quad,\quad 
4\cos^2\left(\frac{\pi}{6}\right)=3$$
$$\ldots$$

(2) When performing a basic construction, we obtain, by trace manipulations on $e_1$:
$$N\notin(1,2)$$

With a double basic construction, we obtain, by trace manipulations on $<e_1,e_2>$:
$$N\notin\left(2,\frac{3+\sqrt{5}}{2}\right)$$

With a triple basic construction, we obtain, by trace manipulations on $<e_1,e_2,e_3>$:
$$N\notin\left(\frac{3+\sqrt{5}}{2},3\right)$$

Thus, we are led to the conclusion in the statement, by a kind of recurrence, involving a certain family of orthogonal polynomials.

\medskip

(3) In practice now, the most elegant way of proving the result is by using the fundamental fact, explained in Theorem 16.19, that that sequence of Jones projections $e_1,e_2,e_3,\ldots\subset B(H)$ generate a copy of the Temperley-Lieb algebra of index $N$:
$$TL_N\subset B(H)$$

With this result in hand, we must prove that such a representation cannot exist in index $N<4$, unless we are in the following special situation:
$$N=4\cos^2\left(\frac{\pi}{n}\right)$$

But this can be proved by using some suitable trace and positivity manipulations on $TL_N$, as in (2) above. For full details here, we refer to \cite{jo1}.
\end{proof}

So long for the subfactor theory basics. As a continuation of the story, the subfactors of index $N\leq4$ can be shown to be classified by the ADE graphs. See \cite{jo5}.

\section*{16c. Planar algebras}

Quite remarkably, the planar algebra structure of $TL_N$, taken in an intuitive sense, of composing diagrams, extends to a planar algebra structure on $P$. In order to discuss this, let us start with axioms for the planar algebras. Following Jones \cite{jo6}, we have:

\index{planar tangle}
\index{planar algebra}

\begin{definition}
The planar algebras are defined as follows:
\begin{enumerate}
\item We consider rectangles in the plane, with the sides parallel to the coordinate axes, and taken up to planar isotopy, and we call such rectangles boxes.

\item A labeled box is a box with $2n$ marked points on its boundary, $n$ on its upper side, and $n$ on its lower side, for some integer $n\in\mathbb N$.

\item A tangle is labeled box, containing a number of labeled boxes, with all marked points, on the big and small boxes, being connected by noncrossing strings.

\item A planar algebra is a sequence of finite dimensional vector spaces $P=(P_n)$, together with linear maps $P_{n_1}\otimes\ldots\otimes P_{n_k}\to P_n$, one for each tangle, such that the gluing of tangles corresponds to the composition of linear maps.
\end{enumerate}
\end{definition}

In this definition we are using rectangles, but everything being up to isotopy, we could have used instead circles with marked points, as in \cite{jo6}. Our choice for using rectangles comes from the main examples that we have in mind, to be discussed below, where the planar algebra structure is best viewed by using rectangles, as above.

\bigskip

Let us also mention that Definition 16.22 is something quite simplified, based on \cite{jo6}. As explained in \cite{jo6}, in order for subfactors to produce planar algebras and vice versa, there are quite a number of supplementary axioms that must be added, and in view of this, it is perhaps better to start with something stronger than Definition 16.22, as basic axioms. However, as before with rectangles vs circles, our axiomatic choices here are mainly motivated by the concrete examples that we have in mind. More on this later.

\bigskip

As a basic example of a planar algebra, we have the Temperley-Lieb algebra:

\begin{theorem}
The Temperley-Lieb algebra $TL_N$, viewed as graded algebra
$$TL_N=(TL_N(n))_{n\in\mathbb N}$$
is a planar algebra, with the corresponding linear maps associated to the planar tangles
$$TL_N(n_1)\otimes\ldots\otimes TL_N(n_k)\to TL_N(n)$$
appearing by putting the various $TL_N(n_i)$ diagrams into the small boxes of the given tangle, which produces a $TL_N(n)$ diagram.
\end{theorem}

\begin{proof}
This is something trivial, which follows from definitions:

\medskip

(1) Assume indeed that we are given a planar tangle $\pi$, as in Definition 16.22, consisting of a box having $2n$ marked points on its boundary, and containing $k$ small boxes, having respectively $2n_1,\ldots,2n_k$ marked points on their boundaries, and then a total of $n+\Sigma n_i$ noncrossing strings, connecting the various $2n+\Sigma 2n_i$ marked points.

\medskip

(2) We want to associate to this tangle $\pi$ a linear map as follows:
$$T_\pi:TL_N(n_1)\otimes\ldots\otimes TL_N(n_k)\to TL_N(n)$$

For this purpose, by linearity, it is enough to construct elements as follows, for any choice of Temperley-Lieb diagrams $\sigma_i\in TL_N(n_i)$, with $i=1,\ldots,k$:
$$T_\pi(\sigma_1\otimes\ldots\otimes\sigma_k)\in TL_N(n)$$

(3) But constructing such an element is obvious, just by putting the various diagrams $\sigma_i\in TL_N(n_i)$ into the small boxes the given tangle $\pi$. Indeed, this procedure produces a certain diagram in $TL_N(n)$, that we can call $T_\pi(\sigma_1\otimes\ldots\otimes\sigma_k)$, as above.

\medskip

(4) Finally, we have to check that everything is well-defined up to planar isotopy, and that the gluing of tangles corresponds to the composition of linear maps. But both these checks are trivial, coming from the definition of $TL_N$, and we are done.
\end{proof}

As a conclusion to all this, $P=TL_N$ is indeed a planar algebra, but of somewhat ``trivial'' type, with the triviality coming from the fact that, in this case, the elements of $P$ are planar diagrams themselves, and so the planar structure appears trivially.

\bigskip

The Temperley-Lieb planar algebra $TL_N$ is however an important planar algebra, because it is the ``smallest'' one, appearing inside the planar algebra of any subfactor. But more on this later, when talking about planar algebras and subfactors.

\bigskip

Moving ahead now, here is our second basic example of a planar algebra, which is also ``trivial'' in the above sense, with the elements of the planar algebra being planar diagrams themselves, but which appears in a bit more complicated way:

\index{Fuss-Catalan algebra}
\index{colored Temperley-Lieb}

\begin{theorem}
The Fuss-Catalan algebra $FC_{N,M}$, which appears by coloring the Temperley-Lieb diagrams with black/white colors, clockwise, as follows 
$$\circ\bullet\bullet\circ\circ\bullet\bullet\circ\ldots\ldots\ldots\circ\bullet\bullet\circ$$
and keeping those diagrams whose strings connect either $\circ-\circ$ or $\bullet-\bullet$, is a planar algebra, with again the corresponding linear maps associated to the planar tangles
$$FC_{N,M}(n_1)\otimes\ldots\otimes FC_{N,M}(n_k)\to FC_{N,M}(n)$$
appearing by putting the various $FC_{N,M}(n_i)$ diagrams into the small boxes of the given tangle, which produces a $FC_{N,M}(n)$ diagram.
\end{theorem}

\begin{proof}
The proof here is nearly identical to the proof of Theorem 16.23, with the only change appearing at the level of the colors. To be more precise:

\medskip

(1) Forgetting about upper and lower sequences of points, which must be joined by strings, a Temperley-Lieb diagram can be thought of as being a collection of strings, say black strings, which compose in the obvious way, with the rule that the value of the circle, which is now a black circle, is $N$. And it is this obvious composition rule that gives the planar algebra structure, as explained in the proof of Theorem 16.23. 

\medskip

(2) Similarly, forgetting about points, a Fuss-Catalan diagram can be thought of as being a collection of strings, which come now in two colors, black and white. These Fuss-Catalan diagrams compose then in the obvious way, with the rule that the value of the black circle is $N$, and the value of the white circle is $M$. And it is this obvious composition rule that gives the planar algebra structure, as before for $TL_N$.
\end{proof}

Getting back now to generalities, and to Definition 16.22, that of a general planar algebra, we have so far two illustrations for it, which, while both important, are both ``trivial'', with the planar structure simply coming from the fact that, in both these cases, the elements of the planar algebra are planar diagrams themselves.

\bigskip

In general, the planar algebras can be more complicated than this, and we will see some further examples in a moment. However, the idea is very simple, namely ``the elements of a planar algebra are not necessarily diagrams, but they behave like diagrams".

\bigskip

In relation now with subfactors, the result, which extends Theorem 16.20, and which was found by Jones in \cite{jo6}, almost 20 years after \cite{jo1}, is as follows:

\index{higher commutant}
\index{planar algebra}

\begin{theorem} 
Given a subfactor $A_0\subset A_1$, the collection $P=(P_n)$ of linear spaces 
$$P_n=A_0'\cap A_n$$
has a planar algebra structure, extending the planar algebra structure of $TL_N$.
\end{theorem}

\begin{proof}
We know from Theorem 16.20 that we have an inclusion as follows, coming from the basic construction, and with $TL_N$ itself being a planar algebra:
$$TL_N\subset P$$

Thus, the whole point is that of proving that the trivial planar algebra structure of $TL_N$ extends into a planar algebra structure of $P$. But this can be done via a long algebraic study, and for the full computation here, we refer to Jones' paper \cite{jo6}.
\end{proof}

As a first illustration for the above result, we have:

\index{Temperley-Lieb}
\index{Fuss-Catalan algebra}

\begin{theorem}
We have the following universality results:
\begin{enumerate}
\item The Temperley-Lieb algebra $TL_N$ appears inside the planar algebra of any subfactor $A_0\subset A_1$ having index $N$.

\item The Fuss-Catalan algebra $FC_{N,M}$ appears inside the planar algebra of any subfactor $A_0\subset A_1$, in the presence of an intermediate subfactor $A_0\subset B\subset A_1$.
\end{enumerate}
\end{theorem}

\begin{proof}
Here the first assertion is something that we already know, from Theorem 16.20, and the second assertion, due to Bisch and Jones, is something quite standard as well, by carefully working out the basic construction for $A_0\subset A_1$, in the presence of an intermediate subfactor $A_0\subset B\subset A_1$. For furher details here, we refer to \cite{jo6}.
\end{proof}

The above results raise the question on whether any planar algebra produces a subfactor. The answer here is yes, but with many subtleties, as follows:

\index{amenable subfactor}
\index{Popa theorem}

\begin{theorem}
We have the following results:
\begin{enumerate}
\item Any planar algebra with positivity produces a subfactor.

\item In particular, we have $TL$ and $FC$ type subfactors.

\item In the amenable case, and with $A_1=R$, the correspondence is bijective.

\item In general, we must take $A_1=L(F_\infty)$, and we do not have bijectivity.

\item The axiomatization of $P$, in the case $A_1=R$, is not known.
\end{enumerate}
\end{theorem}

\begin{proof}
All this is quite heavy, basically coming from the work of Popa in the 90s, using heavy functional analysis, the idea being as follows:

\medskip

(1) As already mentioned in the comments after Definition 16.22, our planar algebra axioms here are something quite simplified, based on \cite{jo6}. However, when getting back to Theorem 16.25, the conclusion is that the subfactor planar algebras there satisfy a number of supplementary ``positivity'' conditions, basically coming from the positivity of the ${\rm II}_1$ factor trace. And the point is that, with these positivity conditions axiomatized, we reach to something which is equivalent to Popa's axiomatization of the lattice of higher relative commutants $A_i'\cap A_j$ of the finite index subfactors, obtained in the 90s via heavy functional analysis. For the full story here, and details, we refer to Jones' paper \cite{jo6}.

\medskip

(2) The existence of the $TL_N$ subfactors, also known as ``$A_\infty$ subfactors'', is something which was known for some time, since some early work of Popa on the subject. As for the existence of the $FC_{N,M}$ subfactors, this can be shown by using the intermediate subfactor picture, $A_0\subset B\subset A_1$, by composing two $A_\infty$ subfactors of suitable indices, $A_0\subset B$ and $B\subset A_1$. For the full story here, we refer as before to Jones \cite{jo6}.

\medskip

(3) This is something fairly heavy, as it is always the case with operator algebra results regarding hyperfiniteness and amenability, due to Popa. For the story here, see \cite{jo6}.

\medskip

(4) This is something a bit more recent, obtained by further building on the above-mentioned constructions of Popa. Again, we refer here to \cite{jo6} and related work.

\medskip

(5) This is the big open question in subfactors. The story here goes back to Jones' original  paper \cite{jo1}, which contains at the end the question, due to Connes, of finding the possible values of the index for the irreducible subfactors of $R$. This question, which certainly looks much easier than (5) in the statement, is in fact still open, now 40 years after its formulation, and with no one having any valuable idea in dealing with it.
\end{proof}

\section*{16d. Group actions}

Getting back to groups, all this machinery is interesting for us. We will need the construction of the tensor and spin planar algebras $\mathcal T_N,\mathcal S_N$. Let us start with:

\index{tensor planar algebra}

\begin{definition}
The tensor planar algebra $\mathcal T_N$ is the sequence of vector spaces 
$$P_k=M_N(\mathbb C)^{\otimes k}$$
with the multilinear maps $T_\pi:P_{k_1}\otimes\ldots\otimes P_{k_r}\to P_k$
being given by the formula
$$T_\pi(e_{i_1}\otimes\ldots\otimes e_{i_r})=\sum_j\delta_\pi(i_1,\ldots,i_r:j)e_j$$
with the Kronecker symbols $\delta_\pi$ being $1$ if the indices fit, and being $0$ otherwise.
\end{definition}

In other words, we put the indices of the basic tensors on the marked points of the small boxes, in the obvious way, and the coefficients of the output tensor are then given by Kronecker symbols, exactly as in the easy quantum group case.

\bigskip

The fact that we have indeed a planar algebra, in the sense that the gluing of tangles corresponds to the composition of linear maps, as required by Definition 16.22, is something elementary, in the same spirit as the verification of the functoriality properties of the correspondence $\pi\to T_\pi$, from easiness, and we refer here to Jones \cite{jo6}. 

\bigskip

Let us discuss now a second planar algebra of the same type, which is important as well for various reasons, namely the spin planar algebra $\mathcal S_N$. This planar algebra appears somehow as the ``square root'' of the tensor planar algebra $\mathcal T_N$. Let us start with:

\begin{definition}
We write the standard basis of $(\mathbb C^N)^{\otimes k}$ in $2\times k$ matrix form,
$$e_{i_1\ldots i_k}=
\begin{pmatrix}i_1 & i_1 &i_2&i_2&i_3&\ldots&\ldots\\
i_k&i_k&i_{k-1}&\ldots&\ldots&\ldots&\ldots 
\end{pmatrix}$$
by duplicating the indices, and then writing them clockwise, starting from top left.
\end{definition}

Now with this convention in hand for the tensors, we can formulate the construction of the spin planar algebra $\mathcal S_N$, also from \cite{jo6}, as follows:

\index{spin planar algebra}

\begin{definition}
The spin planar algebra $\mathcal S_N$ is the sequence of vector spaces 
$$P_k=(\mathbb C^N)^{\otimes k}$$
written as above, with the multiplinear maps $T_\pi:P_{k_1}\otimes\ldots\otimes P_{k_r}\to P_k$ 
being given by
$$T_\pi(e_{i_1}\otimes\ldots\otimes e_{i_r})=\sum_j\delta_\pi(i_1,\ldots,i_r:j)e_j$$
with the Kronecker symbols $\delta_\pi$ being $1$ if the indices fit, and being $0$ otherwise.
\end{definition}

Here are some illustrating examples for the spin planar algebra calculus:

\medskip

(1) The identity $1_k$ is the $(k,k)$-tangle having vertical strings only. The solutions of $\delta_{1_k}(x,y)=1$ being the pairs of the form $(x,x)$, this tangle $1_k$ acts by the identity:
$$1_k\begin{pmatrix}j_1 & \ldots & j_k\\ i_1 & \ldots & i_k\end{pmatrix}=\begin{pmatrix}j_1 & \ldots & j_k\\ i_1 & \ldots & i_k\end{pmatrix}$$

(2) The multiplication $M_k$ is the $(k,k,k)$-tangle having 2 input boxes, one on top of the other, and vertical strings only. It acts in the following way:
$$M_k\left( 
\begin{pmatrix}j_1 & \ldots & j_k\\ i_1 & \ldots & i_k\end{pmatrix}
\otimes\begin{pmatrix}l_1 & \ldots & l_k\\ m_1 & \ldots & m_k\end{pmatrix}
\right)=
\delta_{j_1m_1}\ldots \delta_{j_km_k}
\begin{pmatrix}l_1 & \ldots & l_k\\ i_1 & \ldots & i_k\end{pmatrix}$$

(3) The inclusion $I_k$ is the $(k,k+1)$-tangle which looks like $1_k$, but has one more vertical string, at right of the input box. Given $x$, the solutions of $\delta_{I_k}(x,y)=1$ are the elements $y$ obtained from $x$ by adding to the right a vector of the form $(^l_l)$, and so:
$${I_k}\begin{pmatrix}j_1 & \ldots & j_k\\ i_1 & \ldots & i_k\end{pmatrix}=
\sum_l\begin{pmatrix}j_1 & \ldots & j_k& l\\ i_1 & \ldots & i_k& l\end{pmatrix}$$

(4) The expectation $U_k$ is the $(k+1,k)$-tangle which looks like $1_k$, but has one more string, connecting the extra 2 input points, both at right of the input box:
$$U_k
\begin{pmatrix}j_1 & \ldots &j_k& j_{k+1}\\ i_1 & \ldots &i_k& i_{k+1}\end{pmatrix}=
\delta_{i_{k+1}j_{k+1}}
\begin{pmatrix}j_1 & \ldots & j_k\\ i_1 & \ldots & i_k\end{pmatrix}$$

(5) The Jones projection $E_k$ is a $(0,k+2)$-tangle, having no input box. There are $k$ vertical strings joining the first $k$ upper points to the first $k$ lower points, counting from left to right. The remaining upper 2 points are connected by a semicircle, and the remaining lower 2 points are also connected by a semicircle. We have:
$$E_k(1)=\sum_{ijl}\begin{pmatrix}i_1 & \ldots &i_k&j&j\\ i_1 & \ldots &i_k&l&l\end{pmatrix}$$

The elements $e_k=N^{-1}E_k(1)$ are then projections, and define a representation of the infinite Temperley-Lieb algebra of index $N$ inside the inductive limit algebra $\mathcal S_N$.

\medskip

(6) The rotation $R_k$ is the $(k,k)$-tangle which looks like $1_k$, but the first 2 input points are connected to the last 2 output points, and the same happens at right:
$$R_k=\begin{matrix}
\hskip 0.3mm\Cap \ |\ |\ |\ |\hskip -0.5mm |\cr
|\hskip -0.5mm |\hskip 10.3mm |\hskip -0.5mm |\cr
\hskip -0.3mm|\hskip -0.5mm |\ |\ |\ |\ \hskip -0.1mm\Cup
\end{matrix}$$

The action of $R_k$ on the standard basis is by rotation of the indices, as follows:
$$R_k(e_{i_1i_2\ldots i_k})=e_{i_2\ldots i_ki_1}$$

There are many other interesting examples of $k$-tangles, but in view of our present purposes, we can actually stop here, due to the following fact:

\begin{theorem}
The multiplications, inclusions, expectations, Jones projections and rotations generate the set of all tangles, via the gluing operation.
\end{theorem}

\begin{proof}
This is something well-known and elementary, obtained by ``chopping'' the various planar tangles into small pieces, as in the above list. See \cite{jo6}.
\end{proof}

Finally, in order for our discussion to be complete, we must talk as well about the $*$-structure of the spin planar algebra. This is constructed as follows:
$$\begin{pmatrix}j_1 & \ldots & j_k\\ i_1 & \ldots & i_k\end{pmatrix}^*
=\begin{pmatrix}i_1 & \ldots & i_k\\ j_1 & \ldots & j_k\end{pmatrix}$$

As before, we refer to Jones' paper \cite{jo6} for more on all this. Getting back now to groups, and more specifically to permutations, we have the following result:

\begin{theorem}
Given $G\subset S_N^+$, consider the tensor powers of the associated coaction map on $C(X)$, where $X=\{1,\ldots,N\}$, which are the folowing linear maps:
$$\Phi^k:C(X^k)\to C(X^k)\otimes C(G)$$
$$e_{i_1\ldots i_k}\to\sum_{j_1\ldots j_k}e_{j_1\ldots j_k}\otimes u_{j_1i_1}\ldots u_{j_ki_k}$$
The fixed point spaces of these coactions, which are by definition the spaces
$$P_k=\left\{ x\in C(X^k)\Big|\Phi^k(x)=1\otimes x\right\}$$
are given by $P_k=Fix(u^{\otimes k})$, and form a subalgebra of the spin planar algebra $\mathcal S_N$.
\end{theorem}

\begin{proof}
Since the map $\Phi$ is a coaction, its tensor powers $\Phi^k$ are coactions too, and at the level of fixed point algebras we have the following formula:
$$P_k=Fix(u^{\otimes k})$$

In order to prove now the planar algebra assertion, we will use Theorem 16.31. Consider the rotation $R_k$. Rotating, then applying $\Phi^k$, and rotating backwards by $R_k^{-1}$ is the same as applying $\Phi^k$, then rotating each $k$-fold product of coefficients of $\Phi$. Thus the elements obtained by rotating, then applying $\Phi^k$, or by applying $\Phi^k$, then rotating, differ by a sum of Dirac masses tensored with commutators in $A=C(G)$:
$$\Phi^kR_k(x)-(R_k\otimes id)\Phi^k(x)\in C(X^k)\otimes [A,A]$$

Now let $\int_A$ be the Haar functional of $A$, and consider the conditional expectation onto the fixed point algebra $P_k$, which is given by the following formula:
$$\phi_k=\left(id\otimes\int_A\right)\Phi^k$$

Since $\int_A$ is a trace, it vanishes on commutators. Thus $R_k$ commutes with $\phi_k$:
$$\phi_kR_k=R_k\phi_k$$

The commutation relation $\phi_kT=T\phi_l$ holds in fact for any $(l,k)$-tangle $T$. These tangles are called annular, and the proof is by verification on generators of the annular category. In particular we obtain, for any annular tangle $T$:
$$\phi_kT\phi_l=T\phi_l$$

We conclude from this that the annular category is contained in the suboperad $\mathcal P'\subset\mathcal P$ of the planar operad consisting of tangles $T$ satisfying the following condition, where $\phi =(\phi_k)$, and where $i(.)$ is the number of input boxes:
$$\phi T\phi^{\otimes i(T)}=T\phi^{\otimes i(T)}$$

On the other hand the multiplicativity of $\Phi^k$ gives $M_k\in\mathcal P'$. Now since the planar operad $\mathcal P$ is generated by multiplications and annular tangles, it follows that we have $\mathcal P'=P$. Thus for any tangle $T$ the corresponding multilinear map between spaces $P_k(X)$ restricts to a multilinear map between spaces $P_k$. In other words, the action of the planar operad $\mathcal P$ restricts to $P$, and makes it a subalgebra of $\mathcal S_N$, as claimed.
\end{proof}

As a second result now, further building on Theorem 16.32, we have:

\begin{theorem}
We have a bijection between quantum permutation groups and subalgebras of the spin planar algebra,
$$(G\subset S_N^+)\quad\longleftrightarrow\quad (Q\subset\mathcal S_N)$$
given in one sense by the construction in Theorem 16.32, and in the other sense by a suitable modification of Tannakian duality.
\end{theorem}

\begin{proof}
The idea is that this will follow by applying Tannakian duality to the annular category over $Q$. Let $n,m$ be positive integers. To any element $T_{n+m}\in Q_{n+m}$ we associate a linear map $L_{nm}(T_{n+m}):P_n(X)\to P_m(X)$ in the following way:
$$L_{nm}\left(\begin{matrix}|\ |\ |\\ T_{n+m}\\ |\ |\ |\end{matrix}\right):
\left(\begin{matrix}|\\ a_n\\ |\end{matrix}\right)
\to \left(\begin{matrix}
\hskip 1.5mm |\hskip 3.0mm |\hskip 3.0mm \cap\\
\ \ T_{n+m}\hskip 0.0mm  |\\
\hskip 1.9mm |\hskip 1.2mm |\hskip 3.2mm |\hskip2.2mm |\\
a_n|\hskip 3.2mm |\hskip 2.2mm |\\
\hskip 2.1mm\cup \hskip3.5mm |\hskip 2.2mm |
\end{matrix}\right)$$

That is, we consider the planar $(n,n+m,m)$-tangle having an small input $n$-box, a big input $n+m$-box and an output $m$-box, with strings as on the picture of the right. This defines a certain multilinear map, as follows:
$$P_n(X)\otimes P_{n+m}(X)\to P_m(X)$$

If we put the element $T_{n+m}$ in the big input box, we obtain in this way a certain linear map $P_n(X)\to P_m(X)$, that we call $L_{nm}$. With this convention, let us set:
$$Q_{nm}=\left\{ L_{nm}(T_{n+m}):P_n(X)\to P_m(X)\Big| T_{n+m}\in Q_{n+m}\right\}$$

These spaces form a Tannakian category, so we obtain a certain Woronowicz algebra $(A,u)$, such that the following equalities hold, for any $m,n$:
$$Hom(u^{\otimes m},u^{\otimes n})=Q_{mn}$$

We prove that $u$ is a magic unitary. We have $Hom(1,u^{\otimes 2})=Q_{02}=Q_2$, so the unit of $Q_2$ must be a fixed vector of $u^{\otimes 2}$. But $u^{\otimes 2}$ acts on the unit of $Q_2$ as follows:
\begin{eqnarray*}
u^{\otimes 2}(1)
&=&u^{\otimes 2}\left( \sum_i \begin{pmatrix}i&i\\ i&i\end{pmatrix}\right)\\
&=&\sum_{ikl}\begin{pmatrix}k&k\\ l&l\end{pmatrix}\otimes u_{ki}u_{li}\\
&=&\sum_{kl}\begin{pmatrix}k&k\\ l&l\end{pmatrix}\otimes (uu^t)_{kl}
\end{eqnarray*}

From $u^{\otimes 2}(1)=1\otimes 1$ ve get that $uu^t$ is the identity matrix. Together with the unitarity of $u$, this gives the following formulae:
$$u^t=u^*=u^{-1}$$

Consider the Jones projection $E_1\in Q_3$. After isotoping, $L_{21}(E_1)$ looks as follows:
$$L_{21}\left( \Bigl| \begin{matrix}\cup\\\cap\end{matrix}\right) :
\begin{pmatrix} \,|\ |\\ {\ }^i_j{\ }^i_j\\ \,|\ |\end{pmatrix}\,\to\,
\begin{pmatrix}\hskip -5.8mm |\\ {\ }^i_j{\ }^i_j\supset\\ \hskip -5.8mm |\end{pmatrix}
=\,\delta_{ij}\begin{pmatrix}\,|\\ {\ }^i_i\\ \,|\end{pmatrix}$$

In other words, the linear map $M=L_{21}(E_1)$ is the multiplication $\delta_i\otimes\delta_j\to\delta_{ij}\delta_i$:
$$M\begin{pmatrix}i&i\\ j&j\end{pmatrix}
=\delta_{ij}\begin{pmatrix}i\\ i\end{pmatrix}$$

In order to finish, consider the following element of $C(X)\otimes A$:
$$(M\otimes id)u^{\otimes 2}\left(\begin{pmatrix}i&i\\ j&j\end{pmatrix}\otimes 1\right)
=\sum_k\begin{pmatrix}k\\ k\end{pmatrix}\delta_k\otimes u_{ki}u_{kj}$$

Since $M\in Q_{21}=Hom(u^{\otimes 2},u)$, this equals the following element of $C(X)\otimes A$:
$$u(M\otimes id)\left(\begin{pmatrix}i&i\\ j&j\end{pmatrix}\otimes 1\right)
=\sum_k\begin{pmatrix}k\\ k\end{pmatrix}\delta_k\otimes\delta_{ij}u_{ki}$$

Thus we have $u_{ki}u_{kj}=\delta_{ij}u_{ki}$ for any $i,j,k$, which shows that $u$ is a magic unitary. Now if $P$ is the planar algebra associated to $u$, we have $Hom(1,v^{\otimes n})=P_n=Q_n$, as desired. As for the uniqueness, this is clear from the Peter-Weyl theory.
\end{proof}

All the above might seem a bit technical, but is worth learning, and for good reason, because it is extremely powerful. As an example of application, if you agree with the bijection $G\leftrightarrow Q$ in Theorem 16.33, then $G=S_N^+$ itself, which is the biggest object on the left, must correspond to the smallest object on the right, namely $Q=TL_N$. 

\bigskip

As a useful complement to Theorem 16.33, we have as well the following result:

\index{quantum symmetry group}
\index{2-box}

\begin{theorem}
The planar algebra associated to a quantum symmetry group
$$G^+(X)=C(S_N^+)\Big/\Big<du=ud\Big>$$
is the planar algebra generated by $d$, viewed as $2$-box in the spin planar algebra $\mathcal S_N$.
\end{theorem}

\begin{proof}
We recall from the above that any quantum permutation group $G\subset S_N^+$ produces a subalgebra $P\subset\mathcal S_N$ of the spin planar algebra, given by:
$$P_k=Fix(u^{\otimes k})$$

In our case, the idea is that $G=G^+(X)$ comes via the relation $d\in End(u)$, but we can view this relation, via Frobenius duality, as a relation of the following type:
$$\xi_d\in Fix(u^{\otimes 2})$$

Indeed, let us view the adjacency matrix $d\in M_N(0,1)$ as a 2-box in $\mathcal S_N$, by using the canonical identification between $M_N(\mathbb C)$ and the algebra of 2-boxes $\mathcal S_N(2)$:
$$(d_{ij})\leftrightarrow \sum_{ij} d_{ij}\begin{pmatrix}i&i\\ j&j\end{pmatrix}$$

Let $P$ be the planar algebra associated to $G^+(X)$ and let $Q$ be the planar algebra generated by $d$. The action of $u^{\otimes 2}$ on $d$ viewed as a 2-box is given by:
$$u^{\otimes 2}\left(\sum_{ij} d_{ij}\begin{pmatrix}i&i\\ j&j\end{pmatrix}\right)
=\sum_{ijkl} d_{ij}\begin{pmatrix}k&k\\ l&l\end{pmatrix}\otimes u_{ki}u_{lj}
=\sum_{kl}\begin{pmatrix}k&k\\ l&l\end{pmatrix}\otimes (udu^t)_{kl}$$

Since $v$ is a magic unitary commuting with $d$ we have:
$$udu^t=duu^t=d$$

But this means that $d$, viewed as a 2-box, is in the algebra $P_2$ of fixed points of $u^{\otimes 2}$. Thus $Q\subset P$. As for $P\subset Q$, this follows from the duality found above.
\end{proof}

And with this, end of our basic presentation of planar algebras. As a conclusion to this, planar algebras provide a powerful Tannakian tool for the study of the permutation groups $G$, finite or quantum, and in case you are stuck with difficult questions about $G$, and its actions on graphs $X$, try the planar algebra technology, this might help.

\bigskip

In relation with subfactors, many things can be said too, the idea being that the subfactors associated to the permutation groups $G$, finite or quantum, appear via some standard fixed point procedures, which are something useful and interesting too.

\bigskip

Importantly, all this generalizes to more general actions, such as the adjoint actions of the finite or compact groups $G\subset U_N$, or more generally quantum groups $G\subset U_N^+$, again with fixed point subfactors, and with the corresponding planar algebras being given by $P_k=End(u^{\otimes k})$, as in Principle 16.2, and with a Tannakian correspondence too. 

\bigskip

Finally, and we would like to emphasize this, the subfactor and planar algebra formalism is in fact far more general than this, allowing for indices $N\notin\mathbb N$. Thus, if you ever come upon some strange group-type object, with $|G|\notin N$, say via some theoretical physics computations, try subfactors and planar algebras, these might be the answer.

\section*{16e. Exercises}

Congratulations for having read this book, and no exercises for this final chapter. However, if you want to learn more about the finite groups, and their various generalizations, be them discrete, compact, locally compact, or of quantum type, there are many books available, all containing very interesting mathematics. We have referenced some of them below, and in the hope that you will go this way, advanced group theory.

\baselineskip=14pt

\printindex

\end{document}